\documentclass[10pt]{book}
\usepackage{amsfonts}
\usepackage{amssymb}
\usepackage{amsmath}
\usepackage{tensor}
\usepackage{makeidx}
\usepackage{wasysym}
\usepackage{stmaryrd}
\usepackage{marvosym}
\usepackage{scalerel}
\usepackage{graphicx}
\usepackage[all]{xy}
\usepackage{float}
\usepackage{xparse}

\newtheorem{theorem}{Theorem}[section]
\newtheorem{corollary}[theorem]{Corollary}
\newtheorem{conjecture}[theorem]{Conjecture}
\newtheorem{definition}[theorem]{Definition}
\newtheorem{construction}[theorem]{Construction}
\newtheorem{notation}[theorem]{Notation}
\newtheorem{example}[theorem]{Example}
\newtheorem{proposition}[theorem]{Proposition}
\newtheorem{remark}[theorem]{Remark}
\newtheorem{lemma}[theorem]{Lemma}

\newtheorem{question}[theorem]{Question}
\numberwithin{equation}{chapter}
\newenvironment{proof}[1][Proof]{\textbf{#1.} }{\ \rule{0.5em}{0.5em}}

\makeatletter
\newdimen\rh@wd
\newdimen\rh@hta
\newdimen\rh@htb
\newbox\rh@box
\def\rh@measure#1{\setbox\rh@box=\hbox{$#1$}\rh@wd=\wd\rh@box \rh@hta=\ht\rh@box}

\def\widecheck#1{\rh@measure{#1}%
  \setbox\rh@box=\hbox{$\widehat{\vrule height \rh@hta width\z@ \kern\rh@wd}$}%
  \rh@htb=\ht\rh@box \advance\rh@htb\rh@hta \advance\rh@htb\p@
  \ooalign{$\vrule height \ht\rh@box width\z@ #1$\cr
           \raise\rh@htb\hbox{\scalebox{1}[-1]{\box\rh@box}}\cr}}
\makeatother

\NewDocumentCommand\argcompl{m}{#1}

\NewDocumentCommand\OAOS{d// d-- d[] d() d++}{
  \IfNoValueTF{#1}{
	\IfNoValueTF{#2}{\IfNoValueTF{#3}{
	\IfNoValueTF{#4}{\IfNoValueTF{#5}{\mathfrak{S}}{\mathfrak{S}\left[#5\right]}}{\IfNoValueTF{#5}{OCCHIO!1}{\mathfrak{S}\left[#4,#5\right]}}}
	{\IfNoValueTF{#4}{OCCHIO!2}{\IfNoValueTF{#5}{OCCHIO!3}{\mathfrak{S}\left[#3,#4,#5\right]}}}}
	{\IfNoValueTF{#3}{
	\IfNoValueTF{#4}{\IfNoValueTF{#5}{\mathfrak{S}}{\mathfrak{S}_{#5}}}{\IfNoValueTF{#5}{OCCHIO!4}{\mathfrak{S}_{#4   ,#5}}}}
	{\IfNoValueTF{#4}{OCCHIO!5}{\IfNoValueTF{#5}{OCCHIO!6}{\mathfrak{S}_{#3,#4,#5}}}}}
	}
	{ \IfNoValueTF{#2}{OCCHIO!7}
	{\IfNoValueTF{#3}{
	\IfNoValueTF{#4}{\IfNoValueTF{#5}{OCCHIO!11}{\mathfrak{S}^{#1}_{#5}}}{\IfNoValueTF{#5}{OCCHIO!8}{\mathfrak{S}^{#1}_{#4,#5}}}}
	{\IfNoValueTF{#4}{OCCHIO!9}{\IfNoValueTF{#5}{OCCHIO!10}{\mathfrak{S}^{#1}_{#3,#4,#5}}}}}}
}
\NewDocumentCommand\OAAS{d// d-- d[] d() d++}{
  \IfNoValueTF{#1}{
	\IfNoValueTF{#2}{\IfNoValueTF{#3}{
	\IfNoValueTF{#4}{\IfNoValueTF{#5}{\varsigma}{\varsigma\left[#5\right]}}{\IfNoValueTF{#5}{OCCHIO!1}{\varsigma\left[#4,#5\right]}}}
	{\IfNoValueTF{#4}{OCCHIO!2}{\IfNoValueTF{#5}{OCCHIO!3}{\varsigma\left[#3,#4,#5\right]}}}}
	{\IfNoValueTF{#3}{
	\IfNoValueTF{#4}{\IfNoValueTF{#5}{\varsigma}{\varsigma_{#5}}}{\IfNoValueTF{#5}{OCCHIO!4}{\varsigma_{#4   ,#5}}}}
	{\IfNoValueTF{#4}{OCCHIO!5}{\IfNoValueTF{#5}{OCCHIO!6}{\varsigma_{#3,#4,#5}}}}}
	}
	{ \IfNoValueTF{#2}{OCCHIO!7}
	{\IfNoValueTF{#3}{
	\IfNoValueTF{#4}{\IfNoValueTF{#5}{OCCHIO!11}{\varsigma^{#1}_{#5}}}{\IfNoValueTF{#5}{OCCHIO!8}{\varsigma^{#1}_{#4,#5}}}}
	{\IfNoValueTF{#4}{OCCHIO!9}{\IfNoValueTF{#5}{OCCHIO!10}{\varsigma^{#1}_{#3,#4,#5}}}}}}
}
\NewDocumentCommand\OAASdue{d// d-- d<< d>> d() d++}{
\IfNoValueTF{#1}{
	\IfNoValueTF{#2}{\IfNoValueTF{#3}{
	\IfNoValueTF{#4}{\IfNoValueTF{#5}{\IfNoValueTF{#6}{\varsigma}{\varsigma\left[#6\right]}}{\IfNoValueTF{#6}{OCCHIO!1}{\varsigma\left[#5,#6\right]}}}{OCCHIO!2}}{\IfNoValueTF{#4}{OCCHIO!3}{\IfNoValueTF{#5}{OCCHIO!4}{\IfNoValueTF{#6}{OCCHIO!5}{\varsigma\left[#3,#4,#5,#6\right]}}}}}
	{\IfNoValueTF{#3}{
	\IfNoValueTF{#4}{\IfNoValueTF{#5}{\IfNoValueTF{#6}{\varsigma}{\varsigma_{#6}}}{\IfNoValueTF{#6}{OCCHIO!6}{\varsigma_{#5,#6}}}}{OCCHIO!7}}{\IfNoValueTF{#4}{OCCHIO!8}{\IfNoValueTF{#5}{OCCHIO!9}{\IfNoValueTF{#6}{OCCHIO!10}{\varsigma_{#3,#4,#5,#6}}}}}}}		
	{\IfNoValueTF{#2}{OCCHIO!11}
{\IfNoValueTF{#3}{
	\IfNoValueTF{#4}{\IfNoValueTF{#5}{\IfNoValueTF{#6}{\varsigma^{#1}}{\varsigma^{#1}_{#6}}}{\IfNoValueTF{#6}{OCCHIO!12}{\varsigma^{#1}_{#5,#6}}}}{OCCHIO!13}}{\IfNoValueTF{#4}{OCCHIO!14}{\IfNoValueTF{#5}{OCCHIO!15}{\IfNoValueTF{#6}{OCCHIO!16}{\varsigma^{#1}_{#3,#4,#5,#6}}}}}}}	
}
\NewDocumentCommand\grm{d[]}{
\IfNoValueTF{#1}{VARIA CON IL FUNTORE CONSIDERATO}{ \langle  #1 \rangle}
}

\newcommand{\Catuno}{A}
\newcommand{\objuno}{a}
\newcommand{\arruno}{f}
\newcommand{\Catdue}{B}
\newcommand{\objdue}{b}
\newcommand{\arrdue}{g}
\newcommand{\Cattre}{C}
\newcommand{\objtre}{c}

\newcommand{\Catquattro}{D}
\newcommand{\objquattro}{d}

\newcommand{\Catcinque}{E}

\newcommand{\Catsei}{X}
\newcommand{\objsei}{x}

\newcommand{\Catsette}{Y}
\newcommand{\objsette}{y}
\newcommand{\arrsette}{\upsilon}
\newcommand{\Catindex}{I}
\newcommand{\objindex}{i}

\NewDocumentCommand\terobj{d++}{ % terminal object
  \IfNoValueTF{#1}{\infty}{\infty_{#1}}
}
\NewDocumentCommand\idobj{d++}{ %  identity arrow
  \IfNoValueTF{#1}{1}{1_{#1}}
}                      
\NewDocumentCommand\genfunc{d[]}{
\IfNoValueTF{#1}{\mathbf{F}}{\mathbf{F}_{#1}}
}
\NewDocumentCommand\genfuncuno{d[]}{
\IfNoValueTF{#1}{\mathbf{G}}{\mathbf{G}_{#1}}
}
\NewDocumentCommand\genfuncdue{d[]}{
\IfNoValueTF{#1}{\mathbf{H}}{\mathbf{H}_{#1}}
}
\newcommand{\fundim}{\mathsf{U}}
\NewDocumentCommand\homF{d++ d// d<< d>>}{
  \IfNoValueTF{#1}{\IfNoValueTF{#2}{\IfNoValueTF{#3}{\IfNoValueTF{#4}{\mathsf{hom}}{\mathsf{hom}^{#4}}}
  {\IfNoValueTF{#4}{{^{#3}\mathsf{hom}}}{OCCHIO!1}}}
	{\IfNoValueTF{#3}{\IfNoValueTF{#4}{OCCHIO!2}{H#2^{#4}}}
	{\IfNoValueTF{#4}{{^{#3}H#2}}{OCCHIO!3}}}}
  {\IfNoValueTF{#2}{\IfNoValueTF{#3}{\IfNoValueTF{#4}{\mathsf{hom}_{#1}}{{\mathsf{hom}_{#1}}^{#4}}}
	{\IfNoValueTF{#4}{{^{#3}\mathsf{hom}_{#1}}}{OCCHIO!4}}}
  {\IfNoValueTF{#3}{\IfNoValueTF{#4}{OCCHIO!5}{{{H#2}_{#1}}^{#4}}}
	{\IfNoValueTF{#4}{{^{#3}{H#2}_{#1}}}{OCCHIO!6}}}}
}
\NewDocumentCommand\costfun{d<< d>> d++}{% constant functor
 \IfNoValueTF{#1}{\IfNoValueTF{#2}{\IfNoValueTF{#3}{\mathbf{c}}{OCCHIO1}}{OCCHIO2}}{
\IfNoValueTF{#2}{OCCHIO3}{\IfNoValueTF{#3}{OCCHIO4}{\mathbf{c}\left[#1,#2,#3\right]}}}
}
\newcommand{\Diagcat}{\Delta} % diagonal functor
\NewDocumentCommand\switchcat{d<< d>>}{% funtore permutazione
 \IfNoValueTF{#1}{\IfNoValueTF{#2}{\mathrm{Sw}}{OCCHIO!}}{
	\IfNoValueTF{#2}{\mathrm{Sw}\left[#1\right]}{\mathrm{Sw}\left[#1,#2\right]}}
}
\newcommand{\inclfun}[2]{\mathbf{i}\left[{#1},{#2}\right]} % inclusion functor
\newcommand{\comparrfun}{\circ} % composition of arrows and functors
\newcommand{\nattrasf}{\tau} % natural transformations
\newcommand{\adjnattrasf}{\varphi} % adjunction natural transformations
\newcommand{\sieve}{\mathsf{S}} % generic sieve 
\newcommand{\gensieve}[1]{\sieve_{#1}} % sieve generated by a set of arrow with the same domain 
\NewDocumentCommand\pbsieve{d[] d()}{
\IfNoValueTF{#1}{OCCHIO!1}{\IfNoValueTF{#2}{OCCHIO!2}{{#1}^*\left(#2\right)}}
}
\newcommand{\Grtop}{J} % Grothendieck topology 
\NewDocumentCommand\eqsite{d[]}{
\IfNoValueTF{#1}{\bowtie}{\overset{#1}{\bowtie}}
}
\newcommand{\collarr}{K} % base for a Grothendieck topology
\newcommand{\precollarr}{H} % collection of arrows
\newcommand{\Setcat}{\mathbf{Set}} % category of all sets and set functions
\newcommand{\adjphins}{\varphi} % adjunction between product and hom in the category of sets
\NewDocumentCommand\diradj{m}{\hat{#1}} % aggiunzione su una funzione
\NewDocumentCommand\invadj{m}{\check{#1}} % inversa dell'aggiunzione su una funzione
\NewDocumentCommand\swdiradj{m}{\hat{\ #1}} % aggiunzione su una funzione 
\NewDocumentCommand\swinvadj{m}{\check{\ #1}} % inversa dell'aggiunzione su una funzione 
\newcommand{\Topcat}{\mathbf{Top}} % category of all topological spaces and continuous functions
\newcommand{\GrTc}{T} % Grothendieck topology on \Topcat

\newcommand{\Catmn}{I}
\newcommand{\arrmn}{f}
\newcommand{\Gtmn}{Jv}

\newcommand{\setsymuno}{X}
\newcommand{\elsymuno}{x}
\newcommand{\setsymdue}{Y}
\newcommand{\elsymdue}{y}
\newcommand{\setsymtre}{Z}
\newcommand{\elsymtre}{z}
\newcommand{\setsymquattro}{V}

\newcommand{\elsymquattrobis}{w}
\newcommand{\setsymcinque}{A}
\newcommand{\elsymcinque}{a}
\newcommand{\setsymsei}{B}
\newcommand{\elsymsei}{b}
\newcommand{\setsymsette}{C}

\newcommand{\setsymotto}{O}

\newcommand{\setsymnove}{U}

\newcommand{\setquotfun}{q}
\newcommand{\cardset}{\sharp}
\newcommand{\udenset}{\varnothing} % insieme vuoto
\newcommand{\udenunk}{\varnothing} 
\newcommand{\emptysetfun}{\emptyset}
\newcommand{\sielset}{\ast} % insieme con un singolo elemento
\newcommand{\sielel}{\ast} % elemento generico dell'insieme con un singolo elemento
\newcommand{\partset}[1]{\mathfrak{P}\left({#1}\right)}
\newcommand{\finpartset}[1]{\mathfrak{FP}\left({#1}\right)}
\newcommand{\partiz}{\mathcal{P}}
\newcommand{\realpartset}{\mathbb{P}}
\newcommand{\segnvar}{\ast}

\newcommand{\prefuncaap}[1]{{#1}^{\vartriangle}} % centering a generalized function at a point 
\newcommand{\preffuncaap}[1]{{#1}^{\triangledown}} % basic special set function  
\newcommand{\funcaap}[1]{{#1}^{\blacktriangle}} % centering a pointed generalized continuous function at a point 
\newcommand{\ffuncaap}[1]{{#1}^{\blacktriangledown}} % basic pointed special set function

\newcommand{\setindexuno}{\mathsf{I}}
\newcommand{\topindexuno}{i}
\newcommand{\elindexuno}{\mathsf{i}}
\newcommand{\setindexdue}{\mathsf{J}}

\newcommand{\elindexdue}{\mathsf{j}}

\newcommand{\topindextre}{k}
\newcommand{\elindextre}{\mathsf{k}}

\newcommand{\topindexquattro}{a}
\newcommand{\elindexquattro}{\mathsf{a}}

\newcommand{\topindexcinque}{b}
\newcommand{\elindexcinque}{\mathsf{b}}

\newcommand{\topindexsei}{c}
\newcommand{\elindexsei}{\mathsf{c}}

\newcommand{\topindexotto}{l}
\newcommand{\elindexotto}{\mathsf{l}}
\newcommand{\topindexnove}{m}

\newcommand{\topindexdieci}{n}
\newcommand{\elindexdieci}{\mathsf{n}}
\newcommand{\incloccset}{\varsigma}

\newcommand{\deridxuno}{\lambda} 
\newcommand{\deridxdue}{\sigma}

\newcommand{\ringsym}{\mathbb{A}}
\newcommand{\elringsymuno}{a}
\newcommand{\elringsymdue}{b}
\newcommand{\unoring}{1}
\newcommand{\zeroring}{0}
\newcommand{\neutsumsym}{e}
\newcommand{\neutmultsym}{u}
\newcommand{\modsym}{M}
\newcommand{\algplus}{\oplus} % usato solo in Notation \ref{alg} si usa anche per la somma diretta, meglio cambiare simbolo 
\newcommand{\algper}{\odot} % usato solo in Notation \ref{alg} mettere un simbolo concorde con il nuovo simbolo per \algplus 
\newcommand{\algsclp}{\circledast} % usato solo in Notation \ref{alg} mettere un simbolo concorde con il nuovo simbolo per \algplus 
\newcommand{\VecR}{\mathbf{Vec}_{\mathbb{R}}}
\newcommand{\AlgR}{\mathbf{Alg}_{\mathbb{R}}}
\newcommand{\dirsum}{\oplus}
\newcommand{\bigdirsum}{\bigoplus}
\newcommand{\tensprodR}{\otimes}
\newcommand{\bigtensprodR}{\bigotimes}
\newcommand{\perm}{\sigma}
\newcommand{\innpreul}{\langle}
\newcommand{\innpreur}{\rangle}
\newcommand{\normeul}{\mid}
\newcommand{\normeur}{\mid}
\NewDocumentCommand\tenalgfun{d[]}{% tensor algebra functor 
\IfNoValueTF{#1}{\mathbf{T}}{\mathbf{T}_{#1}}
}
\NewDocumentCommand\unoT{d[]}{
\IfNoValueTF{#1}{u}{u_{#1}}
} % unità tensoriale
\NewDocumentCommand\symalgfun{d[]}{% symmetric algebra functor
\IfNoValueTF{#1}{\mathbf{S}}{\mathbf{S}_{#1}}
}
\NewDocumentCommand\extalgfun{d[]}{% exterior algebra functor
\IfNoValueTF{#1}{\mathbf{E}}{\mathbf{E}_{#1}}
}
\newcommand{\Aug}{\mathbf{A}} % adding the unity to algebras functiors
\NewDocumentCommand\unoA{d[]}{
\IfNoValueTF{#1}{e}{e_{#1}}
} % unità aggiunta
\NewDocumentCommand\inclunoA{d[]}{
\IfNoValueTF{#1}{i}{i_{#1}}
} % inclusione canonica

\newcommand{\ef}{\mathbf{e}} % euclidean frame
\newcommand{\vecuno}{\mathsf{v}}
\newcommand{\vecdue}{\mathsf{w}}

\newcommand{\intuno}{S}
\newcommand{\unkuno}{s}
\newcommand{\intdue}{T}
\newcommand{\unkdue}{t}
\newcommand{\inttre}{R}
\newcommand{\unktre}{r}
\newcommand{\opsspint}{\vDash}
\newcommand{\sqsubseteqdentro}{\overset{\circ}{\sqsubseteq}} 
\newcommand{\domint}{\mathcal{I}}
\newcommand{\att}{\mathcal{A}}
\newcommand{\attfun}{\vartheta}
\newcommand{\pmRn}{\sigma}

\newcommand{\topo}{\mbox{\boldmath$\tau$}}
\newcommand{\compacont}{\subset\subset}
\newcommand{\closure}{\mathsf{Cl}}
\newcommand{\interior}{\mathsf{In}}
\newcommand{\subseteqdentro}{\overset{\circ}{\subseteq}}
\newcommand{\dentro}{\overset{\circ}{\in}} 
\newcommand{\openset}{P}
\newcommand{\SiteSpTop}{\mathcal{O}} % funtore dalla categoria degli spazi topologici a quella dei siti che associa ad ogni spazio topologico il sito dei suoi aperti e delle relative inclsioni (per ora usato solo per gli oggetti)

\newcommand{\Diag}[2]{\Delta\left[{#1},{#2}\right]}

\newcommand{\incl}[2]{\mathrm{i}\left[{#1},{#2}\right]}
\NewDocumentCommand\cost{d<< d>> d++}{
 \IfNoValueTF{#1}{\IfNoValueTF{#2}{\IfNoValueTF{#3}{\mathrm{c}}{OCCHIO1}}{OCCHIO2}}{
\IfNoValueTF{#2}{OCCHIO3}{\IfNoValueTF{#3}{OCCHIO4}{\mathrm{c}\left[#1,#2,#3\right]}}}
}
\NewDocumentCommand\proj{d<< d>>}{% funzione proiezone su un fattore di un prodotto
 \IfNoValueTF{#1}{\IfNoValueTF{#2}{\mathrm{p}}{\mathrm{p}_{#2}}}{
	\IfNoValueTF{#2}{OCCHIO1}{\mathrm{p}_{#2}\left[#1\right]}}
}
\NewDocumentCommand\pointedincl{d<< d>>}{% funzione inclusione su un fattore di un prodotto
 \IfNoValueTF{#1}{\IfNoValueTF{#2}{\mathrm{b}}{OCCHIO1}}{
	\IfNoValueTF{#2}{OCCHIO1}{\mathrm{b}_{#2}\left[#1\right]}}
}
\NewDocumentCommand\switch{d<< d>>}{% permutazione prodotto
 \IfNoValueTF{#1}{\IfNoValueTF{#2}{\mathrm{sw}}{OCCHIO!}}{
	\IfNoValueTF{#2}{\mathrm{sw}\left[#1\right]}{\mathrm{sw}\left[#1,#2\right]}}
}
\NewDocumentCommand\vecsum{d<< d>>}{
\IfNoValueTF{#1}{\IfNoValueTF{#2}{+}{OCCHIO1!}}{\IfNoValueTF{#2}{OCCHIO2!}{+\left[{#1},{#2}\right]}}
}
\newcommand{\vecminus}{-}
\newcommand{\vecprod}{\cdot}
\newcommand{\sectn}{sc}
\newcommand{\proje}[2]{pr\left[{#1},{#2}\right]}
\newcommand{\coord}[2]{co\left[{#1},{#2}\right]}
\NewDocumentCommand\trasl{d[]}{
\IfNoValueTF{#1}{\mathrm{tr}}{\mathrm{tr}_{#1}}
}

\newcommand{\funsum}{+}
\newcommand{\funminus}{-}
\newcommand{\funmult}{\cdot}
\newcommand{\funcomp}{\circ}
\newcommand{\funscalp}{\ast}
\NewDocumentCommand\valass{d||}{%
  \IfNoValueTF{#1}{\mid\,\mid}{\mid\!#1\!\mid}
} 
\newcommand{\Cod}{\mathsf{Cod}}
\newcommand{\Dom}{\mathsf{Dom}}
\newcommand{\Cae}{\mathsf{Ce}}
\newcommand{\Ima}{\mathsf{Im}}

\newcommand{\smder}{\partial} 

\NewDocumentCommand\symalgtvsm{d<< d>>}{% symmetric algebra functor
\IfNoValueTF{#1}{\IfNoValueTF{#2}{\Delta^{\infty}}{OCCHIO1!}}{\IfNoValueTF{#2}{OCCHIO2!}{\Delta^{\infty}_{#1}\left[#2\right]}}
}
\NewDocumentCommand\symjacsm{d<< d>> d[]}{% symmetric algebra functor
\IfNoValueTF{#1}{\IfNoValueTF{#2}{\IfNoValueTF{#3}{d^{\infty}}{OCCHIO!1}}{OCCHIO2!}}{\IfNoValueTF{#2}{OCCHIO3!}{\IfNoValueTF{#3}{d^{\infty}_{#1}\left[#2\right]}}{d^{\infty}_{#1}\left[#2,#3\right]}}
}
\NewDocumentCommand\extalgtvsm{d<< d>>}{% symmetric algebra functor
\IfNoValueTF{#1}{\IfNoValueTF{#2}{\Lambda^{\infty}}{OCCHIO1!}}{\IfNoValueTF{#2}{OCCHIO2!}{\Lambda^{\infty}_{#1}\left[#2\right]}}
}

\newcommand{\smint}{\int}
\newcommand{\intvar}{\lambda}

\NewDocumentCommand\natisprodfibsm{d// d-- d<< d>>}{%
  \IfNoValueTF{#1}{\IfNoValueTF{#2}{\tau^{\infty}\left[#3,#4\right]}{\tau^{\infty}_{#3,#4}\left(#2\right)}}
	{\IfNoValueTF{#2}{OCCHIO1!}{\left(\tau^{\infty}_{#3,#4}\right)^{#1}\left(#2\right)}}
} 
\NewDocumentCommand\natisprodfibinvsm{d// d-- d<< d>>}{
  \IfNoValueTF{#1}{\IfNoValueTF{#2}{\sigma^{\infty}\left[#3,#4\right]}{\sigma^{\infty}_{#3,#4}\left(#2\right)}}
	{\IfNoValueTF{#2}{OCCHIO1!}{\left(\sigma^{\infty}_{#3,#4}\right)^{#1}\left(#2\right)}}
}

\newcommand{\bsfM}{\mbox{\boldmath$\mathsf{M}$}}
\newcommand{\subbdiffsfM}{\mbox{\boldmath$\mathsf{D}$}}
\newcommand{\subintsfM}{\mbox{\boldmath$\mathsf{I}$}}
\newcommand{\subbsfM}{\mbox{\boldmath$\mathsf{O}$}}
\newcommand{\bsfm}{\mbox{\boldmath$\mathsf{m}$}}
\newcommand{\bsfn}{\mbox{\boldmath$\mathsf{n}$}}
\newcommand{\bsfuno}{\mbox{\boldmath$\mathsf{1}$}}
\newcommand{\eumea}{\mathfrak{m}}
\newcommand{\Fint}{\mbox{\boldmath$\varint$}}
\newcommand{\Fpart}{\mathsf{d}}
\newcommand{\Fp}{\mathsf{c}}
\newcommand{\Fq}{\mathsf{p}}
\newcommand{\Fqq}{\mathsf{q}}
\newcommand{\bsfrelsymb}{=}

\newcommand{\abidsp}{\mathcal{X}}
\newcommand{\abgenf}{X}
\newcommand{\abgenfzero}{\abgenf_0}
\newcommand{\abgenfempty}{\abgenf_{\emptyset}}
\newcommand{\abgenfuncomp}{\circledcirc}
\newcommand{\abgenlbound}{\leftmoon}
\newcommand{\abgenrbound}{\rightmoon}
\newcommand{\abgenfunsum}{\oplus}
\newcommand{\bigabgenfunsum}{\scalerel*{\abgenfunsum}{\sum}}
\newcommand{\abgenfunmult}{\odot}
\newcommand{\abgenfunscalp}{\circledast}
\newcommand{\abbsfmu}{\bullet}
\newcommand{\abidsptwo}{\widetilde{\mathcal{X}}}
\newcommand{\abgenftwo}{Y}
\newcommand{\abgenfzerotwo}{\abgenftwo_0}
\newcommand{\abgenfemptytwo}{\abgenfempty}
\newcommand{\abgenfuncomptwo}{\boxcircle}
\newcommand{\abgenlboundtwo}{\llbracket}
\newcommand{\abgenrboundtwo}{\rrbracket}
\newcommand{\abgenfunsumtwo}{\boxplus}
\newcommand{\abgenfunmulttwo}{\boxdot}
\newcommand{\abgenfunscalptwo}{\boxast}
\newcommand{\abbsfmutwo}{{\scriptstyle\blacksquare}}
\newcommand{\IDSMa}{\mathbf{IDS}}
\newcommand{\QIDSMa}{\mathbf{IS}}
\newcommand{\ffIDSMaQIDSMa}{\dot{\Phi}}

\newcommand{\qidscont}{\mathcal{C}^0}
\newcommand{\idssmo}{\mathcal{C}^{\infty}}
\NewDocumentCommand\Cksp{m d() d++}{%
  \IfNoValueTF{#2}{\IfNoValueTF{#3}{C^{#1}}{C^{#1}\left[#3\right]}}{\IfNoValueTF{#3}{OCCHIO!}{C^{#1}\left[#2,#3\right]}}
}
\NewDocumentCommand\Ckspzero{m d() d++}{%
  \IfNoValueTF{#2}{\IfNoValueTF{#3}{C^{#1}_0}{C^{#1}_0\left[#3\right]}}{\IfNoValueTF{#3}{OCCHIO!}{C^{#1}_0\left[#2,#3\right]}}
} 
\newcommand{\smospempty}{\emptyset}
\newcommand{\smosplbound}{\left(}
\newcommand{\smosprbound}{\right)}
\newcommand{\smospbsfmu}{\Join}
\newcommand{\topS}{\mbox{\boldmath$\sigma$}_{\Cksp{\infty}}}
\NewDocumentCommand\Cksploc{m d() d++}{%
  \IfNoValueTF{#2}{\IfNoValueTF{#3}{\overset{\bullet}{C}^{#1}}{\overset{\bullet}{C}^{#1}\left[#3\right]}}{\IfNoValueTF{#3}{OCCHIO!}{\overset{\bullet}{C}^{#1}\left[#2,#3\right]}}
}
\NewDocumentCommand\Cksplocccz{m d() d++}{%
  \IfNoValueTF{#2}{\IfNoValueTF{#3}{\overset{\bullet}{C}^{#1}_{CC}}{\overset{\bullet}{C}^{#1}_CC\left[#3\right]}}{\IfNoValueTF{#3}{OCCHIO!}{\overset{\bullet}{C}^{#1}_CC
	\left[#2,#3\right]}}
	}
\NewDocumentCommand\Ckspquot{m d// d-- d() d++}{
 \IfNoValueTF{#2}{
 \IfNoValueTF{#3}{
	\IfNoValueTF{#4}{\IfNoValueTF{#5}{{\mathbf{C}}^{#1}}{{\mathbf{C}}^{#1}\left[#5\right]}}{\IfNoValueTF{#5}{OCCHIO!1}{{\mathbf{C}}^{#1}\left[#4,#5\right]}}
	}{
	\IfNoValueTF{#4}{\IfNoValueTF{#5}{{\mathbf{C}}^{#1}}{{\mathbf{C}}^{#1}_{#5}}}{\IfNoValueTF{#5}{OCCHIO!2}{{\mathbf{C}}^{#1}_{#4   ,#5}}}
  }}
	{\IfNoValueTF{#3}{
	\IfNoValueTF{#4}{\IfNoValueTF{#5}{OCCHIO!3}{{\mathbf{C}}^{#1}_{#5}}}{\IfNoValueTF{#5}{OCCHIO!4}{{\mathbf{C}}^{#1}_{#4,#5}}}
	}
	{
	\IfNoValueTF{#4}{\IfNoValueTF{#5}{OCCHIO!5}{{\mathbf{C}}^{#1\,#3}_{#5}}}{\IfNoValueTF{#5}{OCCHIO!6}{{\mathbf{C}}^{#1\,#3}_{#4,#5}}}
	}}
}
\NewDocumentCommand\Ckspquotfun{d// d-- d() d++}{ 
 \IfNoValueTF{#1}{
 \IfNoValueTF{#2}{
	\IfNoValueTF{#3}{\IfNoValueTF{#4}{\check{\mathfrak{q}}}{\check{\mathfrak{q}}\left[#4\right]}}{\IfNoValueTF{#4}{OCCHIO!1}{\check{\mathfrak{q}}\left[#3,#4\right]}}
	}{OCCHIO!2}}
	{\IfNoValueTF{#2}{
	\IfNoValueTF{#3}{\IfNoValueTF{#4}{OCCHIO!3}{\check{\mathfrak{q}}_{#4}}}{\IfNoValueTF{#4}{OCCHIO!4}{\check{\mathfrak{q}}_{#3,#4}}}
	}
	{
	\IfNoValueTF{#3}{\IfNoValueTF{#4}{OCCHIO!5}{\check{\mathfrak{q}}^{#2}_{#4}}}{\IfNoValueTF{#4}{OCCHIO!6}{\check{\mathfrak{q}}^{#2}_{#3,#4}}}
	}}
}
\NewDocumentCommand\Ckspquotsf{m d|| d() d++}{\widehat{\mathbf{C}}^{#1}_{#2,#3,#4}}
\NewDocumentCommand\Ckspquotccz{m d// d-- d() d++}{
 \IfNoValueTF{#2}{
 \IfNoValueTF{#3}{
	\IfNoValueTF{#4}{\IfNoValueTF{#5}{{\mathbf{C}}^{#1}_{CC}}{{\mathbf{C}}^{#1}_{CC}\left[#5\right]}}{\IfNoValueTF{#5}{OCCHIO!1}{{\mathbf{C}}^{#1}_{CC}\left[#4,#5\right]}}
	}{OCCHIO!2}}
	{\IfNoValueTF{#3}{
	\IfNoValueTF{#4}{\IfNoValueTF{#5}{OCCHIO!3}{{\mathbf{C}}^{#1}_{CC\, #5}}}{\IfNoValueTF{#5}{OCCHIO!4}{{\mathbf{C}}^{#1}_{CC\, #4,#5}}}}
	{
	\IfNoValueTF{#4}{\IfNoValueTF{#5}{OCCHIO!5}{{\mathbf{C}}^{#1\,#3}_{CC\,#5}}}{\IfNoValueTF{#5}{OCCHIO!6}{{\mathbf{C}}^{#1\,#3}_{CC\,#4,#5}}}
	}}
}
\NewDocumentCommand\Ckspquotcczsf{m d|| d() d++}{\widehat{\mathbf{C}}^{#1}_{CC,#2,#3,#4}}
\NewDocumentCommand\gczsfx{d[]}{
\IfNoValueTF{#1}{\mathbf{f}}{
\wr #1\wr} 
}
\NewDocumentCommand\gczsfy{d[]}{
\IfNoValueTF{#1}{\mathbf{g}}{
\wr #1\wr} 
}
\NewDocumentCommand\gczsfz{d[]}{
\IfNoValueTF{#1}{\mathbf{h}}{
\wr #1\wr} 
}
\newcommand{\classquotempty}{\smospempty}
\NewDocumentCommand\ckfggerm{d[]}{ 
\IfNoValueTF{#1}{\widetilde{\gczsfx}}{
\grm[#1]}
}

\newcommand{\contspempty}{\smospempty} 
\newcommand{\contsplbound}{\smosplbound}
\newcommand{\contsprbound}{\smosprbound}
\newcommand{\contspbsfmu}{\smospbsfmu}
\newcommand{\topC}{\mbox{\boldmath$\sigma$}_{\Cksp{0}}}
\newcommand{\lininclquot}{\mathfrak{li}}
\newcommand{\evalcompquot}{\mathfrak{ev}}
\newcommand{\dimCod}{\widetilde{\Cod}}
\newcommand{\dimDom}{\widetilde{\Dom}}

\newcommand{\bascont}{\beth}

\newcommand{\fremaggen}{G}
\newcommand{\fremag}{\setsymuno}
\newcommand{\fremagcont}{\fremag_C}
\newcommand{\fremagcontcont}{\fremag_{CC}}
\newcommand{\fremagsmooth}{\fremag_S}
\newcommand{\fremagsmoothsmooth}{\fremag_{SS}}
\NewDocumentCommand\canrapp{d()}{
  \IfNoValueTF{#1}{L}{L_{#1}}
} 
\newcommand{\fmx}{\elsymuno}
\newcommand{\fmy}{\elsymdue}
\newcommand{\fmz}{\elsymtre}
\newcommand{\fmw}{\elsymquattrobis}
\newcommand{\compone}{\,\blacklozenge\,}
\newcommand{\codmagone}{\mathbf{C}}
\newcommand{\dommagone}{\mathbf{D}}
\newcommand{\dimcodmagone}{\widetilde{\mathbf{C}}}
\newcommand{\dimdommagone}{\widetilde{\mathbf{D}}}
\newcommand{\inclfremaglev}{\mathbf{i}}
\newcommand{\incllim}{\mbox{\boldmath$\iota$}} 
\newcommand{\lboundone}{\left\bracevert}
\newcommand{\rboundone}{\right\bracevert}

\newcommand{\evalcompone}{\mathbf{ev}}
\newcommand{\lininclone}{\mathbf{li}}
\newcommand{\quotmagone}{\mathbf{q}}
\newcommand{\acmone}{\blacktriangleright}
\newcommand{\symb}{\widehat{\fremaggen}}
\newcommand{\extfremag}{\widehat{\fremag}}
\newcommand{\extx}{\widehat{\fmx}}
\newcommand{\extxtwo}{\widehat{\fmy}}
\newcommand{\extsubstrchar}{\widehat{\mathbf{\fmx}}}
\newcommand{\extlnght}{\mbox{\boldmath$\lambda$}}
\newcommand{\extoccset}{\widehat{\occset}}

\newcommand{\extoccfun}{\widehat{\occfun}}
\newcommand{\extoccfundue}{\widehat{\occfundue}}
\newcommand{\extstrchar}{\widehat{\strchar}}
\newcommand{\extasselstr}{\extoccfun}
\newcommand{\extoccsetstr}{\extoccset}
\newcommand{\occset}{O} %occurrence set
\newcommand{\occsetdue}{P} %occurrence set 2
\newcommand{\occsettre}{Q}
\newcommand{\occfun}{\eta} %occurrence function
\newcommand{\occfundue}{\xi} %occurrence function 2
\newcommand{\strchar}{\mathbf{u}}
\newcommand{\occfunstr}{\occfun}
\newcommand{\occsetstr}{\occset}
\newcommand{\assfun}{\alpha} %associating function
\newcommand{\assfundue}{\beta} %associating function 2
\newcommand{\assfuntre}{\gamma} %associating function 3
\newcommand{\assfunadj}{\textnormal{adj}}
\newcommand{\pmone}{\chi}  %path 
\newcommand{\singsmpth}{\mathsf{S}} %singular set for smooth paths
\newcommand{\ordrelsymb}{\prec}
\newcommand{\relsymbord}{\succ\!\prec}
\newcommand{\relsymb}{\approx}
\newcommand{\fremagfunonsm}{\mathbf{\fmx}}
\newcommand{\fremagfunoncont}{\mathbf{\fmx}}

\newcommand{\magtwo}{W}
\newcommand{\magtwocont}{\magtwo_C}
\newcommand{\magtwocontcont}{\magtwo_{CC}}
\newcommand{\magtwosmooth}{\magtwo_S}
\newcommand{\magtwosmoothsmooth}{\magtwo_{SS}}
\newcommand{\magtwoempty}{\magtwo_{\emptyset}}
\newcommand{\mtx}{u}
\newcommand{\mty}{v}
\newcommand{\mtz}{w}
\newcommand{\comptwo}{\,\lozenge\,}
\newcommand{\comptwocont}{\,\vartriangle_C\,}
\newcommand{\comptwosmooth}{\,\vartriangle_S\,}
\newcommand{\codmagtwo}{\textnormal{C}}
\newcommand{\dommagtwo}{\textnormal{D}}
\newcommand{\dimcodmagtwo}{\widetilde{\textnormal{C}}}
\newcommand{\dimdommagtwo}{\widetilde{\textnormal{D}}}
\newcommand{\lboundtwo}{\arrowvert}
\newcommand{\rboundtwo}{\arrowvert}
\newcommand{\lboundtwocont}{\lboundtwo}
\newcommand{\rboundtwocont}{\rboundtwo_C}
\newcommand{\lboundtwosmooth}{\lboundtwo}
\newcommand{\rboundtwosmooth}{\rboundtwo_S}
\newcommand{\evalcomptwo}{\mathrm{ev}}
\newcommand{\linincltwo}{\mathrm{li}}

\newcommand{\acmtwo}{\rhd}
\newcommand{\acmtwocont}{\,\acmtwo_C\,}
\newcommand{\acmtwosmooth}{\,\acmtwo_S\,}
\newcommand{\pmtwo}{\omega}
\newcommand{\toppathmagtwo}{\mbox{\boldmath$\tau$}}
\newcommand{\smtoppathmagtwo}{\mbox{\boldmath$\sigma$}}

\newcommand{\genfidsp}{\mathcal{\genf}}
\NewDocumentCommand\genf{d() d++}{%
  \IfNoValueTF{#1}{\IfNoValueTF{#2}{\magtwo}{\magtwo\left(#2\right)}}{\IfNoValueTF{#2}{OCCHIO!}{\magtwo\left(#1,#2\right)}}
}
\NewDocumentCommand\genfzero{d() d++}{%
  \IfNoValueTF{#1}{\IfNoValueTF{#2}{\genf_0}{\genf_0\left(#2\right)}}{\IfNoValueTF{#2}{OCCHIO!}{\genf_0\left(#1,#2\right)}}
}
\NewDocumentCommand\genfempty{d() d++}{%
  \IfNoValueTF{#1}{\IfNoValueTF{#2}{\magtwoempty}{\magtwoempty\left(#2\right)}}{\IfNoValueTF{#2}{OCCHIO!}{\magtwoempty\left(#1,#2\right)}}
}
\newcommand{\genfuncomp}{\comptwo}
\newcommand{\lboundgenf}{\lboundtwo}
\newcommand{\rboundgenf}{\rboundtwo}
\newcommand{\genfunsum}{\oplus}
\newcommand{\biggenfunsum}{\scalerel*{\genfunsum}{\sum}}
\newcommand{\genfunmult}{\odot}

\newcommand{\genfunscalp}{\circledast}
\newcommand{\genfbsfmu}{\acmtwo}
\newcommand{\gfx}{\mtx}
\newcommand{\gfy}{\mty}
\newcommand{\gfz}{\mtz}
\newcommand{\gfu}{y} 
\newcommand{\gfw}{x}
\newcommand{\gfv}{z}
\newcommand{\genfidspcont}{\genfidsp_C}
\NewDocumentCommand\genfcont{d() d++}{%
  \IfNoValueTF{#1}{\IfNoValueTF{#2}{\genf_C}{\genf_C\left(#2\right)}}{\IfNoValueTF{#2}{OCCHIO!}{\genf_C\left(#1,#2\right)}}
}
\NewDocumentCommand\genfcontcont{d[] d() d++}{%
 \IfNoValueTF{#1}{\IfNoValueTF{#2}{\IfNoValueTF{#3}{\genf_{CC}}{\genf_{CC}\left(#3\right)}}{\IfNoValueTF{#3}{OCCHIO!}{\genf_{CC}\left(#2,#3\right)}}}{\IfNoValueTF{#2}{\IfNoValueTF{#3}{\genf_{CC, #1}}{\genf_{CC, #1}\left(#3\right)}}{\IfNoValueTF{#3}{OCCHIO!}{\genf_{CC, #1}\left(#2,#3\right)}}}
}
\newcommand{\genfidspsmooth}{\genfidsp_S}
\NewDocumentCommand\genfsmooth{d() d++}{%
  \IfNoValueTF{#1}{\IfNoValueTF{#2}{\genf_{S}}{\genf_{S}\left(#2\right)}}{\IfNoValueTF{#2}{OCCHIO!}{\genf_{S}\left(#1,#2\right)}}
} 
\newcommand{\smofunconfun}{\iota}

\NewDocumentCommand\genfloc{d() d++}{%
  \IfNoValueTF{#1}{\IfNoValueTF{#2}{\overset{\bullet}{\genf}}{\overset{\bullet}{\genf}\left(#2\right)}}{\IfNoValueTF{#2}{OCCHIO!}{\overset{\bullet}{\genf}\left(#1,#2\right)}}
}
\NewDocumentCommand\genflocccz{d() d++}{%
  \IfNoValueTF{#1}{\IfNoValueTF{#2}{\genfloc_{CC}}{\genfloc_{CC}\left(#2\right)}}{\IfNoValueTF{#2}{OCCHIO!}{\genfloc_{CC}\left(#1,#2\right)}}
} 
\NewDocumentCommand\genfloccczz{d() d++}{%
  \IfNoValueTF{#1}{\IfNoValueTF{#2}{\genfloc_{CC0}}{\genfloc_{CC0}\left(#2\right)}}{\IfNoValueTF{#2}{OCCHIO!}{\genfloc_{CC0}\left(#1,#2\right)}}
}

\newcommand{\infneareqrel}{\Bumpeq}
\NewDocumentCommand\genfquot{d// d-- d() d++}{    
 \IfNoValueTF{#1}{
 \IfNoValueTF{#2}{
	\IfNoValueTF{#3}{\IfNoValueTF{#4}{V}{V\left[#4\right]}}{\IfNoValueTF{#4}{OCCHIO!1}{V\left[#3,#4\right]}}
	}{OCCHIO!2}
	}
	{\IfNoValueTF{#2}{
	\IfNoValueTF{#3}{\IfNoValueTF{#4}{OCCHIO!3}{V_{#4}}}{\IfNoValueTF{#4}{OCCHIO!4}{V_{#3,#4}}}
	}
	{
	\IfNoValueTF{#3}{\IfNoValueTF{#4}{OCCHIO!5}{V^{#2}_{#4}}}{\IfNoValueTF{#4}{OCCHIO!6}{V^{#2}_{#3,#4}}}
	}}
} 
\NewDocumentCommand\genfquotfun{d// d-- d() d++}{ 
 \IfNoValueTF{#1}{
 \IfNoValueTF{#2}{
	\IfNoValueTF{#3}{\IfNoValueTF{#4}{\check{\mathfrak{q}}}{\check{\mathfrak{q}}\left[#4\right]}}{\IfNoValueTF{#4}{OCCHIO!1}{\check{\mathfrak{q}}\left[#3,#4\right]}}
	}{OCCHIO!2}}
	{\IfNoValueTF{#2}{
	\IfNoValueTF{#3}{\IfNoValueTF{#4}{OCCHIO!3}{\check{\mathfrak{q}}_{#4}}}{\IfNoValueTF{#4}{OCCHIO!4}{\check{\mathfrak{q}}_{#3,#4}}}
	}
	{
	\IfNoValueTF{#3}{\IfNoValueTF{#4}{OCCHIO!5}{\check{\mathfrak{q}}^{#2}_{#4}}}{\IfNoValueTF{#4}{OCCHIO!6}{\check{\mathfrak{q}}^{#2}_{#3,#4}}}
	}}
}
\NewDocumentCommand\genfquotcccz{d// d-- d() d++}{
 \IfNoValueTF{#1}{
 \IfNoValueTF{#2}{
	\IfNoValueTF{#3}{\IfNoValueTF{#4}{\genfquot_{CC}}{\genfquot_{CC}\left[#4\right]}}{\IfNoValueTF{#4}{OCCHIO!1}{\genfquot_{CC}\left[#3,#4\right]}}
	}{OCCHIO!2}}
	{\IfNoValueTF{#2}{
	\IfNoValueTF{#3}{\IfNoValueTF{#4}{OCCHIO!3}{\genfquot_{CC}_{#4}}}{\IfNoValueTF{#4}{OCCHIO!4}{\genfquot_{CC}_{#3,#4}}}
	}
	{
	\IfNoValueTF{#3}{\IfNoValueTF{#4}{OCCHIO!5}{\genfquot_{CC}^{#2}_{#4}}}{\IfNoValueTF{#4}{OCCHIO!6}{\genfquot_{CC}^{#2}_{#3,#4}}}
	}}
}
\NewDocumentCommand\genfquotccczz{d// d-- d() d++}{ 
 \IfNoValueTF{#1}{
 \IfNoValueTF{#2}{
	\IfNoValueTF{#3}{\IfNoValueTF{#4}{\genfquot_{CC0}}{\genfquot_{CC0}\left[#4\right]}}{\IfNoValueTF{#4}{OCCHIO!1}{\genfquot_{CC0}\left[#3,#4\right]}}
	}{OCCHIO!2}}
	{\IfNoValueTF{#2}{
	\IfNoValueTF{#3}{\IfNoValueTF{#4}{OCCHIO!3}{\genfquot_{CC0}_{#4}}}{\IfNoValueTF{#4}{OCCHIO!4}{\genfquot_{CC0}_{#3,#4}}}
	}
	{
	\IfNoValueTF{#3}{\IfNoValueTF{#4}{OCCHIO!5}{\genfquot_{CC0}^{#2}_{#4}}}{\IfNoValueTF{#4}{OCCHIO!6}{\genfquot_{CC0}^{#2}_{#3,#4}}}
	}}
}
\newcommand{\zeroquot}{0}
\NewDocumentCommand\gggfx{d[]}{
\IfNoValueTF{#1}{\mathbf{\gfx}}{
\wr #1\wr}
} 
\NewDocumentCommand\gggfy{d[]}{
\IfNoValueTF{#1}{\mathbf{\gfy}}{
\wr #1\wr}
}
\NewDocumentCommand\gggfz{d[]}{
\IfNoValueTF{#1}{\mathbf{\gfz}}{
\wr #1\wr}
}
\NewDocumentCommand\gggfu{d[]}{
\IfNoValueTF{#1}{\mathbf{\gfu}}{
\wr #1\wr}
}
\NewDocumentCommand\gggfw{d[]}{
\IfNoValueTF{#1}{\mathbf{\gfw}}{
\wr #1\wr}
}
\NewDocumentCommand\gggfv{d[]}{
\IfNoValueTF{#1}{\mathbf{\gfv}}{
\wr #1\wr}
}
\newcommand{\codgerm}{\overline{\codmagtwo}}
\newcommand{\domgerm}{\overline{\dommagtwo}}
\newcommand{\lboundquotpoint}{\langle}
\newcommand{\rboundquotpoint}{\rangle}
\newcommand{\sqsumquotpoint}{\,\widehat{\boxplus}\,}

\newcommand{\sumquotpoint}{\,\widehat{\oplus}\,}
\newcommand{\bigsumquotpoint}{\widehat{\scalerel*{\oplus}{\sum}}}

\newcommand{\sqmultquotpoint}{\,\widehat{\boxdot}\,}

\newcommand{\multquotpoint}{\,\widehat{\odot}\,}
\newcommand{\bigmultquotpoint}{\widehat{\scalerel*{\odot}{\sum}}}
\newcommand{\compgenfquot}{\,\square\,}
\newcommand{\compquotpoint}{\,\widehat{\circledcirc}\,}
\newcommand{\scalpquotpoint}{\,\widehat{\circledast}\,}
\newcommand{\bsfmuquotpoint}{\,\widehat{\diamond}\,}
\newcommand{\evalcompquotcontcontzero}{\overline{\mathfrak{e}}}
\newcommand{\lboundquotgradpoint}{\left\lfloor }
\newcommand{\rboundquotgradpoint}{\right\rceil}
\newcommand{\sqsumquotgradpoint}{\,\overline{\boxplus}\,}
\newcommand{\sumquotgradpoint}{\,\overline{\oplus}\,}
\newcommand{\sqmultquotgradpoint}{\,\overline{\boxdot}\,}
\newcommand{\multquotgradpoint}{\,\overline{\odot}\,}
\newcommand{\scalpquotgradpoint}{\,\overline{\circledast}\,}
\newcommand{\bsfmuquotgradpoint}{\,\overline{\diamond}\,}
\NewDocumentCommand\vandergf{m}{Z_{#1}}
\NewDocumentCommand\ppvandergf{m}{Z^{\perp}_{#1}}
\newcommand{\pmV}{\nu}
\newcommand{\pmVdue}{\mu}
\newcommand{\ssf}[1]{SSF\left(#1\right)} % special set functions
\NewDocumentCommand\gfggerm{d[]}{ 
\IfNoValueTF{#1}{\mbox{\boldmath$\pmV$}}{
\grm[#1]}
} 
\NewDocumentCommand\gfggermdue{d[]}{ 
\IfNoValueTF{#1}{\mbox{\boldmath$\pmVdue$}}{
\grm[#1]}
}

\newcommand{\genpreder}{\mbox{\boldmath$\partial$}}
\NewDocumentCommand\nullpreder{d// d-- d[] d() d++}{%
 \IfNoValueTF{#1}{
\IfNoValueTF{#2}{\IfNoValueTF{#3}{\IfNoValueTF{#4}{\IfNoValueTF{#5}{\mathcal{N}}{\mathcal{N}\left[#5\right]}}{\IfNoValueTF{#5}{OCCHIO!1}{\mathcal{N}\left[#4,#5\right]}}}{\IfNoValueTF{#4}{OCCHIO!2}{\IfNoValueTF{#5}{OCCHIO!3}{\mathcal{N}\left[#3,#4,#5\right]}}}}
	{\IfNoValueTF{#3}{\IfNoValueTF{#4}{\IfNoValueTF{#5}{\mathcal{N}}{\mathcal{N}_{#5}}}{\IfNoValueTF{#5}{OCCHIO!6}{\mathcal{N}_{#4,#5}}}}{\IfNoValueTF{#4}{OCCHIO!7}{\IfNoValueTF{#5}{OCCHIO!8}{\mathcal{N}_{#3,#4,#5}}}}}}	
{\IfNoValueTF{#2}{OCCHIO!9}{\IfNoValueTF{#3}{\IfNoValueTF{#4}{\IfNoValueTF{#5}{OCCHIO!10}{\mathcal{N}^{#1}_{#5}}}{\IfNoValueTF{#5}{OCCHIO!11}{\mathcal{N}^{#1}_{#4,#5}}}}{\IfNoValueTF{#4}{OCCHIO!12}{\IfNoValueTF{#5}{OCCHIO!13}{\mathcal{N}^{#1}_{#3,#4,#5}}}}}
	}
} 
\NewDocumentCommand\gennullpreder{d// d-- d[] d() d++}{%
 \IfNoValueTF{#1}{
\IfNoValueTF{#2}{\IfNoValueTF{#3}{\IfNoValueTF{#4}{\IfNoValueTF{#5}{G}{G\left[#5\right]}}{\IfNoValueTF{#5}{OCCHIO!1}{G\left[#4,#5\right]}}}{\IfNoValueTF{#4}{OCCHIO!2}{\IfNoValueTF{#5}{OCCHIO!3}{G\left[#3,#4,#5\right]}}}}
	{\IfNoValueTF{#3}{\IfNoValueTF{#4}{\IfNoValueTF{#5}{G\left(#2\right)}{G_{#5}\left(#2\right)}}{\IfNoValueTF{#5}{OCCHIO!4}{G_{#4,#5}\left(#2\right)}}}{\IfNoValueTF{#4}{OCCHIO!5}{\IfNoValueTF{#5}{OCCHIO!6}{G_{#3,#4,#5}\left(#2\right)}}}}}
{\IfNoValueTF{#2}{\IfNoValueTF{#3}{\IfNoValueTF{#4}{\IfNoValueTF{#5}{G}{G\left[#5\right]}}{\IfNoValueTF{#5}{OCCHIO!7}{G\left[#4,#5\right]}}}{\IfNoValueTF{#4}{OCCHIO!8}{\IfNoValueTF{#5}{OCCHIO!9}{G\left[#3,#4,#5\right]}}}}
	{OCCHIO10}}
} 
\newcommand{\setofcores}{\Omega}
\NewDocumentCommand\genprederspace{d// d-- d[] d() d++}{
\IfNoValueTF{#1}{
\IfNoValueTF{#2}{\IfNoValueTF{#3}{\IfNoValueTF{#4}{\IfNoValueTF{#5}{\mathfrak{PDer}}{\mathfrak{PDer}\left[#5\right]}}{\IfNoValueTF{#5}{OCCHIO!1}{\mathfrak{PDer}\left[#4,#5\right]}}}{\IfNoValueTF{#4}{OCCHIO!2}{\IfNoValueTF{#5}{OCCHIO!3}{\mathfrak{PDer}\left[#3,#4,#5\right]}}}}
	{\IfNoValueTF{#3}{\IfNoValueTF{#4}{\IfNoValueTF{#5}{\mathfrak{PDer}}{\mathcal{N}_{#5}}}{\IfNoValueTF{#5}{OCCHIO!6}{\mathfrak{PDer}_{#4,#5}}}}{\IfNoValueTF{#4}{OCCHIO!7}{\IfNoValueTF{#5}{OCCHIO!8}{\mathfrak{PDer}_{#3,#4,#5}}}}}}	
{\IfNoValueTF{#2}{OCCHIO!9}{\IfNoValueTF{#3}{\IfNoValueTF{#4}{\IfNoValueTF{#5}{OCCHIO!10}{\mathfrak{PDer}^{#1}_{#5}}}{\IfNoValueTF{#5}{OCCHIO!11}{\mathfrak{PDer}^{#1}_{#4,#5}}}}{\IfNoValueTF{#4}{OCCHIO!12}{\IfNoValueTF{#5}{OCCHIO!13}{\mathfrak{PDer}^{#1}_{#3,#4,#5}}}}}}
} 
\NewDocumentCommand\genprederspacesm{d// d-- d[] d() d++}{
\IfNoValueTF{#1}{
\IfNoValueTF{#2}{\IfNoValueTF{#3}{\IfNoValueTF{#4}{\IfNoValueTF{#5}{\mathfrak{SPDer}}{\mathfrak{SPDer}\left[#5\right]}}{\IfNoValueTF{#5}{OCCHIO!1}{\mathfrak{SPDer}\left[#4,#5\right]}}}{\IfNoValueTF{#4}{OCCHIO!2}{\IfNoValueTF{#5}{OCCHIO!3}{\mathfrak{SPDer}\left[#3,#4,#5\right]}}}}
	{\IfNoValueTF{#3}{\IfNoValueTF{#4}{\IfNoValueTF{#5}{\mathfrak{SPDer}}{\mathcal{N}_{#5}}}{\IfNoValueTF{#5}{OCCHIO!6}{\mathfrak{SPDer}_{#4,#5}}}}{\IfNoValueTF{#4}{OCCHIO!7}{\IfNoValueTF{#5}{OCCHIO!8}{\mathfrak{SPDer}_{#3,#4,#5}}}}}}	
{\IfNoValueTF{#2}{OCCHIO!9}{\IfNoValueTF{#3}{\IfNoValueTF{#4}{\IfNoValueTF{#5}{OCCHIO!10}{\mathfrak{SPDer}^{#1}_{#5}}}{\IfNoValueTF{#5}{OCCHIO!11}{\mathfrak{SPDer}^{#1}_{#4,#5}}}}{\IfNoValueTF{#4}{OCCHIO!12}{\IfNoValueTF{#5}{OCCHIO!13}{\mathfrak{SPDer}^{#1}_{#3,#4,#5}}}}}}
}
\NewDocumentCommand\Ggenprederspacesm{d[] d() d++}{%
  \IfNoValueTF{#1}{\IfNoValueTF{#2}{\IfNoValueTF{#3}{\mathfrak{GSPDer}}{\mathfrak{GSPDer}\left[#3\right]}}{\IfNoValueTF{#3}{OCCHIO!}{\mathfrak{GSPDer}\left[#2,#3\right]}}}
	{\IfNoValueTF{#2}{OCCHIO!}{\IfNoValueTF{#3}{OCCHIO!}{\mathfrak{GSPDer}\left[#1,#2,#3\right]}}}	
}
\NewDocumentCommand\evspd{d[] d() d++}{%
  \IfNoValueTF{#1}{\IfNoValueTF{#2}{\IfNoValueTF{#3}{\mathfrak{ep}}{\mathfrak{ep}_{#3}}}{\IfNoValueTF{#3}{OCCHIO!}{\mathfrak{ep}_{#2,#3}}}}
	{\IfNoValueTF{#2}{OCCHIO!}{\IfNoValueTF{#3}{OCCHIO!}{\mathfrak{ep}_{#1,#2,#3}}}}	
}
\NewDocumentCommand\dirsumgenpderspace{d// d-- d[] d() d++}{
\IfNoValueTF{#1}{
	\IfNoValueTF{#2}{\IfNoValueTF{#3}{
	\IfNoValueTF{#4}{\IfNoValueTF{#5}{\mathbf{M}}{\mathbf{M}\left[#5\right]}}{\IfNoValueTF{#5}{OCCHIO!1}{\mathbf{M}\left[#4,#5\right]}}}
	{\IfNoValueTF{#4}{OCCHIO!2}{\IfNoValueTF{#5}{OCCHIO!3}{\mathbf{M}\left[#3,#4,#5\right]}}}}
	{\IfNoValueTF{#3}{
	\IfNoValueTF{#4}{\IfNoValueTF{#5}{\mathbf{M}}{\mathbf{M}_{#5}}}{\IfNoValueTF{#5}{OCCHIO!4}{\mathbf{M}_{#4   ,#5}}}}
	{\IfNoValueTF{#4}{OCCHIO!5}{\IfNoValueTF{#5}{OCCHIO!6}{\mathbf{M}_{#3,#4,#5}}}}}
	}
	{ \IfNoValueTF{#2}{OCCHIO!7}
	{\IfNoValueTF{#3}{
	\IfNoValueTF{#4}{\IfNoValueTF{#5}{OCCHIO!11}{\mathbf{M}^{#1}_{#5}}}{\IfNoValueTF{#5}{OCCHIO!8}{\mathbf{M}^{#1}_{#4,#5}}}}
	{\IfNoValueTF{#4}{OCCHIO!9}{\IfNoValueTF{#5}{OCCHIO!10}{\mathbf{M}^{#1}_{#3,#4,#5}}}}}}
} 
\NewDocumentCommand\dirsumgenprediff{d// d-- d[] d() d++}{
\IfNoValueTF{#1}{
	\IfNoValueTF{#2}{\IfNoValueTF{#3}{
	\IfNoValueTF{#4}{\IfNoValueTF{#5}{\mathbf{pd}}{\mathbf{pd}\left[#5\right]}}{\IfNoValueTF{#5}{OCCHIO!1}{\mathbf{pd}\left[#4,#5\right]}}}
	{\IfNoValueTF{#4}{OCCHIO!2}{\IfNoValueTF{#5}{OCCHIO!3}{\mathbf{pd}\left[#3,#4,#5\right]}}}}
	{\IfNoValueTF{#3}{
	\IfNoValueTF{#4}{\IfNoValueTF{#5}{\mathbf{pd}}{\mathbf{pd}_{#5}}}{\IfNoValueTF{#5}{OCCHIO!4}{\mathbf{pd}_{#4   ,#5}}}}
	{\IfNoValueTF{#4}{OCCHIO!5}{\IfNoValueTF{#5}{OCCHIO!6}{\mathbf{pd}_{#3,#4,#5}}}}}}
	{ \IfNoValueTF{#2}{OCCHIO!7}	
	{	\IfNoValueTF{#3}{
	\IfNoValueTF{#4}{\IfNoValueTF{#5}{\mathbf{pd}^{#1}}{\mathbf{pd}^{#1}_{#5}}}{\IfNoValueTF{#5}{OCCHIO!8}{\mathbf{pd}^{#1}_{#4   ,#5}}}}
	{\IfNoValueTF{#4}{OCCHIO!9}{\IfNoValueTF{#5}{OCCHIO!10}{\mathbf{pd}^{#1}_{#3,#4,#5}}}}}}
}
\NewDocumentCommand\nulldirsumgenpderspace{d// d-- d[] d() d++}{
\IfNoValueTF{#1}{
	\IfNoValueTF{#2}{\IfNoValueTF{#3}{
	\IfNoValueTF{#4}{\IfNoValueTF{#5}{\mathcal{M}}{\mathcal{M}\left[#5\right]}}{\IfNoValueTF{#5}{OCCHIO!1}{\mathcal{M}\left[#4,#5\right]}}}
	{\IfNoValueTF{#4}{OCCHIO!2}{\IfNoValueTF{#5}{OCCHIO!3}{\mathcal{M}\left[#3,#4,#5\right]}}}}
	{\IfNoValueTF{#3}{
	\IfNoValueTF{#4}{\IfNoValueTF{#5}{\mathcal{M}}{\mathcal{M}_{#5}}}{\IfNoValueTF{#5}{OCCHIO!4}{\mathcal{M}_{#4   ,#5}}}}
	{\IfNoValueTF{#4}{OCCHIO!5}{\IfNoValueTF{#5}{OCCHIO!6}{\mathcal{M}_{#3,#4,#5}}}}}
	}
	{ \IfNoValueTF{#2}{OCCHIO!7}
	{\IfNoValueTF{#3}{
	\IfNoValueTF{#4}{\IfNoValueTF{#5}{OCCHIO!11}{\mathcal{M}^{#1}_{#5}}}{\IfNoValueTF{#5}{OCCHIO!8}{\mathcal{M}^{#1}_{#4,#5}}}}
	{\IfNoValueTF{#4}{OCCHIO!9}{\IfNoValueTF{#5}{OCCHIO!10}{\mathcal{M}^{#1}_{#3,#4,#5}}}}}}
}     
\NewDocumentCommand\quotdirsumuno{d// d-- d[] d() d++}{
\IfNoValueTF{#1}{
	\IfNoValueTF{#2}{\IfNoValueTF{#3}{
	\IfNoValueTF{#4}{\IfNoValueTF{#5}{\mathsf{q}}{\mathsf{q}\left[#5\right]}}{\IfNoValueTF{#5}{OCCHIO!1}{\mathsf{q}\left[#4,#5\right]}}}
	{\IfNoValueTF{#4}{OCCHIO!2}{\IfNoValueTF{#5}{OCCHIO!3}{\mathsf{q}\left[#3,#4,#5\right]}}}}	
	{\IfNoValueTF{#3}{
	\IfNoValueTF{#4}{\IfNoValueTF{#5}{OCCHIO!4}{\mathsf{q}_{#5}}}{\IfNoValueTF{#5}{OCCHIO!5}{\mathsf{q}_{#4,#5}}}}
	{\IfNoValueTF{#4}{OCCHIO!5}{\IfNoValueTF{#5}{OCCHIO!6}{\mathsf{q}_{#3,#4,#5}}}}}}	
	{ \IfNoValueTF{#2}{OCCHIO!7}{\IfNoValueTF{#3}{
	\IfNoValueTF{#4}{\IfNoValueTF{#5}{OCCHIO!8}{\mathsf{q}^{#1}_{#5}}}{\IfNoValueTF{#5}{OCCHIO!9}{\mathsf{q}^{#1}_{#4,#5}}}}
	{\IfNoValueTF{#4}{OCCHIO!10}{\IfNoValueTF{#5}{OCCHIO!11}{\mathsf{q}^{#1}_{#3,#4,#5}}}}}}
}
\NewDocumentCommand\quotdirsumdue{d// d-- d[] d() d++}{%
 \IfNoValueTF{#1}{
	\IfNoValueTF{#2}{\IfNoValueTF{#3}{
	\IfNoValueTF{#4}{\IfNoValueTF{#5}{\mathsf{p}}{\mathsf{p}\left[#5\right]}}{\IfNoValueTF{#5}{OCCHIO!1}{\mathsf{p}\left[#4,#5\right]}}}
	{\IfNoValueTF{#4}{OCCHIO!2}{\IfNoValueTF{#5}{OCCHIO!3}{\mathsf{p}\left[#3,#4,#5\right]}}}}	
	{\IfNoValueTF{#3}{
	\IfNoValueTF{#4}{\IfNoValueTF{#5}{OCCHIO!4}{\mathsf{p}_{#5}}}{\IfNoValueTF{#5}{OCCHIO!5}{\mathsf{p}_{#4,#5}}}}
	{\IfNoValueTF{#4}{OCCHIO!5}{\IfNoValueTF{#5}{OCCHIO!6}{\mathsf{p}_{#3,#4,#5}}}}}}	
	{ \IfNoValueTF{#2}{OCCHIO!7}{\IfNoValueTF{#3}{
	\IfNoValueTF{#4}{\IfNoValueTF{#5}{OCCHIO!8}{\mathsf{p}^{#1}_{#5}}}{\IfNoValueTF{#5}{OCCHIO!9}{\mathsf{p}^{#1}_{#4,#5}}}}
	{\IfNoValueTF{#4}{OCCHIO!10}{\IfNoValueTF{#5}{OCCHIO!11}{\mathsf{p}^{#1}_{#3,#4,#5}}}}}}
}	
\NewDocumentCommand\nullprederunodue{d// d-- d[] d() d++}{
  \IfNoValueTF{#1}{
	\IfNoValueTF{#2}{\IfNoValueTF{#3}{
	\IfNoValueTF{#4}{\IfNoValueTF{#5}{\mathcal{O}}{\mathcal{O}\left[#5\right]}}{\IfNoValueTF{#5}{OCCHIO!1}{\mathcal{O}\left[#4,#5\right]}}}
	{\IfNoValueTF{#4}{OCCHIO!2}{\IfNoValueTF{#5}{OCCHIO!3}{\mathcal{O}\left[#3,#4,#5\right]}}}}	
	{\IfNoValueTF{#3}{
	\IfNoValueTF{#4}{\IfNoValueTF{#5}{OCCHIO!}{\mathcal{O}_{#5}}}{\IfNoValueTF{#5}{OCCHIO!4}{\mathcal{O}_{#4,#5}}}}
	{\IfNoValueTF{#4}{OCCHIO!5}{\IfNoValueTF{#5}{OCCHIO!6}{\mathcal{O}_{#3,#4,#5}}}}}}	
	{ \IfNoValueTF{#2}{OCCHIO!7}{	
	\IfNoValueTF{#3}{
	\IfNoValueTF{#4}{\IfNoValueTF{#5}{OCCHIO!8}{\mathcal{O}^{#1}_{#5}}}{\IfNoValueTF{#5}{OCCHIO!9}{\mathcal{O}^{#1}_{#4,#5}}}}
	{\IfNoValueTF{#4}{OCCHIO!10}{\IfNoValueTF{#5}{OCCHIO!11}{\mathcal{O}^{#1}_{#3,#4,#5}}}}}}
} 
\NewDocumentCommand\sumgenpreder{d[] d() d++}{%
  \IfNoValueTF{#1}{\IfNoValueTF{#2}{\IfNoValueTF{#3}{+}{+\left[#3\right]}}{\IfNoValueTF{#3}{OCCHIO!}{+\left[#2,#3\right]}}}
	{\IfNoValueTF{#2}{OCCHIO!}{\IfNoValueTF{#3}{OCCHIO!}{+\left[#1,#2,#3\right]}}}
} 
\NewDocumentCommand\bigsumgenpreder{d[] d() d++}{%
  \IfNoValueTF{#1}{\IfNoValueTF{#2}{\IfNoValueTF{#3}{\scalerel*{\sumgenpreder}{\sum}}{\scalerel*{\sumgenpreder}{\sum}\left[#3\right]}}{\IfNoValueTF{#3}{OCCHIO!}{\scalerel*{\sumgenpreder}{\sum}\left[#2,#3\right]}}}
	{\IfNoValueTF{#2}{OCCHIO!}{\IfNoValueTF{#3}{OCCHIO!}{\scalerel*{\sumgenpreder}{\sum}\left[#1,#2,#3\right]}}}
} 
\NewDocumentCommand\scalpgenpreder{d[] d() d++}{%
  \IfNoValueTF{#1}{\IfNoValueTF{#2}{\IfNoValueTF{#3}{\ast}{\ast\left[#3\right]}}{\IfNoValueTF{#3}{OCCHIO!}{\ast\left[#2,#3\right]}}}
	{\IfNoValueTF{#2}{OCCHIO!}{\IfNoValueTF{#3}{OCCHIO!}{\ast\left[#1,#2,#3\right]}}}
}
\NewDocumentCommand\genprediff{d// d-- d[] d() d++}{
\IfNoValueTF{#1}{
	\IfNoValueTF{#2}{\IfNoValueTF{#3}{
	\IfNoValueTF{#4}{\IfNoValueTF{#5}{\mathfrak{pd}}{\mathfrak{pd}\left[#5\right]}}{\IfNoValueTF{#5}{OCCHIO!1}{\mathfrak{pd}\left[#4,#5\right]}}}
	{\IfNoValueTF{#4}{OCCHIO!2}{\IfNoValueTF{#5}{OCCHIO!3}{\mathfrak{pd}\left[#3,#4,#5\right]}}}}
	{\IfNoValueTF{#3}{
	\IfNoValueTF{#4}{\IfNoValueTF{#5}{\mathfrak{pd}}{\mathfrak{pd}_{#5}}}{\IfNoValueTF{#5}{OCCHIO!4}{\mathfrak{pd}_{#4   ,#5}}}}
	{\IfNoValueTF{#4}{OCCHIO!5}{\IfNoValueTF{#5}{OCCHIO!6}{\mathfrak{pd}_{#3,#4,#5}}}}}}
	{ \IfNoValueTF{#2}{OCCHIO!7}	
	{	\IfNoValueTF{#3}{
	\IfNoValueTF{#4}{\IfNoValueTF{#5}{\mathfrak{pd}^{#1}}{\mathfrak{pd}^{#1}_{#5}}}{\IfNoValueTF{#5}{OCCHIO!8}{\mathfrak{pd}^{#1}_{#4   ,#5}}}}
	{\IfNoValueTF{#4}{OCCHIO!9}{\IfNoValueTF{#5}{OCCHIO!10}{\mathfrak{pd}^{#1}_{#3,#4,#5}}}}}}
} 
\NewDocumentCommand\retrofn{d-- d[] d() d++}{
 \IfNoValueTF{#1}{\IfNoValueTF{#2}{OCCHIO!1}{\IfNoValueTF{#3}{OCCHIO!2}{\IfNoValueTF{#4}{OCCHIO!3}{\mathfrak{rpd}_{#2,#3,#4}}}}}
                 {\IfNoValueTF{#2}{OCCHIO!4}{\IfNoValueTF{#3}{OCCHIO!5}{\IfNoValueTF{#4}{OCCHIO!6}{\mathfrak{rpd}_{#2,#3,#4}\left(#1\right)}}}}
}
\newcommand{\pmdspreder}{\mu}

\newcommand{\pmpreder}{\zeta}
\newcommand{\pmprederdue}{\varrho}
\newcommand{\pmpredertre}{\gamma}
\NewDocumentCommand\adsetcat{d() d++}{%
  \IfNoValueTF{#1}{\IfNoValueTF{#2}{\mbox{\boldmath$\Omega$}}{\mbox{\boldmath$\Omega$}_{#2}}}{\IfNoValueTF{#2}{OCCHIO!}{\mbox{\boldmath$\Omega$}_{#1,#2}}}
} 
\NewDocumentCommand\genpredersoc{d() d++}{%
  \IfNoValueTF{#1}{\IfNoValueTF{#2}{\mbox{\boldmath$\omega$}}{\mbox{\boldmath$\omega$}_{#2}}}{\IfNoValueTF{#2}{OCCHIO!}{\mbox{\boldmath$\omega$}_{#1,#2}}}
}  
\NewDocumentCommand\prederFunc{d[] d() d++}{
  \IfNoValueTF{#1}{\IfNoValueTF{#2}{\IfNoValueTF{#3}{\mathfrak{PrDe}}{\mathfrak{PrDe}_{#3}}}{\IfNoValueTF{#3}{OCCHIO!}{\mathfrak{PrDe}_{#2,#3}}}}
	{\IfNoValueTF{#2}{OCCHIO!}{\IfNoValueTF{#3}{OCCHIO!}{\mathfrak{PrDe}_{#1,#2,#3}}}}
}
\NewDocumentCommand\inclcol{d[] d() d++}{%
  \IfNoValueTF{#1}{\IfNoValueTF{#2}{\IfNoValueTF{#3}{\theta}{\theta_{#3}}}{\IfNoValueTF{#3}{OCCHIO!}{\theta_{#2,#3}}}}
	{\IfNoValueTF{#2}{OCCHIO!}{\IfNoValueTF{#3}{OCCHIO!}{\theta_{#1,#2,#3}}}}
} 
\NewDocumentCommand\genprediffNT{d[] d() d++}{%
  \IfNoValueTF{#1}{\IfNoValueTF{#2}{\IfNoValueTF{#3}{\genprediff}{\genprediff_{#3}}}{\IfNoValueTF{#3}{OCCHIO!}{\genprediff_{#2,#3}}}}
	{\IfNoValueTF{#2}{OCCHIO!}{\IfNoValueTF{#3}{OCCHIO!}{\genprediff_{#1,#2,#3}}}}
}
\NewDocumentCommand\genpreinclNT{d<< d>> d() d++}{%
 \IfNoValueTF{#1}{OCCHIO!}{\IfNoValueTF{#2}{OCCHIO!}{\IfNoValueTF{#3}{OCCHIO!}{\IfNoValueTF{#4}{OCCHIO!}{\mathfrak{pi}_{#1,#2,#3,#4}}}}}
 }
\NewDocumentCommand\genpreincl{d<< d>> d() d++}{%
 \IfNoValueTF{#1}{OCCHIO!}{\IfNoValueTF{#2}{OCCHIO!}{\IfNoValueTF{#3}{OCCHIO!}{\IfNoValueTF{#4}{OCCHIO!}{\mathfrak{pi}\left[#1,#2,#3,#4\right]}}}}
 }
\NewDocumentCommand\dirsumgenpreincl{d<< d>> d() d++}{%
 \IfNoValueTF{#1}{OCCHIO!}{\IfNoValueTF{#2}{OCCHIO!}{\IfNoValueTF{#3}{OCCHIO!}{\IfNoValueTF{#4}{OCCHIO!}{\mathbf{pi}\left[#1,#2,#3,#4\right]}}}}
 }

\NewDocumentCommand\generder{d[]}{ 
\IfNoValueTF{#1}{\Game}{\langle #1 \rangle}
}

\NewDocumentCommand\generderspace{d// d-- d[] d() d++}{%
  \IfNoValueTF{#1}{
	\IfNoValueTF{#2}{\IfNoValueTF{#3}{
	\IfNoValueTF{#4}{\IfNoValueTF{#5}{\mathfrak{D}}{\mathfrak{D}\left[#5\right]}}{\IfNoValueTF{#5}{OCCHIO!1}{\mathfrak{D}\left[#4,#5\right]}}}
	{\IfNoValueTF{#4}{OCCHIO!2}{\IfNoValueTF{#5}{OCCHIO!3}{\mathfrak{D}\left[#3,#4,#5\right]}}}}
	{\IfNoValueTF{#3}{
	\IfNoValueTF{#4}{\IfNoValueTF{#5}{\mathfrak{D}}{\mathfrak{D}_{#5}}}{\IfNoValueTF{#5}{OCCHIO!4}{\mathfrak{D}_{#4   ,#5}}}}
	{\IfNoValueTF{#4}{OCCHIO!5}{\IfNoValueTF{#5}{OCCHIO!6}{\mathfrak{D}_{#3,#4,#5}}}}}
	}
	{ \IfNoValueTF{#2}{OCCHIO!7}
	{\IfNoValueTF{#3}{
	\IfNoValueTF{#4}{\IfNoValueTF{#5}{OCCHIO!11}{\mathfrak{D}^{#1}_{#5}}}{\IfNoValueTF{#5}{OCCHIO!8}{\mathfrak{D}^{#1}_{#4,#5}}}}
	{\IfNoValueTF{#4}{OCCHIO!9}{\IfNoValueTF{#5}{OCCHIO!10}{\mathfrak{D}^{#1}_{#3,#4,#5}}}}}}
}
\NewDocumentCommand\generderspacesf{d[] d|| d() d++}{
\IfNoValueTF{#1}{\widehat{\mathfrak{D}}_{#2,#3,#4}}{\widehat{\mathfrak{D}}_{#1,#2,#3,#4}}
}
\NewDocumentCommand\sumgenerder{d-- d[] d() d++}{
  \IfNoValueTF{#1}{\IfNoValueTF{#2}{\IfNoValueTF{#3}{\IfNoValueTF{#4}{+}{+\left[#4\right]}}{\IfNoValueTF{#4}{OCCHIO!1}{+\left[#3,#4\right]}}}{\IfNoValueTF{#3}{OCCHIO!2}{\IfNoValueTF{#4}{+\left[#2,#3\right]}{+\left[#2,#3,#4\right]}}}}
	{\IfNoValueTF{#2}{\IfNoValueTF{#3}{\IfNoValueTF{#4}{+\left[#1\right]}{+\left[#1,#4\right]}}{\IfNoValueTF{#4}{OCCHIO!3}{+\left[#1,#3,#4\right]}}}{\IfNoValueTF{#3}{OCCHIO!4}{\IfNoValueTF{#4}{+\left[#1,#2,#3\right]}{+\left[#1,#2,#3,#4\right]}}}}
}

\NewDocumentCommand\scalpgenerder{d-- d[] d() d++}{%
\IfNoValueTF{#1}{\IfNoValueTF{#2}{\IfNoValueTF{#3}{\IfNoValueTF{#4}{\ast}{\ast\left[#4\right]}}{\IfNoValueTF{#4}{OCCHIO!1}{\ast\left[#3,#4\right]}}}{\IfNoValueTF{#3}{OCCHIO!2}{\IfNoValueTF{#4}{\ast\left[#2,#3\right]}{\ast\left[#2,#3,#4\right]}}}}
	{\IfNoValueTF{#2}{\IfNoValueTF{#3}{\IfNoValueTF{#4}{\ast\left[#1\right]}{\ast\left[#1,#4\right]}}{\IfNoValueTF{#4}{OCCHIO!3}{\ast\left[#1,#3,#4\right]}}}{\IfNoValueTF{#3}{OCCHIO!4}{\IfNoValueTF{#4}{\ast\left[#1,#2,#3\right]}{\ast\left[#1,#2,#3,#4\right]}}}} 
}
\NewDocumentCommand\mmscalpgenerder{d-- d[] d() d++}{%
\IfNoValueTF{#1}{\IfNoValueTF{#2}{\IfNoValueTF{#3}{\IfNoValueTF{#4}{\widehat{\ast}}{\widehat{\ast}\left[#4\right]}}{\IfNoValueTF{#4}{OCCHIO!1}{\widehat{\ast}\left[#3,#4\right]}}}{\IfNoValueTF{#3}{OCCHIO!2}{\IfNoValueTF{#4}{\widehat{\ast}\left[#2,#3\right]}{\widehat{\ast}\left[#2,#3,#4\right]}}}}
	{\IfNoValueTF{#2}{\IfNoValueTF{#3}{\IfNoValueTF{#4}{\widehat{\ast}\left[#1\right]}{\widehat{\ast}\left[#1,#4\right]}}{\IfNoValueTF{#4}{OCCHIO!3}{\widehat{\ast}\left[#1,#3,#4\right]}}}{\IfNoValueTF{#3}{OCCHIO!4}{\IfNoValueTF{#4}{\widehat{\ast}\left[#1,#2,#3\right]}{\widehat{\ast}\left[#1,#2,#3,#4\right]}}}} 
}
\NewDocumentCommand\sqmmscalpgenerder{d-- d[] d() d++}{%
\IfNoValueTF{#1}{\IfNoValueTF{#2}{\IfNoValueTF{#3}{\IfNoValueTF{#4}{\star}{\star\left[#4\right]}}{\IfNoValueTF{#4}{OCCHIO!1}{\star\left[#3,#4\right]}}}{\IfNoValueTF{#3}{OCCHIO!2}{\IfNoValueTF{#4}{\star\left[#2,#3\right]}{\star\left[#2,#3,#4\right]}}}}
	{\IfNoValueTF{#2}{\IfNoValueTF{#3}{\IfNoValueTF{#4}{\star\left[#1\right]}{\star\left[#1,#4\right]}}{\IfNoValueTF{#4}{OCCHIO!3}{\star\left[#1,#3,#4\right]}}}{\IfNoValueTF{#3}{OCCHIO!4}{\IfNoValueTF{#4}{\star\left[#1,#2,#3\right]}{\star\left[#1,#2,#3,#4\right]}}}} 
}

\NewDocumentCommand\gendiff{d// d-- d[] d() d++ d)) d]]}{%
\IfNoValueTF{#7}{
\IfNoValueTF{#6}{
\IfNoValueTF{#1}{
	\IfNoValueTF{#2}{\IfNoValueTF{#3}{
	\IfNoValueTF{#4}{\IfNoValueTF{#5}{\mathbf{d}}{\mathbf{d}\left[#5\right]}}{\IfNoValueTF{#5}{OCCHIO!1}{\mathbf{d}\left[#4,#5\right]}}}
	{\IfNoValueTF{#4}{OCCHIO!2}{\IfNoValueTF{#5}{OCCHIO!3}{\mathbf{d}\left[#3,#4,#5\right]}}}}
	{\IfNoValueTF{#3}{
	\IfNoValueTF{#4}{\IfNoValueTF{#5}{\mathbf{d}}{\mathbf{d}_{#5}}}{\IfNoValueTF{#5}{OCCHIO!4}{\mathbf{d}_{#4   ,#5}}}}
	{\IfNoValueTF{#4}{OCCHIO!5}{\IfNoValueTF{#5}{OCCHIO!6}{\mathbf{d}_{#3,#4,#5}}}}}}
	{ \IfNoValueTF{#2}{OCCHIO!7}	
	{	\IfNoValueTF{#3}{
	\IfNoValueTF{#4}{\IfNoValueTF{#5}{\mathbf{d}^{#1}}{\mathbf{d}^{#1}_{#5}}}{\IfNoValueTF{#5}{OCCHIO!8}{\mathbf{d}^{#1}_{#4   ,#5}}}}
	{\IfNoValueTF{#4}{OCCHIO!9}{\IfNoValueTF{#5}{OCCHIO!10}{\mathbf{d}^{#1}_{#3,#4,#5}}}}}}
}{OCCHIO!11}}{
\IfNoValueTF{#6}{OCCHIO!12}{
\IfNoValueTF{#1}{
	\IfNoValueTF{#2}{\IfNoValueTF{#3}{
	\IfNoValueTF{#4}{\IfNoValueTF{#5}{\mathbf{d}^{#7}_{#6}}{\mathbf{d}\left[#5\right]^{#7}_{#6}}}{\IfNoValueTF{#5}{OCCHIO!1}{\mathbf{d}\left[#4,#5\right]^{#7}_{#6}}}}
	{\IfNoValueTF{#4}{OCCHIO!2}{\IfNoValueTF{#5}{OCCHIO!3}{\mathbf{d}\left[#3,#4,#5\right]^{#7}_{#6}}}}}
	{\IfNoValueTF{#3}{
	\IfNoValueTF{#4}{\IfNoValueTF{#5}{\mathbf{d}^{#7}_{#6}}{\left(\mathbf{d}_{#5}\right)^{#7}_{#6}}}{\IfNoValueTF{#5}{OCCHIO!4}{\left(\mathbf{d}_{#4   ,#5}\right)^{#7}_{#6}}}}
	{\IfNoValueTF{#4}{OCCHIO!5}{\IfNoValueTF{#5}{OCCHIO!6}{\left(\mathbf{d}_{#3,#4,#5}\right)^{#7}_{#6}}}}}}
	{ \IfNoValueTF{#2}{OCCHIO!7}	
	{	\IfNoValueTF{#3}{
	\IfNoValueTF{#4}{\IfNoValueTF{#5}{\left(\mathbf{d}^{#1}\right)^{#7}_{#6}}{\left(\mathbf{d}^{#1}_{#5}\right)^{#7}_{#6}}}{\IfNoValueTF{#5}{OCCHIO!8}{\left(\mathbf{d}^{#1}_{#4   ,#5}\right)^{#7}_{#6}}}}
	{\IfNoValueTF{#4}{OCCHIO!9}{\IfNoValueTF{#5}{OCCHIO!10}{\left(\mathbf{d}^{#1}_{#3,#4,#5}\right)^{#7}_{#6}}}}}}
}
}
}
	
\NewDocumentCommand\generderspacesm{d// d-- d[] d() d++}{ 
\IfNoValueTF{#1}{
	\IfNoValueTF{#2}{\IfNoValueTF{#3}{
	\IfNoValueTF{#4}{\IfNoValueTF{#5}{\mathfrak{SD}}{\mathfrak{SD}\left[#5\right]}}{\IfNoValueTF{#5}{OCCHIO!1}{\mathfrak{SD}\left[#4,#5\right]}}}
	{\IfNoValueTF{#4}{OCCHIO!2}{\IfNoValueTF{#5}{OCCHIO!3}{\mathfrak{SD}\left[#3,#4,#5\right]}}}}
	{\IfNoValueTF{#3}{
	\IfNoValueTF{#4}{\IfNoValueTF{#5}{\mathfrak{SD}}{\mathfrak{SD}_{#5}}}{\IfNoValueTF{#5}{OCCHIO!4}{\mathfrak{SD}_{#4   ,#5}}}}
	{\IfNoValueTF{#4}{OCCHIO!5}{\IfNoValueTF{#5}{OCCHIO!6}{\mathfrak{SD}_{#3,#4,#5}}}}}
	}
	{ \IfNoValueTF{#2}{OCCHIO!7}
	{\IfNoValueTF{#3}{
	\IfNoValueTF{#4}{\IfNoValueTF{#5}{OCCHIO!8}{\mathfrak{SD}^{#1}_{#5}}}{\IfNoValueTF{#5}{OCCHIO!9}{\mathfrak{SD}^{#1}_{#4,#5}}}}
	{\IfNoValueTF{#4}{OCCHIO!10}{\IfNoValueTF{#5}{OCCHIO!11}{\mathfrak{SD}^{#1}_{#3,#4,#5}}}}}}
}
\NewDocumentCommand\generderspacesmsf{d[] d|| d() d++}{
\IfNoValueTF{#1}{\widehat{\mathfrak{SD}}_{#2,#3,#4}}{\widehat{\mathfrak{S}}_{#1,#2,#3,#4}}
}
\NewDocumentCommand\inclgdsmgd{d// d-- d[] d() d++}{ 
  \IfNoValueTF{#1}{
	\IfNoValueTF{#2}{\IfNoValueTF{#3}{
	\IfNoValueTF{#4}{\IfNoValueTF{#5}{s\iota}{s\iota\left[#5\right]}}{\IfNoValueTF{#5}{OCCHIO!1}{s\iota\left[#4,#5\right]}}}
	{\IfNoValueTF{#4}{OCCHIO!2}{\IfNoValueTF{#5}{OCCHIO!3}{s\iota\left[#3,#4,#5\right]}}}}
	{\IfNoValueTF{#3}{
	\IfNoValueTF{#4}{\IfNoValueTF{#5}{s\iota}{s\iota_{#5}}}{\IfNoValueTF{#5}{OCCHIO!4}{s\iota_{#4,#5}}}}
	{\IfNoValueTF{#4}{OCCHIO!5}{\IfNoValueTF{#5}{OCCHIO!6}{s\iota_{#3,#4,#5}}}}}}
	{ \IfNoValueTF{#2}{OCCHIO!7}
	{\IfNoValueTF{#3}{
	\IfNoValueTF{#4}{\IfNoValueTF{#5}{s\iota^{#1}}{s\iota^{#1}_{#5}}}{\IfNoValueTF{#5}{OCCHIO!8}{s\iota^{#1}_{#4,#5}}}}{\IfNoValueTF{#4}{OCCHIO!9}{\IfNoValueTF{#5}{OCCHIO!10}{s\iota^{#1}_{#3,#4,#5}}}}}}
} 
\NewDocumentCommand\gendiffsm{d// d-- d[] d() d++}{
\IfNoValueTF{#1}{
	\IfNoValueTF{#2}{\IfNoValueTF{#3}{
	\IfNoValueTF{#4}{\IfNoValueTF{#5}{\mathbf{sd}}{\mathbf{sd}\left[#5\right]}}{\IfNoValueTF{#5}{OCCHIO!1}{\mathbf{sd}\left[#4,#5\right]}}}
	{\IfNoValueTF{#4}{OCCHIO!2}{\IfNoValueTF{#5}{OCCHIO!3}{\mathbf{sd}\left[#3,#4,#5\right]}}}}
	{\IfNoValueTF{#3}{
	\IfNoValueTF{#4}{\IfNoValueTF{#5}{\mathbf{sd}}{\mathbf{sd}_{#5}}}{\IfNoValueTF{#5}{OCCHIO!4}{\mathbf{sd}_{#4   ,#5}}}}
	{\IfNoValueTF{#4}{OCCHIO!5}{\IfNoValueTF{#5}{OCCHIO!6}{\mathbf{sd}_{#3,#4,#5}}}}}}
	{ \IfNoValueTF{#2}{OCCHIO!7}	
	{	\IfNoValueTF{#3}{
	\IfNoValueTF{#4}{\IfNoValueTF{#5}{\mathbf{sd}^{#1}}{\mathbf{sd}^{#1}_{#5}}}{\IfNoValueTF{#5}{OCCHIO!8}{\mathbf{sd}^{#1}_{#4   ,#5}}}}
	{\IfNoValueTF{#4}{OCCHIO!9}{\IfNoValueTF{#5}{OCCHIO!10}{\mathbf{sd}^{#1}_{#3,#4,#5}}}}}}
} 
\NewDocumentCommand\kersqA{d// d-- d[] d() d++}{ 
  \IfNoValueTF{#1}{
	\IfNoValueTF{#2}{\IfNoValueTF{#3}{
	\IfNoValueTF{#4}{\IfNoValueTF{#5}{\mathbf{kd}}{\mathbf{kd}\left[#5\right]}}{\IfNoValueTF{#5}{OCCHIO!1}{\mathbf{kd}\left[#4,#5\right]}}}
	{\IfNoValueTF{#4}{OCCHIO!2}{\IfNoValueTF{#5}{OCCHIO!3}{\mathbf{kd}\left[#3,#4,#5\right]}}}}
	{\IfNoValueTF{#3}{
	\IfNoValueTF{#4}{\IfNoValueTF{#5}{\mathbf{kd}}{\mathbf{kd}_{#5}}}{\IfNoValueTF{#5}{OCCHIO!4}{\mathbf{kd}_{#4,#5}}}}
	{\IfNoValueTF{#4}{OCCHIO!5}{\IfNoValueTF{#5}{OCCHIO!6}{\mathbf{kd}_{#3,#4,#5}}}}}}
	{ \IfNoValueTF{#2}{OCCHIO!7}
	{	
	\IfNoValueTF{#3}{
	\IfNoValueTF{#4}{\IfNoValueTF{#5}{\mathbf{kd}^{#1}}{\mathbf{kd}^{#1}_{#5}}}{\IfNoValueTF{#5}{OCCHIO!8}{\mathbf{kd}^{#1}_{#4,#5}}}}
	{\IfNoValueTF{#4}{OCCHIO!9}{\IfNoValueTF{#5}{OCCHIO!10}{\mathbf{kd}^{#1}_{#3,#4,#5}}}}}}
}
\NewDocumentCommand\kergenerderspacesm{d// d-- d[] d() d++}{ 
	\IfNoValueTF{#1}{
	\IfNoValueTF{#2}{\IfNoValueTF{#3}{
	\IfNoValueTF{#4}{\IfNoValueTF{#5}{\mathfrak{KSD}}{\mathfrak{KSD}\left[#5\right]}}{\IfNoValueTF{#5}{OCCHIO!1}{\mathfrak{KSD}\left[#4,#5\right]}}}
	{\IfNoValueTF{#4}{OCCHIO!2}{\IfNoValueTF{#5}{OCCHIO!3}{\mathfrak{KSD}\left[#3,#4,#5\right]}}}}
	{\IfNoValueTF{#3}{
	\IfNoValueTF{#4}{\IfNoValueTF{#5}{\mathfrak{KSD}}{\mathfrak{KSD}_{#5}}}{\IfNoValueTF{#5}{OCCHIO!4}{\mathfrak{KSD}_{#4,#5}}}}
	{\IfNoValueTF{#4}{OCCHIO!5}{\IfNoValueTF{#5}{OCCHIO!6}{\mathfrak{KSD}_{#3,#4,#5}}}}}}
	{ \IfNoValueTF{#2}{OCCHIO!7}	
	{\IfNoValueTF{#3}{
	\IfNoValueTF{#4}{\IfNoValueTF{#5}{OCCHIO!8}{\mathfrak{KSD}^{#1}_{#5}}}{\IfNoValueTF{#5}{OCCHIO!9}{\mathfrak{KSD}^{#1}_{#4,#5}}}}
	{\IfNoValueTF{#4}{OCCHIO!10}{\IfNoValueTF{#5}{OCCHIO!11}{\mathfrak{KSD}^{#1}_{#3,#4,#5}}}}}}
}
\NewDocumentCommand\kgdssmar{d// d-- d[] d() d++}{ 
  \IfNoValueTF{#1}{
	\IfNoValueTF{#2}{\IfNoValueTF{#3}{
	\IfNoValueTF{#4}{\IfNoValueTF{#5}{\rho}{\rho\left[#5\right]}}{\IfNoValueTF{#5}{OCCHIO!1}{\rho\left[#4,#5\right]}}}
	{\IfNoValueTF{#4}{OCCHIO!2}{\IfNoValueTF{#5}{OCCHIO!3}{\rho\left[#3,#4,#5\right]}}}}
	{	\IfNoValueTF{#3}{
	\IfNoValueTF{#4}{\IfNoValueTF{#5}{\rho}{\rho_{#5}}}{\IfNoValueTF{#5}{OCCHIO!4}{\rho_{#4,#5}}}}
	{\IfNoValueTF{#4}{OCCHIO!5}{\IfNoValueTF{#5}{OCCHIO!6}{\rho_{#3,#4,#5}}}}}}
	{ \IfNoValueTF{#2}{OCCHIO!7}
	{\IfNoValueTF{#3}{
	\IfNoValueTF{#4}{\IfNoValueTF{#5}{\rho^{#1}}{\rho^{#1}_{#5}}}{\IfNoValueTF{#5}{OCCHIO!8}{\rho^{#1}_{#4,#5}}}}
	{\IfNoValueTF{#4}{OCCHIO!9}{\IfNoValueTF{#5}{OCCHIO!10}{\rho^{#1}_{#3,#4,#5}}}}}} 
}
\NewDocumentCommand\cokergenerderspace{d// d-- d[] d() d++}{
  \IfNoValueTF{#1}{
	\IfNoValueTF{#2}{\IfNoValueTF{#3}{
	\IfNoValueTF{#4}{\IfNoValueTF{#5}{\mathfrak{CD}}{\mathfrak{CD}\left[#5\right]}}{\IfNoValueTF{#5}{OCCHIO!1}{\mathfrak{CD}\left[#4,#5\right]}}}
	{\IfNoValueTF{#4}{OCCHIO!2}{\IfNoValueTF{#5}{OCCHIO!3}{\mathfrak{CD}\left[#3,#4,#5\right]}}}}
	{\IfNoValueTF{#3}{
	\IfNoValueTF{#4}{\IfNoValueTF{#5}{\mathfrak{CD}}{\mathfrak{CD}_{#5}}}{\IfNoValueTF{#5}{OCCHIO!4}{\mathfrak{CD}_{#4   ,#5}}}}
	{\IfNoValueTF{#4}{OCCHIO!5}{\IfNoValueTF{#5}{OCCHIO!6}{\mathfrak{CD}_{#3,#4,#5}}}}}
	}
	{\IfNoValueTF{#2}{OCCHIO!7}	
	{\IfNoValueTF{#3}{
	\IfNoValueTF{#4}{\IfNoValueTF{#5}{\mathfrak{CD}^{#1}}{\mathfrak{CD}^{#1}_{#5}}}{\IfNoValueTF{#5}{OCCHIO!8}{\mathfrak{CD}^{#1}_{#4,#5}}}}
	{\IfNoValueTF{#4}{OCCHIO!9}{\IfNoValueTF{#5}{OCCHIO!10}{\mathfrak{CD}^{#1}_{#3,#4,#5}}}}}}
}
\NewDocumentCommand\ckgdssmar{d// d-- d[] d() d++}{
   \IfNoValueTF{#1}{
	\IfNoValueTF{#2}{\IfNoValueTF{#3}{
	\IfNoValueTF{#4}{\IfNoValueTF{#5}{\epsilon}{\epsilon\left[#5\right]}}{\IfNoValueTF{#5}{OCCHIO!1}{\epsilon\left[#4,#5\right]}}}
	{\IfNoValueTF{#4}{OCCHIO!2}{\IfNoValueTF{#5}{OCCHIO!3}{\epsilon\left[#3,#4,#5\right]}}}}
	{\IfNoValueTF{#3}{
	\IfNoValueTF{#4}{\IfNoValueTF{#5}{\epsilon}{\epsilon_{#5}}}{\IfNoValueTF{#5}{OCCHIO!4}{\epsilon_{#4,#5}}}}
	{\IfNoValueTF{#4}{OCCHIO!5}{\IfNoValueTF{#5}{OCCHIO!6}{\epsilon_{#3,#4,#5}}}}}}
	{ \IfNoValueTF{#2}{OCCHIO!7}		
	{\IfNoValueTF{#3}{
	\IfNoValueTF{#4}{\IfNoValueTF{#5}{\epsilon^{#1}}{\epsilon^{#1}_{#5}}}{\IfNoValueTF{#5}{OCCHIO!8}{\epsilon^{#1}_{#4,#5}}}}
	{\IfNoValueTF{#4}{OCCHIO!9}{\IfNoValueTF{#5}{OCCHIO!10}{\epsilon^{#1}_{#3,#4,#5}}}}}}
}
\NewDocumentCommand\cksqB{d// d-- d[] d() d++}{ 
\IfNoValueTF{#1}{
	\IfNoValueTF{#2}{\IfNoValueTF{#3}{
	\IfNoValueTF{#4}{\IfNoValueTF{#5}{\mathbf{cd}}{\mathbf{cd}\left[#5\right]}}{\IfNoValueTF{#5}{OCCHIO!1}{\mathbf{cd}\left[#4,#5\right]}}}
	{\IfNoValueTF{#4}{OCCHIO!2}{\IfNoValueTF{#5}{OCCHIO!3}{\mathbf{cd}\left[#3,#4,#5\right]}}}}
	{\IfNoValueTF{#3}{
	\IfNoValueTF{#4}{\IfNoValueTF{#5}{\mathbf{cd}}{\mathbf{cd}_{#5}}}{\IfNoValueTF{#5}{OCCHIO!4}{\mathbf{cd}_{#4,#5}}}}
	{\IfNoValueTF{#4}{OCCHIO!5}{\IfNoValueTF{#5}{OCCHIO!6}{\mathbf{cd}_{#3,#4,#5}}}}}}	
	{ \IfNoValueTF{#2}{OCCHIO!7}	
	{\IfNoValueTF{#3}{
	\IfNoValueTF{#4}{\IfNoValueTF{#5}{\mathbf{cd}^{#1}}{\mathbf{cd}^{#1}_{#5}}}{\IfNoValueTF{#5}{OCCHIO!8}{\mathbf{cd}^{#1}_{#4,#5}}}}
	{\IfNoValueTF{#4}{OCCHIO!9}{\IfNoValueTF{#5}{OCCHIO!10}{\mathbf{cd}^{#1}_{#3,#4,#5}}}}}}
} 
\NewDocumentCommand\pmgenerder{d[]}{\mbox{\boldmath$\pmpreder$}}
\NewDocumentCommand\pmgenerderdue{d[]}{\mbox{\boldmath$\pmprederdue$}}
\NewDocumentCommand\pmgenerdertre{d[]}{\mbox{\boldmath$\pmpredertre$}}
\NewDocumentCommand\retezero{d++}{
\IfNoValueTF{#1}{\mathbf{N}}{\mathbf{N}_{#1}}
} 
\NewDocumentCommand\transterezero{d++}{%
  \IfNoValueTF{#1}{\mbox{\boldmath$\eta$}}{\mbox{\boldmath$\eta$}_{#1}}
}
\newcommand{\objrtz}{\mathbf{n}}
\newcommand{\arrtz}{\mathbf{g}}
\NewDocumentCommand\derrtzfun{d-- d[] d() d++}{
 \IfNoValueTF{#1}{\IfNoValueTF{#2}{\IfNoValueTF{#3}{\IfNoValueTF{#4}{\mathbf{F}}{\mathbf{F}_{#4}}}{\IfNoValueTF{#4}{OCCHIO!}{\mathbf{F}_{#3,#4}}}}{\IfNoValueTF{#3}{OCCHIO!}{\IfNoValueTF{#4}{OCCHIO!}{\mathbf{F}_{#2,#3,#4}}}}}
{\IfNoValueTF{#2}{\IfNoValueTF{#3}{\IfNoValueTF{#4}{OCCHIO!}{\mathbf{F}_{#1,#4}}}{\IfNoValueTF{#4}{OCCHIO!}{\mathbf{F}_{#1,#3,#4}}}}{\IfNoValueTF{#3}{OCCHIO!}{\IfNoValueTF{#4}{OCCHIO!}{\mathbf{F}_{#1,#2,#3,#4}}}}}
}
\NewDocumentCommand\quotpreder{d[] d() d++}{ 
  \IfNoValueTF{#1}{\IfNoValueTF{#2}{\IfNoValueTF{#3}{\mathfrak{q}}{\mathfrak{q}_{#3}}}{\IfNoValueTF{#3}{OCCHIO!1}{\mathfrak{q}_{#2,#3}}}}{\IfNoValueTF{#2}{OCCHIO!2}{\IfNoValueTF{#3}{OCCHIO!3}{\mathfrak{q}_{#1,#2,#3}}}}
} 
\NewDocumentCommand\generdiffNT{ d[] d() d++}{%
 \IfNoValueTF{#1}{\IfNoValueTF{#2}{\IfNoValueTF{#3}{\phi}{\phi_{#3}}}{\IfNoValueTF{#3}{OCCHIO!}{\phi_{#2,#3}}}}{\IfNoValueTF{#2}{OCCHIO!}{\IfNoValueTF{#3}{OCCHIO!}{\phi_{#1,#2,#3}}}}
}
\NewDocumentCommand\generinclNT{d<< d>> d() d++}{
  \IfNoValueTF{#1}{OCCHIO!1}{\IfNoValueTF{#2}{OCCHIO!2}{\IfNoValueTF{#3}{OCCHIO!3}{\IfNoValueTF{#4}{OCCHIO!4}{\mathfrak{j}_{#1,#2,#3,#4}}}}}
}
\NewDocumentCommand\generincl{d// d-- d<< d>> d() d++}{
\IfNoValueTF{#1}{
	\IfNoValueTF{#2}{\IfNoValueTF{#3}{
	\IfNoValueTF{#4}{\IfNoValueTF{#5}{\IfNoValueTF{#6}{\mathfrak{i}}{\mathfrak{i}\left[#6\right]}}{\IfNoValueTF{#6}{OCCHIO!1}{\mathfrak{i}\left[#5,#6\right]}}}{OCCHIO!2}}{\IfNoValueTF{#4}{OCCHIO!3}{\IfNoValueTF{#5}{OCCHIO!4}{\IfNoValueTF{#6}{OCCHIO!5}{\mathfrak{i}\left[#3,#4,#5,#6\right]}}}}}
	{\IfNoValueTF{#3}{
	\IfNoValueTF{#4}{\IfNoValueTF{#5}{\IfNoValueTF{#6}{\mathfrak{i}}{\mathfrak{i}_{#6}}}{\IfNoValueTF{#6}{OCCHIO!6}{\mathfrak{i}_{#5,#6}}}}{OCCHIO!7}}{\IfNoValueTF{#4}{OCCHIO!8}{\IfNoValueTF{#5}{OCCHIO!9}{\IfNoValueTF{#6}{OCCHIO!10}{\mathfrak{i}_{#3,#4,#5,#6}}}}}}}		
	{\IfNoValueTF{#2}{OCCHIO!11}
{\IfNoValueTF{#3}{
	\IfNoValueTF{#4}{\IfNoValueTF{#5}{\IfNoValueTF{#6}{\mathfrak{i}^{#1}}{\mathfrak{i}^{#1}_{#6}}}{\IfNoValueTF{#6}{OCCHIO!12}{\mathfrak{i}^{#1}_{#5,#6}}}}{OCCHIO!13}}{\IfNoValueTF{#4}{OCCHIO!14}{\IfNoValueTF{#5}{OCCHIO!15}{\IfNoValueTF{#6}{OCCHIO!16}{\mathfrak{i}^{#1}_{#3,#4,#5,#6}}}}}}}	
}
\NewDocumentCommand\ckgenerincl{d// d-- d<< d>> d() d++}{ %
\IfNoValueTF{#1}{
	\IfNoValueTF{#2}{\IfNoValueTF{#3}{
	\IfNoValueTF{#4}{\IfNoValueTF{#5}{\IfNoValueTF{#6}{\mathfrak{cki}}{\mathfrak{cki}\left[#6\right]}}{\IfNoValueTF{#6}{OCCHIO!}{\mathfrak{cki}\left[#5,#6\right]}}}{OCCHIO!1}}{\IfNoValueTF{#4}{OCCHIO!2}{\IfNoValueTF{#5}{OCCHIO!3}{\IfNoValueTF{#6}{OCCHIO!4}{\mathfrak{cki}\left[#3,#4,#5,#6\right]}}}}}
	{\IfNoValueTF{#3}{
	\IfNoValueTF{#4}{\IfNoValueTF{#5}{\IfNoValueTF{#6}{\mathfrak{i}}{\mathfrak{cki}_{#6}}}{\IfNoValueTF{#6}{OCCHIO!5}{\mathfrak{cki}_{#5,#6}}}}{OCCHIO!6}}{\IfNoValueTF{#4}{OCCHIO!7}{\IfNoValueTF{#5}{OCCHIO!8}{\IfNoValueTF{#6}{OCCHIO!9}{\mathfrak{cki}_{#3,#4,#5,#6}}}}}}}	
	{	\IfNoValueTF{#2}{OCCHIO!10}{\IfNoValueTF{#3}{
	\IfNoValueTF{#4}{\IfNoValueTF{#5}{\IfNoValueTF{#6}{\mathfrak{cki}}{\mathfrak{cki}\left[#6\right]}}{\IfNoValueTF{#6}{OCCHIO!11}{\mathfrak{cki}\left[#5,#6\right]}}}{OCCHIO!12}}{\IfNoValueTF{#4}{OCCHIO!13}{\IfNoValueTF{#5}{OCCHIO!14}{\IfNoValueTF{#6}{OCCHIO!15}{\mathfrak{cki}\left[#3,#4,#5,#6\right]}}}}}}	
}
\NewDocumentCommand\generinclsm{d// d-- d<< d>> d() d++}{%
\IfNoValueTF{#1}{
	\IfNoValueTF{#2}{\IfNoValueTF{#3}{
	\IfNoValueTF{#4}{\IfNoValueTF{#5}{\IfNoValueTF{#6}{\mathfrak{is}}{\mathfrak{is}\left[#6\right]}}{\IfNoValueTF{#6}{OCCHIO!1}{\mathfrak{is}\left[#5,#6\right]}}}{OCCHIO!2}}{\IfNoValueTF{#4}{OCCHIO!3}{\IfNoValueTF{#5}{OCCHIO!4}{\IfNoValueTF{#6}{OCCHIO!5}{\mathfrak{is}\left[#3,#4,#5,#6\right]}}}}}
	{\IfNoValueTF{#3}{
	\IfNoValueTF{#4}{\IfNoValueTF{#5}{\IfNoValueTF{#6}{\mathfrak{is}}{\mathfrak{is}_{#6}}}{\IfNoValueTF{#6}{OCCHIO!6}{\mathfrak{is}_{#5,#6}}}}{OCCHIO!7}}{\IfNoValueTF{#4}{OCCHIO!8}{\IfNoValueTF{#5}{OCCHIO!9}{\IfNoValueTF{#6}{OCCHIO!10}{\mathfrak{is}_{#3,#4,#5,#6}}}}}}}	
	{	\IfNoValueTF{#2}{OCCHIO!11}{	
	\IfNoValueTF{#3}{
	\IfNoValueTF{#4}{\IfNoValueTF{#5}{\IfNoValueTF{#6}{\mathfrak{is}^{#1}}{\mathfrak{is}^{#1}_{#6}}}{\IfNoValueTF{#6}{OCCHIO!11}{\mathfrak{is}^{#1}_{#5,#6}}}}{OCCHIO!12}}{\IfNoValueTF{#4}{OCCHIO!13}{\IfNoValueTF{#5}{OCCHIO!14}{\IfNoValueTF{#6}{OCCHIO!15}{\mathfrak{is}^{#1}_{#3,#4,#5,#6}}}}}}}	
}
\NewDocumentCommand\kergenerinclsm{d// d-- d<< d>> d() d++}{
\IfNoValueTF{#1}{
	\IfNoValueTF{#2}{\IfNoValueTF{#3}{
	\IfNoValueTF{#4}{\IfNoValueTF{#5}{\IfNoValueTF{#6}{\mathfrak{kis}}{\mathfrak{kis}\left[#6\right]}}{\IfNoValueTF{#6}{OCCHIO!1}{\mathfrak{kis}\right[#5,#6\right]}}}{OCCHIO!2}}{\IfNoValueTF{#4}{OCCHIO!3}{\IfNoValueTF{#5}{OCCHIO!4}{\IfNoValueTF{#6}{OCCHIO!5}{\mathfrak{kis}\left[#3,#4,#5,#6\right]}}}}}
	{\IfNoValueTF{#3}{
	\IfNoValueTF{#4}{\IfNoValueTF{#5}{\IfNoValueTF{#6}{\mathfrak{kis}}{\mathfrak{kis}_{#6}}}{\IfNoValueTF{#6}{OCCHIO!6}{\mathfrak{kis}_{#5,#6}}}}{OCCHIO!7}}{\IfNoValueTF{#4}{OCCHIO!8}{\IfNoValueTF{#5}{OCCHIO!9}{\IfNoValueTF{#6}{OCCHIO!10}{\mathfrak{kis}_{#3,#4,#5,#6}}}}}}}
	{	\IfNoValueTF{#2}{OCCHIO!11}{\IfNoValueTF{#3}{
	\IfNoValueTF{#4}{\IfNoValueTF{#5}{\IfNoValueTF{#6}{\mathfrak{kis}^{#1}}{\mathfrak{kis}^{#1}_{#6}}}{\IfNoValueTF{#6}{OCCHIO!6}{\mathfrak{kis}^{#1}_{#5,#6}}}}{OCCHIO!7}}{\IfNoValueTF{#4}{OCCHIO!8}{\IfNoValueTF{#5}{OCCHIO!9}{\IfNoValueTF{#6}{OCCHIO!10}{\mathfrak{kis}^{#1}_{#3,#4,#5,#6}}}}}}}	
}
\NewDocumentCommand\quotsm{d// d-- d[] d() d++}{ % 
\IfNoValueTF{#1}{
	\IfNoValueTF{#2}{\IfNoValueTF{#3}{
	\IfNoValueTF{#4}{\IfNoValueTF{#5}{\mathbf{r}}{\mathbf{r}\left[#5\right]}}{\IfNoValueTF{#5}{OCCHIO!1}{\mathbf{r}\left[#4,#5\right]}}}
	{\IfNoValueTF{#4}{OCCHIO!2}{\IfNoValueTF{#5}{OCCHIO!3}{\mathbf{r}\left[#3,#4,#5\right]}}}}
	{\IfNoValueTF{#3}{
	\IfNoValueTF{#4}{\IfNoValueTF{#5}{\mathbf{r}}{\mathbf{r}_{#5}}}{\IfNoValueTF{#5}{OCCHIO!4}{\mathbf{r}_{#4,#5}}}}
	{\IfNoValueTF{#4}{OCCHIO!5}{\IfNoValueTF{#5}{OCCHIO!6}{\mathbf{r}_{#3,#4,#5}}}}}}	
	{\IfNoValueTF{#2}{OCCHIO!7}{	
	\IfNoValueTF{#3}{
	\IfNoValueTF{#4}{\IfNoValueTF{#5}{\mathbf{r}^{#1}}{\mathbf{r}^{#1}_{#5}}}{\IfNoValueTF{#5}{OCCHIO!8}{\mathbf{r}^{#1}_{#4,#5}}}}
	{\IfNoValueTF{#4}{OCCHIO!9}{\IfNoValueTF{#5}{OCCHIO!10}{\mathbf{r}^{#1}_{#3,#4,#5}}}}}} 
}
\NewDocumentCommand\invquotsm{d// d-- d() d++}{% 
  \IfNoValueTF{#1}{\IfNoValueTF{#2}{\IfNoValueTF{#3}{\IfNoValueTF{#4}{\mathbf{s}}{\mathbf{s}\left[#4\right]}}{\IfNoValueTF{#4}{OCCHIO!1}{\mathbf{s}\left[#3,#4\right]}}}
	{\IfNoValueTF{#3}{\IfNoValueTF{#4}{OCCHIO!2}{\mathbf{s}_{#4}}}{\IfNoValueTF{#4}{OCCHIO!3}{\mathbf{s}_{#3,#4}}}}}	
	{\IfNoValueTF{#2}{OCCHIO!4}{	
	\IfNoValueTF{#3}{\IfNoValueTF{#4}{OCCHIO!2}{\mathbf{s}^{#1}_{#4}}}{\IfNoValueTF{#4}{OCCHIO!3}{\mathbf{s}^{#1}_{#3,#4}}}
	}}
}
\NewDocumentCommand\natisprodfib{d// d-- d[] d() d<< d>>}{  
  \IfNoValueTF{#1}{
	\IfNoValueTF{#2}{\IfNoValueTF{#3}{
	\IfNoValueTF{#4}{\IfNoValueTF{#5}{\IfNoValueTF{#6}{\tau}{OCCHIO!1}}{\IfNoValueTF{#6}{OCCHIO!2}{\tau\left[#5,#6\right]}}}{\IfNoValueTF{#5}{OCCHIO!3}{ 
		\IfNoValueTF{#6}{OCCHIO!4}{\tau\left[#4,#5,#6\right]}}}}
	{\IfNoValueTF{#4}{OCCHIO!5}{\IfNoValueTF{#5}{OCCHIO!6}{\IfNoValueTF{#6}{OCCHIO!7}{\tau\left[#3,#4,#5,#6\right]}}}}}
	{	\IfNoValueTF{#3}{
	\IfNoValueTF{#4}{\IfNoValueTF{#5}{\IfNoValueTF{#6}{\tau}{OCCHIO!8}}
	{\IfNoValueTF{#6}{OCCHIO!9}{\tau_{#5,#6}}}}{
	\IfNoValueTF{#5}{\IfNoValueTF{#6}{\tau_{#4}}{OCCHIO!10}}
	{\IfNoValueTF{#6}{OCCHIO!11}{\tau_{#4,#5,#6}}}}}
	{\IfNoValueTF{#4}{OCCHIO!12}{\IfNoValueTF{#5}{\IfNoValueTF{#6}{\tau_{#3,#4}}{OCCHIO!13}}
	{\IfNoValueTF{#6}{OCCHIO!14}{\tau_{#3,#4,#5,#6}}}}}}}	
	{\IfNoValueTF{#2}{OCCHIO!15}{\IfNoValueTF{#3}{
	\IfNoValueTF{#4}{\IfNoValueTF{#5}{\IfNoValueTF{#6}{\tau^{#1}}{OCCHIO!16}}
	{\IfNoValueTF{#6}{OCCHIO!17}{\tau_{#5,#6}^{#1}}}}{
	\IfNoValueTF{#5}{\IfNoValueTF{#6}{\tau_{#4}^{#1}}{OCCHIO!18}}
	{\IfNoValueTF{#6}{OCCHIO!19}{\tau_{#4,#5,#6}^{#1}}}}}
	{\IfNoValueTF{#4}{OCCHIO!20}{\IfNoValueTF{#5}{\IfNoValueTF{#6}{\tau_{#3,#4}^{#1}}{OCCHIO!21}}
	{\IfNoValueTF{#6}{OCCHIO!22}{\tau_{#3,#4,#5,#6}^{#1}}}}}}}
}
\NewDocumentCommand\natisprodfibinv{d// d-- d[] d() d<< d>>}{ 
\IfNoValueTF{#1}{
	\IfNoValueTF{#2}{\IfNoValueTF{#3}{
	\IfNoValueTF{#4}{\IfNoValueTF{#5}{\IfNoValueTF{#6}{\sigma}{OCCHIO!1}}{\IfNoValueTF{#6}{OCCHIO!2}{\sigma\left[#5,#6\right]}}}{\IfNoValueTF{#5}{OCCHIO!3}{ 
		\IfNoValueTF{#6}{OCCHIO!4}{\sigma\left[#4,#5,#6\right]}}}}
	{\IfNoValueTF{#4}{OCCHIO!5}{\IfNoValueTF{#5}{OCCHIO!6}{\IfNoValueTF{#6}{OCCHIO!7}{\sigma\left[#3,#4,#5,#6\right]}}}}}	
	{\IfNoValueTF{#3}{
	\IfNoValueTF{#4}{\IfNoValueTF{#5}{\IfNoValueTF{#6}{\sigma}{OCCHIO!8}}
	{\IfNoValueTF{#6}{OCCHIO!9}{\sigma_{#5,#6}}}}{
	\IfNoValueTF{#5}{\IfNoValueTF{#6}{\sigma_{#4}}{OCCHIO!10}}
	{\IfNoValueTF{#6}{OCCHIO!11}{\sigma_{#4,#5,#6}}}}}
	{\IfNoValueTF{#4}{OCCHIO!12}{\IfNoValueTF{#5}{\IfNoValueTF{#6}{\sigma_{#3,#4}}{OCCHIO!13}}
	{\IfNoValueTF{#6}{OCCHIO!14}{\sigma_{#3,#4,#5,#6}}}}}}}
	{\IfNoValueTF{#2}{OCCHIO!15}{
	\IfNoValueTF{#3}{
	\IfNoValueTF{#4}{\IfNoValueTF{#5}{\IfNoValueTF{#6}{\sigma^{#1}}{OCCHIO!16}}
	{\IfNoValueTF{#6}{OCCHIO!17}{\sigma_{#5,#6}^{#1}\right)}}}{
	\IfNoValueTF{#5}{\IfNoValueTF{#6}{\sigma_{#4}^{#1}}{OCCHIO!18}}
	{\IfNoValueTF{#6}{OCCHIO!19}{\sigma_{#4,#5,#6}^{#1}}}}}
	{\IfNoValueTF{#4}{OCCHIO!20}{\IfNoValueTF{#5}{\IfNoValueTF{#6}{\sigma_{#3,#4}^{#1}}{OCCHIO!21}}
	{\IfNoValueTF{#6}{OCCHIO!22}{\sigma_{#3,#4,#5,#6}^{#1}}}}}}}
}

\newcommand{\GTopcat}{\mathbb{V}}
\newcommand{\csGT}{\mathcal{T}} 
\newcommand{\GVecR}{\mathbb{VEC}_{\mathbb{R}}}
\newcommand{\GAlgR}{\mathbb{ALG}_{\mathbb{R}}}

\NewDocumentCommand\gdgerm{d[]}{ 
\IfNoValueTF{#1}{\widetilde{\Game}}{
\grm[#1]}
}

\newcommand{\ihtlre}[3]{\phi\left[{#1},{#2},{#3}\right]}
\NewDocumentCommand\corrarrre{d[]}{\IfNoValueTF{#1}{OCCHIO!}{\overset{\star}{#1}}}
\newcommand{\ihtl}[3]{\psi\left[{#1},{#2},{#3}\right]}
\NewDocumentCommand\dspdihtl{d// d-- d[] d() d++}{ 
\IfNoValueTF{#1}{
	\IfNoValueTF{#2}{\IfNoValueTF{#3}{
\IfNoValueTF{#4}{\IfNoValueTF{#5}{\widehat{\psi}}{\widehat{\psi}_{#5}}}{\IfNoValueTF{#5}{OCCHIO!1}{\widehat{\psi}_{#4,#5}}}}
	{\IfNoValueTF{#4}{OCCHIO!2}{\IfNoValueTF{#5}{OCCHIO!3}{\widehat{\psi}_{#3,#4,#5}}}}}
	{\IfNoValueTF{#3}{
	\IfNoValueTF{#4}{\IfNoValueTF{#5}{\widehat{\psi}\left(#2\right)}{\widehat{\psi}_{#5}\left(#2\right)}}{\IfNoValueTF{#5}{OCCHIO!4}{\widehat{\psi}_{#4   ,#5}\left(#2\right)}}}
	{\IfNoValueTF{#4}{OCCHIO!5}{\IfNoValueTF{#5}{OCCHIO!6}{\widehat{\psi}_{#3,#4,#5}\left(#2\right)}}}}}
	{ \IfNoValueTF{#2}{\IfNoValueTF{#3}{
	\IfNoValueTF{#4}{\IfNoValueTF{#5}{OCCHIO!6}{\widehat{\psi}\left[#5\right]}}{\IfNoValueTF{#5}{OCCHIO!7}{\widehat{\psi}\left[#4,#5\right]}}}
	{\IfNoValueTF{#4}{OCCHIO!8}{\IfNoValueTF{#5}{OCCHIO!9}{\widehat{\psi}\left[#3,#4,#5\right]}}}}{OCCHIO!10}} 
}
\NewDocumentCommand\corrarr{d[]}{\IfNoValueTF{#1}{OCCHIO!}{\widetilde{#1}}}

\newcommand{\symls}{\topindexquattro}
\NewDocumentCommand\ls{d>>}{
\IfNoValueTF{#1}{\symls}{\symls_{#1}}
}
\NewDocumentCommand\lss{d>>}{
\IfNoValueTF{#1}{\mathsf{\symls}}{\mathsf{\symls}_{#1}}
}
\newcommand{\symantels}{\topindextre}
\NewDocumentCommand\antels{d>>}{
\IfNoValueTF{#1}{\symantels}{\symantels_{#1}}
}
\NewDocumentCommand\antelss{d>>}{
\IfNoValueTF{#1}{\mathsf{\symantels}}{\mathsf{\symantels}_{#1}}
}
\newcommand{\symddgtv}{\topindexotto}
\newcommand{\elindexddgtv}{\mathsf{l}}
\NewDocumentCommand\ddgtv{d<< d>>}{
\IfNoValueTF{#1}{\IfNoValueTF{#2}{\symddgtv}{\symddgtv_{#2}}}
{\IfNoValueTF{#2}{\mathbf{\symddgtv}}{\mathbf{\symddgtv}_{#2}}}
} 
\newcommand{\symdgtv}{\topindexnove}
\NewDocumentCommand\dgtv{d<< d>>}{
\IfNoValueTF{#1}{\IfNoValueTF{#2}{\symdgtv}{\symdgtv_{#2}}}
{\IfNoValueTF{#2}{\mathbf{\symdgtv}}{\mathbf{\symdgtv}_{#2}}}
}
\newcommand{\symuedgtv}{\topindexdieci}
\NewDocumentCommand\uedgtv{d<< d>>}{
\IfNoValueTF{#1}{\IfNoValueTF{#2}{\symuedgtv}{\symuedgtv_{#2}}}
{\IfNoValueTF{#2}{\mathbf{\symuedgtv}}{\mathbf{\symuedgtv}_{#2}}}
}
\NewDocumentCommand\elindexuedgtv{d>>}{
\IfNoValueTF{#1}{\elindexdieci}{\elindexdieci_{#1}}
}
\NewDocumentCommand\ouol{d<< d>>}{
\IfNoValueTF{#1}{\IfNoValueTF{#2}{\delta}{\delta_{#2}}}
{\IfNoValueTF{#2}{\mbox{\boldmath$\delta$}}{\mbox{\boldmath$\delta$}_{#2}}}
}

\NewDocumentCommand\AQ{d// d-- d++ d<<}{ 
\IfNoValueTF{#1}{
\IfNoValueTF{#2}{
\IfNoValueTF{#3}{\IfNoValueTF{#4}{\mathfrak{A}}{\mathfrak{A}_{#4}}}{\IfNoValueTF{#4}{\mathfrak{A}_{\left(#3\right)}}{\mathfrak{A}_{\left(#3\right),#4}}}
}
{\IfNoValueTF{#3}{\IfNoValueTF{#4}{\mathfrak{A}^{#2}}{\mathfrak{A}^{#2}_{#4}}}{\IfNoValueTF{#4}{\mathfrak{A}^{#2}_{\left(#3\right)}}{\mathfrak{A}^{#2}_{\left(#3\right),#4}}}}
}{
\IfNoValueTF{#2}{
\IfNoValueTF{#3}{\IfNoValueTF{#4}{OCCHIO!4}{OCCHIO5!}}{\IfNoValueTF{#4}{OCCHIO6!}{\mathfrak{A}_{#1,\left(#3\right),#4}}}
}
{\IfNoValueTF{#3}{OCCHIO!8}{\IfNoValueTF{#4}{OCCHIO!9}{\mathfrak{A}^{#2}_{#1,\left(#3\right),#4}}}}
}
}
\NewDocumentCommand\AQdiuni{d// d-- d++ d<<}{ 
\IfNoValueTF{#1}{
\IfNoValueTF{#2}{
\IfNoValueTF{#3}{\IfNoValueTF{#4}{\mathfrak{R}}{\mathfrak{R}_{#4}}}{\IfNoValueTF{#4}{\mathfrak{R}_{\left(#3\right)}}{\mathfrak{R}_{\left(#3\right),#4}}}
}
{\IfNoValueTF{#3}{\IfNoValueTF{#4}{\mathfrak{R}^{#2}}{\mathfrak{R}^{#2}_{#4}}}{\IfNoValueTF{#4}{\mathfrak{R}^{#2}_{\left(#3\right)}}{\mathfrak{R}^{#2}_{\left(#3\right),#4}}}}
}{
\IfNoValueTF{#2}{
\IfNoValueTF{#3}{\IfNoValueTF{#4}{OCCHIO!4}{OCCHIO5!}}{\IfNoValueTF{#4}{OCCHIO6!}{\mathfrak{R}_{#1,\left(#3\right),#4}}}
}
{\IfNoValueTF{#3}{OCCHIO!8}{\IfNoValueTF{#4}{OCCHIO!9}{\mathfrak{R}^{#2}_{#1,\left(#3\right),#4}}}}
}
}
\newcommand{\elAQ}{\elsymuno}
\NewDocumentCommand\consec{d-- d() d++}{ 
\IfNoValueTF{#1}{\IfNoValueTF{#2}{\IfNoValueTF{#3}{\curvearrowright}{OCCHIO!1}}{OCCHIO!2}}
{\IfNoValueTF{#2}{OCCHIO!3}{\IfNoValueTF{#3}{OCCHIO!4}{
\tensor[_{#1}]{\overset{#2}{\curvearrowright}}{_{#3}}
}}}
}
\NewDocumentCommand\anteconsec{d-- d() d++}{ 
\IfNoValueTF{#1}{\IfNoValueTF{#2}{\IfNoValueTF{#3}{\leftrightarrow}{OCCHIO!1}}{OCCHIO!2}}
{\IfNoValueTF{#2}{OCCHIO!3}{\IfNoValueTF{#3}{OCCHIO!4}{
\tensor[_{#1}]{\overset{#2}{\leftrightarrow}}{_{#3}}
}}}
}

\NewDocumentCommand\predcR{d[] d-- d++ d<<}{ 
  \IfNoValueTF{#1}{	
\IfNoValueTF{#2}{
\IfNoValueTF{#3}{\IfNoValueTF{#4}{\widecheck{\mathcal{A}}}{\widecheck{\mathcal{A}}_{#4}}}{\IfNoValueTF{#4}{OCCHIO1!}{\widecheck{\mathcal{A}}_{\left(#3\right),#4}}}
}
{\IfNoValueTF{#3}{\IfNoValueTF{#4}{\widecheck{\mathcal{A}}^{#2}}{\widecheck{\mathcal{A}}^{#2}_{#4}}}{\IfNoValueTF{#4}{OCCHIO!2}{\widecheck{\mathcal{A}}^{#2}_{\left(#3\right),#4}}}}
}{
\widecheck{\mathcal{A}}\left[#1\right]}
}
\NewDocumentCommand\predcRdiuni{d[] d-- d++ d<<}{ 
  \IfNoValueTF{#1}{	
\IfNoValueTF{#2}{
\IfNoValueTF{#3}{\IfNoValueTF{#4}{\widecheck{\mathcal{R}}}{\widecheck{\mathcal{R}}_{#4}}}{\IfNoValueTF{#4}{OCCHIO1!}{\widecheck{\mathcal{R}}_{\left(#3\right),#4}}}
}
{\IfNoValueTF{#3}{\IfNoValueTF{#4}{\widecheck{\mathcal{R}}^{#2}}{\widecheck{\mathcal{R}}^{#2}_{#4}}}{\IfNoValueTF{#4}{OCCHIO!2}{\widecheck{\mathcal{R}}^{#2}_{\left(#3\right),#4}}}}
}{
\widecheck{\mathcal{R}}\left[#1\right]}
}
\newcommand{\predcRel}{\widecheck{\Game}} 
\NewDocumentCommand\nullpredcR{d[] d-- d++ d<< d()}{ 
\IfNoValueTF{#5}{
\IfNoValueTF{#1}{
\IfNoValueTF{#2}{
\IfNoValueTF{#3}{
\IfNoValueTF{#4}{\widecheck{\mathcal{O}}}{\widecheck{\mathcal{O}}_{#4}}
}{\IfNoValueTF{#4}{\widecheck{\mathcal{O}}_{\left(#3\right)}}{
\widecheck{\mathcal{O}}_{\left(#3\right),#4}}}}
{\IfNoValueTF{#3}{\IfNoValueTF{#4}{
\widecheck{\mathcal{O}}_{#2}}{
\widecheck{\mathcal{O}}^{#2}_{#4}
}}{\IfNoValueTF{#4}{OCCHIO!2}{
\widecheck{\mathcal{O}}^{#2}_{\left(#3\right),#4}
}}}
}
{\widecheck{\mathcal{O}}\left[#1\right]}
}{
\IfNoValueTF{#1}{
\IfNoValueTF{#2}{
\IfNoValueTF{#3}{\IfNoValueTF{#4}{\widecheck{\mathcal{O}}_{#5}
}{\widecheck{\mathcal{O}}_{#4,#5}}}{\IfNoValueTF{#4}{
\widecheck{\mathcal{O}}_{\left(#3\right),#5}
}{\widecheck{\mathcal{O}}_{\left(#3\right),#4,#5}}}
}
{\IfNoValueTF{#3}{\IfNoValueTF{#4}{\widecheck{\mathcal{O}}^{#2}_{#5}}{\widecheck{\mathcal{R}}^{#2}_{#4,#5}}}{\IfNoValueTF{#4}{
\widecheck{\mathcal{O}}^{#2}_{\left(#3\right),#5}}{
\widecheck{\mathcal{O}}^{#2}_{\left(#3\right),#4,#5}
}}}
}{\widecheck{\mathcal{O}}_{#5}\left[#1\right]}
}}
\NewDocumentCommand\gennullwpredcR{d++}{ 
\IfNoValueTF{#1}{\setsymuno}{
\setsymuno_{\left(#1\right)}
}}
\NewDocumentCommand\predcRsum{d<< d>>}{
  \IfNoValueTF{#1}{\IfNoValueTF{#2}{\widecheck{+}}{OCCHIO!1}}
{\IfNoValueTF{#2}{OCCHIO!2}{\widecheck{+}\!\left[#1,#2\right]}}
}
\NewDocumentCommand\bigpredcRsum{d<< d>>}{
  \IfNoValueTF{#1}{\IfNoValueTF{#2}{\scalerel*{\widecheck{+}}{\sum}}{OCCHIO!1}}
{\IfNoValueTF{#2}{OCCHIO!2}{\scalerel*{\widecheck{+}}{\sum}\!\left[#1,#2\right]}}
}
\NewDocumentCommand\predcRscalp{d<< d>>}{
  \IfNoValueTF{#1}{\IfNoValueTF{#2}{\,\widecheck{\ast}\,}{OCCHIO!1}}
{\IfNoValueTF{#2}{OCCHIO!2}{\,\widecheck{\ast}\!\left[#1,#2\right]\,}}
}
\NewDocumentCommand\predcRprod{d<< d() d>>}{
 \IfNoValueTF{#2}{\IfNoValueTF{#1}{\IfNoValueTF{#3}{\widecheck{\boxtimes}}{OCCHIO!1}}
{\IfNoValueTF{#3}{OCCHIO!2}{\tensor[_{#1}]{{\widecheck{\boxtimes}}}{_{#3}}}}}
{\tensor[_{#1}]{{\overset{#2}{\widecheck{\boxtimes}}}}{_{#3}}}
}
\NewDocumentCommand\antedcRprod{d<<  d>> d||}{ 
 \IfNoValueTF{#3}{\IfNoValueTF{#1}{\IfNoValueTF{#2}{\widecheck{\boxdot}}{OCCHIO!1}}
{\IfNoValueTF{#2}{OCCHIO!2}{\tensor[_{#1}]{{\widecheck{\boxdot}}}{_{#2}}}}}{
\widecheck{\boxdot}_{#3}
}
}

\NewDocumentCommand\evpdf{d>> d<< d[] d()}{ % valutazione dei germi delle derivazioni composte
\IfNoValueTF{#1}{
\IfNoValueTF{#2}{\IfNoValueTF{#3}{\IfNoValueTF{#4}{EPF}{OCCHIO1!}}{\IfNoValueTF{#4}{EPF_{#3}}{EPF_{#3}\left[#4\right]}}}
{\IfNoValueTF{#3}{\IfNoValueTF{#4}{EPF^{#2}}{EPF^{#2}\left[#4\right]}}{\IfNoValueTF{#4}{EPF^{#2}_{#3}}{EPF^{#2}_{#3}\left[#4\right]}}}
}{
\IfNoValueTF{#2}{\IfNoValueTF{#3}{\IfNoValueTF{#4}{\overset{\bullet}{EPF}}{OCCHIO1!}}{\IfNoValueTF{#4}{\overset{\bullet}{EPF}_{#3}}{\overset{\bullet}{EPF}_{#3}\left[#4\right]}}}
{\IfNoValueTF{#3}{\IfNoValueTF{#4}{\overset{\bullet}{EPF}^{#2}}{\overset{\bullet}{EPF}^{#2}\left[#4\right]}}{\IfNoValueTF{#4}{\overset{\bullet}{EPF}^{#2}_{#3}}{\overset{\bullet}{EPF}^{#2}_{#3}\left[#4\right]}}}}
}
\NewDocumentCommand\PRELOCAT{ d-- d++}{
\IfNoValueTF{#1}{\IfNoValueTF{#2}{\widecheck{\lfloor}}{OCCHIO!2}}{\IfNoValueTF{#2}{OCCHIO!2}{#1 \widecheck{\lfloor}_{#2}}}
}
\NewDocumentCommand\PREDER{ d-- d++}{
\IfNoValueTF{#1}{\IfNoValueTF{#2}{\widecheck{\mathcal{D}}}{OCCHIO!1}}{\IfNoValueTF{#2}{OCCHIO!2}{\widecheck{\mathcal{D}}_{#2}\left(#1\right)}}
}
\NewDocumentCommand\predcRquot{d[] d-- d++ d<< d()}{
\IfNoValueTF{#5}{
\IfNoValueTF{#1}{
\IfNoValueTF{#2}{
\IfNoValueTF{#3}{\IfNoValueTF{#4}{\widecheck{q}}{\widecheck{q}_{#4}}}{\IfNoValueTF{#4}{OCCHIO1!}{\widecheck{q}_{#3,#4}}}
}
{\IfNoValueTF{#3}{OCCHIO!2}{\IfNoValueTF{#4}{OCCHIO!3}{\widecheck{q}_{#2,#3,#4}}}}
}{
\widecheck{q}\left[#1\right]}
}{
\IfNoValueTF{#1}{
\IfNoValueTF{#2}{
\IfNoValueTF{#3}{\IfNoValueTF{#4}{\widecheck{q}_{#5}}{\widecheck{q}_{#4,#5}}}{\IfNoValueTF{#4}{OCCHIO1!}{\widecheck{q}_{#3,#4,#5}}}
}
{\IfNoValueTF{#3}{OCCHIO!2}{\IfNoValueTF{#4}{OCCHIO!3}{\widecheck{q}_{#2,#3,#4,#5}}}}
}{
\widecheck{q}\left[#1\right]}
}
}

\NewDocumentCommand\dcR{d-- d++}{ 
\IfNoValueTF{#1}{
\IfNoValueTF{#2}{\mathcal{R}}{\mathcal{R}_{\left(#2\right)}}
}
{\IfNoValueTF{#2}{\mathcal{R}^{#1}}{\mathcal{R}^{#1}_{\left(#2\right)}}}
}
\NewDocumentCommand\predcRquotfun{d-- d++}{ 
\IfNoValueTF{#1}{
\IfNoValueTF{#2}{q}{q_{\left(#2\right)}}
}
{\IfNoValueTF{#2}{q^{#1}}{q^{#1}_{\left(#2\right)}}}
}

\NewDocumentCommand\dcRsum{d-- d++}{
  \IfNoValueTF{#1}{\IfNoValueTF{#2}{+}{OCCHIO!1}}
{\IfNoValueTF{#2}{OCCHIO!2}{+\!\left[#1,#2\right]}}
}
\NewDocumentCommand\dcRscalp{d-- d++}{
  \IfNoValueTF{#1}{\IfNoValueTF{#2}{\ast}{OCCHIO!1}}
{\IfNoValueTF{#2}{OCCHIO!2}{\ast\!\left[#1,#2\right]}}
}
\NewDocumentCommand\dcRprod{d<< d() d>>}{
 \IfNoValueTF{#2}{\IfNoValueTF{#1}{\IfNoValueTF{#3}{\boxdot}{OCCHIO!1}}
{\IfNoValueTF{#3}{OCCHIO!2}{\boxdot\!\left[#1,#3\right]}}}
{\boxdot\!\left[#1,#2,#3\right]}
}

\NewDocumentCommand\dcRl{d-- d++ d[]}{ 
  \IfNoValueTF{#1}{
\IfNoValueTF{#2}{\IfNoValueTF{#3}{\overset{\leftarrow}{\mathcal{R}}}{\overset{\leftarrow}{\mathcal{R}}\!\left[#3\right]}}{OCCHIO!1}}
{\IfNoValueTF{#2}{OCCHIO!2}{\overset{\leftarrow}{\mathcal{R}}_{#1}\left(#2\right)}}
}
\NewDocumentCommand\intunosk{d>>}{ 
\IfNoValueTF{#1}{\intuno}{\intuno_{#1}}
}
\NewDocumentCommand\unkunosk{d>>}{
\IfNoValueTF{#1}{\unkuno}{\unkuno_{#1}}
}
\NewDocumentCommand\domsk{d>>}{ 
\IfNoValueTF{#1}{\setsymnove}{\setsymnove_{#1}}
}
\NewDocumentCommand\tddmsk{d>>}{ % dimensione dominio degli scheletri, \tddmsk=\ddgtv+\dgtv
\IfNoValueTF{#1}{d}{d_{#1}}
}
\NewDocumentCommand\opcov{d() d[] d-- d++ d<< d>> d||}{ 
\IfNoValueTF{#7}{
\IfNoValueTF{#1}{\IfNoValueTF{#2}{\IfNoValueTF{#3}{\IfNoValueTF{#4}{\IfNoValueTF{#5}{\IfNoValueTF{#6}{\mathcal{O}}{OCCHIO1!}}{OCCHIO2!}}{OCCHIO3!}}{OCCHIO4!}}{\IfNoValueTF{#3}{\IfNoValueTF{#4}{\IfNoValueTF{#5}{\IfNoValueTF{#6}{\mathcal{O}_{#2}}{OCCHIO!5}}{OCCHIO!6}}{OCCHIO!7}}{OCCHIO!8}}}
{\IfNoValueTF{#2}{
\IfNoValueTF{#3}{\IfNoValueTF{#4}{\IfNoValueTF{#5}{\IfNoValueTF{#6}{\widehat{\mathcal{O}}}{OCCHIO!9}}{OCCHIO!10}
}{OCCHIO!11}}{\IfNoValueTF{#4}{OCCHIO!12}{
\IfNoValueTF{#5}{\IfNoValueTF{#6}{\widehat{\mathcal{O}}_{#3,#4}}{OCCHIO!13}}{\IfNoValueTF{#6}{OCCHIO!14}{\widehat{\mathcal{O}}_{#3,#4}^{#5,#6}}}}}}{\widehat{\mathcal{O}}_{#2}}}}{
\IfNoValueTF{#1}{\IfNoValueTF{#2}{\IfNoValueTF{#3}{\IfNoValueTF{#4}{\IfNoValueTF{#5}{\IfNoValueTF{#6}{{_{#7}\mathcal{O}}}{OCCHIO1!}}{OCCHIO2!}}{OCCHIO3!}}{OCCHIO4!}}{\IfNoValueTF{#3}{\IfNoValueTF{#4}{\IfNoValueTF{#5}{\IfNoValueTF{#6}{{_{#7}\mathcal{O}}_{#2}}{OCCHIO!5}}{OCCHIO!6}}{OCCHIO!7}}{OCCHIO!8}}}
{\IfNoValueTF{#2}{
\IfNoValueTF{#3}{\IfNoValueTF{#4}{\IfNoValueTF{#5}{\IfNoValueTF{#6}{{_{#7}\widehat{\mathcal{O}}}}{OCCHIO!9}}{OCCHIO!10}
}{OCCHIO!11}}{\IfNoValueTF{#4}{OCCHIO!12}{
\IfNoValueTF{#5}{\IfNoValueTF{#6}{{_{#7}\widehat{\mathcal{O}}}_{#3,#4}}{OCCHIO!13}}{\IfNoValueTF{#6}{OCCHIO!14}{{_{#7}\widehat{\mathcal{O}}}_{#3,#4}^{#5,#6}}}}}}{{_{#7}\widehat{\mathcal{O}}}_{#2}}}}
}
\NewDocumentCommand\LOCAT{ d-- d++}{
\IfNoValueTF{#1}{\IfNoValueTF{#2}{\lfloor}{OCCHIO!2}}{\IfNoValueTF{#2}{OCCHIO!2}{#1 \lfloor_{#2}}}
}
\NewDocumentCommand\DER{d-- d++}{
\IfNoValueTF{#1}{\IfNoValueTF{#2}{\mathcal{D}}{OCCHIO!1}}{\IfNoValueTF{#2}{OCCHIO!2}{\mathcal{D}_{#2}\left(#1\right)}}
}

\newcommand{\atlas}{\mathcal{A}}
\newcommand{\cl}{\varphi}

\NewDocumentCommand\tenalgtv{d<< d>>}{% tensor algebra functor
\IfNoValueTF{#1}{\IfNoValueTF{#2}{\Theta^{0}}{OCCHIO1!}}{\IfNoValueTF{#2}{OCCHIO2!}{\Theta^{0}_{#1}\left[#2\right]}}
}
\NewDocumentCommand\symalgtv{d<< d>>}{% symmetric algebra functor
\IfNoValueTF{#1}{\IfNoValueTF{#2}{\Delta^{0}}{OCCHIO1!}}{\IfNoValueTF{#2}{OCCHIO2!}{\Delta^{0}_{#1}\left[#2\right]}}
}
\NewDocumentCommand\symjac{d<< d>> d[]}{% symmetric algebra functor
\IfNoValueTF{#1}{\IfNoValueTF{#2}{\IfNoValueTF{#3}{d^{0}}{OCCHIO!1}}{OCCHIO2!}}{\IfNoValueTF{#2}{OCCHIO3!}{\IfNoValueTF{#3}{d^{0}_{#1}\left[#2\right]}}{d^{0}_{#1}\left[#2,#3\right]}}
}
\NewDocumentCommand\extalgtv{d<< d>>}{% symmetric algebra functor
\IfNoValueTF{#1}{\IfNoValueTF{#2}{\Lambda^{0}}{OCCHIO1!}}{\IfNoValueTF{#2}{OCCHIO2!}{\Lambda^{0}_{#1}\left[#2\right]}}
}

\newcommand{\contchart}{LC} % category of continuous local charts
\newcommand{\Gtcch}{J} % Grothendieck topology on \contchart 
\newcommand{\Tmcat}{M} % category of topological manifolds
\newcommand{\smoothchart}{LC^{\infty}} % category of smooth local charts
\newcommand{\Gtsch}{J^{\infty}} % Grothendieck topology on \smoothchart 
\newcommand{\Smcat}{M^{\infty}} % category of smooth manifolds
\newcommand{\ntbf}{\mathbf{0}} % null trivial bundle
\newcommand{\genfib}{\mathbf{F}}
\newcommand{\genfibsm}{\mathbf{SF}}
\newcommand{\kergenfibsm}{\mathbf{KSF}}
\newcommand{\cokergenfib}{\mathbf{CKF}}
\newcommand{\prntgf}{\mathsf{pr}}
\newcommand{\ckprntgf}{\mathsf{cpr}}
\newcommand{\prntgfsm}{\mathsf{prs}}
\newcommand{\prkergenfibsm}{\mathsf{kprs}}
\NewDocumentCommand\gentotfib{d[]}{ 
\IfNoValueTF{#1}{\mathbf{F}_{\segnvar}}{
\mathbf{F}_{#1}}
}
\newcommand{\genbunfunc}{\mathbf{B}}
\newcommand{\smgenbunfunc}{\mathbf{SB}}
\newcommand{\kergenbunfuncsm}{\mathbf{KSB}}
\newcommand{\cokergenbunfunc}{\mathbf{CKB}}
\newcommand{\smoothfib}{\mathbf{F}^{\infty}}
\newcommand{\prntsm}{\mathsf{ps}}
\newcommand{\ptbf}{\mu} %projection trivial fibre bundle
\newcommand{\kerptbf}{\varpi} %projection kernel trivial fibre bundle
\newcommand{\ptbb}{\widetilde{\ptbf}} %projection bundle
\newcommand{\kerptbb}{\widetilde{\kerptbf}} %projection kernel bundle
\newcommand{\itbf}{\nu} %inclusion trivial fibre bundle
\newcommand{\ckitbf}{\varsigma} %inclusion cokernel on trivial fibre bundle
\newcommand{\itbb}{\widetilde{\itbf}} %inclusion  bundle
\newcommand{\ckitbb}{\widetilde{\ckitbf}}%inclusion cokernel bundle
\newcommand{\smoothbunfunc}{\mathbf{B}^{\infty}}
\newcommand{\natisprod}{\widetilde{\tau}}
\newcommand{\natisprodinv}{\widetilde{\sigma}}
\newcommand{\natisprodsm}{\natisprod^{\infty}}
\newcommand{\natisprodinvsm}{\natisprodinv^{\infty}}

\NewDocumentCommand\Ckspquotin{d()}{
\IfNoValueTF{#1}{OCCHIO!1}{\mathbf{j}_{#1}}
}

\NewDocumentCommand\invCkspquot{m d// d-- d() d++}{ % elementi invertibili di \Ckspquot 
 \IfNoValueTF{#2}{
 \IfNoValueTF{#3}{
	\IfNoValueTF{#4}{{\mathbf{A}}^{#1}}{{\mathbf{A}}^{#1}\left[#4\right]}
	}{
	\IfNoValueTF{#4}{OCCHIO!1}{{\mathbf{A}}^{#1}_{#4}}
  }}
	{\IfNoValueTF{#3}{
	\IfNoValueTF{#4}{OCCHIO!2}{{\mathbf{A}}^{#1}_{#4}}
	}
	{
	\IfNoValueTF{#4}{OCCHIO!3}{{\mathbf{A}}^{#1\,#3}_{#4}}
	}}
}

\newcommand{\opdisc}{D}
\newcommand{\circfiss}{\gamma}
\newcommand{\sph}{S}
\newcommand{\upcal}{T_{+}}
\newcommand{\uncal}{T_{-}}
\newcommand{\cor}{C}
\newcommand{\mappn}{f_+}
\newcommand{\mapps}{f_-}
\newcommand{\corcv}{g}
\newcommand{\cv}{\mbox{\boldmath$\vecuno$}}
\newcommand{\cvup}{\mbox{\boldmath$\vecuno$}_+}
\newcommand{\cvun}{\mbox{\boldmath$\vecuno$}_-}

\newcommand{\preisocordue}{\psi}
\newcommand{\isocordue}{\varphi}

\allowdisplaybreaks

\begin{document}

\title{Linear microbundles}
\author{Tommaso Boccellari}
\maketitle

{\scriptsize  \tableofcontents}

\chapter{Introduction \label{Intro}}

The world of manifolds is divided into two parts: the world of smooth manifolds and the world of topological manifolds. The difference origins in the structure of the fibre of the tangent bundle associated to any manifold and in the behavior on the fibre of the differential of a function between manifolds.\newline 
More precisely: in case of smooth manifolds the fibre of the tangent bundle is the topological real vector space of all derivation operators of degree $1$ and the differential of a smooth function is continuous linear when restricted to the fibre; in case of topological manifolds the fibre of the tangent bundle has no linear structure, neither elements belonging to the fibre somehow represent suitably generalized derivation operators of any degree, and therefore the differential of a continuous function has no linear properties, neither the meaning of differential operator when restricted to the fibre.  
We refer to \cite{RMS}\;\;Chap\,12 for a detailed construction of the tangent bundle to smooth manifolds and to \cite{JWM}, \cite{RMS}\;\;Chap\,14 for a detailed construction of the tangent bundle to topological manifolds.\newline
We denote by $\smoothbunfunc_1$ the functor associating to any smooth manifold the corresponding tangent bundle and to any smooth function between manifolds the corresponding differential.\newline
Construction of the tangent bundle functor $\smoothbunfunc_1$ is based on calculus in euclidean vector spaces, and tangent vectors correspond to derivation operators of degree $1$. Derivation operators of any degree on smooth manifolds are obtained starting from $\smoothbunfunc_1$ by applying the functorial construction of symmetric algebra associated to a vector space. Then we obtain the functor $\smoothbunfunc$ associating to any smooth manifold the bundle whose fibre is the $\mathbb{R}$-algebra of all derivation operators, and to any smooth function between manifolds the corresponding differential. Algebraic structure of functor  $\smoothbunfunc$ translates in algebraic terms all relations which hold true in calculus among derivation operators.\newline 
The main result of this work is the extension of the functor $\smoothbunfunc$ to the category of topological manifolds and continuous functions, by preserving its algebraic structure and its link with calculus.
Namely we build a functor $\genbunfunc$ associating to any topological manifold a locally trivial tangent bundle whose fibre is the topological $\mathbb{R}$-algebra of all generalized derivation operators, and to any continuous function between topological manifolds a differential which is a function of topological $\mathbb{R}$-algebras when restricted to the fibre. We also prove that if $\genbunfunc$ is restricted to the category of smooth manifolds and smooth functions then there is a sub-functor $\smgenbunfunc$ of $\genbunfunc$ and a natural epimorphism $\smgenbunfunc\rightarrow \smoothbunfunc$ with a not natural right inverse.\newline
Construction of $\genbunfunc$ follows the construction of $\smoothbunfunc$ where the role of classical calculus on smooth function is played by a new type of calculus on continuous function which we develop specifically for this aim and which we refer to as the generalized calculus.\newline
Generalization of calculus is a universal construction extending the category $\Topcat$ of topological spaces and continuous functions to a category whose arrows are all smooth in a suitable sense which coincides with the classical smoothness on smooth functions of $\Topcat$, namely the smallest category with such property.\newline 
The contribution of this work to the research is outlined below:\newline
1)	The linear structure of the fibre and the linearity of differentials on fibres is new with respect to other extensions of the notion of tangent bundle to topological manifolds.\newline
2)	Linearity allows to extend to topological manifolds those Algebraic Topology techniques which are now available only for smooth manifolds.\newline
For example from the extension of de Raham cohomology and of other deep results, such as \cite{RT}, we think that new facts are to be expected.\newline
3)	Linearity allows to extend tensor calculus to topological manifolds.\newline
4)	The structure of the tangent bundle to smooth manifolds turns out to be richer than it was supposed to be up to now. Actually the fibre $\mathbb{R}^n$ of the smooth tangent bundle to a $n$-dimensional smooth manifold $M$ is a quotient of the fibre of the tangent bundle of $M$ considered as a topological manifold. Then, even on smooth manifolds, there are tensor fields which are still unknown since they cannot be represented neither in the classical setting (\cite{S}) nor in generalized setting available up to now (\cite{KS}, \cite{KS2}). We think that such tensor fields could be of great interest in physics.\newline
This is the latest in a long series of works devoted to extend to geometry and modern analysis the same deep link that differential geometry provides between geometry and classical analysis.\newline  
Since Schwartz introduced Distributions in \cite{LS0} many authors, taking in account \cite{LS}, have worked on the subject in order to solve non linear problems which continually arise in physics and analysis. Some authors enriched the structure of the space of distributions (\cite{JFC}, \cite{EER}) extending to it as many properties of traditional functions as possible. In particular Colombeau Theory of Generalized Functions has proved to be formulated intrinsically on smooth manifolds (\cite{AB}, \cite{CM}, \cite{KSV}), then it provides effective tools to differential geometry (\cite{KS}, \cite{KS2}) and physics (\cite{CVW}, \cite{OGG}). Other authors focused their attention on subspaces contained in the space of distributions, then they showed that such subspaces possess structures richer than that of the whole space of distributions (\cite{ADM}, \cite{ADLM}, \cite{AFP}, \cite{ST}, \cite{AIV}). This approach has led to an extension of differential geometry to elements belonging to a large class of topological spaces properly containing smooth manifolds but not all topological manifolds (\cite{A}, \cite{C}, \cite{F}, \cite{K}).\newline
The outline of the work is the following.\newline
Part \ref{PI} contains the construction of generalized calculus which takes Chapter \ref{DiffSp}-\ref{bascalc}.\newline
In Chapter \ref{DiffSp} we express both classical integro-differential calculus on smooth functions and integral calculus on continuous functions through a formalism which makes easier to understand the extension to continuous functions of the integro-differential calculus. We introduce four notions: integro-differential monoid, integro-differential space, integral monoid, integral space.\newline
The integro-differential monoid $\bsfM$ is the monoid whose elements are all integro-differential operators obtained by the action of partial derivatives, integral operators along coordinate directions and projections of vector valued functions on their coordinates. Monoid $\bsfM$ (Definition \ref{intdiffmon}) is defined through generators and relations: generators are partial derivatives, integral operators along coordinate directions, projections of vector valued functions on their coordinates; relations among generators are relations which hold true among them in the setting of integro-differential calculus on smooth functions.\newline
An integro-differential space (Definition \ref{DiffSpDef}) is a topological space endowed with four continuous internal everywhere defined operations (composition, cartesian product, sum, product) and two continuous everywhere defined external operations (scalar product, $\bsfM$-action) fulfilling conditions which hold true for the corresponding operations on smooth functions.\newline
The set of all smooth functions together with classical composition, cartesian product, sum, product, scalar product, $\bsfM$-action is an integro-differential space denoted by $\idssmo$.\newline
The integral monoid $\subbsfM$ is the sub-monoid of $\bsfM$ generated by integral operators and projections.\newline
An integral space (Definition \ref{QDiffSpDef}) is a topological space endowed with four continuous internal everywhere defined operations (composition, cartesian product, sum, product) and two continuous everywhere defined external operations (scalar product, $\subbsfM$-action) fulfilling conditions which hold true for the corresponding operations on continuous functions.\newline
The set of all continuous functions together with classical composition, cartesian product, sum, product, scalar product, $\bsfM$-action is an integral space denoted by $\qidscont$.\newline
In Chapter \ref{costruct} we define the integro-differential space $\genfidsp$ containing both smooth functions and continuous functions, where composition, cartesian product, sum, product, scalar product are defined and where integro-differential calculus extends from smooth functions to continuous functions (Definition \ref{defgenf}). We call generalized functions all elements belonging to $\genfidsp$. Generalized functions behaves very differently from real smooth functions but any operation in $\genfidsp$ coincide with the corresponding classical operation when performed on smooth functions.\newline
Construction of $\genfidsp$ splits into four steps. \newline
Step 1. For any fixed positive integer $m$ and any open set $\intuno\subseteq  \mathbb{R}^m$ with not empty interior we choose a set $\bascont\left(\intuno\right)$ which can be completed to an algebraic base for the real vector space $C^0\left(\intuno,\mathbb{R}\right)$ by adding elements belonging to $\Cksp{\infty}(\intuno)+\mathbb{R}+$. Elements belonging to sets $\bascont\left(\intuno\right)$ are called generating functions. Construction of $\genfidsp$ turns out to be independent on the choice of generating functions (Remark \ref{remch}-[2]). \newline 
Step 2. We define the free magma $\left(\fremag, \compone   \right)$ (Definition \ref{Cinfty}) as the direct limit of the family of free magmas $\left\{\left(\fremag_i, \compone_i   \right)\right\}_{i\in \mathbb{N}_0}$ defined inductively. Magma $\left(\fremag_0, \compone_0   \right)$ (Definition \ref{basind}) is the base of induction and is defined through elements belonging to sets $\Cksp{\infty}(\intuno)+\mathbb{R}+\cup \bascont\left(\intuno\right)$. The inductive step to obtain magma $\left(\fremag_{i+1}, \compone_{i+1}\right)$ involves both the magma $\left(\fremag_i, \compone_i\right)$ and the integro-differential monoid $\bsfM$.\newline
Related to free magma $\left(\fremag, \compone   \right)$ we define (Propositions \ref{linincl}, \ref{lininclbis},
Definitions \ref{Cinftybis}, \ref{pathinC}): cartesian products $\lboundone\,\rboundone:\fremag^n \rightarrow \fremag$ for any positive
integer $n$; $\bsfM$-action $\acmone: \bsfM \times \fremag \rightarrow \fremag$; domain and co-domain for any element belonging to $\fremag$; the sub-magma $\fremagcont$ containing any element $\fmx$ belonging to $\fremag$ such that no
differential operator acts on a continuous function; the sub-magma $\fremagsmooth$ containing all elements belonging to
$\fremag$ such that no elements belonging to $\bascont\left(\intuno\right)$ occurs in $\fmx$ for any open set $\intuno$;
the set function $\evalcompone$ associating to any $\fmx\in \fremagcont$ the continuous function obtained by evaluating
$\fmx$; the set function $\lininclone$ associating to any continuous function the corresponding unique linear combination of generating functions and smooth functions; a class of set functions $\left(-1,1\right)\rightarrow \fremag$ called paths in $\fremag$; the notion of augmentation of elements belonging to $\fremag$ (Definition \ref{augdef}) formalizing in this abstract setting the well known fact that composition of continuous functions is associative wherever defined.\newline
Through paths in $\fremag$, the function $\evalcompone$ and the notion of augmentation we introduce an equivalence relation $\relsymb$ on $\fremag$ (Definitions \ref{admrelc1}, \ref{admrelc}) which identifies elements suitably linked through special paths called detecting paths. Equivalence relation $\relsymb$ is based on two principles: retention of information (condition \eqref{eqelc}); invariance of interpretation (condition \eqref{eqelc1}). Invariance of interpretation type principle is common to all theories of generalized functions we know (\cite{HAB}, \cite{EER}), on the contrary we never met in literature any use of retention of information principle. Differences between generalized functions and smooth functions origins mainly from the retention of information principle.\newline
We prove that operations $\compone$, $\lboundone\;\rboundone$, $\acmone$ are compatible with $\relsymb$ (Proposition \ref{prpadrel}) and that set function $\evalcompone$ assumes the same value when evaluated on two equivalent elements belonging to $\fremagcont$ (Proposition \ref{augprop}-[6]).\newline
Step 3. We define the magma $\left(\magtwo, \comptwo\right)$ obtained by identifying through $\relsymb$ elements belonging to $\fremag$.\newline
Then (Proposition \ref{alprtwo}, Definition \ref{pathtwodef}): operations $\lboundone\;\rboundone$, $\acmone$ defined on $\fremag$ induce corresponding operations $\lboundtwo\;\rboundtwo$, $\acmtwo$ on $\magtwo$; domain and co-domain for any element belonging to $\fremag$ induce the notion of domain and co-domain for any element belonging to $\magtwo$; to subsets $\fremagcont$, $\fremagsmooth$ of $\fremag$ correspond subsets $\magtwocont$, $\magtwosmooth$ of $\magtwo$; set functions $\evalcompone$, $\lininclone$ induce corresponding set functions $\evalcomptwo$, $\linincltwo$; paths in $\fremag$ induce paths in $\magtwo$.\newline
Through paths in $\magtwo$ we introduce topology in $\magtwo$ (Definition \ref{topmagtwo}-[1]). We prove that $\comptwo$, $\lboundtwo\;\rboundtwo$, $\acmtwo$ are continuous operations (Proposition \ref{topinDprop}-[1-5]) and that $\evalcomptwo$, $\linincltwo$ are continuous functions (Proposition \ref{topinDcontprop}-[7, 8]).\newline
Through paths in $\magtwo$ and the function $\evalcomptwo$ we introduce the set $\magtwoempty$ containing all elements belonging to $\magtwo$ which cannot be evaluated on their whole domain (Definition \ref{topmagtwo}-[2]).\newline
Step 4. We define the integro-differential space $\genfidsp$ obtained by taking magma $\left(\magtwo,\comptwo\right)$, operations $\lboundtwo\;\rboundtwo$, $\acmtwo$ and introducing sum, product, scalar product in $\magtwo$ by exploiting the fact that smooth functions belong to $\magtwo$ by construction (Definitions \ref{defgenf}).\newline
In Chapter \ref{shffgf}, for any fixed non negative integer $m$, we focus our attention to generalized functions whose domains are open neighborhoods of $0\in \mathbb{R}^m$, we call them local generalized functions (Definition \ref{infnearpointdef}). We introduce the equivalence relation $\infneareqrel$ which identifies local generalized functions whenever they coincide locally at $0$ (Definition \ref{releqinfnearpointdef}). We examine compatibility of $\infneareqrel$ with the algebraic structure of $\magtwo$ (Propositions \ref{infnearpointpropbis}, \ref{infnearpointprop}) and we find that it suffices to repeat here, 
in the setting of local generalized function, the construction of germs of local functions which is well known in the setting of smooth and of continuous functions (Remark \ref{genercontfun}). For any fixed $n\in \mathbb{N}$ we call generalized $n$-germ the germ of a local generalized function with fixed co-domain $\mathbb{R}^n$ and $n$-germ the germ of a local function with fixed co-domain $\mathbb{R}^n$, we drop the index $n$ when there is no need of such detail. The space of  generalized germs contains the space of germs. The construction of germs, however, cannot be repeated word by word since local generalized functions behaves not exactly as local functions, in particular we have that: the space of generalized $1$-germs is a commutative weak $\mathbb{R}$-algebra (Notation \ref{alg}-[3]) which is not regular due to the retention of information principle, hence cannot be naturally embedded into an algebra; the space of generalized $n$-germs is a weak $\mathbb{R}$-module (Notation \ref{alg}-[5]).\newline
In Chapter \ref{Derggf}, for any fixed non negative integer $m$, we focus our attention to derivation operators of degree $1$ acting on generalized $1$-germs represented by local generalized functions whose domain is an open neighbprhood of $0 \in \mathbb{R}^m$. We are not able to define directly generalized derivations on generalized $1$-germs and differentials of generalized germs, we are forced to make an intermediate step: define first pre-derivations and pre-differentials. Pre-derivations on generalized $1$-germs are suitable linear combinations of derivation operators belonging to $\bsfM$ (Definition \ref{genderdef}). To any fixed germ of continuous function we associate a linear function between spaces of pre-derivations which generalizes differential of smooth functions (Proposition \ref{algstrprederpoint}-[3]). We organize this material in the functorial language to give a more compact and clear construction of topological real vector spaces of pre-derivations and of pre-differential of germs of continuous functions (Definition \ref{genderdef2}).\newline
Unfortunately construction of pre-derivations and pre-differentials does not provide a functor since the pre-differential of the composition of continuous functions does not coincides with the composition of their pre-differentials due to the fact that composition of generalized functions does not coincides with composition of continuous functions.\newline
However this problem can be solved: we define the real vector space of generalized derivation operators of degree $1$ and the notion of generalized differential of germs of continuous functions by taking the colimit of a suitable functor constructed by employing spaces of pre-derivations and pre-differentials (Definition \ref{coeqdef}). We prove that: such construction does not give the trivial vector space (Proposition \ref{algstrderpoint}); by taking suitable generalized derivations we get a $\mathbb{R}$-vector subspace which is stable by generalized differential of smooth functions and which naturally maps onto $\mathbb{R}^m$ (Proposition \ref{nattosm}).\newline
Sets of germs of generalized functions, sets of pre-derivations operators of order $1$ and sets of generalized derivations operators of order $1$ can not be endowed with classical topology (Remark \ref{siterem2}), topology which is suitable for our purposes is a Grothendieck topology. A site $\left(\GTopcat, \csGT\right)$ is defined, as an extension of the site $\Topcat$ through a suitable universal construction (Definition \ref{augsitedef}, Proposition \ref{augsiteprop}), among whose objects are sets of our interest listed above and among whose arrows are set functions related to them (Definition \ref{defcatD}).\newline
In Chapter \ref{bascalc} we study calculus in $\GTopcat$. We explicitly construct the smallest $\mathbb{R}$-algebra $\dcR$, canonically containing $\mathbb{R}$ as a subfield, such that any arrow of $\GTopcat$ with co-domain $\dcR$ is differentiable infinitely many times in the generalized sense, where differentiation is performed locally with respect to the Grothendieck topology of $\GTopcat$. Then we say that $\dcR$ is the integro-differential closure of $\mathbb{R}$ and that $\GTopcat$ is the integro-differential closure of $\Topcat$. Calculus on topological manifolds and, more generally, the extension of results achieved in the setting of smooth manifolds strictly depends on the knowledge of the structure of $\dcR$. Construction of $\dcR$ reveals that its structure is deeply linked to the local structure of differential equations.\newline 
Part \ref{PII} contains some immediate implications of the theory developed in Part \ref{PI} and takes Chapter \ref{App}-\ref{Qstns}.\newline
In Chapter \ref{App} we define the functor $\genfib$ associating to any local chart of a manifold the corresponding generalized trivial bundle whose fibre is the topological algebra of all generalized derivation operators and to any continuous function between local charts the corresponding generalized differential.
Functor $\genfib$ generalizes functor $\smoothfib$ associating to any local chart of a smooth manifold the corresponding trivial bundle whose fibre is the topological algebra of all derivation operators and to any smooth function between local charts the corresponding differential (Theorem \ref{mainth}-[1]). We also prove that if $\genfib$ is restricted to the category of smooth manifolds and smooth maps then there is a sub-functor $\genfibsm$ which naturally maps onto $\smoothfib$ (Theorem \ref{mainth}-[2, 3, 4]).\newline
Eventually functors $\genbunfunc$ and $\smgenbunfunc$ are defined starting from $\genfib$ and $\genfibsm$ respectively by repeating word by word the argument used to get $\smoothbunfunc$ starting form $\smoothfib$.\newline
In Chapter \ref{chssec} we deal with the problem of finding nowhere vanishing tangent vector fields on spheres. In the classical setting it is well known (see \cite{JWM2}) that nowhere vanishing tangent vector fields do not exists on even dimensional spheres. On the contrary in the generalized setting we exhibit a nowhere vanishing generalized tangent vector field on spheres (Theorem \ref{chsth}). By \cite{JWM2} we have that the projection of such generalized tangent vector field on the classical tangent bundle vanishes somewhere in case of even dimensional spheres (Remark \ref{nullproj}).
The possibility of using in practice the existence of nowhere vanishing tangent vector fields on spheres depends on the wider problem of understanding the physical meaning of the generalized tangent bundle.\newline 
In Chapter \ref{Qstns} we list some questions arising in Mathematics and Physics from this work. \newline
Part \ref{PIII} contains material to recall whenever needed in the work and takes Chapter \ref{relations}, \ref{notat}.\newline
In Chapter \ref{relations} we list those relations among elements belonging to $\fremag$ which are employed in this work. All relations written here are well known to hold true for smooth functions or for continuous functions.\newline
In Chapter \ref{notat} we collect the notation used throughout the work. This section should be consulted only if necessary or whenever it is explicitly mentioned in the work.

\part{Generalized Calculus \label{PI}}

\chapter{Integro-differential spaces\label{DiffSp}}

In this chapter we define the category of integro-differential space and functions.

\bigskip

We introduce the notion of integro-differential monoid. In spite of its abstract appearance an integro-differential monoid is nothing but a structure synthesizing all algebraic properties of integral and differential operators on $\Cksp{\infty}$ functions. 
 
\begin{definition}\label{intdiffmon}
We define the integro-differential monoid $\bsfM$ by setting\\[8pt]
\textbf{Generators:}$\quad \{\Fint_i  \}_{i \in \mathbb{N}} \;\;\cup\;\;\{\Fpart_i  \}_{i \in \mathbb{N}} \;\;\cup\;\;\{\Fp_i  \}_{i \in \mathbb{N}}\;\;\cup\;\;\{\Fq_i  \}_{i \in \mathbb{N}}\;\;\cup\;\;\{\Fqq_i  \}_{i \in \mathbb{N}}$;\\[8pt]
\textbf{Relations:} for any $i,j \in \mathbb{N}$ 
\begin{flalign}
&\begin{array}{ll}
\Fint_i \Fint_j\bsfrelsymb \Fint_{j+1} \Fint_i & \text{if}\; i<j\text{,}\vspace{4pt}\\ 
\end{array}\label{intint}\\[8pt]
&\left\{
\begin{array}{lll}
(i)&\Fpart_{i} \Fpart_j\bsfrelsymb\Fpart_j \Fpart_i \text{,}\vspace{4pt}\\ 
(ii)& \Fpart_{i} \Fint_j\bsfrelsymb\Fint_j \Fpart_i & \text{if}\; i< j \text{,}\vspace{4pt}\\
(iii) &\Fpart_{i+1} \Fint_j\bsfrelsymb \Fint_j \Fpart_i  &\text{if}\; i > j\text{,}\vspace{4pt}\\ 
\end{array}\label{derint}
\right.\\[8pt]
&\left\{
\begin{array}{lll}
(i)& \Fp_i\bsfrelsymb\Fp_1 \Fp_i\text{,}\vspace{4pt}\\ 
(ii)&  \Fp_i\Fint_j\bsfrelsymb\Fint_j\Fp_i\text{,}\vspace{4pt}\\
(iii)&  \Fp_i\Fpart_j\bsfrelsymb\Fpart_j\Fp_i\text{,}\vspace{4pt}\\ 
\end{array}\label{coordint}
\right.\\[8pt]
&\left\{
\begin{array}{llll}
(i)&  \Fq_i\bsfrelsymb\Fpart_{i+1}\Fint_i \text{,}  \vspace{4pt}\\
(ii)&  \Fq_i\Fq_j\bsfrelsymb\Fq_{j+1}\Fq_i & \text{if}\; i\leq j\text{,}  \vspace{4pt}\\
(iii)&  \Fq_i\Fint_{j}\bsfrelsymb\Fint_{j+1}\Fq_i & \text{if}\; i < j\text{,} \vspace{4pt}\\
(iv)&  \Fq_{i}\Fint_j\bsfrelsymb\Fint_{j}\Fq_{i-1} & \text{if}\; i \geq j+2\text{,} \vspace{4pt}\\
(v)&  \Fq_i\Fpart_j\bsfrelsymb\Fpart_{j+1}\Fq_i & \text{if}\; i \leq j\text{,}\vspace{4pt}\\
(vi)&  \Fq_i\Fpart_j\bsfrelsymb\Fpart_{j}\Fq_i & \text{if}\; i > j\text{,}\vspace{4pt}\\
(vii)&  \Fq_i\Fp_j\bsfrelsymb\Fp_{j}\Fq_i \text{,} \vspace{4pt}\\
\end{array}\label{leftproj}
\right.\\[8pt]
&\left\{
\begin{array}{llll}
(i)&  \Fqq_i\bsfrelsymb\Fpart_{i}\Fint_i \text{,}  \vspace{4pt}\\
(ii)&  \Fqq_i\Fqq_j\bsfrelsymb\Fqq_{j+1}\Fqq_i & \text{if}\; i\leq j\text{,}  \vspace{4pt}\\
(iii)&  \Fqq_i\Fint_{j}\bsfrelsymb\Fint_{j+1}\Fqq_i & \text{if}\; i < j\text{,} \vspace{4pt}\\
(iv)&  \Fqq_{i}\Fint_j\bsfrelsymb\Fint_{j}\Fqq_{i-1} & \text{if}\; i \geq j+2\text{,} \vspace{4pt}\\
(v)&  \Fqq_i\Fpart_j\bsfrelsymb\Fpart_{j+1}\Fqq_i & \text{if}\; i < j\text{,}\vspace{4pt}\\
(vi)&  \Fqq_i\Fpart_j\bsfrelsymb\Fpart_{j}\Fqq_i & \text{if}\; i \geq j\text{,} \vspace{4pt}\\
(vii)&  \Fqq_i\Fp_j\bsfrelsymb\Fp_{j}\Fqq_i \text{,} \vspace{4pt}\\
\end{array}\label{rightproj}
\right.\\[8pt]
&\left\{
\begin{array}{llll}
(i)& \Fqq_i\Fq_j\bsfrelsymb\Fq_{j+1}\Fqq_i & \text{if } i<j\text{,}\vspace{4pt}\\
(ii)& \Fqq_i\Fq_j\bsfrelsymb\Fq_{j}\Fqq_{i -1}& \text{if } i>j\text{.}
\end{array} \label{leftrightinter}
 \right.
\end{flalign}
We denote: by $\bsfuno$ the unit of $\bsfM$; by $\bsfm$ the generic element of $\bsfM$; by $\subbdiffsfM$ the sub-monoid of $\bsfM$ generated by the set $\left\{\Fpart_i\right\}_{i\in \mathbb{N}}$; by $\subintsfM$ the sub-monoid of $\bsfM$ generated by the set $\left\{\Fint_i\right\}_{i\in \mathbb{N}}$; by $\subbsfM$ the sub-monoid of $\bsfM$ generated by
$\qquad\{\Fint_i  \}_{i \in \mathbb{N}}  \;\;\cup\;\;\{\Fp_i  \}_{i \in \mathbb{N}}\;\;\cup\;\;\{\Fq_i  \}_{i \in \mathbb{N}}\;\;\cup\;\;\{\Fqq_i  \}_{i \in \mathbb{N}}$.\newline
Monoid $\bsfM$ and its sub-monoids are endowed with the discrete topology (i.e. $\{\bsfm\}$ is an open set for any $\bsfm \in \bsfM $).\newline
With an abuse of language, whenever needed, symbols $\Fint_0$, $\Fpart_{0}$ will mean $\bsfuno$. 
\end{definition}

\begin{remark}\label{intdiffmonrem}\mbox{}
\begin{enumerate}
\item In Definition \ref{intdiffmon} generators $\Fq_i$, $\Fqq_i$ are redundant because of relations \eqref{leftproj}\textnormal{-[}$(i)$\textnormal{]}, \eqref{rightproj}\textnormal{-[}$(i)$\textnormal{]}.
We choose to give such a redundant definition to make more clear constructions and examples.
\item Relation \ref{derint}-[$(i)$] entails that $\subbdiffsfM$ is an abelian monoid with generators $\{\Fpart_i  \}_{i \in \mathbb{N}}$.  
\end{enumerate}
\end{remark}

We introduce the notion of integro-differential spaces. In spite of its abstract appearance an integro-differential space is nothing but a structure synthesizing all algebraic properties of $\Cksp{\infty}$ functions and their relation with integro-differential operators (Example \ref{ExCinf}).\newline
Fix a topological commutative ring $\ringsym$ with unit.

\begin{definition}
\label{DiffSpDef} \mbox{}\\[6pt]
\textnormal{\textbf{Objects.}}\ \ An $\ringsym$-integro-differential space  $\abidsp=\left(\abgenf,\abgenfzero,\abgenfempty ,\abgenfuncomp , \abgenlbound\,\abgenrbound, \abgenfunsum,\abgenfunmult ,\abgenfunscalp , \abbsfmu\right)$ is a
nine-tuple where 
\begin{equation*}
\begin{array}{ccc}
\abgenf\; \text{is a topological space,}& \abgenfzero\text{,}\;\abgenfempty\; \text{are subspaces of}\;\abgenf \text{,}&\abgenfzero\cap \abgenfempty=\udenset  \text{,}
\end{array}
\end{equation*}
\begin{equation*}
\begin{array}{c}
\left.
\begin{array}{l}
\abgenfuncomp, \abgenlbound\,\abgenrbound, \abgenfunsum,\abgenfunmult :\abgenf \!\times\! \abgenf\!\rightarrow \!\abgenf\text{,}\\[4pt]
\abgenfunscalp :\ringsym\!\times
\!\abgenf\!\rightarrow \!\abgenf\text{,}\\[4pt]
\abbsfmu:\bsfM \!\times\! \abgenf\!\rightarrow \!\abgenf
\end{array}
\right\}\qquad\text{are continuous functions,}
\end{array}
\end{equation*}
fulfilling all conditions below $\forall i \in \mathbb{N}$, $\forall \elringsymuno,\elringsymdue\in \ringsym$,  $\forall \bsfm,\bsfn \in \bsfM$, $\forall  
u,v,w,x\in \abgenf$
\begin{flalign}
&\left\{ 
\begin{array}{lll}
(i) & u\abgenfuncomp v \in\abgenfempty\hspace{99pt}\text{if} \;v \in\abgenfempty   \\ 
(ii)& u \abgenfuncomp v \in \abgenfzero \hspace{99pt}\text{if}\;u \in \abgenfzero\text{,}\;u\abgenfuncomp  v\notin \abgenfempty 
\end{array}\label{circle}
\right.\\[8pt]
&\left\{ 
\begin{array}{ll}
(i)&\abgenlbound \abgenlbound u,v   \abgenrbound ,w  \abgenrbound=\abgenlbound u,\abgenlbound v,w   \abgenrbound    \abgenrbound \\
(ii)&\abgenlbound u,v   \abgenrbound =\abgenlbound v,u   \abgenrbound =\abgenfempty\hspace{63pt}\text{if}\; v\in  \abgenfempty \\
(iii)& \abgenlbound u\abgenfuncomp x,v \abgenfuncomp y  \abgenrbound=\abgenlbound u,v   \abgenrbound\abgenfuncomp \abgenlbound  x, y  \abgenrbound \hspace{10pt} \text{if}\;u \abgenfuncomp x,  v \abgenfuncomp y\notin \abgenfempty
\end{array}\label{cartprod}
\right.\\[8pt]
&\left\{ 
\begin{array}{ll}
(i) & \left(u\abgenfunsum v\right)\abgenfunsum w=u\abgenfunsum\left(v\abgenfunsum w\right)   \\ 
(ii) & u\abgenfunsum v=v\abgenfunsum u  \\ 
(iii) & u\abgenfunsum u\in\abgenfempty  \hspace{94pt} \text{if and only if}\;u\in\abgenfempty  \\ 
(iv) & \exists 0_{u}\in \abgenfzero:v\abgenfunsum 0_{u}=v \hspace{47pt} \text{if}\;u \abgenfunsum v\notin \abgenfempty \\
(v) & \left(u\abgenfunsum v\right)\abgenfuncomp w=\left(u\abgenfuncomp w\right)\abgenfunsum\left(v\abgenfuncomp w\right)   \\
\end{array}\label{+}
\right.  \\[8pt]
&\left\{ 
\begin{array}{ll}
(i) & \left(u\abgenfunmult v\right)\abgenfunmult w=u\abgenfunmult  \left(v\abgenfunmult  w\right)   \\ 
(ii) & u\abgenfunmult  v=v \abgenfunmult  u\in\abgenfempty   \hspace{10pt} \text{if}\;v \in \abgenfempty \\ 
(iii)& u\abgenfunmult  v, v\abgenfunmult u \in \abgenfzero \hspace{20pt} \text{if}\;u \in \abgenfzero\text{,\ }u\abgenfunmult  v, v\abgenfunmult  u\notin\abgenfempty\\
(iv) & \left(u\abgenfunmult  v\right)\abgenfuncomp w=\left(u\abgenfuncomp w\right)\abgenfunmult  \left(v\abgenfuncomp w\right)\\
(v) & \left(u\abgenfunsum v\right)\abgenfunmult  w=\left(u\abgenfunmult  w\right)\abgenfunsum\left(v\abgenfunmult  w\right)   \\ 
(vi) & u\abgenfunmult  \left(v\abgenfunsum w\right)=\left(u\abgenfunmult  v\right)\abgenfunsum\left(u\abgenfunmult  w\right) 
\end{array}\label{cdot}
\right. \\[8pt]
&\left\{ 
\begin{array}{ll} 
(i) & \elringsymuno\abgenfunscalp u\in\abgenfempty  \hspace{42pt}\text{if and only if}\;u\in\abgenfempty  \\ 
(ii)& \zeroring\abgenfunscalp u\in \abgenfzero \hspace{42pt}  \text{if}\; u\neq\abgenfempty\\
(iii)& \unoring\genfunscalp u=u  \\
(iv) & \elringsymuno\abgenfunscalp \left(\elringsymdue\abgenfunscalp u\right)=\left(\elringsymuno\elringsymdue\right)\abgenfunscalp  u  \\ 
(v) & \left(\elringsymuno\abgenfunscalp  u\right)\abgenfuncomp v=\elringsymuno\abgenfunscalp  \left(u\abgenfuncomp v\right)\\
(vi) & \elringsymuno\abgenfunscalp\abgenlbound u,v  \abgenrbound=\abgenlbound \elringsymuno\abgenfunscalp u,\elringsymuno\abgenfunscalp v  \abgenrbound\\
(vii) & \left(\elringsymuno+\elringsymuno\right)\abgenfunscalp  u=\left(\elringsymuno\abgenfunscalp  u\right)\abgenfunsum\left(\elringsymuno\abgenfunscalp  u\right)  \\ 
(viii) & \elringsymuno\abgenfunscalp  \left(u\abgenfunsum v\right)=\left(\elringsymuno\abgenfunscalp  u\right)\abgenfunsum\left(\elringsymuno\abgenfunscalp  v\right)  \\ 
(ix) & \left(\elringsymuno\abgenfunscalp  u\right)\abgenfunmult  v=u\abgenfunmult  \left(\elringsymuno\abgenfunscalp  v\right)=\elringsymuno\abgenfunscalp  \left(u\abgenfunmult  v\right)   \\ 
\end{array}\label{star}
\right.\\[8pt]
&\left\{ 
\begin{array}{ll}
\!(i)\! &\! \bsfm \abbsfmu \left(\bsfn \abbsfmu u\right) =\left(\bsfm \bsfn\right) \abbsfmu u\\
\!(ii)\!&\! \bsfuno \abbsfmu u =u\\
\!(iii)\! & \!\bsfm \abbsfmu u \in \abgenfempty \hspace{20pt} \text{if and only if}\;u\in\abgenfempty  \\
\!(iv)\!& \!\bsfm \abbsfmu \left(u\abgenfunsum v\right) =\left(\bsfm \abbsfmu u\right)\abgenfunsum \left(\bsfm \abbsfmu v\right) \\
\!(v)\!& \!\bsfm \abbsfmu \left(a \abgenfunscalp  u\right) =a \abgenfunscalp  \left(\bsfm \abbsfmu u\right) \\
\!(vi)\! & \!\left(\Fpart_i \Fq_i\right) \abbsfmu u ,\left(\Fpart_{i+1} \Fqq_i\right) \abbsfmu u \in \abgenfzero\hspace{59pt} \text{if}\; u \notin\abgenfempty\\
\!(vii)\! &\!\left(\Fp_i \Fp_j\right) \abbsfmu u \in \abgenfzero \hspace{114pt} \text{if}\; u\notin \abgenfempty \text{,}\; i \neq 1\\
\!(viii)\! & \!\left(\Fpart_i \abbsfmu u\right) \abgenfunsum u \notin \abgenfempty \hspace{104pt} \text{if}\;u\notin \abgenfempty\\
\!(ix)\!& \!\Fpart_i \abbsfmu \left(u \abgenfunmult  v\right)= \left(\Fpart_i \abbsfmu u \right)\abgenfunmult  v \abgenfunsum  u \abgenfunmult  \left(\Fpart_i \abbsfmu v\right)\\
\!(x)\!& \!\exists m \in \mathbb{N} :\hspace{10pt} \Fpart_i \abbsfmu \left(u \abgenfuncomp v\right)=\overset{m}{\underset{\mathsf{m}=1}{\bigabgenfunsum}}\left(\left(\Fpart_{\mathsf{m}} \abbsfmu u \right)\abgenfuncomp v\right)\abgenfunmult  \left(\left(\Fp_{\mathsf{m}} \Fpart_i\right)\abbsfmu v\right)\\ 
\!(xi)\!& \!\Fp_i \abbsfmu \left(u \abgenfuncomp v\right)= \left(\Fp_i \abbsfmu u\right) \abgenfuncomp v\\
\!(xii)\!&\!\Fq_{i} \abbsfmu u = \left( \left(\Fint_{i} \Fpart_i\right)\abbsfmu u \right)\abgenfunsum \left(\left(-1\right) \abgenfunscalp \left(\Fqq_i\abbsfmu u\right)\right)\\
\!(xiii)\!& \!\Fq_i \abbsfmu \left(u \abgenfuncomp v\right)=u \abgenfuncomp \left(\Fq_i \abbsfmu v\right)\\
\!(xiv)\! & \!\Fq_i \left(u \abgenfunmult  v\right)= \left(\Fq_i \abbsfmu u\right)\abgenfunmult \left(\Fq_i \abbsfmu v\right)\\
\!(xv)\! &\!\Fqq_i \abbsfmu \left(u \abgenfuncomp v\right)\abgenfunsum u \abgenfuncomp \left(\Fqq_i \abbsfmu\left( -1\abgenfunscalp v\right)\right)\in \abgenfzero
\hspace{10pt} \text{if}\;u\abgenfuncomp v \notin \abgenfempty \\ 
\!(xvi)\! &\! \Fqq_i \left(u \abgenfunmult  v\right)\abgenfunsum \left(\Fqq_i \abbsfmu u\right)\abgenfunmult \left(\Fqq_i \abbsfmu v\right)\in X_0 \hspace{25pt} \text{if}\;u\abgenfuncomp v \notin \abgenfempty\\
\!(xvii)\! &\! \Fint_i \abbsfmu \left(u \abgenfunmult  \left(\Fq_i \abbsfmu v\right)\right)= \left(\Fint_i \abbsfmu u\right) \abgenfunmult  \left(\left(\Fq_i \Fq_i\right) \abbsfmu v\right)\\
\!(xviii) \!&\!\Fint_{i+1} \abbsfmu \left(u \abgenfunmult  \left(\Fqq_i \abbsfmu v\right)\right)= \\
&\hspace{51pt}\left(\Fint_{i+1} \abbsfmu u\right) \abgenfunmult  \left(\left(\Fqq_i\Fqq_i\right) \abbsfmu \left( -1 \genfunscalp  v\right)\right)\\
\!(ixx)\! &\! \left(\Fq_{i}\Fint_i\right)\abbsfmu u = \left(\Fq_{i+1}\Fint_i\right)\abbsfmu u \abgenfunsum  \left(\Fqq_{i+1}\Fint_i\right)\abbsfmu u\\
\!(xx)\! &\! \left(\Fqq_{i}\Fint_i\right)\abbsfmu u \abgenfunsum \left(\Fq_{i+1}\Fint_i\right)\abbsfmu u\in X_0\hspace{44pt} \text{if}\; u\notin \abgenfempty\\
\!(xxi)\! & \!\left(\Fqq_{i}\Fq_i\right)\abbsfmu u \abgenfunsum \left(\Fq_{i}\Fq_i\right)\abbsfmu u\in X_0\hspace{61pt} \text{if}\; u\notin \abgenfempty\text{.}
\end{array}\label{intdiffop}
\right. 
\end{flalign}
We call: $\abgenfempty $ empty element of $\abidsp$; $\abgenfuncomp $ composition; $\abgenlbound \;\abgenrbound$ cartesian product; $\abgenfunsum$ sum; $\abgenfunmult$ matrix product;  $\abgenfunscalp$
product by a scalar; $\abbsfmu$ integro-differential operation.\newline
With an abuse of language, whenever no confusion is possible, we use the
notation $u\in \abidsp$ in place of $u\in \abgenf$ to denote an element of the $
\ringsym$-integro-differential space $\abidsp$.\\[6pt]
\textnormal{\textbf{Arrows.}}\ \ Fix two $\ringsym$-integro-differential spaces 
\begin{equation*}
\abidsp=\left(\abgenf,\abgenfzero,\abgenfempty ,\abgenfuncomp , \abgenlbound\,\abgenrbound, \abgenfunsum,\abgenfunmult  ,\abgenfunscalp  , \abbsfmu\right)\text{,}\qquad 
\abidsptwo=\left(\abgenftwo,\abgenfzerotwo,\abgenfemptytwo,
\abgenfuncomptwo, \abgenlboundtwo\,\abgenrboundtwo, \abgenfunsumtwo,\abgenfunmulttwo,\abgenfunscalptwo, \abbsfmutwo\right)\text{.}
\end{equation*}
A function between $\ringsym$-integro-differential spaces is a continuous function $f :\abgenf\rightarrow \abgenftwo$ fulfilling the system of conditions below
\begin{equation*}
\left\{
\begin{array}{ll}
f\left(\abgenfzero\right)\subseteq \abgenfzerotwo\text{,}\\[4pt]
f\left(\abgenfempty\right)\subseteq \abgenfemptytwo\text{,}\\[4pt]
f\left(u\abgenfuncomp  v\right)=f\left(u\right)\,\abgenfuncomptwo\,f\left(v\right)&\forall u,v\in \abgenf\text{,}\\[4pt]
f \left(\abgenlbound u,  v\abgenrbound \right)=\abgenlboundtwo f\left(u\right),f \left(v\right)\abgenrboundtwo&\forall u,v\in \abgenf\text{,}\\[4pt]
f\left(u\abgenfunsum v\right)=f\left(u\right)\,\abgenfunsumtwo\,f\left(v\right)& \forall u,v\in \abgenf\text{,}\\[4pt]
f\left(u\abgenfunmult  v\right)=f \left(u\right)\,\abgenfunmulttwo\,f\left(v\right)&\forall u,v\in \abgenf\text{,}\\[4pt]
f\left(a\abgenfunscalp  u\right)=a\,\abgenfunscalptwo\,f\left(u\right)&\forall a\in \ringsym\text{,} \;\;\forall u\in \abgenf\text{,}\\[4pt]
f\left(\bsfm \abbsfmu u\right)=\bsfm \,\abbsfmutwo\,f\left(u\right)&     \forall \bsfm \in \bsfM\text{,}\;\;  \forall u,v\in \abgenf\text{.}
\end{array}
\right.
\end{equation*}

\textnormal{\textbf{Category.}}\ \ We denote by $\IDSMa _{\ringsym}$ the category of $\ringsym$
-integro-differential spaces and functions.
\end{definition}

Example \ref{ExCinf} below shows that, in spite of its abstract appearance, Definition \ref{DiffSpDef} is a quite natural one, more precisely: the set of all $\Cksp{\infty}$ functions is an $\mathbb{R}$-integro-differential space in a very natural way. We refer to Notations \ref{ins}-[1, 5, 7], \ref{gentopnot}-[5, 7, 8, 9],  \ref{realvec}-[6], \ref{realfunc}-[3(a), 4], \ref{difrealfunc}-[1, 2, 3].  

\begin{example}
\label{ExCinf} 
We define the eight-tuple $\idssmo\!=\!\left(\Cksp{\infty},\Ckspzero{\infty}, \left\{\emptysetfun\right\} ,\funcomp, \smosplbound\,\smosprbound , \funsum,\funmult, \funscalp ,\smospbsfmu\right)$.
\begin{description}
\item[\textnormal{Definition of} $\Cksp{\infty}$\textnormal{,} $\Ckspzero{\infty}$\textnormal{.}] We set:
\begin{flalign*}
&\Cksp{\infty}=\left\{ f\;:\;f\in \Cksp{\infty}(\intuno)+\mathbb{R}^{n}+\quad \forall m,n\in \mathbb{N}_{0}\quad \forall \intuno\subseteqdentro \mathbb{R}^{m}\right\}\cup\{\emptysetfun \} \text{;}\\[4pt]
&\Ckspzero{\infty}=\left\{ \cost<\intuno<>\mathbb{R}^n>+0+\quad \forall m,n\in \mathbb{N}_{0}\quad \forall \intuno\subseteqdentro \mathbb{R}^{m}\right\} \text{.}
\end{flalign*}
\item[\textnormal{Definition of} $\funcomp$\textnormal{.}] Fix $f,g\in \Cksp{\infty}$. \newline
If $\Ima\left(g\right)\cap  \Dom\left(g\right)=\udenset$, then we set $f\funcomp g=\emptysetfun$.\newline
If $\Ima\left(g\right)\cap \Dom\left(f\right)\neq\udenset$, then we set
\begin{equation*}
\left(f\funcomp g\right)\left(\unkuno\right)= f\left(g\left(\unkuno\right)\right) \quad \forall \unkuno \in g^{-1}\left(\Dom\left(f\right)\right)\text{.}
\end{equation*}
\item[\textnormal{Definition of } $\smosplbound\,\smosprbound$\textnormal{.}] Fix $f,g\in \Cksp{\infty}$.\newline 
If $f = \emptysetfun$ or $g= \emptysetfun$ then we set $\smosplbound f, g\smosprbound=\emptysetfun$.\newline
If $f \neq \emptysetfun$ and $g\neq \emptysetfun$ then we set
\begin{equation*}
\smosplbound f, g\smosprbound \left(\unkuno,\unkdue\right)= \smosplbound f\left(\unkuno\right), g\left(\unkdue\right)\smosprbound  \quad \forall \left(\unkuno,\unkdue\right) \in \Dom\left(f\right)\times \Dom\left(g\right)\text{.}
\end{equation*}
\item[\textnormal{Definition of } $\funsum$\textnormal{.}]Fix $f,g\in \Cksp{\infty}$. \newline
If $f = \emptysetfun$ or $g= \emptysetfun$ or $\Dom\left(f\right)\neq \Dom\left(g\right)$ or $\Cod\left(f\right)\neq\Cod\left(g\right)$,\newline
then we set $f+g=\emptysetfun$.\newline
If $f \neq \emptysetfun$ and $g\neq \emptysetfun$ and $\Dom\left(f\right)=\Dom\left(g\right)$ and $\Cod\left(f\right)=\Cod\left(g\right)=\mathbb{R}^n$, then we set \quad $\left(f+g\right)\left(\unkuno\right)= \vecsum<n<>2>\left(f\left(\unkuno\right),g\left(\unkuno\right)\right) \quad \forall \unkuno \in \Dom\left(f\right)$.
\item[\textnormal{Definition of} $\funmult$\textnormal{.}] Fix $f,g\in \Cksp{\infty}$. \newline
If $f = \emptysetfun$ or $g= \emptysetfun$ or $\Dom\left(f\right)\neq \Dom\left(g\right)$, then we set $f\funmult g=\emptysetfun$.\newline
If $f \neq \emptysetfun$ and $g\neq \emptysetfun$ and $\Dom\left(f\right)=\Dom\left(g\right)$ and $\Cod\left(f\right)= \mathbb{R}^m$ and $\Cod\left(g\right)=\mathbb{R}^n$ then we set $\left(f\!\funmult\! g\right)\left(\unkuno\right)= \vecprod\left[m,n\right]\left(f\left(\unkuno\right),g\left(\unkuno\right)\right) \quad \forall \unkuno \in \Dom\left(f\right)$.
\item[\textnormal{Definition of } $\funscalp$\textnormal{.}] Fix $f\in \Cksp{\infty}$, $a \in \mathbb{R}$. \newline
If $f = \emptysetfun$ then we set $a \funscalp f=\emptysetfun$.\newline
If $f \neq \emptysetfun$ then we set \quad $\left(a \funscalp f\right)\left(\unkuno\right)= a f\left(\unkuno\right) \quad \forall \unkuno \in \Dom\left(g\right)$.
\item[\textnormal{Definition of } $\smospbsfmu$\textnormal{.}] Fix $f\in \Cksp{\infty}$, $i \in \mathbb{N}$.\newline
If $f = \emptysetfun$ then we set $\Fint_i \smospbsfmu f=\Fpart_i \smospbsfmu f=\Fp_i \smospbsfmu f=\Fq_i \smospbsfmu f=\Fqq_i \smospbsfmu f=\emptysetfun$.\newline
If $f \neq \emptysetfun$ then there are $m,n \in \mathbb{N}_0$, $\intuno\subseteqdentro \mathbb{R}^m$ such that $\Dom\left(f\right)=\intuno$ and $\Cod\left(f\right)=\mathbb{R}^n$, hence we set: 
\begin{flalign*}
&\begin{array}{l}
\Fint_i \smospbsfmu f=\left\{
\begin{array}{ll}
\left( \smint_i (f_1),...,\smint_i(f_{\overline{n}})\right)\funcomp\Diag{\domint\left[\intuno,i\right]}{\overline{n}} & \text{if}\; i \leq m\vspace{4pt}\\
\left( \smint_i (g_1),...,\smint_i(g_{\overline{n}})\right)\funcomp\Diag{\domint\left[\intuno,i\right]}{\overline{n}} & \text{if}\; i > m\text{,}
\end{array}\right.\vspace{4pt}\\ 
\text{where}\qquad g_{\mathsf{n}}=f_{\mathsf{n}} \funcomp \proj<\intuno,\mathbb{R}^{i-m}<>1> \quad \forall {\mathsf{n}}\in\{1,...,\overline{n}\}\text{;}
\end{array}\\[8pt]
&\begin{array}{l}
\Fpart_i \smospbsfmu f= \left\{
\begin{array}{ll} 
\left( \smder_i (f_1),...,\smder_i(f_{\overline{n}})\right)\funcomp\Diag{\intuno}{\overline{n}} & \text{if}\; i \leq m\text{,}\vspace{4pt}\\
\cost<\intuno<>\mathbb{R}^n>+0+  & \text{if}\; i > m\text{;}
\end{array}\right.
\end{array}\\[8pt]
&\begin{array}{l}
\Fp_i \smospbsfmu f= \left\{
\begin{array}{ll}
f_i & \text{if}\; i \leq \overline{n}\text{,}\vspace{4pt}\\
\cost<\intuno<>\mathbb{R}>+0+  & \text{if}\; i > \overline{n}\text{;}
\end{array}\right.
\end{array}\\[8pt]
&\begin{array}{l}
\Fq_i \smospbsfmu f= \left\{
\begin{array}{ll}
f \funcomp\proje{m}{i}\funcomp\incl{\domint\left[\intuno,i\right]}{\mathbb{R}^{m+1}}& \text{if}\; i \leq m\text{,}\vspace{4pt}\\
f \funcomp \proj<\intuno, \mathbb{R}^{i-m+1}<>1> & \text{if}\; i > m\text{;}
\end{array}\right.
\end{array}\\[8pt]
&\begin{array}{l}
\Fqq_i \smospbsfmu f= \left\{
\begin{array}{ll}
-f \funcomp\proje{m}{i+1}\funcomp\incl{\domint\left[\intuno,i\right]}{\mathbb{R}^{m+1}}& \text{if}\; i \leq m\text{,}\vspace{4pt}\\
-f \funcomp \proj<\intuno, \mathbb{R}^{i-m+1}<>1> & \text{if}\; i > m\text{.}
\end{array}\right.
\end{array}
\end{flalign*}
\end{description}
It is a long but not difficult exercise to check that conditions \eqref{circle}-\eqref{intdiffop} are fulfilled.\newline
Topology on $\Cksp{\infty}$ is denoted by $\topS$ and is the smallest topology containing all sets listed below 
\begin{multline*}
\setsymotto\left[\setsymcinque,\setsymsei,N\right]=\\
\left\{f\in \Cksp{\infty}\;:\; \setsymcinque\subseteq\left(\smder_{\alpha}\left(f\right)\right)^{-1}\left(\setsymsei\right)\quad \forall \alpha \in \left(\mathbb{N}_0\right)^m\;\text{with}\;\overset{m}{\underset{{\mathsf{m}}=1}{\sum}}\,\alpha_{\mathsf{m}}\leq N\right\}
\end{multline*}
where $m,n,N\in \mathbb{N}_0$, $\setsymcinque\subseteq\mathbb{R}^m$ is a compact set, $\setsymsei\subseteqdentro\mathbb{R}^n$, $\smder_{\alpha}$ is the differential operator given by the composition $\left(\smder_1\right)^{\alpha_1}...\left(\smder_m\right)^{\alpha_m}$.\newline
We list below some facts about topology $\topS$ which will be used in the next sections.\newline
Fix a set $\setsymotto\left[\setsymcinque,\setsymsei,N\right]$, $f \in\setsymotto\left[\setsymcinque,\setsymsei,N\right]$ then $\setsymcinque\subseteq \Cae\left(f\right)$. This entails that   
\begin{equation}
 \begin{array}{l}
\text{for}\hspace{7pt}\text{any}\hspace{7pt}m,n\in\mathbb{N}_0\text{,}\hspace{7pt} \intuno \subseteqdentro\mathbb{R}^m\hspace{7pt}\text{the}\hspace{7pt}\text{subspace}\hspace{7pt}\text{topology}\hspace{7pt}\text{induced}\hspace{7pt}\text{on}\\[4pt]
 \Cksp{\infty}(\intuno)+\mathbb{R}^{n}+\hspace{4pt}
\text{coincides}\hspace{4pt}\text{with}\hspace{4pt}\text{topology}\hspace{4pt}\text{recalled}\hspace{4pt}\text{in}\hspace{4pt}\text{Notation}\hspace{4pt}\text{\ref{difrealfunc}-[2].}
\end{array}\label{extopsmdis2}
\end{equation}
Fix $N \in \mathbb{N}_0$, $f\in \Cksp{\infty}\setminus\left\{\emptysetfun\right\}$, a compact set $\setsymcinque\subseteq \Cae\left(f\right)$ then straightforwardly is $f \in\setsymotto\left[\setsymcinque,\Cod\left(f\right),N\right]$ and $\emptysetfun \notin\setsymotto\left[\setsymcinque,\Cod\left(f\right),N\right]$ since $\emptysetfun\left(\unkuno\right)=\udenunk\notin \Cod\left(f\right)$. This entails that the set $\{\emptysetfun\}$ is a closed subset of $\Cksp{\infty}$.
\begin{flalign}
&\begin{array}{l}
\left\{ 
\begin{array}{l}
\text{Set functions}\vspace{4pt}\\
\hspace{60pt}\funcomp, \smosplbound\,\smosprbound, \funsum,\funmult :\Cksp{\infty}\times \Cksp{\infty}\rightarrow \Cksp{\infty}\vspace{4pt}\\
\hspace{20pt}\funscalp :\mathbb{R}\times\Cksp{\infty}\rightarrow \Cksp{\infty}\hspace{20pt}\smospbsfmu :\bsfM\times \Cksp{\infty}\rightarrow \Cksp{\infty}\vspace{4pt}\\  
\text{are continuous functions.}
\end{array}
\right.
\end{array}\label{extopsmdis3}
\end{flalign}
Topology $\topS$ allows to handle convergence of a sequence of smooth functions $\left\{f_n\right\}_{n \in \mathbb{N}}$ to a smooth function $f$ even if sets $\Cae\left(f_n\right)$ for $n \in \mathbb{N}$, $\Cae\left(f\right)$ do not coincide. Precisely $\underset{n\rightarrow\infty}{\lim}\,f_n=f$ means that for any compact set $\setsymcinque\subseteq \Cae\left(f\right)$ there is $n_{\setsymcinque}$ such that both conditions below are fulfilled:
\begin{flalign*}
&\begin{array}{l}\setsymcinque\subseteq \Cae\left(f_n\right) \quad\forall n \geq n_{\setsymcinque}\text{;}\end{array}\\[6pt]
&\left\{\begin{array}{l}\left\{f_n\funcomp \incl{\setsymcinque}{\Cae\left(f_n\right)}\right\}_{n \geq n_{\setsymcinque}}
\text{uniformly converges to}\; f\funcomp \incl{\setsymcinque}{\Cae\left(f\right)}\\[4pt]\text{with respect to}\;\|\,\|_{k}\;\text{for any}\;k \in \mathbb{N}\text{.}\end{array}\right.
\end{flalign*}
Topology $\topS$ coincides with the final topology on $\Cksp{\infty}$ with respect to the family of all paths in $\Cksp{\infty}$.\newline
The $\mathbb{R}$-integro-differential space $\idssmo$ fulfills many other properties in addition to \eqref{circle}-\eqref{intdiffop}. Such additional properties cannot be fulfilled by an $\mathbb{R}$-integro-differential space if we pretend that it contains all continuous functions. This fact, which was first proven in \cite{LS}, turns out to be intrinsic in construction described in Chapter \ref{costruct}.        
\end{example}

By restricting our attention to the sub-monoid $\subbsfM$ of $\bsfM$ and by dropping suitable conditions in Definition \ref{DiffSpDef} we obtain integral spaces.

\begin{definition}
\label{QDiffSpDef} \mbox{}\\[6pt]
\textnormal{\textbf{Objects.}}\ An $\ringsym$-integral space
is a eight-tuple $\abidsp\!=\!\left(\abgenf,\abgenfzero,\abgenfempty ,\abgenfuncomp ,  \abgenlbound\,\abgenrbound,  \abgenfunsum ,\abgenfunmult  ,\abgenfunscalp, \abbsfmu \right)$
fulfilling word by word all conditions contained in Definition \ref{DiffSpDef} where $\bsfM$ is replaced by $\subbsfM$ and any reference to $\Fpart_i$ for any $i\in\mathbb{N}$ is dropped.\\[6pt]
\textnormal{\textbf{Arrows.}}\ Arrows between $\ringsym$-integral spaces are defined word by word as in Definition 
\ref{DiffSpDef} where $\bsfM$ is replaced by $\subbsfM$.\\[6pt] 
\textnormal{\textbf{Category.}}\ We denote by $\QIDSMa_{\ringsym}$ the category of $\ringsym$
-integral spaces and functions.\newline
We denote by $\ffIDSMaQIDSMa:\IDSMa_{\ringsym}\rightarrow \QIDSMa_{\ringsym}$ the forgetful functor.
\end{definition}

Example \ref{ExCzero} shows that, in spite of its abstract appearance, Definition \ref{DiffSpDef} is a quite natural one, more precisely: the set of all $C^0$ functions is an $\mathbb{R}$-integral space in a very natural way. We refer 
to Notations \ref{ins}-[7], \ref{gentopnot}-[5, 7, 8, 9], \ref{difrealfunc}-[1, 2], Example \ref{ExCinf}.

\begin{example}\label{ExCzero} 
We define the eight-tuple $\qidscont=(\Cksp{0},\Ckspzero{0},\left\{\emptysetfun\right\} ,\funcomp,\contsplbound \,\contsprbound, \funsum ,\funmult,\funscalp, \contspbsfmu)$.
\begin{description}
\item[\textnormal{Definition of} $\Cksp{0}$\textnormal{,} $\Ckspzero{0}$\textnormal{,} $\emptysetfun$\textnormal{.}]
We set:
\begin{flalign*}
&\Cksp{0}=\big\{ f\;:\;f\in \Cksp{0}(\intuno)+\mathbb{R}^{n}+\quad \forall m,n\in \mathbb{N}_{0}\quad \forall \intuno\subseteqdentro \mathbb{R}^{m}\big\}\cup\{\emptysetfun \} \text{;}\\[4pt]
&\Ckspzero{0}=\Ckspzero{\infty}\text{.} 
\end{flalign*}
\item[\textnormal{Definition of} $\funcomp$\textnormal{,}  $\contsplbound \,\contsprbound $\textnormal{,} $\funsum$\textnormal{,} $\funscalp$\textnormal{,} $\funmult$\textnormal{,}  $\contspbsfmu$\textnormal{.}] Set functions $\funcomp$, $\contsplbound \,\contsprbound $, $\funsum$, $\funscalp$, $\funmult$,  $\contspbsfmu$ are defined word by word as in Example \ref{ExCinf} where $f,g \in \Cksp{\infty}$ is replaced by $f,g \in \Cksp{0}$, $\bsfm \in \bsfM$ is replaced by $\bsfm \in \subbsfM$.
\end{description} 
It is a long but not difficult exercise to prove that $\qidscont$ is an integral space by checking conditions \eqref{circle}-\eqref{star}.\newline
Topology on $\Cksp{0}$ is denoted by $\topC$ and is given word by word as in Example \ref{ExCinf} where differential operator $\smder_{\alpha}$ is replaced by $\idobj+\Cksp{0}+$.\newline 
We list below some facts about topology $\topC$ which will be used in the next sections. The proof of these facts is either straightforward or is word by word the same as in Example \ref{ExCinf}.\newline
For any $m,n\in\mathbb{N}_0$, $\intuno \subseteqdentro\mathbb{R}^m$ the subspace topology induced on $C^{0 }\left(\intuno,\mathbb{R}^{n}\right)$ coincides with topology recalled in Notation \ref{difrealfunc}-[2].\newline
\begin{flalign}
&\begin{array}{l}
\text{The set}\;\{\emptysetfun\}\;\text{is a closed subset of}\;\Cksp{0}\text{.}\label{extopcntdis1}
\end{array}\\[8pt]
& \left\{
\begin{array}{l}
\text{Set functions}\\[4pt]
\hspace{80pt}\funcomp, \contsplbound \,\contsprbound ,\funsum,\funmult :\Cksp{0}\times \Cksp{0} \rightarrow \Cksp{0}\text{,} \vspace{4pt}\\
\hspace{55pt}
\funscalp :\mathbb{R}\times\Cksp{0}\rightarrow \Cksp{0}\text{,}\hspace{15pt}\contspbsfmu :\subbsfM\times \Cksp{0}\rightarrow \Cksp{0}\vspace{4pt}\\ 
\text{are continuous functions.}
\end{array}
\label{extopcntdis3}
\right.
\end{flalign}
Inclusion $\incl{\Cksp{\infty}}{\Cksp{0}}$ is a continuous function.\newline
Weierstrass approximation theorem entails that  
\begin{equation}
\Cksp{\infty}\; \text{is dense in}\;\Cksp{0}\text{.} 
\label{extopcntdis5}
\end{equation}
Inclusion $\incl{\Cksp{\infty}}{\Cksp{0}}$ induces a function $\smofunconfun: \ffIDSMaQIDSMa(\idssmo)\rightarrow \qidscont$ between 
integral spaces.\newline
Topology $\topC$ allows to handle convergence of a sequence of smooth functions $\left\{f_n\right\}_{n \in \mathbb{N}}$ to a smooth function $f$ even if sets $\Cae\left(f_n\right)$ for $n \in \mathbb{N}$, $\Cae\left(f\right)$ do not coincide.\newline
Topology $\topC$ coincides with the final topology on $\Cksp{0}$ with respect to the family of all paths in $\Cksp{0}$.\newline
We emphasize that for any $f\in \Cksp{0}$ notion of $\Dom\left(f\right)$ and of $\Cae\left(f\right)$ coincide.
\end{example}

\chapter{Construction of the integro-differential space of generalized functions\label{costruct}}
In this chapter we explicitly construct an integro-differential space $\genfidsp$ which extend in a suitable way the integro-differential space of $ \idssmo $ and the integral space of $\qidscont$.\newline
Throughout this section we will make an extensive use of magmas. We refer to Notation \ref{alg}, Examples \ref{ExCinf}, \ref{ExCzero}, \cite{Bour}\;\;Chap\,I\;\;\textsection\,1\;\;N$^{\circ}$\,1, \cite{Bour}\;\;Chap\,I\;\;\textsection\,7\;\;N$^{\circ}$\,1.\newline
Construction of integro-differential space $\genfidsp$ splits into five steps.

\section{Choice of generating functions \label{cgf}}
For any $m\in \mathbb{N}$, $\intuno\subseteqdentro \mathbb{R}^{m}$ we choose: 
\begin{equation}
\left\{
\begin{array}{l}
\text{an algebraic base}\;\bascont^{\infty }\left(\intuno\right)\text{\ for the vector space }\Cksp{\infty}(\intuno)+\mathbb{R}+\text{;}\vspace{4pt} \\ 
\text{a subset}\;\bascont\left(\intuno\right)\subseteq C^{0}\left(\intuno,\mathbb{R}\right)\;\text{
such that}\;\bascont\left(\intuno\right)\cup \bascont^{\infty }\left(\intuno\right)\;\text{is an}\\ 
\text{algebraic base for the vector space}\;C^{0}\left(\intuno,\mathbb{R}\right)\text{.}
\end{array}
\right.
\label{Bchoice}
\end{equation}
We set: 
\begin{equation}\label{Bchoice2}
\bascont=\underset{m\in\mathbb{N}}{\bigcup}\left(\underset{\intuno\subseteqdentro \mathbb{R}^m}{\bigcup} \bascont\left(\intuno\right)\right)\text{,} \hspace{40pt}
\fremag_{-1}=\left(\Cksp{\infty}\setminus\left\{\smospempty\right\}\right)\cup \bascont\text{.}
\end{equation}

\begin{remark}\label{remch}\mbox{}
\begin{enumerate}
\item We do not consider vector spaces $C^0\left(\intuno,\mathbb{R}^0\right)$, $\Cksp{\infty}(\intuno)+\mathbb{R}^0+$, $C^0(\mathbb{R}^0,\mathbb{R})$, $\Cksp{\infty}(\mathbb{R}^0)+\mathbb{R}+$ in \eqref{Bchoice} since we have that $C^0\left(\intuno,\mathbb{R}^0\right)=\Cksp{\infty}(\intuno)+\mathbb{R}^0+$ and $C^0(\mathbb{R}^0,\mathbb{R})= \Cksp{\infty}(\mathbb{R}^0)+\mathbb{R}+$.    
\item Construction of the integro-differential space $\genfidsp$ does not depend on choice \eqref{Bchoice}. Two different choices of generators $\bascont_1$, $\bascont_2$ in \eqref{Bchoice} give two different integro-differential spaces $\genfidsp_1$, $\genfidsp_2$ respectively which are isomorphic in $\IDSMa_{\mathbb{R}}$: any isomorphism of sets $f:\bascont_1\rightarrow \bascont_2$ fulfilling $f\left(\bascont_1\left(\intuno\right)\right)\subseteq\bascont_2\left(\intuno\right)$ for any $m\in \mathbb{N}$, $\intuno\subseteqdentro \mathbb{R}^{m}$ gives an isomorphism $\genfidsp_1\rightarrow \genfidsp_2$ of integro-differential spaces. This will be clear from the construction.
\item The unconventional choice of algebraic base for $\mathbb{R}$-vector spaces $C^{0}\left(\intuno,\mathbb{R}\right)$ allows to handle a finite number of generators when dealing with a finite number of generalized function. This turns out to be very useful both in handling generalized functions from the algebraic point of view and in defining paths in the space of generalized functions, which are crucial in order to introduce topology.       
\end{enumerate}
\end{remark}

\section{Free magma $(\fremag,\compone)$\label{fremag}}
Free magma $(\fremag,\compone)$ is the direct limit of free magmas 
$\{(\fremag_{n},\compone_{\!n})\}_{n \in \mathbb{N}_0}$ which are defined inductively.\newline
In Definition \ref{basind} below we establish the base of the inductive process.
\begin{definition}\label{basind}
We set $\fremaggen_0=\underset{n\in\mathbb{N}}{\bigcup} \left(\fremag_{-1}\right)^n$. The generic element belonging to $\fremaggen_0$ is denoted by $\fmx$; elements belonging to $\fremag_{-1}$ will be denoted also by $f$; for any fxied $n \in \mathbb{N}$ elements belonging to $\left(\fremag_{-1}\right)^n$ are denoted also by $\left( \fmx_1,...,\fmx_n \right)$.\newline
$(\fremag_0,\compone_{\!0})$ is the free magma generated by the set $\fremaggen_0$.
\end{definition}

\begin{proposition}\label{stepind}
Assume that exists $k \in \mathbb{N}$ such that free magmas $(\fremag_{j}, \compone_{\!j})$ with generators $\fremaggen_j$ and inclusion of magmas $\inclfremaglev_j:\fremag_{j-1}\rightarrow \fremag_{j}$ are defined for any $j \in \{0,...,k-1\}$.\newline
Set $\fremaggen_{k}=\underset{n\in \mathbb{N}}{\bigcup}\left(\bsfM \times \fremag_{k-1}\right)^n$. The generic element belonging to $\fremaggen_k$ is denoted by $\fmx$; for any fixed
$n \in \mathbb{N}$ elements belonging to $\left(\bsfM \times \fremag_{k-1}\right)^n$ are also denoted by $\left( \fmx_1,...,\fmx_n \right)$ where $\fmx_i \in \bsfM \times \fremag_{k-1}$ for any $i \in \left\{1,...,n\right\}$.\newline
Then: 
\begin{flalign*}
&\begin{array}{l}
\text{there is one and only one free magma}\; \left(\fremag_{k},\compone_{\!k}\right)\; \text{generated by}\;\fremaggen_{k}\text{;}
\end{array}\\[8pt]
&\begin{array}{l}
\left\{
\begin{array}{l}
\text{there is one and only one injective function of magmas}\; \inclfremaglev_k:\fremag_{k-1}\rightarrow \fremag_{k}\\[4pt]
\text{defined by setting}
\hspace{10pt}\inclfremaglev_k\left(\fmx\right)= \left(\bsfuno, \fmx\right) \quad \forall \fmx \in \fremaggen_{k-1}\text{.}
\end{array}
\right.
\end{array}
\end{flalign*}
\end{proposition}
\begin{proof}
Construction of a free magma starting from a given set of generators is described in \cite{Bour}\;\;Chap\,I\;\;\textsection\,1\;\;N$^{\circ}$\,1. Set function $\inclfremaglev_k$ is an inclusion of magmas by construction.
\end{proof}

Motivated by Definition \ref{basind}, Proposition \ref{stepind} we introduce the free magma $ \left(\fremag,\compone\right)$.
\begin{definition}\label{Cinfty}Free magma $(\fremag,\compone)$ is the direct limit of the chain
\begin{equation*}
\left(\fremag_0,\compone_0\right) \xrightarrow{\inclfremaglev_1} \left(\fremag_1,\compone_{\!1}\right)\xrightarrow{\inclfremaglev_2}...\xrightarrow{\inclfremaglev_{i-1}}\left(\fremag_{i-1},\compone_{\!i-1}\right) \xrightarrow{\inclfremaglev_i} \left(\fremag_i,\compone_{\!i}\right)\xrightarrow{\inclfremaglev_{i+1}}...
\end{equation*}
Fix $i \in \mathbb{N}_0$, we denote by $\incllim_i :\fremag_i\rightarrow \fremag$ the natural inclusion function of magmas induced by the direct limit. We set $\inclfremaglev_0= \incl{\fremag_{-1}}{\fremag_0}$, $\incllim_{-1}=\incllim_0 \funcomp \inclfremaglev_0$.\newline
The generic element belonging to $\fremag$ will be denoted by $\fmx$.\newline
We define the order of elements belonging to $\fremag$ by saying that $i\in \mathbb{N}_0\cup\left\{-1\right\}$ is the order of $\fmx\in \fremag$ if and only if $\fmx \in \incllim_i\left(\fremag_i\right)$.
\end{definition}

In Proposition \ref{linincl} below we introduce four set functions which are crucial to prove key results. We refer to Notation \ref{realvec}-[1].
\begin{proposition}\label{linincl}\mbox{}
\begin{enumerate}
\item Fix $n \in \mathbb{N}$. Then there is one and only one set function $\lboundone\quad\rboundone:\fremag^n\rightarrow \fremag$ defined by setting
\begin{equation*}
\lboundone \fmx_1,...,\fmx_n \rboundone=\incllim_{\left(\check{k}+1\right)} \left( \left(\bsfuno, \fmz_1\right),...,\left(\bsfuno, \fmz_n\right) \right)\qquad \forall \left(\fmx_1,...,\fmx_n\right)\in \fremag^n\text{,}
\end{equation*}
where: 
\begin{flalign*}
&\begin{array}{l}
\check{k}=\min\left\{k\in \mathbb{N}_0\,:\quad \forall \mathsf{n} \in \left\{1,...,n\right\}\quad\exists \fmy \in \fremag_k\quad\text{with} \quad\incllim_k\left(\fmy\right)=\fmx_{\mathsf{n}}\right\}\text{;}
\end{array}\\[8pt]
&\begin{array}{l}
\left\{
\begin{array}{l}
\left(\fmz_1,...,\fmz_n\right)\;\text{is the unique element belonging to}\;\left(\fremag_{\check{k}}\right)^n\text{and}\vspace{4pt}\\
\text{fulfilling} \;\left(\incllim_{\check{k}}\left(\fmz_1\right),..., \incllim_{\check{k}}\left(\fmz_n\right)\right)=\left(\fmx_1,...,\fmx_n\right)\text{.}
\end{array}
\right.
\end{array}
\end{flalign*}
We emphasize that in case $n=1$ we have $\lboundone\quad\rboundone=\idobj+\fremag+$.
\item  There is one and only one set function $\acmone:\bsfM \times \fremag \rightarrow \fremag$ defined by setting 
$\hspace{10pt}\bsfm \acmone \fmx = \incllim_{\left(\check{k}+1\right)}\left(\left(\bsfm, \fmz\right)\right) \qquad \forall \left(\bsfm , \fmx\right)\in\bsfM \times \fremag$,\newline
where: 
\begin{flalign*}
&\check{k}=\min\left\{k\in \mathbb{N}_0\,:\; \exists \fmy \in \fremag_k\;\text{with} \;\incllim_k\left(\fmy\right)=\fmx\right\}\text{;}\\[4pt]
&\fmz\;\text{is the unique element belonging to}\;\fremag_{\check{k}}\;\text{and fulfilling}\; \incllim_{\check{k}}\left(\fmz\right)=\fmx\text{.}
\end{flalign*}
\item There is one and only one set function $\codmagone:\fremag\rightarrow \realpartset$ fulfilling the system of conditions below
\begin{equation*} 
\left\{
\begin{array}{ll}
\codmagone\left(f\right)=\mathbb{R}\text{,}\vspace{4pt}\\
\codmagone\left(\bsfuno \acmone \fmx\right)=\codmagone\left(\fmx\right)\text{,}\vspace{4pt}\\
\codmagone\left(\Fint_i \acmone \fmx\right)=\codmagone\left(\fmx\right)\text{,}\vspace{4pt}\\
\codmagone\left(\Fpart_i \acmone \fmx\right)=\codmagone\left(\fmx\right)\text{,}\vspace{4pt}\\
\codmagone\left(\Fp_i \acmone \fmx\right)=\mathbb{R}^{0}& \text{if}\; \codmagone\left(\fmx\right)= \mathbb{R}^0\text{,}\vspace{4pt}\\
\codmagone\left(\Fp_i \acmone \fmx\right)=\mathbb{R}& \text{if}\; \codmagone\left(\fmx\right)\neq \mathbb{R}^0\text{,}\vspace{4pt}\\
\codmagone\left(\left(\Fint_i \bsfm\right)\acmone \fmx\right)=\codmagone\left(\bsfm \acmone \fmx\right)\text{,}\vspace{4pt}\\
\codmagone\left(\left(\Fpart_i \bsfm\right)\acmone \fmx\right)=\codmagone\left(\bsfm \acmone \fmx\right)\text{,}\vspace{4pt}\\
\codmagone\left(\left(\Fp_i \bsfm\right)\acmone \fmx\right)=\mathbb{R}^0& \text{if}\; \codmagone\left(\bsfm \acmone \fmx\right)=\mathbb{R}^0\text{,}\vspace{4pt}\\
\codmagone\left(\left(\Fp_i \bsfm\right)\acmone \fmx\right)=\mathbb{R}& \text{if}\; \codmagone\left(\bsfm \acmone \fmx\right)\neq\mathbb{R}^0\text{,}\vspace{4pt}\\
\codmagone\left(\lboundone \fmx_1,...,\fmx_n\rboundone\right)=\overset{n}{\underset{{\mathsf{n}}=1}{\prod}}\codmagone\left(\fmx_{\mathsf{n}}\right)\text{,}\vspace{4pt}\\
\codmagone\left(\left(\fmy\right)\,\compone\, \left(\fmz\right)\right)=\codmagone\left(\fmy\right)\text{,} 
\end{array}
\right.
\end{equation*}
$\forall i \in \mathbb{N}\text{,}\quad  \forall \bsfm \in \bsfM\text{,}\quad \forall f \in \fremag_{-1}\text{,}\quad\forall \fmx, \fmy,  \fmx_1,...,\fmx_n \in \fremag$.
\item There is one and only one set function $\dommagone:\fremag\rightarrow \realpartset$ fulfilling the system of conditions below 
\begin{equation*}
\left\{
\begin{array}{l}
\dommagone\left(f\right)=\Dom\left(f\right)\text{,} \vspace{4pt}\\
\dommagone\left(\bsfuno \acmone \fmx\right)=\dommagone\left(\fmx\right)\text{,} \vspace{4pt}\\
\dommagone\left(\Fint_i \acmone \fmx\right)=\domint\left[\dommagone\left(\fmx\right),i\right]\text{,} \vspace{4pt}\\
\dommagone\left(\Fpart_i \acmone \fmx\right)=\dommagone\left(\fmx\right)\text{,}\vspace{4pt}\\
\dommagone\left(\Fp_i \acmone \fmx\right)=\dommagone\left(\fmx\right)\text{,}\vspace{4pt}\\
\dommagone\left(\left(\Fint_i \bsfm\right)\acmone \fmx\right)=\domint\left[\dommagone\left(\bsfm \acmone \fmx\right),i\right]\text{,}\vspace{4pt}\\
\dommagone\left(\left(\Fpart_i \bsfm\right)\acmone \fmx\right)=\dommagone\left(\bsfm \acmone \fmx\right) \text{,}\vspace{4pt}\\
\dommagone\left(\left(\Fp_i \bsfm\right)\acmone \fmx\right)=\dommagone\left(\bsfm \acmone \fmx\right) \text{,}\vspace{4pt}\\
\dommagone\left(\lboundone \fmx_1,...,\fmx_n\rboundone\right)=\overset{n}{\underset{{\mathsf{n}}=1}{\prod}}\dommagone\left(\fmx_{\mathsf{n}}\right)\text{,} \vspace{4pt}\\
\dommagone\left(\left(\fmy\right)\,\compone\, \left(\fmz\right)\right)=\dommagone\left(\fmz\right)\text{,} 
\end{array}
\right.
\end{equation*} 
$\forall i \in \mathbb{N}\text{,}\quad\forall \bsfm \in \bsfM\text{,}\quad \forall f \in \fremag_{-1}\text{,}\quad\forall \fmx, \fmy,  \fmx_1,...,\fmx_n \in \fremag $.
\end{enumerate}
\end{proposition}
\begin{proof}\mbox{}\newline
\textnormal{\textbf{Proof of statement 1, 2.}}\ \ Statements 1, 2 follow by construction of $\fremag$.\newline
\textnormal{\textbf{Proof of statement 3.}}\ \ We recursively define set functions:
\begin{flalign*}
&\begin{array}{l}
\codmagone_0:\fremag_0 \rightarrow \realpartset\; \text{by setting}\vspace{4pt}\\ 
\left\{ 
\begin{array}{l}
\codmagone_0\left(f\right)=\mathbb{R}\text{,}\vspace{4pt}\\
\codmagone_0\left(\left( \fmx_1,...,\fmx_n\right)\right)=\overset{n}{\underset{{\mathsf{n}}=1}{\prod}}\Cod\left(\fmx_{\mathsf{n}}\right)\text{,}\vspace{4pt}\\
\codmagone_0\left(\left(\fmy\right)\,\compone_{\!0}\, \left(\fmz\right)\right)=\codmagone_0\left(\fmy\right) \text{,} 
\end{array}\right.\vspace{4pt}\\
\forall f \in \fremag_{-1}\text{,}\quad \forall \left( \fmx_1,...,\fmx_n\right) \in \left(\fremag_{-1}\right)^n\text{,}\quad \forall \fmy,\fmz \in \fremag_0\text{;} 
\end{array}\\[8pt]  
&\begin{array}{l}
\codmagone_k:\fremag_{k}\rightarrow \realpartset \;\text{by setting} \vspace{4pt}\\
\left\{
\begin{array}{ll}
\codmagone_k\left(\left(\bsfuno, \fmx\right)\right)=\codmagone_{k-1}(\fmx)\text{,}\vspace{4pt}\\
\codmagone_k\left(\left(\Fint_i, \fmx\right)\right)=\codmagone_{k-1}(\fmx)\text{,}\vspace{4pt}\\
\codmagone_k\left(\left(\Fpart_i, \fmx\right)\right)=\codmagone_{k-1}(\fmx)\text{,}\vspace{4pt}\\
\codmagone_k\left(\left(\Fp_i, \fmx\right)\right)={0}& \text{if}\; \codmagone_{k-1}\left(\fmx\right)= \mathbb{R}^0\text{,}\vspace{4pt}\\
\codmagone_k\left(\left(\Fp_i, \fmx\right)\right)=\mathbb{R}& \text{if}\; \codmagone_{k-1}\left(\fmx\right)\neq \mathbb{R}^0\text{,}\vspace{4pt}\\
\codmagone_k\left(\left(\Fint_i \bsfm, \fmx\right)\right)=\codmagone_{k}\left(\left(\bsfm, \fmx\right)\right)\text{,}\vspace{4pt}\\
\codmagone_k\left(\left(\Fpart_i \bsfm, \fmx\right)\right)=\codmagone_{k}\left(\left(\bsfm, \fmx\right)\right)\text{,}\vspace{4pt}\\
\codmagone_k\left(\left(\Fp_i \bsfm, \fmx\right)\right)=\mathbb{R}^0& \text{if}\; \codmagone_{k}\left(\left(\bsfm, \fmx\right)\right)=\mathbb{R}^0\text{,}\vspace{4pt}\\
\codmagone_k\left(\left(\Fp_i \bsfm, \fmx\right)\right)=\mathbb{R}& \text{if}\; \codmagone_{k}\left(\left(\bsfm, \fmx\right)\right)\neq\mathbb{R}^0\text{,}\vspace{4pt}\\
\codmagone_k\left(\left( \fmx_1,...,\fmx_n\right)\right)=\overset{n}{\underset{{\mathsf{n}}=1}{\prod}}\codmagone_{k}\left(\fmx_{\mathsf{n}}\right)\text{,}\vspace{4pt}\\
\codmagone_k\left(\left(\fmy\right)\,\compone_{\!k}\, \left(\fmz\right)\right)=\codmagone_k\left(\fmy\right)\text{,} 
\end{array}
\right.\vspace{4pt}\\
\forall i \in \mathbb{N}\text{,}\quad\forall \bsfm \in \bsfM\text{,}\quad     \forall \fmx \in \fremag_{k-1}\text{,}\quad \forall \left(\fmx_1,...,\fmx_n\right) \in \fremaggen_k\text{,}\quad \forall \fmy,\fmz \in \fremag_{k}\text{.}
\end{array}
\end{flalign*}
Set functions $\codmagone_k$ are well defined for any $k \in \mathbb{N}$ since: Remark \ref{intdiffmonrem}-[1] holds true; $\codmagone_k\left(\bsfm \left(\fmx\right)\right)=\codmagone_k\left(\bsfn \left(\fmx\right)\right)$ for any $\fmx\in\fremag_{k-1}$ and any relation $\bsfm \bsfrelsymb\bsfn$ chosen among \eqref{intint}-\eqref{leftrightinter}.\newline
Eventually statement follows since $\qquad \codmagone_{k-1}=\codmagone_k \funcomp \inclfremaglev_k \qquad \forall k \in \mathbb{N}_0$.\newline
\textnormal{\textbf{Proof of statement 4.}}\ \ Referring to Notation \ref{realvec}-[6], Remark \ref{rprop} we recursively define set functions:
\begin{flalign*}
&\begin{array}{l}
\dommagone_0:\fremag_0 \rightarrow \realpartset\,\text{by setting}\vspace{4pt}\\
\left\{\begin{array}{l} 
\dommagone_0\left(f\right)=\Dom\left(f\right)\text{,}\vspace{4pt}\\
\dommagone_0\left(\left( \fmx_1,...,\fmx_n\right)\right)=\overset{n}{\underset{{\mathsf{n}}=1}{\prod}}\Dom\left(\fmx_{\mathsf{n}}\right)\text{,}\vspace{4pt}\\
\dommagone_0\left(\left(\fmy\right)\,\compone_{\!0}\, \left(\fmz\right)\right)=\dommagone_0\left(\fmz\right) \text{,} 
\end{array}\right.\vspace{4pt}\\
\forall f \in \fremag_{-1}\text{,}\quad \forall \left( \fmx_1,...,\fmx_n\right) \in \left(\fremag_{-1}\right)^n\text{,}\quad \forall \fmy,\fmz \in \fremag_0\text{;}
\end{array}\\[8pt]
&\begin{array}{l}
\dommagone_k:\fremag_{k}\rightarrow \realpartset\;\text{by setting}\vspace{4pt}\\ 
\left\{
\begin{array}{l}
\dommagone_k\left(\left(\bsfuno, \fmx\right)\right)=\dommagone_{k-1}\left(\fmx\right)\text{,} \vspace{4pt}\\
\dommagone_k\left(\left(\Fint_i, \fmx\right)\right)=\domint\left[\dommagone_{k-1}\left(\fmx\right),i\right]\text{,} \vspace{4pt}\\
\dommagone_k\left(\left(\Fpart_i, \fmx\right)\right)=\dommagone_{k-1}\left(\fmx\right)\text{,}\vspace{4pt}\\
\dommagone_k\left(\left(\Fp_i, \fmx\right)\right)=\dommagone_{k-1}\left(\fmx\right)\text{,}\vspace{4pt}\\
\dommagone_k\left(\left(\Fint_i \bsfm, \fmx\right)\right)=\domint\left[\dommagone_{k}\left(\left(\bsfm, \fmx\right)\right),i\right]\text{,}\vspace{4pt}\\
\dommagone_k\left(\left(\Fpart_i \bsfm, \fmx\right)\right)=\dommagone_{k}\left(\left(\bsfm, \fmx\right)\right) \text{,}\vspace{4pt}\\
\dommagone_k\left(\left(\Fp_i \bsfm, \fmx\right)\right)=\dommagone_k\left(\left(\bsfm, \fmx\right)\right) \text{,}\vspace{4pt}\\
\dommagone_k\left(\left( \fmx_1,...,\fmx_n\right)\right)=\overset{n}{\underset{{\mathsf{n}}=1}{\prod}}\dommagone_{k}\left(\fmx_{\mathsf{n}}\right)\text{,} \vspace{4pt}\\
\dommagone_k\left(\left(\fmy\right)\,\compone_{\!k}\, \left(\fmz\right)\right)=\dommagone_k\left(\fmz\right)\text{,} 
\end{array}
\right.\vspace{4pt}\\
\forall i \in \mathbb{N}\text{,}\quad \forall \bsfm \in \bsfM\text{,}\quad \forall \fmx \in \fremag_{k-1}\text{,}\quad \forall \left(\fmx_1,...,\fmx_n\right) \in \fremaggen_k\text{,}\quad \forall \fmy,\fmz \in \fremag_{k}\text{.}
\end{array}
\end{flalign*} 
Set functions $\dommagone_k$ are well defined for any $k \in \mathbb{N}$ since: Remark \ref{intdiffmonrem}-[1] holds true; $\dommagone_k\left(\left(\bsfm, \fmx\right)\right)=\dommagone_k\left(\left(\bsfn, \fmx\right)\right)$ for any $\fmx\in\fremag_{k-1}$ and any relation $\bsfm \bsfrelsymb\bsfn$ chosen among \eqref{intint}-\eqref{leftrightinter}.\newline
Eventually statement follows since $\qquad \dommagone_{k-1}=\dommagone_k \funcomp \inclfremaglev_k \qquad \forall k \in \mathbb{N}_0$.
\end{proof}

\begin{remark}\label{domcae}
We emphasize that for $\fmx\in \fremag$ the notion of $\Cae\left(\fmx\right)$ cannot be well defined. Precisely it is straightforward to verify that the set $\Cae\left(\fmx\right)$ changes in $\dommagone\left(\fmx\right)$ when evaluating $\fmx$ by replacing elements belonging to $\bascont$ with smooth functions. Similar consideration holds for $\Ima\left(\fmx\right)$ with respect to $\codmagone\left(\fmx\right)$. 
\end{remark}

Motivated by Proposition \ref{linincl}-[3, 4] we introduce set functions giving the dimension of domain and co-domain of elements belonging $\fremag$. 
\begin{definition}\label{dimdomcodfuncontone}
We define set functions
\begin{flalign*}
& \dimcodmagone:\fremag\rightarrow\mathbb{N}_0\quad\text{by setting}  \quad \codmagone\left(\fmx\right)=\mathbb{R}^{\dimcodmagone\left(\fmx \right)}\quad\forall \fmx \in \fremag \text{;}\\[8pt]
&  \dimdommagone:\fremag\rightarrow\mathbb{N}_0\quad\text{by setting}  \quad \dommagone\left(\fmx\right)\subseteqdentro\mathbb{R}^{\dimdommagone\left( \fmx\right)}\quad\forall\fmx \in \fremag\text{.}
\end{flalign*}
\end{definition}

In order to prove key results below we need to handle elements belonging to $\fremag_k$ very accurately. To this aim we consider any fixed $\fmx \in \fremag_k$ as a finite ordered sequence of characters, then we introduce notation and tools which are needed to handle finite ordered sequences of characters with the required accuracy. We refer to Notation \ref{ins}-[2, 5].

\begin{definition}\label{finordseq}\mbox{} 
\begin{enumerate}
\item Fix $k \in \mathbb{N}_0$. We define the set $\symb_k$ by setting:
\begin{equation*}
\begin{array}{c}
\fremag_{-1}\subseteq \symb_k\text{;}\qquad \bsfM \subseteq \symb_k\text{;}\qquad \compone_{\!{\mathsf{k}}} \in \symb_k\qquad \forall {\mathsf{k}} \in\left\{0,...,k\right\}\text{;}\vspace{4pt}\\
( \in \symb_k\text{;}\qquad  ) \in \symb_k\text{;}\qquad  , \in \symb_k\text{.}
\end{array}
\end{equation*}
\item Fix $k \in \mathbb{N}_0$. We set $\extfremag_k=\left\{\extx:\left\{1,...,n\right\}\rightarrow \symb_k\qquad \forall n \in \mathbb{N}\right\} $.
\item Fix $k \in \mathbb{N}_0$. We define the set function $\extlnght_k:\extfremag_k\rightarrow \mathbb{N}$ by setting \begin{equation*}
\extlnght_k\left(\extx\right)=\cardset\left(\Dom\left(\extx\right)\right)\qquad  \forall \extx\in \extfremag_k\text{.}
\end{equation*}
We say that the positive integer $\extlnght_k\left(\extx\right)$ is the length of $\extx$. 
\item Fix $k,l \in \mathbb{N}_0$, $i \in \mathbb{N}$, $\extx \in \extfremag_k$. Assume $i, i+l\in \left\{1,...,\extlnght_k\left(\extx\right)\right\}$.\newline
We define the set function $\extsubstrchar_k\left[\extx,i,l\right]:\left\{1,...,l\right\}\rightarrow \symb_k$ by setting
\begin{equation*}
\extsubstrchar_k\left[\extx,i,l\right]\left({\mathsf{l}}\right)=\extx\left(i+{\mathsf{l}}-1\right) \qquad \forall {\mathsf{l}} \in \left\{1,...,l\right\}\text{.}
\end{equation*} 
\item Fix $k \in \mathbb{N}_0$, $\extx\in \extfremag_{k}$.\newline
Fix $\left(\extxtwo,i \right)\in \extfremag_k\times\mathbb{N}$. We say that  $\left(\extxtwo,i \right)$ is an occurrence pair for $\extx$ if and only if $\extsubstrchar_k\left[\extx,i, \extlnght_k\left(\extxtwo\right)\right]=\extxtwo$. We also say that: $\left(\extxtwo,i \right)$ occurs in $\extx$; $\extxtwo$ occurs in $\extx$ at $i$; $\extxtwo$ occurs in $\extx$, by dropping any reference to $i$ whenever there is no need to specify where $\extxtwo$ occurs in $\extx$; $i$ is an occurrence of $\extxtwo$ in $\extx$.\newline
Fix a set $\extoccset \subseteq \extfremag_k\times\mathbb{N}$. We say that $\extoccset$ is an occurrence set for $\extx$ if and only if either $\extoccset = \udenset$ or the four conditions below are all fulfilled
\begin{equation*}
\left\{
\begin{array}{ll}
(i)&\extoccset\;\text{is a not empty, finite, totally ordered set,}\vspace{6pt}\\
(ii)&i+\extlnght_k\left(\extxtwo\right)\leq \extlnght_k\left(\extx\right)  \hspace{20pt}\forall \left(\extxtwo,i\right) \in \extoccset\text{,} \vspace{6pt}\\
(iii)&i_1+\extlnght_k\left(\extxtwo_1\right)\leq i_2\vspace{4pt}\\ 
&\hspace{15pt}\forall \left(\extxtwo_1,i_1\right), \left(\extxtwo_2,i_2\right) \in \extoccset \;\text{with}\;  \left(\extxtwo_1,i_1\right) \prec\left(\extxtwo_2,i_2\right)\text{,}\vspace{6pt}\\
(iv)&\extsubstrchar_k\left[\extx,i,\extlnght_k\left(\extxtwo\right) \right]=\extxtwo \hspace{20pt} \forall \left(\extxtwo,i\right) \in \extoccset\text{.}
\end{array}
\right.
\end{equation*}
\item Fix $k \in \mathbb{N}_0$, $\extx\in \extfremag_{k}$, an occurrence set $\left\{\left(\extx_{\mathsf{m}},i_{\mathsf{m}}\right)\right\}_{{\mathsf{m}}=1}^m$ for $\extx$, a set function $\extoccfun:\left\{\left(\extx_{\mathsf{m}},i_{\mathsf{m}}\right)\right\}_{{\mathsf{m}}=1}^m\rightarrow \extfremag_k$. Assume $\extoccfun \neq \emptysetfun$. Set:
\begin{flalign*}
&\begin{array}{l}
\deridxuno\left[\extx, \extoccfun\right]=\extlnght_k\left(\extx\right)-\underset{{\mathsf{m}}=1}{\overset{m}{\sum}} \extlnght_k\left(\extx_{\mathsf{m}}\right)+\underset{{\mathsf{m}}=1}{\overset{m}{\sum}} \extlnght_k\left(\extoccfun\left(\extx_{\mathsf{m}},i_{\mathsf{m}}     \right)\right)\text{;}
\end{array}\\[8pt]
&\begin{array}{l}
\deridxdue_{\mathsf{m}}=\left\{
\begin{array}{ll}
i_1-1 &\text{if}\; {\mathsf{m}}=1\text{,}\vspace{4pt}\\
i_{\mathsf{m}} -\underset{a=1}{\overset{{\mathsf{m}}-1}{\sum}} \extlnght_k\left(\extx_a\right)+\underset{a=1}{\overset{{\mathsf{m}}-1}{\sum}} \extlnght_k\left(\extoccfun\left(\extx_a,i_a\right)\right)-1&\text{if}\; 2\leq {\mathsf{m}} \leq m\text{.}
\end{array}
\right.
\end{array}
\end{flalign*}
We define the set function 
$\extstrchar_k\left[\extx, \extoccfun\right]:\left\{1,...,\deridxuno\left[\extx, \extoccfun\right]\right\}\rightarrow \symb_k$ by setting:
\begin{flalign*}
&\left\{
\begin{array}{l}
\extstrchar_k\left[\extx, \extoccfun\right]\left(l \right)=\extx\left(l\right)\hspace{88pt} 
\text{if}\; 1\leq l < i_1\text{,}\vspace{10pt}\\
\extstrchar_k\left[\extx, \extoccfun\right]\left(\deridxdue_{\mathsf{m}}+l \right)=\left(\extoccfun\left(\extx_{\mathsf{m}},i_{\mathsf{m}}\right)\right)\left(l\right)\hspace{20pt} 
\text{if}\; 1\leq l \leq \extlnght_k\left(\extoccfun\left(\extx_{\mathsf{m}},i_{\mathsf{m}}\right)\right)\text{,}\vspace{10pt}\\ 
\extstrchar_k\left[\extx, \extoccfun\right]\left(\deridxdue_{\mathsf{m}}+l \right)=
\extx\left(i_{\mathsf{m}}+l+\extlnght_k\left(\extx_{\mathsf{m}}\right)-\extlnght_k\left(\extoccfun\left(\extx_{\mathsf{m}},i_{\mathsf{m}}\right)\right)-1\right)\\
\text{if}\hspace{20pt} \extlnght_k\left(\extoccfun\left(\extx_{\mathsf{m}},i_{\mathsf{m}}\right)\right)+1\leq l\leq \extlnght_k\left(\extoccfun\left(\extx_{\mathsf{m}},i_{\mathsf{m}}\right)\right)+i_{{\mathsf{m}}+1}-i_{\mathsf{m}}-\extlnght_k\left(\extx_{\mathsf{m}}\right)\\
\text{and} \hspace{9pt}1\leq {\mathsf{m}}\leq m-1\text{,}\vspace{10pt}\\ 
\extstrchar_k\left[\extx, \extoccfun\right]\left(\deridxdue_m+l \right)=
\extx\left(i_m+l+\extlnght_k\left(\extx_m\right)-\extlnght_k\left(\extoccfun\left(\extx_m,i_m\right)\right)-1\right)\\
\text{if}\hspace{20pt}\extlnght_k\left(\extoccfun\left(\extx_m,i_m\right)\right)+1\leq l\\
\text{and}\hspace{9pt}l\leq \extlnght_k\left(\extoccfun\left(\extx_m,i_m\right)\right)+ \extlnght_k\left(\extx\right)  +1 -i_m-\extlnght_k\left(\extx_m\right)\text{.}
\end{array}
\right.
\end{flalign*}
\end{enumerate}
\end{definition}

\begin{remark}\label{remfinordseq}We emphasize that:
\begin{enumerate}
\item $\extsubstrchar_k\left[\extx,i,l\right]$ is the finite ordered sequence obtained from $\extx$ by taking all symbols from the $i$-th to the $(i+l-1)$-th. 
\item $\extstrchar_k\left[\extx, \extoccfun\right]$ is the finite ordered sequence of symbols obtained by replacing any subsequence $\extsubstrchar_k\left[\extx,i_{\mathsf{m}},\extlnght_k\left(\extx_{\mathsf{m}}\right) \right]$ of $\extx$ by the sequence $\extoccfun\left(\extx_{\mathsf{m}},i_{\mathsf{m}}\right)$.
\item Fix $k\in \left\{0,...,l\right\}$. There is an injective set function $\theta\left[k,l\right]:\fremag_k\rightarrow \extfremag_l$ which associates to any fixed $\fmx\in \fremag_k$ the corresponding finite ordered sequence of characters belonging to $\extfremag_l$. With an abuse of language we denote the set $\theta\left[k,l\right]\left(\fremag_k\right)\subseteq \extfremag_l$ again by $\fremag_k$ and we denote $\theta\left[k,l\right]\left(\fmx\right)$ again by $\fmx$ for any $\fmx \in \fremag_k$.
\item With an abuse of language, by referring to Notation \ref{ins}-[7], we denote $\cost<\left\{1\right\}<>\symb_k>+\bsfm+\in \extfremag_k$ again by $\bsfm$ for any $\bsfm \in \bsfM$.
\end{enumerate}
\end{remark}

\begin{proposition}\label{invsost}
Fix $k \in \mathbb{N}_0$, $\extx\in \extfremag_{k}$, an occurrence set $\left\{\left(\extx_{\mathsf{m}},i_{\mathsf{m}}\right)\right\}_{{\mathsf{m}}=1}^m$ for $\extx$, a set function $\extoccfun:\left\{\left(\extx_{\mathsf{m}},i_{\mathsf{m}}\right)\right\}_{{\mathsf{m}}=1}^m\rightarrow \extfremag_k$. Assume $\extoccfun \neq \emptysetfun$.\newline
Then there are: 
\begin{flalign*}
&\text{one and only one occurrence set}\;\;\left\{\left(\extxtwo_{\mathsf{m}},j_{\mathsf{m}}\right)\right\}_{{\mathsf{m}}=1}^m\;\;\text{for}\;\;\extstrchar_k\left[\extx, \extoccfun\right]\text{,} \\[8pt]
&\text{one and only one set function}\;\; \extoccfundue:\left\{\left(\extxtwo_{\mathsf{m}},j_{\mathsf{m}}\right)\right\}_{{\mathsf{m}}=1}^m\rightarrow \extfremag_k  \text{,}
\end{flalign*}
such that $\extx=\extstrchar_k\left[\extstrchar_k\left[\extx, \extoccfun\right], \extoccfundue\right]$. 
\end{proposition} 
\begin{proof}
We set: 
\begin{flalign*}
&\begin{array}{l}
\extxtwo_{\mathsf{m}}=\extoccfun\left(\extx_{\mathsf{m}}, i_{\mathsf{m}} \right)\hspace{149pt} \forall {\mathsf{m}} \in\left\{1,...,m\right\}\text{;}\end{array}\\[12pt]
&\begin{array}{l}
j_{\mathsf{m}}=\left\{
\begin{array}{ll}
i_1 &\text{if}\; {\mathsf{m}}=1\text{,}\vspace{4pt}\\
i_{\mathsf{m}} -\underset{a=1}{\overset{{\mathsf{m}}-1}{\sum}} \extlnght_k\left(\extx_a\right)+\underset{a=1}{\overset{{\mathsf{m}}-1}{\sum}} \extlnght_k\left(\extoccfun\left(\extx_a,i_a\right)\right)&\text{if}\; {\mathsf{m}} \in \left\{2,...,m \right\}\text{;}
\end{array}
\right.
\end{array}\\[12pt]
&\begin{array}{l}
\extoccfundue\left(\extxtwo_{\mathsf{m}}, j_{\mathsf{m}} \right)=\extx_{\mathsf{m}}\hspace{149pt} \forall {\mathsf{m}} \in\left\{1,...,m\right\}\text{.}
\end{array}
\end{flalign*}
Eventually statement follows by computing $\extstrchar_k\left[\extstrchar_k\left[\extx, \extoccfun\right], \extoccfundue\right]$.  
\end{proof}

\begin{proposition}\label{propfinordseq} Fix $l \in \mathbb{N}_0$, $k \in \left\{0,...,l\right\}$, $\fmx\in \fremag_{k}$, an occurrence set $\left\{\left(\fmx_{\mathsf{m}},i_{\mathsf{m}}\right)\right\}_{{\mathsf{m}}=1}^m$ for $\fmx$, a set function $\extoccfun:\left\{\left(\fmx_{\mathsf{m}},i_{\mathsf{m}}\right)\right\}_{{\mathsf{m}}=1}^m\rightarrow \extfremag_l$.
Assume 
\begin{equation*}
\forall {\mathsf{m}} \in \left\{1,...,m\right\}\quad
\exists k_{\mathsf{m}} \in \left\{0,...,k\right\}\;\text{such that}\;\fmx_{\mathsf{m}} \in\fremag_{k_{\mathsf{m}}}\; \text{and}\; \extoccfun\left(\fmx_{\mathsf{m}},i_{\mathsf{m}}\right)\in \fremag_{k_{\mathsf{m}}}\text{.}
\end{equation*}
Then:
\begin{enumerate}
\item $\extstrchar_k\left[\fmx, \extoccfun\right] \in \fremag_k$;
\item there are
\begin{flalign*}
&\text{one and only one occurrence set}\;\;\left\{\left(\fmy_{\mathsf{m}},j_{\mathsf{m}}\right)\right\}_{{\mathsf{m}}=1}^m\;\;\text{for}\;\; \inclfremaglev_{k+1}\left(\fmx\right)\text{,}\\[8pt]
&\text{one and only one set function}\;\; \extoccfundue: \left\{\left( \fmy_{\mathsf{m}},j_{\mathsf{m}}\right)\right\}_{{\mathsf{m}}=1}^m\rightarrow \extfremag_{l+1}\text{,}
\end{flalign*}
such that: 
\begin{flalign*}
&\begin{array}{l}
\fmy_{\mathsf{m}}=\left\{
\begin{array}{ll}
\inclfremaglev_{k+1}\left(\fmx_a\right)&\forall {\mathsf{m}} \in\left\{1,...,m\right\}\;\text{with} \; k_{\mathsf{m}} = k\text{,}\vspace{4pt}\\
\fmx_{\mathsf{m}} & \forall {\mathsf{m}} \in\left\{1,...,m\right\}\;\text{with}\;k_{\mathsf{m}}\neq k\text{;}
\end{array}
\right.
\end{array}\\[8pt]
&\begin{array}{l}
\extoccfundue\left(\fmy_{\mathsf{m}},j_{\mathsf{m}}\right)=
\left\{
\begin{array}{ll}
\inclfremaglev_{k+1}  \left(\extoccfun\left(  \fmx_{\mathsf{m}},i_{\mathsf{m}}\right)\right)&\forall {\mathsf{m}} \in\left\{1,...,m\right\}\;\text{with}\; k_{\mathsf{m}} = k\text{,}\vspace{4pt}\\
\extoccfun\left(  \fmx_{\mathsf{m}},i_{\mathsf{m}}\right)& \forall {\mathsf{m}} \in\left\{1,...,m\right\}\;\text{with}\;k_{\mathsf{m}}\neq k\text{;}
\end{array}
\right.
\end{array}\\[8pt]
&\begin{array}{l}
\inclfremaglev_{k+1}\left(\extstrchar_{l}\left[\fmx, \extoccfun\right]\right)=
\left\{
\begin{array}{ll}
\extstrchar_{l+1}\left[\inclfremaglev_{k+1}  \left(\fmx\right), \extoccfundue\right]&\text{if} \; k = l\text{,}\vspace{4pt}\\
\extstrchar_{l}\left[\inclfremaglev_{k+1}  \left(\fmx\right), \extoccfundue\right]& \text{if}\;k\neq l\text{.}
\end{array}
\right.
\end{array}
\end{flalign*}
\end{enumerate}
\end{proposition}
\begin{proof}
Statements follows by exploiting the following two facts occurring for any $i \in \mathbb{N}_0$: 
$\inclfremaglev_{i+1}  \left(\fremaggen_i\right)\subseteq \fremaggen_{i+1}$; any $z \in \fremag_i$ can be uniquely written in terms of elements belonging to $\fremaggen_i$.      
\end{proof}

\begin{proposition}\label{posrel}
Fix $k \in \mathbb{N}_0$, $\fmx\in \fremag_{k}$, two occurrence pairs $\left(\fmx_1,i_1\right), \left(\fmx_2,i_2\right)$ for $\fmx$. Assume $\quad\forall a \in \left\{1,2\right\}\quad
\exists k_a \in \left\{0,...,k\right\}\;\;\text{such that}\;\;\fmx_a \in\fremag_{k_a}$.\newline
Then one among the following mutually exclusive four cases occurs:
\begin{equation*}
\begin{array}{ll}
i_1+\extlnght_k\left(\fmx_1\right)\leq i_2\text{,} \hspace{20pt}&i_2\leq i_1< i_1+\extlnght_k\left(\fmx_1\right)\leq i_2+\extlnght_k\left(\fmx_2\right)\text{,}\vspace{4pt}\\
i_2+\extlnght_k\left(\fmx_2\right)\leq i_1\text{,}\hspace{20pt}& i_1\leq i_2< i_2+\extlnght_k\left(\fmx_2\right)\leq i_1+\extlnght_k\left(\fmx_1\right)\text{.}
\end{array}
\end{equation*}
\end{proposition}
\begin{proof}
Statement follows by exploiting the following two facts occurring for any $i \in \mathbb{N}_0$: 
$\inclfremaglev_{i+1}  \left(\fremaggen_i\right)\subseteq \fremaggen_{i+1}$; $\fmz$ can be uniquely written in terms of elements belonging to $\fremaggen_i$ for any $\fmz \in \fremag_i$.    
\end{proof}

Motivated by Definition \ref{finordseq}, Proposition \ref{propfinordseq}, Remark \ref{remfinordseq}-[3, 4] we extend to $\fremag$ the approach by sequence of characters developed above for elements belonging to $\fremag_{k}$.

\begin{definition}\mbox{}\label{defpath}
\begin{enumerate}
\item Fix $\fmx\in \fremag$.\newline
Fix $\left(\fmy,i \right)\in \fremag\times\mathbb{N}$. We say that $\left(\fmy,i \right)$ is an occurrence pair for $\fmx$ if and only if exist $k_{\fmx},k_{\fmy}\in \mathbb{N}_0$, $\fmz_{\fmx} \in \fremag_{k_{\fmx}}$, $\fmz_{\fmy} \in \fremag_{k_{\fmy}}$ fulfilling system below:
\begin{equation}
\left\{
\begin{array}{l}
\incllim_{k_{\fmx}} \left(\fmz_{\fmx}\right)=\fmx\text{,}\vspace{4pt}\\
\incllim_{k_{\fmy}} \left(\fmz_{\fmy}\right)=\fmy\text{,}\vspace{4pt}\\
k_{\fmy}\in \left\{0,...,k_{\fmx}\right\}\text{,}\vspace{4pt}\\
\left(\fmz_{\fmy},i \right)\;\text{is an occurrence pair for} \;\fmz_{\fmx}\text{.}
\end{array}\right.\label{sysocc}
\end{equation}
We also say that: $\left(\fmy,i \right)$ occurs in $\fmx$; $\fmy$ occurs in $\fmx$ at $i$; $\fmy$ occurs in $\fmx$ by dropping any reference to $i$ whenever there is no need to specify where $\fmy$ occurs in $\fmx$; $i$ is an occurrence of $\fmy$ in $\fmx$.\newline 
Fix a set $\occset \subseteq \fremag\times\mathbb{N}$, a set function $\occfun:\occset\rightarrow \fremag$.\newline
We say that $\occset$ is an occurrence set for $\fmx$ if and only if there are $k_{\fmx}\in \mathbb{N}_0$, $\fmz_{\fmx} \in \fremag_{k_{\fmx}}$, an occurrence set $\extoccset$ for $\fmz_{\fmx}$ fulfilling system below
\begin{equation}
\left\{
\begin{array}{l}
\incllim_{k_{\fmx}} \left(\fmz_{\fmx}\right)=\fmx\text{,}\vspace{4pt}\\
\forall\left(\fmy,i\right)\in \extoccset \quad
\exists k\left[\fmy,i\right] \in \left\{0,...,k_{\fmx}\right\}\;\;\text{such that}\;\;\fmy \in\fremag_{k\left[\fmy,i\right]}\text{,}\vspace{4pt}\\
\occset=\left\{\left(\incllim_{k\left[\fmy,i\right]} \left(\fmy\right),i\right),\quad  \left(\fmy,i\right)\in \extoccset\right\}\text{.}
\end{array}
\right.\label{hypextelcomp}
\end{equation}
We say that $\occfun$ is an occurrence function on $\occset$ if and only if there is a set function $\extoccfun:\extoccset\rightarrow \extfremag_{k_{\fmx}}$ fulfilling the system of conditions below
\begin{equation}
\left\{
\begin{array}{ll}
\extoccfun\left(\fmy,i\right) \in \fremag_{k\left[\fmy,i\right]}& \forall \left(\fmy,i\right)\in \extoccset\text{,}\\[4pt]
\occfun\left(\incllim_{k\left[\fmy,i\right]} \left(\fmy\right),i\right)=\incllim_{k\left[\fmy,i\right]}\left(\extoccfun\left(\fmy,i\right)\right)&\forall \left(\fmy,i\right)\in \extoccset
\text{.}
\end{array}
\right.
\label{hypasselcomp}
\end{equation}
\item Fix $\bsfm \in \bsfM$, $\fmx\in \fremag$. We say that $\bsfm$ occurs in $\fmx$ if and only if exist $i, k\in \mathbb{N}$, $\fmz \in \fremag_k$ such that:
\begin{equation*}
\incllim_k \left(\fmz\right)=x\text{;}\hspace{20pt}\left(\bsfm, i \right)\,\text{is an occurrence pair for} \;\fmz\text{.}
\end{equation*}
\item Fix $n \in \mathbb{N}$, $\left(\fmx_1,...,\fmx_n\right) \in \fremag^n$, an occurrence set $\occset_{\mathsf{n}}$ for $\fmx_{\mathsf{n}}$, an occurrence function $\occfun_{\mathsf{n}}$ on $\occset_{\mathsf{n}}$ for any ${\mathsf{n}} \in \left\{1,...,n\right\}$.
Define
\begin{equation*}
\begin{array}{l}
N=\left\{{\mathsf{n}} 	\in 	\left\{1,...,n\right\}\,:\; \occfun_{\mathsf{n}}\neq \emptysetfun \right\}\text{,}\vspace{4pt}\\
\occfun: \underset{{\mathsf{n}}\in N}{\coprod}\occset_{{\mathsf{n}}}\rightarrow \fremag\quad\text{by setting}\quad	
\occfun\funcomp \incl{\occset_{\mathsf{m}}}{\underset{{\mathsf{n}}\in N}{\coprod}\occset_{{\mathsf{n}}}}=\occfun_{\mathsf{n}}\quad \forall  {\mathsf{m}}\in N\text{.}
\end{array}
\end{equation*}
We emphasize that if $N=\udenset$ then $\occfun=\emptysetfun$.\newline
We deal with a finite number of elements belonging to $\fremag$ (i.e. $\fmx_1$,...,$\fmx_n$), of occurrence sets (i.e. $\occset_1$,...,$\occset_n$), of occurrence functions (i.e. $\occfun_1$,...,$\occfun_n$), then by Proposition \ref{propfinordseq}-[2] there is no loss of generality by assuming that there are
$k \in \mathbb{N}_0$ , $\fmz_{\mathsf{n}}\in \fremag_{k}$ for any ${\mathsf{n}}\in \left\{1,...,n\right\}$, an occurrence set $\extoccset_{\mathsf{n}}$ for $\fmz_{\mathsf{n}}$ for any ${\mathsf{n}}\in \left\{1,...,n\right\}$, 
a set function $\extoccfun_{\mathsf{n}}:\extoccset_{\mathsf{n}}\rightarrow \extfremag_{k}$ for any ${\mathsf{n}}\in \left\{1,...,n\right\}$, such that: data $\fmx_{\mathsf{n}}$, $\occset_{\mathsf{n}}$, $k$, $\fmz_{\mathsf{n}}$, $\extoccset_{\mathsf{n}}$ fulfill \eqref{hypextelcomp} for any ${\mathsf{n}}\in \left\{1,...,n\right\}$; data $\fmx_{\mathsf{n}}$, $\occset_{\mathsf{n}}$, $\occfun_{\mathsf{n}}$, $k$, $\fmz_{\mathsf{n}}$, $\extoccset_{\mathsf{n}}$, $\extoccfun_{\mathsf{n}}$ fulfill \eqref{hypasselcomp} for any ${\mathsf{n}}\in \left\{1,...,n\right\}$.\newline
Hence we define $\strchar \left[\left(\fmx_1,...,\fmx_n\right), \occfun\right]=\left(\fmy_1,...,\fmy_n\right)\in \fremag^n$ where
\begin{equation*}
\fmy_{\mathsf{n}}=\left\{
\begin{array}{ll}
\incllim_k \left(\extstrchar_k\left[\fmz_{\mathsf{n}}, \extoccfun_{\mathsf{n}}\right]\right)& \text{if}\;{\mathsf{n}} \in N\text{,}\vspace{4pt}\\
\fmx_{\mathsf{n}}& \text{if}\;{\mathsf{n}} \notin N\text{.}
\end{array}
\right.
\end{equation*}
\item Fix $\fmx \in \fremag$. We set 
\begin{flalign*}
&\begin{array}{l}
\occset^{\fmx}=\big\{\left(\incllim_{-1}\left(f\right),i\right)\qquad\forall \left(f,i\right)\in \bascont\times \mathbb{N}\;\text{such that}\\[4pt]
 \hspace{100pt}\left(\incllim_{-1}\left(f\right),i\right) \;\text{is an occurrence pair for}\;\fmx\big\}\text{;}
\end{array}\\[8pt]
& \begin{array}{l}
\occset_{\fmx} =\left\{f \in \bascont : \;\;\incllim_{-1}\left(f\right) \;\text{occurs in}\; \fmx \right\}\text{.}
\end{array}
\end{flalign*}
\end{enumerate}
\end{definition}

\begin{proposition}\label{suboccass}
Fix $\fmx \in\fremag$, an occurrence set $\occset$ for $\fmx$, a subset $\occsetdue\subseteq \occset$, an occurrence function $\occfun$ on $\occset$. Then: $\occsetdue$ is an occurrence set for $\fmx$; $\occfun\funcomp\incl{\occsetdue}{\occset}$ is an occurrence function on $\occsetdue$.
\end{proposition}
\begin{proof}
Statement follows directly from Definition \ref{defpath}-[1].
\end{proof}

In Construction \ref{superevalsm} below we introduce the notion of deduced occurrence set and deduced occurrence function.
\begin{construction}\label{superevalsm}
Fix $\fmx\in \fremag$, an occurrence pair $\left(\fmy,i \right)$ for $\fmx$, an occurrence set $\occset$ for $\fmy$, an occurrence function $\occfun$ on $\occset$.\newline
We denote by 
\begin{flalign*}
&\occsetstr\left[\occset ,\left(\fmy,i\right)\right] \;\text{the occurrence set for}\;\fmy\text{,}\\[8pt]
&\occfunstr\left[\occfun, \left(\fmy,i\right)\right]\;\text{the occurrence function on}\; \occsetstr\left[\left(\fmy,i\right),\occset \right]
\end{flalign*}
defined through the construction below. We say that: $\occsetstr\left[\occset ,\left(\fmy,i\right)\right]$ is the occurrence set deduced on $\fmx$ by $\occset$; $\occfunstr\left[\occfun \left(\fmy,i\right)\right]$ is the occurrence function induced on $\occsetstr\left[\occset, \left(\fmy,i\right)\right]$ by $\occfun$.\\[8pt] 
Since we deal with a finite number of elements belonging to $\fremag$ (i.e. $\fmx$, $\fmy$), of occurrence sets (i.e. $\occset$), of occurrence functions (i.e. $\occfun$), by Proposition \ref{propfinordseq}-[2] there is no loss of generality by assuming that there are
$k_{\fmx},k_{\fmy} \in \mathbb{N}_0$ , $\fmz_{\fmx}\in \fremag_{k_{\fmx}}$, $\fmz_{\fmy}\in \fremag_{k_{\fmy}}$, an occurrence set $\extoccset$ for $\fmz_{\fmy}$, 
a set function $\extoccfun:\extoccset\rightarrow \extfremag_{k_{\fmy}}$ such that: data $\fmx$, $\fmy$, $i$, $k_{\fmx}$, $k_{\fmy}$, $\fmz_{\fmx}$, $\fmz_{\fmy}$ fulfill \eqref{sysocc}; data $\fmy$, $\occset$, $k_{\fmy}$, $\fmz_{\fmy}$, $\extoccset$ fulfill \eqref{hypextelcomp}; data $\fmy$, $\occset$, $\occfun$, $k_{\fmy}$, $\fmz_{\fmy}$, $\extoccset$, $\extoccfun$ fulfill \eqref{hypasselcomp}.\newline
By exploiting Propositions \ref{posrel} we define: 
\begin{flalign*}
&\begin{array}{l}
\extoccsetstr\left[\extoccset,\left(\fmz_{\fmy},i\right) \right]=\left\{\left(\fmz,j+i-1\right)\;:\;\left(\fmz,j\right)\in\extoccset\right\}\text{;}
\end{array}\\[8pt]
&\begin{array}{l}
\text{the set function}\quad\extasselstr\left[\extoccfun,\left(\fmz_{\fmy},i\right) \right]: \extoccsetstr\left[\extoccset,\left(\fmz_{\fmy},i\right) \right]\rightarrow \extfremag_{k_{\fmx}}\quad\text{by setting}\\[4pt] 
\extasselstr\left[\extoccfun,\left(\fmz_{\fmy},i\right),\right]\left(\fmz,j \right)=\extoccfun\left(\fmz,j-i+1 \right)\qquad \forall \left(\fmz,j \right)\in \extoccsetstr\left[\extoccset,\left(\fmz_{\fmy},i\right)\right]\text{.}
\end{array}
\end{flalign*}
If $\extoccsetstr\left[\extoccset,\left(\fmz_{\fmy},i\right) \right]=\udenset$ then we set: $\occsetstr\left[\occset,\left(\fmy,i\right) \right]=\udenset$; $\occfunstr\left[\occfun,\left(\fmy,i\right) \right]=\emptysetfun$.\newline  
If $\extoccsetstr\left[\extoccset,\left(\fmz_{\fmy},i\right) \right]\neq\udenset$ then we set:  
\begin{flalign*}
&\occsetstr\left[\occset,\left(\fmy,i\right) \right]= \left\{\left(\incllim_{k\left[z,j\right]} \left(\fmz\right),j\right)\,:\;
 \left(\fmz,j\right)\in \extoccsetstr\left[\extoccset,\left(\fmz_{\fmy},i\right) \right]\right\}\text{;}\\[8pt]
&\occfunstr\left[\occfun,\left(\fmy,i\right) \right]\left(\incllim_{k\left[\fmz,j\right]} \left(\fmz\right),j \right)=\incllim_{k\left[\fmz,j\right]}\left(
\extasselstr\left[\extoccfun,\left(\fmz_{\fmy},i\right) \right]\left(\fmz,j \right)\right)\hspace{10pt} \forall \left(\fmz,j \right)\in \extoccsetstr\left[\extoccset,\left(\fmz_{\fmy},i\right) \right]\text{.}
\end{flalign*}
\end{construction}

In Construction \ref{subevalsm} below we introduce the notion of induced occurrence set and induced occurrence function.
\begin{construction}\label{subevalsm}
Fix $\fmx\in \fremag$, an occurrence pair $\left(\fmy,i \right)$ for $\fmx$, an occurrence set $\occset$ for $\fmx$, an occurrence function $\occfun$ on $\occset$.\newline
We denote by 
\begin{flalign*}
&\occsetstr\left[\left(\fmy,i\right),\occset \right] \;\text{the occurrence set for}\;\fmy\text{,}\\[8pt]
&\occfunstr\left[\left(\fmy,i\right),\occfun \right]\;\text{the occurrence function on}\; \occsetstr\left[\left(\fmy,i\right),\occset \right]
\end{flalign*}
defined through the construction below. We say that: $\occsetstr\left[\left(\fmy,i\right),\occset \right]$ is the occurrence set induced on $\fmy$ by $\occset$; $\occfunstr\left[\left(\fmy,i\right),\occfun \right]$ is the occurrence function induced on $\occsetstr\left[\left(\fmy,i\right),\occset \right]$ by $\occfun$.\\[8pt] 
Since we deal with a finite number of elements belonging to $\fremag$ (i.e. $\fmx$, $\fmy$), of occurrence sets (i.e. $\occset$), of occurrence functions (i.e. $\occfun$), by Proposition \ref{propfinordseq}-[2] there is no loss of generality by assuming that there are
$k_{\fmx},k_{\fmy} \in \mathbb{N}_0$, $\fmz_{\fmx}\in \fremag_{k_{\fmx}}$, $\fmz_{\fmy}\in \fremag_{k_{\fmy}}$, an occurrence set $\extoccset$ for $\fmz_{\fmx}$, 
a set function $\extoccfun:\extoccset\rightarrow \extfremag_{k_{\fmx}}$ such that: data $\fmx$, $\fmy$, $i$, $k_{\fmx}$, $k_{\fmy}$, $\fmz_{\fmx}$, $\fmz_{\fmy}$ fulfill \eqref{sysocc}; data $\fmx$, $\occset$, $k_{\fmx}$, $\fmz_{\fmx}$, $\extoccset$ fulfill \eqref{hypextelcomp}; data $\fmx$, $\occset$, $\occfun$, $k_{\fmx}$, $\fmz_{\fmx}$, $\extoccset$, $\extoccfun$ fulfill \eqref{hypasselcomp}.\newline
By exploiting Propositions \ref{posrel} we define: 
\begin{flalign*}
&\begin{array}{l}
\extoccsetstr\left[\left(\fmz_{\fmy},i\right),\extoccset \right]=\\[4pt]
\left\{\left(\fmz,j-i+1\right)\,:\;\left(\fmz,j\right) \in\extoccset \text{,}\quad i\leq j<j+ \extlnght_{k_{\fmx}}\left(\fmz\right)\leq i+ \extlnght_{k_{\fmx}} \left(\fmz_{\fmy}\right)\right\}\text{;}
\end{array}\\[8pt]
&\begin{array}{l}
\text{the set function}\quad\extasselstr\left[\left(\fmz_{\fmy},i\right),\extoccfun \right]: \extoccsetstr\left[\left(\fmz_{\fmy},i\right),\extoccset \right]\rightarrow \extfremag_{k_{\fmy}}\quad\text{by setting}\\[4pt] 
\extasselstr\left[\left(\fmz_{\fmy},i\right),\extoccfun \right]\left(\fmz,j \right)=\extoccfun\left(\fmz,j+i-1 \right)\qquad \forall \left(\fmz,j \right)\in \extoccsetstr\left[\left(\fmz_{\fmy},i\right),\extoccset \right]\text{.}
\end{array}
\end{flalign*}
If $\extoccsetstr\left[\left(\fmz_{\fmy},i\right),\extoccset \right]=\udenset$ then we set: $\occsetstr\left[\left(\fmy,i\right),\occset \right]=\udenset$; $\occfunstr\left[\left(\fmy,i\right),\occfun \right]=\emptysetfun$.\newline  
If $\extoccsetstr\left[\left(\fmz_{\fmy},i\right),\extoccset \right]\neq\udenset$ then we set:  
\begin{flalign*}
&\occsetstr\left[\left(\fmy,i\right),\occset \right]= \left\{\left(\incllim_{k\left[z,j\right]} \left(\fmz\right),j\right)\,:\;
 \left(\fmz,j\right)\in \extoccsetstr\left[\left(\fmz_{\fmy},i\right),\extoccset \right]\right\}\text{;}\\[8pt]
&\occfunstr\left[\left(\fmy,i\right),\occfun \right]\left(\incllim_{k\left[\fmz,j\right]} \left(\fmz\right),j \right)=\incllim_{k\left[\fmz,j\right]}\left(
\extasselstr\left[\left(\fmz_{\fmy},i\right),\extoccfun \right]\left(\fmz,j \right)\right)\hspace{10pt} \forall \left(\fmz,j \right)\in \extoccsetstr\left[\left(\fmz_{\fmy},i\right),\extoccset \right]\text{.}
\end{flalign*}
\end{construction}

In Propositions \ref{elpropC1}-\ref{elpropC3} below we study how Construction \ref{subevalsm} specializes in case of elementary structure of $\fmx$, in such cases we also prove a sort of inversion of Construction \ref{subevalsm} (statements 2). 

\begin{proposition}\label{elpropC1}\mbox{}
\begin{enumerate}
\item Fix $ \fmx_1,\fmx_2 \in \fremag$, an occurrence set $\occset$ for $\fmx_1 \compone \fmx_2$, an occurrence function $\occfun$ on $\occset$. Then either $ \occset =\left\{\left(\fmx_1 \compone \fmx_2,1\right)\right\}$ or there are
\begin{equation*}
\begin{array}{ll}
l_1,l_{2} \in \mathbb{N}\text{,}\vspace{4pt}\\
\text{an occurrence set}\;\occset_j\; \text{for}\;\fmx_j & \forall j\in \left\{1,2\right\}\text{,}\vspace{4pt}\\
\text{an occurrence function}\;\occfun_j\;\text{on} \;\occset_j& \forall j\in \left\{1,2\right\}\text{,} 
\end{array}
\end{equation*}
such that all conditions below are fulfilled:
\begin{equation*}
\begin{array}{l}
\occset= \left\{\left(\fmy,i+l_j \right),\quad  \left(\fmy,i\right)\in \occset_j, \quad j\in\left\{1,2\right\}  \right\}\text{,}\vspace{4pt}\\
\occfun_j\left(\fmy,i\right)=\occfun\left(\fmy,i+l_j\right) \hspace{20pt} \forall \left(\fmy,i\right) \in \occset_j\quad \forall j\in\left\{1,2\right\} \text{,}\vspace{4pt}\\
\strchar\left[\fmx_1 \compone \fmx_2,\occfun\right]=
\strchar\left[\fmx_1, \occfun_1\right] \compone \strchar\left[\fmx_2, \occfun_2\right]\text{.}
\end{array}
\end{equation*} 
\item Fix $\fmx, \fmy_1,\fmy_2\in\fremag$, an occurrence set $\occset$ for $\fmx$, an occurrence function $\occfun$ on $\occset$. Assume $ \fmy_1 \compone \fmy_2=\strchar[\fmx,\occfun]$.\newline
Then either $ \occset =\left\{\left(\fmx,1\right)\right\}$ or there are 
\begin{equation*}
\begin{array}{ll}
\fmx_1,\fmx_2\in\fremag\text{,}\vspace{4pt}\\
\text{an occurrence set}\;\occset_j\; \text{for}\;\fmx_j & \forall j\in \left\{1,2\right\}\text{,}\vspace{4pt}\\
\text{an occurrence function}\;\occfun_j\;\text{on}\;\occset_j& \forall j\in \left\{1,2\right\}\text{,} 
\end{array}
\end{equation*}
such that:\hspace{20pt}$\fmx= \fmx_1\compone \fmx_2$;\hspace{20pt} $\fmy_j=\strchar\left[\fmx_j, \occfun_j\right] \quad \forall j\in \left\{1,2\right\}$.
\end{enumerate}
\end{proposition}
\begin{proof}\mbox{}\newline
\textnormal{\textbf{Proof of statement 1.}} \ \ Since we deal with a finite number of elements belonging to $\fremag$ (i.e. $\fmx_1\compone \fmx_2$), of occurrence sets (i.e. $\occset$), of occurrence functions (i.e. $\occfun$), by Proposition \ref{propfinordseq}-[2] there is no loss of generality by assuming that there are
$k \in \mathbb{N}_0$, $\fmz\in \fremag_{k}$, an occurrence set $\extoccset$ for $\fmz$, 
a set function $\extoccfun:\extoccset\rightarrow \extfremag_{k}$ such that: data $\fmx_1\compone \fmx_2$, $\occset$, $k$, $\fmz$, $\extoccset$ fulfill \eqref{hypextelcomp}; data $\fmx_1\compone \fmx_2$, $\occset$, $\occfun$, $k$, $\fmz$, $\extoccset$, $\extoccfun$ fulfill \eqref{hypasselcomp}; $\strchar\left[\fmx_1 \compone \fmx_2,\occfun\right]=\incllim_k \left(\extstrchar\left[\fmz,\extoccfun\right]\right)$.\newline
We set $\check{k}=\min\left\{l\in \mathbb{N}_0\,:\quad \forall j \in \left\{1,2\right\}\quad\exists \fmy \in \fremag_l\quad\text{with} \quad\incllim_l\left(\fmy\right)=\fmx_j \right\}$.\newline 
Then:
\begin{equation*}
\left\{
\begin{array}{lll}
(i)&\exists!\left(\fmz_1,\fmz_2\right) \in \left(\fremag_{\check{k}}\right)^2\,:\;\left(\incllim_{\check{k}}\left(\fmz_1\right),\incllim_{\check{k}}\left(\fmz_2\right)\right)=\left(\fmx_1,\fmx_2\right)\text{,}\vspace{4pt}\\
(ii)&k \geq \check{k}\text{,}\vspace{4pt}\\
(iii)&\fmz=  \fmz_1\compone_{\!\check{k}} \fmz_2 & \text{if}\;k=\check{k} \text{,}\vspace{4pt}\\
(iv)&\fmz=\left(\inclfremaglev_k\funcomp...\funcomp \inclfremaglev_{\check{k}+1} \right) \left(  \fmz_1\compone_{\!\check{k}} \fmz_2 \right)& \text{if}\;k>\check{k} \text{.}
\end{array}
\right.
\end{equation*}
Since\newline
\centerline{$\left(\inclfremaglev_k\funcomp...\funcomp \inclfremaglev_{\check{k}+1} \right) \left(  z_1\compone_{\!\check{k}} \fmz_2 \right)=\left(\left(\inclfremaglev_k\funcomp...\funcomp \inclfremaglev_{\check{k}+1} \right)  \left( \fmz_1\right)\right)\compone_{\!\check{k}}\left(\left(\inclfremaglev_k\funcomp...\funcomp \inclfremaglev_{\check{k}+1} \right)\left( \fmz_2\right)\right)$,}\newline
up to renaming $\left(\inclfremaglev_k\funcomp...\funcomp \inclfremaglev_{\check{k}+1} \right)  \left( \fmz_j\right)$ again by $\fmz_j$ for any $j\in \left\{1,2\right\}$, there is no loss of generality by assuming $k = \check{k}$.\newline
Eventually statement follows by applying Construction \ref{subevalsm} twice: first with data $\fmx_1\compone \fmx_2$, $\left(\fmx_1,1\right)$, $\occset$, $\occfun$, then with data $\fmx_1\compone \fmx_2$, $\left(\fmx_2,\extlnght_k\left(\fmz_1\right)+2\right)$, $\occset$, $\occfun$.\newline 
\textnormal{\textbf{Proof of statement 2.}} \ \ Statement follows by Propositions \ref{invsost}, \ref{propfinordseq} and statement 1.
\end{proof}

\begin{proposition}\label{elpropC2}\mbox{}
\begin{enumerate}
\item  Fix $ \fmx_1,...,\fmx_n \in \fremag$, an occurrence set $\occset$ for $\lboundone \fmx_1,...,\fmx_n\rboundone$, an occurrence function $\occfun$ on $\occset$. Then either $\occset=\left\{\left(\lboundone \fmx_1,...,\fmx_n\rboundone,1\right)\right\}$ or there are
\begin{flalign*}
&l_1,...,l_{n} \in \mathbb{N}\text{,}\\[4pt]
&\text{an occurrence set}\;\occset_{\mathsf{n}}\; \text{for}\;\fmx_{\mathsf{n}}  \hspace{67pt}\forall {\mathsf{n}}\in \left\{1,...,n\right\}\text{,}\\[4pt]
&\text{an occurrence function}\;\occfun_{\mathsf{n}}\;\text{on} \;\occset_{\mathsf{n}} \hspace{40pt}\forall {\mathsf{n}}\in \left\{1,...,n\right\}\text{,} 
\end{flalign*}
such that:
\begin{flalign*}
&\occset= \left\{\left(\fmy,i+l_{\mathsf{n}} \right),\quad  \left(\fmy,i\right)\in \occset_{\mathsf{n}}, \quad {\mathsf{n}}\in\left\{1,...,n\right\}  \right\}\text{,}\\[4pt]
&\occfun_{\mathsf{n}}\left(\fmy,i\right)=\occfun\left(\fmy,i+l_{\mathsf{n}}\right) \hspace{20pt} \forall \left(\fmy,i\right) \in \occset_{\mathsf{n}}\quad \forall {\mathsf{n}}\in\left\{1,...,n\right\} \text{,}\\[4pt]
&\strchar\left[\lboundone \fmx_1,...,\fmx_n\rboundone,\occfun\right]=
\lboundone\strchar\left[\fmx_1, \occfun_1\right],..., \strchar\left[\fmx_n, \occfun_n\right]\rboundone\text{.}
\end{flalign*}
\item Fix $n \in \mathbb{N}$, $\fmx, \fmy_1,...,\fmy_n\in\fremag$, an occurrence set $\occset$ for $\fmx$, an occurrence function $\occfun$ on $\occset$. Assume $\lboundone \fmy_1,...,\fmy_n\rboundone=\strchar[\fmx,\occfun]$.\newline
Then either $\occset=\left\{\left(\fmx, 1\right)\right\}$ or there are 
\begin{flalign*}
&\fmx_1,...,\fmx_n\in\fremag\text{,}\\[4pt]
&\text{an occurrence set}\;\occset_{\mathsf{n}}\; \text{for}\;\fmx_{\mathsf{n}} \hspace{67pt} \forall {\mathsf{n}}\in \left\{1,...,n\right\}\text{,}\\[4pt]
&\text{an occurrence function}\;\occfun_{\mathsf{n}}\;\text{on}\; \occset_{\mathsf{n}}\hspace{40pt} \forall {\mathsf{n}}\in \left\{1,...,n\right\}\text{,} 
\end{flalign*}
such that:\hspace{20pt}$x=\lboundone \fmx_1,...,\fmx_n\rboundone$;\hspace{20pt} $\fmy_{\mathsf{n}}=\strchar\left[\fmx_{\mathsf{n}}, \occfun_{\mathsf{n}}\right]\quad \forall {\mathsf{n}}\in \left\{1,...,n\right\} $.
\end{enumerate}
\end{proposition}
\begin{proof}\mbox{}\newline
\textnormal{\textbf{Proof of statement 1.}} \ \ We refer to Proposition \ref{linincl}-[1] both for the statement and for the notation.\newline
Since we deal with a finite number of elements belonging to $\fremag$ (\!i.e. $\lboundone \fmx_1,...,\fmx_n\rboundone$\!), of occurrence sets (i.e. $\occset$), of occurrence functions (i.e. $\occfun$), by Proposition \ref{propfinordseq}-[2] there is no loss of generality by assuming that there are
$k \in \mathbb{N}_0$, $\fmz\in \fremag_{k}$, an occurrence set $\extoccset$ for $\fmz$, 
a set function $\extoccfun:\extoccset\rightarrow \extfremag_{k}$ such that: data $\lboundone \fmx_1,...,\fmx_n\rboundone$, $\occset$, $k$, $\fmz$, $\extoccset$ fulfill \eqref{hypextelcomp}; data $\lboundone \fmx_1,...,\fmx_n\rboundone$, $\occset$, $\occfun$, $k$, $\fmz$, $\extoccset$, $\extoccfun$ fulfill \eqref{hypasselcomp}; $k\geq \check{k}+1$;  $\strchar\left[\lboundone \fmx_1,...,\fmx_n\rboundone,\occfun\right]=\incllim_k \left(\extstrchar\left[\fmz,\extoccfun\right]\right)$.\newline
By construction we have
\begin{equation*}
\left\{\!
\begin{array}{lll}
(i)\!&\!\fmz= \left( \left(\bsfuno, \fmz_1\right),...,\left(\bsfuno, \fmz_n\right) \right)& \text{if}\;k=\check{k}+1 \text{,}\vspace{4pt}\\
(ii)\!&\!\fmz=\left(\inclfremaglev_k\funcomp...\funcomp \inclfremaglev_{\check{k}+2} \right) \left( \left(\bsfuno, \fmz_1\right),...,\left(\bsfuno, \fmz_n\right) \right)& \text{if}\;k>\check{k}+1 \text{,}\vspace{4pt}\\
(iii)\!&\!\left( \left(\bsfuno, \fmz_1\right),...,\left(\bsfuno, \fmz_n\right) \right)\in \fremaggen_{\check{k}+1}\text{,}\vspace{4pt}\\
(iv)\!&\!\left(\inclfremaglev_l\funcomp...\funcomp \inclfremaglev_{\check{k}+2} \right) \left( \left(\bsfuno, \fmz_1\right),...,\left(\bsfuno, \fmz_n\right) \right)\in \fremaggen_l&\forall l \in \left\{\check{k}+2,...,k\right\}\text{.}
\end{array}
\right.
\end{equation*}
Eventually statement follows by applying Construction \ref{subevalsm} with data\newline
\centerline{$\lboundone \fmx_1,...,\fmx_n\rboundone\text{,}\qquad \left(\fmx_{\mathsf{n}},l_{\mathsf{n}}\right)\text{,}\qquad \occset\text{,}\qquad \occfun$}\newline
for any ${\mathsf{n}} \in \left\{1,...,n\right\}$, where 
\begin{equation*}
l_{\mathsf{n}}=\left\{
\begin{array}{ll}
3\left(k-\check{k}-1\right)+2&\text{if}\;{\mathsf{n}}=1\text{,}\vspace{4pt}\\
\overset{{\mathsf{n}}-1}{\underset{a=1}{\sum}} \extlnght_k\left(\fmz_a\right)+5\left({\mathsf{n}}-1\right)+3\left(k-\check{k}-1\right)+2& \text{if}\;{\mathsf{n}}\in \left\{2,...,n\right\}\text{.}
\end{array}
\right.
\end{equation*}
\textnormal{\textbf{Proof of statement 2.}} \ \ Statement follows by Propositions \ref{invsost}, \ref{propfinordseq} and statement 1.
\end{proof}

\begin{proposition}\label{elpropC3}\mbox{} 
\begin{enumerate}
\item Fix $\fmx\in\fremag$, $\bsfm \in \bsfM$, an occurrence set $\occset$ for $\bsfm \acmone \fmx$, an occurrence function $\occfun$ on $\occset$. Then either $\occset=\left\{\left(\bsfm \acmone \fmx,1\right)\right\}$ or there are
\begin{flalign*}
&l \in \mathbb{N}\text{,}\\[4pt]
&\text{an occurrence set}\;\occset_1\; \text{for}\;\fmx \text{,}\\[4pt]
&\text{an occurrence function}\;\occfun_1\;\text{on}\;\occset_1\text{,} 
\end{flalign*}
such that:
\begin{flalign*}
&\occset= \left\{\left(\fmy,i+l \right),\quad  \left(\fmy,i\right)\in \occset_1 \right\}\text{,}\\[4pt]
&\occfun_1\left(\fmy,i\right)=\occfun\left(\fmy,i+l\right) \hspace{20pt} \forall \left(\fmy,i\right) \in \occset_1 \text{,}\\[4pt]
&\strchar\left[\bsfm \acmone \fmx,\occfun\right]=
\bsfm \acmone\strchar\left[\fmx, \occfun_1\right]\text{.}
\end{flalign*} 
\item Fix $\bsfm \in \bsfM$, $\fmx,\fmy\in\fremag$, an occurrence set $\occset$ for $\fmx$, an occurrence function $\occfun$ on $\occset$. Assume $\bsfm \acmone \fmy=\strchar[\fmx,\occfun]$.\newline
Then either $ \occset =\left\{\left(\fmx,1\right)\right\}$ or there are 
\begin{flalign*}
&\fmz\in\fremag\text{,}\\[4pt]
&\text{an occurrence set}\;\occsetdue\; \text{for}\;\fmz \text{,}\\[4pt]
&\text{an occurrence function}\;\occfundue\;\text{on} \;\occsetdue\text{,} 
\end{flalign*}
such that:\hspace{20pt}$\fmx= \bsfm \acmone \fmz$;\hspace{20pt} $\fmy=\strchar\left[\fmz, \occfundue\right]$.
\end{enumerate}
\end{proposition}
\begin{proof}\mbox{}\newline
\textnormal{\textbf{Proof of statement 1.}} \ \ We refer to Proposition \ref{linincl}-[2] both for the statement and for the notation.\newline
Since we deal with a finite number of elements belonging to $\fremag$ (i.e. $\bsfm \acmone \fmx$), of occurrence sets (i.e. $\occset$), of occurrence functions (i.e. $\occfun$), by Proposition \ref{propfinordseq}-[2] there is no loss of generality by assuming that there are
$k \in \mathbb{N}_0$ , $\fmw\in \fremag_{k}$, an occurrence set $\extoccset$ for $\fmw$, 
a set function $\extoccfun:\extoccset\rightarrow \extfremag_{k}$ such that: data $\bsfm \acmone \fmx$, $\occset$, $k$, $\fmw$, $\extoccset$ fulfill \eqref{hypextelcomp}; data $\bsfm \acmone \fmx$, $\occset$, $\occfun$, $k$, $\fmw$, $\extoccset$, $\extoccfun$ fulfill \eqref{hypasselcomp}; $k\geq \check{k}+1$; $\strchar\left[\bsfm \acmone \fmx,\occfun\right]=\incllim_k \left(\extstrchar\left[\fmw,\extoccfun\right]\right)$.\newline
By construction we have
\begin{equation*}
\left\{
\begin{array}{lll}
(i)&\fmw=  \left(\bsfm, \fmz\right)& \text{if}\;k=\check{k}+1 \text{,}\vspace{4pt}\\
(ii)&\fmz=\left(\inclfremaglev_k\funcomp...\funcomp \inclfremaglev_{\check{k}+2} \right) \left(\bsfm, \fmz\right)& \text{if}\;k>\check{k}+1 \text{,}\vspace{4pt}\\
(iii)&\left(\bsfm, \fmz\right)\in \fremaggen_{\check{k}+1}\text{,}\vspace{4pt}\\
(iv)&\left(\inclfremaglev_l\funcomp...\funcomp \inclfremaglev_{\check{k}+2} \right) \left(\bsfm, \fmz\right)\in \fremaggen_l&\forall l \in \left\{\check{k}+2,...,k\right\}\text{.}
\end{array}
\right.
\end{equation*}
Eventually statement follows by applying Construction \ref{subevalsm} with data $\bsfm \acmone \fmx$, $\left(\fmx,3\left(k-\check{k}\right)\right)$, $\occset$, $\occfun$ for any $i \in \left\{1,...,n\right\}$.\newline
\textnormal{\textbf{Proof of statement 2.}} \ \ Statement follows by Propositions \ref{invsost}, \ref{propfinordseq} and statement 1.
\end{proof}

In Proposition \ref{elpropC4} below we study the structure of sets $\occset_{\fmx}$ in case of elementary structure of $\fmx$. 

\begin{proposition}\label{elpropC4}
Fix $\fmx_1,...,\fmx_n, \fmx, \fmy \in \fremag$, $\bsfm\in \bsfM$.  Assume that $\fmy$ occurs in $\fmx$. Refer to Definition \ref{defpath}-[4]. Then: 
\begin{equation*}
\begin{array}{c}
\occset_{\fmx_1\compone \fmx_2}=\occset_{\fmx_1}\cup  \occset_{\fmx_2}\text{,}\hspace{40pt}\occset_{\lboundone \fmx_1,...,\fmx_n\rboundone}=\overset{n}{\underset{\mathsf{n}=1}{\bigcup}}\,\occset_{\fmx_{\mathsf{n}}} \text{,}\vspace{4pt}\\
\occset_{\mbox{\footnotesize{\bsfm}} \acmone \fmx}=\occset_{\fmx}\text{,}\hspace{40pt} \occset_{\fmy}\subseteq\occset_{\fmx}\text{.}
\end{array}
\end{equation*}
\end{proposition}
\begin{proof}Statement follows by Propositions \ref{elpropC1}, \ref{elpropC2}, \ref{elpropC3}
applied to occurrence sets $\occset^{\fmx}$ for any $\fmx\in \fremag$.
\end{proof}

In Definition \ref{Cinftybis}, Proposition \ref{lininclbis} below we introduce and study two sub-magmas of $\fremag$ which play a central role to establish links between traditional integro-differential calculus on smooth functions and integral calculus on continuous functions. We refer to Definition \ref{intdiffmon}.

\begin{definition}\label{Cinftybis}\mbox{}
\begin{enumerate}
\item We denote by $\fremagcont$ the sub-magma of $\left(\fremag,\compone\right)$ fulfilling all conditions below:
\begin{flalign}
&\left\{
\begin{array}{l}
\fmx\in \fremagcont\;\text{for any}\;\fmx \in \fremag \;\text{fulfilling}\\[4pt]
\hspace{50pt}\begin{array}{ll}
\text{either}& \Fpart_j \;\text{does not occur in}\; \fmx\; \text{for any}\;j \in \mathbb{N}\text{,} \vspace{4pt}\\
\text{or} &f \;\text{does not occur in}\; \fmx\; \text{for any}\; f \in  \bascont\text{;}
\end{array}
\end{array}
\right.\\[12pt]
&\begin{array}{l}
\lboundone \fmx_1,...,\fmx_n\rboundone\in \fremagcont\hspace{32pt} \forall n \in \mathbb{N}\quad \forall \left(\fmx_1,...,\fmx_n\right) \in \fremagcont^n \text{;}
\end{array}\\[12pt]
&\begin{array}{l}
\bsfm\acmone \fmx\in \fremagcont \hspace{57pt}\forall \bsfm \in \subbsfM\quad \forall \fmx \in \fremagcont \text{.}
\end{array}
\end{flalign}
\item 
$\fremagsmooth= \left\{\fmx \in \fremag \,:\; \;f \;\text{does not occur in}\; \fmx\;\text{for any}\; f \in  \bascont\right\}$.
\end{enumerate} 
\end{definition}

\begin{proposition}\label{lininclbis}\mbox{}
\begin{enumerate}
\item $\fremagsmooth$ is a sub-magma of $\fremagcont$.
\item There is one and only one set function $\evalcompone :\fremagcont \rightarrow \Cksp{0}$ recursively defined by setting:
\begin{equation}
\left\{
\begin{array}{lll}
(i)&\evalcompone\left(\fmx\compone \fmy		\right)= \evalcompone\left(\fmx\right)\funcomp \evalcompone\left(\fmy\right)&\forall \fmx,\fmy \in \fremagcont\text{,}\vspace{4pt}\\
(ii)&\evalcompone\left(\lboundone \fmx_1,...,\fmx_n\rboundone\right)=\vspace{4pt}\\
&\hspace{40pt}\left(\evalcompone\left(\fmx_1\right),...,\evalcompone\left(\fmx_n\right)\right)&\forall  \left(\fmx_1,...,\fmx_n \right)\in \fremagcont^n \text{,}\vspace{4pt}\\
(iii)&\evalcompone\left(\bsfm \acmone\fmx\right)=\bsfm \contspbsfmu\evalcompone\left(\fmx\right)&\forall \bsfm\in\subbsfM\quad\forall \fmx \in \fremagcont\text{,}\vspace{4pt}\\
(iv)&\evalcompone\left(\bsfm\acmone \fmx\right)=\bsfm \smospbsfmu \evalcompone\left(\fmx\right)&\forall \bsfm\in\bsfM\quad\forall \fmx \in \fremagsmooth\text{,}\vspace{4pt}\\
(v)&\evalcompone\left(\incllim_{-1}\left(f\right) \right)=f&\forall f \in \fremag_{-1}\text{.}
\end{array}
\right.\label{evcmppropr}
\end{equation}
We set:
\begin{flalign*}
&\fremagcontcont=\left\{\fmx \in\fremagcont\;:\; \dommagone\left(\fmx\right)=\Cae\left(\evalcompone\left(\fmx	\right)\right)\right\}\text{;}\\[4pt]
&\fremagsmoothsmooth=\left\{\fmx \in\fremagsmooth\;:\;  \dommagone\left(\fmx\right)=\Cae\left(\evalcompone\left(\fmx	\right)\right)\right\}.
\end{flalign*}
\item There is a set function $\lininclone : \Cksp{0}\setminus\left\{\contspempty\right\}\rightarrow\fremag$ such that:
\begin{equation}\label{proplininclone}
\left\{
\begin{array}{lll}
(i)&\lininclone\left(f\right)=\incllim_{-1}\left(f\right)& \forall f \in \fremag_{-1}  \text{,}\\[4pt]
(ii)& \lininclone\left(\Cksp{0}  \setminus \{\contspempty \}\right)\subseteq \fremagcontcont \text{,}\\[4pt]
(iii)& \evalcompone \funcomp \lininclone = \idobj+\Cksp{0} \setminus \left\{\contspempty\right\}+\text{.}
\end{array}
\right.
\end{equation} 
\item Fix $n\in\mathbb{N}_0$, $\bsfm \in \bsfM$, $\bsfn \in \subbsfM$, $\left(\fmx_1,...,\fmx_n\right)\in \fremag^n$, $\fmx\in \fremag$.\newline
If $\fmx_1,...,\fmx_n, \fmx\in \fremagcontcont$ then $\lboundone \fmx_1,...,\fmx_n\rboundone, \bsfn \acmone \fmx \in \fremagcontcont$.\newline 
If $\fmx_1,...,\fmx_n, \fmx\in \fremagsmooth$ then $\lboundone \fmx_1,...,\fmx_n\rboundone, \bsfm \acmone \fmx \in \fremagsmooth$.\newline
If $\fmx_1,...,\fmx_n, \fmx\in \fremagsmoothsmooth$ then $\lboundone \fmx_1,...,\fmx_n\rboundone, \bsfm \acmone \fmx \in \fremagsmoothsmooth$. 
\end{enumerate}
\end{proposition}
\begin{proof}
Statements 1, 2, 4 follow directly by definition and construction. Then we prove statement 3 by giving an explicit construction of $\lininclone$.\newline
If $f \in \fremag_{-1}$ then we define $\lininclone\left(f\right)=\incllim_{-1}\left(f\right) $. \newline 
If $f \in \Cksp{0} \setminus \fremag_{-1}$ then, by Remark \ref{remch}, there are $m, n \in \mathbb{N}$, $\intuno \subseteqdentro \mathbb{R}^m$ with 
$\Dom\left(f\right)=\intuno$ and $\Cod\left(f\right)=\mathbb{R}^n$. Definition of $\bascont\left(\intuno\right)$ (see \eqref{Bchoice}) entails that for any ${\mathsf{n}}\in\left\{1,...,n	\right\}$ there are one and only one ordered finite sequence of not zero real numbers $(a_{{\mathsf{n}},{\mathsf{l}}})_{{\mathsf{l}}=1}^{l_{\mathsf{n}}}$, one and only one finite not empty set of functions $\{f_{{\mathsf{n}},{\mathsf{l}}}\}_{{\mathsf{l}}=1}^{l_{\mathsf{n}}} \in (\Cksp{\infty}(\intuno)+\mathbb{R}+ \cup \bascont\left(\intuno\right))$ such that $f_{\mathsf{n}}=\overset{l_{\mathsf{n}}}{\underset{{\mathsf{l}}=1}{\sum}}(a_{{\mathsf{n}},{\mathsf{l}}}\, f_{{\mathsf{n}},{\mathsf{l}}})$. We set:
\begin{flalign*}
&l_0=0\text{;}\\[12pt]
&\begin{array}{l}
g \in \Cksp{\infty}(\mathbb{R}^{l_1+...+l_n})+\mathbb{R}^n+\;\text{defined by setting}\\
g_{\mathsf{n}}\left(\unkuno\right)=\overset{l_{\mathsf{n}}}{\underset{{\mathsf{l}}=1}{\sum}}(a_{{\mathsf{n}},{\mathsf{l}}}\, s_{l_0+...+l_{{\mathsf{n}}-1}+{\mathsf{l}}}) \qquad \forall {\mathsf{n}}\in\left\{1,...,n\right\}\quad \forall \unkuno \in \mathbb{R}^{l_1+...+l_n}\text{;}
\end{array}\\[12pt]
&\begin{array}{l}
\fmx=(\fmx_1,...,\fmx_n)\in \fremaggen_0\;\text{defined by setting}\\
\fmx_{\mathsf{n}}=\incllim_{-1}\left(g_{\mathsf{n}}\right)\qquad \forall {\mathsf{n}}\in\left\{1,...,n\right\}\text{;}
\end{array}\\[12pt]
&\begin{array}{l}
\fmy=(\fmy_1,...,\fmy_{l_1+...+l_n})\in \fremaggen_0\;\text{defined by setting}\\
\fmy_{l_0+...+l_{{\mathsf{n}}-1}+{\mathsf{l}}}=f_{{\mathsf{n}},{\mathsf{l}}}\qquad \forall {\mathsf{n}}\in\left\{1,...,n	\right\}\quad \forall {\mathsf{l}}\in\left\{1,...,l_{\mathsf{n}}\right\}\text{.}
\end{array}
\end{flalign*}
Eventually we define $\lininclone\left(f	\right)= \incllim_{0}\left(\left(\fmx\, \compone_{\!0}\,  \fmy\right) \,\compone_{\!0}\, \Diag{\intuno}{l_1+...+l_n} \right)$.
\end{proof}

Motivated by Propositions \ref{elpropC4}, \ref{lininclbis} we extend Construction \ref{subevalsm} to sets $\occset_{\fmx}$.

\begin{remark}\label{notelpropC4}\mbox{}
\begin{enumerate}
\item Fix $\fmx\in \fremag$, an occurrence pair $\left(\fmy,i \right)$ for $\fmx$. We emphasize that\newline $\occsetstr\left[\left(\fmy,i\right),\occset^{\fmx} \right]=\occset^{\fmy}$.
\item Fix $\fmx \in \fremag$, a set function $\occfun:\occset_{\fmx} \rightarrow \Cksp{0}$.\newline
Then $\occfun$ defines one and only one occurrence function $\occfundue$ on $\occset^{\fmx}$ by setting \begin{equation*}
\occfundue\left(\fmy,i\right)=\lininclone\left(\occfun\left(f\right)\right)\quad \forall \left(\fmy,i\right)\in\occset^{\fmx}\quad  \forall f\in \occset_{\fmx}\;\text{with}\; \fmy=\incllim_{-1}\left(f\right)\text{.}
\end{equation*}
With an abuse of language we denote 
$\strchar \left[\fmx, \occfundue\right]$ by $\strchar \left[\fmx, \occfun\right]$.
\item Fix $\fmx\in \fremag$, an occurrence pair $\left(\fmy,i \right)$ for $\fmx$, a set function $\occfun:\occset_{\fmx}\rightarrow \Cksp{0}$.\newline
We emphasize that $\occset_{\fmy}=\proj<\fremag,\mathbb{N}<>1>\left(\occsetstr\left[\left(\fmy,i\right),\occset^{\fmx} \right] \right)$.\newline
With an abuse of language we adopt the notation\newline
$\occfunstr\left[\left(\fmy,i\right),\occfun \right]=\occfun\funcomp \incl{\occset_{\fmy}}{\occset_{\fmx}}$.
\end{enumerate}
\end{remark}

In Notation \ref{adjointfun}, Definition \ref{pathinC} below we introduce the notion of path in $\fremag$. This establish the first appearance of topology in $\fremag$ which, up to now, was a purely algebraic structure.

\begin{notation}\label{adjointfun}
Fix a set $A$, a set function $\assfun:A \times \left(-1,1\right)\rightarrow \Cksp{0}$. Assume that all conditions below are fulfilled: 
\begin{flalign*}
&(\assfun=\emptysetfun) \Leftrightarrow \left(A=	\udenset \right)\text{;}\\[4pt]
&\assfun(a,\unkdue)\neq \emptysetfun \hspace{133pt} \forall \left(a,\unkdue\right)\in A \times \left(-1,1\right)\text{;}\\[4pt]
&\Dom(\assfun(a,\unkdue_1))=\Dom(\assfun(a,\unkdue_2))\hspace{33pt}\forall \left(a,\unkdue_1\right),(a,\unkdue_2)\in A\times\left(-1,1\right)\text{;}\\[4pt]
&\Cod\left(\assfun\left(a,\unkdue_1\right)\right)=\Cod\left(\assfun\left(a,\unkdue_2\right)\right)\hspace{31pt}\forall \left(a,\unkdue_1\right),\left(a,\unkdue_2\right)\in A\times\left(-1,1\right)\text{.}
\end{flalign*}
Set $B=\left\{f:\Dom\left(\assfun\left(a,0\right)\right)\times\left(-1,1\right)\rightarrow \Cod\left(\assfun\left(a,0\right)\right)\,:\;a\in A\right\}$, we emphasize that elements belonging to $B$ are not necessarily continuous set functions.\newline
We define the set function $\assfunadj\left[\assfun\right]:A\rightarrow B$ by setting 
\begin{equation*}
\left(\assfunadj\left[\assfun\right]\left(a\right)\right)\left(\unkuno,\unkdue\right)=\left(\assfun\left(a,\unkdue\right)\right)\left(\unkuno\right) \qquad  \forall a\in A\quad \forall \left(\unkuno,\unkdue\right)\in \Dom\left(\assfun\left(a,0\right)\right)\times\left(-1,1\right)\text{.}
\end{equation*}
We call $\assfunadj\left[\assfun\right]$ the adjoint of $\assfun$.
\end{notation}

\begin{definition}\label{pathinC}\mbox{}
\begin{enumerate}
\item Fix $n \in \mathbb{N}$, $\left(\fmx_1,...,\fmx_n\right)\in \fremag^n$, a set function $\assfun:\left(\underset{\mathsf{n}=1}{\overset{n}{\bigcup}} \,\occset_{\fmx_{\mathsf{n}}}\right)\times \left(-1,1\right)  \rightarrow \Cksp{0}$. We say that $\assfun$ is an associating function for $\underset{\mathsf{n}=1}{\overset{n}{\bigcup}} \,\occset_{\fmx_{\mathsf{n}}}$ if and only if it fulfills the system of conditions below
\begin{equation}
\left\{
\begin{array}{lll}
\!(i)\!&\!\assfun=\emptysetfun \Leftrightarrow \underset{\mathsf{n}=1}{\overset{n}{\bigcup}} \,\occset_{\fmx_{\mathsf{n}}}=\udenset\text{,}\vspace{4pt}\\
\!(ii)\!&\!\assfun\left(f,\unkdue\right)\neq\emptysetfun   &\forall \left(f,\unkdue\right)\! \in\! \left(\underset{\mathsf{n}=1}{\overset{n}{\cup}} \,\occset_{\fmx_{\mathsf{n}}}\right)\!\times\! \left(-1,1\right)  \text{,}\vspace{4pt}\\
\!(iii)\!&\!\Dom\left(\assfun \left(f,\unkdue\right) \right)\!=\!\Dom\left(f\right)& \forall \left(f,\unkdue\right) \!\in\! \left(\underset{\mathsf{n}=1}{\overset{n}{\cup}} \,\occset_{\fmx_{\mathsf{n}}}\right)\! \times\!\left(-1,1\right) \text{,}\vspace{4pt}\\
\!(iv)\!&\!\Cod\left(\assfun\left(f,\unkdue\right)\right)=\mathbb{R}&\forall \left(f,\unkdue\right) \!\in\!  \left(\underset{\mathsf{n}=1}{\overset{n}{\cup}} \,\occset_{\fmx_{\mathsf{n}}}\right) \!\times\!\left(-1,1\right) \text{,}\vspace{4pt}\\
\!(v)\!&\!\assfunadj[\assfun]\left(f\right)\in \Cksp{0} & \forall f\in \underset{\mathsf{n}=1}{\overset{n}{\cup}} \,\occset_{\fmx_{\mathsf{n}}} \text{.} 
\end{array}
\right.
\label{pathconduro}
\end{equation}
We say that $\assfun$ is a centered associating function for $\underset{\mathsf{n}=1}{\overset{n}{\bigcup}} \,\occset_{\fmx_{\mathsf{n}}}$ if and only if it fulfills both conditions below:
\begin{flalign}
&\begin{array}{l}
\assfun\hspace{4pt}\text{is}\hspace{4pt}\text{an}\hspace{4pt}\text{associating}\hspace{4pt}\text{function}\hspace{4pt}\text{for}\hspace{4pt}\underset{\mathsf{n}=1}{\overset{n}{\bigcup}} \,\occset_{\fmx_{\mathsf{n}}}\text{;}
\end{array}\\[6pt]
&\begin{array}{l}
\assfun\left(f,0\right)=f\hspace{15pt}\forall f\in \underset{\mathsf{n}=1}{\overset{n}{\cup}} \,\occset_{\fmx_{\mathsf{n}}} \text{.}
\end{array}
\end{flalign}
\item Fix $n \in \mathbb{N}$, a set function $\pmone:\left(-1,1\right) \rightarrow \fremag^n$.\newline
We say that $\pmone$ is a path in $\fremag^n$ if and only if there are 
$\left(\fmx_1,...,\fmx_n\right) \in \fremag^n$, an associating function $\assfun$ for $\underset{\mathsf{n}=1}{\overset{n}{\bigcup}} \,\occset_{\fmx_{\mathsf{n}}}$
such that
\begin{equation}
\hspace{-5pt}
\begin{array}{l}
\pmone\left(t\right)=\strchar\bigg[\left(\fmx_1,...,\fmx_n\right),\vspace{4pt}\\
 \left(\left(\lininclone \funcomp \assfun\right) \funcomp \left(\idobj+\underset{\mathsf{n}=1}{\overset{n}{\cup}} \,\occset_{\fmx_{\mathsf{n}}}+ ,\cost<\underset{\mathsf{n}=1}{\overset{n}{\cup}} \,\occset_{\fmx_{\mathsf{n}}}<>(-1,1)>+\unkdue+\right)\right)\funcomp \Diag{\underset{\mathsf{n}=1}{\overset{n}{\cup}} \,\occset_{\fmx_{\mathsf{n}}}}{2}\bigg] \vspace{4pt}\\
\hspace{245pt}\forall   \unkdue \in \left(-1,1\right)\text{.} 
\end{array}
\label{pathconduro2}
\end{equation}
We say that: $\left(\fmx_1,...,\fmx_n\right)$ is a core of $\pmone$; $\assfun$ is an associating function of $\pmone$; the pair $\left(\left(\fmx_1,...,\fmx_n\right),\assfun\right)$ is a skeleton of $\pmone$.
\item Fix $n \in \mathbb{N}$, a path $\pmone$ in $\fremag^n$, a subset $\singsmpth \subseteq \left(-1,1\right)$.\newline 
We say that $\left(\pmone,\singsmpth\right)$ is a smooth path in $\fremag^n$ if and only if
\begin{flalign}
&\hspace{-10pt}\begin{array}{l}
\text{exists a skeleton}\;\; \left(\left(\fmx_1,...,\fmx_n\right),\assfun\right)\;\; \text{of}\;\; \pmone\;\;\text{such that:} \vspace{4pt}\\
\left\{
\begin{array}{ll}
\!(i)\!&\!0\in \singsmpth\text{,}\vspace{4pt}\\
\!(ii)\!&\!\singsmpth\setminus (-\varepsilon, \varepsilon)\;\; \text{is a finite (maybe empty) set} \;\;\forall \varepsilon>0\text{,}\vspace{4pt}\\
\!(iii)\!&\!\left(\assfunadj[\assfun]\right)\left(f\right)\funcomp \incl{\Dom\left(f\right)\times B}{\Dom\left(f\right)\times\left(-1,1\right)} \in \Cksp{\infty}\vspace{4pt}\\
 &  \forall f\in \underset{\mathsf{n}=1}{\overset{n}{\bigcup}} \,\occset_{\fmx_{\mathsf{n}}}\text{,}\quad \forall \;\text{path component}\;\; B\;\; \text{of} \;\; \left(-1,1\right)\setminus \singsmpth\text{.}
\end{array}
\right.
\end{array}\label{smptass}
\end{flalign}
We say that: $\singsmpth$ is the singular set of $\left(\pmone, \singsmpth\right)$; $\left(\left(\fmx_1,...,\fmx_n\right),\assfun\right)$ is a skeleton of $\left(\pmone,\singsmpth\right)$.\newline
We emphasize that nothing is assumed about $\pmone\left(\unkdue\right)$ when $\unkdue \in \singsmpth$.
\item Fix $n \in \mathbb{N}$, $\left(\fmx_1,...,\fmx_n\right)\in \fremag^n$, a path $\pmone$ in $\fremag^n$.\newline
We say that $\pmone$ is a path in $\fremag^n$ through $\left(\fmx_1,...,\fmx_n\right)$ if and only if we have $\pmone\left(0\right)=\left(\fmx_1,...,\fmx_n\right)$.
\item Fix $\fmx\in \fremag$, a smooth path $\left(\pmone,\left\{0\right\} \right)$ in $\fremag$ through $\fmx$.\newline
We say that $\left(\pmone,\left\{0\right\} \right)$ is a path in $\fremag$ detecting $\fmx$ if and only if there is a centered associating function $\assfun$ for $\occset_{\fmx}$ such that pair $\left(\fmx, \assfun \right)$ is a skeleton of $\pmone$.\newline
We say that the pair $\left(\fmx,\assfun\right)$ is a detecting skeleton of $\left(\pmone,\left\{0\right\} \right)$.\newline 
With an abuse of language we drop any reference to the pair and we refer to $\left(\pmone,\left\{0\right\} \right)$ simply by $\pmone$.
\end{enumerate}
\end{definition}

\begin{remark}\label{notunsk}
We emphasize that:
\begin{enumerate}
\item
Any pair $\left(\left(\fmx_1,...,\fmx_n\right),\assfun\right)$ fulfilling \eqref{pathconduro} defines one and only one path in $\fremag^n$ by setting \eqref{pathconduro2}.  
\item
For any path $\pmone$ in $\fremag$ there could be more than one pair $\left(\left(\fmx_1,...,\fmx_n\right),\assfun\right)$ fulfilling \eqref{pathconduro} which defines $\pmone$ by setting \eqref{pathconduro2}.
\end{enumerate}
\end{remark}

In Propositions \ref{occpath2}-\ref{pathoneprop1} below we study how and to what extent fixed paths and corresponding skeletons can be modified in order to meet specific needs. We refer to \eqref{Bchoice}, \eqref{Bchoice2}.

\begin{proposition}\label{occpath2} Fix $n \in \mathbb{N}$, a path $\pmone$ in $\fremag^n$, a finite (possibly empty) subset $E \subseteq \bascont$.
\begin{enumerate}
\item For any skeleton $\left(\left(\fmx_1,...,\fmx_n\right),\assfun\right)$ of $\pmone$ there is an injective set function $\theta: \underset{\mathsf{n}=1}{\overset{n}{\coprod}}\,\occset^{\fmx_{\mathsf{n}}} \rightarrow \bascont$, a skeleton $\left(\left(\fmy_1,...,\fmy_n\right),\assfundue\right)$ of $\pmone$ such that
\begin{flalign}
&\left(\fmy_1,...,\fmy_n\right) = \strchar\left[\left(\fmx_1,...,\fmx_n\right),\incllim_{-1} \funcomp \theta\right]\text{;}\label{skelprop1}\\[8pt]
&\assfundue\left(\theta\left(f,i\right),t\right)=\assfun\left(f,t\right)\qquad \forall \left(f,i\right) \in \underset{\mathsf{n}=1}{\overset{n}{\coprod}}\,\occset^{\fmx_{\mathsf{n}}}	\qquad \forall t \in \left(-1,1\right) \text{;}\label{skelprop3}\\[8pt]
&\left(\underset{\mathsf{n}=1}{\overset{n}{\bigcup}}\occset_{\fmy_{\mathsf{n}}}\right) \cap E=\udenset\text{.}\label{skelprop2}
\end{flalign}
\item Fix $I\subseteq \left\{1,...,n\right\}$. For any skeleton $\left(\left(\fmx_1,...,\fmx_n\right),\assfun\right)$ of $\pmone$ with $\fmx_{\mathsf{n}}\in \fremagcont$ for any $\mathsf{n} \in I$ then there is an injective set function $\theta: \underset{\mathsf{n}=1}{\overset{n}{\coprod}}\,\occset^{\fmx_{\mathsf{n}}} \rightarrow \bascont$, a skeleton $\left(\left(\fmy_1,...,\fmy_n\right),\assfundue\right)$ of $\pmone$ fulfilling all conditions below:
\begin{equation*}
\eqref{skelprop1}\text{;}\hspace{40pt}\eqref{skelprop3}\text{;}\hspace{40pt}\eqref{skelprop2}\text{;}\hspace{40pt}\fmy_{\mathsf{n}}\in \fremagcont\qquad \forall \mathsf{n} \in I\text{.}
\end{equation*}
\end{enumerate}
\end{proposition}
\begin{proof}\mbox{}\newline 
\textnormal{\textbf{Proof of statement 1.}}\ \ We refer to Definition \ref{defpath}-[4], Proposition \ref{elpropC4}, Remark \ref{notelpropC4}-[2].\newline
We choose a set function $\theta: \underset{\mathsf{n}=1}{\overset{n}{\coprod}}\,\occset^{\fmx_{\mathsf{n}}} \rightarrow \bascont$ fulfilling the system of conditions below
\begin{equation*}
\left\{
\begin{array}{l}
\theta\;\text{is an injective set function,}\vspace{4pt}\\
\theta\left(f,i\right)\in \bascont\left(\Dom\left(f\right)\right)\setminus E\qquad\forall \left(f,i\right) \in \underset{\mathsf{n}=1}{\overset{n}{\coprod}}\,\occset^{\fmx_{\mathsf{n}}}\text{.}
\end{array}
\right.
\end{equation*}
This choice is possible since: $\bascont\left(\intuno\right)$ is an infinite set for any $\intuno \subseteqdentro \mathbb{R}^m$ and any $ m \in \mathbb{N}$; $\underset{\mathsf{n}=1}{\overset{n}{\coprod}}\,\occset^{\fmx_{\mathsf{n}}}$ is a finite set by construction of $\fremag$; $E$ is a finite set by assumption.\newline  
Then we define pair $\left(\left(\fmy_1,...,\fmy_n\right),\assfundue\right)$ by setting \eqref{skelprop1}, \eqref{skelprop3}. Eventually statement 1 follows by checking that pair $\left(\left(\fmy_1,...,\fmy_n\right),\assfundue\right)$ fulfills \eqref{pathconduro}, \eqref{pathconduro2}.\newline
\textnormal{\textbf{Proof of statement 2.}}\ \ The proof coincides word by word with the proof of statement 1. Eventually we check that if $\fmx_{\mathsf{n}} \in \fremagcont$ for any $\mathsf{n} \in I$ then $\fmy_{\mathsf{n}}\in \fremagcont$ for any $\mathsf{n}\in I$. 
\end{proof}

\begin{proposition}\label{occpath3} Fix $n \in \mathbb{N}$, a path  $\pmone_{\mathsf{n}}$ in $\fremag^{m_{\mathsf{n}}}$ for any $\mathsf{n}\in\left\{1,...,n\right\}$, a finite (possibly empty) subset $E \subseteq \bascont$.
\begin{enumerate}
\item For any $\mathsf{n} \in \left\{1,...,n\right\}$ there is a skeleton $\left(\left(\fmx_{\mathsf{n},1},...,\fmx_{\mathsf{n},m_{\mathsf{n}}}\right),\assfun_{\mathsf{n}}\right)$ of $\pmone_\mathsf{n}$ such that
\begin{equation*}
\left(\underset{\mathsf{m}=1}{\overset{m_i}{\bigcup}}\occset_{\fmx_{i,\mathsf{m}}}\right)\cap \left(E\;\cup\;\underset{\mathsf{m}=1}{\overset{m_j}{\bigcup}}\occset_{\fmx_{j,\mathsf{m}}}\right)=\udenset \qquad\forall i,j \in \left\{1,...,n\right\}\;\text{with}\;i\neq j\text{.}
\end{equation*}
\item Fix $I\subseteq \left\{1,...,n\right\}$, $J_{\mathsf{n}}\subseteq \left\{1,...,m_{\mathsf{n}}\right\}$ for any $\mathsf{n} \in I$. If there is a core $\left(\fmy_{\mathsf{n},1},...,\fmy_{\mathsf{n}, m_{\mathsf{n}}}\right)$ of $\pmone_{\mathsf{n}}$ with $\fmy_{\mathsf{n},\mathsf{m}}\in \fremagcont$ for $\mathsf{n} \in I$ and any $\mathsf{m} \in J_{\mathsf{n}}$ then we can choose $\left(\fmx_{\mathsf{n},1},...,\fmx_{\mathsf{n},m_{\mathsf{n}}}\right) \in \fremag^{m_{\mathsf{n}}}$ fulfilling $\fmx_{\mathsf{n},\mathsf{m}}\in \fremagcont$ for any $\mathsf{n} \in I$ and any $\mathsf{m} \in J_{\mathsf{n}}$.
\end{enumerate}
\end{proposition}
\begin{proof}
We argue by induction on $n$. If $n=1$ then statements 1, 2 follows directly by Proposition \ref{occpath2}. If $n=l+1$ with $l\geq 1$ then, by induction, statements 1, 2 hold true with data $l \in \mathbb{N}$, paths  $\pmone_{\mathsf{l}}$ in $\fremag^{m_{\mathsf{l}}}$ for $i\in\left\{1,...,l\right\}$. Eventually statements 1, 2 follow by applying Proposition \ref{occpath2} with data $m_{l+1}$, $\pmone_{l+1}$, $E\;\cup\;\overset{l}{\underset{\mathsf{l}=1}{\bigcup}}\left(\underset{\mathsf{m}=1}{\overset{m_{\mathsf{l}}}{\bigcup}}\occset_{\fmx_{\mathsf{l},\mathsf{m}}}\right)$.
\end{proof}

\begin{proposition}\label{pathoneprop1}
Fix $n_1,n_2 \in \mathbb{N}$, a path $\pmone_1$ in $\fremag^{n_1}$, a path $\pmone_2$ in $\fremag^{n_2}$. 
\begin{enumerate}
\item  The set function $\pmone=\left(\pmone_1,\pmone_2\right)\funcomp \Diag{\left(-1,1\right)}{2}$ is a path in $\fremag^{n_1+n_2}$.
\item If there are subsets $\singsmpth_1,\singsmpth_2 \subseteq \left(-1,1\right)$ such that pair $\left(\pmone_1,\singsmpth_1\right)$ is a smooth path in $\fremag^{n_1}$ and pair $\left(\pmone_2,\singsmpth_2\right)$ is a smooth path in $\fremag^{n_2}$ then pair\newline
$\left(\pmone,\singsmpth_1\cup\singsmpth_2 \right)$ is a smooth path in $\fremag^{n_1 + n_2}$.
\item Fix $I \subseteq \left\{1,2\right\}$, $J_i \subseteq \left\{1,...,n_i\right\}$ and a skeleton $\left(\left(\fmx_{i,1},...,\fmx_{i,n_i}\right), \assfun_i\right)$ of $\pmone_i$ for any $i \in I$, a skeleton $\left(\left(\fmy_1,...,\fmy_{n_1+n_2}\right), \assfun_{\fmy}\right)$ of $\pmone$.\newline
Assume $\hspace{10pt} \fmx_{i,\mathsf{n}} \in \fremagcont\quad \forall i\in I\quad \forall \mathsf{n}\in J_i$.\hspace{10pt}Set $n_0=0$.\newline
Then exists a skeleton $\left(\left(\fmx_1,...,\fmx_{n_1+n_2}\right), \assfun\right)$ of $\pmone$ such that: 
\begin{flalign*}
&\fmx_k=y_k \hspace{20pt}\text{if} \quad k  \neq n_{i-1}+\mathsf{n}\quad \forall i \in I \quad \forall \mathsf{n}\in J_i \text{;}\\[4pt]
&\fmx_k\in \fremagcont \hspace{18pt}\text{if}\quad \exists i \in I \quad \exists \mathsf{n}\in J_i\quad \text{such that}\quad k  = n_{i-1}+\mathsf{n}\text{.}
\end{flalign*}
\end{enumerate}
\end{proposition}
\begin{proof}\mbox{}\newline
\textnormal{\textbf{Proof of statements 1, 2.}}\newline
Choose a skeleton $\left(\left(\fmz_{1,1},...,\fmz_{1,n_1}\right), \assfuntre_1\right)$ of $\fremag_1$, a skeleton $\left(\left(\fmz_{2, 1},...,\fmz_{2,n_2}\right), \assfuntre_2\right)$ of $\fremag_2$ fulfilling 
$ \left(\underset{\mathsf{n}=1}{\overset{n_1}{\bigcup}} \,\occset_{\fmz_{1,\mathsf{n}}}\right)\cap\left(\underset{\mathsf{n}=1}{\overset{n_2}{\bigcup}} \,\occset_{\fmz_{2,\mathsf{n}}}\right)=\udenset$. This choice is made possible by Proposition \ref{occpath3}. We define:\newline
$\left(\fmx_1,...,\fmx_{n_1+n_2}\right)\in \fremag^{n_1+n_2}$ by setting
\begin{equation*}
\fmx_{\mathsf{n}} =\left\{
\begin{array}{ll}
\fmz_{1,\mathsf{n}}&\text{if} \; \mathsf{n} \leq n_1\text{,}\vspace{4pt}\\
\fmz_{2,\mathsf{n}-n_1}& \text{if}\;\mathsf{n}>n_1\text{;}
\end{array}
\right.
\end{equation*}
the set function $\assfun: \left(\overset{n_1+n_2}{\underset{\mathsf{n}=1}{\bigcup}} \occset_{\fmx_{\mathsf{n}}}\right)  \times \left(-1,1\right)\rightarrow \Cksp{0}$ by setting
\begin{equation*}
\assfun(f,t)=\left\{
\begin{array}{ll}
\assfun_1(f,t)&\text{if}\;f \in \overset{n_1}{\underset{\mathsf{n}=1}{\bigcup}} \occset_{\fmz_{1,\mathsf{n}}}\text{,}\vspace{4pt}\\
\assfun_2(f,t)&\text{if}\;f \in \overset{n_2}{\underset{\mathsf{n}=1}{\bigcup}} \occset_{\fmz_{2,\mathsf{n}}}\text{.}
\end{array}
\right.
\end{equation*}
Eventually statements 1, 2 follow by checking that pair $\left(\left(\fmx_1,...,\fmx_{n_1+n_2}\right),\assfun\right)$ fulfills \eqref{pathconduro} and \eqref{smptass} respectively, and defines the set function $\pmone$ through \eqref{pathconduro2}.\newline
\textnormal{\textbf{Proof of statements 3.}}\ \ For any $i \in I$ we denote by $\pmone_i$ the path in $\fremag^{n_i}$ defined  by the pair $\left(\left(\fmx_{i,1},...,\fmx_{i,n_i}\right), \assfun_i\right)$ (see Remark \ref{notunsk}). For any $i \notin I$ we choose a path $\pmone_i$ in $\fremag^{n_i}$.\newline
By Proposition \ref{occpath3} applied with data 2, paths $\pmone_1$, $\pmone_2$, $\underset{\mathsf{n}=1}{\overset{n_1+n_2}{\bigcup}}\occset_{\fmy_{\mathsf{n}}}$ we get a skeleton $\left(\left(\fmz_{i,1},...,\fmz_{i,n_i}\right), \assfuntre_i\right)$ of $\pmone_i$ for any $i \in \left\{1,2\right\}$ fulfilling
\begin{equation*}
\left(\underset{\mathsf{n}=1}{\overset{n_i}{\bigcup}}\occset_{\fmz_{i,\mathsf{n}}}\right)\cap \left(\underset{\mathsf{n}=1}{\overset{n_1+n_2}{\bigcup}}\occset_{\fmy_{\mathsf{n}}}\;\cup\;\underset{\mathsf{n}=1}{\overset{n_k}{\bigcup}}\occset_{\fmz_{k,\mathsf{n}}}\right)=\udenset \qquad\forall i,k \in \left\{1,2\right\}\text{.}
\end{equation*}
Then we define:\newline
$\left(\fmx_1,...,\fmx_{n_1+n_2}\right)\in \fremag^{n_1+n_2}$ by setting 
\begin{equation*}
\fmx_k=\left\{
\begin{array}{ll}
\fmy_k &\text{if} \quad k  \neq n_{i-1}+\mathsf{n}\quad \forall i \in I \quad \forall \mathsf{n}\in J_i \text{,}\vspace{4pt}\\
\fmz_{i,\mathsf{n}}\in \fremagcont &\text{if}\quad \exists i \in I \quad \exists \mathsf{n}\in J_i\quad \text{such that}\quad k  = n_{i-1}+\mathsf{n}\text{;}
\end{array}
\right.
\end{equation*}
the set function $\assfun:\left(\underset{\mathsf{n}=1}{\overset{n_1+n_2}{\bigcup}}\occset_{\fmx_{\mathsf{n}}}\right)\times\left(-1,1\right)\rightarrow \Cksp{0}$ by setting
\begin{equation*}
\assfun\left(f,t\right)=\left\{
\begin{array}{ll}
\assfun_{\fmy}\left(f,t\right) &\forall\left(f,t\right)\in \left(\underset{\mathsf{n}=1}{\overset{n_1+n_2}{\bigcup}}\occset_{\fmy_{\mathsf{n}}}\right)\times\left(-1,1\right) \text{,}\vspace{4pt}\\
\assfuntre_i\left(f,t\right) &\forall i \in\left\{1,2\right\}\quad \forall \left(f,t\right)\in\left( \underset{\mathsf{n}=1}{\overset{n_i}{\bigcup}}\occset_{\fmz_{i,\mathsf{n}}}   \right)\times\left(-1,1\right)\text{.}
\end{array}
\right.
\end{equation*}
Eventually statement follows by checking axioms and exploiting Proposition \ref{occpath3}-[2].  
\end{proof}

In Proposition \ref{domcodpathone} below we prove that domain and co-domain are constant along a path in $\fremag$.   

\begin{proposition}\label{domcodpathone}
Fix a path $\pmone$ in $\fremag$, a skeleton $\left(\fmx,\assfun\right)$ of $\pmone$. Then
\begin{equation*}
\codmagone\left(\pmone\left(\unkdue\right)\right)=\codmagone\left(x\right)\text{,}\qquad \dommagone\left(\pmone\left(\unkdue\right)\right)=\dommagone\left(x\right)\qquad \forall \unkdue \in \left(-1,1\right)\text{.}
\end{equation*}
\end{proposition}
\begin{proof}
Statement follows by \eqref{pathconduro}-[$(iii)$, $(iv)$].
\end{proof}

In Proposition \ref{indassfun} below we prove that associating functions for a set $\occset_{\fmx}$ induce an associating function for any subset $\occset_{\fmy}$ of $\occset_{\fmx}$. 

\begin{proposition}
\label{indassfun} 
Fix $\fmx,\fmy\in \fremag$, an associating function $\assfun$ for $\occset_{\fmx}$.\newline
Assume $\occset_{\fmy}\subseteq \occset_{\fmx}$.
\begin{enumerate}
\item Set function $\assfun\funcomp\left(\incl{\occset_{\fmy}}{\occset_{\fmx}}, \idobj+\left(-1,1\right)+\right)$
is an associating function for $\occset_{\fmy}$.\newline
If $\assfun$ is a centered associating function for $\occset_{\fmx}$ then $\assfun\funcomp\left(\incl{\occset_{\fmy}}{\occset_{\fmx}}, \idobj+\left(-1,1\right)+\right)$ is a centered associating function for $\occset_{\fmy}$. 
\item If there is a subset $\singsmpth \subseteq \left(-1,1\right)$ such that data $1$, $\assfun$, $\singsmpth$ fulfill \eqref{smptass} then data $1$, $\assfun\funcomp\left(\incl{\occset_{\fmy}}{\occset_{\fmx}}, \idobj+\left(-1,1\right)+\right)$, $\singsmpth$ fulfill \eqref{smptass} as well.
\end{enumerate}
\end{proposition}
\begin{proof}
Statement is a straightforward consequence of Definition \ref{pathinC}. 
\end{proof}

Motivated by Propositions \ref{elpropC4}, \ref{indassfun} we introduce the notion of induced path.

\begin{definition}\label{subpath}\mbox{}
\begin{enumerate}
\item  Fix a path $\pmone$ in $\fremag$, a skeleton $\left(\fmx, \assfun\right)$ of $\pmone$, $\fmy\in \fremag$. Assume $\occset_{\fmy}\subseteq \occset_{\fmx}$.
We denote by $\pmone\left[\fmy, \left(\fmx, \assfun\right)\right]$ the path in $\fremag$ defined by $\left(\fmy, \assfun\funcomp\left(\incl{\occset_{\fmy}}{\occset_{\fmx}}, \idobj+\left(-1,1\right)+\right)\right)$ through \eqref{pathconduro2}. We say that $\pmone\left[\fmy, \left(\fmx, \assfun\right)\right]$ is the path induced on $\fmy$ by $\left(\fmx, \assfun\right)$. By Remark \ref{notunsk}-[2] we say that $\pmone\left[\fmy, \left(\fmx, \assfun\right)\right]$ is a path induced on $\fmy$ by $\pmone$.\newline
If $\fmy$ occurs in $\fmx$ at some $i \in \mathbb{N}$ then we set $\pmone\left[\left(\fmy,i\right), \left(\fmx, \assfun\right)\right]=\pmone\left[\fmy, \left(\fmx, \assfun\right)\right]$ and we say that $\pmone\left[\left(\fmy,i\right), \left(\fmx, \assfun\right)\right]$ is the path induced on $\fmy$ by $\left(\fmx, \assfun\right)$ at $i$. By Remark \ref{notunsk}-[2] we say that $\pmone\left[\left(\fmy,i\right), \left(\fmx, \assfun\right)\right]$ is a path induced on $\fmy$ by $\pmone$ at $i$. We drop any reference to the occurrence $i$ whenever no confusion is possible.
\item  Fix a smooth path $\left(\pmone, \singsmpth\right)$ in $\fremag$, a skeleton $\left(\fmx, \assfun\right)$ of $\left(\pmone, \singsmpth\right)$, $\fmy\in \fremag$. Assume $\occset_{\fmy}\subseteq \occset_{\fmx}$.
We denote by $\left(\pmone\left[\fmy, \left(\fmx, \assfun\right)\right],\singsmpth\right)$ the path in $\fremag$ defined by $\left(\fmy, \assfun\funcomp\left(\incl{\occset_{\fmy}}{\occset_{\fmx}}, \idobj+\left(-1,1\right)+\right)\right)$ through \eqref{pathconduro2}. We say that $\left(\pmone\left[\fmy, \left(\fmx, \assfun\right)\right], \singsmpth\right)$ is the smooth path induced on $\fmy$ by $\left(\fmx, \assfun\right)$. By Remark \ref{notunsk}-[2] we say that $\left(\pmone\left[\fmy, \left(\fmx, \assfun\right)\right], \singsmpth\right)$ is a smooth path induced on $\fmy$ by $\left(\pmone, \singsmpth\right)$.\newline
If $\fmy$ occurs in $\fmx$ at some $i \in \mathbb{N}$ then we set 
\begin{equation*}
\left(\pmone\left[\left(\fmy,i\right), \left(\fmx, \assfun\right)\right], \singsmpth\right)=\left(\pmone\left[\fmy, \left(\fmx, \assfun\right)\right], \singsmpth\right)
\end{equation*}
and we say that $\left(\pmone\left[\left(\fmy,i\right), \left(\fmx, \assfun\right)\right], \singsmpth\right)$ is the path induced on $\fmy$ by $\left(\fmx, \assfun\right)$ at $i$.
By Remark \ref{notunsk}-[2] we say that $\left(\pmone\left[\left(\fmy,i\right), \left(\fmx, \assfun\right)\right], \singsmpth\right)$ is a path induced on $\fmy$ by $\left(\pmone, \singsmpth\right)$ at $i$. We drop any reference to the occurrence $i$ whenever no confusion is possible.
\end{enumerate}
\end{definition}

In Proposition \ref{indindassfun} below we study the behavior of iterating the notion of induced paths. 
\begin{proposition}\mbox{}
\label{indindassfun} 
\begin{enumerate}
\item Fix a path $\pmone$ in $\fremag$, a skeleton $\left(\fmx, \assfun\right)$ of $\pmone$, an occurrence pair $\left(\fmy,i\right)$ for $\fmx$, an occurrence pair $\left(\fmz,j\right)$ for $\fmy$. \newline
Then $\pmone\left[\left(\fmz,j+i-1\right), \left(\fmx, \assfun\right)\right]=\pmone\left[\left(\fmz,j\right), \left(\fmy, \assfun\funcomp\left(\incl{\occset_{\fmy}}{\occset_{\fmx}}, \idobj+\left(-1,1\right)+\right)\right)\right]$.
\item Fix a smooth path $\left(\pmone, \singsmpth\right)$ in $\fremag$, a skeleton $\left(\fmx, \assfun\right)$ of $\left(\pmone, \singsmpth\right)$, an occurrence pair $\left(\fmy,i\right)$ for $\fmx$, an occurrence pair $\left(\fmz,j\right)$ for $\fmy$. Then
\begin{multline*}
\left(\pmone\left[\left(\fmz,j+i-1\right), \left(\fmx, \assfun\right)\right], \singsmpth\right)=\\
\left(\pmone\left[\left(\fmz,j\right), \left(\fmy, \assfun\funcomp\left(\incl{\occset_{\fmy}}{\occset_{\fmx}}, \idobj+\left(-1,1\right)+\right)\right)\right], \singsmpth\right)\text{.}
\end{multline*}
\item Fix $i \in \mathbb{N}$, $\fmx,\fmy \in \fremag$, a path $\pmone$ in $\fremag$ detecting $\fmx$, a detecting skeleton $\left(\fmx,\assfun\right)$ of $\pmone$. Assume that $\left(\fmy,i\right)$ is an occurrence pair for $\fmx$. Then:\newline
$\pmone\left[\left(\fmy,i\right), \left(\fmx, \assfun\right)\right]$ is a path in $\fremag$ detecting $\fmy$;\newline
$\left(\fmy, \assfun\funcomp\left(\incl{\occset_{\fmy}}{\occset_{\fmx}}, \idobj+\left(-1,1\right)+\right)\right)$ is a detecting skeleton of $\pmone\left[\left(\fmy,i\right), \left(\fmx, \assfun\right)\right]$. 
\end{enumerate} 
\end{proposition}
\begin{proof}
Statement is a straightforward consequence of Proposition \ref{indassfun}. 
\end{proof}

In Proposition \ref{incastrpath} we study how paths can be glued together to get a new path. We refer to Construction \ref{subevalsm}, Remark \ref{notelpropC4}-[3]. 

\begin{proposition}\label{incastrpath}
Fix $\fmx \in \fremag$, an occurrence set $\occset$ for $\fmx$, a path $\pmone\left[\fmy,i\right]$ in $\fremag$ through $\fmy$ for any $\left(\fmy,i\right) \in \occset$.
\begin{enumerate}
\item The set function  $\occfundue_t:\occset\rightarrow \fremag$ defined by 
\begin{equation*}
\occfundue_{\unkdue}\left(\fmy,i\right)= \pmone\left[\fmy,i\right]\left(\unkdue\right) \quad \forall \left(\fmy,i\right)\in \occset\quad \forall \unkdue \in\left(-1,1\right)
\end{equation*}
is an occurrence function for any $\unkdue \in\left(-1,1\right)$.\newline
The set function $\pmone\left[\fmx,\left\{\pmone\left[\fmy,i\right]\right\}_{\left(\fmy,i\right)\in \occset}\right]:\left(-1,1\right) \rightarrow \fremag$ defined by setting 
\begin{equation*}
\pmone\left[\fmx,\left\{\pmone\left[\fmy,i\right]\right\}_{\left(\fmy,i\right)\in \occset}\right]\left(\unkdue\right)=\strchar\left[\fmx,\occfundue_{\unkdue}\right]\qquad \forall \unkdue \in \left(-1,1\right)
\end{equation*}
is a path in $\fremag$ through $\fmx$.
\item Fix a skeleton $\left(\widehat{\fmy}\left[\fmy,i\right], \assfun\left[\fmy,i\right]\right)$ of $\pmone\left[\fmy,i\right]$ for any $\left(\fmy,i\right) \in \occset$. Assume:
\begin{flalign}
&\left\{\begin{array}{l}
\text{for any}\; \left(\incllim_{-1}\left(f\right),j	\right)\in \occset^{\fmx} \;\text{with}\; f\in 
\underset{\left(\fmy,i\right)\in \occset}{\bigcup}\,
\occset_{\widehat{\fmy}\left[\fmy,i\right]}\\[4pt]
\text{there is}\; \left(\fmy,i	\right)\in \occset\;\text{such that}\; \left(\incllim_{-1}\left(f\right),j-i+1	\right)\in \occset^{\widehat{\fmy}\left[\fmy,i\right]}\text{;}
\end{array}\label{eqcore10}
\right.\\[8pt]
&
\left\{\begin{array}{l}
\assfun\left[\fmy_1,i_1\right]\left(f,\unkdue\right)=\assfun\left[\fmy_2,i_2\right]\left(f,\unkdue\right)\\[4pt]
 \hspace{93pt}\forall \left(f,\unkdue\right)\in \left(\occset_{\widehat{\fmy}\left[\fmy_1,i_1\right]}\cap \occset_{\widehat{\fmy}\left[\fmy_2,i_2\right]}\right)\times\left(-1,1\right)\\[4pt] 
\hspace{93pt}\forall \left(\fmy_1,i_1\right), \left(\fmy_2,i_2\right)\in \occset \text{.}\label{eqcore20}
\end{array}
\right.
\end{flalign}
The set function $\occfundue:\occset\rightarrow \fremag$ defined by setting
\begin{equation*}
\occfundue\left(\fmy,i\right)=\widehat{\fmy}\left[\fmy,i\right] \qquad \forall \left(\fmy,i\right) \in \occset
\end{equation*}
is an occurrence function on $\occset$.\newline
Define:\newline
$\widehat{\fmx}\in \fremag$ by setting $\widehat{\fmx} =\strchar\left[\fmx, \occfundue \right] $;\newline
the set function $\assfun:\occset_{\widehat{\fmx}}\times\left(-1,1\right)\rightarrow \Cksp{0}$ by setting
\begin{equation*}
\assfun\left( f  ,\unkdue\right)=\left\{
\begin{array}{l}
\left(\assfun\left[\fmy,i\right]\right)\left(f  ,\unkdue\right)\hspace{16pt}\forall \left(\fmy,i\right)\in \occset\quad \forall \left( f  ,\unkdue\right)\in \occset_{\widehat{\fmy}\left[\fmy,i\right]}\times\left(-1,1\right)\text{,}\\[4pt]
f \hspace{49pt}\text{if}\;\left( f  ,\unkdue\right) \in\left(\occset_{\widehat{\fmx}}\setminus\left(\underset{\left(\fmy,i\right)\in\occset}{\bigcup}\, \occset_{\widehat{\fmy}\left[\fmy,i\right]}\right)\right)\times\left(-1,1\right)\text{.}
\end{array}
\right.
\end{equation*}
Then pair $\left(\widehat{\fmx}, \assfun \right)$ is a skeleton of $\pmone\left[\fmx,\left\{\pmone\left[\fmy,i\right]\right\}_{\left(\fmy,i\right)\in \occset}\right]$.\newline
If $\assfun\left[\fmy,i\right]$ is a centered associating function for $\occset_{\widehat{\fmy}\left[\fmy,i\right]}$ for any $\left(\fmy,i\right) \in \occset$ then $\assfun$ is a centered associating function for $\occset_{\widehat{\fmx}}$.
\item If $\fmx\in \fremagcont$ and for any $\left(\fmy,i\right) \in \occset$ there is a core $\widehat{\fmy}\left[\fmy,i\right]\in \fremagcont$ of $\pmone\left[\fmy,i\right]$ then there is a core $\widehat{\fmx}\in \fremagcont$ of $\pmone\left[\fmx,\left\{\pmone\left[\fmy,i\right]\right\}_{\left(\fmy,i\right)\in \occset}\right]$.  
\end{enumerate}
\end{proposition}
\begin{proof}\mbox{}\newline
\textnormal{\textbf{Proof of statement 1.}}\newline
By Propositions \ref{occpath2} for any $\left(\fmy,i\right) \in \occset$ there is a skeleton $\left(\widehat{\fmy}\left[\fmy,i\right],\assfun\left[\fmy,i\right]\right)$ of $\pmone\left[\fmy,i\right]$ such that conditions below are all fulfilled:
\begin{flalign}
&\begin{array}{l}
\occset_{\widehat{\fmy}\left[\fmy,i\right]}\cap \occset_{\fmx}=\udenset\hspace{15pt}\forall 
\left(\fmy,i\right)\in \occset\text{;}
\end{array}\label{eqcore12}\\[8pt]
&\left\{\begin{array}{l}
\occset_{\widehat{y}\left[\fmy_1,i_1\right]}\cap \occset_{\widehat{\fmy}\left[\fmy_2,i_2\right]}=\udenset\\[4pt] 
\hspace{60pt}\forall \left(\fmy_1,i_1\right), \left(\fmy_2,i_2\right)\in \occset\;\text{with}\; \left(\fmy_1,i_1\right)\neq\left(\fmy_2,i_2\right) \text{.}\label{eqcore22}
\end{array}\right.
\end{flalign}
Then Propositions \ref{invsost}, \ref{propfinordseq}-[1] entail that the set function $\occfundue:\occset\rightarrow \fremag$ defined by setting
\begin{equation*}
\occfundue\left(\fmy,i\right)=\widehat{\fmy}\left[\fmy,i\right] \qquad \forall \left(\fmy,i\right) \in \occset
\end{equation*}
is an occurrence function on $\occset$.\newline
Define $\widehat{\fmx}\in \fremag$ by setting $\widehat{\fmx} =\strchar\left[\fmx, \occfundue \right]$.
Then Proposition \ref{elpropC4}, \eqref{eqcore12}, \eqref{eqcore22} entail that
\begin{equation}
\occset_{\widehat{\fmx}}=\left(\occset_{\fmx}\setminus\left(\underset{\left(\fmy,i\right)\in\occset}{\bigcup}\, \occset_{\fmy}\right)\right)\coprod
\left(\underset{\left(\fmy,i\right) \in \occset}{\coprod}\,\occset_{\widehat{\fmy}\left[\fmy,i\right]}\right)\text{.}\label{genfcore}
\end{equation} 
Define the set function $\assfun:\occset_{\widehat{\fmx}}\times\left(-1,1\right)\rightarrow \Cksp{0}$ by setting
\begin{equation*}
\assfun\left( f  ,\unkdue\right)=\left\{
\begin{array}{l}
\left(\assfun\left[\fmy,i\right]\right)\left(f  ,\unkdue\right)\hspace{16pt}\forall \left(\fmy,i\right)\in \occset\quad \forall \left( f  ,\unkdue\right)\in \occset_{\widehat{\fmy}\left[\fmy,i\right]}\times\left(-1,1\right)\text{,}\\[4pt]
f \hspace{56pt}\text{if}\;\left( f  ,\unkdue\right) \in\left(\occset_{\widehat{\fmx}}\setminus\left(\underset{\left(\fmy,i\right)\in\occset}{\bigcup}\, \occset_{\fmy}\right)\right)\times\left(-1,1\right)\text{.}
\end{array}
\right.
\end{equation*} 
Eventually statement 1 follows by checking that pair $\left(\widehat{\fmx}, \assfun\right)$ defines the set function $\pmone\left[\fmx,\left\{\pmone\left[\fmy,i\right]\right\}_{\left(\fmy,i\right)\in \occset}\right]$ through \eqref{pathconduro2}.\newline
\textnormal{\textbf{Proof of statement 2.}}\newline
Propositions \ref{invsost}, \ref{propfinordseq}-[1] entail that the set function $\occfundue$ is an occurrence function on $\occset$
Proposition \ref{elpropC4}, \eqref{eqcore10}, \eqref{eqcore20} entaail that
\begin{equation}
\occset_{\widehat{\fmx}}=\occset_{\fmx}\cup\left(\underset{\left(\fmy,i\right)\in\occset}{\bigcup}\, \occset_{\widehat{\fmy}\left[\fmy,i\right]}\right)\text{.}\label{genfcoredue}
\end{equation}
Since $\assfun\left[\fmy,i\right]$ is an associating function for $\occset_{\fmy}$ for any $\left(\fmy,i\right)\in \occset$ and \eqref{genfcoredue} holds true, we have that set function $\assfun$ is an associating function for $\occset_{\widehat{\fmx}}$.\newline
Eventually statements 2 follows by exploiting \eqref{evcmppropr}-[$(v)$], \eqref{proplininclone}-[$(i)$] to check that pair $\left(\widehat{\fmx}, \assfun \right)$
defines $\pmone\left[\fmx,\left\{\pmone\left[\fmy,i\right]\right\}_{\left(\fmy,i\right)\in \occset}\right]$ through \eqref{pathconduro2}.\newline
\textnormal{\textbf{Proof of statement 3.}}\newline
Statement 3 follows by checking that $\widehat{\fmx}$ defined in statements 1 and 2 belongs to $\fremagcont$.
\end{proof}

In Proposition \ref{incastrsmpath} below we study how smooth paths can be glued together to get a new smooth path. We refer to  Construction \ref{subevalsm}, Proposition \ref{incastrpath}, Remark \ref{notelpropC4}-[3]. 

\begin{proposition}\label{incastrsmpath}
Fix $\fmx \in \fremag$, an occurrence set $\occset$ for $\fmx$, a smooth path $\left(\pmone\left[\fmy,i\right], \singsmpth\left[\fmy,i\right]\right)$ in $\fremag$ through $\fmy$ for any $\left(\fmy,i\right) \in \occset$. Assume 
\begin{equation*}
\forall \left(\fmz,j	\right)\in \occset^{\fmx} \quad \exists \left(\fmy,i	\right)\in \occset\,:\; \left(\fmz,j-i+1	\right)\in \occset^{\fmy}\text{.}
\end{equation*}
\begin{enumerate}
\item The pair 
\begin{equation*}
\left(\pmone\left[\fmx, \left\{\pmone\left[\fmy,i\right] \right\}_{\left(\fmy,i\right)\in \occset}\right], \underset{\left(\fmy,i\right)\in \occset}{\bigcup}\,\singsmpth\left[\fmy,i\right] \right)
\end{equation*}
is a smooth path in $\fremag$ through $\fmx$. 
\item Fix a skeleton $\left(\widehat{\fmy}\left[	\fmy,i\right], \assfun\left[\fmy,i\right]\right)$ of $\left(\pmone\left[\fmy,i\right], \singsmpth\left[	\fmy,i\right]\right)$ for any $\left(\fmy,i\right) \in \occset$. Assume
\begin{equation*}
\left\{
\begin{array}{l}
\assfun\left[\fmy_1,i_1\right]\left(f,\unkdue\right)=\assfun\left[\fmy_2,i_2\right]\left(f,\unkdue\right)\\[4pt]
\hspace{93pt} \forall \left(f,\unkdue\right)\in \left(\occset_{\widehat{\fmy}\left[\fmy_1,i_1\right]}\cap \occset_{\widehat{\fmy}\left[\fmy_2,i_2\right]}\right)\times\left(-1,1\right)\\[4pt]
 \hspace{93pt}\forall \left(\fmy_1,i_1\right), \left(\fmy_2,i_2\right)\in \occset \text{.}
\end{array}
\right.
\end{equation*}
The set function $\occfundue:\occset\rightarrow \fremag$ defined by setting
\begin{equation*}
\occfundue\left(\fmy,i\right)=\widehat{\fmy}\left[\fmy,i\right] \qquad \forall \left(\fmy,i\right) \in \occset
\end{equation*}
is an occurrence function on $\occset$.\newline
Define:\newline
$\widehat{\fmx}\in \fremag$ by setting $\widehat{\fmx} =\strchar\left[\fmx, \occfundue \right] $;\newline
the set function $\assfun:\occset_{\widehat{\fmx}}\times\left(-1,1\right)\rightarrow \Cksp{0}$ by setting
\begin{equation*}
\assfun\left( \fmz  ,\unkdue\right)=
\assfun\left[\fmy,i\right]\left(\fmz  ,\unkdue\right)\qquad\forall \left(\fmy,i\right)\in \occset\quad \forall \left( \fmz  ,\unkdue\right)\in \occset_{\widehat{\fmy}\left[\fmy,i\right]}\times\left(-1,1\right)\text{.}
\end{equation*}
Then pair $\left(\widehat{\fmx}, \assfun \right)$ is a skeleton of $\left(\pmone\left[\fmx, \left\{\pmone\left[\fmy,i\right] \right\}_{\left(\fmy,i\right)\in \occset}\right], \underset{\left(\fmy,i\right)\in \occset}{\bigcup}\,\singsmpth\left[\fmy,i\right] \right)$.
\end{enumerate}
\end{proposition}
\begin{proof} Statements follow by corresponding statements in Proposition \ref{incastrpath} by checking that \eqref{pathconduro2} is fulfilled.
\end{proof}

In Proposition \ref{incastrindpath} below we study the link between induced paths and glued paths. 

\begin{proposition}\label{incastrindpath}
Fix $\fmx \in \fremag$, an associating function $\assfun$ for $\occset_{\fmx}$, an occurrence pair $\left(\fmy, i\right)$ for $\fmx$, an occurrence set $\occset$ for $\fmy$.\newline
Assume that\hspace{20pt}  
$\forall \left(\fmw,l	\right)\in \occset^{\fmy} \quad \exists \left(\fmz,j	\right)\in \occset\,:\; \left(\fmw,l+i-j	\right)\in \occset^{\fmz}$.\newline
Then $\pmone\left[\fmy, \left\{\pmone\left[\left(\fmz,j+i-1\right), \left(\fmx,\assfun\right)\right] \right\}_{\left(\fmz,j\right)\in \occset}\right] =\pmone\left[\left(\fmy,i\right),\left(\fmx,\assfun\right)\right]$.
\end{proposition}
\begin{proof} Statement follows by exploiting Propositions \ref{indindassfun}, \ref{incastrpath}.
\end{proof}

In Propositions \ref{pathoneprop2}-\ref{pathoneprop5} below we specialize general situations described in Propositions \ref{incastrpath}, \ref{incastrsmpath} to the particular case of combining paths in $\fremag$ through operations $\compone$, $\lboundone\,\rboundone$, $\acmone$ to get new paths in $\fremag$. Propositions \ref{pathoneprop2}-\ref{pathoneprop5} could be proved by referring to Propositions \ref{incastrpath}, but it is easier to give explicit proofs.          

\begin{proposition}\label{pathoneprop2}
Fix two paths $\pmone_1, \pmone_2$ in $\fremag$.
\begin{enumerate}
\item  The set function $\pmone_1 \compone \pmone_2:\left(-1,1\right)\rightarrow\fremag$ defined by setting
\begin{equation*}
\left(\pmone_1 \compone \pmone_2\right)\left(\unkdue\right)=\pmone_1\left(\unkdue\right)\compone\pmone_2\left(t\right)\qquad \forall \unkdue\in \left(-1,1\right)
\end{equation*}
is a path in $\fremag$.
\item If there are subsets $\singsmpth_1,\singsmpth_2 \subseteq \left(-1,1\right)$ such that pairs $\left(\pmone_1,\singsmpth_1\right)$, $\left(\pmone_2,\singsmpth_2\right)$ are smooth paths in $\fremag$ then pair $\left(\pmone_1 \compone \pmone_2,\singsmpth_1\cup\singsmpth_2 \right)$ is a smooth path in $\fremag$.
\item If there is a skeleton $\left(\fmx_1, \assfun_1\right)$ of $\pmone_1$, a skeleton $\left(\fmx_2, \assfun_2\right)$ of $\pmone_2$ with $\fmx_1,\fmx_2 \in \fremagcont$ then there is a skeleton $\left(\fmx, \assfun\right)$ of $\pmone_1 \compone \pmone_2$ with $\fmx \in \fremagcont$.
\end{enumerate}
\end{proposition}
\begin{proof}\mbox{}\newline
\textnormal{\textbf{Proof of statements 1, 2.}}\newline
Choose a skeleton $\left(\fmz_1, \assfuntre_1\right)$ of $\pmone_1$, a skeleton $\left(\fmz_2, \assfuntre_2\right)$ of $\pmone_2$ fulfilling $\occset_{\fmz_1}\cap\occset_{\fmz_2}=\udenset$. This choice is made possible by Proposition \ref{occpath3}. We define:\newline
$\fmx\in \fremag$ by setting $\fmx=\fmz_1\compone \fmz_2$;\newline
the set function $\assfun: \occset_{\fmx}  \times \left(-1,1\right)\rightarrow \Cksp{0}$ by setting
\begin{equation*}
\assfun\left(f,\unkdue\right)=\assfuntre_i\left(f,\unkdue\right) \qquad  \forall i \in \left\{1,2\right\}\qquad\forall \left(f,\unkdue\right)  \in \occset_{\fmz_i}\times \left(-1,1\right)\text{.}
\end{equation*}
Eventually statements 1, 2 follow by checking that pair $\left(\fmx,\assfun\right)$ fulfills \eqref{pathconduro} and \eqref{smptass} respectively, and defines the set function $\pmone_1 \compone \pmone_2$ through \eqref{pathconduro2}.\newline
\textnormal{\textbf{Proof of statement 3.}}\ \ By Proposition \ref{occpath3} there is skeleton $\left(\fmz_1, \assfuntre_1\right)$ of $\pmone_1$, a skeleton $\left(\fmz_2, \assfuntre_2\right)$ of $\pmone_2$ fulfilling both conditions: $\occset_{\fmz_1}\cap\occset_{\fmz_2}=\udenset$; $\fmz_1, \fmz_2 \in \fremagcont$.
Then the proof coincides word by word with the proof of statements 1, 2. Eventually statement follows by checking that $\fmx \in \fremagcont$.
\end{proof}

\begin{proposition}\label{pathoneprop4} Fix $n \in \mathbb{N}$, an ordered $n$-tuple $\left(\pmone_1,...,\pmone_n\right)$ of paths in $\fremag$.
\begin{enumerate}
\item The set function $\lboundone\pmone_1,...,\pmone_n\rboundone:\left(-1,1\right)\rightarrow\fremag$ defined by setting
\begin{equation*}\lboundone\pmone_1,...,\pmone_n\rboundone\left(\unkdue\right)=\lboundone\pmone_1\left(\unkdue\right),...,\pmone_n\left(\unkdue\right)\rboundone\qquad \forall \unkdue\in \left(-1,1\right)
\end{equation*}
is a path in $\fremag$.
\item If there is a subset $\singsmpth_{\mathsf{n}} \subseteq \left(-1,1\right)$ such that the pair $\left(\pmone_{\mathsf{n}},\singsmpth_{\mathsf{n}}\right)$ is a smooth path in $\fremag$ for any $\mathsf{n} \in \left\{1,...,n\right\}$ then pair $\left(\lboundone\pmone_1,...,\pmone_n\rboundone,\underset{\mathsf{n}=1}{\overset{n}{\bigcup}}\singsmpth_{\mathsf{n}} \right)$ is a smooth path in $\fremag$.
\item If there is a skeleton $\left(\fmx_{\mathsf{n}}, \assfun_{\mathsf{n}}\right)$ of $\pmone_{\mathsf{n}}$ with $\fmx_{\mathsf{n}} \in \fremagcont$ for any $\mathsf{n} \in \left\{1,...,n\right\}$ then there is a skeleton $\left(\fmx, \assfun\right)$ of $\lboundone\pmone_1,...,\pmone_n\rboundone$ with $\fmx \in \fremagcont$.
\end{enumerate}
\end{proposition}
\begin{proof}\mbox{}\newline
\textnormal{\textbf{Proof of statements 1, 2.}}\ \ Choose a skeleton $\left(\fmz_{\mathsf{n}}, \assfuntre_{\mathsf{n}}\right)$ of $\pmone_{\mathsf{n}}$ for any $\mathsf{n}\in\left\{1,...,n\right\}$ fulfilling $\hspace{10pt} \occset_{\fmz_i}\cap\occset_{\fmz_j}=\udenset \qquad \forall i,j\in\left\{1,...,n\right\}$.\newline
This choice is made possible by Proposition \ref{occpath3}.\newline
We define (see Definition \ref{defpath}-[$4$], Proposition \ref{elpropC2}-[1]):\newline
$\fmx\in \fremag$ by setting $\fmx= \lboundone \fmz_1,...,\fmz_n\rboundone$;\newline 
the set function $\assfun: \occset_{\fmx}  \times \left(-1,1\right)\rightarrow \Cksp{0}$ by setting
\begin{equation*}
\assfun(f,t)=\left\{
\begin{array}{ll}
\assfuntre_1(f,\unkdue)&\text{if}\;\left(f,\unkdue\right) \in \occset_{\fmz_1}\times \left(-1,1\right)\text{,}\vspace{4pt}\\
\assfuntre_2(f,\unkdue)&\text{if}\;\left(f,\unkdue\right) \in \occset_{\fmz_2}\times \left(-1,1\right)\text{.}
\end{array}
\right.
\end{equation*}
Eventually statements 1, 2 follow by checking that pair $\left(\fmx,\assfun\right)$ fulfills \eqref{pathconduro} and \eqref{smptass} respectively, and defines the set function $\pmone$ through \eqref{pathconduro2}.\newline
\textnormal{\textbf{Proof of statements 3.}}\ \ The proof is word by word the same of statements 1, 2 by choosing  $\left(\fmz_{\mathsf{n}}, \assfuntre_{\mathsf{n}}\right)=\left(\fmx_{\mathsf{n}}, \assfun_{\mathsf{n}}\right)$ for any $\mathsf{n}\in\left\{1,...,n\right\}$ and exploiting Proposition \ref{occpath3}-[2].
\end{proof}

\begin{proposition}\label{pathoneprop5}
Fix a path $\pmone$ in $\fremag$, $\bsfm\in\bsfM$.
\begin{enumerate}
\item The set function $\bsfm\acmone\pmone:\left(-1,1\right)\rightarrow\fremag$ defined by setting
\begin{equation*}
\left(\bsfm\acmone\pmone\right)\left(\unkdue\right)=\bsfm\acmone\pmone\left(\unkdue\right)\qquad \forall \unkdue\in \left(-1,1\right)
\end{equation*}
is a path in $\fremag$.
\item If there is a subset $\singsmpth \subseteq \left(-1,1\right)$ such that the pair $\left(\pmone,\singsmpth\right)$ is a smooth path in $\fremag$ then pair $\left(\bsfm\acmone\pmone,\singsmpth \right)$ is a smooth path in $\fremag$.\newline
\item If $\bsfm\in\subbsfM$ and there is a skeleton $\left(\fmy, \assfun_{\fmy}\right)$ of $\pmone$ with $\fmy \in \fremagcont$ then there is a skeleton $\left(\fmx, \assfun\right)$ of $\bsfm\acmone\pmone$ with $\fmx \in \fremagcont$.
\end{enumerate}
\end{proposition}
\begin{proof} \mbox{}\newline
\textnormal{\textbf{Proof of statements 1, 2.}}\ \ Choose a skeleton $\left(\fmz, \assfuntre\right)$ of $\pmone$. Statement follows by checking that pair $\left(\bsfm \acmone \fmz, \assfuntre\right)$ fulfills \eqref{pathconduro} and \eqref{smptass} respectively, and defines the set function $\bsfm\acmone\pmone$ through \eqref{pathconduro2}.\newline
\textnormal{\textbf{Proof of statements 3.}}\ \ The proof is word by word the same of statements 1, 2 by choosing  $\left(\fmz, \assfuntre\right)=\left(\fmy, \assfun_{\fmy}\right)$.
\end{proof}

In Proposition \ref{funonsm} below we study the behavior of elements $\fmx \in \fremag$ considered as functions with $\Cksp{\infty}$-valued unknowns. We refer to Notations \ref{ins}-[2], \ref{difrealfunc}-[2], Remark \ref{notelpropC4}.

\begin{proposition}\label{funonsm}
Fix $\fmx \in \fremag$. 
Assume $\cardset\left(\occset_{\fmx}\right)\geq 1$. Fix an order on $\occset_{\fmx}$.\newline
Set:\hspace{15pt}
$k=\cardset\left(\occset_{\fmx}\right)\text{,} \hspace{15pt} \occset_{\fmx}=\left\{f_1,...,f_k\right\}\text{,}\hspace{15pt}
\Cksp{\infty}\left[\fmx\right]=\underset{\mathsf{k}=1}{\overset{k}{\prod}}\Cksp{\infty}(\Dom\left(f_{\mathsf{k}}\right))+\mathbb{R}+.$\newline
Define the set function $\occfun\left[g_1,...,g_k\right]:\occset_{\fmx}\rightarrow \Cksp{\infty}$ for any $\left(g_1,...,g_k\right)  \in \Cksp{\infty}\left[\fmx\right]$ by setting $\qquad\occfun\left[g_1,...,g_k\right]\left(f_{\mathsf{k}}\right)=g_{\mathsf{k}} \qquad \forall \mathsf{k} \in \left\{1,...,k\right\}$.\newline
Then:
\begin{enumerate}
\item there is a set function $\fremagfunonsm\left[\fmx\right] : \Cksp{\infty}\left[\fmx\right]\rightarrow \Cksp{\infty}$ defined by setting
\begin{equation*}
\fremagfunonsm\left[\fmx\right]\left(g_1,...,g_k\right)=\evalcompone\left( \strchar\left[\fmx,\lininclone \funcomp \occfun\left[g_1,...,g_k\right]\right] \right)\quad\forall \left(g_1,...,g_k\right)\in \Cksp{\infty}\left[\fmx\right] \text{;}
\end{equation*}
\item  the set function $\fremagfunonsm\left[\fmx\right]\funcomp \incl{\setsymcinque\left[\fmx\right]}{\Cksp{\infty}\left[\fmx\right]}$ is a continuous function where:
\begin{flalign*}
&\setsymcinque\left[\fmx\right]=\left\{\left(g_1,...,g_k\right)\in\Cksp{\infty}\left[\fmx\right]\;:\; \fremagfunonsm\left[\fmx\right]\left(g_1,...,g_k\right)\neq\smospempty\right\}\text{;}\\[4pt]
& \setsymcinque\left[\fmx\right]\;\text{is endowed with the subspace topology with respect to}\;\left(\Cksp{\infty}\right)^k\text{.}
\end{flalign*}
\end{enumerate} 
\end{proposition} 
\begin{proof}\mbox{}\newline
\textnormal{\textbf{Proof of statement 1.}}\ \ Statement follows directly by definition.\newline 
\textnormal{\textbf{Proof of statement 2.}}\ \  We set $L=\left\{\fmx \in \fremag: \fmx\;\text{fulfills statement 2} \right\}$, then we prove that $L=\fremag$ by a recursive argument.\newline
\emph{\textbf{Step 1.}}\ \ Inclusion $\fremaggen_0 \subseteq L$ straightforwardly follows by \eqref{extopsmdis2}.\newline
\emph{\textbf{Step 2.}}\ \ Fix $n \in \mathbb{N}$, $\left(\fmx_1,...,\fmx_n\right)\in L^n$. Proof that $\lboundone \fmx_1,...,\fmx_n\rboundone \in L$ is straightforward.\newline
\emph{\textbf{Step 3.}}\ \ Fix $\bsfm \in \bsfM$, $\fmx \in L$. Proof that $\bsfm\acmone \fmx \in L$ follows since $\smospbsfmu$ is a continuous function (see \eqref{extopsmdis3}).\newline
\emph{\textbf{Step 4.}}\ \ Fix $\fmy,\fmz \in L$. We prove that $\fmy\compone \fmz\in L$.\newline
Denote $\fmy\compone \fmz$ by $\fmx$. We emphasize that limits of smooth functions are computed in $\idssmo$ endowed with topology defined in Example \ref{ExCinf}. \newline 
Fix: an order on $\occset_{\fmx}$; a $k$-tuple $\left(g_1,...,g_k \right)\in \setsymcinque\left[\fmx\right]$; a sequence of $k$-tuples $\left\{\left({_{i}g_1},...,{_{i}g_k} \right)\right\}_{i\in \mathbb{N}}\subseteq \setsymcinque\left[\fmx\right]$.\newline
Set $\occset_{\fmx}=\left\{f_1,...,f_k\right\}$. Assume
\begin{equation}
\underset{i\rightarrow+\infty}{lim}\,{_{i}g_{\mathsf{k}}}=g_{\mathsf{k}}\; \text{in}\;\Cksp{\infty}(\Dom\left(f_{\mathsf{k}}\right))+\mathbb{R}++\qquad\forall \mathsf{k} \in \left\{1,...,k\right\}\text{.}
\label{asscomp}
\end{equation}
By Definition \ref{defpath}-[$4$], Proposition \ref{elpropC1}-[1] the order on $\occset_{\fmx}$ induces order on both $\occset_{\fmy}$ and $\occset_{\fmz}$.\newline
We set $k_{\fmy}=\cardset\left(\occset_{\fmy}\right)$. If $k_{\fmy} \geq 1$ then inclusion of ordered sets $\incl{\occset_{\fmy}}{\occset_{\fmx}}$ induces a set function $\incloccset_{\fmy} :\left\{1,...,k_{\fmy}\right\}\rightarrow\left\{1,...,k\right\}$ fulfilling
\begin{equation*}
\incl{\occset_{\fmy}}{\occset_{\fmx}}\left(f_{\mathsf{k}}\right)=f_{\incloccset_{\fmy}\left(\mathsf{k}\right)} \qquad \forall \mathsf{k} \in \left\{1,...,k_{\fmy}\right\}\text{.}
\end{equation*} 
We set $k_{\fmz}=\cardset\left(\occset_{\fmz}\right)$. If $k_{\fmz} \geq 1$ then inclusion of ordered sets $\incl{\occset_{\fmz}}{\occset_{\fmx}}$ induces a set function $\incloccset_{\fmz} :\left\{1,...,k_{\fmz}\right\}\rightarrow\left\{1,...,k\right\}$ fulfilling
\begin{equation*}
\incl{\occset_{\fmz}}{\occset_{\fmx}}\left(f_{\mathsf{k}}\right)=f_{\incloccset_{\fmz}\left(\mathsf{k}\right)} \qquad \forall \mathsf{k}\in \left\{1,...,k_{\fmz}\right\}\text{.}
\end{equation*} 
By \eqref{evcmppropr}-[$(i)$] we have
\begin{flalign*}
&\fremagfunonsm\left[\fmx\right]\left(g_1,...g_k\right)=\fremagfunonsm\left[\fmy\right]\left(g_{\incloccset_{\fmy}\left(1\right)},...g_{\incloccset_{\fmy}\left(k_{\fmy}\right)}\right)\funcomp \fremagfunonsm\left[\fmz\right]\left(g_{\incloccset_{\fmz}\left(1\right)},...,g_{\incloccset_{\fmz}\left(k_{\fmz}\right)}\right)\text{,}\\[8pt]
&\fremagfunonsm\left[\fmx\right]\left({_{i}g_1},...,{_{i}g_k}\right)=\\[4pt]
&\qquad\fremagfunonsm\left[\fmy\right]\left({_{i}g_{\incloccset_{\fmy}\left(1\right)}},...,{_{i}g_{\incloccset_{\fmy}\left(k_{\fmy}\right)}}\right)\funcomp\fremagfunonsm\left[\fmz\right]\left({_{i}g_{\incloccset_{\fmz}\left(1\right)}},...,{_{i}g_{\incloccset_{\fmz}\left(k_{\fmz}\right)}}\right)\qquad \forall i \in\mathbb{N}\text{,}
\end{flalign*}
hence
\begin{flalign}
&\fremagfunonsm\left[\fmy\right]\left(g_{\incloccset_{\fmy}\left(1\right)},...g_{\incloccset_{\fmy}\left(k_{\fmy}\right)}\right)\neq \smospempty\text{,}\label{asscomp21}\\[4pt]
&\fremagfunonsm\left[\fmy\right]\left({_{i}g_{\incloccset_{\fmy}\left(1\right)}},...,{_{i}g_{\incloccset_{\fmy}\left(k_{\fmy}\right)}}\right)\neq \smospempty\hspace{30pt} \forall i \in\mathbb{N}\text{,}\label{asscomp22}\\[4pt]
&\fremagfunonsm\left[\fmz\right]\left(g_{\incloccset_{\fmz}\left(1\right)},...,g_{\incloccset_{\fmz}\left(k_{\fmz}\right)}\right)\neq \smospempty\text{,}\label{asscomp23}\\[4pt]
&\fremagfunonsm\left[\fmz\right]\left({_{i}g_{\incloccset_{\fmz}\left(1\right)}},...,{_{i}g_{\incloccset_{\fmz}\left(k_{\fmz}\right)}}\right)\neq \smospempty\hspace{30pt} \forall i \in\mathbb{N}\text{.}\label{asscomp24}
\end{flalign}
Eventually statement is achieved since continuity of $\funcomp$ stated in \eqref{extopsmdis3}, \eqref{asscomp}, \eqref{asscomp21}-\eqref{asscomp24}, recursive assumption on $\fmy$, $\fmz$ entail that
\begin{flalign*}
&\underset{i\rightarrow+\infty}{lim}\,\fremagfunonsm\left[\fmx\right]\left({_{i}g_1},...,{_{i}g_k}\right)
=\\[4pt]
&\underset{i\rightarrow+\infty}{lim}\,\fremagfunonsm\left[\fmy\right]\left({_{i}g_{\incloccset_{\fmy}\left(1\right)}},...,{_{i}g_{\incloccset_{\fmy}\left(k_{\fmy}\right)}}\right)
\funcomp \underset{i\rightarrow+\infty}{lim}\,\fremagfunonsm\left[\fmz\right]\left({_{i}g_{\incloccset_{\fmz}\left(1\right)}},...,{_{i}g_{\incloccset_{\fmz}\left(k_{\fmz}\right)}}\right)
=\\[4pt] 
&\fremagfunonsm\left[\fmy\right]\left(g_{\incloccset_{\fmy}\left(1\right)},...g_{\incloccset_{\fmy}\left(k_{\fmy}\right)}\right)
\funcomp \fremagfunonsm\left[\fmz\right]\left(g_{\incloccset_{\fmz}\left(1\right)},...,g_{\incloccset_{\fmz}\left(k_{\fmz}\right)}\right)=\fremagfunonsm\left[\fmz\right]\left(g_{\incloccset_{\fmz}\left(1\right)},...,g_{\incloccset_{\fmz}\left(k_{\fmz}\right)}\right)\text{.}
\end{flalign*}
\end{proof}

In Proposition \ref{funoncont} below we study the behavior of elements $\fmx \in \fremagcont$ considered as functions with $\Cksp{0}$-valued unknowns. We refer to Notation \ref{ins}-[2], Remark \ref{notelpropC4}.

\begin{proposition}\label{funoncont}
Fix $\fmx \in \fremagcont$. 
Assume $\cardset\left(\occset_{\fmx}\right)\geq 1$. Fix an order on $\occset_{\fmx}$.\newline
Set:\hspace{20pt}$k=\cardset\left(\occset_ {\fmx}\right)\text{,} \qquad \occset_{\fmx}=\left\{f_1,...,f_k\right\}\text{,}\qquad
\Cksp{0}\left[\fmx\right]=\underset{\mathsf{k}=1}{\overset{k}{\prod}}\Cksp{0}\left(\Dom\left(f_{\mathsf{k}}\right),\mathbb{R}\right)$.\newline
Define the set function $\occfun\left[g_1,...,g_k\right]:\occset_{\fmx}\rightarrow \Cksp{0}$ for any $\left(g_1,...,g_k\right)  \in \Cksp{0}\left[\fmx\right]$ by setting $\qquad\occfun\left[g_1,...,g_k\right]\left(f_{\mathsf{k}}\right)=g_{\mathsf{k}} \qquad \forall \mathsf{k} \in \left\{1,...,k\right\}$.\newline
Then:
\begin{enumerate}
\item there is a set function $\fremagfunoncont\left[\fmx\right] : \Cksp{0}\left[\fmx\right]\rightarrow \Cksp{0}$ defined by setting
\begin{equation*}
\fremagfunoncont\left[\fmx\right]\left(g_1,...g_k\right)=\evalcompone\left( \strchar\left[\fmx,\lininclone \funcomp \occfun\left[g_1,...,g_k\right]\right] \right)\quad\forall \left(g_1,...,g_k\right)\in \Cksp{0}\left[\fmx\right] \text{;}
\end{equation*}
\item  the set function $\fremagfunoncont\left[\fmx\right]\funcomp \incl{\setsymcinque\left[\fmx\right]}{\Cksp{0}\left[\fmx\right]}$ is a continuous function where:
\begin{flalign*}
&\setsymcinque\left[\fmx\right]=\left\{\left(g_1,...,g_k\right)\in\Cksp{0}\left[\fmx\right]\;:\; \fremagfunoncont\left[\fmx\right]\left(g_1,...,g_k\right)\neq\contspempty\right\}\text{;}\\[4pt]
& \setsymcinque\left[\fmx\right]\;\text{is endowed with the subspace topology with respect to}\;\left(\Cksp{0}\right)^k\text{.}
\end{flalign*}
\end{enumerate} 
\end{proposition} 
\begin{proof}
The proof coincides word by word with proof of Proposition \ref{funonsm} by replacing $\bsfM$ by $\subbsfM$ and $\Cksp{\infty}$ by $\Cksp{0}$. 
\end{proof}

In Definition \ref{admrelc1} below we introduce a relation establishing an order among elements belonging to $\fremag$ based on the behavior of their detecting paths. We refer to Notation \ref{ins}-[5].

\begin{definition}\label{admrelc1}\mbox{}
\begin{enumerate}
\item Fix $\fmx_1,\fmx_2\in \fremag$. We say that $\fmx_1 \ordrelsymb \fmx_2$ if and only if conditions below are all fulfilled: 
\begin{flalign}
&\begin{array}{l}
\dommagone\left(\fmx_1\right)=\dommagone\left(\fmx_2\right)\text{;}
\end{array}\label{eqelc0}\\[8pt]
&\begin{array}{l}
\occset_{\fmx_1}=\occset_{\fmx_2}\text{;}\label{eqelc}
\end{array}\\[8pt]
&\left\{
\begin{array}{l}
\text{for any choice of a path} \;\pmone_i\;\text{in}\; \fremag\;\text{detecting}\; \fmx_i\text{, of a dete-}\\[4pt]
\text{cting skeleton}\;\left(\fmx_i, \assfun_i\right) \;\text{of}\;\pmone_i\;\text{for any}\;i \in \left\{1,2\right\}\; \text{fulfilling}\\[4pt]
(i)\hspace{8pt}
\assfun_1\left(f,\unkdue\right)=\assfun_2\left(f,\unkdue\right)\quad\forall \left(f,\unkdue\right)\in\occset_{\fmx_1}\times  \left(-1,1\right)\text{,}\\[4pt]
\text{we}\hspace{4pt}\text{have}\\[4pt]
(ii)\hspace{3pt}\left(\evalcompone\left(\pmone_1\left(\unkdue\right)\right)\right)\left(\unkuno\right)\neq\udenunk\; \Rightarrow \\[4pt]
\hspace{91pt}
\left(\evalcompone\left(\pmone_1\left(\unkdue\right)\right)\right)\left(\unkuno\right)=\left(\evalcompone\left(\pmone_2\left(\unkdue\right)\right)\right)\left(\unkuno\right)\\[4pt] 
\hspace{122pt}\forall \unkuno \in\dommagone\left(\fmx_1\right)\;\forall \unkdue \in\left(-1,1\right)\setminus \left\{0\right\}\text{.}
\end{array}
\right.\label{eqelc1}
\end{flalign}
\item Fix $\fmx_1,\fmx_2\in \fremag$. We say that $\fmx_1 \relsymbord \fmx_2$ if and only if  $\fmx_1 \ordrelsymb \fmx_2$ and  $\fmx_2 \ordrelsymb \fmx_1$.
\end{enumerate}
\end{definition}

In Proposition \ref{prpadrel1} below we study the nature of relation $\ordrelsymb$ and we examine its compatibility with operations $\compone$, $\lboundone\;\rboundone$, $\acmone$. We refer to Notation \ref{alg}-[2].

\begin{proposition}\label{prpadrel1}\mbox{}
\begin{enumerate}
\item $\ordrelsymb$ is a reflexive and transitive relation in $\fremag$.
\item $\ordrelsymb$ is compatible with operations $\compone$, $\lboundone\;\rboundone$, $\acmone$.
\item Fix $\fmx,\fmy,\fmz \in \fremag$. Then $\left(\fmx\compone \fmy\right)\compone \fmz \ordrelsymb \fmx\compone \left(\fmy\compone \fmz\right)$.
\end{enumerate}
\end{proposition}
\begin{proof}
Proof of statements 1, 2 follow directly by checking axioms and exploiting Proposition \ref{incastrindpath}.\newline
We prove statement 3. Set $\fmx_1=\left(\fmx\compone \fmy\right)\compone \fmz$, $\fmx_2=\fmx\compone \left(\fmy\compone \fmz\right)$. Then \eqref{eqelc0}, \eqref{eqelc} straightforwardly hold true. We prove \eqref{eqelc1} arguing as follows. Fix a path $\pmone_i$ in $\fremag$ detecting $\fmx_i$, a detecting skeleton $\left(\fmx_i, \assfun_i\right)$ of $\pmone_i$ for any $i \in \left\{1,2\right\}$ fulfilling $\assfun_1\left(f,\unkdue\right)=\assfun_2\left(f,\unkdue\right)\quad\forall \left(f,\unkdue\right)\in  \occset_{x_1}\times\left(-1,1\right)$.\newline
Propositions \ref{indassfun}, \ref{indindassfun}, \ref{incastrindpath}, \ref{pathoneprop2} entail that there are a path $\pmone_{\fmx}$ in $\fremag$ detecting $\fmx$, a path $\pmone_{\fmy}$ in $\fremag$ detecting $\fmy$, a path $\pmone_{\fmz}$ in $\fremag$ detecting $\fmz$ such that\newline
\centerline{$\pmone_1\left(\unkdue\right) = \left(\left(\pmone_{\fmx} \compone \pmone_{\fmy}\right)\compone \pmone_{\fmz}\right)\left(\unkdue\right)\hspace{20pt} \forall \unkdue \in \left(-1,1\right)$.}\newline
Eventually statement follows by exploiting Proposition \ref{pathoneprop2} to check that:
\begin{flalign*}
&\pmone_{\fmx} \compone \left(\pmone_{\fmy}\compone \pmone_{\fmz}\right)\;\text{is a path in}\; \fremag\;\text{detecting}\; \fmx \compone \left(\fmy\compone \fmz\right)\text{;}\\[8pt]
&\pmone_2\left(\unkdue\right)=\left(\pmone_{\fmx} \compone \left(\pmone_{\fmy}\compone \pmone_{\fmz}\right)\right)\left(\unkdue\right)\hspace{20pt} \forall \unkdue \in \left(-1, 1\right)\setminus \left\{0\right\}\text{.}
\end{flalign*}
\end{proof}

Motivated by Proposition \ref{prpadrel1}-[3] we introduce the notion of augmentation of elements belonging to $\fremag$. We refer to Definition \ref{defpath}.

\begin{definition}\label{augdef}\mbox{}
\begin{enumerate}
\item Fix $\fmx_1,\fmx_2 \in \fremag$. We say that $\fmx_2$ is an elementary augmentation of $\fmx_1$ if and only if there are $i \in \mathbb{N}$, $\fmy_1, \fmy_2, \fmy_3  \in \fremag$ such that both conditions below are fulfilled:
\begin{flalign*}
&\begin{array}{l}\left(\left(\fmy_1\compone\fmy_2\right)\compone\fmy_3,i\right) \;\text{occurs in}\; \fmx_1\text{;}\end{array}\\[8pt]
&\left\{\begin{array}{l}
\fmx_2=\strchar \left[\fmx_1, \occfun\right]\\[4pt]
\text{where} \;\occfun\; \text{is the occurrence function on}
\left\{\left(\left(\fmy_1\compone\fmy_2\right)\compone\fmy_3,i\right)\right\}\\[4pt]
\text{defined by setting} \;\occfun \left(\left(\fmy_1\compone\fmy_2\right)\compone\fmy_3,i\right)=\fmy_1\compone\left(\fmy_2\compone\fmy_3\right)\text{.}
\end{array}
\right.
\end{flalign*}
We say that: $\left(\left(\fmy_1\compone\fmy_2\right)\compone\fmy_3,i\right)$ is the occurrence pair relative to the elementary augmentation $\fmx_2$ of $\fmx_1$; $\left\{\left(\left(\fmy_1\compone\fmy_2\right)\compone\fmy_3,i\right)\right\}$ is the occurrence set relative to the elementary augmentation $\fmx_2$ of $\fmx_1$; $\occfun$ is the occurrence function relative to the elementary augmentation $\fmx_2$ of $\fmx_1$;
\item Fix $\fmx_1,\fmx_2 \in \fremag$. We say that $\fmx_2$ is an augmentation of $\fmx_1$ if and only if there are $n \in \mathbb{N}$, $\fmz_1, ..., \fmz_n \in \fremag$ such that all conditions below are fulfilled 
\begin{flalign}
&\begin{array}{l}\fmz_1=\fmx_1\text{;}\end{array}\label{aug1}\\[8pt]
&\begin{array}{l}\fmz_n=\fmx_2\text{;}\end{array}\label{aug2}\\[8pt]
&\left\{\begin{array}{l}
\fmz_{\mathsf{n}+1}\;\text{is an elementary augmentation of}\; \fmz_{\mathsf{n}}\;
\text{for any}\\[4pt]\mathsf{n} \in \left\{1,...,n-1\right\} \text{.}
\end{array}\label{aug3}\right.
\end{flalign}
We say that the ordered $n$-tuple $\left(\fmz_1,....,\fmz_n\right)$ is the elementary chain from $\fmx_1$ to $\fmx_2$. 
\item Fix $\fmx_1,\fmx_2 \in \fremag$. We say that $\fmx_2$ is a maximal augmentation of $\fmx_1$ if and only if 
$\fmx_2$ is an augmentation of $\fmx_1$ and there are no elementary augmentations of $\fmx_2$.
\item Fix $\fmx \in \fremag$. We say that $\fmx$ is a maximal element if and only if there is no elementary augmentation of $\fmx$. 
\end{enumerate}
\end{definition}

In Propositions \ref{augpropl1}-\ref{augprop} below we study the behavior of augmentation with respect to operations $\compone$, $\lboundone\;\rboundone$, $\acmone$ and to relation $\ordrelsymb$. 

\begin{proposition}\label{augpropl1} 
Fix $i, n \in\mathbb{N}$, $\fmx, \fmx_1,...,\fmx_n,\fmy_1,\fmy_2,\fmy_3, \fmz\in \fremag$. Assume that all conditions below are fulfilled:
$\fmx= \fmx_1\compone...\compone\fmx_n$ up to associativity; $\fmz$ is an elementary augmentation of $\fmx$; $\left(\left(\fmy_1\compone\fmy_2\right)\compone\fmy_3,i\right)$ is the occurrence pair relative to the elementary augmentation $\fmz$ of $\fmx$.\newline
Then only two mutually exclusive possibilities occur:
\begin{flalign}
&\left\{\begin{array}{l} \text{there are}\;\mathsf{n}_1,\mathsf{n}_2, \mathsf{n}_3, \mathsf{n}_4 \in \left\{1,...,n\right\}\;\text{such that}\;\\[4pt]
(i)\hspace{11pt}\mathsf{n}_1<\mathsf{n}_2< \mathsf{n}_3< \mathsf{n}_4\text{,}\\[4pt]
(ii)\hspace{8pt}\fmy_i=\fmx_{\mathsf{n}_i}\compone...\compone\fmx_{\mathsf{n}_{i+1}-1}\;\text{up to associativity}\quad\forall i \in \left\{1,2,3\right\}\text{,}\\[4pt]
 (iii)\hspace{4pt}\fmx=\fmx_1\compone...\compone\fmx_{\mathsf{n}_1-1}\compone \fmy\compone \fmx_{\mathsf{n}_4}\compone...\compone\fmx_{n}\;\text{up to associativity}\\[4pt]
\hspace{25pt}\text{where}\;\fmy=\left(\left(\fmy_1\compone\fmy_2\right)\compone\fmy_3\right)\text{;}
\end{array}\label{p1}
\right.\\[8pt]
&\left\{\begin{array}{l} \text{there are}\; l \in \mathbb{N}\text{,}\; \mathsf{n}\in\left\{1,...,n\right\}\;\text{such that}\; \left(\left(\fmy_1\compone\fmy_2\right)\compone\fmy_3,i-l\right)\\[4pt]
\text{occurs in}\;\fmx_{\mathsf{n}}\text{.}
\end{array}\label{p2}
\right.
\end{flalign}
\end{proposition}
\begin{proof}
Statement follows by iterated application of Proposition \ref{elpropC1}-[1].
\end{proof}

\begin{proposition}\label{augpropl2} 
Fix $i, n \in\mathbb{N}$, $\fmx_1,...,\fmx_n,\fmy_1,\fmy_2,\fmy_3, \fmz\in \fremag$.
Assume that both conditions below are fulfilled:
$\fmz$ is an elementary augmentation of $\lboundone\fmx_1,...,\fmx_n\rboundone$;
$\left(\left(\fmy_1\compone\fmy_2\right)\compone\fmy_3,i\right)$ is the occurrence pair relative to the elementary augmentation $\fmz$ of $\lboundone\fmx_1,...,\fmx_n\rboundone$.\newline
Then there are $l\in\mathbb{N}$, $\mathsf{n}\in \left\{1,...,n\right\}$ such that both conditions below are fulfilled:
\begin{equation*}
\begin{array}{l}
\left(\left(\fmy_1\compone\fmy_2\right)\compone\fmy_3,i-l\right)\;\text{occurs in}\;\fmx_{\mathsf{n}}\text{;}\\[8pt]
\fmz=\lboundone\fmx_1,...,\strchar \left[\fmx_{\mathsf{n}}, \occfun\right],...,\fmx_n\rboundone\text{,}\\
\text{where}\;\occfun\;\text{is the occurrence function on}\;\left\{\left(\left(\fmy_1\compone\fmy_2\right)\compone\fmy_3,i-l\right)\right\}\\
\text{defined by setting}\;\occfun \left(\left(\fmy_1\compone\fmy_2\right)\compone\fmy_3,i-l\right)=\fmy_1\compone\left(\fmy_2\compone\fmy_3\right)\text{.}
\end{array}
\end{equation*}
\end{proposition}
\begin{proof}
Statement follows by Proposition \ref{elpropC2}-[1] since $\left(\left(\fmy_1\compone\fmy_2\right)\compone\fmy_3,i\right)\neq \left(\lboundone\fmx_1,...,\fmx_n\rboundone,1\right)$.
\end{proof}

\begin{proposition}\label{augpropl3} 
Fix $i \in\mathbb{N}$, $\bsfm\in \bsfM$, $\fmx,\fmy_1,\fmy_2,\fmy_3, \fmz\in \fremag$. Assume that both conditions below are fulfilled:
$\fmz$ is an elementary augmentation of $\bsfm\acmone\fmx$; $\left(\left(\fmy_1\compone\fmy_2\right)\compone\fmy_3,i\right)$ is the occurrence pair relative to the elementary augmentation $\fmz$ of $\bsfm\acmone\fmx$.\newline
Then there is $l\in\mathbb{N}$  such that both conditions below are fulfilled:
\begin{equation*}
\begin{array}{ll}
\left(\left(\fmy_1\compone\fmy_2\right)\compone\fmy_3,i-l\right)\;\text{occurs in}\;\fmx\text{;}\\[8pt]
\fmz=\bsfm\acmone\strchar \left[\fmx, \occfun\right]\text{,}\\
\text{where}\;\occfun\;\text{is the occurrence function on}
\left\{\left(\left(\fmy_1\compone\fmy_2\right)\compone\fmy_3,i-l\right)\right\}\\
\text{defined by setting} \;\occfun \left(\left(\fmy_1\compone\fmy_2\right)\compone\fmy_3,i-l\right)=\fmy_1\compone\left(\fmy_2\compone\fmy_3\right)\text{.}
\end{array}
\end{equation*}
\end{proposition}
\begin{proof}
Statement follows by Proposition \ref{elpropC3}-[1] since $\left(\left(\fmy_1\compone\fmy_2\right)\compone\fmy_3,i\right)\neq \left(\bsfm\acmone\fmx,1\right)$.
\end{proof}

\begin{proposition}\label{augprop}\mbox{}
\begin{enumerate}
\item Fix $\fmx_1,\fmx_2 \in \fremag$. Assume that $\fmx_2$ is an augmentation of $\fmx_1$. Then:
\begin{flalign}
& \codmagone\left(\fmx_1\right)=\codmagone\left(\fmx_2\right)\text{;}\nonumber\\[4pt]
& \dommagone\left(\fmx_1\right)=\dommagone\left(\fmx_2\right)\text{;}\nonumber\\[4pt]
& \fmx_1\ordrelsymb\fmx_2\label{augord}\text{.}
\end{flalign}
\item Fix $\fmx \in \fremag$. Then there is one and only one $ \fmz\in \fremag$ such that $\fmz$ is a maximal augmentation of $\fmx$.\newline
If $\fmx$ has order $0$ then $\fmz=\fmx$.\newline
If $\fmx$ has order $k \in \mathbb{N}_0$ then there are $m \in \mathbb{N}$, $\bsfm_{\mathsf{m},1},...,\bsfm_{\mathsf{m},n_{\mathsf{m}}}\in \bsfM$, $\fmx_{\mathsf{m},1},...,\fmx_{\mathsf{m},n_{\mathsf{m}}}\in \incllim_{k-1}\left(\fremag_{k-1}\right)$, $\fmz_{\mathsf{m},1},...,\fmz_{\mathsf{m},n_{\mathsf{m}}}\in \incllim_{k-1}\left(\fremag_{k-1}\right)$ for any $\mathsf{m }\in \left\{1,...,m\right\}$ such that all conditions below are fulfilled:
\begin{flalign*}
&\left\{\begin{array}{l}
\fmx=\fmx_1\compone...\compone\fmx_m\;\text{up to associativity, where}\\[4pt]
\fmx_{\mathsf{m}}=\lboundone\bsfm_{\mathsf{m},1}\acmone\fmx_{\mathsf{m},1},...,\bsfm_{\mathsf{m},n_{\mathsf{m}}}\acmone\fmx_{\mathsf{m},n_{\mathsf{m}}}\rboundone \hspace{20pt} \forall 
\mathsf{m }\in \left\{1,...,m\right\}\text{;}
\end{array}\right.\\[8pt]
&\left\{\begin{array}{l}\fmz_{\mathsf{m},\mathsf{n}}\;\text{is the maximal augmentation of}\;\fmx_{\mathsf{m},\mathsf{n}} \\[4pt]
\hspace{134pt}\forall \mathsf{m }\in \left\{1,...,m\right\}  \quad \forall\mathsf{n }\in \left\{1,...,n_{\mathsf{m}}\right\}\text{;}
\end{array}\right.\\[8pt]
&\left\{\begin{array}{l}
\fmz=\fmz_1\compone\left(\fmz_2\compone\left(...\compone\left(\fmz_{m-1}\compone\fmz_m\right)...\right)\right)\text{,}\quad \text{where}\\[4pt]
\fmz_{\mathsf{m}}=\lboundone\bsfm_{\mathsf{m},1}\acmone\fmz_{\mathsf{m},1},...,\bsfm_{\mathsf{m},n_{\mathsf{m}}}\acmone\fmz_{\mathsf{m},n_{\mathsf{m}}}\rboundone  
 \hspace{25pt}\forall\mathsf{m }\in \left\{1,...,m\right\}\text{.}
\end{array}\right.
\end{flalign*}
\item Fix $\fmx_{1},\fmx_{2}, \fmz_{1},\fmz_{2} \in \fremag$. Assume that all the following conditions are fulfilled: $\fmx_2$ is an augmentation of $\fmx_1$; $\fmz_i$ is the maximal augmentation of $\fmx_i$ for any $i \in \left\{1,2\right\}$. Then $\fmz_1=\fmz_2$.
\item Fix $\fmx_{1,1},\fmx_{2,1},\fmx_{1,2},\fmx_{2,2}, \fmz_1,\fmz_2 \in \fremag$. Assume that all the following conditions are fulfilled:
\begin{flalign}
&\fmx_{i,j}\;\text{is a maximal element}\hspace{111pt} \forall i,j \in \left\{1,2\right\}\text{;}\label{rappcomp1}\\[4pt]
&\fmx_{i,1}\ordrelsymb\fmx_{i,2}\hspace{179pt} \forall i \in \left\{1,2\right\}\text{;}\label{rappcomp}\\[4pt]
&\fmz_j\;\text{is the maximal augmentation of }\;\fmx_{1,j}\compone\fmx_{2,j}\hspace{25pt} \forall j \in \left\{1,2\right\}\text{.}
\end{flalign}
 Then $\fmz_1\ordrelsymb \fmz_2$.
\item Fix $m_1,m_2,m_3\in \mathbb{N}_0$, $\intuno_i\subseteqdentro \mathbb{R}^{m_i}$ for any $i \in \left\{1,2,3\right\}$, $\fmx_{1},\fmx_{2},\fmx_{3} \in \fremag$. Assume $\intuno_i\subseteq \dommagone\left(\fmx_i\right)$ for any $i \in \left\{1,2,3\right\}$. Then 
\begin{multline*}
\fmx_1\compone\left( \incl{\intuno_1}{\dommagone\left(\fmx_1\right)}\compone\left(\fmx_2\compone \left(\incl{\intuno_2}{\dommagone\left(\fmx_2\right)}\compone\left(\fmx_3\compone \incl{\intuno_3}{\dommagone\left(\fmx_3\right)}\right)\right)\right)\right)\ordrelsymb\\
\fmx_1\compone\left(\fmx_2\compone\left(\fmx_3\compone \incl{\intuno_3}{\dommagone\left(\fmx_3\right)}\right)\right)\text{.}
\end{multline*}
\item Fix $\fmx_1,\fmx_2\in \fremagcont$. Assume that $\fmx_2$ is an augmentation of $\fmx_1$. Then $\evalcompone\left(\fmx_1\right)=\evalcompone\left(\fmx_2\right)$.
\end{enumerate}
\end{proposition}
\begin{proof}\mbox{}\newline
\textnormal{\textbf{Proof of statement 1.}}\ \ Statement straightforwardly follows by definition of augmentation and by Proposition \ref{prpadrel1}-[3].\newline
\textnormal{\textbf{Proof of statements 2.}}\ \ We prove the statement by a three step induction argument based on the order of $\fmx$ (see Definition \ref{Cinfty}).\newline
\emph{\textbf{Step 1.}}\ \ If $-1$ is the order of $\fmx$ then statement straightforwardly holds true since the maximal augmentation of $\fmx$ coincides with $\fmx$ itself.\newline
\emph{\textbf{Step 2.}}\ \ If $k\geq 0$ is the order of $\fmx$ and $\fmx\in \incllim_k\left(\fremaggen_{k}\right)$ then there are $m\in \mathbb{N}$, $\bsfm_1,...,\bsfm_m \in \bsfM$, $\fmx_1,...,\fmx_m\in \incllim_{k-1}\left(\fremag_{k-1}\right)$ such that\newline
\centerline{$\fmx=\lboundone\bsfm_1\acmone\fmx_1,...,\bsfm_m\acmone\fmx_m\rboundone$.}
\newline
Fix a maximal augmentation $\fmz$ of $\fmx$. Iterated application of Propositions \ref{augpropl2}, \ref{augpropl3} entail that there are $\fmz_1,...,\fmz_m\in \incllim_{k-1}\left(\fremag_{k-1}\right)$ such that\newline
\centerline{$\fmz=\lboundone\bsfm_1\acmone\fmz_1,...,\bsfm_m\acmone\fmz_m\rboundone$.}\newline
Since $\fmz$ is a maximal augmentation of $\fmx$ we have that $\fmz_{\mathsf{m}}$ is a maximal augmentation of $\fmx_{\mathsf{m}}$ for any $\mathsf{m}\in \left\{1,...,m\right\}$. In fact, on the contrary, any elementary augmentation of an element $\fmz_{\mathsf{m}}$ would provide an elementary augmentation of $\fmz$ by Construction \ref{superevalsm}. Then statement follows since the maximal augmentation of $\fmx_{\mathsf{m}}$ is unique by induction for any $\mathsf{m}\in \left\{1,...,m\right\}$.\newline
\emph{\textbf{Step 3.}}\ \ If $k\geq 0$ is the order of $\fmx$ and $\fmx\in \incllim_k\left(\fremag_{k}\right)$ then there are $\fmx_1,....,\fmx_m\in \incllim_k\left(\fremaggen_{k}\right)$ such that $\fmx=\fmx_1\compone...\compone\fmx_m$ up to associativity. Fix a maximal augmentation $\fmz$ of $\fmx$. Iterated application of Propositions \ref{augpropl1}, \ref{augpropl2}, \ref{augpropl3} entails that that there are $\fmz_1,...,\fmz_m\in \incllim_k\left(\fremaggen_{k}\right)$, fulfilling both conditions below
\begin{flalign*}
&\begin{array}{l}
\fmz_{\mathsf{m}}\;\text{is a maximal augmentation of}\;\fmx_{\mathsf{m}}\hspace{20pt} \forall\mathsf{m}\in \left\{1,...,m\right\}\text{;} 
\end{array}\\[8pt]
&\begin{array}{l}
\fmz=\fmz_1\compone\left(\fmz_2\compone\left(...\compone\left(\fmz_{m-1}\compone\fmz_m\right)...\right)\right)\text{.}
\end{array}
\end{flalign*}
In fact, on the contrary, any arrangement of elements $\fmz_1,...,\fmz_m$ in $\fmz$ of type $\left(\fmw_1\compone\fmw_2\right)\compone\fmw_3$, or any elementary augmentation of an element $\fmz_{\mathsf{m}}$ would provide an elementary augmentation of $\fmz$ by Construction \ref{superevalsm}. Then statement follows since the maximal augmentation of $\fmx_{\mathsf{m}}$ is unique by step 2 for any $\mathsf{m}\in \left\{1,...,m\right\}$.\newline
\textnormal{\textbf{Proof of statements 3.}}\ \ Statement follows arguing as in the proof of Statement 2 by induction on the degree of elements belonging to $\fremag$.\newline
\textnormal{\textbf{Proof of statement 4.}}\ \  Assumptions straightforwardly entail that data $\fmz_1$, $\fmz_2$ fulfill conditions \eqref{eqelc0}, \eqref{eqelc}. We prove \eqref{eqelc1} by arguing as follows. Fix a path $\pmone_j$ in $\fremag$ detecting $\fmz_j$, a detecting skeleton $\left(\fmz_j, \assfun_j\right)$ of $\pmone_j$ for any $j \in \left\{1,2\right\}$ fulfilling \eqref{eqelc1}-[$(i)$].
Fix $i,j \in \left\{1,2\right\}$, construction of $\fremag$ entails that there is $k_{i,j}\in \mathbb{N}_0\cup\left\{-1\right\}$ with $\fmx_{i,j}\in \incllim_{k_{i,j}}\left(\fremag_{k_{i,j}}\right)$. Propositions \ref{augpropl2}, \ref{augpropl3}, statement 2, assumption \eqref{rappcomp1} entail that there are $\fmx_{i,j,1},...,\fmx_{i,j,k_{i,j}}\in \incllim_{k_{i,j}-1}\left(\fremaggen_{k_{i,j}}\right)$ for any $i,j \in \left\{1,2\right\}$ fulfilling both conditions below:
\begin{flalign*}
&\fmx_{i,j,\mathsf{k}} \;\text{is a maximal element}\hspace{20pt} \forall i,j \in \left\{1,2\right\}\quad\forall \mathsf{k}\in \left\{1,...,k_{i,j}\right\}\text{;}\\[8pt]
&\fmx_{i,j}=\fmx_{i,j,1}\compone\left(\fmx_{i,j,2}\compone\left(...\compone\left(\fmx_{i,j,k_{i,j}-1}\compone\fmx_{i,j,k_{i,j}}\right)...\right)\right)\hspace{20pt}  \forall i,j \in \left\{1,2\right\}\text{.}
\end{flalign*}
Propositions \ref{indassfun}, \ref{indindassfun}, \ref{incastrindpath}, \ref{pathoneprop2} entail that there are a path $\pmone_{i,j}$  detecting $\fmx_{i,j}$, a detecting skeleton $\left(\fmx_{i,j}, \assfun_{i,j}\right)$ of $\pmone_{i,j}$, a path $\pmone_{i,j,\mathsf{k}}$ detecting $\fmx_{i,j, \mathsf{k}}$, a detecting skeleton $\left(\fmx_{i,j, \mathsf{k}}, \assfun_{i,j, \mathsf{k}}\right)$ of $\pmone_{i,j,\mathsf{k}}$ for any $i,j \in \left\{1,2\right\}$, $\mathsf{k}\in \left\{1,...,k_{i,j}\right\}$ fulfilling both conditions below:
\begin{flalign}
&\hspace{-10pt}\begin{array}{l}
\pmone_{i,j}=\pmone_{i,j,1}\compone\left(\pmone_{i,j,2}\compone\left(...\compone\left(\pmone_{i,j,k_{i,j}-1}\compone\pmone_{i,j,k_{i,j}}\right)...\right)\right)\hspace{12pt}\forall i,j \in \left\{1,2\right\}\text{;}\label{compopt1}
\end{array}\\[8pt]
&\hspace{-10pt}\begin{array}{l}
\pmone_{j}=\pmone_{1,j,1}\compone\left(\pmone_{1,j,2}\compone\left(...\compone\left(\pmone_{1,j,k_{i,j}-1}\compone\left(\pmone_{1,j,k_{i,j}}\compone\pmone_{2,j}\right)\right)...\right)\right)\\
\hspace{255pt}\forall j \in \left\{1,2\right\}\text{.}\label{compopt2}
\end{array}
\end{flalign} 
By construction data $\pmone_{i,1}$, $\left(\fmx_{i,1}, \assfun_{i,1}\right)$, $\pmone_{i,2}$, $\left(\fmx_{i,2}, \assfun_{i,2}\right)$ fulfill condition \eqref{eqelc1}-[$(i)$] for any $i \in \left\{1,2\right\}$.\newline
Then assumption \eqref{rappcomp} entails that data $\pmone_{i,1}$, $\left(\fmx_{i,1}, \assfun_{i,1}\right)$, $\pmone_{i,2}$, $\left(\fmx_{i,2}, \assfun_{i,2}\right)$ fulfill condition \eqref{eqelc1}-[$(ii)$] for any $i \in \left\{1,2\right\}$. \newline
Fix $\unkuno \in \dommagone\left(\fmz_1\right)$, $\unkdue \in \left(-1,1\right)\setminus\left\{0\right\}$ with $\left(\evalcompone\left(\pmone_1\left(\unkdue\right)\right)\right)\left(\unkuno\right)\neq\smospempty$.\newline
Then \eqref{compopt1}, \eqref{compopt2} entail that both conditions below are fulfilled:
\begin{flalign}
&\left(\evalcompone\left(\pmone_{2,1}\left(\unkdue\right)\right)\right)\left(\unkuno\right)\neq\udenunk\text{;}\label{que1}\\[8pt]
&\left(\evalcompone\left(\pmone_{1,1}\left(\unkdue\right)\right)\right)\left(\left(\evalcompone\left(\pmone_{2,1}\left(\unkdue\right)\right)\right)\left(\unkuno\right)\right)\neq\udenunk\text{.}\label{que2}
\end{flalign}
Eventually condition \eqref{eqelc1}-[$(ii)$], assumption \eqref{rappcomp}, identity \eqref{compopt2}, inequalities \eqref{que1}, \eqref{que2} entail that 
\begin{multline*}
\left(\evalcompone\left(\pmone_{1}\left(\unkdue\right)\right)\right)\left(\unkuno\right)=\left(\evalcompone\left(\pmone_{1,1}\left(\unkdue\right)\right)\right)\left(\left(\evalcompone\left(\pmone_{2,1}\left(\unkdue\right)\right)\right)\left(\unkuno\right)\right)=\\
\left(\evalcompone\left(\pmone_{1,2}\left(\unkdue\right)\right)\right)\left(\left(\evalcompone\left(\pmone_{2,2}\left(\unkdue\right)\right)\right)\left(\unkuno\right)\right)=\left(\evalcompone\left(\pmone_{2}\left(\unkdue\right)\right)\right)\left(\unkuno\right)\text{.}
\end{multline*}  
\textnormal{\textbf{Proof of statement 5.}}\ \ Statement follows straightforwardly by checking axioms.\newline
\textnormal{\textbf{Proof of statement 6.}}\ \ Fix $\fmx,\fmy,\fmz\in\fremagcont$. Composition defined in Examples \ref{ExCinf}, \ref{ExCzero} and \eqref{evcmppropr}-[$(i)$] together entail that $\evalcompone\left(\left(\fmx\compone\fmy\right)\compone\fmz\right)=\evalcompone\left(\fmx\compone\left(\fmy\compone\fmz\right)\right)$. Then statement follows by Propositions \ref{augpropl1}-\ref{augpropl3}. 
\end{proof}

In Definition \ref{admrelc} below we introduce an equivalence relation which identifies elements belonging to $\fremag$ whenever paths detecting their maximal augmentations suitably coincides.

\begin{definition}\label{admrelc}
Fix $\fmx_1,\fmx_2 \in \fremag$. We say that $\fmx_1 \relsymb \fmx_2$ if and only if $\fmz_1 \relsymbord \fmz_2$ where $\fmz_i$ is the maximal augmentation of $\fmx_i$ for any $i \in \left\{1,2\right\}$.  
\end{definition}

In Proposition \ref{prpadrel} below we study the nature of relation $\relsymb$ and we examine its compatibility with operations $\compone$, $\lboundone\;\rboundone$, $\acmone$.

\begin{proposition}\label{prpadrel}
$\relsymb$ is an equivalence relation in $\fremag$ and is compatible with operations $\compone$, $\lboundone\;\rboundone$, $\acmone$. 
\end{proposition}
\begin{proof}
Statement follows straightforwardly by Propositions \ref{prpadrel1}-[1], \ref{augprop}-[2, 4].
\end{proof}

\section{Topological magma $\left(\magtwo,\comptwo\right)$\label{magtwosec}}
First we define magma $\left(\magtwo,\comptwo\right)$ as an algebraic object and study its basic properties, then we introduce topology on $\left(\magtwo,\comptwo\right)$ and study its basic properties.\newline

Motivated by \eqref{eqelc}, Definition \ref{admrelc}, Proposition \ref{prpadrel} and referring to Notation \ref{ins}-[13], we give Definition \ref{magintdef} below.

\begin{definition}\label{magintdef}
Magma $(\magtwo,\comptwo)$ is the quotient magma between $(\fremag,\compone)$ and equivalence relation $\relsymb$.\newline
We denote by $\quotmagone :(\fremag,\compone) \rightarrow (\magtwo,\comptwo)$ the quotient function of magmas.\newline
We set: 
\begin{equation*}
\linincltwo = \quotmagone \funcomp \lininclone\text{,}
\end{equation*}
\begin{equation*}
\magtwocont=\quotmagone (\fremagcont)\text{,} \qquad\magtwocontcont=\quotmagone \left(\fremagcontcont\right)\text{,}\qquad \magtwosmooth=\quotmagone (\fremagsmooth)\text{,}\qquad\magtwosmoothsmooth=\quotmagone (\fremagsmoothsmooth)\text{.}
\end{equation*}
The generic element belonging to $\magtwo$ will be denoted by $\mtx$.\newline
For any $\mtx\in \magtwo$ we denote by $\occset_{\mtx}$ any set $\occset_{\fmx}$ with $\mtx=\quotmagone\left(\fmx\right)$.\newline
For any $\mtx\in \magtwo$, $\fmx\in \fremag$ with $\fmx$ maximal element and $\mtx=\quotmagone\left(\fmx\right)$ we say that $\fmx$ is a maximal representative of $\mtx$.\newline
\end{definition}

\begin{remark} We emphasize that $\magtwo$ is not the single point set.\newline 
In fact, referring to \eqref{Bchoice2}, Definition \ref{Cinfty} we have that condition \eqref{eqelc} entails\newline
\centerline{$ \quotmagone\left(\inclfremaglev_0\left(\fmx\right)\right)\neq \quotmagone\left(\inclfremaglev_0\left(\fmy\right)\right)\qquad\forall \fmx,\fmy \in \fremag_{-1}$.}
\end{remark}

In Proposition \ref{alprtwo} below we state and prove basic algebraic properties of $\magtwo$.

\begin{proposition}\label{alprtwo}\mbox{}
\begin{enumerate}
\item Fix $n \in \mathbb{N}$. There is one and only one set function $\lboundtwo\quad\rboundtwo:\magtwo^n\rightarrow \magtwo$ defined by setting 
$\quad\lboundtwo\quad \rboundtwo\funcomp \quotmagone^n =\quotmagone\funcomp \lboundone\quad \rboundone$.\newline
We emphasize that in case $n=1$ we have $\lboundtwo\quad\rboundtwo=\idobj+\magtwo+$.
\item There is one and only one set function 
$\acmtwo:\bsfM \times \magtwo \rightarrow \magtwo$ defined by setting 
$ \qquad\acmtwo \funcomp \left(\idobj+\bsfM+, \quotmagone\right)= \quotmagone \funcomp \acmone$.\newline
Set function $\acmtwo$ fulfills the system of conditions below:
\begin{equation}
\left\{
\begin{array}{lll}
(i)&(\bsfm\bsfn) \acmtwo \mtx =\bsfm  \acmtwo( \bsfn \acmtwo \mtx)& \forall \bsfm ,\bsfn \in \bsfM\;\; \forall \mtx \in \magtwo\text{,} \vspace{4pt}\\
(ii)& \bsfm \acmtwo \mtx\in\magtwocont & \forall \bsfm \in \subbsfM\;\; \forall \mtx \in \magtwocont\text{,}
\vspace{4pt}\\
(iii)& \bsfm \acmtwo \mtx\in\magtwocontcont & \forall \bsfm \in \subbsfM\;\; \forall \mtx \in \magtwocontcont\text{,} \vspace{4pt}\\
(iv)& \bsfm \acmtwo \mtx\in\magtwosmooth & \forall \bsfm \in \bsfM\;\; \forall \mtx \in \magtwosmooth\text{,}
\vspace{4pt}\\
(v)& \bsfm \acmtwo \mtx\in\magtwosmoothsmooth & \forall \bsfm \in \bsfM\;\; \forall \mtx \in \magtwosmoothsmooth
\text{.} 
\end{array}
\right.\label{evalcomptwopropr}
\end{equation}
\item There is one and only one set function $\dommagtwo : \magtwo  \rightarrow \mathbf{P}$ fulfilling $\dommagone = \dommagtwo\funcomp\quotmagone $.
\item There is one and only one set function $\codmagtwo : \magtwo  \rightarrow \mathbf{P}$ fulfilling $\codmagone = \codmagtwo\funcomp\quotmagone $.
\item $\magtwocont$ is a sub-magma of $\magtwo$, $\magtwosmooth$ is a sub-magma of $\magtwocont$.
\item There is one and only one set function $\evalcomptwo : \magtwocont  \rightarrow \Cksp{0}$ defined by setting
\begin{equation}
\left\{
\begin{array}{l}
\evalcomptwo\left(\mtx\right)=\evalcompone\left(\fmz\right) \hspace{10pt}\forall \mtx\in \magtwocont\hspace{10pt} \forall\fmz\in\fremagcont\hspace{10pt} \text{with}\\[4pt] 
\mtx=\quotmagone\left(\fmz\right)\text{,}\hspace{20pt}\fmz	\; \text{is a maximal element.}
\end{array}
\right.\label{defevctw}
\end{equation}
Set function $\evalcomptwo$ fulfills the system of conditions below:
\begin{equation}\label{propevaltwo}
\left\{
\begin{array}{lll}
(i)& \evalcomptwo \left(\mtx\comptwo \mty\right)= \evalcomptwo \left(\mtx\right)\funcomp \evalcomptwo \left(\mty\right)& \forall \mtx,\mty \in \magtwocont\text{,}\\[4pt]
(ii)&\evalcomptwo\left(\lboundtwo \mtx_1,...,\mtx_n\rboundtwo\right)=\vspace{4pt}\\
&\hspace{40pt}\left(\evalcomptwo\left(\mtx_1\right),...,\evalcomptwo\left(\mtx_n\right)\right)&\forall  \left(\mtx_1,...,\mtx_n \right)\in \magtwocont^n\text{,}\vspace{4pt}\\
(iii)&\evalcomptwo\left(\bsfm \acmtwo\mtx\right)=\bsfm\contspbsfmu\evalcomptwo\left(\mtx\right)&\forall \bsfm\in\subbsfM\quad\forall \mtx \in \magtwocont \text{,}\vspace{4pt}\\
(iv)&\evalcomptwo\left(\bsfm \acmtwo\mtx\right)=\bsfm \smospbsfmu \evalcomptwo\left(\mtx\right)& \forall \bsfm\in\bsfM\quad\forall \mtx \in \magtwosmooth\text{,}\vspace{4pt}\\
(v)& \evalcomptwo \left(\mtx\right)\in \Cksp{\infty} & \forall  \mtx \in \magtwosmooth\text{.}
\end{array}
\right.
\end{equation} 
\end{enumerate}
\end{proposition}
\begin{proof}\mbox{}\newline
\textnormal{\textbf{Proof of statements 1, 2.}}\ \ Statements follow by compatibility of $\relsymb$ with operations $\lboundone\;\rboundone$, $\acmone$ (see Proposition \ref{prpadrel}). \newline
\textnormal{\textbf{Proof of statements 3, 4.}}\ \ Statements follow since for any $\fmx_1,\fmx_2\in \fremag$ with $\fmx_1 \relsymb\fmx_2$ we have also $\dommagone\left(\fmx_1\right)=\dommagone\left(\fmx_2\right)$, $\codmagone\left(\fmx_1\right)=\codmagone\left(\fmx_2\right)$.\newline
\textnormal{\textbf{Proof of statement 5.}}\ \ Statements follow by compatibility of $\relsymb$ with the operation $\compone$ (see Proposition \ref{prpadrel}). \newline 
\textnormal{\textbf{Proof of statement 6.}}\ \ Set function $\evalcomptwo$ is everywhere defined on $\magtwocont$ by Proposition \ref{augprop}-[2] and it is well defined on $\magtwocont$ since \eqref{extopcntdis5}, Proposition \ref{funoncont}, relation $\relsymb$ together entail that $\evalcompone\left(\fmx_1\right)=\evalcompone\left(\fmx_2\right)$ for any $\fmx_1,\fmx_2 \in \fremagcont$ with $\fmx_1$, $\fmx_2$ maximal elements and $\quotmagone\left(\fmx_1\right)=\quotmagone\left(\fmx_2\right)$. Eventually \eqref{propevaltwo} follows by \eqref{evcmppropr},  Proposition \ref{augprop}-[2, 6].
\end{proof}

Motivated by Proposition \ref{alprtwo}-[3, 4] we introduce set functions giving the dimension of domain and co-domain of elements belonging $\magtwo$. We refer to Definition \ref{dimdomcodfuncontone}.
\begin{definition}\label{dimdomcodfuncont}
We define set functions
\begin{flalign*}
& \dimcodmagtwo:\magtwo\rightarrow\mathbb{N}_0\quad\text{by setting}  \quad \dimcodmagtwo\left(\quotmagone\left(\fmx\right)\right)=\dimcodmagone\left(\fmx\right)\quad\forall \fmx \in \fremag \text{;}\\[8pt]
&  \dimdommagtwo:\magtwo\rightarrow\mathbb{N}_0\quad\text{by setting}  \quad \dimdommagtwo\left(\quotmagone\left(\fmx\right)\right)=\dimdommagone\left(\fmx\right)\quad\forall\fmx \in \fremag\text{.}
\end{flalign*}
\end{definition}

In Definition \ref{pathtwodef} below we introduce the notion of path in sets $\magtwo^n$. We refer to Definition \ref{pathinC}.

\begin{definition}\label{pathtwodef}\mbox{}
\begin{enumerate}
\item Fix $n \in \mathbb{N}$, a set function $\pmtwo:\left(-1,1\right) \rightarrow \magtwo^n$.\newline
We say that $\pmtwo$ is a path in $\magtwo^n$ if and only if there is a path $\pmone$ in $\fremag^n$ such that $\pmtwo = \quotmagone^n \funcomp \pmone$.\newline
With an abuse of language we say that: any skeleton $\left(\left(\fmx_1,...,\fmx_n\right),\assfun\right)$ of $\pmone$ is a skeleton of $\pmtwo$; $\left(\fmx_1,...,\fmx_n\right)$ is a core of $\pmtwo$; $\assfun$ is an associating function of $\pmtwo$.
\item Fix $n \in \mathbb{N}$, a path $\pmtwo$ in $\magtwo^n$, a subset $\singsmpth \subseteq (-1,1)$.\newline 
We say that $\left(\pmtwo,\singsmpth\right)$ is a smooth path in $\magtwo^n$ if and only if there is a smooth path $\left(\pmone,\singsmpth\right)$ in $\fremag^n$ such that $\pmtwo = \quotmagone^n \funcomp \pmone$.\newline
With an abuse of language we say that:
$\singsmpth$ is the singular set of $\left(\pmtwo, \singsmpth\right)$;
any skeleton $\left(\left(\fmx_1,...,\fmx_n\right),\assfun\right)$ of $\left(\pmone,\singsmpth\right)$ is a skeleton of $\left(\pmtwo,\singsmpth\right)$.\newline
We emphasize that nothing is assumed about $\pmtwo\left(\unkdue\right)$ when $\unkdue \in \singsmpth$.
\item Fix $n \in \mathbb{N}$, $\left(\mtx_1,...,\mtx_n\right)\in \magtwo^n$, a path $\pmtwo$ in $\magtwo^n$.\newline
We say that $\pmtwo$ is a path in $\magtwo^n$ through $\left(\mtx_1,...,\mtx_n\right)$ if and only if we have $\pmtwo(0)=\left(\mtx_1,...,\mtx_n\right)$.
\item Fix $\mtx\in \magtwo$, a path $\pmtwo$ in $\magtwo$ through $\mtx$.
We say that $\pmtwo$ is a path in $\magtwo$ detecting $\mtx$ if and only if 
there is $x\in \fremag$, a path $\pmone$ in $\fremag$ detecting $x$ such that $\pmtwo = \quotmagone \funcomp \pmone$.\newline
With an abuse of language we say that any detecting skeleton $\left(\fmx,\assfun\right)$ of $\pmone$ is a detecting skeleton of $\pmtwo$.
\end{enumerate}
\end{definition}

In Propositions \ref{domcodpathtwo}, \ref{pathprop1}, \ref{pathprop2}, \ref{pathprop3}, \ref{pathprop4}, \ref{pathprop5}, \ref{pathprop5bis}, \ref{pathprop8} below we study the link between algebraic structure of $\magtwo$ and paths in $\magtwo$.

\begin{proposition}\label{domcodpathtwo}
Fix a path $\pmtwo$ in $\magtwo$, a skeleton $\left(\fmx,\assfun\right)$ of $\pmtwo$. Then
\begin{equation*}
\codmagtwo\left(\pmtwo\left(\unkdue\right)\right)=\codmagone\left(\fmx\right)\text{,}\qquad \dommagtwo\left(\pmtwo\left(\unkdue\right)\right)=\dommagone\left(\fmx\right)\qquad \forall \unkdue \in \left(-1,1\right)\text{.}
\end{equation*}
\end{proposition}
\begin{proof}
Statement follows by  Definition \ref{pathtwodef}, Propositions \ref{domcodpathone}, \ref{augprop}-[1], \ref{alprtwo}-[3, 4].
\end{proof}

Motivated by Proposition \ref{domcodpathtwo} we extend to continuous functions with path connected domain the notions of domain and co-domain.
 
\begin{definition}\label{domcodfuncont}
Fix a path connected topological space $\setsymuno$, a continuous function $f:\setsymuno\rightarrow \magtwo$.  We define: 
\begin{flalign*}
& \codmagtwo\left(\pmtwo\right)=\codmagtwo\left(\pmtwo\left(\elsymuno\right)\right)\hspace{10pt}\elsymuno\in\setsymuno\text{;}\\[6pt]
& \dommagtwo\left(\pmtwo\right)=\dommagtwo\left(\pmtwo\left(\elsymuno\right)\right)\hspace{10pt}\elsymuno\in\setsymuno\text{;}\\[6pt]
& \dimcodmagtwo\left(\pmtwo\right)=\codmagtwo\left(\pmtwo\left(\elsymuno\right)\right)\hspace{10pt}\elsymuno\in\setsymuno\text{;}\\[6pt]
& \dimdommagtwo\left(\pmtwo\right)=\dommagtwo\left(\pmtwo\left(\elsymuno\right)\right)\hspace{10pt}\elsymuno\in\setsymuno\text{.}
\end{flalign*}
We say that: $\codmagtwo\left(\pmtwo\right)$ is the co-domain of $\pmtwo$; $\dommagtwo\left(\pmtwo\right)$ is the domain of $\pmtwo$; $\dimcodmagtwo\left(\pmtwo\right)$ is the dimension of the co-domain of $\pmtwo$; $\dimdommagtwo\left(\pmtwo\right)$ is the dimension of the domain of $\pmtwo$.
\end{definition}

Motivated by \eqref{eqelc}, Proposition \ref{indassfun}, in Definition \ref{subpathmtw} we introduce the notion of induced paths.

\begin{definition}\label{subpathmtw}\mbox{}
\begin{enumerate}
\item  Fix $\mtx, \mty\in \magtwo$, a path $\pmtwo$ in $\magtwo$, a skeleton $\left(\fmx, \assfun\right)$ of $\pmtwo$, $\fmy\in \fremag$.\newline
Assume that all the following conditions hold true:
\begin{equation*}
 \occset_{\mty}\subseteq \occset_{\mtx}\text{;}\hspace{20pt} \mtx=\quotmagone\left(\fmx\right)\text{;}\hspace{20pt}\mty=\quotmagone\left(\fmy\right)\text{.}
\end{equation*}
We denote by $\pmtwo\left[\mty, \left(\mtx, \assfun\right)\right]$ the path in $\magtwo$ obtained by composing with $\quotmagone$ the path in $\fremag$ defined by $\left(\fmy, \assfun\funcomp\left(\incl{\occset_{\mty}}{\occset_{\mtx}}, \idobj+\left(-1,1\right)+\right)\right)$ through \eqref{pathconduro2}.\newline
We say that $\pmtwo\left[\mty, \left(\mtx, \assfun\right)\right]$ is the path induced on $\mty$ by $\left(\mtx, \assfun\right)$. By Remark \ref{notunsk}-[2] we say that $\pmtwo\left[\mty, \left(\mtx, \assfun\right)\right]$ is a path induced on $\mty$ by $\pmtwo$. 
\item Fix $\mtx, \mty\in \magtwo$, a smooth path $\left(\pmtwo, \singsmpth\right)$ in $\magtwo$, a skeleton $\left(\fmx, \assfun\right)$ of $\pmtwo$, $\fmy\in \fremag$.\newline
Assume that all the following conditions hold true: 
\begin{equation*}
 \occset_{\mty}\subseteq \occset_{\mtx}\text{;}\hspace{20pt} \mtx=\quotmagone\left(\fmx\right)\text{;}\hspace{20pt}\mty=\quotmagone\left(\fmy\right)\text{.}
\end{equation*}
We denote by $\left(\pmtwo\left[\mty, \left(\mtx, \assfun\right)\right], \singsmpth\right)$ the smooth path in $\magtwo$ obtained by composing with $\quotmagone$ the smooth path in $\fremag$ defined by $\left(\fmy, \assfun\funcomp\left(\incl{\occset_{\mty}}{\occset_{\mtx}}, \idobj+\left(-1,1\right)+\right)\right)$ through \eqref{pathconduro2}.\newline
We say that $\left(\pmtwo\left[\mty, \left(\mtx, \assfun\right)\right], \singsmpth\right)$ is the smooth path induced on $\mty$ by $\left(\mtx, \assfun\right)$. By Remark \ref{notunsk}-[2] we say that $\left(\pmtwo\left[\mty, \left(\mtx, \assfun\right)\right], \singsmpth\right)$ is a smooth path induced on $\mty$ by $\pmtwo$.
\end{enumerate}
\end{definition}

\begin{proposition}\label{pathprop1}
Fix $n_1,n_2 \in \mathbb{N}$, a path $\pmtwo_1$ in $\magtwo^{n_1}$, a path $\pmtwo_2$ in $\magtwo^{n_2}$. 
\begin{enumerate}
\item  The set function $\pmtwo=\left(\pmtwo_1,\pmtwo_2\right)\funcomp \Diag{\left(-1,1\right)}{2}$ is a path in $\magtwo^{n_1+n_2}$.
\item If there are subsets $\singsmpth_1,\singsmpth_2 \subseteq \left(-1,1\right)$ such that pair $\left(\pmtwo_1,\singsmpth_1\right)$ is a smooth path in $\magtwo^{n_1}$ and pair $\left(\pmtwo_2,\singsmpth_2\right)$ is a smooth path in $\magtwo^{n_2}$ then the pair $\left(\left(\pmtwo_1,\pmtwo_2\right)\funcomp \Diag{\left(-1,1\right)}{2},\singsmpth_1\cup\singsmpth_2 \right)$ is a smooth path in $\magtwo^{n_1 + n_2}$.
\item Fix $I \subseteq \left\{1,2\right\}$, $J_i \in \left\{1,...,n_i\right\}$ and a skeleton $\left(\left(\fmx_{i,1},...,\fmx_{i,n_i}\right), \assfun_i\right)$ of $\pmtwo_i$ for any $i \in I$, a skeleton $\left(\left(\fmy_1,...,\fmy_{n_1+n_2}\right), \assfun_{\fmy}\right)$ of $\pmtwo$.\newline
Assume $\quad \fmx_{i,\mathsf{n}} \in \fremagcont\quad \forall i\in I\quad \forall \mathsf{n}\in J_i$.\quad Set $n_0=0$.\newline
Then exists a skeleton $\left(\left(\fmx_1,...,\fmx_{n_1+n_2}\right), \assfun\right)$ of $\pmtwo$ such that: 
\begin{flalign*}
&\fmx_k=\fmy_k \hspace{20pt}\text{if} \quad k  \neq n_{i-1}+\mathsf{n}\quad \forall i \in I \quad \forall \mathsf{n}\in J_i \text{;}\\[4pt]
&\fmx_k\in \fremagcont \hspace{18pt}\text{if}\quad \exists i \in I \quad \exists \mathsf{n}\in J_i\quad \text{such that}\quad k  = n_{i-1}+\mathsf{n}\text{.}
\end{flalign*}
\end{enumerate}
\end{proposition}
\begin{proof}
Statements follow by corresponding statements in Proposition \ref{pathoneprop1}.  
\end{proof}

\begin{proposition}\label{pathprop2}
Fix two paths $\pmtwo_1, \pmtwo_2$ in $\magtwo$.
\begin{enumerate}
\item  The set function $\pmtwo_1\comptwo \pmtwo_2:\left(-1,1\right)\rightarrow\magtwo$ defined by setting \begin{equation*}
\left(\pmtwo_1\comptwo \pmtwo_2\right)\left(\unkdue\right)=\pmtwo_1\left(\unkdue\right)\comptwo\pmtwo_2\left(\unkdue\right)\qquad \forall \unkdue\in \left(-1,1\right)
\end{equation*}
is a path in $\magtwo$.
\item If there are subsets $\singsmpth_1,\singsmpth_2 \subseteq \left(-1,1\right)$ such that pairs $\left(\pmtwo_1,\singsmpth_1\right)$, $\left(\pmtwo_2,\singsmpth_2\right)$ are smooth paths in $\magtwo$ then pair $\left(\pmtwo_1\comptwo \pmtwo_2,\singsmpth_1\cup\singsmpth_2 \right)$ is a smooth path in $\magtwo$.
\item If there is a skeleton $\left(\fmx_1, \assfun_1\right)$ of $\pmtwo_1$, a skeleton $\left(\fmx_2, \assfun_2\right)$ of $\pmtwo_2$ with $\fmx_1,\fmx_2 \in \fremagcont$ then there is a skeleton $\left(\fmx, \assfun\right)$ of $\pmtwo_1\comptwo \pmtwo_2$ with $\fmx \in \fremagcont$.
\end{enumerate}
\end{proposition}
\begin{proof}
Statements follow by corresponding statements in Proposition \ref{pathoneprop2}.
\end{proof}

\begin{proposition}\label{pathprop3}
Fix $\mtx,\mty \in \magtwo$, a path $\pmtwo$ in $\magtwo$ detecting $\mtx\comptwo \mty$.\newline
Then there are a path $\pmtwo_{\mtx}$ in $\magtwo$ detecting $\mtx$, a path $\pmtwo_{\mty}$ in $\magtwo$ detecting $\mty$
such that\hspace{20pt}$\pmtwo\left(\unkdue\right) = \left(\pmtwo_{\mtx}\comptwo \pmtwo_{\mty}\right)\left(\unkdue\right)\qquad \forall \unkdue \in \left(-1,1\right)\setminus\left\{0\right\}$.\newline
We emphasize that \hspace{20pt}$\left(\evalcomptwo\left(\pmtwo\left(\unkdue\right)\right)\right)\left(\unkuno\right)\neq \udenunk$\hspace{20pt} entails:\newline
\centerline{$\left(\evalcomptwo\left( \pmtwo_{\mty}\left(\unkdue\right)\right)\right)\left(\unkuno\right)\neq \udenunk\text{;}\quad\left(\evalcomptwo\left(\pmtwo_{\mtx}\left(\unkdue\right)\right)\right)\left(\left(\evalcomptwo\left( \pmtwo_{\mty}\left(\unkdue\right)\right)\right)\left(\unkuno\right)\right)\neq \udenunk$.}
\end{proposition}
\begin{proof} Choose: $\fmx,\fmy \in\fremag$ with $\quotmagone\left(\fmx\right)=\mtx$, $\quotmagone\left(\fmy\right)=\mty$; a detecting skeleton $\left(\fmz, \assfun\right)$ of $\pmtwo$. We have $\fmz \relsymb \fmx\compone \fmy$ then by Proposition \ref{indassfun}, relation \eqref{(R 17)} there is a path $\pmtwo_1$ in $\magtwo$ detecting $\mtx \comptwo \mty$ fulfilling both conditions below:
\begin{flalign*}
&\left(\fmx\compone \fmy, \assfun\right)\;\text{is a detecting skeleton of}\;\pmtwo_1\text{;}\\[8pt]
&\pmtwo\left(\unkdue\right)= \pmtwo_1\left(\unkdue\right)\qquad \forall \unkdue \in \left(-1,1\right)\setminus\left\{0\right\}\text{.}
\end{flalign*}
By Proposition \ref{elpropC4} we set:\newline
\centerline{
$\assfun_{\fmx}=\assfun \funcomp \left(\incl{\occset_{\fmx}}{\occset_{\fmz}}, \idobj+\left(-1,1	\right)+\right)\text{,} \qquad\assfun_{\fmy}=\assfun \funcomp \left(\incl{\occset_{\fmy}}{\occset_{\fmz}}, \idobj+\left(-1,1\right)+\right)$.}\newline
Eventually statement follows by Definition \ref{pathtwodef}-[4], Propositions \ref{indindassfun}-[3], \ref{incastrindpath} if we define paths $\pmtwo_{\mtx}$, $\pmtwo_{\mty}$ through \eqref{pathconduro2} applied with data $\left(\fmx, \assfun_{\fmx}\right)$, $\left(\fmy, \assfun_{\fmy}\right)$ respectively.
\end{proof}

\begin{proposition}\label{pathprop4} Fix $n \in \mathbb{N}$, an ordered $n$-tuple $\left(\pmtwo_1,...,\pmtwo_n\right)$ of paths in $\magtwo$.
\begin{enumerate}
\item The set function $\lboundtwo\pmtwo_1,...,\pmtwo_n\rboundtwo:\left(-1,1\right)\rightarrow\magtwo$ defined by setting
\begin{equation*}\lboundtwo\pmtwo_1,...,\pmtwo_n\rboundtwo\left(\unkdue\right)=\lboundtwo\pmtwo_1\left(\unkdue\right),...,\pmtwo_n\left(\unkdue\right)\rboundtwo\qquad \forall \unkdue\in \left(-1,1\right)
\end{equation*}
is a path in $\magtwo$.
\item If there is a subset $\singsmpth_{\mathsf{n}} \subseteq \left(-1,1\right)$ such that the pair $\left(\pmtwo_{\mathsf{n}},\singsmpth_{\mathsf{n}}\right)$ is a smooth path in $\magtwo$ for any $\mathsf{n} \in \left\{1,...,n\right\}$ then pair $\left(\lboundtwo\pmtwo_1,...,\pmtwo_n\rboundtwo,\underset{\mathsf{n}=1}{\overset{n}{\bigcup}}\singsmpth_{\mathsf{n}} \right)$ is a smooth path in $\magtwo$.
\item If there is a skeleton $\left(\fmx_{\mathsf{n}}, \assfun_{\mathsf{n}}\right)$ of $\pmtwo_{\mathsf{n}}$ with $\fmx_{\mathsf{n}} \in \fremagcont$ for any $\mathsf{n} \in \left\{1,...,n\right\}$ then there is a skeleton $\left(\fmx, \assfun\right)$ of $\lboundtwo\pmtwo_1,...,\pmtwo_n\rboundtwo$ with $\fmx \in \fremagcont$.
\end{enumerate}
\end{proposition}
\begin{proof}
Statements follow by corresponding statements in Proposition \ref{pathoneprop4}.
\end{proof}

\begin{proposition}\label{pathprop4bis}
Fix $n \in \mathbb{N}$, $\left(\mtx_1,...,\mtx_n\right) \in \magtwo^n$, a path $\pmtwo$ in $\magtwo$ detecting $\lboundtwo \mtx_1,...,\mtx_n\rboundtwo$.\newline 
Then there is a path $\pmtwo_{\mathsf{n}}$ in $\magtwo$ detecting $\mtx_{\mathsf{n}}$ for any $\mathsf{n} \in\left\{1,...,n\right\}$ such that\hspace{20pt}$ \pmtwo\left(\unkdue\right)=\lboundtwo\pmtwo_1,...,\pmtwo_n\rboundtwo\left(\unkdue\right)\qquad \forall \unkdue \in \left(-1,1\right)\setminus\left\{0\right\}$.\newline
We emphasize that $\left(\evalcomptwo\left(\pmtwo\left(\unkdue\right)\right)\right)\left(\unkuno\right)\neq \udenunk$ entails 
$\left(\evalcomptwo\left(\pmtwo_{\mathsf{n}}\left(\unkdue\right)\right)\right)\left(\unkuno\right)\neq \udenunk$ for any $\mathsf{n}\in \left\{1,...,n\right\}$.
\end{proposition}
\begin{proof} Choose: $\left(\fmx_{1},...,\fmx_{n}\right) \in\fremag^n$ with $\quotmagone\left(\fmx_{\mathsf{n}}\right)=\mtx_{\mathsf{n}}$ for any $\mathsf{n} \in \left\{1,...,n\right\}$; a detecting skeleton $\left(\fmz, \assfun\right)$ of $\pmtwo$. We have $\fmz \relsymb \lboundtwo \fmx_1,...,\fmx_n\rboundtwo$ then by Proposition \ref{indassfun}, relation \eqref{(R 17)} there is path $\pmtwo_1$ in $\magtwo$ detecting $\lboundtwo \mtx_1,...,\mtx_n\rboundtwo$ fulfilling both conditions below:
\begin{flalign*}
&\left(\lboundtwo \fmx_1,...,\fmx_n\rboundtwo, \assfun\right)\;\text{is a detecting skeleton of}\;\pmtwo_1\text{;}\\[8pt]
&\pmtwo\left(\unkdue\right)= \pmtwo_1\left(\unkdue\right)\qquad \forall \unkdue \in \left(-1,1\right)\setminus\left\{0\right\}\text{.}
\end{flalign*}
By Proposition \ref{elpropC4} we set
$\assfun_{i}=\assfun \funcomp \left(\incl{\occset_{\fmx_{\mathsf{n}}}}{\occset_{\fmz}}, \idobj+\left(-1,1\right)+\right)$ for any $\mathsf{n} \in\left\{1,...,n\right\}$.\newline
Eventually statement follows by Definition \ref{pathtwodef}-[4], Propositions \ref{indindassfun}-[3], \ref{incastrindpath} if we define path $\pmtwo_{\mathsf{n}}$ through \eqref{pathconduro2} applied with datum $\left(\fmx_{\mathsf{n}}, \assfun_{\mathsf{n}}\right)$ for any $\mathsf{n} \in\left\{1,...,n\right\}$.
\end{proof}

\begin{proposition}\label{pathprop5}
Fix a path $\pmtwo$ in $\magtwo$, $\bsfm\in\bsfM$.
\begin{enumerate}
\item The set function $\bsfm \acmtwo \pmtwo:\left(-1,1\right)\rightarrow\magtwo$ defined by setting
\begin{equation*}\left(\bsfm \acmtwo \pmtwo\right)\left(\unkdue\right)=\bsfm\acmtwo\pmtwo\left(\unkdue\right)\qquad \forall \unkdue\in \left(-1,1\right)
\end{equation*}
is a path in $\magtwo$.\newline
\item If there is a subset $\singsmpth \subseteq \left(-1,1\right)$ such that the pair $\left( \pmtwo,\singsmpth\right)$ is a smooth path in $\magtwo$ then pair $\left(\bsfm \acmtwo \pmtwo,\singsmpth \right)$ is a smooth path in $\magtwo$.\newline
\item If $\bsfm\in\subbsfM$ and there is a skeleton $\left(\fmy, \assfun_{\fmy}\right)$ of $\pmtwo$ with $\fmy \in \fremagcont$ then there is a skeleton $\left(\fmx, \assfun\right)$ of $\bsfm \acmtwo \pmtwo$ with $\fmx \in \fremagcont$.
\end{enumerate}
\end{proposition}
\begin{proof}
Statements follow by corresponding statements in Proposition \ref{pathoneprop5}.
\end{proof}

\begin{proposition}\label{pathprop5bis}
Fix $\bsfm \in \bsfM$, $\mtx \in \magtwo$, a path $\pmtwo$ in $\magtwo$ detecting $\bsfm \acmtwo \mtx$. \newline
Then there is a path $\pmtwo_{\mtx}$ in $\magtwo$ detecting $\mtx$ such that\newline
\centerline{$\pmtwo\left(\unkdue\right)=\left(\bsfm\acmone\pmtwo_{\mtx}\right)\left(\unkdue\right)\qquad \forall \unkdue \in \left(-1,1\right)\setminus\left\{0\right\}$.}\newline
We emphasize that $\evalcomptwo\left(\pmtwo\left(\unkdue\right)\right)\neq \contspempty$ entails
$\evalcomptwo\left(\pmtwo_{\mtx}\left(\unkdue\right)\right)\neq \contspempty$.
\end{proposition}
\begin{proof} Choose: $\fmx_1 \in\fremag$ with $\quotmagone\left(\fmx_1\right)=\mtx$; a detecting skeleton $\left(\fmz, \assfun\right)$ of $\pmtwo$.
We have $\fmz \relsymb \bsfm\acmone \fmx$ then by Proposition \ref{indassfun}, relation \eqref{(R 17)} there is a path $\pmtwo_1$ in $\magtwo$ detecting $\bsfm \acmtwo \mtx$ fulfilling both conditions below::
\begin{flalign*}
&\left(\bsfm\acmone \fmx, \assfun\right)\;\text{is a detecting skeleton of}\;\pmtwo_1\text{;}\\[8pt]
&\pmtwo\left(\unkdue\right)= \pmtwo_1\left(\unkdue\right)\qquad \forall \unkdue \in \left(-1,1\right)\setminus\left\{0\right\}\text{.}
\end{flalign*}
By Proposition \ref{elpropC4} we set
$\assfun_{\fmx}=\assfun \funcomp \left(\incl{\occset_{\fmx}}{\occset_{\fmz}}, \idobj+\left(-1,1\right)+\right)$.\newline
Eventually statement follows by Definition \ref{pathtwodef}-[4], Propositions \ref{indindassfun}-[3], \ref{incastrindpath} if we define path $\pmtwo_{\mtx}$ through \eqref{pathconduro2} applied with datum $\left(\fmx, \assfun_{\fmx}\right)$. 
\end{proof}

\begin{proposition}\label{pathprop8}
Fix $n \in \mathbb{N}$, a path $\pmtwo$ in $\magtwo^n$ and a set $\singsmpth\subseteq [0,1)$. Assume that all conditions below are fulfilled:
\begin{equation}
0\in \singsmpth\text{,}\hspace{30pt} \singsmpth \cap (\varepsilon ,1)\;\;\text{is a finite set for any}\;\;\varepsilon >0\text{.}
\label{assclaimg1}
\end{equation} 
Then there is a smooth path $(\pmtwo_1,\singsmpth)$ in $\magtwo^n$ fulfilling  
\begin{equation}
\pmtwo_1\left(\unkdue\right)=\pmtwo\left(\unkdue\right) \qquad \forall \unkdue \in \singsmpth\text{.}
\label{claimg1}
\end{equation} 
\end{proposition}
\begin{proof}\mbox{}\newline 
Referring to Definitions \ref{pathinC}, \ref{pathtwodef}-[1] we choose a skeleton $\left(\left(\fmx_1,...,\fmx_n\right), \assfun\right)$ of $\pmtwo$. Then there is a set function $\assfun_1: \left(\underset{\mathsf{n}=1}{\overset{n}{\bigcup}} \,\occset_{\fmx_{\mathsf{n}}}\right)\times \left(-1,1\right)  \rightarrow \Cksp{0}$ fulfilling all conditions below:
\begin{flalign*}
&\text{\eqref{pathconduro};} \hspace{40pt}\text{\eqref{smptass}-[}(iii)\text{] with data}\; \assfun_1\text{,}\;\singsmpth\text{;}\\[4pt]
& \assfun_1\left(f,\unkdue\right)=\assfun\left(f,\unkdue\right)\qquad \forall \left(f,\unkdue\right)\in  \left(\underset{\mathsf{n}=1}{\overset{n}{\bigcup}} \,\occset_{\fmx_{\mathsf{n}}}\right)\times \singsmpth\text{.}
\end{flalign*}
Existence of $\assfun_1$ follows by Weierstrass approximation theorem and the path connectivity of $\Cksp{0}(\intuno)+\mathbb{R}^n+$ (refer to Notation \ref{difrealfunc}-[2]).\newline
Statement is achieved by setting  $\pmtwo_1 = \quotmagone^n \funcomp \pmone_1$ where $\pmone_1$ is the unique path in $\fremag^n$ defined by the pair $\left(\left(\fmx_1,...,\fmx_n\right), \assfun_1\right)$ through \eqref{pathconduro2}.
\end{proof}

\begin{conjecture}\label{pathconj7}
Fix $\mtx \in\magtwocontcont$. Then there is a path $\pmtwo$ detecting $\mtx$ such that
\hspace{20pt}$\left(\evalcomptwo\left(\pmtwo \left(\unkdue\right)\right)\right)\left(\unkuno\right)\neq \udenunk \hspace{10pt} \forall \unkdue \in \left(-1,1\right)\setminus\left\{0\right\}\quad \forall \unkuno\in \dommagtwo\left(\mtx\right)$.
\end{conjecture}
Conjecture \ref{pathconj7} seems to be quite natural, however we were not able to prove it. Fortunately for our aims (Chapter \ref{App}) we do not need to prove Conjecture \ref{pathconj7} but a weaker version of it (Proposition \ref{pathprop7}).\\[12pt]

In Definition \ref{topmagtwo} below we give topology to sets $\magtwo^n$ through paths introduced in Definition \ref{pathtwodef}. We refer to Notations \ref{ins}-[5], \ref{gentopnot}-[5].

\begin{definition}\label{topmagtwo}\mbox{}
\begin{enumerate}
\item Fix $n \in \mathbb{N}$. Set $\magtwo^n$ is endowed with the final topology with respect to the family of all paths in $\magtwo^n$. We say that such topology is the the path topology on $\magtwo^n$.
\item We set 
\begin{multline*}
\magtwoempty =\Big\{\mtx \in \magtwo :\;\; \forall\;\text{path}\;\pmtwo\; 
\text{in} \;\magtwo\;\text{detecting}\; \mtx\text{,}\hspace{10pt}\forall \delta>0\\[4pt]
\exists \unkdue\in \left(-\delta,\delta\right)\setminus\left\{0\right\}\quad \exists\unkuno\in \dommagtwo\left(\mtx\right)\,:\;
\left(\evalcomptwo\left( \pmtwo\left(\unkdue\right)\right)\right)\left(\unkuno\right)=\udenunk
\Big\}\text{.}
\end{multline*}
\end{enumerate}
\end{definition}

In Proposition \ref{topinDprop} below we state and prove basic properties of the path topology introduced in Definition \ref{topmagtwo}. We refer to Notation \ref{gentopnot}-[4, 5], Definitions \ref{intdiffmon}, \ref{admrelc}.

\begin{proposition}\label{topinDprop}\mbox{}
\begin{enumerate}
\item Fix $m,n \in \mathbb{N}$. Then path topology on $\magtwo^{m+n}$ coincides with the product of path topologies on 
$\magtwo^m$ and on $\magtwo^n$.
\item Fix $n \in \mathbb{N}$. Then path topology on $\magtwo^n$ coincides with the final topology of $\magtwo^n$ with respect to all smooth paths in $\magtwo^n$.
\item $\comptwo:\magtwo\times \magtwo \rightarrow \magtwo$ is a continuous function.
\item Fix $n \in \mathbb{N}$. Then the set function $\lboundtwo \quad\rboundtwo : \magtwo^n \rightarrow \magtwo$ introduced in Proposition \ref{alprtwo}-[1] is a continuous function. 
\item The set function $\acmtwo:\bsfM \times \magtwo \rightarrow \magtwo$ introduced in Proposition \ref{alprtwo}-[2] is a continuous function.
\item Fix $\bsfm \in\bsfM$, $\mtx,\mty \in \magtwo$. Then:
\begin{flalign}
&\mtx\comptwo \mty \notin \magtwoempty\quad \Rightarrow \quad \mty \notin \magtwoempty\text{,}\label{compempty}\\[4pt]
&\lboundtwo \mtx, \mty  \rboundtwo  \notin \magtwoempty\quad \Rightarrow \quad \left(\mtx \notin \magtwoempty\quad \text{and}\quad \mty \notin \magtwoempty\right)\text{,}\label{cartprodempty}\\[4pt]
&\bsfm\acmtwo \mtx  \notin \magtwoempty\quad \Leftrightarrow \quad \mtx \notin \magtwoempty\nonumber\text{.}
\end{flalign}
\end{enumerate}
\end{proposition}    
\begin{proof}\mbox{}\newline
\textnormal{\textbf{Proof of statement 1.}}\ \ We denote by $\toppathmagtwo_k$ the path topology on $\magtwo^{k}$ for any $k \in \mathbb{N}$. Inclusion $\toppathmagtwo_{m+n}\subseteq \toppathmagtwo_m\times \toppathmagtwo_n$ immediately follows by checking that $\proj<\magtwo^{m},\magtwo^n<>1> \funcomp \pmtwo$ is a path in $\magtwo^m$ and $\proj<\magtwo^{m},\magtwo^n<>2>\funcomp \pmtwo$ is a path in $\magtwo^n$ for any fixed path $\pmtwo$ in $\magtwo^{m+n}$.\newline
We are left to prove $ \toppathmagtwo_{m}\times \toppathmagtwo_n\subseteq \toppathmagtwo_{m+n}$. Fix a set function $\pmtwo: (-1,1)\rightarrow \magtwo^{m+n}$. Set $\pmtwo_1=\proj<\magtwo^{m},\magtwo^m<>1> \funcomp \pmtwo$ and $\pmtwo_2=\proj<\magtwo^m,\magtwo^n<>2> \funcomp \pmtwo$. Assume that $\pmtwo_1$ is a path in $\magtwo^m$ and $\pmtwo_2$ is a path in $\magtwo^n$ (i.e. $\pmtwo$ is continuous with respect to Euclidean topology on $(-1,1)$ and to $\toppathmagtwo_m\times \toppathmagtwo_n$). 
Proposition \ref{pathprop1} entails that $(\pmtwo_1,\pmtwo_2)\funcomp \Diag{\left(-1,1\right)}{2}$ is a path in $\magtwo^{m+n}$. Eventually statement follows since $\pmtwo=(\pmtwo_1,\pmtwo_2)\funcomp \Diag{\left(-1,1\right)}{2}$ by construction.\newline
\textnormal{\textbf{Proof of statement 2.}}\ \ We denote by $\toppathmagtwo_n$ the path topology on $\magtwo^{n}$ and by $\smtoppathmagtwo_n$ the final topology on $\magtwo^n$ with respect to the family of all smooth paths in $\magtwo^n$. Inclusion $\toppathmagtwo_n\subseteq \smtoppathmagtwo_n$ is straightforward. To prove $\smtoppathmagtwo_n\subseteq \toppathmagtwo_n$ fix $\openset\in\smtoppathmagtwo_n$ and a path $\pmtwo$ in $\magtwo^n$. Assume $\pmtwo\left(0\right)\in \openset$ and $\openset \notin \toppathmagtwo_n$ by contradiction, then there would be a sequence $\{\unkdue_i\}_{i\in \mathbb{N}}$ fulfilling all conditions below:
\begin{equation} 
0<\unkdue_{i+1}<\unkdue_i<1\quad\forall i\in\mathbb{N}\text{,}\qquad \underset{i\rightarrow \infty}{lim}\,\unkdue_i=0\text{,}\qquad 
\pmtwo\left(\unkdue_i\right)\notin \openset\quad\forall i\in \mathbb{N}\text{.}\label{contr1}
\end{equation}
Set $\singsmpth=\{0\}\cup\{\unkdue_i\}_{i\in \mathbb{N}}$, since the pair $\left(\pmtwo,\singsmpth\right)$ fulfills assumptions \eqref{assclaimg1} then there is a smooth path $\left(\pmtwo_1,\singsmpth\right)$ fulfilling \eqref{claimg1}. Since $\openset \in \smtoppathmagtwo_n$ then $\exists \varepsilon>0\; :\;  \pmtwo_1\left(\unkdue\right)\in\openset \;\; \forall \unkdue\in \left(-\varepsilon, \varepsilon\right)$, hence $\pmtwo\left(\unkdue\right)=\pmtwo_1\left(\unkdue\right)\in \openset \;\; \forall \unkdue \in \singsmpth\cap (-\varepsilon, \varepsilon)$ contradicting \eqref{contr1}.\newline
\textnormal{\textbf{Proof of statements 3.}}\ \ Statement follows since $\comptwo$ sends paths in $\magtwo\times \magtwo$ to paths in $\magtwo$ by Propositions \ref{pathprop1}, \ref{pathprop2}.\newline
\textnormal{\textbf{Proof of statements 4.}}\ \ Statement follows since $\lboundtwo \quad\rboundtwo : \magtwo^n \rightarrow \magtwo$ sends paths in $\magtwo^n$ to paths $\magtwo$  by Proposition \ref{pathprop4}.\newline
\textnormal{\textbf{Proof of statements 5.}}\ \ Statement follows since $\bsfm\acmtwo\pmtwo$ is a path in $\magtwo$ for any $\bsfm \in \bsfM$ and any path $\pmtwo$ in $\magtwo$ by Proposition \ref{pathprop5}.\newline
\textnormal{\textbf{Proof of statement 6.}}\ \ Statement \eqref{compempty} follow by Proposition \ref{pathprop3}. Statement \eqref{cartprodempty} follows by Proposition \ref{pathprop4bis}.\newline
Statement $\quad\bsfm\acmtwo \mtx  \notin \magtwoempty\quad \Leftarrow \quad \mtx \notin \magtwoempty\quad$ can be straightforwardly verified by direct computation.\newline
Statement $\quad\bsfm\acmtwo \mtx  \notin \magtwoempty\quad \Rightarrow \quad \mtx \notin \magtwoempty\quad$ follows by Proposition \ref{pathprop5bis}.
\end{proof}

\begin{remark}
Assumption \eqref{smptass}-[(ii)] is crucial in order to prove Proposition \ref{topinDprop}-[2] and cannot be weakened.
In fact if we replace \eqref{smptass}-[(ii)] by
\begin{equation*}
\begin{array}{ll}
(ii)'&\singsmpth\;\; \text{is a finite set,}
\end{array}
\end{equation*}
then we are able to prove just that an open set for the path topology is open also for the final topology with respect to all smooth paths in $\magtwo^n$, but not the converse.      
\end{remark}

In Definition \ref{pathtwocontdef} below we introduce the notion of path in sets $\magtwocont^n$. 

\begin{definition}\label{pathtwocontdef}\mbox{}
\begin{enumerate}
\item Fix $n \in \mathbb{N}$, a set function $\pmtwo:(-1,1) \rightarrow \magtwocont^n$.\newline
We say that $\pmtwo$ is a path in $\magtwocont^n$ if and only if both conditions below are fulfilled:
\begin{flalign*}
&\incl{\magtwocont^n}{\magtwo^n}\funcomp \pmtwo\;\text{is a path in}\; \magtwo^n\text{;}
\\[8pt]
&\text{there is a core}\; \left(\fmx_1,...,\fmx_n\right)\; \text{of}\; \incl{\magtwocont^n}{\magtwo^n}\funcomp\pmtwo\;
\text{with}\;\left(\fmx_1,...,\fmx_n\right)\in \fremagcont^n\text{.}
\end{flalign*}
With an abuse of language we say that: any skeleton $\left(\left(\fmx_1,...,\fmx_n\right),\assfun\right)$ of $\incl{\magtwocont^n}{\magtwo^n}\funcomp \pmtwo$ is a skeleton of $\pmtwo$; $\left(\fmx_1,...,\fmx_n\right)$ is a core of $\pmtwo$; $\assfun$ is an associating function of $\pmtwo$.
\item Fix $n \in \mathbb{N}$, a path $\pmtwo$ in $\magtwocont^n$, a subset $\singsmpth \subseteq \left(-1,1\right)$.\newline
We say that $\left(\pmtwo,\singsmpth\right)$ is a smooth path in $\magtwocont^n$ if and only if the pair defined by $\left(\incl{\magtwocont^n}{\magtwo^n}\funcomp\pmtwo,\singsmpth\right)$ is a smooth path in $\magtwo^n$.\newline
With an abuse of language we say that: $\singsmpth$ is the singular set of
$\left(\pmtwo, \singsmpth\right)$; any skeleton $\left(\left(\fmx_1,...,\fmx_n\right)\!,\assfun\right)$ of $\left(\!\incl{\magtwocont^n}{\magtwo^n}\!\funcomp\!\pmtwo,\singsmpth\right)$ is a skeleton of  $\left(\!\pmtwo,\singsmpth\right)$.\newline
We emphasize that nothing is assumed about $\pmtwo\left(\unkdue\right)$ when $\unkdue \in \singsmpth$.
\item Fix $n \in \mathbb{N}$, $\left(\mtx_1,...,\mtx_n\right)\in \magtwocont^n$, a path $\pmtwo$ in $\magtwocont^n$.\newline
We say that $\pmtwo$ is a path in $\magtwocont^n$ through $(\mtx_1,...,\mtx_n)$ if and only if we have $\pmtwo\left(0\right)=\left(\mtx_1,...,\mtx_n\right)$.
\end{enumerate}
\end{definition}

In Definition \ref{topmagtwocont} below we give topology to sets $\magtwocont^n$ through paths introduced in Definition \ref{pathtwocontdef}. We refer to Notation \ref{gentopnot}-[5].

\begin{definition}\label{topmagtwocont} Fix $n \in \mathbb{N}$. Set $\magtwocont^n$ is endowed with the final topology with respect to the family of all paths in $\magtwocont^n$. We say that such topology is the path topology on $\magtwocont^n$.
\end{definition}

In Proposition \ref{topinDcontprop} below we state and prove basic properties of the path topology introduced in Definition \ref{topmagtwocont}. We refer to Notation \ref{gentopnot}-[4, 5], Definition \ref{intdiffmon}.

\begin{proposition}\label{topinDcontprop}\mbox{} 
\begin{enumerate}
\item Fix $m,n \in \mathbb{N}$. Then path topology on $\magtwocont^{m+n}$ coincides with the product of path topologies on $\magtwocont^m$ and on $\magtwocont^n$.  
\item Fix $n \in \mathbb{N}$. Then path topology on $\magtwocont^n$ coincides with the final topology on $\magtwocont^n$ with respect to the family of all smooth paths in $\magtwocont^n$.
\item Inclusion $\incl{\magtwocont}{\magtwo}$ is a continuous function.
\item There is one and only one continuous function $\comptwocont:\magtwocont\times \magtwocont \rightarrow \magtwocont$ fulfilling $ \comptwo \funcomp \incl{\magtwocont}{\magtwo}^2= \incl{\magtwocont}{\magtwo} \funcomp \comptwocont$.\newline
With an abuse of language, from now on, we denote $\comptwocont$ by $\comptwo$. 
\item Fix $n \in \mathbb{N}$. There is one and only one continuous function $\lboundtwocont\;\rboundtwocont:\magtwocont^n \rightarrow \magtwocont$ defined by setting $ \lboundtwocont\;\rboundtwocont \funcomp \incl{\magtwocont}{\magtwo}^n= \incl{\magtwocont}{\magtwo} \funcomp \lboundtwocont\;\rboundtwocont$. \newline
With an abuse of language, from now on, we denote $\lboundtwocont\;\rboundtwocont$ by $\lboundtwo\;\rboundtwo$. 
\item There is one and only one continuous function $\acmtwocont:\subbsfM \times \magtwocont \rightarrow \magtwocont$ defined by setting $ \acmtwo \funcomp \left(\idobj+\subbsfM+, \incl{\magtwocont}{\magtwo}\right)= \incl{\magtwocont}{\magtwo} \funcomp \acmtwocont$.\newline
With an abuse of language, from now on, we denote $\acmtwocont$ by $\acmtwo$.
\item Fix a path $\pmtwo$ in $\magtwocont$, $\unkdue_0 \in \left(-1,1\right)$, $\intdue\subseteq \left(-1,1\right)$.\newline 
If $\evalcomptwo\left(\pmtwo\left(t_0\right)\right)=\contspempty$ then $\evalcomptwo \funcomp \pmtwo$ is a continuous function at $\unkdue=\unkdue_0$.\newline
If $\evalcomptwo\left(\pmtwo\left(\unkdue\right)\right)\neq\contspempty \quad \forall \unkdue\in \intdue$ then $\left(\evalcomptwo \funcomp \pmtwo\right)\funcomp \incl{\intdue}{\left(-1,1\right)}$ is a continuous function where $\intdue$ is endowed with the subspace topology with respect to $\left(-1,1\right)$.
\item The set function $\linincltwo: \Cksp{0}\setminus \{\contspempty \}\rightarrow \magtwo$ fulfills the system of conditions below
\begin{equation}\label{proplinincltwo}
\left\{
\begin{array}{lll}
(i)& \linincltwo\left(\Cksp{0}  \setminus \{\contspempty \}\right)\subseteq \magtwocontcont\setminus\magtwoempty \text{,}\vspace{4pt} \\
(ii)& \evalcomptwo \funcomp \linincltwo = \idobj+\Cksp{0} \setminus \{\emptysetfun\}+\text{,}\vspace{4pt} \\
(iii)&\linincltwo\; \text{is a continuous function}\\
&\text{(} \Cksp{0}\setminus \{\contspempty \} \;\text{is endowed with the subspace topology}\\
&\text{with respect to}\; \Cksp{0}\text{).}
\end{array}
\right.
\end{equation}
\item Fix $\mtx\in \magtwocont$, a path $\pmtwo$ in $\magtwo$ detecting $\mtx$. Then $\pmtwo$ is a path in $\magtwocont$.
\end{enumerate}
\end{proposition}    
\begin{proof}\mbox{}\newline
\textnormal{\textbf{Proof of statement 1.}}\ \ The proof coincides word by word with the proof of Proposition \ref{topinDprop}-[1] by replacing $\fremag$ by $\fremagcont$, $\magtwo$ by $\magtwocont$ and by choosing skeleta $\left(\left(\mtx_1,...,\mtx_m\right), \assfun_1\right)$ with $\left(\mtx_1,...,\mtx_m\right)\in \fremagcont^m$, $\left(\left(\mty_1,...,\mty_n\right), \assfun_2\right)$ with $\left(\mty_1,...,\mty_n\right)\in \fremagcont^n$.\newline
\textnormal{\textbf{Proof of statement 2.}}\ \ The proof coincides word by word with the proof of Proposition \ref{topinDprop}-[2] by replacing $\fremag$ by $\fremagcont$, $\magtwo$ by $\magtwocont$.\newline
\textnormal{\textbf{Proof of statement 3.}}\ \ Statement follows since any path in $\magtwocont$ is also a path in $\magtwo$.\newline
\textnormal{\textbf{Proof of statement 4.}}\ \ Statement follows by Propositions \ref{alprtwo}-[5], \ref{pathprop2}-[3].\newline
\textnormal{\textbf{Proof of statement 5.}}\ \ Statement follows by Proposition \ref{pathprop4}-[3].\newline
\textnormal{\textbf{Proof of statement 6.}}\ \ Statement follows by Proposition \ref{pathprop5}-[3].\newline
\textnormal{\textbf{Proof of statement 7.}}\ \ If $\evalcomptwo\left( \pmtwo\left(\unkdue_0\right)\right)= \contspempty$ then continuity of $\evalcomptwo\funcomp \pmtwo$ at $\unkdue_0$ follows since $\Cksp{0}$ is the only open set belonging to the topology $\topC$ and containing $\contspempty$. If $\evalcomptwo\left(\pmtwo\left(\unkdue\right)\right)\neq\contspempty \quad \forall \unkdue\in \intdue$ then statement is straightforward if there is a core $\fmx$ of $\pmtwo$ with $\occset_{\fmx}=\udenset$, otherwise statement follows by Proposition \ref{funoncont}.\newline
\textnormal{\textbf{Proof of statement 8.}}\ \ Statement \eqref{proplinincltwo}-[$(i)$] is proved by a two steps argument. The first part of the statement, say $\linincltwo\left(\Cksp{0}  \setminus \{\contspempty \}\right)\subseteq \magtwocontcont$, follows by \eqref{proplininclone}-[$(ii)$]. Then we prove $\linincltwo\left(\Cksp{0}  \setminus \{\contspempty \}\right)\cap \magtwoempty=\udenset$ by arguing as follows. Fix $f \in \Cksp{0}\setminus \left\{\contspempty\right\}$. By Weierstrass approximation theorem there is at least one associating function $\assfun$ for $\occset_{\lininclone\left(f\right)}$ fulfilling both conditions below: 
\begin{flalign*}
&\begin{array}{l} \assfun\left(g,\unkdue\right)\neq \contspempty\hspace{20pt}\forall \left(g,\unkdue\right) \in \occset_{\lininclone\left(f\right)}\times \left(-1,1\right)\text{,}
\end{array}\\[8pt]
&\begin{array}{l} \assfun\left(g,0\right) =g\hspace{20pt}\forall g \in \occset_{\lininclone\left(f\right)}\text{.}
\end{array}
\end{flalign*}
Then construction of $\lininclone$ in Proposition \ref{lininclbis}-[3] entails that pair $\left(\lininclone\left(f\right),\assfun\right)$ defines through \eqref{pathconduro2} a path $\pmtwo$ in $\magtwo$ detecting $\linincltwo\left(f\right)$ with\newline
\centerline{$\pmtwo\left(\unkdue\right)\notin \magtwoempty \hspace{20pt} \forall \unkdue \in \left(-1,1\right)\setminus \left\{0\right\}$.}\newline
Statement \eqref{proplinincltwo}-[$(ii)$] follows by \eqref{proplininclone}-[$(iii)$].\newline
Statement \eqref{proplinincltwo}-[$(iii)$] is proved as follows. Fix a path $\pmtwo$ in $\Cksp{0}$. By \eqref{extopcntdis1} the set $\pmtwo^{-1}\left(\Cksp{0}\setminus\{\contspempty\}\right)$ is an open subset of $\left(-1,1\right)$ then there is no loss of generality by assuming that\hspace{20pt}$\pmtwo\left(\unkdue\right) \neq \contspempty \qquad \forall \unkdue\in \left(-1,1\right)$.\newline
Eventually statement follows since $\linincltwo\funcomp\pmtwo $ is a path in $\magtwocont$ by definition.\newline
\textnormal{\textbf{Proof of statement 9.}}\ \ Statement follows straightforwardly by \eqref{eqelc}, Definitions \ref{pathtwodef}-[4], \ref{pathtwocontdef}-[1]. 
\end{proof}

In Definition \ref{pathtwocontcontdef} below we introduce the notion of path in sets $\magtwocontcont^n$.

\begin{definition}\label{pathtwocontcontdef}\mbox{}
\begin{enumerate}
\item Fix $n \in \mathbb{N}$, a set function $\pmtwo:(-1,1) \rightarrow \magtwocontcont^n$.\newline
We say that $\pmtwo$ is a path in $\magtwocontcont^n$ if and only if $\incl{\magtwocontcont^n}{\magtwocont^n}\funcomp \pmtwo$ is a path in $\magtwocont^n$.\newline
With an abuse of language we say that: any skeleton $\left(\left(\fmx_1,...,	\fmx_n\right),\assfun\right)$ of $\incl{\magtwocontcont^n}{\magtwocont^n}\funcomp \pmtwo$ is a skeleton of $\pmtwo$; $\left(\fmx_1,...,\fmx_n\right)$ is a core of $\pmtwo$; $\assfun$ is an associating function of $\pmtwo$.
\item Fix $n \in \mathbb{N}$, a path $\pmtwo$ in $\magtwocontcont^n$, a subset $\singsmpth \subseteq \left(-1,1\right)$.\newline
We say that $\left(\pmtwo,\singsmpth\right)$ is a smooth path in $\magtwocontcont^n$ if and only if the pair defined by $\left(\incl{\magtwocontcont^n}{\magtwocont^n}\funcomp\pmtwo,\singsmpth\right)$ is a smooth path in $\magtwocont^n$.\newline
With an abuse of language we say that: $\singsmpth$ is the singular set of $\left(\pmtwo, \singsmpth\right)$; any skeleton $\left(\left(\fmx_1,...,\fmx_n\right),\assfun\right)$ of $\left(\incl{\magtwocontcont^n}{\magtwocont^n}\funcomp\pmtwo,\singsmpth\right)$ is a skeleton of $\left(\pmtwo,\singsmpth\right)$.\newline
We emphasize that nothing is assumed about $\pmtwo\left(\unkdue\right)$ when $\unkdue \in \singsmpth$.
\item Fix $n \in \mathbb{N}$, $\left(\mtx_1,...,\mtx_n\right)\in \magtwocontcont^n$, a path $\pmtwo$ in $\magtwocontcont^n$.\newline
We say that $\pmtwo$ is a path in $\magtwocontcont^n$ through $\left(\mtx_1,...,\mtx_n\right)$ if and only if we have $\pmtwo\left(0\right)=\left(\mtx_1,...,\mtx_n\right)$.
\end{enumerate}
\end{definition}

In Definition \ref{topmagtwocontcont} below we give topology to sets $\magtwocontcont^n$ through paths introduced in Definition \ref{pathtwocontcontdef}. We refer to  Notation \ref{gentopnot}-[5].

\begin{definition}\label{topmagtwocontcont} Fix $n \in \mathbb{N}$. Set $\magtwocontcont^n$ is endowed with the final topology with respect to the family of all paths in $\magtwocontcont^n$. We say that such topology is the the path topology on $\magtwocontcont^n$.
\end{definition}

In Proposition \ref{topinDcontcontprop} below we state and prove basic properties of the path topology introduced in Definition \ref{topmagtwocontcont}. We refer to Notation \ref{gentopnot}-[4], Definition \ref{intdiffmon}.

\begin{proposition}\label{topinDcontcontprop}\mbox{} 
\begin{enumerate}
\item Fix $m,n \in \mathbb{N}$. Then path topology on $\magtwocontcont^{m+n}$ coincides with the product of paths topologies on $\magtwocontcont^m$ and on $\magtwocontcont^n$.  
\item Inclusion $\incl{\magtwocontcont}{\magtwocont}$ is a continuous function.
\item Fix paths $\pmtwo_1$, $\pmtwo_2$ in $\magtwocontcont$. Assume that 
\begin{equation*}
\pmtwo_1\left(t\right) \comptwo \pmtwo_2\left(t\right) \notin \magtwoempty \qquad \forall t \in \left(-1,1\right)\text{.}
\end{equation*} 
Then the set function $\pmtwo:\left(-1,1\right)\rightarrow\magtwocontcont$ defined by setting
\begin{equation*}
\pmtwo\left(t\right)=\pmtwo_1\left(t\right) \comptwo \pmtwo_2\left(t\right)\qquad \forall t \in \left(-1,1\right),
\end{equation*}
is a path in $\magtwocontcont$.
\item Fix $n \in \mathbb{N}$, paths $\pmtwo_1$, ..., $\pmtwo_n$ in $\magtwocontcont$.\newline
Then the set function $\pmtwo:\left(-1,1\right)\rightarrow\magtwocontcont$ defined by setting
\begin{equation*}
\pmtwo\left(t\right)=\lboundtwo\pmtwo_1\left(t\right),.., \pmtwo_n\left(t\right)\rboundtwo\qquad \forall t \in \left(-1,1\right),
\end{equation*}
is a path in $\magtwocontcont$.
\item  Fix $\bsfm \in\subbsfM$, a path $\pmtwo$ in $\magtwocontcont$.\newline
Then the set function $\pmtwo_1:\left(-1,1\right)\rightarrow\magtwocontcont$ defined by setting
\begin{equation*}
\pmtwo_1\left(t\right)=\bsfm\acmtwo\pmtwo\left(t\right)\qquad \forall t \in \left(-1,1\right),
\end{equation*}
is a path in $\magtwocontcont$.
\end{enumerate}
\end{proposition}
\begin{proof}\mbox{}\newline
\textnormal{\textbf{Proof of statement 1.}}\ \ The proof coincides word by word with the proof of Proposition \ref{topinDprop}-[1] by replacing $\fremag$ by $\fremagcontcont$, $\magtwo$ by $\magtwocontcont$ and by choosing skeletons $\left(\left(\mtx_1,...,\mtx_m\right), \assfun_1\right)$ with $\left(\mtx_1,...,\mtx_m\right)\in \fremagcont^m$, $\left(\left(\mty_1,...,\mty_n\right), \assfun_2\right)$ with $\left(\mty_1,...,\mty_n\right)\in \fremagcont^n$.\newline
\textnormal{\textbf{Proof of statement 2.}}\ \ Statement follows since any path in $\magtwocontcont$ is also a path in $\magtwocont$.\newline
\textnormal{\textbf{Proof of statement 3.}}\ \ Statement follows by \eqref{propevaltwo}-[$(i)$], Proposition \ref{topinDcontprop}-[4].\newline
\textnormal{\textbf{Proof of statement 4.}}\ \ Statement follows by \eqref{propevaltwo}-[$(ii)$], Proposition \ref{topinDcontprop}-[5].\newline
\textnormal{\textbf{Proof of statement 5.}}\ \ Statement follows by \eqref{propevaltwo}-[$(iii)$], Proposition \ref{topinDcontprop}-[6].
\end{proof}

In Definition \ref{pathtwosmoothdef} below we introduce the notion of path in sets $\magtwosmooth^n$.

\begin{definition}\label{pathtwosmoothdef}\mbox{}
\begin{enumerate}
\item Fix $n \in \mathbb{N}$, a set function $\pmtwo:\left(-1,1\right) \rightarrow \magtwosmooth^n$.\newline
We say that $\pmtwo$ is a path in $\magtwosmooth^n$ if and only if $\incl{\magtwosmooth^n}{\magtwo^n}\funcomp\pmtwo$ is a path in $\magtwo^n$.\newline
With an abuse of language we say that: any skeleton $\left(\left(\fmx_1,...,\fmx_n\right),\assfun\right)$ of $\incl{\magtwosmooth^n}{\magtwo^n}\funcomp\pmtwo$ is a skeleton of $\pmtwo$; $\left(\fmx_1,...,\fmx_n\right)$ is a core of $\pmtwo$; $\assfun$ is an associating function of $\pmtwo$.
\item Fix $n \in \mathbb{N}$, a path $\pmtwo$ in $\magtwosmooth^n$, a subset $\singsmpth \subseteq \left(-1,1\right)$.\newline
We say that $\left(\pmtwo,\singsmpth\right)$ is a smooth path in $\magtwosmooth^n$ if and only if the pair defined by $\left(\incl{\magtwosmooth^n}{\magtwo^n}\funcomp\pmtwo,\singsmpth\right)$ is a smooth path in $\magtwo^n$.\newline
With an abuse of language we say that: $\singsmpth$ is the singular set of $\left(\pmtwo, \singsmpth\right)$; any skeleton $\left(\left(\fmx_1,...,\fmx_n\right)\!,\assfun\right)$ of $\left(\incl{\magtwosmooth^n}{\magtwo^n}\funcomp\pmtwo,\singsmpth\right)$ is a skeleton of $\left(\pmtwo,\singsmpth\right)$.\newline
We emphasize that nothing is assumed about $\pmtwo\left(\unkdue\right)$ when $\unkdue \in \singsmpth$.
\item Fix $n \in \mathbb{N}$, $\left(\mtx_1,...,\mtx_n\right)\in \magtwosmooth^n$, a path $\pmtwo$ in  $\magtwosmooth^n$.\newline
We say that $\pmtwo$ is a path in $\magtwosmooth^n$ through $\left(\mtx_1,...,\mtx_n\right)$ if and only if we have $\pmtwo\left(0\right)=\left(\mtx_1,...,\mtx_n\right)$.
\end{enumerate}
\end{definition}

In Proposition \ref{smpathtocont} below we prove that any path in $\magtwosmooth^n$ is a path in $\magtwocont^n$.

\begin{proposition}\label{smpathtocont}
Fix $n \in \mathbb{N}$, a path $\pmtwo$ in $\magtwosmooth^n$.
\begin{enumerate}
\item Then $\incl{\magtwosmooth^n}{\magtwocont^n}\funcomp \pmtwo$ is a path in $\magtwocont^n$.
\item If there is a subset $\singsmpth \subseteq \left(-1,1\right)$ such that the pair $\left(\pmtwo,\singsmpth\right)$ is a smooth path in $\magtwosmooth^n$ then pair $\left(\incl{\magtwosmooth^n}{\magtwocont^n}\funcomp\pmtwo,\singsmpth \right)$ is a smooth path in $\magtwocont^n$.
\end{enumerate}
\end{proposition}
\begin{proof}\mbox{}\newline
\textnormal{\textbf{Proof of statement 1.}}\ \ Choose a skeleton $\left(\left(\fmz_1,...,\fmz_n\right), \assfundue\right)$ of $\pmtwo$. Denote by $\pmone$ the path in $\fremag$ defined by $\left(\left(\fmz_1,...,\fmz_n\right), \assfundue\right)$ through \eqref{pathconduro2}.Referring to Definition \ref{intdiffmon}, relations \eqref{derint}, \eqref{coordint}-[$(iii)$], \eqref{leftproj}-[$(i)$], \eqref{rightproj}-[$(i)$], \eqref{(R 9.1)}, \eqref{(R 11)}, \eqref{(R 12)}, \eqref{(R 17)}, \eqref{(R 17 bis)} entail that there are 
\begin{flalign*}
&\text{an element}\;\left(\fmx_1,...,\fmx_n\right) \in \fremagcont^n\text{,}\\[4pt]
&\text{a set function}\;\varphi:\overset{n}{\underset{\mathsf{n}=1}{\bigcup}}\,\occset^{\fmx_{\mathsf{n}}} \rightarrow \subbdiffsfM\times\left(\overset{n}{\underset{\mathsf{n}=1}{\bigcup}}\,\occset_{\fmz_{\mathsf{n}}}\right)\text{,}\\[4pt]
&\text{a set function}\;\vartheta :\overset{n}{\underset{\mathsf{n}=1}{\bigcup}}\,\occset_{\fmx_{\mathsf{n}}}\rightarrow \overset{n}{\underset{\mathsf{n}=1}{\bigcup}}\,\occset^{\fmx_{\mathsf{n}}}
\end{flalign*}
fulfilling all conditions below:
\begin{flalign*}
& \vartheta \;\text{is one-to-one}\text{,}\\[8pt]
&\left(\incllim_{1}\!\funcomp\! \incl{\subbdiffsfM\!\times\! \left(\overset{n}{\underset{\mathsf{n}=1}{\bigcup}}\,\occset_{\fmz_{\mathsf{n}}}\right)}{\fremag_{1}}\right)\!\funcomp \!\varphi\;\text{is an occurrence function on}\;\overset{n}{\underset{\mathsf{n}=1}{\bigcup}}\,\occset^{\fmx_\mathsf{n}}\text{,}\\[8pt]
&\quotmagone^n\left(\fmz_1,...,\fmz_n\right)=\quotmagone^n\left(\strchar\!\left[\left(\fmx_1,...,\fmx_n\right),  \left(\incllim_{1}\!\funcomp\! \incl{\subbdiffsfM\!\times\!\left(\overset{n}{\underset{\mathsf{n}=1}{\bigcup}}\,\occset_{\fmz_{\mathsf{n}}}\right)
}{\fremag_{1}}\right)\!\funcomp\! \varphi \right]\right)\text{.}
\end{flalign*}
Since \eqref{pathconduro}-[$(ii)$] holds true and $\pmone\left(\left(-1,1\right)\right)\subseteq \fremagsmooth$, we have that\newline
\centerline{$\assfundue\left(f,\unkdue\right)\in \Cksp{\infty}\setminus\left\{\smospempty\right\} \qquad \forall \left(f,\unkdue\right)\in\left(\overset{n}{\underset{\mathsf{n}=1}{\bigcup}}\,\occset_{\fmz_{\mathsf{n}}}\right)\times\left(-1,1\right)$,}\newline
then we define the set function $\assfun:\left(\overset{n}{\underset{\mathsf{n}=1}{\bigcup}}\,\occset_{\fmx_{\mathsf{n}}}\right)\times\left(-1,1\right)\rightarrow \Cksp{\infty}$ by setting
\begin{multline*}
\assfun\left(f,\unkdue\right)=\proj<\subbdiffsfM, \left(\overset{n}{\underset{\mathsf{n}=1}{\bigcup}}\,\occset_{\fmz_{\mathsf{n}}}\right)<>1>\left(\varphi\left(\vartheta\left(f\right)\right)\right)\acmone\\[4pt]\assfundue\left(\proj<\subbdiffsfM,\left(\overset{n}{\underset{\mathsf{n}=1}{\bigcup}}\,\occset_{\fmz_{\mathsf{n}}}\right)<>2>\left(\varphi\left(\vartheta\left(f\right)   \right)\right),\unkdue\right)\\[4pt]
\forall \left(f,\unkdue\right)\in \left(\overset{n}{\underset{\mathsf{n}=1}{\bigcup}}\,\occset_{\fmx_{\mathsf{n}}}\right)\times\left(-1,1\right)\text{.}
\end{multline*}
Eventually statement follows directly by checking that the pair $\left(\left(\fmx_1,...,\fmx_n\right),\assfun\right)$ is a skeleton of $\pmtwo$.\newline
\textnormal{\textbf{Proof of statement 2.}}\ \ The proof of coincides word by word with the proof of statement 1 by replacing path $\pmtwo$ by the smooth path $\left(\pmtwo,\singsmpth\right)$.
\end{proof}

In Definition \ref{topmagtwosmooth} below we give topology to sets $\magtwosmooth^n$ through paths introduced in Definition \ref{pathtwosmoothdef}. We refer to Notation \ref{gentopnot}-[5].

\begin{definition}\label{topmagtwosmooth} Fix $n \in \mathbb{N}$. Set $\magtwosmooth^n$ is endowed with the final topology with respect to the family of all paths in $\magtwosmooth^n$. We say that such topology is the path topology on $\magtwosmooth^n$.
\end{definition}

In Proposition \ref{topinDsmoothprop} below we state and prove basic properties of the path topology introduced in Definition \ref{topmagtwosmooth}. We refer to Notation \ref{gentopnot}-[4, 5], Definition \ref{intdiffmon}.

\begin{proposition}\label{topinDsmoothprop}\mbox{}
\begin{enumerate}
\item Fix $m,n \in \mathbb{N}$. Then path topology on $\magtwosmooth^{m+n}$ coincides with the product of path topologies on $\magtwosmooth^{m}$ and on $\magtwosmooth^{n}$. 
\item Fix $n \in \mathbb{N}$. Then path topology on $\magtwosmooth^{n}$ coincides with the final topology on $\magtwosmooth^{n}$ with respect to the family of all smooth functions in $\magtwosmooth^{n}$.
\item Inclusion $\incl{\magtwosmooth}{\magtwocont}$ is a continuous function.
\item There is one and only one continuous function $\comptwosmooth:\magtwosmooth\times \magtwosmooth \rightarrow \magtwosmooth$ defined by setting $ \comptwo \funcomp \incl{\magtwosmooth}{\magtwocont}^2= \incl{\magtwosmooth}{\magtwocont} \funcomp \comptwosmooth$. \newline
With an abuse of language, from now on, we denote $\comptwosmooth$ by $\comptwo$. 
\item Fix $n \in \mathbb{N}$. There is one and only one continuous function $\lboundtwosmooth\;\rboundtwosmooth:\magtwosmooth^n \rightarrow \magtwosmooth$ defined by setting $ \lboundtwo\;\rboundtwo \funcomp \incl{\magtwosmooth}{\magtwocont}^n= \incl{\magtwosmooth}{\magtwocont} \funcomp \lboundtwosmooth\;\rboundtwosmooth$. \newline
With an abuse of language, from now on, we denote $\lboundtwosmooth\;\rboundtwosmooth$ by $\lboundtwo\;\rboundtwo$.
\item There is one and only one continuous function $\acmtwosmooth:\bsfM \times \magtwosmooth \rightarrow \magtwosmooth$ defined by setting $ \acmtwo \funcomp \left(\idobj+\bsfM+, \incl{\magtwosmooth}{\magtwo}\right)= \incl{\magtwosmooth}{\magtwo} \funcomp \acmtwosmooth$.\newline
With an abuse of language, from now on, we denote $\acmtwosmooth$ by $\acmtwo$.
\end{enumerate}
\end{proposition}    
\begin{proof}\mbox{}\newline
\textnormal{\textbf{Proof of statement 1.}}\ \ The proof coincides word by word with the proof of Proposition \ref{topinDprop}-[1] by replacing $\fremag$ by $\fremagsmooth$, $\magtwo$ by $\magtwosmooth$.\newline
\textnormal{\textbf{Proof of statement 2.}}\ \ The proof coincides word by word with the proof of Proposition \ref{topinDprop}-[2]  by replacing $\fremag$ by $\fremagcont$, $\magtwo$ by $\magtwosmooth$.\newline
\textnormal{\textbf{Proof of statement 3.}}\ \ Statement follows by Proposition \ref{smpathtocont} which entails that any path in $\magtwosmooth$ is also a path in $\magtwocont$.\newline
\textnormal{\textbf{Proof of statement 4.}}\ \ Statement follows by Propositions \ref{alprtwo}-[5], \ref{pathprop2}, statement 3.\newline
\textnormal{\textbf{Proof of statement 5.}}\ \ Statement follows by Propositions \ref{lininclbis}-[4], \ref{alprtwo}-[1], \ref{pathprop4}, statement 3.\newline
\textnormal{\textbf{Proof of statement 6.}}\ \ Statement follows by Propositions \ref{lininclbis}-[4], \ref{alprtwo}-[2], \ref{pathprop5}-[3], \ref{smpathtocont}.
\end{proof}

\section{The integro-differential space $\genfidsp$}
We define the integro-differential space $\genfidsp$ by contribution of both algebraic structure and topology of magma 
$\left(\magtwo,\comptwo\right)$.\vspace{10pt}

In Definition \ref{defgenf} we introduce the missing ingredients for the construction of $\genfidsp$. We refer to Notations \ref{ins}-[7, 8], \ref{realfunc}-[3(a), 4].

\begin{definition}\label{defgenf}We define: 
\begin{enumerate}
\item The set $\genfzero=\left\{\cost<\dommagtwo\left(\mtx\right)<>\codmagtwo\left(\mtx\right)>+0+\comptwo \mtx \,:\mtx  \in \magtwo\setminus\magtwoempty \right\}$.
\item The set function $\genfunsum:\magtwo\times \magtwo\rightarrow \magtwo$ by setting
\begin{equation*}
\mtx_1\genfunsum \mtx_2 = \vecsum<n<>2>\comptwo \left(\lboundtwo \mtx_1,\mtx_2\rboundtwo\comptwo\Diag{\intuno}{2}\right)
\end{equation*}
where 
\begin{equation*}
 n=\max\left\{n_1,n_2\right\}\text{,}\hspace{20pt}
\intuno=\left\{
\begin{array}{ll}
\dommagtwo\left(\mtx_1\right)&\text{if}\;m_1>m_2\text{,}\\
\dommagtwo\left(\mtx_2\right)&\text{if}\;m_1<m_2\text{,}\\
\dommagtwo\left(\mtx_1\right)\cup\dommagtwo\left(\mtx_2\right)&\text{if}\;m_1=m_2\text{,}
\end{array} \right.
\end{equation*}
for any $m_1,m_2,n_1,n_2\in \mathbb{N}_0$, $\mtx_1,\mtx_2\in \magtwo$ fulfilling all conditions below: 
\begin{equation*}
\codmagtwo\left(\mtx_1\right)=\mathbb{R}^{n_1}\text{,}\quad \codmagtwo\left(\mtx_2\right)=\mathbb{R}^{n_2}\text{,}\quad\dommagtwo\left(\mtx_1\right)\subseteqdentro\mathbb{R}^{m_1}\text{,}\quad\dommagtwo\left(\mtx_1\right)\subseteqdentro\mathbb{R}^{m_1}\text{.}
\end{equation*}
\item The set function $\genfunmult:\magtwo\times \magtwo\rightarrow \magtwo$ by setting
\begin{equation*}
\mtx_1\genfunmult \mtx_2 = \vecprod\left[n_1,n_2\right]\comptwo \left(\lboundtwo \mtx_1,\mtx_2\rboundtwo\comptwo\Diag{\intuno}{2}\right)
\end{equation*}
where 
\begin{equation*}
 \intuno=\left\{
\begin{array}{ll}
\dommagtwo\left(\mtx_1\right)&\text{if}\;m_1>m_2\text{,}\\
\dommagtwo\left(\mtx_2\right)&\text{if}\;m_1<m_2\text{,}\\
\dommagtwo\left(\mtx_1\right)\cup\dommagtwo\left(\mtx_2\right)&\text{if}\;m_1=m_2\text{,}
\end{array} \right.
\end{equation*}
for any $m_1,m_2,n_1,n_2\in \mathbb{N}_0$, $\mtx_1,\mtx_2\in \magtwo$ fulfilling all conditions below: 
\begin{equation*}
\codmagtwo\left(\mtx_1\right)=\mathbb{R}^{n_1}\text{,}\quad \codmagtwo\left(\mtx_2\right)=\mathbb{R}^{n_2}\text{,}\quad\dommagtwo\left(\mtx_1\right)\subseteqdentro\mathbb{R}^{m_1}\text{,}\quad\dommagtwo\left(\mtx_1\right)\subseteqdentro\mathbb{R}^{m_1}\text{.}
\end{equation*}
\item The set function $\genfunscalp:\mathbb{R}\times \magtwo\rightarrow \magtwo$ by setting
\begin{equation*}
\unkuno \genfunscalp  \mtx= \vecprod\left[1,n\right] \comptwo\left(\lboundtwo\cost<\dommagtwo\left(\mtx\right)<>\mathbb{R}>+\unkuno+,\mtx\rboundtwo\comptwo \Diag{\dommagtwo\left(\mtx\right)}{2}\right)
\end{equation*}
for any $n\in \mathbb{N}_0$, $\unkuno\in \mathbb{R}$, $\mtx\in \magtwo$ with $\codmagtwo\left(\mtx\right)=\mathbb{R}^{n}$.
\end{enumerate}
\end{definition}

\begin{theorem}\label{genstrth}\mbox{}
\begin{enumerate}
\item $\genfidsp=\left(\genf, \genfzero, \genfempty, \genfuncomp,  \lboundgenf\,\rboundgenf   ,\genfunsum , \genfunmult , \genfunscalp , \genfbsfmu\right)$ is an integro-differential space.
\item $\genfidspcont=\left(\genfcont, \genfzero, \genfempty,\genfuncomp,  
\lboundgenf\,\rboundgenf,   \genfunsum , \genfunmult ,  \genfunscalp, \genfbsfmu \right)$ is an integral space.\newline
Inclusion $\incl{\genfcont}{\genf}$ induces a function $\genfidspcont \rightarrow  \ffIDSMaQIDSMa\left(\genfidsp\right)$ between integral spaces.
\item $\genfidspsmooth=\left(\genfsmooth,\genfzero,\genfempty,\genfuncomp,  \lboundgenf\,\rboundgenf,\genfunsum , \genfunmult ,  \genfunscalp,\genfbsfmu\right)$ is an integro-differential space.\newline
Inclusion $\incl{\genfsmooth}{\genfcont}$ induces inclusion $\incl{\ffIDSMaQIDSMa\left(\genfidspsmooth\right)}{\genfidspcont}$ of integral spaces.\newline
Inclusion $\incl{\genfsmooth}{\genf}$ induces inclusion $\incl{\genfidspsmooth}{\genfidsp}$ of integro-differential spaces.
\end{enumerate}
\end{theorem}
\begin{proof}
Statements are achieved by carefully checking axioms in Definitions \ref{DiffSpDef}, \ref{QDiffSpDef}. It is a routine, long but not difficult exercise which can be carried out by exploiting relations \eqref{(R 1)}-\eqref{(R 17 bis)}, Propositions \ref{alprtwo}-[2, 5], \ref{topinDprop}-[6], \ref{topinDcontprop}-[3, 4, 5, 6], \ref{topinDsmoothprop}-[3, 4, 5, 6].\newline
We prove \eqref{+}-[$(iv)$], \eqref{star}-[$(iii)$] whose proof is not straightforward.\newline
\textnormal{\textbf{Proof of}} \eqref{+}-[$(iv)$]\textnormal{\textbf{.}}\ \ Fix $\gfx,\gfy\in \genf$. By definition of $\genfunsum$ in Definition \ref{defgenf} we have that $\gfx \genfunsum \gfy\notin \genfempty$ if and only if all the following conditions are fulfilled: $\gfx,\gfy\in \genf\setminus\genfempty $; $\dommagtwo\left(\gfx\right)=\dommagtwo\left(\gfy\right)$; $\codmagtwo\left(\gfx\right)=\codmagtwo\left(\gfy\right)$. 
We prove that $0_{\gfx} = \cost<\dommagtwo\left(\gfx\right)<>\codmagtwo\left(\gfx\right)>+0+$. 
Choose $n\in\mathbb{N}_0$, $f \in \Cksp{\infty}\left(\dommagtwo\left(\gfx\right), \codmagtwo\left(\gfx\right)\right)$ with $\codmagtwo\left(\gfy\right)=\mathbb{R}^n$. Then the following chain of equalities holds true 
\begin{equation*}
\begin{array}{l}
\gfy \genfunsum \cost<\dommagtwo\left(\gfx\right)<>\codmagtwo\left(\gfx\right)>+0+\overset{(1)}{=}\\[4pt]
  \vecsum<n<>2>\comptwo \left(\lboundtwo \gfy, \cost<\dommagtwo\left(\gfx\right)<>\codmagtwo\left(	\gfx\right)>+0+ \rboundtwo\comptwo \Diag{\dommagtwo\left(\gfy\right)}{2} \right)\overset{(2)}{=}\\[4pt]
 \vecsum<n<>2>\comptwo \left(\left(
\lboundtwo \idobj+\codmagtwo\left(\gfy\right)+,  \cost<\codmagtwo\left(\gfx\right)<>\codmagtwo\left(\gfx\right)>+0+ \rboundtwo\comptwo
\lboundtwo \gfy,f \rboundtwo\right)\comptwo \Diag{\dommagtwo\left(\gfy\right)}{2} \right)\overset{(3)}{=}\\[4pt]
  \left(\left(\vecsum<n<>2>\comptwo 
\lboundtwo \idobj+\codmagtwo\left(\gfy\right)+, \cost<\codmagtwo\left(\gfx\right)<>\codmagtwo\left(\gfx\right)>+0+ \rboundtwo\right)\comptwo
\lboundtwo \gfy,f \rboundtwo\right)\comptwo \Diag{\dommagtwo\left(\gfy\right)}{2} \overset{(4)}{=}\\[4pt]
  \left(\proj<\codmagtwo\left(\gfx\right), \codmagtwo\left(f\right)<>1>\comptwo
\lboundtwo \gfy,f \rboundtwo\right)\comptwo \Diag{\dommagtwo\left(\gfy\right)}{2} \overset{(5)}{=}
\gfy\text{,}
\end{array}
\end{equation*}
where: equality (1) follows by definition of $\genfunsum$ in Definition \ref{defgenf}; equality (2) follows by relations \eqref{(S 0)}, \eqref{(R 7)}, \eqref{(R 17)}; equality (3) follows by relation \eqref{(R 7 quater)}; equality (4) follows by relation \eqref{(R 17)}; equality (5) follows by relations \eqref{(R 1 bis)}, \eqref{(R 17)}.\newline
\textnormal{\textbf{Proof of}} \eqref{star}-[$(iii)$]\textnormal{\textbf{.}}\ \ Fix $\gfx\in \genf$. Choose $n \in\mathbb{N}_0$,
$f \in \Cksp{\infty}(\argcompl{\dommagtwo\left(\gfx\right)})+\codmagtwo\left(\gfx\right)+$ with $\codmagtwo\left(\gfx\right)=\mathbb{R}^n$. Then the following chain of equalities holds true 
\begin{equation*}
\begin{array}{l}
1 \genfunscalp \gfx \overset{(1)}{=}
 \vecprod\left[1,n\right] \comptwo\left(\lboundtwo \cost<\dommagtwo\left(\gfx\right)<>\mathbb{R}>+1+,\gfx\rboundtwo\comptwo \Diag{\dommagtwo\left(\gfx\right)}{2}
\right)\overset{(2)}{=}\\[4pt]
  \vecprod\left[1,n\right]\comptwo \left(\left(
\lboundtwo \cost<\codmagtwo\left(\gfx\right)<>\codmagtwo\left(\gfx\right)>+1+,\idobj+\codmagtwo\left(\gfx\right)+ \rboundtwo\comptwo
\lboundtwo f,\gfx \rboundtwo\right)\comptwo \Diag{\dommagtwo\left(\gfx\right)}{2} \right)\overset{(3)}{=}\\[4pt]
  \left(\left(\vecprod\left[1,n\right]\comptwo 
\lboundtwo \cost<\codmagtwo\left(\gfx\right)<>\codmagtwo\left(\gfx\right)>+1+,\idobj+\codmagtwo\left(\gfx\right)+ \rboundtwo\right)\comptwo
\lboundtwo f,\gfx \rboundtwo\right)\comptwo \Diag{\dommagtwo\left(\gfx\right)}{2}\overset{(4)}{=}\\[4pt]
 \left(\proj<\codmagtwo\left(f\right), \codmagtwo\left(\gfx\right)<>2>\comptwo
\lboundtwo f,\gfx \rboundtwo\right)\comptwo \Diag{\dommagtwo\left(\gfx\right)}{2}\overset{(5)}{=}\gfx\text{,}
\end{array}
\end{equation*}
where: equality (1) follows by definition of $\genfunscalp$ in Definition \ref{defgenf}; equality (2) follows by relations \eqref{(S 0)}, \eqref{(R 7)}, \eqref{(R 17)}; equality (3) follows by relation \eqref{(R 7 quater)}; equality (4) follows by relation \eqref{(R 17)}; equality (5) follows by relations \eqref{(R 1 bis)}, \eqref{(R 17)}.\newline
\end{proof}

\begin{notation}\label{magtwopartic}\mbox{}
\begin{enumerate}
\item Elements belonging to $\genf$ are called generalized functions.\newline
Elements belonging to $\genfcont$ are called generalized continuous functions.
\item Fix $m,n\in\mathbb{R}^n$, $\intuno \subseteqdentro \mathbb{R}^m$. We set: 
\begin{flalign*}
&\genf(\intuno)+\mathbb{R}^n+=\left\{\gfx\in \genf\;:\quad\dommagtwo\left(\gfx\right)=\intuno\text{,}\quad\codmagtwo\left(\gfx\right)=\mathbb{R}^n\right\}\text{;}\\[4pt]
&\genfzero(\intuno)+\mathbb{R}^n+=\genf(\intuno)+\mathbb{R}^n+\cap\genfzero\text{;}\\[4pt]
&\genfempty(\intuno)+\mathbb{R}^n+=\genf(\intuno)+\mathbb{R}^n+\cap\genfempty\text{;}\\[4pt]
&\genfcont(\intuno)+\mathbb{R}^n+=\genf(\intuno)+\mathbb{R}^n+\cap\genfcont\text{;}\\[4pt]
&\genfcontcont(\intuno)+\mathbb{R}^n+=\genf(\intuno)+\mathbb{R}^n+\cap\magtwocontcont\text{;}\\[4pt]
&\genfsmooth(\intuno)+\mathbb{R}^n+=\genf(\intuno)+\mathbb{R}^n+\cap\magtwosmooth\text{.}
\end{flalign*}
\end{enumerate}
\end{notation}

\begin{remark}\label{Cinfdentro}\mbox{}
\begin{enumerate}
\item We have that $\magtwosmoothsmooth$ gives rise to a copy of $\idssmo$ in  $\genfidsp$. Then we will no longer distinguish between  
$\magtwosmoothsmooth$ and $\Cksp{\infty}$ in the notation. In particular we will denote elements belonging to $\Cksp{\infty}$ and $\magtwosmoothsmooth$ by the same symbol without making use of functions $\evalcomptwo$ or $\linincltwo$ any more.
\item Operation $\genfunscalp$ defines a bilateral action of $\mathbb{R}$ on $\genf$.
\end{enumerate}
\end{remark}

\section{Structure of generalized continuous functions\label{strcgf}}
In this section we study the structure of generalized continuous function and their behavior with respect to derivatives. We refer to Section \ref{fremag}.\vspace{12pt}

In Definition \ref{Ldef} below we introduce a special class of element belonging to $\fremagcontcont$.
We refer to Definition \ref{augdef}-[4].
\begin{definition}\label{Ldef} 
We denote by $\canrapp$ the subset of $\fremagcontcont$ defined recursively as follows:
\begin{flalign*}
&\begin{array}{l}
\text{any}\hspace{4pt}\fmx\in \canrapp\hspace{4pt}\text{is}\hspace{4pt}\text{a}\hspace{4pt}\text{maximal}\hspace{4pt}\text{element;}
\end{array}\\[8pt]
&\begin{array}{l}
\incllim_{-1}\left(f\right)\in \canrapp\hspace{51pt} \forall f \in \fremag_{-1} \text{;}
\end{array}\\[8pt]
&\begin{array}{l}
 \fmx_1 \compone  \fmx_2\in \canrapp \hspace{50pt}  \forall \fmx_1, \fmx_2 \in \canrapp \hspace{4pt}\text{with}\hspace{4pt}\fmx_1 \compone  \fmx_2\in\fremagcontcont \text{;}
\end{array}\\[8pt]
&\left\{
\begin{array}{l}
\lboundone \fmx_1,...,\fmx_n\rboundone\in \canrapp \hspace{20pt} \forall n \in \mathbb{N}\quad \forall\left(\fmx_1,...,\fmx_n\right)\in \fremag^n\\[4pt]
\text{with}\;\fmx_1 \in \lininclone\left(\Cksp{0}\setminus\left\{\contspempty\right\}\right)\cup\left\{\Fint_i\acmone\fmy \;:\; i\in\mathbb{N} \text{,}\;\fmy\in \canrapp \right\} \text{,}\\[4pt]
\left(\fmx_2,...,\fmx_n\right) \in \left(\lininclone\left(\Cksp{\infty}\setminus\left\{\smospempty\right\}\right)\right)^{n-1} \text{.}
\end{array}
\right.
\end{flalign*}
\end{definition} 

In Proposition \ref{Lprop} below we prove that any generalized continuous function admits a representative in $\canrapp$.
\begin{proposition}\label{Lprop}
Fix $\gfx \in \genfcontcont$. Then there is $\fmx \in \canrapp$ with $\gfx=\quotmagone\left(\fmx\right)$. 
\end{proposition} 
\begin{proof}
Statement follows by relations \eqref{(R 5)}-\eqref{(R 7 ter)}, \eqref{(R 17)}.
\end{proof}

In Definition \ref{Lorddef} below we introduce the notion of order which distinguishes elements belonging to $\canrapp$ with respect to the nesting level of integral operators.
\begin{definition}\label{Lorddef} 
We define the order of elements belonging to $\canrapp$ recursively as follows:
\begin{enumerate}
\item We say that $0$ is the order of $\fmx\in\canrapp$ if and only if there is an occurrence set $\occsettre_0=\left\{\left(\fmy_1,j_1\right),...,\left(\fmy_{\topindexcinque},j_{\topindexcinque}\right)\right\}$ for $\fmx$ fulfilling all conditions below:
\begin{flalign}
& \begin{array}{l} \fmx=\fmy_1\compone\left(\fmy_2\compone\left(...\compone\left( \fmy_{\topindexcinque-1}\compone\fmy_{\topindexcinque}\right)...\right)\right)\text{;}\end{array}\label{ordcond04}\\[8pt]
&\left\{
\begin{array}{l}
\forall \left(\fmy,j\right)\in \occsettre_0\quad \exists \topindexsei \in \mathbb{N}\quad\exists \left( {_1\fmy},...,{_{\topindexsei}\fmy}\right)\in \fremag^{\topindexsei} \quad\text{such that}\\[4pt]
\begin{array}{ll}
(i)&\fmy=\lboundone{_1\fmy},...,{_{\topindexsei}\fmy}\rboundone\text{,}\\[4pt]
(ii)&{_1\fmy}\in\lininclone\left(\Cksp{0}\setminus\left\{\contspempty\right\}\right)\text{,}\\[4pt]
(iii)&\left({_2\fmy},...,{_{\topindexsei}\fmy}\right)\in\left(\lininclone\left(\Cksp{\infty}\setminus\left\{\smospempty\right\}\right)\right)^{\topindexsei-1}\text{.}
\end{array}
\end{array}\label{ordcond05}
\right.
\end{flalign}
\item Fix $\topindexotto \in \mathbb{N}$. We say that $\topindexotto$ is the order of $\fmx\in \canrapp$ if and only if there are $\fmx_0\in \canrapp$ of order $0$, an occurrence set $\occsetdue=\left\{\left(\fmx_1,i_1\right),...,\left(\fmx_{\topindexquattro},i_{\topindexquattro}\right)\right\}$ for $\fmx$, two occurrence sets $\occsetdue_0=\left\{\left(\fmx_{0,1},i_{0,1}\right),...,\left(\fmx_{0,\topindexquattro},i_{0,\topindexquattro}\right)\right\}$, $\occsettre_0=\left\{\left(\fmy_1,j_1\right),...,\left(\fmy_{\topindexcinque},j_{\topindexcinque}\right)\right\}$ for $\fmx_0$, a set function $\theta:\occsetdue\rightarrow \mathbb{N}$, an injective set function $\rho_0 :\occsetdue_0\rightarrow\occsettre_0$ fulfilling all conditions below:
\begin{flalign}
&\begin{array}{l} 
\text{the}\hspace{4pt}\text{order}\hspace{4pt}\text{of}\hspace{4pt} \fmx_{\elindexquattro}\hspace{4pt}\text{belongs}\hspace{4pt}\text{to}\hspace{4pt} \left\{0,...,\topindexotto-1\right\}\quad \forall \elindexquattro \in\left\{1,...,\topindexquattro\right\}\text{;}
\end{array}\label{ordcond1}
\\[8pt]
&\begin{array}{l} 
\exists \elindexquattro \in\left\{1,...,\topindexquattro\right\}\;:\; \text{the}\hspace{4pt}\text{order}\hspace{4pt}\text{of}\hspace{4pt}\fmx_{\elindexquattro}\;\text{is}\;\topindexotto-1\; \text{;}
\end{array}\\[8pt]
&\left\{\begin{array}{l} 
\fmx=\strchar\left[\fmx_0, \varphi  \right]\quad
\text{where}\hspace{4pt}\varphi \hspace{4pt}\text{is}\hspace{4pt}\text{the}\hspace{4pt}\text{occurrence}\hspace{4pt}\text{function}\hspace{4pt}\text{on}\hspace{4pt} \occsetdue_0 \\[4pt]
\text{defined}\hspace{4pt}\text{by}\hspace{4pt}\text{setting}\\[4pt] \hspace{20pt}\varphi\left(\fmx_{0,\elindexquattro},i_{0,\elindexquattro}\right)=\Fint_{\theta\left(\fmx_{\elindexquattro},i_{\elindexquattro}\right)}\acmone\fmx_{\elindexquattro}\qquad \forall \elindexquattro \in \left\{1,...,\topindexquattro\right\}\text{;} 
\end{array}\label{ordcond3}\right.\\[8pt]
& \begin{array}{l} \fmx_0=\fmy_1\compone\left(\fmy_2\compone\left(...\compone\left( \fmy_{\topindexcinque-1}\compone\fmy_{\topindexcinque}\right)...\right)\right)\text{;}\end{array}\label{ordcond4}\\[8pt]
&\left\{
\begin{array}{l}
\forall \left(\fmy,j\right)\in \occsettre_0\quad \exists \topindexsei \in \mathbb{N}\quad\exists \left( {_1\fmy},...,{_{\topindexsei}\fmy}\right)\in \fremag^{\topindexsei} \quad\text{such}\hspace{4pt}\text{that}\\[4pt]
\begin{array}{ll}
(i)&\fmy=\lboundone{_1\fmy},...,{_{\topindexsei}\fmy}\rboundone\text{,}\\[4pt]
(ii)&{_1\fmy}=\fmx_{0,\elindexquattro} \hspace{49pt}\text{if}\hspace{4pt}\rho_0\left(\fmx_{0,\elindexquattro},i_{0,\elindexquattro}\right)=\left(\fmy,j\right)\text{,}\\[4pt]
(iii)&{_1\fmy}\in\lininclone\left(\Cksp{0}\setminus\left\{\contspempty\right\}\right)\hspace{10pt}
 \text{if}\;\left(\fmy,j\right)\in   \occsettre_0 \setminus \rho_0\left( \occsetdue_0\right)\text{,}\\[4pt]
(iv)&\left({_2\fmy},...,{_{\topindexsei}\fmy}\right)\in\left(\lininclone\left(\Cksp{\infty}\setminus\left\{\smospempty\right\}\right)\right)^{\topindexsei-1}\text{;}
\end{array}
\end{array}\label{ordcond5}
\right.\\[8pt]
&\fmx_{0,\elindexquattro}=\lininclone\left(\evalcompone\left(\Fint_{\theta\left(\fmx_{\elindexquattro},i_{\elindexquattro}\right)}\acmone\fmx_{\elindexquattro}\right)\right)\qquad \forall \elindexquattro \in \left\{1,...,\topindexquattro\right\}\text{.}\label{ordcond6}
\end{flalign} 
\end{enumerate}
We set:
\begin{flalign*}
&\canrapp(\topindexotto)=\left\{\fmx \in \canrapp\;:\;  \topindexotto\hspace{4pt}\text{is}\hspace{4pt}\text{the}\hspace{4pt}\text{order}\hspace{4pt}\text{of}\hspace{4pt}\fmx\right\}\hspace{15pt}\forall \topindexotto \in \mathbb{N}_0 \text{;}\\[6pt]
&\canrapp(-1)=\incllim_{-1}\left(\bascont\right)\text{.}
\end{flalign*}
\end{definition}

\begin{remark}\label{unlpertuttirem}
Relation \eqref{(R 1)} and choices made in Section \ref{cgf} entail that there is no loss of generality by improving \eqref{ordcond05}-[$(ii)$] to ${_1\fmy}\in\lininclone\left( \fremag_{-1}\right)$.\newline
Relations \eqref{(R 1)}, \eqref{(R 10)} entail that there is no loss of generality by adding the condition $\codmagone\left({_1\fmy}\right)=\mathbb{R}$ to \eqref{ordcond5}. 
\end{remark}

In Proposition \ref{unlpertutti} below we prove that any element belonging to $\canrapp$ has an order.
\begin{proposition}\label{unlpertutti}
Fix $\fmx\in \canrapp$. Then there is $\topindexotto \in \mathbb{N}_0$ such that $\fmx$ has order $\topindexotto$. 
\end{proposition}
\begin{proof}
Statement follows by Proposition \ref{invsost}, conditions \eqref{ordcond1}-\eqref{ordcond3}.
\end{proof}

Motivated by Proposition \ref{Lprop} we give Definition \ref{Worddef} below where we introduce the notion of order in $\genfcontcont$ which distinguishes its elements with respect to the nesting level of integral operators. 
\begin{definition}\label{Worddef} 
Fix $\topindexotto \in \mathbb{N}_0\cup\left\{-1\right\}$, $\gfx \in \genfcontcont$. We say that $\topindexotto$ is the order of $\gfx$ if and only if there is $\fmx\in \canrapp(\topindexotto)$ such that $\gfx=\quotmagone\left(\fmx\right)$.
We set:
\begin{flalign*}
&\genfcontcont[\topindexotto]=  \quotmagone\left(\canrapp(\topindexotto)\right)\hspace{15pt}\forall \topindexotto \in \mathbb{N}_0 \text{;}\\[6pt]
&\genfcontcont[-1]=\quotmagone\left(\canrapp(-1)\right)\text{.}
\end{flalign*}
\end{definition}

\begin{remark}\label{nontuttiinl}
Because of identifications induced by equivalence relation $\relsymb$ in Definition \ref{magintdef} it is possible that set $\genfcontcont[\topindexotto]$ contains elements of order lower than $l$.
\end{remark}

\chapter{Localization of generalized functions\label{shffgf}} 
In this chapter we study the behavior of generalized functions in a suitably small neighborhood of a fixed point belonging to their domain.

\section{Local behavior of generalized functions\label{lagf}}

In Lemma \ref{funintapprox} below we investigate the effect of integral operators in performing composition of generalized functions. We refer to Notations \ref{gentopnot}-[7, 8], \ref{realvec}-[4, 6].

\begin{lemma}\label{funintapprox}
Fix $m, n \in \mathbb{N}_0$, $i \in \mathbb{N}$, $ \intdue_0 \subseteqdentro \mathbb{R}^n$, $\fmx\in \fremagcont$, $\intuno_1\sqsubseteqdentro  \domint\left[\dommagone\left(\fmx\right), i\right]$, $\inttre_1\subseteqdentro \dommagone\left(\fmx\right)$, a path $\pmone$ in $\fremagcont$ detecting $\fmx$. Assume that all conditions below are fulfilled:
\begin{flalign}
& \inttre_1\compacont\dommagone\left(\fmx\right)\sqsubseteqdentro \mathbb{R}^m\text{;}\label{intcomp1}\\[8pt]
& \domint\left[\inttre_1,i\right]\compacont \intuno_1\compacont\domint\left[\dommagone\left(\fmx\right), i\right]\text{;}\nonumber\\[8pt]
&\exists \delta>0\;:\;
\left(\evalcompone\left(\pmone\left(\unkdue\right)\right)\right)\left(\unkuno\right)\neq \udenunk \quad \forall \unkdue\in \left(-\delta,\delta\right)\quad \forall\unkuno\in\dommagone\left(\fmx\right)\text{;}\nonumber\\[8pt]
&\evalcompone\left(\incl{\intdue_0}{\mathbb{R}^n}\compone\left(\left(\Fint_i \acmone \fmx\right)\compone\ \incl{\intuno_1}{\mathbb{R}^{m+1}}   \right)\right)\neq \contspempty\text{.}\nonumber
\end{flalign} 
Then there is $\delta_1 \in \left(0,\delta\right]$ such that
\begin{multline}
\left(\evalcompone\left(
\incl{\intdue_0}{\mathbb{R}^n}\compone
\left(
\Fint_i \acmone \left(\left(\pmone \left(\unkdue\right)\right)
\compone \incl{\inttre_1}{\mathbb{R}^{m}}\right)
\right)
\right)\right)\left(\unkuno\right)\neq \udenunk\\
 \forall \unkdue \in \left(-\delta_1,\delta_1\right)\quad \forall\unkuno\in\domint\left[\inttre_1,i\right]\text{.}
\label{intnz}
\end{multline}
\end{lemma}
\begin{proof}
Statement follows by continuity of integral operators on $\qidscont$ (see \eqref{extopcntdis3}) since $\left(\evalcompone\left(\Fint_i \acmone \fmx\right)\right)\left(\intuno_1\right)\compacont\intdue_0$ and $\intdue_0$ is open in $\mathbb{R}^n$.
\end{proof}

In Proposition \ref{pathprop7} below we prove that elements belonging to $\magtwocontcont$ are suitably linked, at least locally, with elements not belonging to $\magtwoempty$. We refer to Notations \ref{ins}-[5], \ref{gentopnot}-[7, 8].

\begin{proposition}\label{pathprop7}
Fix $\mtx \in\magtwocontcont$, $\unkuno_0\in\dommagtwo\left(\mtx\right)$. Then there is $\intuno\subseteqdentro\dommagtwo\left(\mtx\right)$ fulfilling both conditions below:
\begin{flalign}
&\begin{array}{l}\unkuno_0 \in \intuno \compacont\dommagtwo\left(\mtx\right)\text{;}\label{condcont1}
\end{array}\\[8pt]
&\left\{
\begin{array}{l} 
\text{for any path}\;\pmtwo\;\text{in}\;\magtwocont\;\text{detecting}\;\mtx\; \text{there is}\;\delta>0\;\text{such that, for}\\[4pt]
\text{any}\hspace{6pt} \text{skeleton}\hspace{6pt}
\left(\fmx,\hspace{1pt} \assfun\right)\hspace{6pt}\text{of}\hspace{6pt}\pmtwo\text{,}\hspace{6pt} \text{pair}\; \left(\fmx\compone\incl{\intuno}{\dommagtwo\left(\mtx\right)}, \hspace{1pt}\assfun\right)\hspace{6pt}\text{defines,}\\[4pt]
\text{through}\hspace{6pt} \text{equality}\hspace{6pt} \text{\eqref{pathconduro2},}\hspace{6pt} \text{a}\hspace{6pt} \text{skeleton}\hspace{6pt} \text{of}\hspace{6pt} \text{a}\hspace{6pt} \text{path}\hspace{6pt}\pmtwo_{\intuno}\hspace{6pt}\text{detecting}\\[4pt]
\mtx\comptwo \incl{\intuno}{\dommagtwo\left(\mtx\right)}\;\text{fulfilling}\\[4pt]
\hspace{53pt}\left(\evalcomptwo\left( \pmtwo_{\intuno}\left(\unkdue\right)\right)\right)\left(\unkuno\right)\neq\udenunk\quad \forall \unkdue\in \left(-\delta,\delta\right)\quad \forall\unkuno\in \intuno\text{.}
\end{array}\label{condcont2}
\right.
\end{flalign}
\end{proposition}
\begin{proof}
Proposition \ref{Lprop} entails that there are $\fmx \in \canrapp$ with $\mtx=\quotmagone\left(\fmx\right)$.\newline
Then statement follows by \eqref{augord}, Proposition \ref{augprop}-[5], \eqref{defevctw}, once we have proven the claim below:
\begin{equation}
\left\{
\begin{array}{l}
\text{Fix}\hspace{4pt} \fmx \in \canrapp\text{,}\hspace{4pt}\unkuno_0\in \dommagone\left(\fmx\right)\text{.}\hspace{5pt}\text{Then}\hspace{5pt}\text{there}\hspace{5pt}\text{is}\hspace{5pt}\intuno\subseteqdentro\dommagone\left(\fmx\right)\hspace{5pt}\text{with}\hspace{5pt}\unkuno_0\in\intuno\text{,}\\[4pt]
\intuno\compacont\dommagone\left(\fmx\right)\hspace{4pt}
\text{such}\hspace{4pt}\text{that}\hspace{4pt}\text{for}\hspace{4pt}\text{any}\hspace{4pt}\pmone\hspace{4pt}\text{in}\hspace{4pt}\fremagcont\hspace{4pt}\text{detecting}\hspace{4pt}\fmx\text{,}\hspace{4pt}\text{any}\hspace{4pt}\text{dete-}\\[4pt]
\text{cting}\hspace{5pt}\text{pair}\hspace{5pt}\left(\fmx,\assfun\right)\hspace{5pt}\text{for}\hspace{5pt}\pmone\text{,}\hspace{5pt}\text{there}\hspace{5pt}\text{are}\hspace{5pt}\text{an}\hspace{5pt}\text{occurrence}\hspace{5pt}\text{set}\hspace{4pt} \occset\hspace{4pt}\text{for} \hspace{4pt}\fmx\text{,}\\[4pt]
\text{an}\hspace{4pt}\text{occurrence}\hspace{4pt}\text{function}\hspace{4pt}\occfun\hspace{4pt}
\text{on}\hspace{4pt} \occset\text{,}\hspace{4pt} \delta >0\hspace{4pt}\text{fulfilling}\hspace{4pt}\text{all}\hspace{4pt}\text{conditions:}\\[4pt]
\begin{array}{ll}
(i)&\fmy\in \incllim_{-1}\left(\fremag_{-1}\right) 
\qquad\forall \left(\fmy,i\right)\in \occset\text{;}
\\[8pt]
(ii)&\forall \left(\fmy,i\right)\in \occset \quad \exists \inttre \sqsubseteqdentro \dommagone\left(\fmy\right)\;:\\[4pt]
&\hspace{75pt}\inttre \compacont\dommagone\left(\fmy\right) \text{,}\quad
\occfun\left(\fmy,i\right)=\fmy\compone \incl{\inttre}{\dommagone\left(\fmy\right)}\text{;}
\\[8pt]
(iii)&\text{pair}\;\left(\strchar\left[\fmx,\occfun\right], \assfun\right)\;\text{is a detecting skeleton of a path}\;\pmone\left[\occfun\right]\\[4pt]
&\text{in}\;\fremagcont	\;\text{detecting}\;\strchar\left[\fmx,\occfun\right]\text{;}
\\[8pt]
(iv)&\intuno\compacont \dommagone\left(\strchar\left[\fmx,\occfun\right]\right)\text{;}\\[8pt]
(v)&\left(\evalcompone\left(\pmone \left[\occfun\right]\left(\unkdue\right)\right)\right)\left(\unkuno\right)\neq\udenunk
 \quad \forall \unkdue \in \left(-\delta,\delta\right)\quad \forall\unkuno\in \intuno\text{.}
\end{array}
\end{array}
\right.\label{claimcond}
\end{equation}
We prove claim \eqref{claimcond} by induction on the nesting level of operators $\Fint_i$ in $\fmx \in \canrapp$.\newline
If $x\in \canrapp$ has order $0$ then claim \ref{ordcond3} follows by exploiting the following facts: conditions \eqref{ordcond04}, \eqref{ordcond05}; the image of a point through a continuous function is a point; the domain of any generalized function is an open set by definition (see Notation \ref{difrealfunc}).\newline
If $x\in \canrapp$ has order $l\geq 1$ then we argue as follows. Since $\fmx_0$ has order $0$, $\occsettre_0$ is a finite set, $\unkuno_0\in \dommagone\left(\fmx\right)$, the image of a compact set through continuous functions is a compact set, the domain of any generalized function is an open set by definition (see Notation \ref{difrealfunc}), then we have that
\begin{flalign}
&\left\{
\begin{array}{l}
\text{for any}\;\elindexcinque\in \left\{1,...,\topindexcinque\right\}\;\text{there are}\; \intuno_{\elindexcinque}\sqsubseteqdentro \dommagone\left(\fmy_{\elindexcinque}\right)$\text{,}\; $\intdue_{\elindexcinque}\sqsubseteqdentro \intuno_{\elindexcinque}\;\text{fulfilling}\\[4pt]
\text{all conditions below:}\\[4pt]
\begin{array}{ll}
(i)&\intdue_{\elindexcinque}\compacont\intuno_{\elindexcinque}\compacont \dommagone\left(\fmy_{\elindexcinque}\right)\text{,}\\[4pt]
(ii)& \unkuno_0\in \intdue_{\topindexcinque}\,\text{,}\\[4pt]
(iii)&\left(\evalcompone\left(\fmy_{\elindexcinque}\right)\right)\left(\intuno_{\elindexcinque}\right)\subseteq \intdue_{\elindexcinque-1} \quad\text{if}\hspace{4pt}\elindexcinque>1\text{;}
\end{array}
\end{array}
\right.\label{primcond}\\[8pt]
&\left\{
\begin{array}{l}
\text{for any}\;\elindexquattro\in \left\{1,...,\topindexquattro \right\}$, $\elindexcinque\in \left\{1,...,\topindexcinque\right\}\;\text{with}\; \rho_0\left(\fmx_{0,\elindexquattro},i_{0,\elindexquattro}\right)=\left(\fmy_{\elindexcinque},j_{\elindexcinque}\right)\\[4pt]
\text{there are}\;{{_1}\inttre_{\elindexcinque}}\sqsubseteqdentro \dommagone\left(\fmx_{\elindexquattro}\right)\text{,}\;
{{_2}\inttre_{\elindexcinque}}\sqsubseteqdentro \dommagone\left(\lboundone{{_2}\fmy_{\elindexcinque}},...,{{_{\topindexsei}}\fmy_{\elindexcinque}}\rboundone\right)
\;\text{fulfilling all}\\[4pt]
\text{conditions below:}\\[4pt]
\begin{array}{ll}
(i)&{{_1}\inttre_{\elindexcinque}}\compacont \dommagone\left(\fmx_{\elindexquattro}\right)\text{,}\\[8pt]
(ii)&{{_2}\inttre_{\elindexcinque}}\compacont \dommagone\left(\lboundone{{_2}\fmy_{\elindexcinque}},...,{{_{\topindexsei}}\fmy_{\elindexcinque}}\rboundone\right)\text{,}\\[8pt]
(iii)&\intdue_{\elindexcinque}\compacont\domint\left[{{_1}\inttre_{\elindexcinque}},\,\theta\left(\fmx_{\elindexquattro},i_{\elindexquattro}\right)\right]\times{{_2}\inttre_{\elindexcinque}}\compacont\intuno_{\elindexcinque}\compacont \dommagone\left(\fmy_{\elindexcinque}\right)\text{.}
\end{array}
\end{array}
\right.\nonumber
\end{flalign}
We set:
\begin{flalign*}
&\widetilde{\fmx}= \lboundone\fmx_0,\fmx_1,...,\fmx_{\topindexquattro}\rboundone\text{;}\\[8pt]
&\widetilde{i}_{\elindexquattro}\hspace{4pt}\text{the}\hspace{4pt}\text{occurrence}\hspace{4pt}\text{of}\hspace{4pt}\fmx_{\elindexquattro}\hspace{4pt}\text{in}\hspace{4pt}\widetilde{\fmx}\hspace{4pt}\text{for}\hspace{4pt}\text{any}\hspace{4pt}\elindexquattro\in \left\{0,...,\topindexquattro \right\}\text{;}\\[8pt] 
&\setsymcinque_0=\intdue_{\elindexcinque}\,\text{;}\\[8pt] 
&\setsymcinque_{\elindexquattro}={{_1}\inttre_{\elindexcinque}}\qquad\text{if}\;\rho_0\left(\fmx_{0,\elindexquattro},i_{0,\elindexquattro}\right)=\left(\fmy_{\elindexcinque},j_{\elindexcinque}\right)\text{.}
\end{flalign*} 
We choose a path $\widetilde{\pmone}$ in $\fremagcont$ detecting $\widetilde{\fmx}$, a detecting skeleton $\left(\widetilde{\fmx},\widetilde{\assfun} \right)$ of $\widetilde{\pmone}$ such that  
$\quad\widetilde{\assfun}\left(f,\unkdue\right)=\assfun\left(f,\unkdue\right) \quad\forall \left(f,t\right)\in \left(\occset_{\fmx}\cap \occset_{\widetilde{\fmx}}\right) \times \left(-1,1\right)$.\newline   
We emphasize that $ \occset_{\widetilde{\fmx}}\setminus\occset_{\fmx}\subseteq \overset{\topindexquattro}{\underset{\elindexquattro=1}{\bigcup}}\occset_{\fmx_{0,\elindexquattro}}$.\newline
Fix $\elindexquattro\in \left\{0,...,\topindexquattro\right\}$. 
Referring to Definition \ref{subpath} we set $\widetilde{\pmone}_{\elindexquattro}=\pmone\left[\left(\fmx_{\elindexquattro},\widetilde{i}_{\elindexquattro}\right), \left(\widetilde{\fmx}, \widetilde{\assfun}\right)\right]$ the path in $\fremagcont$ defined by the pair $\left(\fmx_{\elindexquattro}, \widetilde{\assfun}\funcomp\left(\incl{\occset_{\fmx_{\elindexquattro}}}{\occset_{\widetilde{\fmx}}}, \idobj+\left(-1,1\right)+\right)\right)$ through \eqref{pathconduro2}.
By construction $\fmx_{\elindexquattro}\in \canrapp$ and the degree of $\fmx_{\elindexquattro}$ belongs to $\left\{0,...,l-1\right\}$ then, by induction, claim \eqref{claimcond} applies with data $\fmx_{\elindexquattro}$, $\setsymcinque_{\elindexquattro}$, $\widetilde{\pmone}_{\elindexquattro}$, $\widetilde{\assfun}\funcomp\left(\incl{\occset_{\fmx_{\elindexquattro}}}{\occset_{\widetilde{\fmx}}}, \idobj+\left(-1,1\right)+\right)$ and provides an occurrence set $\occset_{\elindexquattro}$ for $\fmx_{\elindexquattro}$, an occurrence function $\occfun_{\elindexquattro}$ on $\occset_{\elindexquattro}$, $\delta_{\elindexquattro}>0$ fulfilling \eqref{claimcond}-[$(i)$-$(v)$].\newline
If $\elindexquattro\in \left\{1,...,\topindexquattro\right\}$ then, by choosing data $\elindexcinque\in \left\{1,...,\topindexcinque \right\}$, $m_{\elindexquattro},m_{\elindexcinque-1}\in \mathbb{N}_0 $ with $\rho_{0}\left(\fmx_{0,\elindexquattro}, i_{0,\elindexquattro}\right)=\left(\fmy_{\elindexcinque}, j_{\elindexcinque}\right)$, $\dommagone\left(\fmy_{\elindexcinque-1}\right)\sqsubseteqdentro \mathbb{R}^{m_{\elindexcinque-1}}$, $\dommagone\left(\fmx_{\elindexquattro}\right)\sqsubseteqdentro \mathbb{R}^{m_{\elindexquattro}}$, by \eqref{ordcond5}-[$(ii)$], \eqref{ordcond6}, \eqref{primcond}-[$(iii)$], referring to Definition \ref{dimdomcodfuncontone} we have that Lemma \ref{funintapprox} applies with data
\begin{equation*}
\begin{array}{l}
m_{\elindexquattro}\text{,}\hspace{23pt} 
m_{\elindexcinque-1}\text{,}\hspace{23pt}
\theta\left(\fmx_{\elindexquattro}, i_{\elindexquattro} \right)\text{,}\hspace{23pt}
\proj<\codmagone\left(\fmx_{\elindexquattro}\right),\mathbb{R}^{m_{\elindexcinque-1}-\dimcodmagone\left(\fmx_{\elindexquattro}\right)}<>1>\left(\intdue_{\elindexcinque-1}\right)\text{,}\hspace{23pt}
\fmx_{\elindexquattro}\text{,}\\[4pt]
\proj<\mathbb{R}^{\dimdommagone\left(\fmx_{\elindexquattro}\right)+1}, \mathbb{R}^{\dommagone\left(\fmy_{\elindexcinque}\right)-\dimdommagone\left(\fmx_{\elindexquattro}\right)-1}<>1>\left(\intuno_{\elindexcinque}\right)
\text{,}\hspace{23pt}
{{_1}\inttre_{\elindexcinque}}\text{,}\hspace{23pt}
\widetilde{\pmone}_{\elindexquattro}\text{,} 
\end{array}
\end{equation*}
and provides $\delta_{1, \elindexquattro}\in \left(0,\delta_{\elindexquattro}\right]$ fulfilling \eqref{intnz}.\newline
Since we deal with a finite number of elements belonging to $\fremag$ (i.e. $\fmx_{0}$,...,$\fmx_{\topindexquattro}$, $\fmx_{0,1}$,...,$\fmx_{0,\topindexquattro}$ ), of occurrence sets (i.e. $\occset_{0}$,...,$\occset_{\topindexquattro}$), of occurrence functions (i.e. $\occfun_{0}$,...,$\occfun_{\topindexquattro}$), by Proposition \ref{propfinordseq}-[2] there is no loss of generality by assuming that for any $\elindexquattro\in \left\{0,...,\topindexquattro\right\}$ there are
$k_{\elindexquattro}, k_{0,\elindexquattro} \in \mathbb{N}_0$, $\fmz_{\elindexquattro}\in \fremag_{k_{\elindexquattro}}$ $\fmz_{0,\elindexquattro}\in \fremag_{k_{0,\elindexquattro}}$, an occurrence set $\extoccset_{\elindexquattro}$ for $\fmz_{\elindexquattro}$, 
a set function $\extoccfun_{\elindexquattro}:\extoccset_{\elindexquattro}\rightarrow \extfremag_{k_{\elindexquattro}}$ such that: $k_{\elindexquattro}\leq k_0$; $k_{0,\elindexquattro}\leq k_0$; data $\fmx_0$, $\fmx_{0,\elindexquattro}$, $i_{0,\elindexquattro}$, $k_{0}$, $k_{0,\elindexquattro}$, $\fmz_{0}$, $\fmz_{0,\elindexquattro}$ 
fulfill \eqref{sysocc} if $\elindexquattro\geq 1$; data $\fmx_{\elindexquattro}$, $\occset_{\elindexquattro}$, $k_{\elindexquattro}$, $\fmz_{\elindexquattro}$, $\extoccset_{\elindexquattro}$ fulfill \eqref{hypextelcomp}; data $\fmx_{\elindexquattro}$, $\occset_{\elindexquattro}$, $\occfun_{\elindexquattro}$, $k_{\elindexquattro}$, $\fmz_{\elindexquattro}$, $\extoccset_{\elindexquattro}$, $\extoccfun_{\elindexquattro}$ fulfill \eqref{hypasselcomp}. \newline
Proposition \ref{posrel}, \eqref{claimcond}-[$(i)$] entail that for any $\left(\fmz,j\right) \in\extoccset_{0}$ one and only one among the following mutually exclusive cases occurs:
\begin{flalign*}
&\text{there}\hspace{4pt}\text{is}\hspace{4pt} \elindexquattro\in\left\{1,...,\topindexquattro\right\} \hspace{4pt}\text{such}\hspace{4pt}\text{that}\\[4pt]
&\hspace{87pt}i_{0,\elindexquattro}\leq j< j+ \extlnght_{k_{\fmz}}\left(\fmz\right) \leq i_{0,\elindexquattro}+\extlnght_{k_{0,\elindexquattro}}\left(\fmz_{0,\elindexquattro}\right)\; \text{where}\;\fmz \in\extfremag_{k_{\fmz}}\text{;}\\[8pt]
&\text{for}\hspace{4pt}\text{any}\hspace{4pt} \elindexquattro\in\left\{1,...,\topindexquattro\right\} \hspace{4pt} \text{either}\hspace{4pt} i_{0,\elindexquattro}+\extlnght_{k_{0,\elindexquattro}}\left(\fmz_{0,\elindexquattro}\right)\leq j \\[4pt]
&\hspace{179pt}\text{or}\hspace{4pt}j+ \extlnght_{k_{\fmz}}\left(\fmz\right)\leq i_{0,\elindexquattro} \hspace{4pt} \text{where}\hspace{4pt}\fmz \in\extfremag_{k_{\fmz}}\text{.}
\end{flalign*}
Then we define:
\begin{flalign*}
&\begin{array}{l}
{_{0,0}\extoccsetstr}=
\left\{\left(\fmz,j\right)\,:\;\left(\fmz,j\right) \in\extoccset_{0}\hspace{4pt}\text{with}\hspace{4pt} j<i_{0,1}\right\}\text{;}
\end{array}\\[8pt]
&\begin{array}{l}
{_{0,\elindexquattro}\extoccsetstr}=\bigg\{\left(\fmz,j+\overset{\elindexquattro}{\underset{\elindexuno=1}{\sum}}\left(\extlnght_{k_{\elindexuno}}\left(\fmz_{\elindexuno}\right)-\extlnght_{k_{0,\elindexuno}}\left(\fmz_{0,\elindexuno}\right)\right)+2\elindexquattro\right)\,:\\[4pt]
\hspace{14pt}\left(\fmz,j\right) \in\extoccset_{0}\quad\text{with}\quad i_{0,\elindexquattro }+\extlnght_{k_{0,\elindexquattro }}\left(\fmz_{0,\elindexquattro }\right)\leq j<i_{0,\elindexquattro +1}\bigg\}\hspace{15pt}\text{for}\hspace{4pt}\elindexquattro\in \left\{1,...,\topindexquattro-1\right\}\text{;}
\end{array}\\[8pt]
&\begin{array}{l}
{_{0,\topindexquattro}\extoccsetstr}=\bigg\{\left(\fmz,j+\overset{\topindexquattro}{\underset{\elindexuno=1}{\sum}}\left(\extlnght_{k_{\mathsf{i}}}\left(\fmz_{\elindexuno}\right)-\extlnght_{k_{0,\elindexuno}}\left(\fmz_{0,\elindexuno}\right)\right)+2\topindexquattro\right)\,:\\[4pt]
\hspace{146pt}\left(\fmz,j\right) \in\extoccset_{0}\quad\text{with}\quad i_{0,\topindexquattro}+\extlnght_{k_{0,\topindexquattro}}\left(\fmz_{0,\topindexquattro}\right)\leq j\bigg\}\text{;}
\end{array}\\[8pt]
&\begin{array}{l}
{_{1}\extoccsetstr}=
\left\{\left(\fmz,j+i_{0,1}+1\right)\,:\;\left(\fmz,j\right) \in\extoccset_{1}\right\}\text{;}
\end{array}\\[8pt]
&\begin{array}{l}
{_{\elindexquattro}\extoccsetstr}=
\bigg\{\left(\fmz,j+i_{0,\elindexquattro}+\overset{\elindexquattro-1}{\underset{\elindexuno=1}{\sum}}\left(\extlnght_{k_{\elindexuno}}\left(\fmz_{\elindexuno}\right)-\extlnght_{k_{0,\elindexuno}}\left(\fmz_{0,\elindexuno}\right)\right)+2\elindexquattro-1\right)\,:\\[4pt]
\hspace{197pt}\left(\fmz,j\right) \in\extoccset_{\elindexquattro}\bigg\}\hspace{15pt} \text{for}\hspace{4pt}\elindexquattro\in \left\{2,...,\topindexquattro\right\}\text{;}
\end{array}\\[8pt]
&\begin{array}{l}
\extoccsetstr= \left(\overset{\topindexquattro}{\underset{\elindexquattro=0}{\bigcup}}\hspace{4pt}{_{0,\elindexquattro}\extoccsetstr}\right)\hspace{4pt}\cup\hspace{4pt}\left(\overset{\topindexquattro}{\underset{\elindexquattro=1}{\bigcup}}\hspace{4pt}{_{\elindexquattro}\extoccsetstr}\right)\text{;}
\end{array}\\[8pt]
&\begin{array}{l}
\text{the}\hspace{4pt}\text{set}\hspace{4pt}\text{function}\hspace{4pt}\extasselstr: \extoccsetstr\rightarrow \extfremag_{k_{0}}\hspace{4pt}\text{by}\hspace{4pt}\text{setting} \\[4pt] 
\extasselstr\left(\fmz,j \right)=\left\{
\begin{array}{l}
\extoccfun_{0}\left(\fmz,j\right)\hspace{116pt}\text{if}\hspace{4pt}\left(\fmz,j\right)\in {_{0,0}\extoccsetstr}\text{,}\\[6pt]
\extoccfun_{0}\left(\fmz,j-
\overset{\elindexquattro}{\underset{\elindexuno=1}{\sum}}\left(\extlnght_{k_{\elindexuno}}\left(\fmz_{\elindexuno}\right)-\extlnght_{k_{0,\elindexuno}}\left(\fmz_{0,\elindexuno}\right)\right)-2\elindexquattro
\right)\\[4pt]
\hspace{150pt}\text{if}\hspace{4pt}\left(\fmz,j\right)\in {_{0,\elindexquattro}\extoccsetstr}\text{,} \hspace{4pt}\elindexquattro\in\left\{1,...,\topindexquattro\right\}\text{,}\\[6pt]
\extoccfun_{1}\left(\fmz,j-i_{0,1}-1
\right)\hspace{742pt}\text{if}\hspace{4pt}\left(\fmz,j\right)\in {_{1}\extoccsetstr}\text{,}\\[6pt]
\extoccfun_{\elindexquattro}\left(\fmz,j-i_{0,\elindexquattro}-
\overset{\elindexquattro-1}{\underset{\elindexuno=1}{\sum}}\left(\extlnght_{k_{\elindexuno}}\left(\fmz_{\elindexuno}\right)-\extlnght_{k_{0,\elindexuno}}\left(\fmz_{0,\elindexuno}\right)\right)-2\elindexquattro+1
\right)\\[4pt]
\hspace{150pt}\text{if}\hspace{4pt}\left(\fmz,j\right)\in {_{0,\elindexquattro}\extoccsetstr}\text{,} \hspace{4pt}\elindexquattro\in\left\{1,...,\topindexquattro\right\}\text{.}
\end{array}
\right.
\end{array}
\end{flalign*}
If $\extoccsetstr=\udenset$ then we set: $\intuno=\setsymcinque_0$, $\occsetstr=\udenset$; $\occfunstr=\emptysetfun$; $\delta= \min\left\{\delta_0, \delta_{1,1},...,\delta_{1,\topindexquattro}\right\}$.\newline  
If $\extoccsetstr\neq\udenset$ then we set:  
\begin{flalign*}
&\intuno=\setsymcinque_0\text{;}\\[4pt]
&\occsetstr=
 \left\{\left(\incllim_{k\left[z,j\right]} \left(\fmz\right),j\right)\,:\;
 \left(\fmz,j\right)\in \extoccsetstr\right\}\text{;}\\[4pt]
&\occfunstr\left(\incllim_{k\left[\fmz,j\right]} \left(\fmz\right),j \right)=\incllim_{k\left[\fmz,j\right]}\left(
\extasselstr\left(\fmz,j \right)\right)\hspace{10pt} \forall \left(\fmz,j \right)\in \extoccsetstr\text{;}\\[4pt]
&\delta= \min\left\{\delta_0, \delta_{1,1},...,\delta_{1,\topindexquattro}\right\}\text{.}
\end{flalign*}
Eventually claim \eqref{claimcond} follows directly by checking conditions $(i)$-$(v)$.
\end{proof}

In Defintion \ref{infnearpointdef} below we focus our attention to subsets of $\genf$ containing elements whose domain is an open neighborhood of $0$.

\begin{definition}\mbox{}\label{infnearpointdef}
\begin{enumerate}
\item We set:
\begin{flalign*}
&\genfloc=\left\{\gfx\in \genf :\; 0\in \dommagtwo\left(\gfx\right) \right\}\text{;} \\[4pt] 
&\genflocccz=\left\{\gfx\in \magtwocontcont :\quad 0\in \dommagtwo\left(\gfx\right),\quad \left(\evalcomptwo\left(u\right)\right)\left(0\right)=0\right\}\text{.}
\end{flalign*}
Elements belonging to $\genfloc$ are called local generalized functions.\newline
Elements belonging to $\genflocccz$ are called local generalized continuous pointed functions.\newline
\item Fix $m,n\in\mathbb{R}^n$, $\intuno \subseteqdentro \mathbb{R}^m$ with $0\in \intuno$. We set: 
\begin{flalign*}
&\genfloc(\intuno)+\mathbb{R}^n+=\left\{\gfx\in \genfloc\;:\quad\dommagtwo\left(\gfx\right)=\intuno\text{,}\quad\codmagtwo\left(\gfx\right)=\mathbb{R}^n\right\}\text{;}\\[4pt]
&\genflocccz(\intuno)+\mathbb{R}^n+=\left\{\gfx\in \genflocccz\;:\quad\dommagtwo\left(\gfx\right)=\intuno\text{,}\quad\codmagtwo\left(\gfx\right)=\mathbb{R}^n\right\}\text{.}
\end{flalign*}
\end{enumerate}
\end{definition}

In Remark \ref{genercontfun0} below we describe the link between $\genfloc$ and the localization process in the setting of smooth and continuous functions. We refer to Examples \ref{ExCinf}, \ref{ExCzero}.
\begin{remark}\label{genercontfun0}
Construction of $\genfloc$ extends to $\genfidsp$ the construction of spaces of local smooth functions and of local continuous functions which is well known in $\idssmo$ and in $\qidscont$ respectively (see \cite{HU}).\newline
We denote by $\Cksploc{\infty}$ the space of local smooth functions, by $\Cksplocccz{\infty}$ the space of local smooth pointed functions, by $\Cksploc{0}$ the space of local continuous functions, by $\Cksplocccz{0}$ the space of local continuous pointed functions.
\end{remark}

In Definition \ref{dualnot0}, Remark \ref{dualrem} below we introduce the notion of localization of generalized functions. We refer to Notations \ref{magtwopartic}, \ref{realfunc}-[6], Definition \ref{magintdef}.

\begin{definition}\label{dualnot0}\mbox{}
\begin{enumerate}
\item Fix $m,n\in \mathbb{N}_0$, a path connected topological space $\setsymuno$, $\intuno\subseteqdentro\mathbb{R}^m$, a continuous function $\pmtwo:\setsymuno\rightarrow\genf(\intuno)+\mathbb{R}^n+$.\newline
We define the continuous function $\prefuncaap{\pmtwo}:\intuno\times\setsymuno\rightarrow\genf(\mathbb{R}^m)+\mathbb{R}^n+$ by setting
\begin{equation*}
\prefuncaap{\pmtwo}\left(\unkuno,\elsymuno\right)= \pmtwo\left(\elsymuno\right)\genfuncomp
\trasl[\unkuno]\hspace{15pt}\forall \left(\unkuno,\elsymuno\right)\in \intuno\times\setsymuno\text{.}
\end{equation*} 
We say that $\prefuncaap{\pmtwo}$ is the localization of $\pmtwo$. \newline
If $\setsymuno$ is the single point set then $\pmtwo$ can be identified with an element $\gfx$ belonging to $\genf$, hence we write $\prefuncaap{\gfx}$ in place of $\prefuncaap{\pmtwo}$.
\item Fix $m,n\in \mathbb{N}_0$, a path connected topological space $\setsymuno$, $\intuno\subseteqdentro\mathbb{R}^m$, a continuous function $\pmtwo:\setsymuno\rightarrow\genf(\intuno)+\mathbb{R}^n+$.\newline
We define the continuous function $\funcaap{\pmtwo}:\intuno\times\setsymuno\rightarrow\genf(\mathbb{R}^m)+\mathbb{R}^n+$ by setting
\begin{multline*}
\funcaap{\pmtwo}\left(\unkuno,\elsymuno\right)= \prefuncaap{\pmtwo}\left(\unkuno,\elsymuno\right)\genfunsum 
 -1 \genfunscalp\left(\prefuncaap{\pmtwo}\left(\unkuno,\elsymuno\right)\genfuncomp \cost<\mathbb{R}^{m}<>\mathbb{R}^{m}>+0+\right)\\[4pt]
\forall \left(\unkuno,\elsymuno\right)\in \intuno\times\setsymuno\text{.}
\end{multline*} 
We say that $\funcaap{\pmtwo}$ is the pointed localization of $\pmtwo$.\newline
If $\setsymuno$ is the single point set then $\pmtwo$ can be identified with an element $\gfx$ belonging to $\genf$, hence we write $\funcaap{\gfx}$ in place of $\funcaap{\pmtwo}$.
\end{enumerate}
\end{definition}

\begin{remark}\mbox{}\label{dualrem}
\begin{enumerate}
\item If $\setsymuno$ is not path connected then $\prefuncaap{\pmtwo}$ and $\funcaap{\pmtwo}$ are defined for any path component of $\setsymuno$ separately.
\item Localization and pointed localization of a generalized functions extend to $\genfidsp$ the corresponding construction which is well known in $\idssmo$ and in $\qidscont$. Localization and pointed localization of continuous functions are denote by symbols used for corresponding operations on generalized functions. 
\end{enumerate}
\end{remark}

In Propositions \ref{comederint1} below we study localizing generalized functions behaves with respect to operations $\genfuncomp$, $\lboundgenf\,\rboundgenf$, $\genfbsfmu$.  
We refer to Notation \ref{realvec}-[6], Definition \ref{intdiffmon}.

\begin{proposition}\label{comederint1}\mbox{}
\begin{enumerate}
\item Fix $m_1,m_2,n\in \mathbb{N}_0$, path connected topological spaces $\setsymuno_1$, $\setsymuno_2$, $\intuno_1\subseteqdentro\mathbb{R}^m_1$, $\intuno_2\subseteqdentro\mathbb{R}^m_2$, continuous functions $\pmtwo_1:\setsymuno_1\rightarrow\genf(\intuno_1)+\mathbb{R}^{m_2}+$, $\pmtwo_2:\setsymuno_2\rightarrow\genf(\intuno_2)+\mathbb{R}^{n}+$.\newline
Define the continuous functions $\pmtwo_3:\intuno_2\times\setsymuno_1\times\setsymuno_2\rightarrow\genf(\intuno_1)+\mathbb{R}^{n}+$ by setting
\begin{multline*}
\pmtwo_3\left(\unkuno_2,\elsymuno_1, \elsymuno_2\right)=\pmtwo_2\left(\elsymuno_2\right)\genfuncomp\left(\trasl[\unkuno_2]\genfuncomp\pmtwo_1\left(\elsymuno_1\right)\right)\\[4pt]
\forall \left(\unkuno_2,\elsymuno_1, \elsymuno_2\right)\in \intuno_2\times\setsymuno_1\times\setsymuno_2\text{.}
\end{multline*}
Then 
\begin{multline*}
\prefuncaap{\pmtwo_2}\left(\unkuno_2,\elsymuno_2\right)\genfuncomp\prefuncaap{\pmtwo_1}\left(\unkuno_1,\elsymuno_1 \right)=
\prefuncaap{\pmtwo_3}\left(\unkuno_1, \unkuno_2,\elsymuno_1,\elsymuno_2 \right) \\[4pt]
 \forall\left( \unkuno_1, \unkuno_2,\elsymuno_1,\elsymuno_2 \right)\in\intuno_1\times\intuno_2\times\setsymuno_1\times\setsymuno_2\text{.}
\end{multline*}
\item Fix $\topindextre\in \mathbb{N}$, $m_{\elindextre}, n_{\elindextre}\in \mathbb{N}$, path connected topological space $\setsymuno_{\elindextre}$, $\intuno_{\elindextre}\subseteqdentro\mathbb{R}^{m_{\elindextre}}$, a continuous functions $\pmtwo_{\elindextre}:\setsymuno_{\elindextre}\rightarrow\genf(\intuno_{\elindextre})+\mathbb{R}^{n_{\elindextre}}+$ for any $\elindextre\in \left\{1,...,\topindextre\right\}$.\newline
Define the continuous functions $\pmtwo:\underset{\elindextre=1}{\overset{\topindextre}{\prod}}\setsymuno_{\elindextre}\rightarrow\genf(\underset{\elindextre=1}{\overset{\topindextre}{\prod}}\intuno_{\elindextre})+\underset{\elindextre=1}{\overset{\topindextre}{\prod}}\mathbb{R}_{n_{\elindextre}}+$ by setting
\begin{equation*}
\pmtwo\left(\elsymuno_1, ...,\elsymuno_{\topindextre}\right)=\lboundgenf \pmtwo_1\left(\elsymuno_1\right),...,\pmtwo_{\topindextre}\left(\elsymuno_{\topindextre}\right)\rboundgenf\hspace{15pt}
\forall \left(\elsymuno_1, ...,\elsymuno_{\topindextre}\right)\in\underset{\elindextre=1}{\overset{\topindextre}{\prod}}\setsymuno_{\elindextre} \text{.}
\end{equation*}
Then 
\begin{multline*}
\lboundgenf \prefuncaap{\pmtwo_1}\left(\unkuno_1,\elsymuno_1\right),...,\prefuncaap{\pmtwo}_{\topindextre}\left(\unkuno_{\topindextre},\elsymuno_{\topindextre}\right)\rboundgenf=\prefuncaap{\pmtwo}\left(\unkuno_1,...,\unkuno_{\topindextre},\elsymuno_1, ...,\elsymuno_{\topindextre}\right)\\[4pt]
 \forall\left(\unkuno_1,...,\unkuno_{\topindextre},\elsymuno_1, ...,\elsymuno_{\topindextre}  \right)\in \left(\underset{\elindextre=1}{\overset{\topindextre}{\prod}}\intuno_{\elindextre}\right) \times	\left(\underset{\elindextre=1}{\overset{\topindextre}{\prod}}\setsymuno_{\elindextre}\right)\text{.}
\end{multline*}
\item Fix $m,n\in \mathbb{N}_0$, $\topindexuno\in \mathbb{N}$ a path connected topological space $\setsymuno$, $\intuno\subseteqdentro\mathbb{R}^m$, a continuous function $\pmtwo:\setsymuno\rightarrow\genf(\intuno)+\mathbb{R}^n+$, $\bsfm\in \left\{\Fint_i, \Fpart_i, \Fp_i, \Fq_i, \Fqq_i \right\}$.\newline
Define : 
\begin{flalign*}
&\left\{
\begin{array}{l}
\text{the}\hspace{4pt}\text{continuous}\hspace{4pt}\text{function}\hspace{4pt}\pmtwo_1:\setsymuno\rightarrow\genf\hspace{4pt}\text{by}\hspace{4pt}\text{setting}\\[4pt]
\pmtwo_1\left(\elsymuno\right)=\bsfm \genfbsfmu \left(\pmtwo\left(\elsymuno\right)\right)\hspace{15pt}\forall \elsymuno \in \setsymuno\text{;}
\end{array}
\right.\\[8pt]
&\left\{
\begin{array}{l}
\text{the}\hspace{4pt}\text{smooth}\hspace{4pt}\text{function}\hspace{4pt} h: \intuno\rightarrow \domint\left[\intuno,\topindexuno\right]  \hspace{4pt}\text{by}\hspace{4pt}\text{setting}\\[4pt]
h\left(\unkuno_1,...,\unkuno_{m}\right)=\left\{
\begin{array}{ll}
\left(\unkuno_1,...,\unkuno_{\topindexuno-1},\unkuno_{\topindexuno},\unkuno_{\topindexuno},\unkuno_{\topindexuno+1},...,\unkuno_{m}\right)&\text{if}\hspace{4pt}\topindexuno\leq m\text{,}\\[4pt]
\left(\unkuno_1,...,\unkuno_{m},...,\unkuno_{\topindexuno-1},\unkuno_{\topindexuno},\unkuno_{\topindexuno}\right)&\text{if}\hspace{4pt}\topindexuno > m\text{.}
\end{array}
\right.
\end{array}
\right.
\end{flalign*}
Then 
\begin{equation*}
\bsfm\genfbsfmu \left(\prefuncaap{\pmtwo}\left(\unkuno, \elsymuno\right)\right)=\left\{
\begin{array}{ll}
\left(\prefuncaap{\pmtwo_1}\right)\left(\unkuno,\elsymuno \right) & \forall\left( \unkuno,\elsymuno\right)\intuno\times\setsymuno\hspace{15pt}\text{if}\hspace{4pt}\bsfm\in \left\{\Fpart_i, \Fp_i\right\}\text{;}\\[8pt]
\left(\prefuncaap{\pmtwo_1}\right)\left(h\left(\unkuno\right),\elsymuno\right)& \forall\left(\unkuno,\elsymuno \right)\intuno\times\setsymuno\hspace{15pt}\text{if}\hspace{4pt}\bsfm\in \left\{\Fint_i, \Fq_i, \Fqq_i \right\}\text{.}
\end{array}
\right.
\end{equation*}
\end{enumerate}
\end{proposition}
\begin{proof} 
Statements follow straightforwardly by direct computation.
\end{proof}

\begin{remark}
Propositions \ref{comederint1} hold true word by word in $\idssmo$ and in $\qidscont$. 
\end{remark}

In Proposition \ref{valutpuntuale} below we prove that the notion of local vanishing of generalized functions is independent on the choice of detecting paths.\newline
We refer to Notation \ref{magtwopartic}-[2].
\begin{proposition}\label{valutpuntuale}
Fix $m,n\in \mathbb{N}_0$, an open neighborhood $\intuno\subseteq \mathbb{R}^{m}$ of $0$, $\gfx\in\genf(\intuno)+\mathbb{R}^n+$. Then both statements below are equivalent:
\begin{flalign}
&\left\{\begin{array}{l}
\text{there}\hspace{5pt}\text{is}\hspace{5pt}\text{an}\hspace{5pt}\text{open}\hspace{5pt}
\text{neighborhood}\hspace{6pt}\intuno_{0}\subseteq\intuno\hspace{6pt}\text{of}\hspace{6pt}0\hspace{6pt}\text{such}\hspace{5pt}\text{that}
\hspace{5pt}\text{for}\hspace{5pt}\text{any}\\[4pt]
\text{path}\hspace{4pt}
\pmtwo\hspace{4pt}\text{in}\hspace{4pt}\genf\hspace{4pt}\text{detecting}\hspace{4pt}\gfx\text{,}\hspace{4pt}\unkdue\in \left(-1,1\right)\setminus\left\{0\right\}\text{,}\hspace{4pt}\unkuno \in\intuno_0\hspace{4pt}\text{occurs} \hspace{4pt}
\text{either}\\[4pt]
\left(\pmtwo\left(\unkdue\right)\right)\left(\unkuno\right)=\udenunk\hspace{4pt}
\text{or}\hspace{4pt}\left(\pmtwo\left(\unkdue\right)\right)\left(\unkuno\right)=0\text{;}
\end{array}\label{tesiintorno}
\right.\\[8pt]
&\left\{\begin{array}{l}
\text{for}\hspace{4pt}\text{any}\hspace{4pt}\text{path}\hspace{4pt}\pmtwo\hspace{4pt}\text{in}\hspace{4pt}\genf\hspace{4pt}\text{detecting}\hspace{4pt}\gfx\hspace{4pt}\text{there}\hspace{4pt}\text{is}\hspace{4pt}\text{an}\hspace{4pt}\text{open}\hspace{4pt}\text{neighborhood}\\[4pt]\intuno_{0}\subseteq \intuno\hspace{5pt}\text{of}\hspace{5pt}0\hspace{5pt}
\text{such}\hspace{4pt}\text{that}
\hspace{4pt}\text{for}\hspace{4pt}\text{any}\hspace{5pt}
\unkdue\in \left(-1,1\right)\setminus\left\{0\right\}\text{,}\hspace{5pt}\unkuno \in\intuno_0\hspace{5pt}\text{occurs}\\[4pt]
\text{either} \hspace{4pt}
\left(\pmtwo\left(\unkdue\right)\right)\left(\unkuno\right)=\udenunk\hspace{4pt}
\text{or}\hspace{4pt}\left(\pmtwo\left(\unkdue\right)\right)\left(\unkuno\right)=0
\text{.}
\end{array}\label{ipotesiintorno}
\right.
\end{flalign}
\end{proposition}
\begin{proof}\mbox{}\newline
Implication $\eqref{tesiintorno} \Rightarrow \eqref{ipotesiintorno}$ is straightforward then we prove 
$\eqref{ipotesiintorno} \Rightarrow \eqref{tesiintorno}$.\newline
Fix a local base $\left\{\setsymcinque_{\topindexuno}\right\}_{\topindexuno \in \mathbb{N}}$ of open neighborhoods of $0\in\intuno$. Arguing by contradiction and referring to Notation \ref{ins}-[5] we have that
\begin{equation}
\left\{
\begin{array}{l}
\text{for}\hspace{4pt}\text{any}\hspace{4pt}\topindexuno \in \mathbb{N}\hspace{4pt}\text{there}\hspace{4pt}\text{are}\hspace{4pt}\text{a}\hspace{4pt}\text{path}\hspace{4pt}\pmtwo_{\topindexuno}\hspace{4pt}\text{in}\hspace{4pt}\genf\hspace{4pt}\text{detecting}\hspace{4pt}
\gfx\text{,}\hspace{4pt}\text{a}\hspace{4pt}\text{point}\\[4pt]
\left(\unkdue_{\topindexuno}, \unkuno_{\topindexuno}\right)\in \left(\left(-1,1\right)\times\setsymcinque_{\topindexuno}\right)\setminus\left(\left\{0\right\}\times\setsymcinque_{\topindexuno}\right)\hspace{4pt}\text{such}\hspace{4pt}\text{that}\hspace{4pt}
\left(\pmtwo_{\topindexuno}\left(\unkdue_{\topindexuno}\right)\right)\left(\unkuno_{\topindexuno}\right)\neq 0
\text{.}
\end{array}
\right.\label{clunouno}
\end{equation}
Referring to Definitions \ref{magintdef}, \ref{pathtwodef} fix a maximal representative $\fmz$ of $\gfx$, a centered associating function $\assfun_{\topindexuno}$ of $\pmtwo_{\topindexuno}$ for any $\topindexuno \in \mathbb{N}$. \newline
Since  $\left\{\setsymcinque_{\topindexuno}\right\}_{\topindexuno \in \mathbb{N}}$ is a local base of open neighborhoods of $0 \in\intuno$ we have
\begin{equation}
\underset{\topindexuno\rightarrow +\infty}{\lim}\,\unkuno_{\topindexuno}=0\text{.}\label{limszero}
\end{equation}
Since $\left(-1,1\right)$ is contained in the compact subset $\left[-1,1\right]$ of $\mathbb{R}$ then there is no loss of generality by assuming that is $\unkdue_0 \in \left[-1,1\right]$ fulfilling
\begin{equation*}
\underset{\topindexuno\rightarrow +\infty}{\lim}\,\unkdue_{\topindexuno}=\unkdue_0\text{.}
\end{equation*}
Fix a sequence $\left\{\unkdue_{1,\topindexuno}\right\}_{\topindexuno\in \mathbb{N}}$ fulfilling both conditions below:
\begin{flalign*}
&\begin{array}{l}
0<\unkdue_{1,\topindexuno}<\unkdue_{1,\topindexuno+1} \hspace{15pt}\forall \topindexuno\in \mathbb{N}\text{;}
\end{array}\\[6pt]
&\begin{array}{l}
\underset{\topindexuno\rightarrow +\infty}{\lim}\,\unkdue_{1,\topindexuno}=1\text{.}
\end{array}
\end{flalign*}
Path connectivity of $\Cksp{0}(\Dom\left(f\right))+\mathbb{R}+$ for any $f\in \occset_{\gfx}$ (refer to Notation \ref{difrealfunc}-[2]) entail that there is a centered associating function $\assfun$ for $\occset_{\gfx}$ fulfilling
\begin{equation*}
\assfun\left(f,\unkdue_{1,\topindexuno}\right)=\assfun_{\topindexuno}\left(f,\unkdue_{\topindexuno}\right)\hspace{10pt} \forall f \in \occset_{\gfx}\text{,}\hspace{10pt} \forall \topindexuno \in \mathbb{N}\text{.}
\end{equation*}
We denote by $\pmtwo$ the path in $\magtwo$ detecting $\gfx$ whose detecting skeleton is $\left(\fmz,\assfun\right)$.\newline
Eventually statement follows since \eqref{ipotesiintorno} fails for path $\pmtwo$ because \eqref{limszero} holds true and 
\hspace{15pt}$\left(\pmtwo\left(\unkdue_{\topindexuno}\right)\right)\left(\unkuno_{\topindexuno}\right)=\left(\pmtwo_{\topindexuno}\left(\unkdue_{1,\topindexuno}\right)\right)\left(\unkuno_{\topindexuno}\right)\neq 0\hspace{15pt}\forall \topindexuno 	\in\mathbb{N}$. 
\end{proof}

Motivated by Remark \ref{valutpuntuale} in Definition \ref{releqinfnearpointdef} below we introduce a relation which identifies elements belonging to $\genfloc$ whenever they coincide locally at $0$.\newline
We refer to Notation \ref{realvec}, Definition \ref{admrelc}.

\begin{definition}\label{releqinfnearpointdef} Fix $\gfx_1,\gfx_2 \in \genfloc$. We say that $\gfx_1 \infneareqrel \gfx_2$ if and only if
there are $m \in \mathbb{N}_0$, $\intuno \subseteqdentro\mathbb{R}^m$ fulfilling all conditions below:
\begin{flalign}
& 0 \in \intuno\text{;}\nonumber\\[4pt]
& \intuno \subseteqdentro\dommagtwo\left(\gfx_i\right)\subseteqdentro \mathbb{R}^m\hspace{20pt} \forall i \in \left\{1,2\right\}\text{;}\label{iddomcond}\\[4pt]
& \gfx_1\genfuncomp\incl{\intuno}{\dommagtwo\left(\gfx_1\right)}=\gfx_2\genfuncomp\incl{\intuno}{\dommagtwo\left(\gfx_2\right)}\label{thcond}\text{.}
\end{flalign}
\end{definition}

In Proposition \ref{infnearpointpropbis} below we prove that $\infneareqrel$ is an equivalence relation and we examine its compatibility with the evaluation $\evalcomptwo$. We refer to Definition \ref{magintdef}.

\begin{proposition} \label{infnearpointpropbis}\mbox{}
\begin{enumerate}
\item $\infneareqrel$ is an equivalence relation.
\item Fix $\gfx_1,\gfx_2\in \genfloc$. Assume $\gfx_1 \infneareqrel \gfx_2$. Then:
\begin{equation*}
\occset_{\gfx_1} = \occset_{\gfx_2}\text{;}\hspace{20pt} \codmagtwo\left(\gfx_1\right) = \codmagtwo\left(\gfx_2\right)\text{;}\hspace{20pt}\dimdommagtwo\left(\gfx_1\right) = \dimdommagtwo\left(\gfx_2\right)\text{.}
\end{equation*}
\item Fix $\gfx_1,\gfx_2\in \genflocccz$. Assume that $\gfx_1 \infneareqrel \gfx_2$.\newline
Then there is 
$\intuno \subseteqdentro \dommagtwo\left(\gfx_1\right)\cap \dommagtwo\left(\gfx_2\right)$ fulfilling both conditions below:
\begin{equation*}
0 \dentro \intuno\text{,}\hspace{50pt}
\evalcomptwo\left(\gfx_1 \genfuncomp  \incl{\intuno}{\dommagtwo\left(\gfx_1\right)}\right) =\evalcomptwo\left(\gfx_2\genfuncomp  \incl{\intuno}{\dommagtwo\left(\gfx_2\right)}\right)\text{.}
\end{equation*}  
\end{enumerate}
\end{proposition}
\begin{proof}
Statement 1 straightforwardly follows by definition of $\infneareqrel$. Statement 2 follows by conditions \eqref{iddomcond}, \eqref{thcond}. Statement 3 follows by \eqref{extopcntdis5}, definition of $\infneareqrel$, Proposition \ref{funoncont}. 
\end{proof}

In Definition \ref{loccost0} below we introduce the set of locally vanishing generalized continuous functions.

\begin{definition}\label{loccost0}\mbox{}
\begin{enumerate}
\item We set:
\begin{multline*}
\genfloccczz=\Big\{\gfx\in \genflocccz\; :\; \exists \intuno\subseteqdentro \dommagtwo\left(\gfx\right)\;\text{with}\\[4pt]
 0\in \intuno\text{,}\qquad\evalcomptwo\left(\gfx\genfuncomp\incl{\intuno}{\dommagtwo\left(\gfx\right)}\right)=\cost<\intuno<>\codmagtwo\left(\gfx\right)>+0+ \Big\}\text{.}
\end{multline*}
Elements belonging to $\genfloccczz$ are called locally vanishing generalized continuous functions.
\item Fix $m,n\in\mathbb{R}^n$, $\intuno \subseteqdentro \mathbb{R}^m$ with $0\in \intuno$. We set: 
\begin{equation*}
\genfloccczz(\intuno)+\mathbb{R}^n+=\left\{\gfx\in \genfloccczz\;:\quad\dommagtwo\left(\gfx\right)=\intuno\text{,}\quad\codmagtwo\left(\gfx\right)=\mathbb{R}^n\right\}\text{.}
\end{equation*}
\end{enumerate}
\end{definition}

In Proposition \ref{infnearpointprop} below we examine compatibility of equivalence relation $\infneareqrel$ with the algebraic structure induced on $\genfloc$ by $\genf$.

\begin{proposition} \label{infnearpointprop}\mbox{}
\begin{enumerate}
\item Fix $\gfx_1,\gfx_2 \in \genfloc$. Assume that $\gfx_1 \infneareqrel \gfx_2$. Then:
\begin{equation*}
\gfx_1\in \genflocccz\Leftrightarrow \gfx_2\in \genflocccz\text{,}\hspace{50pt} \gfx_1\in \genfloccczz\Leftrightarrow \gfx_2\in \genfloccczz\text{.}
\end{equation*}
\item Fix $\gfx \in \genfloc$, $\gfy  \in \genflocccz$. Assume that both the following conditions are fulfilled: $\gfx \notin \genfempty$; $\Ima\left(\evalcomptwo\left(\gfy\right)\right)\subseteq \dommagtwo\left(\gfx\right)$.\newline
Then there is $\gfy_1\in\genflocccz$ fulfilling both conditions: $\gfy_1\infneareqrel \gfy$; $\gfx\genfuncomp \gfy_1 \notin \genfempty$. 
\item Fix $\gfx  \in \genf$, $ \gfy_1,\gfy_2 \in \genfloc$. Assume $\gfy_1\infneareqrel \gfy_2$. Then $\gfx\genfuncomp \gfy_1\infneareqrel \gfx\genfuncomp \gfy_2$.
\item Fix $\gfx_1,\gfx_2 \in \genfloc$, $\gfy_1,\gfy_2  \in \genflocccz$. Assume that both the following conditions are fulfilled:  
$\gfx_1\infneareqrel \gfx_2$, $\gfy_1\infneareqrel \gfy_2$. Then  $\gfx_1\genfuncomp \gfy_1\infneareqrel \gfx_2\genfuncomp \gfy_2$.
\item Fix $k\in \mathbb{N}$, four finite ordered sequences $\left(\gfx_{1},...,\gfx_{k}\right), \left(\gfy_{1},...,\gfy_{k}\right)\subseteq \genfloc^k$. Assume $\gfx_{\mathsf{k}}\infneareqrel \gfy_{\mathsf{k}}$ for any $\mathsf{k}\in \left\{1,...,k\right\}$.
Then $\lboundtwo \gfx_{1},...,\gfx_{k} \rboundtwo \infneareqrel \lboundtwo \gfy_{1},...,\gfy_{k} \rboundtwo$.
\item Fix $\bsfm \in \bsfM$, $\gfx_1, \gfx_2  \in \genfloc$. Assume $\gfx_1 \infneareqrel \gfx_2$. Then $\bsfm \genfbsfmu \gfx_1 \infneareqrel \bsfm \genfbsfmu \gfx_2$.
\end{enumerate}
\end{proposition}
\begin{proof} \mbox{}\newline
\textnormal{\textbf{Proof of statement 1.}}\ \ Statement follows by condition \eqref{thcond}.\newline
\textnormal{\textbf{Proof of statement 2.}}\newline
First we construct $\gfy_1$.\newline
Set $\gfy_0=\incl{\dommagtwo\left(\gfx\right)}{\codmagtwo\left(\gfy\right)}\genfuncomp\gfy$. Then $\gfy_0\in \magtwocontcont$ by assumption $\Ima\left(\evalcomptwo\left(\gfy\right)\right)\subseteq \dommagtwo\left(\gfx\right)$. Then we apply Proposition \ref{pathprop7} to data $\gfy_0$, $0 \in \dommagtwo\left(\gfy_0\right)$ obtaining a set $\intuno\subseteqdentro\dommagtwo\left(\gfy_0\right)$ such that data $\gfy_0$, $0$, $\intuno$ fulfill both \eqref{condcont1} and \eqref{condcont2}. We set $\gfy_1=\gfy\genfuncomp\incl{\intuno}{\codmagtwo\left(\gfy\right)}$. By construction we have $\gfy_1\infneareqrel \gfy$.\newline
Then we prove that $\gfx\genfuncomp \gfy_1 \notin \genfempty$.\newline
Fix a path $\pmtwo_{\gfx}$ in $\genf$ detecting $\gfx$, $\delta_{\gfx}>0$ such that\newline
\centerline{$\left(\evalcomptwo\left(\pmtwo_{\gfx}\left(\unkdue\right)\right)\right)\left(\unkuno\right)\neq\udenunk\quad\forall\unkdue \in  \left(-\delta_{\gfx},\delta_{\gfx}\right)\setminus\left\{0\right\} \quad \forall \unkuno\in \dommagtwo\left(\gfx\right)$,}\newline
we emphasize that such a detecting path exists by assumption $\gfx \notin \genfempty$.\newline
Fix a detecting skeleton $\left(\fmx_{\gfx},\assfun_{\gfx}\right)$ of $\pmtwo_{\gfx}$, then we are able to choose $\fmx_{\gfy_1}\in \fremagcontcont$ with $\gfy_1=\quotmagone\left(\fmx_{\gfy_1}\right)$, a centered associating function $\assfun_{\gfy_1}$ for $\occset_{\gfy_1}$ such that\newline
\centerline{$\assfun_{\gfx}\left(f,\unkdue\right)=\assfun_{\gfy_1}\left(f,\unkdue\right) \quad \forall\left(f,\unkdue\right) \in \left(\occset_{\gfx}\cap \occset_{\gfy_1}\right)\times \left(-1,1\right)$.}\newline
Pair $\left(\fmx_{\gfy_1},\assfun_{\gfy_1}\right)$ defines, through \eqref{pathconduro2}, a path $\pmtwo_{\gfy_1}$ in $\magtwocont$ detecting $\gfy_1$. Then construction of $\gfy_1$ entails that there is $\delta_{\gfy_1}>0$ fulfilling both conditions below:
\begin{flalign*}
&\begin{array}{l}
\left(\evalcomptwo\left(\pmtwo_{\gfy_1}\left(\unkdue\right)\right)\right)\left(\unkuno\right)\neq\udenunk\quad\forall\unkdue \in  \left(-\delta_{\gfy_1},\delta_{\gfy_1}\right) \quad \forall \unkuno\in \dommagtwo\left(\gfy_1\right)\text{;}
\end{array}\\[8pt]
&\begin{array}{l}
\Ima\left(\evalcomptwo\left(\pmtwo_{\gfy_1}\left(\unkdue\right)\right)\right)\subseteq \dommagtwo\left(\gfx\right)\quad\forall\unkdue \in  \left(-\delta_{\gfy_1},\delta_{\gfy_1}\right)\text{.}
\end{array}
\end{flalign*} 
By referring to Proposition \ref{pathprop2} we set $\pmtwo=\pmtwo_{\gfx} \comptwo \pmtwo_{\gfy_1}$, $\delta=\min\left\{\delta_{\gfx}, \delta_{\gfy_1}\right\}$ then by construction $\pmtwo$ is a path in $\magtwo$ detecting $\gfx\genfuncomp \gfy_1$ such that\newline
\centerline{$\left(\evalcomptwo\left( \pmtwo\left(\unkdue\right)\right)\right)\left(\unkuno\right)\neq\udenunk\quad
\forall \unkdue\in \left(-\delta,\delta\right)\setminus\left\{0\right\}\quad \forall\unkuno\in \dommagtwo\left(\gfx\genfuncomp \gfy_1\right)$.}\newline
\textnormal{\textbf{Proof of statements 3, 4, 5, 6.}}\ \ Statements follows straightforwardly by definition of $\infneareqrel$ and statement 2.
\end{proof}

In Definition \ref{pathgenfcent} below we introduce the notion of path in sets $\genfloc^n$.

\begin{definition}\label{pathgenfcent}\mbox{}
\begin{enumerate}
\item Fix $n \in \mathbb{N}$, a set function $\pmtwo:\left(-1,1\right) \rightarrow \genfloc^n$.\newline
We say that $\pmtwo$ is a path in $\genfloc^n$ if and only if $\incl{\genfloc^n}{\genf^n}\funcomp \pmtwo$ is a path in $\genf^n$.\newline
With an abuse of language we say that: any skeleton $\left(\left(\fmx_1,...,\fmx_n\right),\assfun\right)$ of $\incl{\genfloc^n}{\genf^n}\funcomp \pmtwo$ is a skeleton of $\pmtwo$; $\left(\fmx_1,...,\fmx_n\right)$ is a core of $\pmtwo$; $\assfun$ is an associating function of $\pmtwo$.
\item Fix $n \in \mathbb{N}$, a path $\pmtwo$ in $\genfloc^n$, a subset $\singsmpth \subseteq \left(-1,1\right)$.\newline
We say that pair $\left(\pmtwo,\singsmpth\right)$ is a smooth path in $\genfloc^n$ if and only if pair defined by $\left(\incl{\genfloc^n}{\genf^n}\funcomp \pmtwo,\singsmpth\right)$ is a smooth path in $\genf^n$.\newline
With an abuse of language we say that: $\singsmpth$ is the singular set of $\pmtwo$; any skeleton $\left(\left(\fmx_1,...,\fmx_n\right),\assfun\right)$ of $\left(\incl{\genfloc^n}{\genf^n}\funcomp \pmtwo,\singsmpth\right)$ is a skeleton of $\left(\pmtwo,\singsmpth\right)$.\newline
We emphasize that nothing is assumed about $\pmtwo\left(\unkdue\right)$ when $\unkdue \in \singsmpth$.
\item Fix $n \in \mathbb{N}$, $\left(\gfx_1,...,\gfx_n\right)\in \genfloc^n$, a path $\pmtwo$ in $\genfloc^n$.\newline
We say that $\pmtwo$ is a path in $\genfloc^n$ through $\left(\gfx_1,...,\gfx_n\right)$ if and only if we have $\pmtwo\left(0\right)=\left(\gfx_1,...,\gfx_n\right)$. 
\item Fix $\gfx\in \genfloc$, a path $\pmtwo$ in $\genfloc$ through $\gfx$.
We say that $\pmtwo$ is a path in $\genfloc$ detecting $\gfx$ if and only if $\incl{\genfloc}{\genf}\funcomp \pmtwo$ is a path in $\genf$ detecting $u$.
With an abuse of language we say that any detecting skeleton $\left(\fmx,\assfun\right)$ of $\incl{\genfloc}{\genf}\funcomp\pmtwo$ is a detecting skeleton of $\pmtwo$.
\end{enumerate}
\end{definition}
 
In Definition \ref{topgenfcent} below we give topology to sets $\genfloc^n$ through paths introduced in Definition \ref{pathgenfcent}. We refer to Notation \ref{gentopnot}-[5].  

\begin{definition}\label{topgenfcent}\mbox{}
Fix $n \in \mathbb{N}$. Set $\genfloc^n$ is endowed with the final topology with respect to the family of all paths in $\genfloc^n$. We say that such topology is the path topology on $\genfloc^n$.
\end{definition}

In Definition \ref{pathgenfcentcontcontzero} below we introduce the notion of path in sets $\genflocccz^n$.

\begin{definition}\label{pathgenfcentcontcontzero}\mbox{}
\begin{enumerate}
\item Fix $n \in \mathbb{N}$, a set function $\pmtwo:\left(-1,1\right) \rightarrow \genflocccz^n$.\newline
We say that $\pmtwo$ is a path in $\genflocccz^n$ if and only if $\incl{\genflocccz^n}{\genfcont^n}\funcomp \pmtwo$ is a path in $\genfcont^n$.\newline
With an abuse of language we say that: any skeleton $\left(\left(\fmx_1,...,\fmx_n\right),\assfun\right)$ of $\incl{\genflocccz^n}{\genfcont^n}\funcomp \pmtwo$ is a skeleton of $\pmtwo$; $\left(\fmx_1,...,\fmx_n\right)$ is a core of $\pmtwo$; $\assfun$ is an associating function of $\pmtwo$.
\item Fix $n \in \mathbb{N}$, a path $\pmtwo$ in $\genflocccz^n$, a subset $\singsmpth \subseteq \left(-1,1\right)$.\newline
We say that path $\left(\pmtwo,\singsmpth\right)$ is a smooth path in $\genflocccz^n$ if and only if path $\left(\incl{\genflocccz^n}{\genfcont^n}\funcomp \pmtwo,\singsmpth\right)$ is a smooth path in $\genfcont^n$.\newline
With an abuse of language we say that: $\singsmpth$ is the singular set of $\pmtwo$; any skeleton $\left(\left(\fmx_1,...,\fmx_n\right),\assfun\right)$ of $\left(\incl{\genflocccz^n}{\genfcont^n}\funcomp \pmtwo,\singsmpth\right)$ is a skeleton of $\left(\pmtwo,\singsmpth\right)$.\newline
We emphasize that nothing is assumed about $\pmtwo\left(\unkdue\right)$ when $\unkdue \in \singsmpth$.
\item Fix $n \in \mathbb{N}$, $\left(\gfx_1,...,\gfx_n\right)\in \genflocccz^n$, a path $\pmtwo$ in $\genflocccz^n$.\newline
We say that $\pmtwo$ is a path in $\genflocccz^n$ through $\left(\gfx_1,...,\gfx_n\right)$ if and only if we have $\pmtwo\left(0\right)=\left(\gfx_1,...,\gfx_n\right)$.\newline
\end{enumerate}
\end{definition}

In Definition \ref{topgenfcentcontcontzero} below we give topology to sets $\genflocccz^n$ through paths introduced in Definition \ref{pathgenfcentcontcontzero}. We refer to Notation \ref{gentopnot}-[5]. 

\begin{definition}\label{topgenfcentcontcontzero}\mbox{}
Fix $n \in \mathbb{N}$. Set $\genflocccz^n$ is endowed with the final topology with respect to the family of all paths in $\genflocccz^n$. We say that such topology is the path topology on $\genflocccz^n$.
\end{definition}

In Definition \ref{pathgenfcentcontcontzerozero} below we introduce the notion of path in sets $\genfloccczz^n$.

\begin{definition}\label{pathgenfcentcontcontzerozero}\mbox{}
\begin{enumerate}
\item Fix $n \in \mathbb{N}$, a set function $\pmtwo:\left(-1,1\right) \rightarrow \genfloccczz^n$.\newline
We say that $\pmtwo$ is a path in $\genfloc^n$ if and only if $\incl{\genfloccczz^n}{\genfcont^n}\funcomp \pmtwo$ is a path in $\genfcont^n$.\newline
With an abuse of language we say that: any skeleton $\left(\left(\fmx_1,...,\fmx_n\right),\assfun\right)$ of $\incl{\genfloccczz^n}{\genfcont^n}\funcomp \pmtwo$ is a skeleton of $\pmtwo$; $\left(\fmx_1,...,\fmx_n\right)$ is a core of $\pmtwo$; $\assfun$ is an associating function of $\pmtwo$.
\item Fix $n \in \mathbb{N}$, a path $\pmtwo$ in $\genfloccczz^n$, a subset $\singsmpth \subseteq \left(-1,1\right)$.\newline
We say that path $\left(\pmtwo,\singsmpth\right)$ is a smooth path in $\genfloccczz^n$ if and only if path $\left(\incl{\genfloccczz^n}{\genfcont^n}\funcomp \pmtwo,\singsmpth\right)$ is a smooth path in $\genfcont^n$.\newline
With an abuse of language we say that: $\singsmpth$ is the singular set of $\pmtwo$; any skeleton $\left(\left(\fmx_1,...,\fmx_n\right),\assfun\right)$ of $\left(\incl{\genfloccczz^n}{\genfcont^n}\funcomp \pmtwo,\singsmpth\right)$ is a skeleton of $\left(\pmtwo,\singsmpth\right)$.\newline
We emphasize that nothing is assumed about $\pmtwo\left(\unkdue\right)$ when $\unkdue \in \singsmpth$.
\item Fix $n \in \mathbb{N}$, $\left(\gfx_1,...,\gfx_n\right)\in \genfloccczz^n$, a path $\pmtwo$ in  $\genfloccczz^n$.\newline
We say that $\pmtwo$ is a path in $\genfloccczz^n$ through $\left(\gfx_1,...,\gfx_n\right)$ if and only if we have $\pmtwo\left(0\right)=\left(\gfx_1,...,\gfx_n\right)$. 
\end{enumerate}
\end{definition}

In Definition \ref{topgenfcentcontcontzero} below we give topology to sets $\genfloccczz^n$ through paths introduced in Definition \ref{pathgenfcentcontcontzerozero}. We refer to Notation \ref{gentopnot}-[5].

\begin{definition}\label{topgenfcentcontcontzerozero}\mbox{}
Fix $n \in \mathbb{N}$. Set $\genfloccczz^n$ is endowed with the final topology with respect to the family of all paths in $\genfloccczz^n$. We say that such topology is the path topology on $\genfloccczz^n$.
\end{definition}

\section{Germs of generalized functions\label{ggf}}

Motivated by condition \eqref{iddomcond}, Propositions \ref{infnearpointpropbis}, \ref{infnearpointprop} we introduce germs of local generalized functions and domain and co-domain of them. We refer to Notations \ref{ins}-[13], \ref{gentopnot}-[8], Definition \ref{magintdef}.

\begin{definition} \label{genpointdef} \mbox{}
We set:
\begin{flalign*}
&\genfquot=\genfloc / \infneareqrel\text{;}\\[6pt]
&\genfquotfun : \genfloc \rightarrow \genfquot\hspace{4pt}\text{the}\hspace{4pt}\text{quotient}\hspace{4pt}\text{function;}\\[6pt]
&\genfquotcccz=\genfquotfun\left(\genflocccz\right)\text{;}\\[6pt]
&\genfquotccczz=\genfquotfun\left(\genfloccczz\right)\text{.}
\end{flalign*}
Elements belonging to $\genfquot$ are called germs of local generalized functions, or generalized germs for short.\newline
Elements belonging to $\genfquotcccz$ are called germs of local generalized continuous pointed functions, or continuous pointed generalized germs for short.\newline
Elements belonging to $\genfquotccczz$ are called germs of locally vanishing generalized continuous functions, or locally vanishing continuous generalized germs for short.\newline
Generalized germs are equivalence classes. Generalized germs are denoted by $\gggfx$ or by $\gggfx[\gfx]$ whenever we need to emphasize a representative $\gfx\in \genfloc$ of the equivalence class.\newline
With an abuse of language we denote any element $\gggfx[\gfx]$ with $\gfx \infneareqrel 0\genfunscalp \gfx$ by $\zeroquot$ and we say that $\gggfx[\gfx]$ vanishes.\newline
Motivated by Proposition \ref{infnearpointpropbis}-[2], we are allowed to denote by $\occset_{\gggfx}$ any set $\occset_{\gfx}$ with $ \gggfx =\genfquotfun\left(\gfx \right)$.\newline
With an abuse of language we say that a path in $\genf$ detecting a representative of a generalized germ $\gggfx$ is a path in $\genf$ detecting $\gggfx$.  
\end{definition}

In Proposition \ref{dimcoddomgerm} below we prove that set functions $\dimcodmagtwo$, $\dimdommagtwo$ factors through $\genfquotfun$.
\begin{proposition}\label{dimcoddomgerm}
There are set functions: 
\begin{flalign*}
& \codgerm:\genfquot\rightarrow\mathbb{N}_0\quad\text{by setting}  \quad \codgerm\left(\gggfx[\gfx]\right)=\dimcodmagtwo\left(\gfx\right)\quad\forall\gggfx[\gfx]\in \genfquot \text{;}\\[8pt]
&  \domgerm:\genfquot\rightarrow\mathbb{N}_0\quad\text{by setting}  \quad \domgerm\left(\gggfx[\gfx]\right)=\dimdommagtwo\left(\gfx\right)\quad\forall\gggfx[\gfx]\in \genfquot\text{.}
\end{flalign*}
With an abuse of language we denote $\codgerm$ again by $\dimcodmagtwo$ and $\domgerm$ again by $\dimdommagtwo$.
\end{proposition}
\begin{proof}
Statement follows straightforwardly by Proposition \ref{infnearpointpropbis}-[2].
\end{proof}

\begin{remark}\mbox{}\label{genercontfun}
\begin{enumerate}
\item The quotient topology of $\genfquot$ with respect to $\genfquotfun$ is the trivial topology. The proof coincides with the proof of the corresponding result which is well known to hold true for germs of local continuous functions.
\item We refer to Remark \ref{genercontfun0}, Examples \ref{ExCinf}, \ref{ExCzero}.\newline
Construction of $\genfquot$ extends to $\genfidsp$ the construction of spaces of germs of local smooth functions and of local continuous functions which is well known in $\idssmo$ and in $\qidscont$ respectively (see \cite{HU}).\newline
We denote by $\Ckspquot{\infty}$ the space of germs of local smooth functions, by $\Ckspquotccz{\infty}$ the space of germs of local smooth pointed functions, by $\Ckspquot{0}$ the space of germs of local continuous functions, by $\Ckspquotccz{0}$ the space of germs of local continuous pointed functions. Quotient functions are denoted by the same symbol used in Definition \ref{genpointdef}. In both cases germs are equivalence classes. Germs are denoted by $\gczsfx$ or by $\gczsfx[f]$ whenever we need to emphasize a representative $f \in \Cksploc{0}$ of the equivalence class, the equivalence class represented by $\smospempty$ is denoted simply by $\classquotempty$. Algebraic, intgro-differential, topological structures existing on $\idssmo$, $\qidscont$ induce corresponding well known structures on $\Ckspquot{\infty}$, $\Ckspquot{0}$ respectively. With an abuse of language we use symbols denoting operations in $\idssmo$, $\qidscont$ to denote the corresponding operations induced on $\Ckspquot{\infty}$, $\Ckspquot{0}$ respectively.\newline
Inclusion $\incl{\Cksploc{\infty}}{\Cksploc{0}}$ induces an inclusion $\incl{\Ckspquot{\infty}}{\Ckspquot{0}}$.\newline
\end{enumerate}
\end{remark}

In Proposition \ref{genpointprop} below we study the algebraic structure induced on $\genfquot$ by $\genfloc$. We refer to Notation \ref{realfunc}-[3(a), 4].

\begin{proposition}\label{genpointprop}\mbox{}
\begin{enumerate}
\item There is one and only one set function
\begin{multline*}
\compgenfquot:\genf \times \genfquot \rightarrow \genfquot\;\; \text{defined by setting}\\[4pt] 
 \gfx \,\compgenfquot\gggfy[\gfy]  =\gggfx[\gfx \genfuncomp \gfy]
\qquad \forall \left( \gfx, \gggfy[\gfy] \right) \in \genf \times \genfquot \text{.}
\end{multline*}
We emphasize that $\gfx \,\compgenfquot\gggfy \in \genfquotcccz$ for any $\left(\gfx ,\gggfy\right)\in\genflocccz\times\genfquotcccz$.
\item There is one and only one set function
\begin{multline*}
\compquotpoint:\genfquot \times \genfquotcccz \rightarrow \genfquot\;\; \text{defined by setting}\\[4pt] 
\gggfx[\gfx] \, \compquotpoint\,\gggfy[\gfy]  =\gggfx[\gfx \genfuncomp \gfy]
\qquad\forall \left(\gggfx[\gfx], \gggfy[\gfy] \right) \in \genfquot \times \genfquotcccz \text{.}
\end{multline*}
We emphasize that: 
\begin{flalign*}
&\gggfx\compquotpoint\gggfy \in \genfquotcccz\hspace{23pt}
\forall\left(\gggfx,\gggfy\right)\in\genfquotcccz\times\genfquotcccz\text{;}\\[6pt]
&\gggfx\compquotpoint \gggfy \in \genfquotccczz\hspace{20pt}
\forall\left(\gggfx,\gggfy\right)\in\genfquotcccz\times\genfquotccczz\cup\genfquotccczz\times\genfquotcccz\text{.}
\end{flalign*}
\item Fix $n \in \mathbb{N}$. There is one and only one set function 
\begin{multline*}
\lboundquotpoint\quad \rboundquotpoint:\genfquot^n\rightarrow \genfquot\;\text{defined by setting}\\[4pt]
\lboundquotpoint \gggfx[\gfx_1] ,...,\gggfx[\gfx_n] \rboundquotpoint=\gggfx[\lboundtwo \gfx_i,...,\gfx_n \rboundtwo] \qquad \forall \left(\gfx_i,...,\gfx_n\right)\in \genfquot^n\text{.}
\end{multline*}
We emphasize that: 
\begin{flalign*}
&\lboundquotpoint\gggfx_1,...,\gggfx_n\rboundquotpoint\in \genfquotcccz\hspace{24pt} \forall \left(\gggfx_1,...,\gggfx_n\right)\in\genfquotcccz^n\text{;}\\[4pt]
&\lboundquotpoint\gggfx_1,..., \gggfx_n \rboundquotpoint\in \genfquotccczz\hspace{20pt} \forall \left(\gggfx_1,...,\gggfx_n\right)\in\genfquotccczz^n\text{.}
\end{flalign*}
\item There is one and only one set function
\begin{multline*}
\sqsumquotpoint:\genfquot \times \genfquot \rightarrow \genfquot\quad \text{by setting}\quad
\gggfx[\gfx]\,\sqsumquotpoint\,\gggfy[\gfy] =
\vecsum<n<>2>\compgenfquot\lboundquotpoint
\gggfx[\gfx], \gggfy[\gfy] 
 \rboundquotpoint\\[4pt]
\forall n \in \mathbb{N}_0 \quad \forall (\gggfx[\gfx], \gggfy[\gfy] ) \in \genfquot \times \genfquot\quad\text{with}\quad\dimcodmagtwo\left(\gggfy[\gfy]\right)=n\text{.}
\end{multline*}
We emphasize that:\newline
$\gggfx \sqsumquotpoint\gggfy \in \genfquotcccz$ for any $\left(\gggfx,\gggfy\right)\in\genfquotcccz\times\genfquotcccz$ with $\dimcodmagtwo\left(\gggfx\right)=\dimcodmagtwo\left(\gggfy\right)$;\newline
$\gggfx\sqsumquotpoint \gggfy \in \genfquotccczz$ for any $\left(\gggfx,\gggfy\right)\in\genfquotccczz\times\genfquotccczz$ with $\dimcodmagtwo\left(\gggfx\right)=\dimcodmagtwo\left(\gggfy\right)$.
\item There is one and only one set function $
\sumquotpoint:\genfquot \times \genfquot \rightarrow \genfquot$ defined by setting:\newline
$\gggfx[\gfx]\,\sumquotpoint\,\gggfy[\gfy]  =\gggfx[\gfx \genfunsum \gfy]$ for any $\left(\gggfx[\gfx],\gggfy[\gfy]\right)\in\genfquot\times\genfquot$ with $\dimdommagtwo\left(\gggfx[\gfx]\right)\neq\dimdommagtwo\left(\gggfy[\gfy]\right)$,\newline
$\gggfx[\gfx]\,\sumquotpoint\,\gggfy[\gfy]  =\gggfx[\argcompl{\left(\gfx\genfuncomp\incl{\dommagtwo\left(\gfx\right)\cap\dommagtwo\left(\gfy\right)}{\dommagtwo\left(\gfx\right)}\right)\genfunsum \left(\gfy\genfuncomp\incl{\dommagtwo\left(\gfx\right)\cap\dommagtwo\left(\gfy\right)}{\dommagtwo\left(\gfy\right)}\right)}]$ for any $\left(\gggfx[\gfx],\gggfy[\gfy]\right)\in\genfquot\times\genfquot$ with $\dimdommagtwo\left(\gggfx[\gfx]\right)=\dimdommagtwo\left(\gggfy[\gfy]\right)$.\newline
We emphasize that:\newline
$\gggfx \sumquotpoint\gggfy \in \genfquotcccz$ for any $\left(\gggfx,\gggfy\right)\in\genfquotcccz\times\genfquotcccz$ with $\dimdommagtwo\left(\gggfx\right)=\dimdommagtwo\left(\gggfy\right)$, $\dimcodmagtwo\left(\gggfx\right)=\dimcodmagtwo\left(\gggfy\right)$;\newline
$\gggfx\sumquotpoint \gggfy \in \genfquotccczz$ for any $\left(\gggfx,\gggfy\right)\in\genfquotccczz\times\genfquotccczz$ with $\dimdommagtwo\left(\gggfx\right)=\dimdommagtwo\left(\gggfy\right)$, $\dimcodmagtwo\left(\gggfx\right)=\dimcodmagtwo\left(\gggfy\right)$.
\item There is one and only one set function
\begin{multline*}
\sqmultquotpoint:\genfquot \times \genfquot \rightarrow \genfquot\quad \text{by setting}\quad
\gggfx[\gfx]\,\sqmultquotpoint\,\gggfy[\gfy]  =
\vecprod\left[m,n\right]\compgenfquot\lboundquotpoint
\gggfx[\gfx], \gggfy[\gfy] 
 \rboundquotpoint\\[4pt]
\forall (\gggfx[\gfx], \gggfy[\gfy] ) \in \genfquot \times \genfquot\quad\text{with}\quad \dimcodmagtwo\left(\gggfx[\gfx]\right)=m\text{,}\quad\dimcodmagtwo\left(\gggfy[\gfy]\right)=n\text{.}
\end{multline*}
We emphasize that:
\begin{flalign*}
&\gggfx \sqmultquotpoint\gggfy\in \genfquotcccz\hspace{24pt}\forall\left(\gggfx,\gggfy\right)\in\genfquotcccz\times\genfquotcccz \text{;}\\[4pt]
&\gggfx\sqmultquotpoint \gggfy\in \genfquotccczz\hspace{20pt}\forall \left(\gggfx,\gggfy\right)\in\genfquotcccz\times\genfquotccczz\cup\genfquotccczz\times\genfquotcccz\text{.}
\end{flalign*}
\item There is one and only one set function $
\multquotpoint:\genfquot \times \genfquot \rightarrow \genfquot$ defined by setting:\newline
$\gggfx[\gfx]\,\multquotpoint\,\gggfy[\gfy]  =\gggfx[\gfx \genfunmult \gfy]$ for any $\left(\gggfx[\gfx],\gggfy[\gfy]\right)\in\genfquot\times\genfquot$ with $\dimdommagtwo\left(\gggfx[\gfx]\right)\neq\dimdommagtwo\left(\gggfy[\gfy]\right)$,\newline
$\gggfx[\gfx]\,\multquotpoint\,\gggfy[\gfy]  =\gggfx[\argcompl{\left(\gfx\genfuncomp\incl{\dommagtwo\left(\gfx\right)\cap\dommagtwo\left(\gfy\right)}{\dommagtwo\left(\gfx\right)}\right)\genfunmult \left(\gfy\genfuncomp\incl{\dommagtwo\left(\gfx\right)\cap\dommagtwo\left(\gfy\right)}{\dommagtwo\left(\gfy\right)}\right)}]$ for any $\left(\gggfx[\gfx],\gggfy[\gfy]\right)\in\genfquot\times\genfquot$ with $\dimdommagtwo\left(\gggfx[\gfx]\right)=\dimdommagtwo\left(\gggfy[\gfy]\right)$.\newline
We emphasize that:\newline
$\gggfx \multquotpoint\gggfy\in \genfquotcccz$ for any $\left(\gggfx,\gggfy\right)\in\genfquotcccz\times\genfquotcccz$ with $\dimdommagtwo\left(\gggfx\right)=\dimdommagtwo\left(\gggfy\right)$;\newline
$\gggfx\multquotpoint \gggfy\in \genfquotccczz$ for any $\left(\gggfx,\gggfy\right)\in\genfquotcccz\times\genfquotccczz\cup\genfquotccczz\times\genfquotcccz$ with $\dimdommagtwo\left(\gggfx\right)=\dimdommagtwo\left(\gggfy\right)$.
\item There is one and only one set function
\begin{multline*}
\scalpquotpoint:\mathbb{R} \times \genfquot \rightarrow \genfquot\;\text{defined by setting}\\[4pt]
a\,\scalpquotpoint\,\gggfx[\gfx]  =\gggfx[\argcompl{a \genfunscalp \gfx}] \qquad \forall \left(a, \gggfx[\gfx]\right)  \in \mathbb{R} \times \genfquot\text{.}
\end{multline*}
We emphasize that: 
\begin{flalign*}
&a\scalpquotpoint\gggfx\in \genfquotcccz\hspace{25pt}\forall \gggfx\in\genfquotcccz\text{;}\\[4pt]
&a\scalpquotpoint\gggfx \in \genfquotccczz\hspace{20pt}\forall \gggfx\in\genfquotccczz\text{.}
\end{flalign*}
\item There is one and only one set function
\begin{multline*}
\bsfmuquotpoint: \bsfM \times \genfquot \rightarrow \genfquot\;\text{defined by setting}\\[4pt]
\bsfm\,\bsfmuquotpoint\,\gggfx[\gfx]  =\gggfx[\argcompl{\bsfm \genfbsfmu \gfx}]\qquad \forall (\bsfm, \gggfx[\gfx] ) \in \bsfM \times \genfquot\text{.}
\end{multline*}
We emphasize that: 
\begin{flalign*}
&\bsfm\bsfmuquotpoint\gggfx\in \genfquotcccz\hspace{24pt}\forall\left(\bsfm,\gggfx\right)\in \subbsfM\times \genfquotcccz\text{;}\\[4pt]
&\bsfm\bsfmuquotpoint\gggfx \in \genfquotccczz\hspace{20pt}\forall\left(\bsfm,\gggfx\right)\in \subbsfM\times\genfquotccczz\text{.}
\end{flalign*}
\item There is one and only one set function $\lininclquot :\Ckspquot{\infty}\rightarrow\genfquot$ defined by setting 
$\qquad\lininclquot\left(\gczsfx[f]\right)= \gggfx[f]\qquad \forall f \in\Cksploc{\infty}$. 
\item There is one and only one set function $\evalcompquotcontcontzero :\genfquotcccz \rightarrow \Ckspquotccz{0}$ defined by setting $\qquad\evalcompquotcontcontzero\left( \gggfx[\gfx]\right) =\gczsfx[\evalcomptwo\left(\gfx\right)] \qquad \forall\gfx \in\genflocccz$.
\end{enumerate}
\end{proposition}
\begin{proof}
Statements follow by Definition \ref{defgenf}, Propositions \ref{infnearpointpropbis}, \ref{infnearpointprop}.
\end{proof}

\begin{remark}\label{existoperquot}\mbox{}
\begin{enumerate}
\item Propositions \ref{infnearpointprop}, \ref{genpointprop} entail that for any combination of finite ordered sequences of elements $\left(a_1,...,a_k\right)\in \mathbb{R}^k$, $\left(\bsfm_1,...,\bsfm_l\right)\in \bsfM^l$, $\left( \gfx_1 ,..., \gfx_m \right)\in\genf^m$,  $\left(\gggfy[\gfy_1] ,...,\gggfy[\gfy_n] \right)\in\genfquot^n$ obtained by a finite application of operations $\compgenfquot$, $\compquotpoint$, $\lboundquotpoint\quad \rboundquotpoint$,  $\sqsumquotpoint$, $\sumquotpoint$, $\sqmultquotpoint$, $\multquotpoint$, $\scalpquotpoint$, $\bsfmuquotpoint$ which admits a representative not belonging to $\magtwoempty$, there is no loss of generality by assuming that such representative is obtained by performing the same sequence of operations on data $\left(a_1,...,a_k\right)\in \mathbb{R}^k$, $\left(\bsfm_1,...,\bsfm_l\right)\in \bsfM^l$, $\left( \gfx_1 ,..., \gfx_m \right)\in\genf^m$,  $\left( \gfy_1 ,...,\gfy_n \right)\in\genf^n$.
\item We refer to Notation \ref{realfunc}-[1], Remark \ref{genercontfun}-[2].\newline
Proposition \ref{dimcoddomgerm} holds true for germs of local smooth functions and for germs of local continuous functions by replacing  
$\genfquot$ by $\Ckspquot{\infty}$ and $\Ckspquot{0}$,  $\dimcodmagtwo$ by $\dimCod$, $\dimdommagtwo$ by $\dimDom$.\newline
Proposition \ref{genpointprop} holds true word by word by replacing $\genf$, $\genfquot$, $\genfquotcccz$, respectively by $\Cksp{\infty}$, $\Ckspquot{\infty}$, $\Ckspquotccz{\infty}$, or by $\Cksp{0}$, $\Ckspquot{0}$, $\Ckspquotccz{0}$.\newline 
Set functions obtained by specializing Proposition \ref{genpointprop} to germs of smooth and continuous functions  will be denoted by same symbols used to denote corresponding set functions constructed in Proposition \ref{genpointprop}.\newline
Statement corresponding to Remark \ref{existoperquot}-[1] holds true for finite sequences of operations in $\Ckspquot{\infty}$, $\Ckspquot{0}$.\newline   
Topology of germs of local generalized, smooth and continuous functions is matter of Section \ref{Grtopss}.   
\end{enumerate}
\end{remark}

Motivated by condition \eqref{iddomcond}, Proposition \ref{infnearpointpropbis}-[2], in Definition \ref{genpointdefmnloc} below we restrict our attention to generalized functions with the same domain or co-domain.

\begin{definition}\mbox{}\label{genpointdefmnloc}
Fix $m,n \in \mathbb{N}_0$. We set:
\begin{flalign*}
&\begin{array}{l}
\genfquot(m)+n+= \left\{
\gggfx \;:\quad \gggfx \in \genfquot\text{,}\quad \dimdommagtwo\left(\gggfx\right)=m\text{,}\quad 
\dimcodmagtwo\left(\gggfx\right)=n\right\} \text{;}
\end{array}\\[8pt]
&\begin{array}{l}
\genfquot(\segnvar)+n+= \left\{
\gggfx \;:\quad \gggfx\in \genfquot\text{,}\quad 
\dimcodmagtwo\left(\gggfx\right)=n\right\} \text{.}
\end{array}
\end{flalign*}
\end{definition}

In Propositon \ref{operAn} below we study algebraic structure induced on $\genfquot(m)+n+$ and $\genfquot(\segnvar)+n+$ by $\genfquot$. 

\begin{proposition}\label{operAn}\mbox{}
\begin{enumerate}
\item Fix $a \in \mathbb{N}$, $m_1,...,m_a,n_1,...,n_a\in \mathbb{N}_0$.
Then there is one and only one set function 
\begin{multline*}
\lboundquotgradpoint\;\rboundquotgradpoint: \overset{a}{\underset{\mathsf{a}=1}{\prod}}\,\genfquot(m_{\mathsf{a}})+n_{\mathsf{a}}+\rightarrow\genfquot(\overset{a}{\underset{\mathsf{a}=1}{\sum}}\,m_{\mathsf{a}})+\overset{a}{\underset{\mathsf{a}=1}{\sum}}\,n_{\mathsf{a}}+\qquad\text{fulfilling}\\[4pt]
\lboundquotpoint\;\rboundquotpoint \funcomp \incl{\overset{a}{\underset{\mathsf{a}=1}{\prod}}\,\genfquot(m_{\mathsf{a}})+n_{\mathsf{a}}+}{\genfquot^a}= \incl{\genfquot(\overset{a}{\underset{\mathsf{a}=1}{\sum}}\,m_{\mathsf{a}})+\overset{a}{\underset{\mathsf{a}=1}{\sum}}\,n_{\mathsf{a}}+}{\genfquot} \funcomp \lboundquotgradpoint\;\rboundquotgradpoint\text{.}
\end{multline*}
With an abuse of language, from now on, we denote $\lboundquotgradpoint\;\rboundquotgradpoint$ again by $\lboundquotpoint\;\rboundquotpoint$. 
\item Fix $m_1,m_2, n \in \mathbb{N}_0$. There is one and only one set function
\begin{multline*}
\sqsumquotgradpoint :\genfquot(m_1)+n+\times \genfquot(m_2)+n+ \rightarrow \genfquot(m_1+m_2)+n+\qquad\text{fulfilling}\\[4pt]
 \sqsumquotpoint \funcomp \incl{\genfquot(m_1)+n+\times \genfquot(m_2)+n+}{\genfquot^2}= \incl{\genfquot(m_1+m_2)+n+}{\genfquot} \funcomp \sqmultquotgradpoint\text{.}
\end{multline*}
With an abuse of language, from now on, we denote $\sqsumquotgradpoint$ again by $\sqsumquotpoint$. 
\item Fix $m,n \in \mathbb{N}_0$.
Then there is one and only one set function 
\begin{multline*}
\sumquotgradpoint :\genfquot(m)+n+\times \genfquot(m)+n+ \rightarrow \genfquot(m)+n+\qquad \text{fulfilling}\\[4pt]
\sumquotpoint \funcomp \incl{\left(\genfquot(m)+n+\right)^2}{\genfquot^2}= \incl{\genfquot(m)+n+}{\genfquot} \funcomp \sumquotgradpoint\text{.}
\end{multline*}
With an abuse of language, from now on, we denote $\sumquotgradpoint$ again by $\sumquotpoint$.
\item Fix $m_1,m_2, n_1,n_2 \in \mathbb{N}_0$. There is one and only one set function
\begin{multline*}
\sqmultquotgradpoint :\genfquot(m_1)+n_1+\times \genfquot(m_2)+n_2+ \rightarrow \genfquot(m_1+m_2)+n_1 n_2+\qquad\text{fulfilling}\\[4pt]
 \sqmultquotpoint \funcomp \incl{\genfquot(m_1)+n_1+\times \genfquot(m_2)+n_2+}{\genfquot^2}= \incl{\genfquot(m_1+m_2)+n_1 n_2+}{\genfquot} \funcomp \sqmultquotgradpoint\text{.}
\end{multline*}
With an abuse of language, from now on, we denote $\sqmultquotgradpoint$ again by $\sqmultquotpoint$.
\item Fix $m, n_1,n_2 \in \mathbb{N}_0$. There is one and only one set function
\begin{multline*}
\multquotgradpoint :\genfquot(m)+n_1+\times \genfquot(m)+n_2+ \rightarrow \genfquot(m)+n_1 n_2+\qquad\text{fulfilling}\\[4pt]
 \multquotpoint \funcomp \incl{\genfquot(m)+n_1+\times \genfquot(m)+n_2+}{\genfquot^2}= \incl{\genfquot(m)+n_1 n_2+}{\genfquot} \funcomp \multquotgradpoint\text{.}
\end{multline*}
With an abuse of language, from now on, we denote $\multquotgradpoint$ again by $\multquotpoint$.
\item Fix $m,n \in \mathbb{N}_0$.
There is one and only one set function
\begin{multline*}
\scalpquotgradpoint :\mathbb{R}\times \genfquot(m)+n+ \rightarrow \genfquot(m)+n+\qquad\text{fulfilling}\\[4pt]
\scalpquotpoint \funcomp \left(\idobj+\mathbb{R}+,\incl{\genfquot(m)+n+}{\genfquot}\right)= \incl{\genfquot(m)+n+}{\genfquot} \funcomp \scalpquotgradpoint\text{.}
\end{multline*}
With an abuse of language, from now on, we denote $\scalpquotgradpoint$ again by $\scalpquotpoint$.
\item Fix $m,n \in \mathbb{N}_0$, $\bsfm \in \bsfM$. Refer to Remark \ref{intdiffmonrem}-[1]. \newline
If there is $ i \in \mathbb{N}$ such that $\bsfm = \Fpart_i$ or $\bsfm = \Fp_i $ then 
\begin{equation*}
\begin{array}{l}
\text{there is a set function}\;\bsfmuquotgradpoint\left[\bsfm\right] :\genfquot(m)+n+ \rightarrow \genfquot(m)+n+\\[4pt]
\text{defined by setting}\quad
 \bsfmuquotgradpoint\left[\bsfm\right] \left(\gggfx\right)= \bsfm \bsfmuquotpoint\gggfx\quad \forall  \gggfx \in \genfquot(m)+n+\text{.}
\end{array}
\end{equation*}
If there is $ i \in \mathbb{N}$ such that $\bsfm = \Fint_i$  then
\begin{equation*}
\begin{array}{l}
\text{there is a set function}\;\bsfmuquotgradpoint\left[\bsfm\right] :\genfquot(m)+n+ \rightarrow \genfquot(m+1)+n+\\[4pt]
\text{defined by setting}\quad
 \bsfmuquotgradpoint\left[\bsfm\right] \left(\gggfx\right)= \bsfm \bsfmuquotpoint \gggfx \quad \forall  \gggfx \in \genfquot(m)+n+\text{.}
\end{array}
\end{equation*}
With an abuse of language, from now on, we denote $\bsfmuquotgradpoint\left[\bsfm\right]$ by $\bsfm\bsfmuquotpoint$.
\item The 4-tuple $\left(\genfquot(\segnvar)+1+, \sumquotpoint, \sqmultquotpoint,  \scalpquotpoint\right)$ is a graded weak $\mathbb{R}$-algebra.
\item Fix $n\in \mathbb{N}_0$. The triplet $\left(\genfquot(\segnvar)+n+, \sumquotpoint, \sqmultquotpoint\right)$ is a graded weak $\genfquot(\segnvar)+1+$-module.
\item The 4-tuple $\left(\genfquot(\segnvar)+1+, \sumquotpoint, \multquotpoint,  \scalpquotpoint\right)$ is a graded commutative weak $\mathbb{R}$-algebra.
\item Fix $n\in \mathbb{N}_0$. The triplet $\left(\genfquot(\segnvar)+n+, \sumquotpoint, \multquotpoint\right)$ is a graded weak $\genfquot(\segnvar)+1+$-module.
\item Fix $m\in \mathbb{N}_0$. The 4-tuple $\left(\genfquot(m)+1+, \sumquotpoint, \multquotpoint,  \scalpquotpoint\right)$ is a commutative weak $\mathbb{R}$-algebra.
\item Fix $m,n\in \mathbb{N}_0$. The triplet $\left(\genfquot(m)+n+, \sumquotpoint, \multquotpoint\right)$ is a weak $\genfquot(m)+1+$-module.
\end{enumerate}
\end{proposition}
\begin{proof}
Statements follows directly by Proposition \ref{genpointprop}.
\end{proof}

In Definitions \ref{genpointdefmnCC0loc}, \ref{genpointdefmnCC00loc}, Remarks \ref{extremAnCC0}, \ref{extremAnCC00} below we specialize Propositions
\ref{genpointprop}, \ref{operAn} in case of generalized germs, of continuous pointed generalized germs and of locally vanishing continuous generalized germs.

\begin{definition}\mbox{}\label{genpointdefmnCC0loc} 
\begin{enumerate}
\item Fix $m,n \in \mathbb{N}_0$. We set:
\begin{flalign*}
&\genfquotcccz(m)+n+=\genfquot(m)+n+\cap \genfquotcccz\text{;}\\[4pt]
&\genfquotcccz(\segnvar)+n+=\genfquot(\segnvar)+n+\cap \genfquotcccz\text{.}
\end{flalign*}
\item  Fix $\topindexquattro,\topindexcinque\in \mathbb{N}$, $\left(m_{\elindexquattro}\right)_{\elindexquattro=1}^{\topindexquattro}\in {\mathbb{N}_0}^{\topindexquattro}$, $\left(n_{\elindexcinque}\right)_{\elindexcinque=1}^{\topindexcinque}\in {\mathbb{N}_0}^{\topindexcinque}$. 
Set $m=\underset{\elindexquattro=1}{\overset{\topindexquattro}{\sum}}\,m_{\elindexquattro}$, $n=\underset{\elindexcinque=1}{\overset{\topindexcinque}{\sum}}\,n_{\elindexcinque}$.\newline
We write $\genfquotcccz(\argcompl{\left(m_{\elindexquattro}\right)_{\elindexquattro=1}^{\topindexquattro}})+\argcompl{\left(n_{\elindexcinque}\right)_{\elindexcinque=1}^{\topindexcinque}}+$ in place of $\genfquotcccz(m)+n+$ whenever identifications $\mathbb{R}^m=\underset{\elindexquattro=1}{\overset{\topindexquattro}{\prod}}\,\mathbb{R}^{m_{\elindexquattro}}$, $\mathbb{R}^n=\underset{\elindexcinque=1}{\overset{\topindexcinque}{\prod}}\,\mathbb{R}^{n_{\elindexcinque}}$ are to be considered.
\end{enumerate}
\end{definition}

\begin{remark}\label{extremAnCC0}\mbox{}
\begin{enumerate}
\item Proposition \ref{genpointprop} holds true word by word for sets $\genfquotcccz(m)+n+$ by replacing $\genfquotcccz$ by $\genfquotcccz(m)+n+$, $\genfquotcccz^k$ by $\overset{k}{\underset{\mathsf{k}=1}{\prod}}\,\genfquotcccz(m_i)+n_i+$, 
$\genfquot$ by $\genfquot(m)+n+$, $\genfquot^k$ by $\overset{k}{\underset{\mathsf{k}=1}{\prod}}\,\genfquot(m_i)+n_i+$;\newline
Proposition \ref{operAn} hold true word by word by replacing symbol $\genfquot(m)+n+$ by $\genfquotcccz(m)+n+$, $\genfquot(\segnvar)+n+$ by $\genfquotcccz(\segnvar)+n+$, $\genfquot$ by $\genfquotcccz$.
\item Definition \ref{genpointdefmnloc} extends to spaces $\Ckspquot{\infty}$, $\Ckspquotccz{\infty}$, $\Ckspquot{0}$,  $\Ckspquotccz{0}$ by replacing symbol $\genfquot$ by $\Ckspquot{\infty}$, $\Ckspquotccz{\infty}$, $\Ckspquot{0}$,  $\Ckspquotccz{0}$ respectively, symbol $\dimcodmagtwo$ by $\dimCod$, symbol $\dimdommagtwo$ by $\dimDom$.\newline
Evaluation of a function at a point extends to germs of local smooth functions and to germs of local continuous functions by defining the set function
\begin{equation*}
\evalcompquot_n: \Ckspquot{0}(m)+n+\rightarrow \mathbb{R}^n\quad\text{by setting}\quad \evalcompquot_n\left(\gczsfx[f]\right)=f(0) \quad \forall \gczsfx[f]\in \Ckspquot{0}(m)+n+\text{.}
\end{equation*}
\end{enumerate}
\end{remark}

\begin{definition}\label{genpointdefmnCC00loc}
Fix $m,n \in \mathbb{N}_0$. We set $\genfquotccczz(m)+n+=\genfquot(m)+n+\cap \genfquotccczz$.
\end{definition}

\begin{remark}\label{extremAnCC00}
Proposition \ref{genpointprop} holds true word by word for sets $\genfquotcccz(m)+n+$ by replacing $\genfquotcccz$ by $\genfquotccczz(m)+n+$, $\genfquotccczz^k$ by $\overset{k}{\underset{\mathsf{k}=1}{\prod}}\,\genfquotccczz(m_i)+n_i+$, 
$\genfquot$ by $\genfquot(m)+n+$, $\genfquot^k$ by $\overset{k}{\underset{\mathsf{k}=1}{\prod}}\,\genfquot(m_i)+n_i+$;\newline
Proposition \ref{operAn} hold true word by word for sets $\genfquotccczz(m)+n+$ by replacing symbol $\genfquot(m)+n+$ by $\genfquotccczz(m)+n+$, $\genfquot$ by $\genfquotccczz$.
\end{remark}

In Definition \ref{socdef} and Proposition \ref{socprop} below we introduce a class of subsets of $\genfquotcccz(l)+m+$ which will play a central role in Section \ref{pdggf}.  

\begin{definition} \label{socdef}
Fix $l,m\in\mathbb{N}_0$, a subset $ \setofcores \subseteq \genfquotcccz(l)+m+$. We sat that $\setofcores$ is an admissible subset of $\genfquotcccz(l)+m+$ if and only if 
\begin{equation*}
\occset_{\gggfz_1 }=\occset_{\gggfz_2 }
\qquad \forall \gggfz_1, \gggfz_2 \in \setofcores\text{.}
\end{equation*}
\end{definition}

\begin{proposition}\label{socprop}
Fix $l,m,n\in\mathbb{N}_0$, an admissible subset $ \setofcores$ of $\genfquotcccz(l)+m+$, $\gggfx\in\genfquotcccz(m)+n+$. Then set $\quad\gggfx \compquotpoint\setofcores=\left\{\gggfx\compquotpoint  \gggfz\;:\;\gggfz\in \setofcores\right\}\quad$ is an admissible subset of  $\genfquotcccz(l)+n+$.
\end{proposition}
\begin{proof} The proof follows by definition of $\compquotpoint$ and of admissible subset. 
\end{proof}

\chapter{Localization of generalized derivations\label{Derggf}}

In this chapter we study invariance of differential operators with respect to local change of coordinates.

\section{Pre-derivations of germs of generalized functions\label{pdggf}}

In Definition \ref{genderdef} below we introduce the notion of pre-derivation: a suitable continuous $\mathbb{R}$-linear $\genfquot(\segnvar)+1+$-valued function on $\genfquot(m)+1+$. We refer to Definitions  \ref{genpointdefmnCC0loc}, \ref{socdef}, Remark \ref{Cinfdentro}-[2].

\begin{definition}\label{genderdef}\mbox{}
\begin{enumerate}
\item Fix $m\in\mathbb{N}_0$, a set function $\genpreder:\genfquot(m)+1+\rightarrow \genfquot(\segnvar)+1+$.\newline
We say that $\genpreder$ is a pre-derivation if and only if there are $l\in \mathbb{N}_0$, a vector $\vecuno\in \mathbb{R}^{l}$, $\gggfz \in \genfquotcccz(l)+m+$,  such that:
\begin{equation}
\left\{
\begin{array}{l}
(i)\hspace{7pt}\genpreder\left(\gggfw\right)=\underset{\mathsf{l}=1}{\overset{l}{\bigsumquotpoint}}\, \left( \Fpart_{\mathsf{l}} \bsfmuquotpoint \left(\gggfw\compquotpoint \gggfz\right)\right)\scalpquotpoint  \vecuno_{\mathsf{l}}
\hspace{10pt} \forall \gggfw \in \genfquot(m)+1+
\hspace{15pt}\text{if}\;l\neq 0\text{,}\\[8pt]
(ii)\hspace{5pt}\genpreder\left(\gggfw\right)= \left(\gggfw\compquotpoint \gggfz\right)\scalpquotpoint 0 \hspace{52pt} \forall \gggfw \in \genfquot(m)+1+\hspace{15pt} \text{if}\;l= 0\text{.}
\end{array}
\right.\label{predereq}
\end{equation}
We say that: $\gggfz$ is a core of $\genpreder$; $\vecuno$ is an associating direction of $\genpreder$; $\left(\gggfz, \vecuno\right)$ is a skeleton of $\genpreder$.\newline
To emphasize the dependence of a pre-derivation on a specific skeleton we write $\genpreder\left[\gggfz, \vecuno \right]$ whenever needed.
\item  Fix $l,m\in\mathbb{N}_0$, a $\mathbb{R}$-vector subspace $\modsym$ of $\mathbb{R}^l$, an admissible subset $ \setofcores$ of $\genfquotcccz(l)+m+$.\newline
We define the additive abelian group $\left(\genprederspace[\modsym](\setofcores)+m+,\sumgenpreder[\modsym](\setofcores)+m+\right)$ by setting\\[8pt]
\textbf{Generators:}  
\begin{equation*}
\big\{\genpreder\;:\hspace{7pt}\forall \gggfz\in  \setofcores\hspace{7pt}\exists \vecuno \in \modsym \hspace{5pt}\text{such}\hspace{5pt}\text{that}\hspace{5pt}\left(\gggfz,\vecuno\right)\hspace{5pt}\text{is}\hspace{5pt}\text{a}\hspace{5pt}\text{skeleton}\hspace{5pt}\text{of}\hspace{5pt}\genpreder\big\}\text{;}
\end{equation*}
\textbf{Relations:}
\begin{equation}
\genpreder\left[\gggfz, \vecuno_{1} \right]\;\sumgenpreder[\modsym](\setofcores)+m+\;\genpreder\left[\gggfz, \vecuno_{2} \right]=\genpreder\left[\gggfz, \vecuno_{1}+\vecuno_{2}\right]\text{.}\label{relpreder}
\end{equation}
\end{enumerate}
\end{definition}

\begin{remark}\label{geosignpreder}
Here we describe the geometric meaning of pre-derivations. A pre-derivation $\genpreder\left[\gggfz, \vecuno \right]$ is the derivation operator with respect to what can be considered a generalized vector of $\mathbb{R}^m$, that is the image of $\vecuno$ through $\gggfz$. 
\end{remark}

In Proposition \ref{algstrprederpoint} below we prove algebraic properties of pre-derivations and we associate to any continuous pointed generalized germ a suitable notion of differential. We refer to Definitions \ref{genpointdef}, \ref{socdef}, Proposition \ref{socprop}.

\begin{proposition}\label{algstrprederpoint}\mbox{}
\begin{enumerate}
\item Fix $m\in\mathbb{N}_0$, two continuous pointed generalized germs $\gggfz_1$, $\gggfz_2$, a pre-derivation $\genpreder:\genfquot(m)+1+\rightarrow \genfquot(\segnvar)+1+$. If $\gggfz_1$ and $\gggfz_2$ are both cores of $\genpreder$ then
\begin{equation*}
\occset_{\gggfz_1}=\occset_{\gggfz_2}\text{,} \qquad
\dimdommagtwo\left(\gggfz_1 \right)=\dimdommagtwo\left(\gggfz_2 \right)\text{.}
\end{equation*}
We emphasize that this result justifies Definition \ref{socdef} of admissible subset.
\item Fix $l,m\in\mathbb{N}_0$, a $\mathbb{R}$-vector subspace $\modsym$ of $\mathbb{R}^l$, an admissible subset $\setofcores$ of  $\genfquotcccz(l)+m+$. Then: 
\begin{flalign*}
&\left\{
\begin{array}{l}
\text{there is a set function} \\
\hspace{47pt}\scalpgenpreder[\modsym](\setofcores)+m+:\mathbb{R} \times\genprederspace[\modsym](\setofcores)+m+ \rightarrow \genprederspace[\modsym](\setofcores)+m+\\[4pt]
\text{defined by setting}\\[4pt]
 \hspace{10pt}a\;\scalpgenpreder[\modsym](\setofcores)+m+\; \genpreder\left[\gggfz, \vecuno \right]=
\genpreder\left[\gggfz, a\vecuno\right]\\[4pt] 
\hspace{101pt}\forall a \in\mathbb{R}\text{,}\quad\forall\genpreder\left[\gggfz, \vecuno \right] \in \genprederspace[\modsym](\setofcores)+m+\text{;}
\end{array}
\right.\\[8pt]
&\begin{array}{l}
\left(\genprederspace[\modsym](\setofcores)+m+, \sumgenpreder[\modsym](\setofcores)+m+, \scalpgenpreder[\modsym](\setofcores)+m+\right)\;\text{is an $\mathbb{R}$-vector space.}
\end{array}
\end{flalign*} 
With an abuse of notation we denote $\sumgenpreder[\modsym](\setofcores)+m+$ simply by $\sumgenpreder$, $\scalpgenpreder[\modsym](\setofcores)+m+$ simply by $\scalpgenpreder$, whenever no confusion is possible.
\item Fix $l,m,n\in\mathbb{N}_0$, a $\mathbb{R}$-vector subspace $\modsym$ of $\mathbb{R}^l$, $\gggfx \in \genfquotcccz(m)+n+$, an admissible subset $\setofcores$ of  $\genfquotccczz(l)+m+$. Then there is an $\mathbb{R}$-linear function
\begin{equation*}
\genprediff[\modsym](\setofcores)+\gggfx+ : \genprederspace[\modsym](\setofcores)+m+\rightarrow \genprederspace[\modsym]( \gggfx  \setofcores)+n+
\end{equation*}
defined by setting
\begin{equation*}
\genprediff[\modsym](\setofcores)+\gggfx+\left(\genpreder\left[\gggfz, \vecuno\right]\right)\!=\!
\genpreder\left[\gggfx \compquotpoint\gggfz, \vecuno\right]\hspace{10pt}
 \forall\genpreder\left[\gggfz,\vecuno \right]\in \genprederspace[\modsym](\setofcores)+m+\text{.}
\end{equation*}
We say that $\genprediff[\modsym](\setofcores)+\gggfx+$ is the pre-differential of $\gggfx$ at $\setofcores$.
\item Fix $l,m,n\in\mathbb{N}_0$, two $\mathbb{R}$-vector subspaces $\modsym_1$, $\modsym_2$ of $\mathbb{R}^l$, $\gggfx \in \genfquotcccz(m)+n+$, two admissible subsets $\setofcores_{1}$, $\setofcores_{2}$ of  $\genfquotccczz(l)+m+$. Assume that both the following conditions hold true:  $\modsym_1\subseteq\modsym_2$; $\setofcores_{1} \subseteq \setofcores_{2}$.
Then:
\begin{flalign*}
&\begin{array}{l} \genprederspace[\modsym_1](\setofcores_{2})+m+\subseteq  \genprederspace[\modsym_2](\setofcores_{1})+m+\text{;}\end{array}\\[8pt]
&\begin{array}{l} \incl{\genprederspace[\modsym_1](\setofcores_{2})+m+}{\genprederspace[\modsym_2](\setofcores_{1})+m+} \;\text{is an $\mathbb{R}$-linear function;}\end{array}\\[8pt]
&\begin{array}{l}  \genprediff[\modsym_2](\setofcores_{1})+\gggfx+\funcomp \incl{\genprederspace[\modsym_1](\setofcores_{2})+m+}{\genprederspace[\modsym_2](\setofcores_{1})+m+}= \\
\hspace{45pt}\incl{\genprederspace[\modsym_1](\gggfx  \setofcores_{2})+n+}{\genprederspace[\modsym_2](\gggfx \setofcores_{1})+n+}\funcomp\genprediff[\modsym_1](\setofcores_{2})+\gggfx+\text{.}
\end{array}
\end{flalign*}
\end{enumerate}
\end{proposition}
\begin{proof} Statements 2, 3, 4 follow directly by checking axioms. We prove statement 1 by arguing as follows. Fix vectors $\vecuno_{1},\vecuno_{2}\in \mathbb{R}^{l}$ such that $\genpreder=\genpreder\left[\gggfz_{i}, \vecuno_{i} \right]$ for any $i \in \left\{1,2\right\}$. Referring to Proposition \ref{genpointprop}-[10] statement follows by exploiting equality\newline
\centerline{$\genpreder\left[\gggfz_{1}, \vecuno_{1} \right]\left(\lininclquot\left(\gczsfx[f]\right)\right)=\genpreder\left[\gggfz_{2}, \vecuno_{2} \right]\left(\lininclquot\left(\gczsfx[f]\right)\right)$}\newline
for $f$ real valued smooth function defined in an open neighborhood of $0 \in \mathbb{R}^m$ and fulfilling $f\left(0\right)=0$.
\end{proof}

Motivated by Proposition \ref{algstrprederpoint}, in Definition \ref{genderdef2} below we organize the material introduced above in the functorial language in order to make easier constructions by exploiting co-limits. We refer to Notations \ref{catfun}-[1, 4], \ref{alg}-[10].

\begin{definition} \label{genderdef2}\mbox{}
\begin{enumerate}
\item Fix $l,m \in \mathbb{N}_0$. We denote by $\adsetcat(l)+m+$ the category defined as follows.\newline
Objects of $\adsetcat(l)+m+$ are admissible subsets of $\genfquotcccz(l)+m+$.\newline
Arrows of $\adsetcat(l)+m+$ are inclusions of admissible subsets.
\item Fix $l, m, n \in \mathbb{N}_0$,
$\gggfx \in \genfquotcccz(m)+n+$. We denote by $\genpredersoc(l)+\gggfx+:\adsetcat(l)+m+\rightarrow\adsetcat(l)+n+$ the functor defined by setting 
\begin{flalign*}
&\begin{array}{l}
\genpredersoc(l)+\gggfx+\left(\setofcores\right)=\gggfx \compquotpoint  \setofcores \hspace{4pt}\text{for any object}\;\setofcores\;\text{of}\;\adsetcat(l)+m+\text{;} 
\end{array}\\[4pt]
&\begin{array}{l}
\genpredersoc(l)+\gggfx+\left(\incl{\setofcores_1}{\setofcores_2}\right)= \incl{\gggfx \compquotpoint\setofcores_1}{\gggfx \compquotpoint\setofcores_2}\hspace{4pt}
\text{for any arrow}\;\incl{\setofcores_1}{\setofcores_2}\;\text{of}\;\adsetcat(l)+m+ \text{.}
\end{array}
\end{flalign*}
\item Fix $l,m\in \mathbb{N}_0$, a $\mathbb{R}$-vector subspace $\modsym$ of $\mathbb{R}^l$.\newline
We denote by $\prederFunc[\modsym](l)+m+:\adsetcat(l)+m+\rightarrow \VecR$ the contravariant functor defined by setting 
\begin{flalign*}
&\begin{array}{l}
\prederFunc[\modsym](l)+m+\left(\setofcores\right)=
\genprederspace[\modsym](\setofcores)+m+ \hspace{4pt}\text{for}\hspace{4pt}\text{any}\hspace{4pt}\text{object}\hspace{4pt}\setofcores\hspace{4pt}\text{of}\hspace{4pt}\adsetcat(l)+m+\text{,} 
\end{array}\\[4pt]
&\left\{
\begin{array}{l}
\prederFunc[\modsym](l)+m+\left(\incl{\setofcores_1}{\setofcores_2}\right)= \incl{\genprederspace[\modsym](\setofcores_2)+m+}{\genprederspace[\modsym](\setofcores_1)+m+}\\[4pt]
\text{for any arrow}\;\incl{\setofcores_1}{\setofcores_2}\;\text{of}\;\adsetcat(l)+m+ \text{.}
\end{array}
\right.
\end{flalign*}
\item Fix $l,m,n\in\mathbb{N}_0$, a $\mathbb{R}$-vector subspace $\modsym$ of $\mathbb{R}^l$, $\gggfx \in \genfquotcccz(m)+n+$.\newline
We denote by $\genprediffNT[\modsym](l)+\gggfx+: \prederFunc[\modsym](l)+m+\rightarrow \left(\prederFunc[\modsym](l)+n+\funcomp\genpredersoc(l)+\gggfx+\right)$ the natural arrow defined by setting
\begin{equation*}
\genprediffNT[\modsym](l)+\gggfx+\left(\setofcores\right) =\genprediff[\modsym](\setofcores)+\gggfx+ \hspace{4pt} \text{for}\hspace{4pt}\text{any}\hspace{4pt}\text{object}\hspace{4pt}\setofcores\hspace{4pt}\text{of}\hspace{4pt}\adsetcat(l)+m+\text{.}
\end{equation*}
\item Fix $l,m\in\mathbb{N}_0$, $\mathbb{R}$-vector subspaces $\modsym_1$, $\modsym_2$ of $\mathbb{R}^l$. Assume $\modsym_1\subseteq \modsym_2$.\newline
We denote by $\genpreinclNT<\modsym_1<>\modsym_2>(l)+m+: \prederFunc[\modsym_1](l)+m+\rightarrow \prederFunc[\modsym_2](l)+m+$ the natural arrow defined by setting
\begin{multline*}
\genpreinclNT<\modsym_1<>\modsym_2>(l)+m+\left(\setofcores\right) = \\
\incl{\genprederspace[\modsym_1](\setofcores)+m+}{\genprederspace[\modsym_2](\setofcores)+m+} \hspace{4pt} \text{for}\hspace{4pt}\text{any}\hspace{4pt}\text{object}\hspace{4pt}\setofcores\hspace{4pt}\text{of}\hspace{4pt}\adsetcat(l)+m+\text{.}
\end{multline*}
\item Fix $l,m\in \mathbb{N}_0$, a $\mathbb{R}$-vector subspace $\modsym$ of $\mathbb{R}^l$.\newline
We define the pair $\left(  \genprederspace[\modsym](l)+m+, \inclcol[\modsym](l)+m+\right)$ as the co-limit of the functor $\prederFunc[\modsym](l)+m+$. We say that $\genprederspace[\modsym](l)+m+$ is the space of $m$-dimensional $\left(\modsym, l\right)$-pre-derivations. Elements belonging to $\genprederspace[\modsym](l)+m+$ are denoted by $\genpreder$ and called $m$-dimensional $\left(\modsym, l\right)$-pre-derivations, with an abuse of language. We drop any reference to $m$, $\modsym$ or $l$ whenever no confusion is possible.\newline
\item Fix $l, m, n \in \mathbb{N}_0$, a $\mathbb{R}$-vector subspace $\modsym$ of $\mathbb{R}^l$,
$\gggfx \in \genfquotcccz(m)+n+$.\newline
We define the arrow $\genprediff[\modsym](l)+\gggfx+ : \genprederspace[\modsym](l)+m+\rightarrow \genprederspace[\modsym](l)+n+$ as the unique arrow induced on co-limts by the natural arrow $\genprediffNT[\modsym](l)+\gggfx+$. We say that $\genprediff[\modsym](l)+\gggfx+$ is the $\left(\modsym, l\right)$-pre-differential of $\gggfx$ at $0$.
\item Fix $l,m\in\mathbb{N}_0$, $\mathbb{R}$-vector subspaces $\modsym_1$, $\modsym_2$ of $\mathbb{R}^l$. Assume $\modsym_1\subseteq \modsym_2$.\newline
We define the arrow $\genpreincl<\modsym_1<>\modsym_2>(l)+m+ : \genprederspace[\modsym_1](l)+m+\rightarrow \genprederspace[\modsym_2](l)+m+$ as the unique arrow induced on co-limts by the natural arrow $\genpreinclNT<\modsym_1<>\modsym_2>(l)+m+$.
\end{enumerate}
\end{definition}

In Remark \ref{pathinlpreder} below we explicitly describe the structure of the co-limit pair $\left(  \genprederspace[\modsym](l)+m+, \inclcol[\modsym](l)+m+\right)$ introduced in Definition \ref{genderdef2}. We emphasize that we are working in the category $\VecR$ of $\mathbb{R}$-vector space and linear functions and we refer to Notation \ref{gentopnot}-[6], Definition \ref{genpointdefmnCC0loc}, Proposition \ref{genpointprop}.

\begin{remark}\label{pathinlpreder}\mbox{} 
\begin{enumerate}
\item $\genprederspace[\modsym](l)+m+$ is a $\mathbb{R}$-vector space which is not trivial since Proposition \ref{algstrprederpoint}-[1] entails that there are at least two different pre-derivations defined through \eqref{predereq}. Sum and scalar product of $\genprederspace[\modsym](l)+m+$ are denoted by $\sumgenpreder[\modsym](l)+m+$ and $\scalpgenpreder[\modsym](l)+m+$ respectively. With an abuse of notation we denote $\sumgenpreder[\modsym](l)+m+$ simply by $\sumgenpreder$, $\scalpgenpreder[\modsym](l)+m+$ simply by $\scalpgenpreder$, whenever no confusion is possible.
\item $\genprederspace[\modsym](l)+m+=\frac{\dirsumgenpderspace[\modsym](l)+m+}{\nulldirsumgenpderspace[\modsym](l)+m+}$ where: 
\begin{flalign}
&\begin{array}{l}
\dirsumgenpderspace[\modsym](l)+m+=\underset{\gggfz\in  \genfquotcccz(l)+m+}{\bigdirsum}\genprederspace[\modsym](\left\{ \gggfz \right\})+m+\text{;}
\end{array}\label{Mdef}\\[8pt]
&\left\{\begin{array}{l}
\nulldirsumgenpderspace[\modsym](l)+m+\hspace{3pt}\text{is}\hspace{3pt}\text{the}\hspace{3pt} \mathbb{R}\text{-vector}\hspace{3pt}\text{subspace}\hspace{3pt}\text{of}\hspace{3pt}\dirsumgenpderspace[\modsym](l)+m+\;\text{generated}\\[4pt]
\text{by}\hspace{3pt}\text{the}\hspace{3pt}\text{set}\\[4pt]
\big\{
\genpreder\left[\gggfz_1, \vecuno_1 \right]-\genpreder\left[\gggfz_2, \vecuno_2 \right]\hspace{4pt}:\\[4pt]
\hspace{20pt} \forall \gggfz_1, \gggfz_2\in  \genfquotcccz(l)+m+\text{,}\hspace{4pt}\forall\vecuno_1,\vecuno_2 \in \modsym \hspace{4pt}\text{with}\\[4pt]
\hspace{30pt}\genpreder\left[\gggfz_1, \vecuno_1 \right]\left(\gggfw\right)=\genpreder\left[\gggfz_2, \vecuno_2 \right]\left(\gggfw\right)
\quad  \forall\gggfw \in \genfquot(m)+1+
\big\}\text{.}
\end{array}
\right.\label{Mcorsdef}
\end{flalign}
We denote by $\quotdirsumuno[\modsym](l)+m+ :\dirsumgenpderspace[\modsym](l)+m+\rightarrow \genprederspace[\modsym](l)+m+$ the quotient function. We drop any reference to $\modsym$, $l$ or $m$ whenever there is no risk of confusion or there is no need of such a detail.
\item Natural arrow $\inclcol[\modsym](l)+m+:\prederFunc[\modsym](l)+m+\rightarrow \costfun<\argcompl{\adsetcat(l)+m+}<>\VecR>+\argcompl{\genprederspace[\modsym](l)+m+}+$ consists of $\mathbb{R}$-linear functions $\inclcol[\modsym](l)+m+\left(\setofcores\right): \prederFunc[\modsym](l)+m+\left(\setofcores\right)\rightarrow \genprederspace[\modsym](l)+m+$ for any object $\setofcores$ of $\adsetcat(l)+m+$.
\item The $\mathbb{R}$-linear function $\genprediff[\modsym](l)+\gggfx+$ is the unique function induced on quotient spaces by the $\mathbb{R}$-linear function
\begin{equation*}
\dirsumgenprediff[\modsym](l)+\gggfx+=\underset{\gggfz\in  \genfquotcccz(l)+m+}{\bigdirsum}\genprediff[\modsym](\left\{\gggfz\right\})+\gggfx+\text{.}
\end{equation*}
\item The $\mathbb{R}$-linear function $\genpreincl<\modsym_1<>\modsym_2>(l)+m+$ is the unique function induced on quotient spaces by the $\mathbb{R}$-linear function 
\begin{multline*}
\dirsumgenpreincl<\modsym_1<>\modsym_2>(l)+m+=\\
\underset{\gggfz\in  \genfquotcccz(l)+m+}{\bigdirsum}
\incl{\genprederspace[\modsym_1](\left\{\gggfz\right\})+m+}{\genprederspace[\modsym_2](\left\{\gggfz\right\})+m+}
\text{.}
\end{multline*}
\item Fix a pre-derivation $\genpreder\in \genprederspace[\modsym](l)+m+$. Then there are a finite set $\setindexdue$, a pair $\left(  \gggfz_{\elindexdue},\vecuno_{\elindexdue}\right)\in\genfquotcccz(l)+m+\times\modsym$ for any $\elindexdue\in\setindexdue$ fulfilling 
\begin{equation*}
\genpreder=\underset{\elindexdue\in  \setindexdue}{\sum} \left(\inclcol[\modsym](l)+m+\left(\left\{\gggfz_{\elindexdue}\right\}\right)\right)\left(\genpreder\left[\gggfz_{\elindexdue},\vecuno_{\elindexdue} \right]\right)\text{.}
\end{equation*} 
We set:
\begin{flalign*}
&\begin{array}{l}
\left( \gggfz_{\setindexdue},\vecuno_{\setindexdue}\right)=\left\{\left( \gggfz_{\elindexdue},\vecuno_{\elindexdue}\right)\right\}_{\elindexdue\in  \setindexdue}\text{;}
\end{array}\\[6pt]
&\begin{array}{l}
\genpreder\left[\gggfz_{\setindexdue},\vecuno_{\setindexdue} \right]=\left\{\genpreder\left[\gggfz_{\elindexdue},\vecuno_{\elindexdue} \right]\right\}_{\elindexdue\in  \setindexdue}\text{;}
\end{array}\\[6pt]
&\begin{array}{l}
\genpreder_{\elindexdue}=\left(\inclcol[\modsym](l)+m+\left(\left\{\gggfz_{\elindexdue}\right\}\right)\right)\left(\genpreder\left[\gggfz_{\elindexdue},\vecuno_{\elindexdue} \right]\right)\hspace{10pt}\forall \elindexdue \in  \setindexdue\text{;}
\end{array}\\[6pt]
&\begin{array}{l}
\genpreder_{\setindexdue}=\left\{\genpreder_{\elindexdue}\right\}_{\elindexdue\in  \setindexdue}\text{.}
\end{array}
\end{flalign*} 
We say that:\newline
$\left( \gggfz_{\setindexdue},\vecuno_{\setindexdue}\right)$ is a skeleton of $\genpreder$;\newline
$\genpreder\left[\gggfz_{\setindexdue},\vecuno_{\setindexdue} \right]$ is a realization of $\genpreder$;\newline
$\genpreder_{\setindexdue}$ is a decomposition of $\genpreder$.\newline
If $\left( \gggfz_{\setindexdue},\vecuno_{\setindexdue}\right)$ fulfills also condition
\begin{equation*}
\begin{array}{l}\genpreder\left[\gggfz_{\elindexdue_1},\vecuno_{\elindexdue_1} \right]\neq \genpreder\left[\gggfz_{\elindexdue_2},\vecuno_{\elindexdue_2} \right]\hspace{15pt}\forall \elindexdue_1, \elindexdue_2 \in \setindexdue \hspace{4pt}\text{with}\hspace{4pt}\elindexdue_1\neq\elindexdue_2 
\end{array}
\end{equation*}
then we say that:\newline
$\left( \gggfz_{\setindexdue},\vecuno_{\setindexdue}\right)$ is a canonical skeleton of $\genpreder$;\newline
$\genpreder\left[\gggfz_{\setindexdue},\vecuno_{\setindexdue} \right]$ is a canonical realization of $\genpreder$;\newline
$\genpreder_{\setindexdue}$ is a canonical decomposition of $\genpreder$.\newline
Whenever will be convenient we will choose 
\begin{equation*}
\setindexdue = \genfquotcccz(l)+m+\hspace{15pt}\text{and}\hspace{15pt}\gggfz_{\elindexdue}=\elindexdue\hspace{10pt}\forall\elindexdue\in \setindexdue\text{.}
\end{equation*}
\item Fix a pre-derivation $\genpreder \in \genprederspace[\modsym](l)+m+$, a skeleton $\left( \gggfz_{\setindexdue},\vecuno_{\setindexdue}\right)$ of $\genpreder$, $\gggfw\in \genfquot(m)+1+$ we set
\begin{equation*}
\genpreder\left(\gggfw\right)=\underset{\elindexdue\in  \setindexdue}{\bigsumquotpoint} \genpreder\left[\gggfz_{\elindexdue},\vecuno_{\elindexdue} \right]\left(\gggfw\right)\text{.}
\end{equation*} 
We emphasize that $\genpreder\left(\gggfw\right)\in \genfquot(l)+1+$ is independent on the chosen skeleton of $\genpreder$.
\end{enumerate}
\end{remark}

In Definition \ref{gpdssm} below we introduce the $\mathbb{R}$-subspace of $\genprederspace[\modsym](l)+m+$ containing all pre-derivations with a smooth core.

\begin{definition}\label{gpdssm} 
Fix $l,m \in\mathbb{N}_0$, a $\mathbb{R}$-vector subspace $\modsym$ of $\mathbb{R}^l$. We denote by $\genprederspacesm[\modsym](l)+m+$ the $\mathbb{R}$-vector subspace of $\genprederspace[\modsym](l)+m+$ generated by the set
\begin{equation*}
\Ggenprederspacesm[\modsym](l)+m+\!=\!\left\{\quotdirsumuno[\modsym](l)+m+\left(\genpreder\left[ \lininclquot\left(\gczsfx[f] \right)  ,\vecuno\right]\right)\;:\; \hspace{4pt}\forall\vecuno\in\modsym\text{,}\hspace{4pt} \forall\gczsfx[f] \in \Ckspquot{\infty}(l)+m+\right\}\text{.}
\end{equation*}
\end{definition}

In Proposition \ref{gpdssmprop} below we prove the existence of a unique $\mathbb{R}$-linear function evaluating pre-derivations with smooth core at $0$.    
\begin{proposition}\label{gpdssmprop}\mbox{}
\begin{enumerate}
\item Fix $l,m \in\mathbb{N}_0$, a $\mathbb{R}$-vector subspace $\modsym$ of $\mathbb{R}^l$. There is one and only one $\mathbb{R}$-linear function $\evspd[\modsym](l)+m+:\genprederspacesm[\modsym](l)+m+\rightarrow \mathbb{R}^m$ defined by setting
\begin{equation}
\left\{
\begin{array}{l}
\evspd[\modsym](l)+m+\left(\quotdirsumuno\left(\genpreder\left[ \lininclquot\left(\gczsfx[f] \right)  ,\vecuno\right]\right)\right)=\left(\smder_{\vecuno}f\right)\left(0\right)\\[4pt]
\text{for}\hspace{4pt}\text{any}\hspace{4pt}\vecuno\in\modsym\text{,}\hspace{4pt} \gczsfx[f]  \in \Ckspquot{\infty}(l)+m+\text{.}
\end{array}
\right.\label{defevspd}
\end{equation}
\item Fix $l,m \in\mathbb{N}_0$, two $\mathbb{R}$-vector subspaces $\modsym_1$, $\modsym_2$ of $\mathbb{R}^l$. Assume $\modsym_1 \subseteq \modsym_2$. Then $\evspd[\modsym_2](l)+m+\funcomp\genpreincl<\modsym_1<>\modsym_2>(l)+m+=\evspd[\modsym_1](l)+m+$.
\end{enumerate}
\end{proposition}
\begin{proof}\newline
Statement 1  follows straightforwardly by proving that $\evspd[\modsym](l)+m+$ is well defined on generators by \eqref{defevspd}. This holds true since, by definition, any element belonging to $\Ggenprederspacesm[\modsym](l)+m+$ defines the same pre-derivation $\genfquot(m)+1+\rightarrow \genfquot(\segnvar)+1+$ regardless of skeleton.\newline
Statement 2 follows by direct computation.   
\end{proof}

\begin{remark}\label{noM}
If $\modsym=\mathbb{R}^l$ then symbol $\modsym$ will be dropped in the notation.
\end{remark}

In Proposition \ref{struttint} below characterize skeletons defining the same pre-derivation.
We refer to Notations \ref{magtwopartic}-[2], \ref{difrealfunc}-[3], Definitions \ref{pathtwodef}-[1], \ref{infnearpointdef}-[2], Proposition \ref{algstrprederpoint}-[1].

\begin{proposition}\label{struttint}
Fix $l,m\in\mathbb{N}_0$, vectors $\vecuno_1,\vecuno_2\in \mathbb{R}^{l}$, an open neighborhood $\intuno\subseteq\mathbb{R}^l$ of $0 \in \mathbb{R}^l$, $\gfz_1,\gfz_2\in \genflocccz(\intuno)+\mathbb{R}^m+$,  paths $\pmtwo_1$, $\pmtwo_2$ in $\genf$ detecting $\gfz_1$, $\gfz_2$ respectively.\newline
Assume that both conditions below hold true:
\begin{flalign}
&\begin{array}{l}
\genpreder\left[\gggfz[\gfz_1], \vecuno_1 \right]=\genpreder\left[\gggfz[\gfz_2], \vecuno_2 \right]\text{;}\label{ipoteqprediff}
\end{array}\\[6pt]
&\begin{array}{l}
\pmtwo_1\text{,}\hspace{4pt}\pmtwo_2\hspace{4pt}\text{have}\hspace{4pt}\text{the}\hspace{4pt}\text{same}\hspace{4pt}\text{associating}\hspace{4pt}\text{function.}\label{ipoteqprediffdue}
\end{array}
\end{flalign}
Then there are $\delta>0$, an open neighborhood $\setsymnove \sqsubseteq  \intuno$ of $0 \in \intuno$ fulfilling all conditions below:
\begin{flalign}
&\begin{array}{l}
\left(\evalcomptwo\left(\pmtwo_1\left(\unkdue\right)\right)\right)\left(\unkuno\right)\neq\udenunk\hspace{15pt}\forall\left(\unkdue,\unkuno\right)\in \left(-\delta,\delta\right)\times\setsymnove\setminus\left\{0\right\}\times\setsymnove\text{;}\label{novuouno}
\end{array}\\[6pt]
&\begin{array}{l}
\left(\evalcomptwo\left(\pmtwo_2\left(\unkdue\right)\right)\right)\left(\unkuno\right)\neq\udenunk\hspace{15pt}\forall\left(\unkdue,\unkuno\right)\in \left(-\delta,\delta\right)\times\setsymnove\setminus\left\{0\right\}\times\setsymnove\text{;}\label{novuodue}
\end{array}
\\[6pt]
&\left\{
\begin{array}{l}
\left(\evalcomptwo\left(\pmtwo_1\left(\unkdue\right)\right)\right)\left(\unkuno\right)=\left(\evalcomptwo\left(\pmtwo_2\left(\unkdue\right)\right)\right)\left(\unkuno\right)\Rightarrow\\[4pt]
\hspace{85pt}\left(\smder_{\vecuno_1}\left(
\evalcomptwo\left(\pmtwo_1\left(\unkdue\right)\right)
\right)\right)\left(\unkuno\right)=\left(\smder_{\vecuno_2}\left(\evalcomptwo\left(\pmtwo_2\left(\unkdue\right)\right)\right)\right)\left(\unkuno\right)\\[4pt]
\forall\left(\unkdue,\unkuno\right)\in \left(-\delta,\delta\right)\times\setsymnove\setminus\left\{0\right\}\times\setsymnove\text{;}\label{tesitteq}
\end{array}
\right.\\[6pt]
&\left\{
\begin{array}{l}
\left(\evalcomptwo\left(\pmtwo_1\left(\unkdue\right)\right)\right)\left(\unkuno\right)\neq\left(\evalcomptwo\left(\pmtwo_2\left(\unkdue\right)\right)\right)\left(\unkuno\right)\Rightarrow\\[4pt]
\hspace{65pt}\left(\smder_{\vecuno_1}\left(\evalcomptwo\left(\pmtwo_1\left(\unkdue\right)\right)\right)\right)\left(\unkuno\right)=\left(\smder_{\vecuno_2}\left(\evalcomptwo\left(\pmtwo_2\left(\unkdue\right)\right)\right)\right)\left(\unkuno\right)=0\\[4pt]
\forall\left(\unkdue,\unkuno\right)\in \left(-\delta,\delta\right)\times\setsymnove\setminus\left\{0\right\}\times\setsymnove\text{.}\label{tesittdiv}
\end{array}
\right.
\end{flalign}
\end{proposition}
\begin{proof} 
Proposition \ref{pathprop7} entails that there are $\delta>0$, an open neighborhood $\setsymnove \sqsubseteq  \intuno$ of $0 \in \intuno$ fulfilling both \eqref{novuouno} and \eqref{novuodue}.\newline
Arguing by contradiction we prove that there are $\delta>0$, an open neighborhood $\setsymnove \sqsubseteq  \intuno$ of $0 \in \intuno$ fulfilling also both \eqref{tesitteq} and \eqref{tesittdiv}.\newline
\textnormal{\textbf{Proof of \eqref{tesitteq}.}}\ \ Arguing by contradiction fix a sequence $\left\{\left(\unkdue_{\topindexuno},\unkuno_{\topindexuno}\right)\right\}_{\topindexuno \in \mathbb{N}}\subseteq \left(-1,1\right)\times\intuno\setminus\left\{0\right\}\times\intuno$ fulfilling one and only one among the two systems of conditions below:
\begin{flalign}
&\left\{
\begin{array}{ll}
(i)&\underset{\topindexuno\rightarrow +\infty}{\lim}\,\left(\unkdue_{\topindexuno},\unkuno_{\topindexuno}\right)=\left(0,0\right)\text{,}\\[4pt]
(ii)&\left(\evalcomptwo\left(\pmtwo_{\topindextre}\left(\unkdue_{\topindexuno}\right)\right)\right)\left(\unkuno_{\topindexuno}\right) \neq \udenunk \hspace{15pt}\forall \topindextre\in\left\{1,2\right\}\text{,}\hspace{10pt}\forall \topindexuno\in\mathbb{N}\text{,}\\[4pt]
(iii)&\left(\evalcomptwo\left(\pmtwo_1\left(\unkdue_{\topindexuno}\right)\right)\right)\left(\unkuno_{\topindexuno}\right)=0\hspace{15pt}\forall \topindexuno\in \mathbb{N}\text{,}\\[4pt]
(iv)&\left(\smder_{\vecuno_1}\left(\evalcomptwo\left(\pmtwo_1\left(\unkdue_{\topindexuno}\right)\right)\right)\right)\left(\unkuno_{\topindexuno}\right)\neq\left(\smder_{\vecuno_2}\left(\evalcomptwo\left(\pmtwo_2\left(\unkdue_{\topindexuno}\right)\right)\right)\right)\left(\unkuno_{\topindexuno}\right)\hspace{15pt}\forall \topindexuno\in \mathbb{N}\text{;}
\end{array}\label{tt0}
\right.\\[6pt]
&\left\{
\begin{array}{ll}
(i)&\underset{\topindexuno\rightarrow +\infty}{\lim}\,\left(\unkdue_{\topindexuno},\unkuno_{\topindexuno}\right)=\left(0,0\right)\text{,}\\[4pt]
(ii)&\left(\evalcomptwo\left(\pmtwo_{\topindextre}\left(\unkdue_{\topindexuno}\right)\right)\right)\left(\unkuno_{\topindexuno}\right) \neq \udenunk \hspace{62pt}\forall \topindextre\in\left\{1,2\right\}\text{,}\hspace{10pt}\forall \topindexuno\in\mathbb{N}\text{,}\\[4pt]
(iii)&\valass|\left(\evalcomptwo\left(\pmtwo_1\left(\unkdue_{\topindexuno+1}\right)\right)\right)\left(\unkuno_{\topindexuno+1}\right)|\,<\,\valass|\left(\evalcomptwo\left(\pmtwo_1\left(\unkdue_{\topindexuno}\right)\right)\right)\left(\unkuno_{\topindexuno}\right)|\hspace{35pt}\forall \topindexuno\in \mathbb{N}\text{,}\\[4pt]
(iv)&\left(\smder_{\vecuno_1}\left(\evalcomptwo\left(\pmtwo_1\left(\unkdue_{\topindexuno}\right)\right)\right)\right)\left(\unkuno_{\topindexuno}\right)\neq\left(\smder_{\vecuno_2}\left(\evalcomptwo\left(\pmtwo_2\left(\unkdue_{\topindexuno}\right)\right)\right)\right)\left(\unkuno_{\topindexuno}\right)\hspace{15pt}\forall \topindexuno\in \mathbb{N}\text{.}
\end{array}\label{monotonadec}
\right.
\end{flalign}
If \eqref{tt0} occurs then it is possible to choose a function $f\in \Cksp{\infty}(\mathbb{R}^m)+\mathbb{R}+$ fulfilling
\begin{multline}
\left(\smder_{\vecuno}\left(f\right)\right)\left(0\right)\neq \left(\smder_{\vecdue}\left(f\right)\right)\left(0\right) \hspace{15pt}\forall\vecuno,\vecdue\in \left\{\left(\smder_{\vecuno_1}\left(\evalcomptwo\left(\pmtwo_1\left(\unkdue_{\topindexuno}\right)\right)\right)\right)\left(\unkuno_{\topindexuno}\right)\right\}_{\topindexuno\in \mathbb{N}}\cup\\[4pt]
\left\{\left(\smder_{\vecuno_2}\left(\evalcomptwo\left(\pmtwo_2\left(\unkdue_{\topindexuno}\right)\right)\right)\right)\left(\unkuno_{\topindexuno}\right)\right\}_{\topindexuno\in \mathbb{N}}\hspace{4pt}\text{with}\hspace{4pt}\vecuno\neq\vecdue\text{.}\label{diversi1}
\end{multline} 
Fix representatives $\gfy_1$ of $\genpreder\left[\gggfz[\gfz_1], \vecuno_1 \right]\left(\lininclquot\left(\gczsfx[f] \right)\right)$, $\gfy_2$ of $\genpreder\left[\gggfz[\gfz_2], \vecuno_2 \right]\left(\lininclquot\left(\gczsfx[f] \right)\right)$.\newline
Since $f$ is smooth it is possible to choose a path $\pmtwo_{\gfy_1}$ in $\genf$ detecting $\gfy_1$ with the same associating function of $\pmtwo_1$, a path $\pmtwo_{\gfy_2}$ in $\genf$ detecting $\gfy_2$ with the same associating function of $\pmtwo_2$, then assumption \eqref{tt0}-[(i)] entails that
\begin{equation}
\begin{array}{l}
\hspace{-10pt}\exists\overline{\topindexuno}\in \mathbb{N}\;:\;
\left\{\hspace{-5pt}
\begin{array}{ll}
\left(\evalcomptwo\left(\pmtwo_{\gfy_1}\left(\unkdue_{\topindexuno}\right)\right)\right)\left(\unkuno_{\topindexuno}\right)=\left(\smder_{\left(\smder_{\vecuno_1}\left(\evalcomptwo\left(\pmtwo_1\left(\unkdue_{\topindexuno}\right)\right)\right)\right)\left(\unkuno_{\topindexuno}\right)}\left(f\right)\right)\left(0\right)&\forall \topindexuno>\overline{\topindexuno}\text{,}\\[4pt]
\left(\evalcomptwo\left(\pmtwo_{\gfy_2}\left(\unkdue_{\topindexuno}\right)\right)\right)\left(\unkuno_{\topindexuno}\right)=\left(\smder_{\left(\smder_{\vecuno_2}\left(\evalcomptwo\left(\pmtwo_2\left(\unkdue_{\topindexuno}\right)\right)\right)\right)\left(\unkuno_{\topindexuno}\right)}\left(f\right)\right)\left(0\right)&\forall \topindexuno>\overline{\topindexuno}\text{.}
\end{array}
\right.
\end{array}\label{stscs}
\end{equation}
Eventually the claimed contradiction with \eqref{ipoteqprediff} follows by Definition \ref{releqinfnearpointdef}, \eqref{predereq}, \eqref{ipoteqprediffdue}, \eqref{diversi1}, \eqref{stscs}.\newline
If \eqref{monotonadec} occurs then it is possible to choose a function $f\in \Cksp{\infty}(\mathbb{R}^m)+\mathbb{R}+$ fulfilling
\begin{multline}
\left(\smder_{\left(\smder_{\vecuno_1}\left(\evalcomptwo\left(\pmtwo_1\left(\unkdue_{\topindexuno}\right)\right)\right)\right)\left(\unkuno_{\topindexuno}\right)}\left(f\right)\right)\left(\left(\evalcomptwo\left(\pmtwo_1\left(\unkdue_{\topindexuno}\right)\right)\right)\left(\unkuno_{\topindexuno}\right)\right)\neq\\[4pt] \left(\smder_{\left(\smder_{\vecuno_2}\left(\evalcomptwo\left(\pmtwo_2\left(\unkdue_{\topindexuno}\right)\right)\right)\right)\left(\unkuno_{\topindexuno}\right)}\left(f\right)\right)\left(\left(\evalcomptwo\left(\pmtwo_2\left(\unkdue_{\topindexuno}\right)\right)\right)\left(\unkuno_{\topindexuno}\right)\right) \hspace{15pt}\forall\topindexuno\in \mathbb{N}\text{.}\label{diversi2}
\end{multline}
Fix representatives $\gfy_1$ of $\genpreder\left[\gggfz[\gfz_1], \vecuno_1 \right]\left(\lininclquot\left(\gczsfx[f] \right)\right)$, $\gfy_2$ of $\genpreder\left[\gggfz[\gfz_2], \vecuno_2 \right]\left(\lininclquot\left(\gczsfx[f] \right)\right)$.\newline
Since $f$ is smooth it is possible to choose a path $\pmtwo_{\gfy_1}$ in $\genf$ detecting $\gfy_1$ with the same associating function of $\pmtwo_1$, a path $\pmtwo_{\gfy_2}$ in $\genf$ detecting $\gfy_2$ with the same associating function of $\pmtwo_2$, then assumption \eqref{tt0}-[(i)] entails that
\begin{equation}
\begin{array}{l}
\hspace{-10pt}\exists\overline{\topindexuno}\in \mathbb{N}\;:\;
\left\{\hspace{-5pt}
\begin{array}{ll}
\left(\evalcomptwo\left(\pmtwo_{\gfy_1}\left(\unkdue_{\topindexuno}\right)\right)\right)\left(\unkuno_{\topindexuno}\right)=\\[4pt]
\hspace{25pt}\left(\smder_{\left(\smder_{\vecuno_1}\left(\evalcomptwo\left(\pmtwo_1\left(\unkdue_{\topindexuno}\right)\right)\right)\right)\left(\unkuno_{\topindexuno}\right)}\left(f\right)\right)\left(\left(\evalcomptwo\left(\pmtwo_1\left(\unkdue_{\topindexuno}\right)\right)\right)\left(\unkuno_{\topindexuno}\right)\right)&\forall \topindexuno>\overline{\topindexuno}\text{,}\\[18pt]
\left(\evalcomptwo\left(\pmtwo_{\gfy_2}\left(\unkdue_{\topindexuno}\right)\right)\right)\left(\unkuno_{\topindexuno}\right)=\\[4pt]
\hspace{25pt}\left(\smder_{\left(\smder_{\vecuno_2}\left(\evalcomptwo\left(\pmtwo_2\left(\unkdue_{\topindexuno}\right)\right)\right)\right)\left(\unkuno_{\topindexuno}\right)}\left(f\right)\right)\left(\left(\evalcomptwo\left(\pmtwo_2\left(\unkdue_{\topindexuno}\right)\right)\right)\left(\unkuno_{\topindexuno}\right)\right)&\forall \topindexuno>\overline{\topindexuno}\text{.}
\end{array}
\right.
\end{array}\label{stscs1}
\end{equation}
Eventually the claimed contradiction with \eqref{ipoteqprediff} follows by Definition \ref{releqinfnearpointdef}, \eqref{predereq}, \eqref{ipoteqprediffdue}, 
\eqref{diversi2}, \eqref{stscs1}.\newline
\textnormal{\textbf{Proof of \eqref{tesittdiv}.}}\ \ Arguing by contradiction fix a sequence $\left\{\left(\unkdue_{\topindexuno},\unkuno_{\topindexuno}\right)\right\}_{\topindexuno \in \mathbb{N}}\subseteq \left(-1,1\right)\times\intuno\setminus\left\{0\right\}\times\intuno$ fulfilling one and only one among the two systems of conditions below: 
\begin{flalign}
&\left\{
\begin{array}{ll}
(i)&\valass|\unkdue_{\topindexuno+1}|<\valass|\unkdue_{\topindexuno}|\text{,}\hspace{8pt}\valass|\unkuno_{\topindexuno+1}|<\valass|\unkuno_{\topindexuno}|\hspace{15pt}\forall \topindexuno\in\mathbb{N}\text{,}\\[6pt]
(ii)&\underset{\topindexuno\rightarrow +\infty}{\lim}\,\left(\unkdue_{\topindexuno},\unkuno_{\topindexuno}\right)=\left(0,0\right)\text{,}\\[6pt]
(iii)&\left(\evalcomptwo\left(\pmtwo_{\topindextre}\left(\unkdue_{\topindexuno}\right)\right)\right)\left(\unkuno_{\topindexuno}\right) \neq \udenunk \hspace{15pt}\forall \topindextre\in\left\{1,2\right\}\text{,}\hspace{10pt}\forall \topindexuno\in\mathbb{N}\text{,}\\[6pt]
(iv)&\left(\evalcomptwo\left(\pmtwo_{\topindextre}\left(\unkdue_{\topindexuno}\right)\right)\right)\left(\unkuno_{\topindexuno}\right) \neq 0 \hspace{15pt}\forall \topindextre\in\left\{1,2\right\}\text{,}\hspace{10pt}\forall \topindexuno\in\mathbb{N}\text{,}\\[6pt]
(v)& \left(\evalcomptwo\left(\pmtwo_1\left(\unkdue_{\topindexuno}\right)\right)\right)\left(\unkuno_{\topindexuno}\right) \neq \left(\evalcomptwo\left(\pmtwo_2\left(\unkdue_{\topindexuno}\right)\right)\right)\left(\unkuno_{\topindexuno}\right) \hspace{15pt}\forall \topindexuno\in\mathbb{N}\text{,}\\[6pt]
(vi)&\max\left\{\valass|\left(\evalcomptwo\left(\pmtwo_1\left(\unkdue_{\topindexuno+1}\right)\right)\right)\left(\unkuno_{\topindexuno+1}\right)|, \valass|\left(\evalcomptwo\left(\pmtwo_2\left(\unkdue_{\topindexuno+1}\right)\right)\right)\left(\unkuno_{\topindexuno+1}\right)|\right\}<\\[4pt]
&\hspace{52pt}\min\left\{\valass|\left(\evalcomptwo\left(\pmtwo_1\left(\unkdue_{\topindexuno}\right)\right)\right)\left(\unkuno_{\topindexuno}\right)|, \valass|\left(\evalcomptwo\left(\pmtwo_2\left(\unkdue_{\topindexuno}\right)\right)\right)\left(\unkuno_{\topindexuno}\right)|\right\}\hspace{15pt}\forall i\in\mathbb{N}\text{,}\\[6pt]
(vi)&\left(\smder_{\vecuno_1}\left(\evalcomptwo\left(\pmtwo_1\left(\unkdue_{i}\right)\right)\right)\right)\left(\unkuno_{\topindexuno}\right)\neq 0\hspace{15pt}\forall \topindexuno \in \mathbb{N}\text{;}
\end{array}
\right.\label{Puno}\\[8pt]
&\left\{
\begin{array}{ll}
(i)&\valass|\unkdue_{\topindexuno+1}|<\valass|\unkdue_{\topindexuno}|\text{,}\hspace{8pt}\valass|\unkuno_{\topindexuno+1}|<\valass|\unkuno_{\topindexuno}|\hspace{15pt}\forall \topindexuno\in\mathbb{N}\text{,}\\[6pt]
(ii)&\underset{\topindexuno\rightarrow +\infty}{\lim}\,\left(\unkdue_{\topindexuno},\unkuno_{\topindexuno}\right)=\left(0,0\right)\text{,}\\[6pt]
(iii)&\left(\evalcomptwo\left(\pmtwo_{\topindextre}\left(\unkdue_{\topindexuno}\right)\right)\right)\left(\unkuno_{\topindexuno}\right) \neq \udenunk \hspace{15pt}\forall \topindextre\in\left\{1,2\right\}\text{,}\hspace{10pt}\forall \topindexuno\in\mathbb{N}\text{,}\\[6pt]
(iv)&0<\valass|\left(\evalcomptwo\left(\pmtwo_1\left(\unkdue_{\topindexuno+1}\right)\right)\right)\left(\unkuno_{\topindexuno+1}\right)|<\valass|\left(\evalcomptwo\left(\pmtwo_1\left(\unkdue_{\topindexuno}\right)\right)\right)\left(\unkuno_{\topindexuno}\right)|\hspace{15pt}\forall \topindexuno\in\mathbb{N}\text{,}\\[6pt]
(v)&\text{for}\hspace{4pt}\text{any}\hspace{4pt}\topindexuno \in \mathbb{N}\hspace{4pt}\text{there}\hspace{4pt}\text{is}\hspace{4pt}\text{an}\hspace{4pt}\text{open}\hspace{4pt}\text{neighborhood}\hspace{4pt}\setsymcinque_{2,\topindexuno}\subseteq\intuno\hspace{4pt}\text{of}\\[4pt]
&\unkuno_{\topindexuno}\hspace{4pt}\text{with}\hspace{4pt}\pmtwo_2\left(\unkdue_{\topindexuno}\right)\genfuncomp\incl{\setsymcinque_{2,\topindexuno}}{\intuno}=\cost<\setsymcinque_{2,\topindexuno}<>\mathbb{R}^m>+0+\text{,}\\[6pt]
(vi)&\left(\smder_{\vecuno_1}\left(\evalcomptwo\left(\pmtwo_1\left(\unkdue_{\topindexuno}\right)\right)\right)\right)\left(\unkuno_{\topindexuno}\right)\neq 0 \hspace{15pt}\forall \topindexuno \in \mathbb{N}\text{.}
\end{array}
\right.\label{Pdue} 
\end{flalign}
First we prove that it is possible to fix a sequence fulfilling either \eqref{Puno} or \eqref{Pdue}. Up to choosing a suitable subsequence of $\left\{\left(\unkdue_{\topindexuno},\unkuno_{\topindexuno}\right)\right\}_{\topindexuno \in \mathbb{N}}$ and relabeling paths $\pmtwo_1$, $\pmtwo_2$, either \eqref{Puno} is fulfilled or 
\begin{equation}
\left\{
\begin{array}{ll}
(i)&\valass|\unkdue_{\topindexuno+1}|<\valass|\unkdue_{\topindexuno}|\text{,}\hspace{8pt}\valass|\unkuno_{\topindexuno+1}|<\valass|\unkuno_{\topindexuno}|\hspace{15pt}\forall \topindexuno\in\mathbb{N}\text{,}\\[6pt]
(ii)&\underset{\topindexuno\rightarrow +\infty}{\lim}\,\left(\unkdue_{\topindexuno},\unkuno_{\topindexuno}\right)=\left(0,0\right)\text{,}\\[6pt]
(iii)&\left(\evalcomptwo\left(\pmtwo_{\topindextre}\left(\unkdue_{\topindexuno}\right)\right)\right)\left(\unkuno_{\topindexuno}\right) \neq \udenunk \hspace{15pt}\forall \topindextre\in\left\{1,2\right\}\text{,}\hspace{10pt}\forall \topindexuno\in\mathbb{N}\text{,}\\[6pt]
(iv)&0<\valass|\left(\evalcomptwo\left(\pmtwo_1\left(\unkdue_{\topindexuno+1}\right)\right)\right)\left(\unkuno_{\topindexuno+1}\right)|<\valass|\left(\evalcomptwo\left(\pmtwo_1\left(\unkdue_{\topindexuno}\right)\right)\right)\left(\unkuno_{\topindexuno}\right)|\hspace{15pt}\forall \topindexuno\in\mathbb{N}\text{,}\\[6pt]
(v)&\left(\evalcomptwo\left(\pmtwo_2\left(\unkdue_{\topindexuno}\right)\right)\right)\left(\unkuno_{\topindexuno}\right)=0\hspace{15pt}\forall \topindexuno\in\mathbb{N}\text{,}\\[6pt]
(vi)&\left(\smder_{\vecuno_1}\left(\evalcomptwo\left(\pmtwo_1\left(\unkdue_{\topindexuno}\right)\right)\right)\right)\left(\unkuno_{\topindexuno}\right)\neq 0 \hspace{15pt}\forall \topindexuno \in \mathbb{N}\text{.}
\end{array}
\right.\label{zeroconddue}
\end{equation}
If \eqref{zeroconddue} occurs then, up to choosing a suitable subsequence of $\left\{\left(\unkdue_{\topindexuno},\unkuno_{\topindexuno}\right)\right\}_{\topindexuno\in \mathbb{N}}$, one and only one among the following two mutually exclusives cases occurs:
\begin{flalign}
&\left\{
\begin{array}{l}
\text{for}\hspace{4pt}\text{any}\hspace{4pt}\topindexuno \in \mathbb{N}\hspace{4pt}\text{there}\hspace{4pt}\text{is}\hspace{4pt}\unkuno_{\topindexuno,1}\in \intuno
\hspace{4pt}\text{fulfilling}\hspace{4pt}\text{all}\hspace{4pt}\text{conditions}\hspace{4pt}\text{below:}\\[4pt]
\begin{array}{ll}
(i)&\valass|\unkuno_{\topindexuno}-\unkuno_{\topindexuno,1}|<
\left\{
\begin{array}{ll}
\frac{1}{2}\valass|\unkuno_{2}-\unkuno_{1}| &\text{if}\hspace{4pt}\topindexuno=1\text{,}\\[4pt]
\frac{1}{2}\min\left\{\valass|\unkuno_{\topindexuno+1}-\unkuno_{\topindexuno}|,\valass|\unkuno_{\topindexuno}-\unkuno_{\topindexuno-1}|\right\}&\text{if}\hspace{4pt}\topindexuno>1\text{,}
\end{array}
\right.\\[4pt]
(ii)&\left(\evalcomptwo\left(\pmtwo_{\topindextre}\left(\unkdue_{\topindexuno}\right)\right)\right)\left(\unkuno_{\topindexuno,1}\right) \neq \udenunk \hspace{15pt}\forall \topindextre\in\left\{1,2\right\}\text{,}\hspace{10pt}\forall \topindexuno\in\mathbb{N}\text{,}\\[4pt]
(iii)&\left(\evalcomptwo\left(\pmtwo_{\topindextre}\left(\unkdue_{\topindexuno}\right)\right)\right)\left(\unkuno_{\topindexuno,1}\right)\neq 0\hspace{15pt}\forall \topindextre \in \left\{1,2\right\}\text{,}\hspace{10pt}\forall \topindexuno\in\mathbb{N}\text{,}\\[4pt]
(iv)&\left(\evalcomptwo\left(\pmtwo_1\left(\unkdue_{\topindexuno}\right)\right)\right)\left(\unkuno_{\topindexuno,1}\right)\neq \left(\evalcomptwo\left(\pmtwo_2\left(\unkdue_{\topindexuno}\right)\right)\right)\left(\unkuno_{\topindexuno,1}\right)\text{,}\\[4pt]
(v)&\left(\smder_{\vecuno_1}\left(\evalcomptwo\left(\pmtwo_1\left(\unkdue_{\topindexuno}\right)\right)\right)\right)\left(\unkuno_{\topindexuno,1}\right)\neq 0\text{;}
\end{array}
\end{array}\label{zerointdue}
\right.\\[6pt]
&\left\{
\begin{array}{l}
\text{for}\hspace{4pt}\text{any}\hspace{4pt}\topindexuno \in \mathbb{N}\hspace{4pt}\text{there}\hspace{4pt}\text{is}\hspace{4pt}\text{an}\hspace{4pt}\text{open}\hspace{4pt}\text{neighborhood}\hspace{4pt}\setsymcinque_{2,\topindexuno}\subseteq\intuno\hspace{4pt}\text{of}\hspace{4pt}\unkuno_{\topindexuno}\\[4pt]
\text{with}\hspace{4pt}\pmtwo_2\left(\unkdue_{\topindexuno}\right)\genfuncomp\incl{\setsymcinque_{2,\topindexuno}}{\intuno}=\cost<\setsymcinque_{2,\topindexuno}<>\mathbb{R}^m>+0+\text{.}
\end{array}
\right.\label{zerointuno}
\end{flalign}
If \eqref{zerointdue} occurs then by taking a suitable subsequence of $\left\{\left(\unkdue_{\topindexuno},\unkuno_{\topindexuno,1}\right)\right\}_{\topindexuno \in \mathbb{N}}$ we obtain a sequence fulfilling \eqref{Puno}.\newline 
If \eqref{zerointuno} occurs then we have that sequence $\left\{\left(\unkdue_{\topindexuno},\unkuno_{\topindexuno}\right)\right\}_{\topindexuno \in \mathbb{N}}$ fulfills \eqref{Pdue}.\newline
If \eqref{Puno} occurs then we argue as follows. By exploiting Euclidean topology of $\mathbb{R}^m$ and referring to Notation \ref{gentopnot}-[7], there are an open neighborhood $\setsymnove_{\topindextre,\topindexuno}\subseteq \mathbb{R}^m$ of $\left(\evalcomptwo\left(\pmtwo_{\topindextre}\left(\unkdue_{\topindexuno}\right)\right)\right)\left(\unkuno_{\topindexuno}\right)$ for any $\topindextre \in \left\{1,2\right\}$, $\topindexuno \in \mathbb{N}$, $f\in C^{\infty}\left(\mathbb{R}^m,\mathbb{R}\right)$ fulfilling all conditions below: 
\begin{flalign*}
&\begin{array}{l}
\closure\left[\setsymnove_{\topindextre_1,\topindexuno_1},\mathbb{R}^m \right]\cap\closure\left[\setsymnove_{\topindextre_2,\topindexuno_2},\mathbb{R}^m \right]=\udenset\hspace{15pt}
\forall \topindextre_1,\topindextre_2 \in \left\{1,2\right\}\text{,}\hspace{10pt} \forall \topindexuno_1,\topindexuno_2 \in \mathbb{N}\text{;}
\end{array}\\[6pt]
&\begin{array}{l}
f\left(\unkuno\right)=0\hspace{15pt} \forall \unkuno \in \setsymnove_{2,\topindexuno}\text{,}\hspace{10pt} \forall \topindexuno \in \mathbb{N}\text{;}
\end{array}\\[6pt]
&\begin{array}{l}
\left(\smder_{\vecuno_1}\left(f\funcomp\left(\evalcomptwo\left(\pmtwo_1\left(\unkdue_{\topindexuno}\right)\right)\right)\right)\right)\left(\unkuno_{\topindexuno}\right)\neq 0 \hspace{15pt}\forall \topindexuno \in \mathbb{N}\text{.}
\end{array}
\end{flalign*}
Then \eqref{ipoteqprediff} fails since $\genpreder\left[\gggfz[\gfz_1], \vecuno_1 \right]\left(\lininclquot\left(\gczsfx[f]\right)\right)\neq\genpreder\left[\gggfz[\gfz_2], \vecuno_2 \right]\left(\lininclquot\left(\gczsfx[f]\right)\right)$.\newline  
If \eqref{Pdue} occurs then we argue as follows. By exploiting Euclidean topology of $\mathbb{R}^m$ and referring to Notation \ref{gentopnot}-[7], there are an open neighborhood $\setsymnove_{\topindexuno}\subseteq \mathbb{R}^m$ of $\left(\evalcomptwo\left(\pmtwo_1\left(\unkdue_{\topindexuno}\right)\right)\right)\left(\unkuno_{\topindexuno}\right)$ for any $\topindexuno \in \mathbb{N}$, $f\in C^{\infty}\left(\mathbb{R}^m,\mathbb{R}\right)$ fulfilling all conditions below: 
\begin{flalign*}
&\begin{array}{l}
\closure\left[\setsymnove_{\topindexuno_1},\mathbb{R}^m \right]\cap\closure\left[\setsymnove_{\topindexuno_2},\mathbb{R}^m \right]=\udenset\hspace{15pt}
 \forall \topindexuno_1,\topindexuno_2 \in \mathbb{N}\text{;}
\end{array}\\[6pt]
&\begin{array}{l}
f\left(0\right)=0\text{;}
\end{array}\\[6pt]
&\begin{array}{l}
\left(\smder_{\vecuno_1}\left(f\funcomp\left(\pmtwo_1\left(\unkdue_{\topindexuno}\right)\right)\right)\right)\left(\unkuno_{\topindexuno}\right)\neq 0 \hspace{15pt}\forall \topindexuno \in \mathbb{N}\text{.}
\end{array}
\end{flalign*}
Eventually \eqref{ipoteqprediff} fails since $\genpreder\left[\gggfz[\gfz_1], \vecuno_1 \right]\left(\lininclquot\left(\gczsfx[f]\right)\right)\neq\genpreder\left[\gggfz[\gfz_2], \vecuno_2 \right]\left(\lininclquot\left(\gczsfx[f]\right)\right)$.
\end{proof}

In Corollary \ref{zerosucomp} below we study where pre-derivations of germs of local generalized functions vanish. We Refer to Example \ref{ExCzero}.

\begin{corollary}\label{zerosucomp}
Fix $l,m\in\mathbb{N}_0$, vectors $\vecuno_1,\vecuno_2\in \mathbb{R}^{l}$, an open neighborhood $\intuno\subseteq\mathbb{R}^l$ of $0 \in \mathbb{R}^l$, $\setsymsette \compacont \intuno$, $\gfz_1,\gfz_2\in \genflocccz(\intuno)+\mathbb{R}^m+$,  paths $\pmtwo_1$, $\pmtwo_2$ in $\genf$ detecting $\gfz_1$, $\gfz_2$ respectively. Assume that data fulfill both \eqref{ipoteqprediff} and \eqref{ipoteqprediffdue}.\newline
Define: 
\begin{flalign}
&\begin{array}{l}
\intuno_{\neq}=\left\{\unkuno\in \intuno\;:\; \left(\evalcomptwo\left(\gfz_1\right)\right)\left(\unkuno\right)\neq\left(\evalcomptwo\left(\gfz_2\right)\right)\left(\unkuno\right)\right\}\text{;}
\end{array}\\[8pt]
&\left\{
\begin{array}{l}
\pmtwo:\left(-1,1\right)\setminus\left\{0\right\}\rightarrow  \Cksp{0}\hspace{4pt}\text{by}\hspace{4pt}\text{setting}\\[4pt]
\pmtwo\left(\unkdue\right)= \contsplbound\smder_{\vecuno_1}\left(\evalcomptwo\left(\pmtwo_1\left(\unkdue\right)\right)\right) , \smder_{\vecuno_2}\left(\evalcomptwo\left(\pmtwo_2\left(\unkdue\right)\right)\right)\contsprbound \funcomp \Diag{\intuno}{2} \hspace{15pt}\forall \unkdue \in\left(-1,1\right)\setminus\left\{0\right\}\text{.}
\end{array}
\right.
\end{flalign}
Then there are sequence $\left\{\delta_{\topindexuno}\right\}_{\topindexuno\in \mathbb{N}}$, an open neighborhood $\setsymnove \sqsubseteq  \intuno$ of $0 \in \intuno$ fulfilling all conditions below: 
\begin{flalign}
&\begin{array}{l}
0<\delta_{\topindexuno+1}<\delta_{\topindexuno} \hspace{15pt}\forall\topindexuno\in \mathbb{N}\text{;}\label{decdelta}
\end{array}\\[6pt]
&\begin{array}{l}
\underset{\topindexuno\rightarrow \infty}{\lim}\,\delta_{\topindexuno}=0\text{;}\label{decdeltaazero}
\end{array}\\[6pt]
&\begin{array}{l}
\left(\pmtwo\left(\unkdue\right)\right)\left(\unkuno\right)\neq \udenunk \hspace{15pt}\left(\unkdue, \unkuno\right)\in\left( \left(-\delta_{1},\delta_1\right)\times\setsymnove\right)\setminus \left(\left\{0\right\}\times\setsymnove\right)\text{;}\label{noevalvuoto}
\end{array}\\[6pt]
&\begin{array}{l}
\setsymnove\cap\intuno_{\neq}\cap\setsymsette\cap\left(\underset{\unkdue\in\left[-\delta_{\topindexuno},\delta_{\topindexuno}\right]\setminus\left\{	0\right\}}{\bigcap} \left(\pmtwo\left(\unkdue\right)\right)^{-1}\left(0\right)\right)\subseteq\\[4pt]
\hspace{100pt}  \setsymnove\cap\intuno_{\neq}\cap\setsymsette\cap\left(\underset{\unkdue\in\left[-\delta_{\topindexuno+1},\delta_{\topindexuno+1}\right]\setminus\left\{	0\right\}}{\bigcap} \left(\pmtwo\left(\unkdue\right)\right)^{-1}\left(0\right)\right)\hspace{15pt}\forall\topindexuno\in \mathbb{N}\text{;}\label{nonzerodentro}
\end{array}\\[6pt]
&\begin{array}{l}
\setsymnove\cap\setsymsette\cap\left(\underset{\topindexuno\in \mathbb{N}}{\bigcup} \left(
\underset{\unkdue\in\left[-\delta_{\topindexuno},\delta_{\topindexuno}\right]\setminus\left\{	0\right\}}{\bigcap} \left(\pmtwo\left(\unkdue\right)\right)^{-1}\left(0\right)
\right)\right)= 
\setsymnove\cap\intuno_{\neq}\cap\setsymsette\text{.}\label{nonzerolim}
\end{array}
\end{flalign}
\end{corollary}
\begin{proof}\mbox{}\newline
If $l=0$ then statement follows straightforwardly by \eqref{predereq}-[(ii)].\newline
If $l\neq 0$ then we argue as follows.\newline
Statement will be proved by inductively defining the sequence $\left\{\delta_{\topindexuno}\right\}_{\topindexuno\in \mathbb{N}}$ fulfilling \eqref{decdelta}-\eqref{nonzerolim}.\newline
Fix $\delta$ and $\setsymnove$ obtained by applying Proposition \ref{struttint}.\newline
Since $\setsymsette \compacont \intuno$ we are able to fix a sequence of sets $\left\{\setsymsette_{\elindexuno}\right\}_{\elindexuno\in \mathbb{N}}$ fulfilling both conditions below 
\begin{flalign}
&\begin{array}{l}
\setsymsette_{\elindexuno} \compacont \setsymsette_{\elindexuno+1} \compacont \setsymnove \cap\intuno_{\neq} \cap \setsymsette\hspace{15pt}\forall \elindexuno\in \mathbb{N}\text{;}\label{seqcompacont}
\end{array}\\[4pt]
&\begin{array}{l}
\underset{\elindexuno\in \mathbb{N}}{\bigcup}\setsymsette_{\elindexuno}=\setsymnove \cap\intuno_{\neq}\cap\setsymsette\text{.}\label{esaurimento}
\end{array}
\end{flalign}   
Since $\gfz_1,\gfz_2\in \genflocccz(\intuno)+\mathbb{R}^m+$ there is $\delta_1\in\left(0,\delta\right)$ fulfilling
\begin{multline}
\left(\evalcomptwo\left(\pmtwo_1\left(\unkdue\right)\right)\right)\left(\unkuno\right)\neq\left(\evalcomptwo\left(\pmtwo_2\left(\unkdue\right)\right)\right)\left(\unkuno\right)\\[4pt]
\hspace{110pt}
\forall\left(\unkdue,\unkuno\right) \in \left(\left(-\delta_{1},\delta_{1}\right)\times\setsymsette_{1}\right)\setminus\left(\left\{0\right\}\times\setsymsette_{1}\right)\text{,}\label{zerodali}
\end{multline}
then Proposition \ref{struttint}, \eqref{zerodali} entail that 
\begin{equation*}
\left(\evalcomptwo\left(\pmtwo\left(\unkdue\right)\right)\right)\left(\unkuno\right)=0\hspace{15pt}\forall\left(\unkdue,\unkuno\right) \in \left(\left(-\delta_{1},\delta_{1}\right)\times\setsymsette_{1}\right)\setminus\left(\left\{0\right\}\times\setsymsette_{1}\right)\text{.}
\end{equation*}
By induction assume that there are $\topindexuno\in\mathbb{N}$, a finite sequence $\left\{\delta_{\elindexuno}\right\}_{\elindexuno=1}^{\topindexuno}$ of positive real numbers fulfilling all conditions below:
\begin{flalign}
&\begin{array}{l}
\delta_{\elindexuno+1} <\frac{\delta_{\elindexuno}}{2}\hspace{15pt}\forall \elindexuno\in\left\{1,...,\topindexuno-1\right\}\text{;}\label{dminm}
\end{array}\\[6pt]
&\left\{
\begin{array}{l}
\left(\evalcomptwo\left(\pmtwo_1\left(\unkdue\right)\right)\right)\left(\unkuno\right)\neq\left(\evalcomptwo\left(\pmtwo_2\left(\unkdue\right)\right)\right)\left(\unkuno\right)\\[4pt]
\text{for}\hspace{4pt}\text{any}\hspace{4pt}\elindexuno\in\left\{1,...,\topindexuno\right\}\text{,}\hspace{4pt}\left(\unkdue,\unkuno\right) \in \left(\left(-\delta_{\elindexuno},\delta_{\elindexuno}\right)\times\setsymsette_{\elindexuno}\right)\setminus\left(\left\{0\right\}\times\setsymsette_{\elindexuno}\right)\text{;}\label{zerodaliw1w2}
\end{array}
\right.\\[6pt]
&\left\{
\begin{array}{l}
\left(\evalcomptwo\left(\pmtwo\left(\unkdue\right)\right)\right)\left(\unkuno\right)=0\\[4pt]
\text{for}\hspace{4pt}\text{any}\hspace{4pt}\elindexuno\in\left\{1,...,\topindexuno\right\}\text{,}\hspace{4pt}\left(\unkdue,\unkuno\right) \in \left(\left(-\delta_{\elindexuno},\delta_{\elindexuno}\right)\times\setsymsette_{\elindexuno}\right)\setminus\left(\left\{0\right\}\times\setsymsette_{\elindexuno}\right)\text{.}\label{zeroinC}
\end{array}
\right.
\end{flalign}
Since $\gfz_1,\gfz_2\in \genflocccz(\intuno)+\mathbb{R}^m+$ there is $\delta_{\topindexuno+1}\in\left(0,1\right)$ fulfilling
\begin{flalign}
&\begin{array}{l}
\delta_{\topindexuno+1} < \frac{\delta_{\topindexuno}}{2}\text{;}
\end{array}\\[6pt]
&\begin{array}{l}
\left(\evalcomptwo\left(\pmtwo_1\left(\unkdue\right)\right)\right)\left(\unkuno\right)\neq\left(\evalcomptwo\left(\pmtwo_2\left(\unkdue\right)\right)\right)\left(\unkuno\right)\\[4pt]
\hspace{90pt}
\forall\left(\unkdue,\unkuno\right) \in \left(\left(-\delta_{\topindexuno+1},\delta_{\topindexuno+1}\right)\times\setsymsette_{\topindexuno+1}\right)\setminus\left(\left\{0\right\}\times\setsymsette_{\topindexuno+1}\right)\text{,}\label{zerodali1}
\end{array}
\end{flalign}
then Proposition \ref{struttint}, \eqref{zerodali1} entail that
\begin{equation*}
\left(\evalcomptwo\left(\pmtwo\left(\unkdue\right)\right)\right)\left(\unkuno\right)=0\hspace{15pt}\forall\left(\unkdue,\unkuno\right) \in \left(\left(-\delta_{\topindexuno+1},\delta_{\topindexuno+1}\right)\times\setsymnove_{\topindexuno+1}\right)\setminus\left(\left\{0\right\}\times\setsymnove_{\topindexuno+1}\right)\text{.}
\end{equation*}
Eventually statement follows, namely: \eqref{dminm} entails both \eqref{decdelta} and \eqref{decdeltaazero}; \eqref{seqcompacont} and \eqref{zeroinC} together entail \eqref{nonzerodentro}; \eqref{esaurimento} entails \eqref{nonzerolim}.
\end{proof}

Motivated by Proposition \ref{struttint}, Corollary \ref{zerosucomp} we give Definition \ref{zerinsggfdef} below.
We refer to Notation \ref{realvec}-[2].

\begin{definition}\label{zerinsggfdef}\mbox{}
\begin{enumerate}
\item Fix $l,m\in\mathbb{N}_0$, $\gfx \in \genflocccz(l)+m+$. We define sets:
\begin{flalign*}
&\vandergf{\gfx}=\left\{
\begin{array}{ll}
\left\{\vecuno \in \mathbb{R}^l\;:\; \underset{\mathsf{l}=1}{\overset{l}{\biggenfunsum}}\, \left( \Fpart_{\mathsf{l}} \genfbsfmu \gfx\right)\genfunscalp \vecuno_{\mathsf{l}}
 =0\scalpquotpoint\gfx
\right\}& \text{if}\hspace{4pt}l\neq 0\text{,}\\[6pt]
\left\{0\right\}& \text{if}\hspace{4pt}l=0\text{;}
\end{array}
\right.\\[8pt]
&\ppvandergf{\gfx}=\left\{\vecuno \in \mathbb{R}^l\;:\; \innpreul\vecuno,\vecdue\innpreur_{m} =0\hspace{10pt}\forall \vecdue\in \vandergf{\gfx}\right\}\text{.}
\end{flalign*}
\item Fix $l,m\in\mathbb{N}_0$, $\gggfx\in \genfquotcccz(l)+m+$. We define sets:
\begin{flalign*}
&\vandergf{\gggfx}=\left\{
\begin{array}{ll}
\left\{\vecuno \in \mathbb{R}^l\;:\; 
 \underset{\mathsf{l}=1}{\overset{l}{\bigsumquotpoint}}\, \left( \Fpart_{\mathsf{l}} \bsfmuquotpoint \gggfx\right)\scalpquotpoint  \vecuno_{\mathsf{l}}
 =0\scalpquotpoint\gggfx\right\}& \text{if}\hspace{4pt}l\neq 0\text{,}\\[6pt]
\left\{0\right\}& \text{if}\hspace{4pt}l=0\text{;}
\end{array}
\right.\\[8pt]
&\ppvandergf{\gggfx}=\left\{\vecuno \in \mathbb{R}^l\;:\; \innpreul \vecuno,\vecdue\innpreur_{m} =0\hspace{10pt}\forall \vecdue\in \vandergf{\gggfx}\right\}\text{.}
\end{flalign*}
\end{enumerate}
\end{definition}

In Proposition \ref{zerinsggfprop} below we study properties of sets introduced in Definition \ref{zerinsggfdef}.
We refer to Proposition \ref{genpointprop}.

\begin{proposition}\label{zerinsggfprop}\mbox{}
\begin{enumerate}
\item Fix $l,m\in\mathbb{N}_0$, $\gggfx \in \genfquotcccz(l)+m+$. Then $\vandergf{\gggfx}$ is a vector subspace of $\mathbb{R}^l$.
\item Fix $l,m\in\mathbb{N}_0$, $\gggfx \in \genfquotcccz(l)+m+$. Then there is $\gfx \in \genflocccz$ fulfilling both conditions below: 
\begin{equation*}
\gggfx=\gggfx[\gfx]\text{;}\hspace{50pt}
\vandergf{\gggfx}=\vandergf{\gfx}\text{.}
\end{equation*}
\item Fix $\gfx\in\genf$, $\gggfy\in\genfquot$. Then  $\vandergf{\gggfy}\subseteq\vandergf{\gfx \,\compgenfquot\gggfy}$.
\item Fix $\gggfy_1\in \genfquot$, $\gggfy_2 \in \genfquotcccz$. Then $\vandergf{\gggfy_2}\subseteq\vandergf{\gggfy_1\compquotpoint\gggfy_2}$.
\item Fix $\gggfy_1, \gggfy_2\in \genfquot$. Then $\vandergf{\gggfy_1} \dirsum\vandergf{\gggfy_2}=\vandergf{\lboundquotpoint\gggfy_1, \gggfy_2 \rboundquotpoint}$.
\item Fix $a\in \mathbb{R}$, $\gggfy\in \genfquot$. Then $\vandergf{\gggfy}=\vandergf{a\scalpquotpoint \gggfy}$.
\item Fix $l,m\in\mathbb{N}_0$, $\vecuno \in \mathbb{R}^l$, $\gggfz \in \genfquotcccz(l)+m+$, $\genpreder\in\genprederspace(l)+m+$. Then the orthogonal projection $\vecdue$
of $\vecuno$ on $\vandergf{\gggfz}$ is the unique vector belonging to $\vandergf{\gggfz}$ fulfilling $\quotdirsumuno\left(\genpreder\left[\gggfz, \vecuno \right]\right)=\quotdirsumuno\left(\genpreder\left[\gggfz, \vecdue \right]\right)$. 
\end{enumerate}
\end{proposition}
\begin{proof}\mbox{}\newline
\textnormal{\textbf{Proof of statement 1.}}\ \ Statement 1 straightforwardly follows by bilinearity of $\genfunscalp$.\newline
\textnormal{\textbf{Proof of statement 2.}}\ \ Fix a generating set $\bascont$ of $\vandergf{\gggfx}$. Since $\vandergf{\gggfz}\subseteq \mathbb{R}^l$ then there is no loss of generality by assuming that $\bascont$ is a finite set, hence there is $\gfx \in \genflocccz$ fulfilling both conditions below: 
\begin{flalign*}
&\begin{array}{l}
\gggfx=\gggfx[\gfx]\text{;}
\end{array}\\[8pt]
&\begin{array}{l}
\underset{\mathsf{l}=1}{\overset{l}{\biggenfunsum}}\, \left( \Fpart_{\mathsf{l}} \genfbsfmu \gfx\right)\genfunscalp \vecuno_{\mathsf{l}}
 =0\scalpquotpoint\gfx\hspace{15pt}\forall \vecuno\in \bascont\text{.}
\end{array}
\end{flalign*} 
Then bilinearity of $\genfunscalp$ entails that $\vandergf{\gggfx}\subseteq\vandergf{\gfx}$. Eventually $\vandergf{\gfx}\subseteq\vandergf{\gggfx}$ follows since any $\vecuno \in\vandergf{\gfx}$ belongs to $\vandergf{\gggfx}$ by definition of $\vandergf{\gggfx}$.\newline
\textnormal{\textbf{Proof of statements 3 - 6.}}\ \ Statements follows by properties of operators $\left\{\Fpart_{\topindexuno}\right\}_{\topindexuno \in \mathbb{N}}$\newline
\textnormal{\textbf{Proof of statement 7.}}\ \ Fix $\vecuno_1, \vecuno_2\in \mathbb{R}^l$ such that $\quotdirsumuno\left(\genpreder\left[\gggfz, \vecuno_1 \right]\right)=\quotdirsumuno\left(\genpreder\left[\gggfz, \vecuno_2 \right]\right)$.\newline
By exploiting statement 2 fix an open neighborhood $\intuno\subseteq\mathbb{R}^l$ of $0 \in \mathbb{R}^l$, $\gfz\in \genflocccz(\intuno)+\mathbb{R}^m+$ fulfilling both $\gggfz=\gggfz[\gfz]$ and $\vandergf{\gggfz}=\vandergf{\gfz}$.\newline
Eventually statement follows since Proposition \ref{struttint} applied in case $\gfz_1=\gfz_2=\gfz$ entails that $\vecuno_1-\vecuno_2\in \vandergf{\gggfz}$. 
\end{proof}

\section{Topology of pre-derivations \label{Grtopss}}

In this section we construct a site containing the site $\Topcat$ as a subsite and all diagrams involving arrows introduced in Sections \ref{ggf}, \ref{Derggf}. Later we will prove that any reasonable topology on pre-derivations cannot be given in terms of open sets, but a Grothendieck topology is needed (Remark \ref{siterem2}-[2]).\newline
We refer to Notation \ref{catfun} and to \cite{MM} Chap. III.\vspace{12pt}

In Definition \ref{augsitedef} we introduce the notion of augmentation of a site by adding selected objects and arrows.

\begin{definition}\label{augsitedef}
Fix a category $\Catuno$, a site $\left(\Catdue,\Grtop\right)$ together with a forgetful functor $\fundim\left[\Catdue,\Catuno\right]$ forgetting the Grothendieck topology $\Grtop$, a subcategory $\Cattre$ of $\Catuno$, a category $\Catquattro$ whose objects are categories $\Catsei$ each one equipped with a forgetful functor $\fundim\left[\Catsei,\Catuno\right]$, a category $\Catcinque$ whose objects are adjunctions $\left(\Catsei, \Catsette, \genfuncuno, \genfuncdue, \adjnattrasf \right)$ where $\Catsei$ is equipped with a forgetful functor $\fundim\left[\Catsei,\Catuno\right]$ and $\Catsette$ is equipped with a forgetful functor $\fundim\left[\Catsette,\Catuno\right]$,  a function $\precollarr$ assigning to any object $\objtre$ of $\Cattre$ a set $\precollarr\left(\objtre\right)$ of sets of arrows of $\Catuno$ whose domains are objects of $\Catdue$ and co-domains all coincide with $\objtre$.\newline 
A site $\left(\Catdue_1,\Grtop_1\right)$ is called augmentation of $\left(\Catdue,\Grtop\right)$ through $\left(\Catuno,\Cattre,\Catquattro,\Catcinque, \precollarr\right)$, or a $\left(\Catuno, \Cattre, \Catquattro, \Catcinque, \precollarr\right)$-augmentation of $\left(\Catdue,\Grtop\right)$, if and only if it fulfills all conditions below:
\begin{flalign}
&\begin{array}{l}
\text{there is a forgetful functor}\;\fundim\left[\Catdue_1,\Catuno\right] \text{;}
\end{array}\nonumber\\[6pt]
&\begin{array}{l}
\Catdue\;\text{and}\;\Cattre \;\text{are both subcategories of}\; \Catdue_1\text{;}
\end{array}\nonumber\\[6pt]
&\left\{\begin{array}{l}
\text{for}\hspace{4pt}\text{any}\hspace{4pt}\text{object}\hspace{4pt}\objdue\hspace{4pt}\text{of}\;\Catdue\text{,}\hspace{4pt}\text{sieve} \hspace{4pt}\sieve_1\hspace{4pt}\text{on}\hspace{4pt}\objdue\hspace{4pt}\text{in}\hspace{4pt}\Catdue_1\hspace{4pt}\text{such}\hspace{4pt}\text{that}\hspace{4pt}\sieve\subseteq\sieve_1\\
\text{for}\hspace{4pt}\text{some}\hspace{4pt}\text{sieve}\hspace{4pt}\sieve\in\Grtop\left(\objdue\right)\text{,}\hspace{4pt}\text{we}\hspace{4pt}\text{have} \hspace{4pt}\sieve_1\in\Grtop_1\left(\objdue\right)\text{;}
\end{array}
\right.\nonumber\\[6pt]
&\left\{\begin{array}{l}
\text{for}\hspace{4pt}\text{any}\hspace{4pt}\text{object}\hspace{4pt}\objtre\hspace{4pt}\text{of}\;\Cattre\text{,}\hspace{4pt}\text{sieve} \hspace{4pt}\sieve_1\hspace{4pt}\text{on}\hspace{4pt}\objtre\hspace{4pt}\text{in}\hspace{4pt}\Catdue_1\hspace{4pt}\text{such}\hspace{4pt}\text{that}\hspace{4pt}\sieve\subseteq\sieve_1\\
\text{for}\hspace{4pt}\text{some}\hspace{4pt}\text{set}\hspace{4pt}\sieve\in\precollarr\left(\objtre\right)\text{,}\hspace{4pt}\text{we}\hspace{4pt}\text{have} \hspace{4pt}\sieve_1\in\Grtop_1\left(\objdue\right)\text{;}
\end{array}
\right.\nonumber\\[6pt]
&\left\{
\begin{array}{l}
\text{any limit diagram}\;\objuno\overset{\nattrasf}{\rightarrow}\genfunc\;\text{in an object}\;\Catsei\;\text{of}\;\Catquattro\;\text{is a  diagram}\\
\text{in}\;\Catdue_1\;\text{if the functor}\;\fundim\left[\Catsei,\Catuno\right]\comparrfun
\genfunc\;\text{factors through}\;\fundim\left[\Catdue_1,\Catuno\right]\text{;}
\end{array}
\right.\nonumber\\[6pt]
&\left\{
\begin{array}{l}
\text{any co-limit diagram}\;\genfunc\overset{\nattrasf}{\rightarrow}\objuno\;\text{in an object}\;\Catsei\;\text{of}\;\Catquattro\;\text{is a diagram}\\
\text{in}\;\Catdue_1\;\text{if the functor}\;\fundim\left[\Catsei,\Catuno\right]\comparrfun
\genfunc\;\text{factors through}\;\fundim\left[\Catdue_1,\Catuno\right]\text{;}
\end{array}
\right.\nonumber\\[6pt]
&\left\{
\begin{array}{l}
\text{for any object}\;\left(\Catsei, \Catsette, \genfuncuno, \genfuncdue, \adjnattrasf \right)\;\text{of}\; \Catcinque\text{,}\;\text{object}\;\objsei\;\text{of}\; \Catsei\text{, object}\;\objsette\;\text{of}\; \Catsette\text{,}\\[4pt]
\arrsette\in\homF+\Catsette+\left(\genfuncuno\left(\objsei\right),\objsette\right)\;\text{all condition below are fulfilled:}\\[4pt]
\hspace{5pt}
\begin{array}{ll}
(i)&\text{if}\;\objsei\;\text{is an object of}\; \Catdue_1\;\text{then}\;\genfuncuno\left(\objsei\right)\;\text{is an object of}\;\Catdue_1\text{,}\\[4pt]
(ii)&\text{if}\;\objsette\;\text{is an object of}\; \Catdue_1\;\text{then}\;\genfuncdue\left(\objsette\right)\;\text{is an object of}\;\Catdue_1\text{,}\\[4pt]
(iii)&\arrsette\hspace{4pt}\text{is}\hspace{3pt}\text{an}\hspace{3pt}\text{arrow}\hspace{3pt}\text{of}\hspace{4pt}\Catdue_1\hspace{4pt}\text{if}\hspace{3pt}\text{and}\hspace{3pt}\text{only}\hspace{3pt}\text{if}\hspace{4pt} \adjnattrasf_{\objsei,\objsette}\left(\arrsette\right) \hspace{4pt}\text{is}\hspace{3pt}\text{an}\hspace{3pt}\text{arrow}\\
&\text{of}\hspace{4pt}\Catdue_1\text{,}\\[4pt]
(iv)&\text{if}\hspace{4pt}\objsette\hspace{4pt}\text{is}\hspace{4pt}\text{an}\hspace{4pt}\text{object}\hspace{4pt}\text{of}\hspace{4pt}\Catdue_1\hspace{4pt}\text{then}\hspace{4pt}\text{a}\hspace{4pt}\text{set}\hspace{4pt}\setsymnove\hspace{4pt}\text{of}\hspace{4pt}\text{arrows}\hspace{4pt}\text{of}\hspace{4pt}\Catdue_1\\[4pt]
&\text{fulfilling}\hspace{5pt}\text{condition}\\[4pt]
&\hspace{5pt}\text{for}\hspace{5pt}\text{any}\hspace{5pt} \arruno\in \setsymnove\hspace{5pt}\text{exists}\hspace{5pt}\text{an}\hspace{5pt}\text{object}\hspace{6pt} \objsei_{\arruno}\hspace{6pt}\text{of}\hspace{6pt}\Catsei\\[4pt]
&	\hspace{110pt}\text{with}\hspace{4pt} \arruno\in\homF+\Catsette+\left(\genfuncuno\left(\objsei_{\arruno}\right),\objsette\right)\\[4pt]
&\text{is}\hspace{5pt}\text{a}\hspace{5pt}\text{base}\hspace{5pt}\text{component}\hspace{5pt}\text{set}\hspace{5pt}\text{for}\hspace{6pt}
\Grtop_1\left(\objsette\right) \hspace{6pt}\text{if}\hspace{5pt}\text{and}\hspace{5pt}\text{only}\hspace{5pt}\text{if}\\[4pt]
&\left\{\adjnattrasf_{\objsei_{\arruno},\objsette}\left(\arruno\right)\,:\,\arruno\in \setsymnove\right\}\hspace{6pt}
\text{is}\hspace{6pt}\text{a}\hspace{6pt}\text{base}\hspace{6pt}
\text{component}\hspace{6pt}\text{set}\hspace{6pt}\text{for}\\[4pt]
&\Grtop_1\left(\genfuncdue\left(\objsette\right)\right)\text{;}
\end{array}
\end{array}
\right.\label{bcs}\\[6pt]
&\left\{\!
\begin{array}{l}
\text{an arrow}\; \objuno_1\overset{\arruno}{\rightarrow}\objuno_2\;\text{of}\;\Catuno\; \text{is an arrow of}\;\Catdue_1\;\text{if both conditions below}\\
\text{are fulfilled:}\\[4pt]
\hspace{3pt}(i) \hspace{6pt}\objuno_1\text{,}\;\objuno_2\; \text{are objects of}\; \Catdue_1\text{,}\\[4pt]
\hspace{2pt}(ii)\hspace{3pt} \text{for any}\; \sieve_2\in\Grtop_1\left(\objuno_2\right)\;\text{we have}\;
\left\{\arrdue\;:\;\arruno\comparrfun \arrdue\in \sieve_2\right\}\in\Grtop_1\left(\objuno_1\right)\text{.} 
\end{array}
\right.\nonumber
\end{flalign}
An augmentation $\left(\Catdue_0,\Grtop_0\right)$ of $\left(\Catdue,\Grtop\right)$ through $\left(\Catuno,\Cattre,\Catquattro, \Catcinque,\precollarr\right)$ is called canonical if and only if for any augmentation $\left(\Catdue_1,\Grtop_1\right)$ of $\left(\Catdue,\Grtop\right)$ through $\left(\Catuno,\Cattre,\Catquattro, \Catcinque,\precollarr\right)$ there is one and only one functor $\genfunc:\Catdue_0\rightarrow\Catdue_1$ fulfilling all conditions below:
\begin{flalign}
&\begin{array}{l}
\genfunc\comparrfun\inclfun{\Catdue}{\Catdue_0}=\inclfun{\Catdue}{\Catdue_1}\text{;}
\end{array}\\[6pt]
&\begin{array}{l}
\genfunc\comparrfun\inclfun{\Cattre}{\Catdue_0}=\inclfun{\Cattre}{\Catdue_1}\text{;}
\end{array}\\[6pt]
&\begin{array}{l}
\genfunc\left(\arruno\right)=\arruno\;\text{for any arrow} \;\arruno\in \precollarr\text{.}
\end{array}
\end{flalign}
\end{definition}

In Proposition \ref{augsiteprop} below we prove that the canonical augmentation of a site always exists. 
\begin{proposition}\label{augsiteprop}
Fix a category $\Catuno$, a site $\left(\Catdue,\Grtop\right)$ together with a forgetful functor $\fundim\left[\Catdue,\Catuno\right]$ forgetting the Grothendieck topology $\Grtop$, a subcategory $\Cattre$ of $\Catuno$, a category $\Catquattro$ whose objects are categories $\Catsei$ each one equipped with a forgetful functor $\fundim\left[\Catsei,\Catuno\right]$, a category $\Catcinque$ whose objects are adjunctions $\left(\Catsei, \Catsette, \genfuncuno, \genfuncdue, \adjnattrasf \right)$ where $\Catsei$ is equipped with a forgetful functor $\fundim\left[\Catsei,\Catuno\right]$ and $\Catsette$ is equipped with a forgetful functor $\fundim\left[\Catsette,\Catuno\right]$,  a function $\precollarr$ assigning to any object $\objtre$ of $\Cattre$ a set $\precollarr\left(\objtre\right)$ of sets of arrows of $\Catuno$ whose domains are objects of $\Catdue$ and co-domains all coincide with $\objtre$.\newline  
Then there is a canonical augmentation $\left(\Catdue_0,\Grtop_0\right)$ of the site $\left(\Catdue,\Grtop\right)$ through $\left(\Catuno,\Cattre,\Catquattro, \Catcinque,\precollarr\right)$.
\end{proposition}
\begin{proof}
Site $\left(\Catdue_0,\Grtop_0\right)$ is defined recursively, in particular Grothendieck topology $\Grtop_0$ on $\Catdue_0$ is given by recursively defining a base $\collarr_0$. The fact that $\left(\Catdue_0,\Grtop_0\right)$ is canonical will be evident from the construction.  
Precisely we set:
\begin{flalign}
&\left\{
\begin{array}{l} 
\text{for}\hspace{4pt}\text{any}\hspace{4pt}\text{object}\hspace{4pt}\objdue\hspace{4pt}\text{of}\hspace{4pt}\Catdue_0\text{,}\hspace{4pt} \text{any}\hspace{4pt}
\text{base}\hspace{4pt}\text{component}\hspace{4pt}\text{set}\hspace{4pt}\setsymnove\hspace{4pt}\text{for}\hspace{4pt}\Grtop_0\left(\objdue\right)\\[4pt]
\text{we}\hspace{4pt}\text{have}\hspace{4pt}\setsymnove\in\collarr_0\left(\objdue\right) \text{;}
\end{array}
\right.\label{bGt0}\\[6pt]
&\begin{array}{l}
\Catdue\;\text{and}\;\Cattre \;\text{are both subcategories of}\; \Catdue_0\text{;}
\end{array}\label{bGt1}
\\[6pt]
&\left\{
\begin{array}{l}
\text{any limit diagram}\;\objuno\overset{\nattrasf}{\rightarrow}\genfunc\;\text{in an object}\;\objquattro\;\text{of}\;\Catquattro\;\text{is a diagram in}\\
\Catdue_0\;\text{if the functor}\;\inclfun{\objquattro}{\Catuno}\comparrfun
\genfunc\;\text{factors through}\;\inclfun{\Catdue_0}{\Catuno}\text{;}
\end{array}
\right.\label{bGt2}
\\[6pt]
&\left\{
\begin{array}{l}
\text{any co-limit diagram}\;\genfunc\overset{\nattrasf}{\rightarrow}\objuno\;\text{in an object}\;\objquattro\;\text{of}\;\Catquattro\;\text{is a diagram in}\\
\Catdue_0\;\text{if the functor}\;\inclfun{\objquattro}{\Catuno}\comparrfun
\genfunc\;\text{factors through}\;\inclfun{\Catdue_0}{\Catuno}\text{;}
\end{array}
\right.\label{bGt2bis}
\\[6pt]
&\left\{
\begin{array}{l}
\text{for any object}\;\left(\Catsei, \Catsette, \genfuncuno, \genfuncdue, \adjnattrasf \right)\;\text{of}\; \Catcinque\text{,}\;\text{object}\;\objsei\;\text{of}\; \Catsei\text{, object}\;\objsette\;\text{of}\; \Catsette\text{,}\\[4pt]
\arrsette\in\homF+\Catsette+\left(\genfuncuno\left(\objsei\right),\objsette\right)\;\text{all condition below are fulfilled:}\\[4pt]
\hspace{5pt}
\begin{array}{ll}
(i)&\text{if}\;\objsei\;\text{is an object of}\; \Catdue_0\;\text{then}\;\genfuncuno\left(\objsei\right)\;\text{is an object of}\;\Catdue_0\text{,}\\[4pt]
(ii)&\text{if}\;\objsette\;\text{is an object of}\; \Catdue_0\;\text{then}\;\genfuncdue\left(\objsette\right)\;\text{is an object of}\;\Catdue_0\text{,}\\[4pt]
(iii)&\arrsette\hspace{4pt}\text{is}\hspace{3pt}\text{an}\hspace{3pt}\text{arrow}\hspace{3pt}\text{of}\hspace{4pt}\Catdue_0\hspace{4pt}\text{if}\hspace{3pt}\text{and}\hspace{3pt}\text{only}\hspace{3pt}\text{if}\hspace{4pt} \adjnattrasf_{\objsei,\objsette}\left(\arrsette\right) \hspace{4pt}\text{is}\hspace{3pt}\text{an}\hspace{3pt}\text{arrow}\\
&\text{of}\hspace{4pt}\Catdue_0\text{,}\\[4pt]
(iv)&\text{if}\hspace{5pt}\objsette\hspace{4pt}\text{is}\hspace{4pt}\text{an}\hspace{4pt}\text{object}\hspace{4pt}\text{of}\hspace{5pt}\Catdue_0\hspace{5pt}\text{then}\hspace{4pt}\text{a}\hspace{4pt}\text{set}\hspace{5pt}\setsymnove\hspace{5pt}\text{of}\hspace{4pt}\text{arrows}\hspace{4pt}\text{of}\hspace{5pt}\Catdue_0\\[4pt]
&\text{fulfilling}\hspace{5pt}\text{condition}\\[4pt]
&\hspace{5pt}\text{for}\hspace{5pt}\text{any}\hspace{5pt} \arruno\in \setsymnove\hspace{5pt}\text{exists}\hspace{5pt}\text{an}\hspace{5pt}\text{object}\hspace{6pt} \objsei_{\arruno}\hspace{6pt}\text{of}\hspace{6pt}\Catsei\\[4pt]
&	\hspace{115pt}\text{with}\hspace{4pt} \arruno\in\homF+\Catsette+\left(\genfuncuno\left(\objsei_{\arruno}\right),\objsette\right)\\[4pt]
&\text{belongs}\hspace{3pt}\text{to}\hspace{3pt}\collarr_0\left(\objsette\right) \hspace{3pt}\text{if}\hspace{3pt}\text{and}\hspace{3pt}\text{only}\hspace{3pt}\text{if}\hspace{3pt}\text{set}\hspace{3pt}
\left\{\adjnattrasf_{\objsei_{\arruno},\objsette}\left(\arruno\right)\,:\,\arruno\!\in\! \setsymnove\right\}\\[4pt]
&\text{belongs}\hspace{4pt}\text{to}\hspace{4pt}\collarr_0\left(
\genfuncdue\left(\objsette\right)\right) \text{;}
\end{array}
\end{array}
\right.\label{bGt2ter}
\\[6pt]
&\left\{\!
\begin{array}{l}
\text{an arrow}\; \objuno_1\overset{\arruno}{\rightarrow}\objuno_2\;\text{of}\;\Catuno\; \text{is an arrow of}\;\Catdue_0\;\text{if both conditions below}\\
\text{are fulfilled:}\\[4pt]
\hspace{3pt}(i) \hspace{6pt}\objuno_1\text{,}\;\objuno_2\; \text{are objects of}\; \Catdue_0\text{,}\\[4pt]
\hspace{2pt}(ii)\hspace{3pt} \text{for any}\; \setsymnove_2\in\collarr_0\left(\objuno_2\right)\;\text{there is}\;\setsymnove_1\in\collarr_0\left(\objuno_1\right) \;\text{such that}\\
\hspace{20pt}\text{for any}\; \arrdue_1\in \setsymnove_1\;\text{there are arrows}\;\arrdue_2\in \setsymnove_2\text{,}\;\widetilde{\arruno}\;\text{of}\;\Catdue_0\;\text{with}\\
 \hspace{20pt}\arruno \comparrfun\arrdue_1=\arrdue_2\comparrfun\widetilde{\arruno}\text{;} 
\end{array}
\right.\label{bGt3}
\\[6pt]
&\left\{
\begin{array}{l}
\text{for}\hspace{4pt}\text{any}\hspace{4pt}\text{object}\hspace{4pt}\objdue\hspace{4pt}\text{of}\;\Catdue\text{,}\hspace{4pt}\text{sieve} \hspace{4pt}\sieve_0\hspace{4pt}\text{on}\hspace{4pt}\objdue\hspace{4pt}\text{in}\hspace{4pt}\Catdue_0\hspace{4pt}\text{such}\hspace{4pt}\text{that}\hspace{4pt}\sieve\subseteq\sieve_0\\
\text{for}\hspace{4pt}\text{some}\hspace{4pt}\text{sieve}\hspace{4pt}\sieve\in\Grtop\left(\objdue\right)\text{,}\hspace{4pt}\text{we}\hspace{4pt}\text{have} \hspace{4pt}\sieve_0\in\collarr_0\left(\objdue\right)\text{;}
\end{array}
\right.\label{bGt4}\\[6pt]
&\left\{
\begin{array}{l}
\text{for}\hspace{4pt}\text{any}\hspace{4pt}\text{object}\hspace{4pt}\objtre\hspace{4pt}\text{of}\;\Cattre\text{,}\hspace{4pt}\text{sieve} \hspace{4pt}\sieve_0\hspace{4pt}\text{on}\hspace{4pt}\objtre\hspace{4pt}\text{in}\hspace{4pt}\Catdue_0\hspace{4pt}\text{such}\hspace{4pt}\text{that}\hspace{4pt}\setsymnove\subseteq\sieve_0\\
\text{for}\hspace{4pt}\text{some}\hspace{4pt}\text{set}\hspace{4pt}\setsymnove\in\precollarr\left(\objtre\right)\text{,}\hspace{4pt}\text{we}\hspace{4pt}\text{have} \hspace{4pt}\sieve_0\in\collarr_0\left(\objtre\right)\text{;}
\end{array}
\right.\label{bGt5}\\[6pt]
&\begin{array}{l}
\text{if}\;\objdue_1\overset{\arrdue}{\rightarrow}\objdue_2\;\text{is an isomorphism of}\;\Catdue_0\;\text{then}\;\left\{\arrdue\right\}\in\collarr_0\left(\objdue_2\right)\text{;}
\end{array}\label{bGt6}
\\[6pt]
&\left\{\!
\begin{array}{l}
\text{for any category}\;\Catindex\text{, functor} \;\Catindex\overset{\genfunc}{\rightarrow}\Catdue_0\text{, co-limit diagram} \;\genfunc\overset{\nattrasf}{\rightarrow}\objdue\;
\text{in an}\\
\text{object}\hspace{4pt}\text{of}\hspace{4pt}\Catquattro\text{,} \hspace{5pt}\text{function}\hspace{4pt} \varphi\hspace{4pt}\text{associating}\hspace{4pt}\text{to}\hspace{4pt}\text{any}\hspace{4pt}\text{object} \hspace{4pt}\objindex\hspace{4pt}\text{of} \;\Catindex\; \text{one}\hspace{4pt}\text{and}\\
\text{only one set in} \;
\collarr_0\left(\genfunc\left(\objindex\right)\right)\;
\text{we have}\\
\hspace{54pt}\underset{\objindex\;\text{object of}\;\Catindex}{\bigcup}\left\{\nattrasf\left(\objindex\right)\comparrfun \arrdue\;:\;\arrdue \in \varphi\left(\objindex\right) \right\}\in \collarr_0\left(\objdue\right)\text{;}
\end{array}
\right.\label{bGt7bis}
\\[6pt]
&\left\{
\begin{array}{l}
\text{for any arrow}\;\objdue_1\overset{\arrdue}{\rightarrow}\objdue_2\; \text{of}\;\Catdue_0\text{,}\;\setsymnove_2\in \collarr_0\left(\objdue_2\right)\;\text{we have that}\\[4pt]
\big\{\arrdue_1\;\text{arrow of}\;\Catdue_0\;:\;\text{there are arrows}\;\arrdue_2\in \setsymnove_2\text{,}\;\widetilde{\arrdue}\;\text{of}\;\Catdue_0\\
\hspace{148pt}\text{with}\;\arrdue \comparrfun\arrdue_1=\arrdue_2\comparrfun\widetilde{\arrdue}\big\}\in\collarr_0\left(\objdue_1\right)\text{;}
\end{array}
\right.\label{bGt7}
\\[6pt]
&\left\{\!
\begin{array}{l}
\text{for any object}\;\objdue\;\text{of}\;\Catdue_0\text{,}\; \text{set}\;\left\{\objdue_i\overset{\arrdue_i}{\rightarrow}\objdue\;:\;i \in I\right\}\in \collarr_0\left(\objdue\right)\text{,}\\
\text{set}\;\left\{\objdue_{i,j}\overset{\arrdue_{i,j}}{\rightarrow}\objdue_i\;:\;j \in J\right\}\in \collarr_0\left(\objdue_i\right)\text{ for any}\;i \in I\text{,}\\
\text{we have}\;\left\{\arrdue_i\comparrfun\arrdue_{i,j}\;:\;i \in I\text{,}\;j\in J\right\}\in \collarr_0\left(\objdue\right)\text{.}
\end{array}
\right.\label{bGt8}
\end{flalign}
To prove that $\left(\Catdue_0,\Grtop_0\right)$ is an augmentation of $\left(\Catdue,\Grtop\right)$ through $\left(\Catuno,\Cattre,\Catquattro, \Catcinque,\precollarr\right)$ we must show that $\Catdue_0$ is a category and that $\collarr_0$ is a base for the Grothendieck topology $\Grtop$ on $\Catdue_0$. Proof splits in six steps. \newline
\emph{\textbf{Step 1.}}\ \ We prove that $\idobj+\objdue+$ is an arrow of $\Catdue_0$ for any object $\objdue$ of $\Catdue_0$. This follows by condition \eqref{bGt3}.\newline
\emph{\textbf{Step 2.}}\ \ We prove that $\arrdue_2\comparrfun\arrdue_1$ is an arrow of $\Catdue_0$ for any pair of arrows $\objdue_1\overset{\arrdue_1}{\rightarrow}\objdue_2$, $\objdue\overset{\arrdue_2}{\rightarrow}\objdue_3$ in $\Catdue_0$. Arrows $\arrdue_1$, $\arrdue_2$ both fulfill condition \eqref{bGt7} by assumption, then $\arrdue_2\comparrfun\arrdue_1$ fulfills \eqref{bGt3}-[(i), (ii)], hence $\arrdue_2\comparrfun\arrdue_1$ is an arrow of $\Catdue_0$.\newline
\emph{\textbf{Step 3.}}\ \ Base first axiom. This axiom follows since condition \eqref{bGt6} entail that $\collarr_0$ fulfills axiom (i') of \cite{MM}\;\;Chap\,III\;\;\textsection\,2\,\;Definition 2.\newline
\emph{\textbf{Step 4.}}\ \ Base second axiom. This axiom follows since condition \eqref{bGt7bis} entail that $\collarr_0$ fulfills axiom (ii'') of \cite{MM}\;\;Chap\,III\;\;Exercise 3.\newline
\emph{\textbf{Step 5.}}\ \ Base third axiom. This axiom follows since condition \eqref{bGt8} entails that $\collarr_0$ fulfills axiom (iii') of \cite{MM}\;\;Chap\,III\;\;\textsection\,2\,\;Definition 2.\newline
\emph{\textbf{Step 6.}}\ \ We prove that $\left(\Catdue_0,\Grtop_0\right)$ fulfills \eqref{bcs}-[4]. Fix an object\;$\left(\Catsei, \Catsette, \genfuncuno, \genfuncdue, \adjnattrasf \right)$ of $\Catcinque$, an object $\objsette$  of $\Catsette$, a set $\setsymnove$ of arrows of $\Catdue_0$ fulfilling condition 
\begin{equation*}
\text{for any}\; \arruno\in \setsymnove\; \text{exists an object}\; \objsei_{\arruno}\;\text{of}\;\Catsei\;
\text{with}\; \arruno\in\homF+\Catsette+\left(\genfuncuno\left(\objsei_{\arruno}\right),\objsette\right)\text{.}
\end{equation*}
Then condition \eqref{bcs}-[4] is fulfilled since:\newline
$\setsymnove$ is a base component set of $\Grtop_0\left(\objsette \right) $ if and only if $\setsymnove \in \collarr_0\left(\objsette \right)$ by \eqref{bGt0};\newline
$\setsymnove \in \collarr_0\left(\objsette \right)$ if and only if $\left\{\adjnattrasf_{\objsei_{\arruno},\objsette}\left(\arruno\right)\;: \;\arruno\in \setsymnove\right\}\in \collarr\left(
\genfuncdue\left(\objsette\right)\right)$ by \eqref{bGt2ter}-[iv];\newline
$\left\{\adjnattrasf_{\objsei_{\arruno},\objsette}\left(\arruno\right)\;:\;\arruno\in \setsymnove\right\}$ is a base component set of $\collarr\left(
\genfuncdue\left(\objsette\right)\right)$ if and only if we have $\left\{\adjnattrasf_{\objsei_{\arruno},\objsette}\left(\arruno\right)\;:\;\arruno\in \setsymnove\right\}\in \collarr\left(\genfuncdue\left(\objsette\right)\right)$ by \eqref{bGt0}. 
\end{proof}

In Notation \ref{locappa} below we establish a short terminology to denote the behavior of arrows with respect to covering sieves. We refer to Notation \ref{catfun}-[9].

\begin{notation}\label{locappa}
Fix a category $\Catuno$, a site $\left(\Catdue,\Grtop\right)$ together with a forgetful functor $\fundim\left[\Catdue,\Catuno\right]$ forgetting the Grothendieck topology $\Grtop$, a subcategory $\Cattre$ of $\Catuno$, a category $\Catquattro$ whose objects are categories $\Catsei$ each one equipped with a forgetful functor $\fundim\left[\Catsei,\Catuno\right]$, a category $\Catcinque$ whose objects are adjunctions $\left(\Catsei, \Catsette, \genfuncuno, \genfuncdue, \adjnattrasf \right)$ where $\Catsei$ is equipped with a forgetful functor $\fundim\left[\Catsei,\Catuno\right]$ and $\Catsette$ is equipped with a forgetful functor $\fundim\left[\Catsette,\Catuno\right]$,  a function $\precollarr$ assigning to any object $\objtre$ of $\Cattre$ a set $\precollarr\left(\objtre\right)$ of sets of arrows of $\Catuno$ whose domains are objects of $\Catdue$ and co-domains all coincide with $\objtre$, a canonical augmentation $\left(\Catdue_0,\Grtop_0\right)$ of the site $\left(\Catdue,\Grtop\right)$ through $\left(\Catuno,\Cattre,\Catquattro, \Catcinque,\precollarr\right)$.
\begin{enumerate}
\item We locally characterize arrows of $\Catdue_0$ as follows.\newline
Fix a base $\collarr_0$ for the Grothendieck topology $\Grtop_0$, two objects $\objdue_1$, $\objdue_2$ of $\Catdue_0$, an arrow $\objdue_1\overset{\arruno}{\rightarrow}\objdue_2$ of $\Catuno$, a covering sieve $\sieve \in \Grtop_0\left(\objdue_2\right)$, a set $\setsymnove\in\collarr_0\left(\objdue_2\right)$.\newline
We say that $\arruno$ locally belongs to $\sieve$ if and only if the pull back sieve $\arruno^*\left(\sieve\right)$ is a covering sieve belonging to $\Grtop_0\left(\objdue_1\right)$.\newline 
We say that $\arruno$ locally factors through $\setsymnove$ if and only if the pull back sieve $\arruno^*\left(\sieve_{\setsymnove}\right)$ is a covering sieve belonging to $\Grtop_0\left(\objdue_1\right)$.
\item Fix an object $\objtre$ of $\Cattre$, an arrow $\arruno$ belonging to a set belonging to $\precollarr\left(\objtre\right)$. We say that $\arruno$ is a special arrow of $\Catuno$.
\end{enumerate} 
\end{notation}

In Definitions \ref{specfunV}, \ref{specfunVCC} below we introduce classes of set functions which are crucial to define Grothendieck topology of the site where spaces of generalized derivations of degree $1$ live in. We refer to Definition \ref{dualnot0}, Proposition \ref{genpointprop}, Remark \ref{dualrem}-[1].

\begin{definition}\label{specfunV}
We define special set functions recursively as follows:
\begin{enumerate}
\item Fix $\gfx\in\genf$. The set function $\preffuncaap{\gfx}=\genfquotfun\funcomp\prefuncaap{\gfx}$ is a special set function.
\item  Fix a topological space $\setsymuno$, $\gfx \in \genf$, a special set function $\pmV:\setsymuno\rightarrow\genfquot$.\newline  
Set function $\pmV_1:\setsymuno\rightarrow\genfquot$ defined by setting
\begin{equation*}
\pmV_1\left(\elsymuno\right)=\gfx \compgenfquot \pmV\left(\elsymuno\right)\hspace{10pt}\forall \elsymuno\in \setsymuno
\end{equation*}
is a special set function.
\item Fix two topological spaces $\setsymuno_1$, $\setsymuno_2$, a special set function $\pmV_i:\setsymuno_i\rightarrow\genfquot$ for any $i\in\left\{1,2\right\}$ with $\pmV_2$ factoring through $\incl{\genfquotcccz}{\genfquot}$.\newline 
Set function $\pmV:\setsymuno_1\times\setsymuno_2\rightarrow\genfquot$ defined by setting 
\begin{equation*}
\pmV\left(\elsymuno_1,\elsymuno_2\right)= \pmV_1\left(\elsymuno_1\right)\compquotpoint\pmV_2\left(\elsymuno_2\right)\hspace{10pt}\forall \left(\elsymuno_1,\elsymuno_2\right)\in \setsymuno_1\times\setsymuno_2
\end{equation*} 
is a special set function.
\item Fix $\topindextre\in\mathbb{N}$, a topological space $\setsymuno_{\elindextre}$, a special set function $\pmV_{\elindextre}:\setsymuno_{\elindextre}\rightarrow\genfquot$ for any $\elindextre\in \left\{1,...,\topindextre\right\}$.\newline
Set function $\pmV:\overset{\topindextre}{\underset{\elindextre=1}{\prod}}\setsymuno_{\elindextre}\rightarrow\genfquot$ defined  by setting 
\begin{equation*}
\pmV\left(\elsymuno_1,...,\elsymuno_{\topindextre}\right)=\lboundquotpoint\pmV_1\left(\elsymuno_1\right),...,\pmV_{\topindextre}\left(\elsymuno_{\topindextre}\right) \rboundquotpoint\hspace{10pt}\forall \left(\elsymuno_1,...,\elsymuno_{\topindextre}\right) \in \overset{\topindextre}{\underset{\elindextre=1}{\prod}}\setsymuno_{\elindextre}
\end{equation*} 
is  a special set function.
\end{enumerate}
\end{definition}

\begin{remark}\label{specfunrem} Definition \ref{specfunV} can be repeated word by word in the setting of $\Ckspquot{\infty}$ and $\Ckspquot{0}$, taking in account Remark \ref{existoperquot}-[2], to get the notion of special set functions in $\Ckspquot{\infty}$ and $\Ckspquot{0}$ respectively.
Special set functions in $\Ckspquot{\infty}$ and $\Ckspquot{0}$ are denoted by symbols used for corresponding notions on generalized functions.
\end{remark}

\begin{definition}\label{specfunVCC}
We define pointed special set functions recursively as follows:
\begin{enumerate}
\item Fix $\gfx \in \genfcontcont$. Then $\ffuncaap{\gfx}=\genfquotfun\funcomp\funcaap{\gfx}$ is a pointed special set function.
\item  Fix a topological space $\setsymuno$, $\gfx \in \genfcontcont$, a pointed special set function $\pmV:\setsymuno\rightarrow\genfquotcccz$.\newline  
Set function $\pmV_1:\setsymuno\rightarrow\genfquotcccz$ defined by setting
\begin{equation*}
\pmV_1\left(\elsymuno\right)=\gfx \compgenfquot \pmV\left(\elsymuno\right)\hspace{10pt}\forall \elsymuno\in \setsymuno
\end{equation*}
is a pointed special set function.
\item Fix two topological spaces $\setsymuno_1$, $\setsymuno_2$, a pointed special set function $\pmV_i:\setsymuno_i\rightarrow\genfquotcccz$ for any $i\in\left\{1,2\right\}$.\newline 
Set function $\pmV:\setsymuno_1\times\setsymuno_2\rightarrow\genfquotcccz$ defined by setting 
\begin{equation*}
\pmV\left(\elsymuno_1,\elsymuno_2\right)= \pmV_1\left(\elsymuno_1\right)\compquotpoint\pmV_2\left(\elsymuno_2\right)\hspace{10pt}\forall \left(\elsymuno_1,\elsymuno_2\right)\in \setsymuno_1\times\setsymuno_2
\end{equation*} 
is a pointed special set function.
\item Fix $\topindextre\in\mathbb{N}$, a topological space $\setsymuno_{\elindextre}$, a pointed special set function $\pmV_{\elindextre}:\setsymuno_{\elindextre}\rightarrow\genfquotcccz$ for any $\elindextre\in \left\{1,...,\topindextre\right\}$.\newline
Set function $\pmV:\overset{\topindextre}{\underset{\elindextre=1}{\prod}}\setsymuno_{\elindextre}\rightarrow\genfquotcccz$ defined  by setting 
\begin{equation*}
\pmV\left(\elsymuno_1,...,\elsymuno_{\topindextre}\right)=\lboundquotpoint\pmV_1\left(\elsymuno_1\right),...,\pmV_{\topindextre}\left(\elsymuno_{\topindextre}\right) \rboundquotpoint\hspace{10pt}\forall \left(\elsymuno_1,...,\elsymuno_{\topindextre}\right) \in \overset{\topindextre}{\underset{\elindextre=1}{\prod}}\setsymuno_{\elindextre}
\end{equation*} 
is a pointed special set function.
\end{enumerate}
\end{definition}

\begin{remark}\label{specfunremVCC} Definition \ref{specfunVCC} can be repeated word by word in the setting of $\Ckspquot{\infty}$ and $\Ckspquot{0}$, taking in account Remark \ref{existoperquot}-[2], to get the notion of pointed special set functions in $\Ckspquotccz{\infty}$ and $\Ckspquotccz{0}$ respectively.
Pointed special set functions in $\Ckspquotccz{\infty}$ and $\Ckspquotccz{0}$ are denoted by symbols used for corresponding notions on generalized functions.
\end{remark}

In Definition \ref{defcatD} below we define the site where spaces of generalized derivation of degree $1$ live in.
We refer to Definitions \ref{genpointdef}, \ref{specfunV}, \ref{specfunVCC}, \ref{domcodspfdef}, Remarks \ref{genercontfun}, \ref{existoperquot}-[2]. 

\begin{definition}\label{defcatD}\mbox{}
\begin{enumerate}
\item The category $\Cattre$ is the category containing all diagrams existing in some abelian category of sets involving objects and arrows listed below:\newline
Objects: all subsets of $\genf$;  all $\mathbb{R}$-vector subspaces of $\mathbb{R}$, $\genfquot$, $\Ckspquot{\infty}$, $\Ckspquot{0}$, $\dirsumgenpderspace[\modsym](l)+m+$ for any $l,m \in \mathbb{N}_0$, $\mathbb{R}$-vector subspace $\modsym$ of $\mathbb{R}^l$.\newline
Arrows: all inclusion functions between objects listed above; all functions between sets $\Ckspquot{\infty}$, $\Ckspquot{0}$, $\genfquot$ defined in Proposition \ref{genpointprop}, Remark \ref{existoperquot}-[2]; sum and scalar product of $\dirsumgenpderspace[\modsym](l)+m+$; quotient functions $\quotdirsumuno[\modsym](l)+m+$ defined in Remark \ref{pathinlpreder}-[2]; evaluation functions $\evspd[\modsym](l)+m+$ introduced in Proposition \ref{gpdssmprop}; functions $\dirsumgenprediff[\modsym](l)+\gggfx[\gfx]+$ introduced in Remark \ref{pathinlpreder}-[4].  
\item The category $\Catquattro$ is the category whose objects are categories of sets.
\item The category $\Catcinque$ is the category whose objects are adjunctions between categories of sets.
\item The function $\precollarr$ is defined by setting:\newline
\begin{flalign}
&\begin{array}{l}
\precollarr\left(\genf\right)=\left\{\udenset\right\}\text{.}\label{discgenf}
\end{array}\\[10pt]
&\begin{array}{l}
\precollarr\left(\mathbb{R}\right)\;\text{is the set of all sieves induced by Euclidean topology.}
\end{array}\\[10pt]
& \begin{array}{l}
\precollarr\left(\genfquot\right)\hspace{5pt}\text{contains}\hspace{6pt}\text{only}\hspace{6pt}\text{the}\hspace{6pt}\text{set}\hspace{5pt}\ssf{\genfquot}\hspace{5pt}\text{of}\hspace{6pt}\text{all}\hspace{6pt}\text{special}\hspace{6pt}\text{set}\\
\text{functions}\hspace{4pt}\text{introduced in Definition \ref{specfunV}.}
\end{array}\label{spfunV}\\[10pt]
& \begin{array}{l}
\precollarr\left(\genfquotcccz\right)\hspace{5pt}\text{contains}\hspace{5pt}\text{only}\hspace{5pt}\text{the}\hspace{5pt}\text{set}\hspace{5pt}\ssf{\genfquotcccz}\hspace{5pt}\text{of}\hspace{5pt}\text{all}\hspace{5pt}\text{pointed}\\
\text{special}\hspace{4pt}\text{set}\hspace{4pt}
\text{functions}\hspace{4pt}\text{introduced}\hspace{4pt}\text{in}\hspace{4pt}\text{Definition}\hspace{4pt}\text{\ref{specfunVCC}.}
\end{array}\label{spfunVCC}\\[10pt]
& \begin{array}{l}
\precollarr\left(\Ckspquot{\infty}\right)\hspace{4pt}\text{contains}\hspace{4pt}\text{only}\hspace{4pt}\text{the}\hspace{4pt}\text{set}\hspace{4pt}\ssf{\Ckspquot{\infty}}\hspace{4pt}\text{of}\hspace{4pt}\text{all}\hspace{4pt}\text{special}\hspace{4pt}\text{set}\\
\text{functions}\hspace{4pt}\text{introduced}\hspace{6pt}\text{in}\hspace{6pt}\text{Remark}\hspace{4pt}\text{\ref{specfunrem}.}
\end{array}\\[10pt]
& \begin{array}{l}
\precollarr\left(\Ckspquotccz{\infty}\right)\hspace{4pt}\text{contains}\hspace{4pt}\text{only}\hspace{4pt}\text{the}\hspace{4pt}\text{set}\hspace{4pt}\ssf{\Ckspquotccz{\infty}}\hspace{4pt}\text{of}\hspace{4pt}\text{all}\hspace{4pt}\text{pointed}\\
\text{special}\hspace{4pt}\text{set}\hspace{4pt}
\text{functions}\hspace{4pt}\text{introduced}\hspace{6pt}\text{in}\hspace{6pt}\text{Remark}\hspace{4pt}\text{\ref{specfunremVCC}.}
\end{array}\\[10pt]
& \begin{array}{l}
\precollarr\left(\Ckspquot{0}\right)\hspace{4pt}\text{contains}\hspace{4pt}\text{only}\hspace{4pt}\text{the}\hspace{4pt}\text{set}\hspace{4pt}\ssf{\Ckspquot{0}}\hspace{4pt}\text{of}\hspace{4pt}\text{all}\hspace{4pt}\text{special}\hspace{4pt}\text{set}\\
\text{functions}\hspace{4pt}\text{introduced}\hspace{6pt}\text{in}\hspace{6pt}\text{Remark}\hspace{6pt}\text{\ref{specfunrem}.}
\end{array}\\[10pt]
& \begin{array}{l}
\precollarr\left(\Ckspquotccz{0}\right)\hspace{4pt}\text{contains}\hspace{4pt}\text{only}\hspace{4pt}\text{the}\hspace{4pt}\text{set}\hspace{4pt}\ssf{\Ckspquotccz{0}}\hspace{4pt}\text{of}\hspace{4pt}\text{all}\hspace{4pt}\text{pointed}\\
\text{special}\hspace{4pt}\text{set}\hspace{4pt}
\text{functions}\hspace{4pt}\text{introduced}\hspace{6pt}\text{in}\hspace{6pt}\text{Remark}\hspace{4pt}\text{\ref{specfunremVCC}.}
\end{array}\\[10pt]
&\begin{array}{l}
\precollarr\left(\dirsumgenpderspace[\modsym](l)+m+\right)\hspace{3pt}\text{contains}\hspace{3pt}\text{only}\hspace{3pt}\text{the}\hspace{3pt}\text{set}\hspace{3pt}\ssf{\dirsumgenpderspace[\modsym](l)+m+}\hspace{3pt}\text{of}\hspace{3pt}\text{all}\\
\text{set}\hspace{3pt}\text{functions}\hspace{4pt}
\pmdspreder:\setsymuno\rightarrow\dirsumgenpderspace[\modsym](l)+m+\hspace{3pt}\text{such}\hspace{3pt}\text{that}\hspace{4pt}\setsymuno\hspace{4pt}\text{is}\hspace{3pt}\text{a}\hspace{3pt}\text{topological}\\
\text{space}\hspace{4pt}\text{and}\hspace{4pt}
\text{there}\hspace{4pt}\text{are}\hspace{4pt}\text{a}\hspace{4pt}\text{finite}\hspace{4pt}\text{set}\hspace{5pt}\setindexdue\text{,}\hspace{5pt}\text{a}\hspace{4pt}\text{pointed}\hspace{4pt}\text{special}\hspace{4pt}\text{set}\hspace{4pt}\text{funct-}\\
\text{ion}\hspace{6pt}
\pmV_{\elindexdue}:\setsymuno\rightarrow
\genfquotcccz(l)+m+\text{,}\hspace{6pt}
\text{a}\hspace{6pt}\text{continuous}\hspace{6pt}\text{function}\hspace{6pt}f_{\elindexdue}:\setsymuno\rightarrow\modsym\\
\text{for}\hspace{4pt}
\text{any}\hspace{4pt}\elindexdue\in\setindexdue\hspace{4pt}\text{such}\hspace{4pt}\text{that}\hspace{4pt}
\pmdspreder\left(\elsymuno\right)=\underset{\elindexdue\in  \setindexdue}{\sum} \genpreder\left[\pmV_{\elindexdue}\left(\elsymuno\right),f_{\elindexdue}\left(\elsymuno\right) \right]\text{.} \\[4pt]
\text{We}\hspace{3pt}\text{say}\hspace{3pt}\text{that:}\hspace{4pt}\left\{\pmV_{\elindexdue}\right\}_{\elindexdue\in\setindexdue}\hspace{4pt}\text{is}\hspace{4pt}\text{a}\hspace{4pt}\text{core}\hspace{4pt}\text{of}\hspace{4pt}\pmdspreder\text{;}\hspace{4pt}\left\{f_{\elindexdue}\right\}_{\elindexdue\in\setindexdue}\hspace{4pt}\text{is}\hspace{4pt}\text{an}\hspace{4pt}\text{associating}\\
\text{direction}\hspace{6pt}\text{of}\hspace{6pt}\pmdspreder\text{;}\hspace{6pt}\left\{\left(\pmV_{\elindexdue}, f_{\elindexdue}\right)\right\}_{\elindexdue\in\setindexdue}\hspace{6pt}\text{is}\hspace{6pt}\text{a}\hspace{6pt}\text{skeleton}\hspace{5pt}\text{of}\hspace{5pt}\pmdspreder\text{.}\\
\text{With}\hspace{5pt}\text{an}\hspace{5pt}\text{abuse}\hspace{5pt}\text{of}\hspace{5pt}\text{language}\hspace{5pt}\text{we}\hspace{5pt}\text{denote}\hspace{5pt}\pmdspreder\hspace{5pt}\text{by}\hspace{4pt}\genpreder\left[\pmV, f\right] \hspace{4pt}\text{and}\hspace{5pt}	\text{say}\\
\text{that}\hspace{4pt}\text{elements}\hspace{4pt}\text{belonging}\hspace{4pt}\text{to}\hspace{6pt}\ssf{\dirsumgenpderspace[\modsym](l)+m+}\hspace{6pt}\text{are}\hspace{4pt}\text{special}\hspace{4pt}\text{set}\\
\text{functions.}
\end{array}\label{pcM}\\[10pt]
& \begin{array}{l}
\precollarr\left(\objtre\right)=\udenset\hspace{5pt}\text{for}\hspace{4pt}\text{any}\hspace{4pt}\text{object}\hspace{4pt}\objtre\hspace{4pt}\text{of}\hspace{5pt}\Cattre\hspace{5pt}\text{different}\hspace{4pt}\text{form}\hspace{6pt}\mathbb{R}\text{,}\hspace{6pt}\Ckspquot{\infty}\text{,}\\
\Ckspquot{0}\text{,}\hspace{6pt}
\genfquot\text{,}\hspace{6pt}\dirsumgenpderspace[\modsym](l)+m+\text{.}\nonumber
\end{array}
\end{flalign}
\item We denote by $\left(\GTopcat, \csGT\right)$ the site obtained by augmentation of $\left(\Topcat,\GrTc\right)$ through $\left(\Setcat, \Cattre,\Catquattro, \Catcinque, \precollarr\right)$ constructed in Proposition \ref{augsiteprop}.\newline
We denote by $\terobj+\GTopcat+$ the terminal object of $\GTopcat$, that is the single point set.
\end{enumerate} 
\end{definition}

\begin{conjecture}
Grothndieck topology $\csGT$ is subcanonical.\newline
If this conjecture is true then the formulation of calculus in $\GTopcat$ (Chapter \ref{bascalc}) would be considerably easier.  
\end{conjecture}

In Proposition \ref{pointspecfun} below we study the Grothendieck topology of $\genfquot$. We refer to Notation \ref{ins}-[7], Definitions \ref{dualnot0}, \ref{genpointdef}, Proposition \ref{genpointprop}-[5, 8].

\begin{proposition}\label{pointspecfun}\mbox{}
\begin{enumerate}
\item Fix $m,n \in \mathbb{N}_0$, an object $\setsymuno$ of $\GTopcat$, $\gggfx\in\genfquot(m)+n+$.\newline
Then $\cost<\setsymuno<>\argcompl{\genfquot(m)+n+}>+\gggfx+\in\homF+\GTopcat+\left(\setsymuno,\genfquot\right)$.
\item  Fix $m,n \in \mathbb{N}_0$, $\pmV \in\homF+\GTopcat+\left(\setsymuno,\genfquot(m)+n+\right)$.\newline
Then the set function $\pmV_1:\setsymuno\rightarrow\genfquot(m)+n+$ defined by setting 
\begin{equation*}
\pmV_1\left(\elsymuno\right)=\pmV\left(\elsymuno\right)\sumquotpoint \left(-1 \scalpquotpoint\pmV\left(\elsymuno\right)  \compquotpoint \lininclquot\left(\gczsfx[\argcompl{\cost<\mathbb{R}^m<>\mathbb{R}^m>+0+}]\right)\right)
\hspace{15pt}\forall \elsymuno\in \setsymuno
\end{equation*}
is an arrow of $\GTopcat$.
\item  The pull back sieve of any covering sieve $\sieve\in \csGT\left(\genfquot\right)$ contains $\ssf{\genfquotcccz}$.
\end{enumerate}
\end{proposition}
\begin{proof}\mbox{}\newline
\textnormal{\textbf{Proof of statement 1.}}\ \  Fix an open neighborhood $\intuno\subseteq \mathbb{R}^m$ of 0, $\gfx\in \genfloc(\intuno)+\mathbb{R}^n+$ with $\genfquotfun\left(\gfx\right)=\gggfx$. 
Statement follows by Definition \ref{specfunV}-[1] since \eqref{bGt5} entail that the set function $\preffuncaap{\gfx} \funcomp\cost<\setsymuno<>\intuno>+0+$ is an arrow of $\GTopcat$.\newline
\textnormal{\textbf{Proof of statement 2.}}\ \ Statement follows by statement 1 since $\pmV_1$ is obtained by suitably assembling special set functions $\pmV$, $\cost<\setsymuno<>\argcompl{\genfquot(m)+n+}>+\argcompl{\lininclquot\left(\gczsfx[\argcompl{\cost<\mathbb{R}^m<>\mathbb{R}^n>+0+}]\right)}+$ and continuous function $\Diag{\setsymuno}{2}$ defined in Notation \ref{ins}-[8] by exploiting Definition \ref{specfunV}-[2, 3, 4].\newline
\textnormal{\textbf{Proof of statement 3.}}\ \ Statement follows by statement 2.
\end{proof}

In Remark \ref{siterem} below we explicitly describe some features of the site $\left(\GTopcat, \csGT\right)$ which are evident from its construction. We refer to Notations \ref{gentopnot}-[7], \ref{realvec}-[2], \ref{realfunc}-[6]. 

\begin{remark}\label{siterem}\mbox{}
\begin{enumerate}
\item Many objects of $\GTopcat$ share the same underlying set. Few examples are listed below. The set $\genf$ is endowed with both the corresponding discrete topology (see condition \eqref{discgenf}) and the topology given in Chapter \ref{magtwosec}. Sets $\Ckspquot{\infty}$, $\Ckspquotccz{\infty}$, $\Ckspquot{0}$, $\Ckspquotccz{0}$ are endowed with both topologies described in Examples \ref{ExCinf}, \ref{ExCzero} and topologies given by special and pointed special set functions.
\item Proposition \ref{pointspecfun} entails that $\incl{\genfquotcccz}{\genfquot}$ is in $\GTopcat$ by \eqref{bGt3}, then it is redundant to require that $\incl{\genfquotcccz}{\genfquot}$ is an arrow of $\Cattre $ in Definition \ref{defcatD}-[1].\newline
The corresponding statement hold true for arrows $\incl{\Ckspquotccz{\infty}}{\Ckspquot{\infty}}$, $\incl{\Ckspquotccz{0}}{\Ckspquot{0}}$, $\incl{\Ckspquot{\infty}}{\Ckspquot{0}}$, $\incl{\Ckspquotccz{\infty}}{\Ckspquotccz{0}}$ when $\Ckspquot{\infty}$, $\Ckspquotccz{\infty}$, $\Ckspquot{0}$, $\Ckspquotccz{0}$ are endowed with both topologies described in Examples \ref{ExCinf}, \ref{ExCzero} or topologies given by special and pointed special set functions. 
\item In category $\GTopcat$ exist all diagrams giving limits and co-limits in objects of $\Catquattro$. We adopt the same notation to denote limit and co-limit diagrams existing in objects of $\Catquattro$ and the corresponding diagrams in $\GTopcat$.\newline
Since all abelian categories of sets are objects of the category $\Catquattro$, then vector spaces and linear functions introduced in Definition \ref{genderdef2}, Remark \ref{pathinlpreder} are all objects and arrows of $\GTopcat$.\newline
For any objects $\setsymuno_1$, $\setsymuno_2$, $\setsymuno_3$, $\setsymuno_4$ of $\GTopcat$, arrow $\bot:\setsymuno_2\times\setsymuno_3\rightarrow \setsymuno_4$ of $\GTopcat$ there are uniquely defined arrows: 
\begin{flalign*}
&\setsymuno_2\times\homF+\GTopcat+\left(\setsymuno_1,\setsymuno_3\right)\rightarrow \homF+\GTopcat+\left(\setsymuno_1,\setsymuno_4\right)
\text{;}\\[6pt]
&\homF+\GTopcat+\left(\setsymuno_1,\setsymuno_3\right)\times\setsymuno_3\rightarrow \homF+\GTopcat+\left(\setsymuno_1,\setsymuno_4\right)\text{;}\\[6pt]
& \homF+\GTopcat+\left(\setsymuno_1,\setsymuno_2\right)\times\homF+\GTopcat+\left(\setsymuno_1,\setsymuno_3\right)\rightarrow \homF+\GTopcat+\left(\setsymuno_1,\setsymuno_4\right)\text{.}
\end{flalign*}
With an abuse of language we denote both by the same symbol $\bot$ being clear, from the unknowns which it acts on, which is the arrows we are dealing with.
\item In category $\GTopcat$ all matches exist between arrows given by adjunctions which are objects of $\Catcinque$. We adopt the same notation to denote adjunctions of $\Catcinque$ and their corresponding in $\GTopcat$. 
\item Fix an object $\setsymuno$ of $\GTopcat$.\newline
Conditions \eqref{bGt3}, \eqref{bGt5}, \eqref{spfunV} together entail that for any arrow $\setsymuno\overset{\pmV}{\rightarrow}\genfquot$ of $\GTopcat$, any covering sieve $\sieve\in \csGT\left(\genfquot\right)$ we have that $\pmV$ locally belongs to $\sieve$, in particular $\pmV$ locally factors through special set functions (see Definition \ref{specfunV}).\newline
Conditions \eqref{bGt3}, \eqref{bGt5}, \eqref{spfunVCC} and Proposition \ref{pointspecfun} together entail that for any arrow $\setsymuno\overset{\pmV}{\rightarrow}\genfquotcccz$ of $\GTopcat$, any covering sieve $\sieve\in \csGT\left(\genfquotcccz\right)$ we have that $\pmV$ locally belongs to $\sieve$, in particular $\pmV$ locally factors through pointed special set functions (see Definition \ref{specfunVCC}).\newline
This entails that set functions defined in Proposition \ref{genpointprop} are all continuous functions.\newline
In particular: function $\compgenfquot$ defined in Proposition \ref{genpointprop}-[1] is continuous with respect to discrete topology on $\genf$ (see 1) and product topology on $\genf\times\genfquotcccz$; function $\bsfmuquotpoint$ defined in Proposition \ref{genpointprop}-[9] is continuous with respect to discrete topology on $\bsfM$ (see 1) and product topology on $\bsfM\times\genfquotcccz$.\newline
The corresponding statement hold true for $\Ckspquotccz{\infty}$, $\Ckspquotccz{0}$, in particular set functions defined in Proposition \ref{genpointprop}-[10, 11] are continuous functions.
\item Fix two objects $\setsymuno$, $\setsymdue$ of $\Topcat$, an arrow $\setsymuno\overset{f}{\rightarrow}\setsymdue$ of $\GTopcat$. Then $f$ is an arrow of $\Topcat$. In fact $f$ is an arrow of $\GTopcat$ if and only if it locally belongs to any sieve of $\csGT\left(\setsymdue\right)$, in particular it locally factors through arrows of $\Topcat$ by \eqref{bGt4}.
\item Fix $l,m\in \mathbb{N}_0$,  a $\mathbb{R}$-vector subspace $\modsym$ of $\mathbb{R}^l$. We set
\begin{multline}
\ssf{\genprederspace[\modsym](l)+m+}=\\
\bigg\{\quotdirsumuno[\modsym](l)+m+\funcomp \pmdspreder\;:\;\pmdspreder\in\ssf{\dirsumgenpderspace[\modsym](l)+m+}\bigg\}\text{.}
\label{specpredersieve0}
\end{multline}
We say that elements belonging to set $\ssf{\genprederspace[\modsym](l)+m+}$ are special set functions.\newline
By \eqref{bGt7bis}, \eqref{pcM} and referring to Notation \ref{catfun}-[9] we have that set $\ssf{\genprederspace[\modsym](l)+m+}$ is a base component set for $\csGT\left(\genprederspace[\modsym](l)+m+\right)$.
\end{enumerate}
\end{remark}

In Proposition \ref{domcodspstfun} below we study the behavior of arrows of $\GTopcat$ with path connected domain.

\begin{proposition}\label{domcodspstfun}
Fix a path connected topological space $\setsymuno$, $\pmV \in\homF+\GTopcat+\left(\setsymuno,\genfquot\right)$.
Then there are $l,m\in\mathbb{N}_0$ such that
\begin{equation*}
\dimcodmagtwo\left(\pmV\left(\elsymuno\right)\right)=m\text{,}\hspace{20pt}
\dimdommagtwo\left(\pmV\left(\elsymuno\right)\right)=l\hspace{20pt}\forall\elsymuno\in \setsymuno\text{.}\\[4pt]
\end{equation*}
\end{proposition}
\begin{proof}
Statement follows by referring to Remark \ref{siterem}-[5] and exploiting recursive definition of special and pointed special set functions.
\end{proof}

Motivated by Proposition \ref{domcodspstfun} we define domain and co-domain of the image of special and pointed special set functions.

\begin{definition}\label{domcodspfdef} Fix a path connected topological space $\setsymuno$, $\pmV \in\homF+\GTopcat+\left(\setsymuno,\genfquot\right)$.
We define 
\begin{flalign*}
&\dimcodmagtwo\left(\pmV\right)=\dimcodmagtwo\left(\pmV\left(\elsymuno\right)\right)\hspace{10pt}\elsymuno\in \setsymuno\text{;}\\[4pt]
&\dimdommagtwo\left(\pmV\right)=\dimdommagtwo\left(\pmV\left(\elsymuno\right)\right)\hspace{10pt}\elsymuno\in \setsymuno\text{.}
\end{flalign*}
\end{definition}

\begin{remark}\label{siteremdue}
Proposition \ref{domcodspstfun} and Definition \ref{domcodspfdef} can be repeated word by word by replacing $\genfquot$ with $\genfquotcccz$, $\Ckspquotccz{\infty}$ or $\Ckspquotccz{0}$.
\end{remark}

In Proposition \ref{comelocV} below we study the local behavior of pointed special set functions.\newline
We refer to Notation \ref{gentopnot}-[8], Definition \ref{genpointdef}, Proposition \ref{dimcoddomgerm}.
\begin{proposition}\label{comelocV}
Fix $m,n \in \mathbb{N}_0$, a special set function $\pmV :\setsymuno\rightarrow\genfquotcccz(m)+n+$, $\setsymcinque\subseteq\mathbb{R}^{m}$, $\setsymsei\subseteq\mathbb{R}^{n}$, $\setsymuno_0\subseteq \setsymuno$. Assume that
\begin{flalign}
&\begin{array}{l}
0\in \mathbb{R}^{m}\hspace{4pt}\text{is}\hspace{4pt}\text{adherent}\hspace{4pt}\text{to}\hspace{4pt}\setsymcinque\text{;}\label{ad0A}
\end{array}\\[6pt]
&\left\{
\begin{array}{l}
\text{for}\hspace{4pt}\text{any}\hspace{4pt}\elsymuno\in \setsymuno_0\hspace{4pt}\text{there}\hspace{4pt}\text{are}\hspace{4pt}\text{an}\hspace{4pt}\text{open}\hspace{4pt}\text{neighborhood}\hspace{4pt}\intuno_{\elsymuno} \subseteq\mathbb{R}^{m}\\[4pt]
\text{of}\hspace{4pt} 0\text{,}\hspace{4pt}
\text{a}\hspace{4pt}\text{generalized}\hspace{4pt}\text{function}\hspace{4pt}\gfy_{\elsymuno}\in \genfcontcont(\intuno_{\elsymuno})+\argcompl{\mathbb{R}^n}+\hspace{6pt}\text{fulfilling}\\[4pt]
\text{both}\hspace{4pt}\text{conditions}\hspace{4pt}\text{below:}\\[4pt]
\begin{array}{ll}
\pmV\left(\elsymuno\right)= \genfquotfun\left(\gfy_{\elsymuno}\right)\text{,}\\[4pt]
\left(\evalcomptwo\left(\gfy_{\elsymuno}\right)\right)\left(\unkuno\right)\in \setsymsei&\forall \unkuno\in \intuno_{\elsymuno}\cap\setsymcinque\text{.}
\end{array}
\end{array}
\right.\label{assppp}
\end{flalign}
Then there is a covering sieve $\sieve\in \csGT\left(\setsymuno\right)$ such that for any $g:\setsymdue\rightarrow\setsymuno$ belonging to $\sieve$ there are an open neighborhood $\intuno\subseteq\mathbb{R}^{m}$ of $0$ with $\setsymuno_0\cap g\left(\setsymdue\right)\neq\udenset$, $\pmtwo\in\homF+\GTopcat+\left(\setsymdue,\genfcontcont(\intuno)+\mathbb{R}^{n}+\right)$ fulfilling both conditions below:
\begin{flalign*}
&\begin{array}{l}
\pmV\funcomp g=\genfquotfun\funcomp\pmtwo\text{;}
\end{array}\\[6pt]
&\begin{array}{l}
\left(\evalcomptwo\left(\pmtwo\left(\elsymdue\right)\right)\right)\left(\unkuno\right)\in \setsymsei\hspace{15pt}\forall \left(\elsymdue,\unkuno\right)\in \setsymdue\times\intuno\text{.}
\end{array}
\end{flalign*}
\end{proposition}
\begin{proof} We prove the statement by following one by one all five steps of Definition \ref{specfunVCC}.\newline
Refer to step 1 of Definition \ref{specfunVCC}. In this case statement follows straightforwardly.\newline
Refer to step 2 of Definition \ref{specfunVCC}. Set $l=\dimcodmagtwo\left(\pmV\right)=\dimdommagtwo\left(\gfx\right)$.\newline 
Fix $\elsymuno\in \setsymuno_0$. Assumption \eqref{assppp} entails that there are an open neighborhood $\intuno_{1,\elsymuno} \subseteq\mathbb{R}^{m}$ of 0, $\gfy_{1,\elsymuno}\in \genfcontcont(\intuno_{1,\elsymuno})+\mathbb{R}^n+$ fulfilling both conditions below:
\begin{flalign}
&\begin{array}{l}
\pmV_1\left(\elsymuno\right)= \genfquotfun\left(\gfy_{1,\elsymuno}\right)\text{;}\nonumber
\end{array}\\[6pt]
&\begin{array}{l}
\left(\evalcomptwo\left(\gfy_{1,\elsymuno}\right)\right)\left(\unkuno\right)\in \setsymsei\hspace{15pt}\forall \unkuno\in \intuno_{1,\elsymuno}\cap \setsymcinque\text{.}
\end{array}\label{dentroA}
\end{flalign}
Fix an open neighborhood $\intuno_{\elsymuno}\subseteq\intuno_{1,\elsymuno}$ of $0$, $\gfy_{\elsymuno}\in \genfcontcont(\intuno_{\elsymuno})+\mathbb{R}^l+$ with $\pmV\left(\elsymuno\right)=\genfquotfun\left(\gfy_{\elsymuno}\right)$.\newline
Referring to Notation \ref{ins}-[5] set $\setsymsei_{\gfx}=\left(\evalcomptwo\left(\gfx\right)\right)^{-1}\left(\setsymsei\right)$.\newline
Since $\genfquotfun\left(\gfx \genfuncomp\gfy_{\elsymuno}\right)=\genfquotfun\left(\gfy_{1,\elsymuno}\right)$ condition \eqref{dentroA} and locality of the argument entail that there is no loss of generality by assuming
\begin{equation*}
\left(\evalcomptwo\left(\gfy_{\elsymuno}\right)\right)\left(\unkuno\right)\in \setsymsei_{\gfx}\hspace{15pt}\forall \unkuno\in \intuno_{\elsymuno}\cap \setsymcinque\text{.}
\end{equation*} 
Then induction hypothesis entail that there are an open neighborhood $\setsymdue\subseteq\setsymuno$ of $\elsymuno$, an open neighborhood $\intuno\subseteq \intuno_{\elsymuno}$ of $0$, an arrow $\pmtwo\in\homF+\GTopcat+\left(\setsymdue,\genfcontcont(\intuno)+\mathbb{R}^l+\right)$ fulfilling both conditions below:
\begin{flalign}
&\begin{array}{l}
\pmV\left(\elsymdue\right)=\genfquotfun\left(\pmtwo\left(\elsymdue\right)\right)\hspace{15pt}\forall\elsymdue\in \setsymdue\cap\setsymcinque\text{;}\label{rappY}
\end{array}\\[6pt]
&\begin{array}{l}
\left(\evalcomptwo\left(\pmtwo\left(\elsymdue\right)\right)\right)\left(\unkuno\right)\in \setsymsei_{\gfx}\hspace{15pt}\forall \left(\elsymdue,\unkuno\right)\in \setsymdue\times\intuno\text{.}\label{inAu}
\end{array}
\end{flalign}
Referring to Remark \ref{siterem}-[3] define $\pmtwo_{1}:\setsymdue\rightarrow\genfcontcont(\intuno)+\mathbb{R}^n+$ by setting $\pmtwo_{1}=\gfx \genfuncomp \pmtwo$.\newline
Eventually Proposition \ref{genpointprop}-[1], \eqref{rappY}, \eqref{inAu} entail that
\begin{flalign*}
&\begin{array}{l}
\pmV\left(\elsymdue\right)=\genfquotfun\left(\pmtwo_1\left(\elsymdue\right)\right)\hspace{15pt}\forall\elsymdue\in \setsymdue\text{;}
\end{array}\\[6pt]
&\begin{array}{l}
\left(\evalcomptwo\left(\pmtwo_1\left(\elsymdue\right)\right)\right)\left(\unkuno\right)\in \setsymsei\hspace{15pt}\forall \left(\elsymdue,\unkuno\right)\in \setsymdue\times\intuno\text{.}
\end{array}
\end{flalign*} 
Refer to step 3 of Definition \ref{specfunVCC}. Up to taking its path components one by one there is no loss of generality by assuming that $\setsymuno$ is a path connected topological space, then by Proposition \ref{domcodspstfun} we are able to set $\dimcodmagtwo\left(\pmV_2\right)=\dimdommagtwo\left(\pmV_1\right)=l$.\newline
Fix $\left(\elsymuno_1,\elsymuno_2\right)\in \setsymuno_0$. Assumption \eqref{assppp} entails that there are an open neighborhood $\intuno_{\elsymuno_1,\elsymuno_2} \subseteq\mathbb{R}^{m}$ of 0, $\gfy_{\elsymuno_1,\elsymuno_2}\in \genfcontcont(\intuno_{\elsymuno_1,\elsymuno_2})+\argcompl{\mathbb{R}^{n}}+$ fulfilling both conditions below:
\begin{flalign}
&\begin{array}{l}
\pmV\left(\elsymuno_{1},\elsymuno_{2} \right)= \genfquotfun\left(\gfy_{\elsymuno_1,\elsymuno_2}\right)\text{;}\nonumber
\end{array}\\[6pt]
&\begin{array}{l}
\left(\evalcomptwo\left(\gfy_{\elsymuno_1,\elsymuno_2}\right)\right)\left(\unkuno\right)\in \setsymsei\hspace{15pt}\forall \unkuno\in \intuno_{\elsymuno_1,\elsymuno_2}\cap\setsymcinque\text{.}
\end{array}\label{dentroA1}
\end{flalign}
Fix an open neighborhood $\intuno_{\elsymuno_1}\subseteq\mathbb{R}^{l}$ of $0 \in\mathbb{R}^{l} $, an open neighborhood $\intuno_{\elsymuno_2}\subseteq\mathbb{R}^{m}$ of $0 \in \mathbb{R}^{m}$, $\gfy_{\elsymuno_1}\in \genfcontcont(\intuno_{\elsymuno_1})+\argcompl{\mathbb{R}^{n}}+$ with $\pmV_1\left(\elsymuno_{1}\right)=\genfquotfun\left(\gfy_{\elsymuno_1}\right)$, $\gfy_{\elsymuno_2}\in \genfcontcont(\intuno_{\elsymuno_2})+\argcompl{\mathbb{R}^{l}}+$ with $\pmV_2\left(\elsymuno_{2}\right)=\genfquotfun\left(\gfy_{\elsymuno_2}\right)$.\newline 
Since 
\begin{equation}
\genfquotfun\left(\gfy_{\elsymuno_1} \genfuncomp\gfy_{\elsymuno_2}\right)=\genfquotfun\left(\gfy_{\elsymuno_1, \elsymuno_2}\right)\label{quotugu12}
\end{equation}
condition \eqref{dentroA1} entail that there is no loss of generality by assuming that for any $\left(\elsymuno_1,\elsymuno_2\right)\in \setsymuno_1\times\setsymuno_2$ 
there is an open neighborhood $\intuno_{\elsymuno_1,\elsymuno_2,2}\subseteq\intuno_{\elsymuno_1,\elsymuno_2}\cap\intuno_{\elsymuno_2}$ of $0\in\mathbb{R}^m$ fulfilling
\begin{equation}
\left(\evalcomptwo\left(\gfy_{\elsymuno_1}\right)\right)\left(\left(\evalcomptwo\left(\gfy_{\elsymuno_2}\right)\right)\left(\unkuno\right)\right)\in \setsymsei\hspace{15pt}\forall\unkuno\in \intuno_{\elsymuno_1,\elsymuno_2,2}\cap\setsymcinque\text{.}\label{dentroA12}
\end{equation}   
Referring to Notation \ref{ins}-[5] set
\begin{equation}
\setsymcinque_{\elsymuno_{2}}=\left(\evalcomptwo\left(\gfy_{\elsymuno_{2}}\right)\right)\left(\intuno_{\elsymuno_{2}}\cap\setsymcinque\right) \hspace{15pt}\forall \elsymuno_2\in \setsymuno_2\text{.} 
\end{equation}
Assumption \eqref{ad0A}, continuity of  $\evalcomptwo\left(\gfy_{\elsymuno_{2}}\right)$ entail that at least one among the two conditions below holds true:
\begin{flalign}
&\begin{array}{l}
\setsymcinque_{\elsymuno_{2}}=\left\{0\right\}\text{;}\label{adAzero}
\end{array}\\[6pt]
&\begin{array}{l}
0\in \mathbb{R}^{l}\hspace{4pt}\text{is}\hspace{4pt}\text{adherent}\hspace{4pt}\text{to}\hspace{4pt}\setsymcinque_{\elsymuno_{2}}\text{.}\label{adlA}
\end{array}
\end{flalign}
Assume that condition \eqref{adAzero} holds true for any $\elsymuno_2\in \proj<\left\{\setsymuno_1, \setsymuno_2\right\}<>2>\left(\setsymuno_0\right)$.\newline
Fix $\left(\elsymuno_{1,0},\elsymuno_{2,0}\right)\in \setsymuno_0$. Then induction hypothesis applied to $\pmV_2$ entails that there are an open neighborhood $\setsymdue_2\subseteq\setsymuno_2$ of $\elsymuno_{2,0}$, an open neighborhood $\intuno_{2}\subseteq \intuno_{\elsymuno_{2,0}}$ of $0\in  \mathbb{R}^l$,  $\pmtwo_{2}\in\homF+\GTopcat+\left(\setsymdue_2,\genfcontcont(\intuno_2)+\mathbb{R}^{n}+\right)$ fulfilling 
\begin{equation}
\left(\evalcomptwo\left(\pmtwo_2\left(\elsymdue\right)\right)\right)\left(\unkuno\right)=0\hspace{15pt}\forall \left(\elsymdue,\unkuno\right)\in \setsymdue_2\times\intuno_2\text{.}\label{S1in0}
\end{equation}
Fix an open neighborhood $\setsymdue_1\subseteq\setsymuno_1$ of $\elsymuno_{1,0}$, an open neighborhood $\intuno_{1}\subseteq \intuno_{\elsymuno_1}$ of $0\in  \mathbb{R}^l$,  $\pmtwo_{1}\in\homF+\GTopcat+\left(\setsymdue_1,\genfcontcont(\intuno_1)+\mathbb{R}^{n}+\right)$ with
\begin{equation}
\pmV\left(\elsymdue\right)=\genfquotfun\left(\pmtwo_1\left(\elsymdue\right)\right)\hspace{15pt}\forall\elsymdue\in \setsymdue_1\text{.}\label{rappglob}
\end{equation}
Define $\pmtwo:\setsymdue_1\times\setsymdue_2\rightarrow \genfcontcont(\intuno_2)+\mathbb{R}^{n}+$ by setting
\begin{equation*}
\pmtwo\left(\elsymdue_1,\elsymdue_2\right)=\pmtwo_1\left(\elsymdue_1\right) \genfuncomp\pmtwo_2\left(\elsymdue_2\right)\hspace{15pt}\forall \left(\elsymdue_{1},\elsymdue_{2}\right)\in \setsymdue_1\times\setsymdue_2\text{.}
\end{equation*}
Eventually Proposition \ref{genpointprop}-[2], \eqref{quotugu12}, \eqref{S1in0}, \eqref{rappglob} entail that
\begin{flalign*}
&\begin{array}{l}
\pmV\left(\elsymdue_1,\elsymdue_2\right)=\genfquotfun\left(\pmtwo\left(\elsymdue_1,\elsymdue_2\right)\right)\hspace{15pt}\forall\left(\elsymdue_{1},\elsymdue_{2}\right)\in \setsymdue_1\times\setsymdue_2\text{;}
\end{array}\\[6pt]
&\begin{array}{l}
\left(\evalcomptwo\left(\pmtwo\left(\elsymdue_1,\elsymdue_2\right)\right)\right)\left(\unkuno\right)=0\in \setsymsei\hspace{15pt}\forall \left(\elsymdue_1,\elsymdue_2,\unkuno\right)\in \setsymdue_1\times\setsymdue_2\times\intuno_1\text{.}
\end{array}
\end{flalign*} 
Assume that there is $\left(\elsymuno_{1,0},\elsymuno_{2,0}\right)\in \setsymuno_0$ fulfilling condition \eqref{adlA}.\newline
Conditions \eqref{dentroA12}, \eqref{adlA}, $\evalcomptwo\left(\gfy_{\elsymuno_{2,0}}\right)\in \Cksp{0}(\intuno_{\elsymuno_{2,0}})+\argcompl{\mathbb{R}^{l}}+$ entail that for any $\elsymuno_1\in \proj<\left\{\setsymuno_1, \setsymuno_2\right\}<>1>\left(\setsymuno_0\right)$ there is an open neighborhood $\intuno_{\elsymuno_1,\elsymuno_{2,0},1}\subseteq \intuno_{\elsymuno_1}$ of $0\in \mathbb{R}^l$ such that 
\begin{equation*}
\left(\evalcomptwo\left(\gfy_{\elsymuno_1}\right)\right)\left(\unkuno\right)\in \setsymsei\hspace{15pt}\forall\unkuno\in \intuno_{\elsymuno_1,\elsymuno_{2,0},1}\cap
\setsymcinque_{\elsymuno_{2,0}}
\text{.}
\end{equation*}
Then induction hypothesis applied to $\pmV_1$ entails that there are an open neighborhood $\setsymdue_1\subseteq\setsymuno_1$ of $\elsymuno_{1,0}$, an open neighborhood $\intuno_{1}\subseteq \intuno_{\elsymuno_1}$ of $0\in  \mathbb{R}^l$,  $\pmtwo_{1}\in\homF+\GTopcat+\left(\setsymdue_1,\genfcontcont(\intuno_1)+\mathbb{R}^{n}+\right)$ fulfilling 
\begin{equation}
\left(\evalcomptwo\left(\pmtwo_1\left(\elsymdue\right)\right)\right)\left(\unkuno\right)\in \setsymsei\hspace{15pt}\forall \left(\elsymdue,\unkuno\right)\in \setsymdue_1\times\intuno_1\text{.}\label{S1inB}
\end{equation}
Since is $\intuno_{1}$ an open neighborhood of $0\in  \mathbb{R}^l$ and $\evalcomptwo\left(\gfy_{\elsymuno_{2}}\right)\in \Cksp{0}(\intuno_{\elsymuno_2})+\argcompl{\mathbb{R}^{l}}+$ for any $\elsymuno_2\in \setsymuno_2$ then there is an open neighborhood $\intuno_{\elsymuno_{2},2}\subseteq \intuno_{\elsymuno_2}$ of $0\in \mathbb{R}^l$ such that 
\begin{equation*}
\left(\evalcomptwo\left(\gfy_{\elsymuno_2}\right)\right)\left(\unkuno\right)\in \intuno_{1}\hspace{15pt}\forall\unkuno\in \intuno_{\elsymuno_{2},2}
\text{.}
\end{equation*}
Then induction hypothesis applied to $\pmV_2$ entails that there are an open neighborhood $\setsymdue_2\subseteq\setsymuno_2$ of $\elsymuno_{2,0}$, an open neighborhood $\intuno_{2}\subseteq \intuno_{\elsymuno_{2,0},2}$ of $0\in  \mathbb{R}^m$,  $\pmtwo_{2}\in\homF+\GTopcat+\left(\setsymdue_2,\genfcontcont(\intuno_2)+\mathbb{R}^{l}+\right)$ fulfilling 
\begin{equation}
\left(\evalcomptwo\left(\pmtwo_2\left(\elsymdue\right)\right)\right)\left(\unkuno\right)\in \intuno_{1}\hspace{15pt}\forall \left(\elsymdue,\unkuno\right)\in \setsymdue_2\times\intuno_2\text{.}\label{S2inS1}
\end{equation}
Define $\pmtwo:\setsymdue_1\times\setsymdue_2\rightarrow \genfcontcont(\intuno_2)+\mathbb{R}^{n}+$ by setting
\begin{equation*}
\pmtwo\left(\elsymdue_1,\elsymdue_2\right)=\pmtwo_1\left(\elsymdue_1\right) \genfuncomp\pmtwo_2\left(\elsymdue_2\right)\hspace{15pt}\forall \left(\elsymdue_{1},\elsymdue_{2}\right)\in \setsymdue_1\times\setsymdue_2\text{.}
\end{equation*}
Eventually Proposition \ref{genpointprop}-[2], \eqref{quotugu12}, \eqref{S1inB}, \eqref{S2inS1} entail that
\begin{flalign*}
&\begin{array}{l}
\pmV\left(\elsymdue_1,\elsymdue_2\right)=\genfquotfun\left(\pmtwo\left(\elsymuno_1,\elsymuno_2\right)\right)\hspace{15pt}\forall\left(\elsymdue_{1},\elsymdue_{2}\right)\in \setsymdue_1\times\setsymdue_2\text{;}
\end{array}\\[6pt]
&\begin{array}{l}
\left(\evalcomptwo\left(\pmtwo\left(\elsymdue_1,\elsymdue_2\right)\right)\right)\left(\unkuno\right)\in \setsymsei\hspace{15pt}\forall \left(\elsymdue_1,\elsymdue_2,\unkuno\right)\in \setsymdue_1\times\setsymdue_2\times\intuno_1\text{.}
\end{array}
\end{flalign*} 
Refer to step 4 of Definition \ref{specfunVCC}. For any $\elindexdieci\in \left\{1,...,\topindextre\right\}$ there are $m_{\elindextre}, n_{\elindextre} \in \mathbb{N}_0$ such that $\dimdommagtwo\left(\pmV_{\elindextre}\right)=m_{\elindextre}$, $\dimcodmagtwo\left(\pmV_{\elindextre}\right)=n_{\elindextre}$.\newline
Fix $\left(\elsymuno_1,...,\elsymuno_{\topindextre}\right)\in \setsymuno_0$. Assumption \eqref{assppp} entails that there are an open neighborhood $\intuno_{\elsymuno_1,...,\elsymuno_{\topindexdieci}}\subseteq \mathbb{R}^{l_1+...+l_{\topindexdieci}}$ of $0$,
 $\gfy_{\elsymuno_1,...,\elsymuno_{\topindextre}}\in \genfcontcont(\intuno_{\elsymuno_1,...,\elsymuno_{\topindextre}})+\argcompl{\mathbb{R}^{n_1+...+n_{\topindextre}}}+$ fulfilling both conditions below:
\begin{flalign}
&\begin{array}{l}
\pmV\left(\elsymuno_{1},...,\elsymuno_{\topindextre}\right)= \genfquotfun\left(\gfy_{\elsymuno_1,...,\elsymuno_{\topindextre}}\right)\text{;}\nonumber
\end{array}\\[6pt]
&\begin{array}{l}
\left(\evalcomptwo\left(\gfy_{\elsymuno_1,...,\elsymuno_{\topindextre}}\right)\right)\left(\unkuno\right)\in \setsymsei\hspace{15pt}\forall \unkuno\in\intuno_{\elsymuno_1,...,\elsymuno_{\topindextre}}\cap\setsymcinque\text{.}
\end{array}\label{dentroAn}
\end{flalign}
Fix an open neighborhood $\intuno_{\elsymuno_{\elindextre}}\subseteq \mathbb{R}^{m_{\elindextre}}$ of $0 \in \mathbb{R}^{m_{\elindextre}}$, a generalized function
 $\gfy_{\elsymuno_{\elindextre}}\in \genfcontcont(\intuno_{\elsymuno_{\elindextre}})+\argcompl{\mathbb{R}^{n_{\elindextre}}}+$ with $\pmV_{\elindextre}\left(\elsymuno_{\elindextre}\right)=\genfquotfun\left(\gfy_{\elsymuno_{\elindextre}}\right)$ for any $\elindextre \in \left\{1,...,\topindextre\right\}$.\newline
Since $\genfquotfun\left(\lboundgenf\gfy_{\elsymuno_{1}},...,\gfy_{\elsymuno_{\topindextre}} \rboundgenf \right)=\genfquotfun\left(\gfy_{\elsymuno_1,...,\elsymuno_{\topindextre}}\right)$
condition \eqref{dentroAn} and locality of the argument entail that there is no loss of generality by assuming that both conditions below:
\begin{flalign}
&\begin{array}{l}
\overset{\topindextre}{\underset{\elindextre=1}{\prod}}\intuno_{\elsymuno_{\elindextre}}\subseteq\intuno_{\elsymuno_1,...,\elsymuno_{\topindextre}}\hspace{15pt}\forall\left(\elsymuno_1,...,\elsymuno_{\topindextre}\right)\in\overset{\topindextre}{\underset{\elindextre=1}{\prod}}\setsymuno_{\elindextre}  \text{;}\label{intdomv}
\end{array}\\[6pt]
&\begin{array}{l}
\left(\left(\evalcomptwo\left(\gfy_{\elsymuno_{1}}\right)\right)\left(\unkuno_1\right),...,\left(\evalcomptwo\left(\gfy_{\elsymuno_{\topindextre}}\right)\right)\left(\unkuno_{\topindextre}\right)\right)\in \setsymsei\\[4pt]
\hspace{150pt}
\forall \left( \unkuno_1,...,\unkuno_{\topindextre}\right)\in \left(\overset{\topindextre}{\underset{\elindextre=1}{\prod}}\intuno_{\elsymuno_{\elindextre}}\right)\cap\setsymcinque  \text{.}
\end{array}\label{prodAinA}
\end{flalign}
Referring to Notation \ref{ins}-[5, 10] set 
\begin{equation*}
\setsymsei_{\elindextre}=\underset{\elsymuno_{\elindextre}\in \setsymuno_{\elindextre}}{\bigcup}\left(\evalcomptwo\left(\gfy_{\elsymuno_{\elindextre}}\right)\right)\left(\intuno_{\elsymuno_{\elindextre}}\cap  \proj<\left\{\mathbb{R}^{m_{\elindextre}}\right\}_{\elindextre=1}^{\topindextre}<>\elindextre>\left(\setsymcinque\right)\right)\hspace{15pt}\forall\elindextre \in \left\{1,...,\topindextre\right\}\text{.}
\end{equation*}
Condition \eqref{prodAinA} entail that  
\begin{equation*}
\overset{\topindextre}{\underset{\elindextre=1}{\prod}}\setsymsei_{\elindextre}\subseteq\setsymsei\text{.} 
\end{equation*}
Induction hypothesis entail that there are 
an open neighborhood $\setsymdue_{\elindextre}\subseteq\setsymuno_{\elindextre}$ of $\elsymuno_{\elindextre}$, an open neighborhood $\intuno_{\elindextre}\subseteq \intuno_{\elsymuno_{\elindextre}}$ of $0\in \intuno_{\elsymuno_{\elindextre}}$,
an arrow $\pmtwo_{\elindextre}\in\homF+\GTopcat+\left(\setsymdue_{\elindextre},\genfcontcont(\intuno_{\elindextre})+\argcompl{\mathbb{R}^{m_{\elindextre}}}+\right)$ for any $\elindextre \in \left\{1,...,\topindextre\right\}$ fulfilling all conditions below:
\begin{flalign*}
&\begin{array}{l}
\pmV_{\elindextre}\left(\elsymdue\right)=\genfquotfun\left(\pmtwo_{\elindextre}\left(\elsymdue\right)\right)\hspace{4pt}\text{for}\hspace{4pt}\text{any}\hspace{4pt}\elsymdue\in \setsymdue_{\elindextre}\text{,}\hspace{4pt}\elindextre \in \left\{1,...,\topindextre\right\}\text{;}
\end{array}\\[6pt]
&\left\{
\begin{array}{l}
\left(\evalcomptwo\left(\pmtwo_{\elindextre}\left(\elsymdue\right)\right)\right)\left(\unkuno\right)\in \setsymsei_{\elindextre}\\[4pt]
\text{for}\hspace{4pt}\text{any}\hspace{4pt} \left(\elsymdue,\unkuno\right)\in \setsymdue_{\elindextre}\times\left(\intuno_{\elindextre}\cap  \proj<\left\{\mathbb{R}^{m_{\elindextre}}\right\}_{\elindextre=1}^{\topindextre}<>\elindextre>\left(\setsymcinque\right) \right)\text{,}\hspace{4pt}\elindextre \in \left\{1,...,\topindextre\right\}\text{.}
\end{array}
\right.
\end{flalign*}
Referring to Remark \ref{siterem}-[3] define $\pmtwo:\overset{\topindextre}{\underset{\elindextre=1}{\prod}}\setsymdue_{\elindextre}\rightarrow\genfcontcont(\argcompl{\overset{\topindextre}{\underset{\elindextre=1}{\prod}}\intuno_{\elindextre}})+\argcompl{\mathbb{R}^{n_1+...+n_{\topindextre}}}+$ by setting $\pmtwo=\lboundgenf\pmtwo_1,...,\pmtwo_{\topindextre}\rboundgenf$.\newline
Eventually Proposition \ref{genpointprop}-[3], \eqref{dentroAn}, \eqref{intdomv},\eqref{prodAinA} entail that
\begin{flalign*}
&\begin{array}{l}
\pmV\left(\elsymdue_{1},...,\elsymdue_{\topindextre}\right)=\genfquotfun\left(\pmtwo\left(\elsymdue_{1},...,\elsymdue_{\topindextre}\right)\right)\hspace{15pt}\forall\left(\elsymdue_{1},...,\elsymdue_{\topindextre}\right) \in \overset{\topindextre}{\underset{\elindextre=1}{\prod}}\setsymdue_{\elindextre}\text{;}
\end{array}\\[6pt]
&\begin{array}{l}
\left(\evalcomptwo\left(\pmtwo\left(\elsymdue_{1},...,\elsymdue_{\topindextre}\right)\right)\right)\left(\unkuno_{1},...,\unkuno_{\topindextre}\right)\in \setsymsei\hspace{15pt}\forall \left(\elsymdue_{1},...,\elsymdue_{\topindextre},\unkuno_{1},...,\unkuno_{\topindextre}\right)\in \overset{\topindextre}{\underset{\elindextre=1}{\prod}}\setsymdue_{\elindextre}\times\overset{\topindextre}{\underset{\elindextre=1}{\prod}}\intuno_{\elindextre}\text{.}
\end{array}
\end{flalign*}
\end{proof}

In Corollary \ref{zeroloc} below we list a straightforward consequence of Proposition \ref{comelocV}.

\begin{corollary}\label{zeroloc}
Fix $m,n\in \mathbb{N}_0$, $\pmV \in\homF+\GTopcat+\left(\setsymuno,\genfquotcccz(m)+n+\right)$, $\elsymuno_0 \in \setsymuno$.\newline
Assume that 
\begin{equation}
\evalcompquotcontcontzero\left(\pmV\left(\elsymuno_0\right)\right)=\gczsfx[\argcompl{\cost<\mathbb{R}^m<>\mathbb{R}^n>+0+}]\text{.}
\end{equation}
Then there is a covering sieve $\sieve\in \csGT\left(\setsymuno\right)$ such that for any $g:\setsymdue\rightarrow\setsymuno$ belonging to $\sieve$ with $\elsymuno_0 \in g\left(\setsymdue\right)$ there are an open neighborhood $\intuno\subseteq\mathbb{R}^m$ of $0$, $\pmtwo\in\homF+\GTopcat+\left(\setsymdue,\genfcontcont(\intuno)+\mathbb{R}^n+\right)$ fulfilling both conditions below:
\begin{flalign}
&\begin{array}{l}
\pmV\funcomp g=\genfquotfun\funcomp\pmtwo\text{;}
\end{array}\\[6pt]
&\begin{array}{l}
\evalcomptwo\left(\pmtwo\left(\elsymdue\right)\right)=\cost<\intuno<>\mathbb{R}^n>+0+\text{.}
\end{array}
\end{flalign}
\end{corollary}

In Definition \ref{alssitesdef}, Remark \ref{alssitesrem} below we introduce categories of real vector spaces, algebras and related continuous functions existing in the setting of the site $\left(\GTopcat, \csGT\right)$. 
We refer to Notation \ref{alg}-[10-15], Definition \ref{defcatD}-[5], Remarks \ref{siterem}, \ref{siteremdue}.

\begin{definition}\mbox{}\label{alssitesdef}
\begin{enumerate}
\item We denote by $\GVecR$ the subcategory of $\GTopcat$ whose objects are all $\mathbb{R}$-vector spaces in $\GTopcat$ and arrows are $\mathbb{R}$-linear arrows in $\GTopcat$.
\item We denote by $\GAlgR$ the subcategory of $\GTopcat$ whose objects are all $\mathbb{R}$-algebras in $\GTopcat$ and arrows are $\mathbb{R}$-algebra arrows in $\GTopcat$.
\end{enumerate}
\end{definition}

\begin{remark}\mbox{}\label{alssitesrem}
\begin{enumerate}
\item Categories $\GVecR$, $\GAlgR$ admits all limits, co-limits and adjunctions existing in $\VecR$, $\AlgR$ by definition of site $\left(\GTopcat,\csGT\right)$ and of categories $\Catquattro$, $\Catcinque$ in Definition \ref{defcatD}. With an abuse of language we denote limits, co-limits and adjunctions in $\GVecR$, $\GAlgR$ by the same symbols used in categories $\VecR$, $\AlgR$. 
\item Functors $\tenalgfun$, $\extalgfun$, $\symalgfun$, $\Aug$ introduced in Notation \ref{alg}-[12, 13, 14, 16] factor through forgetful functors $ \GVecR\rightarrow \VecR$, $ \GAlgR\rightarrow \AlgR$. With an abuse of language we denote such factorization again by $\tenalgfun$, $\extalgfun$, $\symalgfun$, $\Aug$ respectively. 
\end{enumerate}
\end{remark}

\section{Generalized derivations of germs of generalized functions\label{dggf}}

In Definition \ref{coeqdef} below we introduce notion of generalized $\left(\setsymuno,\modsym,l\right)$-derivation of degree $1$.  A $\left(\setsymuno, \modsym, l\right)$-derivation of degree $1$ is an equivalence class of $\left(\modsym,l\right)$-pre-derivations valued arrows of $\GTopcat$ with domain $\setsymuno$. Construction of $\left(\setsymuno, \modsym, l\right)$-derivations will follow by organizing material introduced in Definition \ref{genderdef2}-[6, 7] in the functorial language, then taking co-limits. If $\setsymuno$ is the terminal object of $\GTopcat$, i.e. the single point, and $\modsym=\mathbb{R}^l$ then we obtain the notion of generalized $l$-derivation of degree $1$ extending the classical notion of derivation of degree $1$. We refer to Notation \ref{catfun}-[1, 4], Definition \ref{genpointdefmnCC0loc}, Proposition \ref{genpointprop}-[11].

\begin{definition}\label{coeqdef}\mbox{}
\begin{enumerate}
\item Fix $m \in \mathbb{N}_0$. We denote by $\retezero+m+$ the category defined as follows.\newline
Objects of $\retezero+m+$ are listed below:\newline
symbol $\terobj+m+$;\newline
triplets $\left(k, \gggfx[\gfx_1], \gggfx[\gfx_2]\right)$ for any $k \in\mathbb{N}_0$, $\gggfx[\gfx_1], \gggfx[\gfx_2]\in\genfquotcccz(k)+m+$ with $\evalcompquotcontcontzero\left(\gggfx[\gfx_1]\right)=\evalcompquotcontcontzero\left(\gggfx[\gfx_2]\right)$.\newline
Arrows of $\retezero+m+$ are listed below:\newline
for any triplet $\left(k, \gggfx[\gfx_1], \gggfx[\gfx_2]\right)$ there are exactly two arrows from $\left(k, \gggfx[\gfx_1],\gggfx[\gfx_2]\right)$ to $\terobj+m+$, we denote such arrows by $\arrtz\left[k, \gggfx[\gfx_1],\gggfx[\gfx_2]\right]_i$ for $i \in \left\{1,2\right\}$;\newline
for any object $\objrtz$ of $\retezero+m+$ the identity $\idobj+\objrtz+$ is the only arrow from $\objrtz$ to $\objrtz$;\newline 
for any pair of objects $\objrtz_1$, $\objrtz_2$ of $\retezero+m+$ there is no arrow from $\objrtz_1$ to $\objrtz_2$ unless either $\objrtz_1=\objrtz_2$ or $\objrtz_2=\terobj+m+$.
\item Fix $m, n \in \mathbb{N}_0$,
$\gggfx[\gfx] \in \genfquotcccz(m)+n+$. We denote by $\transterezero+\gggfx[\gfx]+:\retezero+m+\rightarrow\retezero+n+$ the functor defined by setting:
\begin{flalign*}
&\left\{
\begin{array}{l}
\transterezero+\gggfx[\gfx]+\left(k, \gggfx[\gfx_1], \gggfx[\gfx_2]\right)=\left(k, \gggfx[\gfx]\compquotpoint\gggfx[\gfx_1], \gggfx[\gfx]\compquotpoint\gggfx[\gfx_2]\right)\text{,}\\[4pt]
\transterezero+\gggfx[\gfx]+\left(\terobj+m+\right)=\terobj+n+\text{;}
\end{array}
\right.\\[8pt]
&\left\{
\begin{array}{l}
\transterezero+\gggfx[\gfx]+\left(\arrtz\left[k, \gggfx[\gfx_1], \gggfx[\gfx_2]\right]_i\right)=\arrtz\left[k, \gggfx[\gfx]\compquotpoint\gggfx[\gfx_1], \gggfx[\gfx]\compquotpoint\gggfx[\gfx_2]\right]_i\text{,}\\[4pt]
\transterezero+\gggfx[\gfx]+\left(\idobj+\objrtz+\right)=\idobj+\argcompl{\transterezero+\gggfx[\gfx]+\left(\objrtz\right)}+\text{.}
\end{array}
\right.
\end{flalign*}
\item Fix $l,m \in \mathbb{N}_0$, a $\mathbb{R}$-vector subspace $\modsym$ of $\mathbb{R}^l$.\newline
We denote by $\derrtzfun[\modsym](l)+m+: \retezero+m+\rightarrow \VecR$ the functor defined by setting: 
\begin{flalign*}
&\left\{
\begin{array}{ll}
\derrtzfun[\modsym](l)+m+\left(k, \gggfx[\gfx_1], \gggfx[\gfx_2]\right)=\genprederspace[\modsym](l)+k+\text{,}\\[4pt]
\derrtzfun[\modsym](l)+m+\left(\terobj+m+\right)=\genprederspace[\modsym](l)+m+\text{;}
\end{array}
\right.\\[8pt]
&\left\{
\begin{array}{l}
\derrtzfun[\modsym](l)+m+\left(\arrtz\left[k, \gggfx[\gfx_1], \gggfx[\gfx_2]\right]_i\right)=\genprediff[\modsym](l)+\gggfx[\gfx_i]+\text{,}\\[4pt]
\derrtzfun[\modsym](l)+m+\left(\idobj+\objrtz+
\right)=\idobj+\argcompl{\derrtzfun[\modsym](l)+m+\left(\objrtz\right)}+\text{.}
\end{array}
\right.
\end{flalign*}
\item Fix $l,m,n\in\mathbb{N}_0$, a $\mathbb{R}$-vector subspace $\modsym$ of $\mathbb{R}^l$, $\gggfx[\gfx] \in \genfquotcccz(m)+n+$.\newline
We denote by $\generdiffNT[\modsym](l)+\gggfx[\gfx]+:\derrtzfun[\modsym](l)+m+\rightarrow \left(\derrtzfun[\modsym](l)+n+\comparrfun\transterezero+\gggfx[\gfx]+\right)$ the natural arrow defined by setting:
\begin{equation*}
\begin{array}{l}
\generdiffNT[\modsym](l)+\gggfx[\gfx]+\left(k, \gggfx[\gfx_1], \gggfx[\gfx_2]\right) =\idobj+\argcompl{\genprederspace[\modsym](l)+k+}+\text{;}\\[4pt]
\generdiffNT[\modsym](l)+\gggfx[\gfx]+\left(\terobj+m+\right)=\genprediff[\modsym](l)+\gggfx[\gfx]+\text{.}
\end{array}
\end{equation*}
\item Fix $l,m\in\mathbb{N}_0$, two $\mathbb{R}$-vector subspaces $\modsym_1$, $\modsym_2$ of $\mathbb{R}^l$. Assume $\modsym_1 \subseteq\modsym_2$.\newline
We denote by $\generinclNT<\modsym_1<>\modsym_2>(l)+m+: \derrtzfun[\modsym_1](l)+m+\rightarrow \derrtzfun[\modsym_2](l)+m+$ the natural arrow defined by setting
\begin{equation*}
\begin{array}{l}
\generinclNT<\modsym_1<>\modsym_2>(l)+m+\left(k, \gggfx[\gfx_1], \gggfx[\gfx_2]\right)=\genpreincl<\modsym_1<>\modsym_2>(l)+k+\text{;}\\[4pt]
\generinclNT<\modsym_1<>\modsym_2>(l)+m+\left(\terobj+m+\right)=\genpreincl<\modsym_1<>\modsym_2>(l)+m+\text{.}
\end{array}
\end{equation*}
\item Fix $l,m\in \mathbb{N}_0$, a $\mathbb{R}$-vector subspace $\modsym$ of $\mathbb{R}^l$.\newline
We define the pair $\left(\generderspace[\modsym](l)+m+, \quotpreder[\modsym](l)+m+\right)$ as the co-limit of the functor $\derrtzfun[\modsym](l)+m+$. We say that $\generderspace[\modsym](l)+m+$ is the space of $m$-dimensional generalized $\left(\modsym, l\right)$-derivations of degree $1$, or $m$-dimensional generalized $\left(\modsym, l\right)$-tangent vectors. Elements belonging to $\generderspace[\modsym](l)+m+$ are equivalence classes denoted by $\generder$ or by $\generder[\genpreder]$ whenever we need to emphasize a representative $\genpreder$ of the equivalence class.
We say that $\generder$ is a $m$-dimensional generalized $\left(\modsym,l\right)$-derivation of degree $1$, or a $m$-dimensional generalized $\left(\modsym,l\right)$-tangent vector. We drop any reference to $m$, $\modsym$ or $l$ whenever no confusion is possible.
\item Fix $l,m, n \in \mathbb{N}_0$, a $\mathbb{R}$-vector subspace $\modsym$ of $\mathbb{R}^l$, $\gggfx \in \genfquotcccz(m)+n+$.\newline
We define the arrow $\gendiff[\modsym](l)+\gggfx+ : \generderspace[\modsym](l)+m+\rightarrow \generderspace[\modsym](l)+n+$ as the unique arrow induced on co-limits by the natural arrow $\generdiffNT[\modsym](l)+\gggfx+$. We say that $\gendiff[\modsym](l)+\gggfx+$ is the generalized $\left(\modsym, l\right)$-differential of degree $1$ of $\gggfx$ at $0$. 
\item Fix $l,m\in\mathbb{N}_0$, two $\mathbb{R}$-vector subspaces $\modsym_1$, $\modsym_2$ of $\mathbb{R}^l$. Assume $\modsym_1 \subseteq\modsym_2$.\newline
We define the arrow $\generincl<\modsym_1<>\modsym_2>(l)+m+: \generderspace[\modsym_1](l)+m+\rightarrow \generderspace[\modsym_2](l)+m+$ as the unique arrow induced on co-limits by the natural arrow $\generinclNT<\modsym_1<>\modsym_2>(l)+m+$.
\end{enumerate}
\end{definition}

In Remark \ref{pathinlgenerder} below we explicitly describe the structure of the co-limit pair $\left(\generderspace[\modsym](l)+m+, \quotpreder[\modsym](l)+m+\right)$ introduced in Definition \ref{coeqdef}. We refer to Notation \ref{catfun}-[9, 10], Proposition \ref{genpointprop}-[11], Remarks \ref{genercontfun}, \ref{existoperquot}-[2], \ref{extremAnCC0}-[2], \ref{pathinlpreder}, \ref{siterem}-[3, 7], \eqref{spfunV}, \eqref{spfunVCC}.

\begin{remark}\label{pathinlgenerder}\mbox{}
\begin{enumerate}
\item Fix $l,m\in \mathbb{N}_0$,  a $\mathbb{R}$-vector subspace $\modsym$ of $\mathbb{R}^l$.\newline
Then $\generderspace[\modsym](l)+m+$ is a $\mathbb{R}$-vector space in $\GTopcat$ which is not trivial by Proposition \ref{algstrderpoint} below.\newline
Sum and scalar product of $\generderspace[\modsym](l)+m+$ are denoted respectively by $\sumgenerder[\modsym](l)+m+$ and $\scalpgenerder[\modsym](l)+m+$. With an abuse of notation we denote operation $\sumgenerder[\modsym](l)+m+$ simply by $\sumgenerder$, $\scalpgenerder[\modsym](l)+m+$
simply by $\scalpgenerder$, whenever no confusion is possible.
\item Fix $l,m\in \mathbb{N}_0$,  a $\mathbb{R}$-vector subspace $\modsym$ of $\mathbb{R}^l$.\newline
Then $\generderspace[\modsym](l)+m+=\frac{\genprederspace[\modsym](l)+m+}{\nullpreder[\modsym](l)+m+}$ where $\nullpreder[\modsym](l)+m+$ is the $\mathbb{R}$-vector subspace of $\genprederspace[\modsym](l)+m+$ generated by  the set 
\begin{multline}\label{quotgenset}
\gennullpreder[\modsym](l)+m+=\big\{
\genprediff[\modsym](l)+\gggfx[\gfx_{1}]+ \left(\genpreder\right)\,\sumgenpreder\, -1\scalpgenpreder
\genprediff[\modsym](l)+\gggfx[\gfx_{2}]+\left(\genpreder\right)\;:\\
\genpreder \in \genprederspace[\modsym](l)+k+\text{,}\;
 \gggfx[\gfx_{1}], \gggfx[\gfx_{2}]\in \genfquotcccz(k)+m+\;\text{with}\\
 \evalcompquotcontcontzero\left(\gggfx[\gfx_1]\right)=\evalcompquotcontcontzero\left(\gggfx[\gfx_2]\right)
\big\}\text{.}
\end{multline}
We denote by $\quotdirsumdue[\modsym](l)+m+ :\genprederspace[\modsym](l)+m+\rightarrow \generderspace[\modsym](l)+m+$ the quotient function.\newline
We denote by $\nullprederunodue[\modsym](l)+m+$ the kernel of $\quotdirsumdue[\modsym](l)+m+\funcomp\quotdirsumuno[\modsym](l)+m+$.\newline
We drop any reference to $\modsym$, $l$ or $m$ whenever there is no risk of confusion or there is no need of such a detail.
\item Fix $l,m\in \mathbb{N}_0$,  a $\mathbb{R}$-vector subspace $\modsym$ of $\mathbb{R}^l$. We set
\begin{multline}
\ssf{\generderspace[\modsym](l)+m+}=\\
\bigg\{\quotdirsumdue[\modsym](l)+m+\funcomp f\;:\;f\in\ssf{\genprederspace[\modsym](l)+m+}\bigg\}\text{.}
\label{specpredersieve1}
\end{multline}
We say that elements belonging to set $\ssf{\generderspace[\modsym](l)+m+}$ are special set functions.\newline
By \eqref{bGt7bis}, \eqref{pcM} and referring to Notation \ref{catfun}-[9] we have that set $\ssf{\generderspace[\modsym](l)+m+}$ is a base component set for $\csGT\left(\generderspace[\modsym](l)+m+\right)$.
\item Fix $l,m\in \mathbb{N}_0$,  a $\mathbb{R}$-vector subspace $\modsym$ of $\mathbb{R}^l$.\newline
Natural arrow $\quotpreder[\modsym](l)+m+:\derrtzfun[\modsym](l)+m+\rightarrow \costfun<\argcompl{\retezero+m+}<>\VecR>+\argcompl{\generderspace[\modsym](l)+m+}+$ consists of linear function $\quotpreder[\modsym](l)+m+\left(\objrtz\right): \derrtzfun[\modsym](l)+m+\left(\objrtz\right)\rightarrow \generderspace[\modsym](l)+m+$ for any object $\objrtz$ of $\retezero+m+$.
\item Fix $l,m\in \mathbb{N}_0$,  a $\mathbb{R}$-vector subspace $\modsym$ of $\mathbb{R}^l$.\newline
The $\mathbb{R}$-linear function $\gendiff[\modsym](l)+\gggfx[\gfx]+$ is the unique linear function induced on quotient spaces by the $\mathbb{R}$-linear function $\genprediff[\modsym](l)+\gggfx[\gfx]+$.
\item Fix $l,m,n\in \mathbb{N}_0$, a $\mathbb{R}$-vector subspace $\modsym$ of $\mathbb{R}^l$, $\gggfx[\gfx_1], \gggfx[\gfx_2] \in \genfquotcccz(m)+n+$. If $\evalcompquotcontcontzero\left(\gggfx[\gfx_1]\right)=\evalcompquotcontcontzero\left(\gggfx[\gfx_2]\right)$ then straightforwardly $\gendiff[\modsym](l)+\gggfx[\gfx_1]+=\gendiff[\modsym](l)+\gggfx[\gfx_2]+$. \newline
We emphasize that:
\begin{flalign}
&\left\{
\begin{array}{l}
\text{equality}\hspace{4pt}\gendiff[\modsym](l)+\argcompl{\gggfx \sumquotpoint \gggfy}+ = \gendiff[\modsym](l)+\gggfx+  \sumgenerder  \gendiff[\modsym](l)+\gggfy+\hspace{4pt}\text{is}\hspace{4pt}\text{not}\\[4pt]
\text{fulfilled}\hspace{4pt}\text{for}\hspace{4pt}\text{any}\hspace{4pt}\gggfx, \gggfy \in \genfquotcccz(m)+n+\text{;}
\end{array}
\right.\\[8pt]
&\left\{
\begin{array}{l}
\text{equality}\hspace{4pt}\gendiff[\modsym](l)+\argcompl{a \scalpquotpoint \gggfx}+ = a  \scalpgenerder  \gendiff[\modsym](l)+\gggfx+\hspace{4pt}\text{is}\hspace{4pt}\text{not}\hspace{4pt}\text{fulfilled}\\[4pt]
\text{for}\hspace{4pt}\text{any}\hspace{4pt}a\in \mathbb{R}\text{,}\hspace{4pt}\gggfx\in \genfquotcccz(m)+n+\text{.}
\end{array}
\right.
\end{flalign}
Fix $l,m,n\in \mathbb{N}_0$, a $\mathbb{R}$-vector subspace $\modsym$ of $\mathbb{R}^l$, $ \gczsfx[f] \in\Ckspquotccz{0}(m)+n+$. We denote by $\gendiff[\modsym](l)+\gczsfx[f]+$ the natural arrow $\gendiff[\modsym](l)+\gggfx[\gfx]+$ where $\gggfx[\gfx] \in \genfquotcccz(m)+n+$ fulfills $\evalcompquotcontcontzero\left(\gggfx[\gfx]\right)=\gczsfx[f]$.
\item Fix $l,m\in \mathbb{N}_0$, a $\mathbb{R}$-vector subspace $\modsym$ of $\mathbb{R}^l$. Then $\generincl<\modsym<>\modsym>(l)+m+=\idobj+\argcompl{\generderspace[\modsym](l)+m+}+$.\newline
Fix $l,m\in \mathbb{N}_0$, three $\mathbb{R}$-vector subspaces $\modsym_1$, $\modsym_2$, $\modsym_3$ of $\mathbb{R}^l$ with $\modsym_1\subseteq\modsym_2\subseteq\modsym_3$. Then $\generincl<\modsym_2<>\modsym_3>(l)+m+\funcomp\generincl<\modsym_1<>\modsym_2>(l)+m+=\generincl<\modsym_1<>\modsym_3>(l)+m+$.\newline
Fix $l,m\in \mathbb{N}_0$, two $\mathbb{R}$-vector subspaces $\modsym_1$, $\modsym_2$ of $\mathbb{R}$-vector subspaces of $\mathbb{R}^l$, $\gggfx[\gfx] \in \genfquotcccz(m)+n+$.\newline
Then $\gendiff[\modsym_2](l)+\gggfx[\gfx]+\funcomp \generincl<\modsym_1<>\modsym_2>(l)+m+=\generincl<\modsym_1<>\modsym_2>(l)+n+ \funcomp \gendiff[\modsym_1](l)+\gggfx[\gfx]+$.
\item Fix $l,m\in \mathbb{N}_0$,  a $\mathbb{R}$-vector subspace $\modsym$ of $\mathbb{R}^l$, an object $\setsymuno$ of $\GTopcat$, an arrow $\pmpreder\in \homF+\GTopcat+\left(\setsymuno,\genprederspace[\modsym](l)+m+\right)$, a covering sieve $\sieve\in \csGT\left(\setsymuno\right)$ with $\sieve \subseteq\pbsieve[\pmpreder](\argcompl{\gensieve{\ssf{\genprederspace[\modsym](l)+m+}}})$.\newline
Then for any arrow $\setsymdue\overset{g}{\rightarrow}\setsymuno$ belonging to $\sieve$ there are a finite set $\setindexdue_{g}$, a pair $\left( \pmV_{\elindexdue},f_{\elindexdue}\right)\in\left(\homF+\GTopcat+\left(\setsymdue,\genfquotcccz(l)+m+\right)\cap\gensieve{\ssf{\genfquot}}
\right)\times\homF+\GTopcat+\left(\setsymdue,\modsym\right)$ for any $\elindexdue\in\setindexdue_{g}$ such that $\quotdirsumuno[\modsym](l)+m+\funcomp\underset{\elindexdue\in  \setindexdue_g}{\sum} \genpreder\left[\pmV_{\elindexdue},f_{\elindexdue} \right]=\pmpreder\funcomp g$ for any arrow $\setsymdue\overset{g}{\rightarrow}\setsymuno$ belonging to $\sieve$.\newline  
We set:
\begin{flalign*}
&\begin{array}{l}
\setindexdue_{\sieve}=\left\{\setindexdue_{g}\right\}_{g \in \sieve}\text{;}
\end{array}\\[6pt]
&\left\{
\begin{array}{l}
\left( \pmV_{\setindexdue_{g}},f_{\setindexdue_{g}}\right)=\left\{\left( \pmV_{\elindexdue},f_{\elindexdue}\right)\right\}_{\elindexdue\in  \setindexdue_{g}}\hspace{15pt}\forall g \in \sieve\text{,}\\[4pt]
\left( \pmV_{\setindexdue_{\sieve}},f_{\setindexdue_{\sieve}}\right)=\left\{\left( \pmV_{\setindexdue_{g}},f_{\setindexdue_{g}}\right)\right\}_{g\in  \sieve}\text{;}
\end{array}
\right.\\[6pt]
&\left\{
\begin{array}{l}
\genpreder\left[\pmV_{\setindexdue_{g}},f_{\setindexdue_{g}} \right]=\left\{\genpreder\left[\pmV_{\elindexdue},f_{\elindexdue} \right]\right\}_{\elindexdue\in  \setindexdue_{g}}\hspace{15pt}\forall g \in \sieve\text{,}\\[4pt]
\genpreder\left[\pmV_{\setindexdue_{\sieve}},f_{\setindexdue_{\sieve}} \right]=\left\{\genpreder\left[\pmV_{\setindexdue_{g}},f_{\setindexdue_{g}}  \right]\right\}_{g \in \sieve}\text{;}
\end{array}
\right.\\[6pt]
&\left\{
\begin{array}{l}
\pmpreder_{\elindexdue}=\quotdirsumuno[\modsym](l)+m+\funcomp\genpreder\left[\pmV_{\elindexdue},f_{\elindexdue} \right]\hspace{15pt}\forall \elindexdue \in  \setindexdue_{g}\text{,}\\[4pt]
\pmpreder_{\setindexdue_{g}}=\left\{\pmpreder_{\elindexdue}\right\}_{\elindexdue\in  \setindexdue_{g}}\hspace{15pt}\forall \elindexdue \in  \setindexdue_{g}\text{,}\\[4pt]
\pmpreder_{\setindexdue_{\sieve}}=\left\{\pmpreder_{\setindexdue_{g}}\right\}_{g\in  \sieve}\text{.}
\end{array}
\right.
\end{flalign*}
We say that:\newline
$\left( \pmV_{\setindexdue_{g}},f_{\setindexdue_{g}}\right)$ is a skeleton of $\pmpreder\funcomp g$ for any $g\in \sieve$;\newline
$\left( \pmV_{\setindexdue_{\sieve}},f_{\setindexdue_{\sieve}}\right)$ is a skeleton of $\pmpreder$ relative to the covering sieve $\sieve$;\newline
$\genpreder\left[\pmV_{\setindexdue_{g}},f_{\setindexdue_{g}} \right]$ is a realization of $\pmpreder\funcomp g$ for any $g\in \sieve$;\newline
$\genpreder\left[\pmV_{\setindexdue_{\sieve}},f_{\setindexdue_{\sieve}} \right]$ is a realization of $\pmpreder$ relative to the covering sieve $\sieve$;\newline
$\pmpreder_{\setindexdue_{g}}$ is a decomposition of $\pmpreder\funcomp g$ for any $g\in \sieve$;\newline
$\pmpreder_{\setindexdue_{\sieve}}$ is a decomposition of $\pmpreder$ relative to the covering sieve $\sieve$.\newline
If $\sieve$ is the maximal sieve then we only consider skeleton, realization and decomposition relative to arrow $\idobj+\setsymuno+\in \sieve$ and we drop any reference to both $\idobj+\setsymuno+$ and $\sieve$ in the notation.
\item Fix $l,m\in \mathbb{N}_0$,  a $\mathbb{R}$-vector subspace $\modsym$ of $\mathbb{R}^l$, an object $\setsymuno$ of $\GTopcat$, an arrow $\pmgenerder\in \homF+\GTopcat+\left(\setsymuno,\generderspace[\modsym](l)+m+\right)$, a covering sieve $\sieve\in \csGT\left(\setsymuno\right)$ with $\sieve \subseteq\pbsieve[\pmgenerder](\argcompl{\gensieve{\ssf{\generderspace[\modsym](l)+m+}}})$.\newline
Then for any arrow $\setsymdue\overset{g}{\rightarrow}\setsymuno$ belonging to $\sieve$ there are a finite set $\setindexdue_{g}$, a pair $\left( \pmV_{\elindexdue},f_{\elindexdue}\right)\in\left(\homF+\GTopcat+\left(\setsymdue,\genfquotcccz(l)+m+\right)\cap\gensieve{\ssf{\genfquot}}
\right)\times\modsym$ for any $\elindexdue\in\setindexdue_{g}$ such that $
\quotdirsumdue[\modsym](l)+m+\funcomp\left(\quotdirsumuno[\modsym](l)+m+\funcomp\left(\underset{\elindexdue\in  \setindexdue_g}{\sum} \genpreder\left[\pmV_{\elindexdue},f_{\elindexdue} \right]\right)\right)=\pmgenerder\funcomp g$ for any arrow $\setsymdue\overset{g}{\rightarrow}\setsymuno$ belonging to $\sieve$.\newline 
We set:
\begin{flalign*}
&\begin{array}{l}
\setindexdue_{\sieve}=\left\{\setindexdue_{g}\right\}_{g \in \sieve}\text{;}
\end{array}\\[6pt]
&\left\{
\begin{array}{l}
\left( \pmV_{\setindexdue_{g}},f_{\setindexdue_{g}}\right)=\left\{\left( \pmV_{\elindexdue},f_{\elindexdue}\right)\right\}_{\elindexdue\in  \setindexdue_{g}}\hspace{15pt}\forall g \in \sieve\text{,}\\[4pt]
\left( \pmV_{\setindexdue_{\sieve}},f_{\setindexdue_{\sieve}}\right)=\left\{\left( \pmV_{\setindexdue_{g}},f_{\setindexdue_{g}}\right)\right\}_{g\in  \sieve}\text{;}
\end{array}
\right.\\[6pt]
&\left\{
\begin{array}{l}
\genpreder\left[\pmV_{\setindexdue_{g}},f_{\setindexdue_{g}} \right]=\left\{\genpreder\left[\pmV_{\elindexdue},f_{\elindexdue} \right]\right\}_{\elindexdue\in  \setindexdue_{g}}\hspace{15pt}\forall g \in \sieve\text{,}\\[4pt]
\genpreder\left[\pmV_{\setindexdue_{\sieve}},f_{\setindexdue_{\sieve}} \right]=\left\{\genpreder\left[\pmV_{\setindexdue_{g}},f_{\setindexdue_{g}}  \right]\right\}_{g \in \sieve}\text{;}
\end{array}
\right.\\[6pt]
&\left\{
\begin{array}{l}
\pmpreder_{\elindexdue}=\quotdirsumuno[\modsym](l)+m+\funcomp\genpreder\left[\pmV_{\elindexdue},f_{\elindexdue} \right]\hspace{15pt}\forall \elindexdue \in  \setindexdue_{g}\text{,}\\[4pt]
\pmpreder_{\setindexdue_{g}}=\left\{\pmpreder_{\elindexdue}\right\}_{\elindexdue\in  \setindexdue_{g}}\hspace{15pt}\forall \elindexdue \in  \setindexdue_{g}\text{,}\\[4pt]
\pmpreder_{\setindexdue_{\sieve}}=\left\{\pmpreder_{\setindexdue_{g}}\right\}_{g\in  \sieve}\text{.}
\end{array}
\right.
\end{flalign*}
We say that:\newline
$\left( \pmV_{\setindexdue_{g}},f_{\setindexdue_{g}}\right)$ is a skeleton of $\pmgenerder\funcomp g$ for any $g\in \sieve$;\newline
$\left( \pmV_{\setindexdue_{\sieve}},f_{\setindexdue_{\sieve}}\right)$ is a skeleton of $\pmgenerder$ relative to the covering sieve $\sieve$;\newline
$\genpreder\left[\pmV_{\setindexdue_{g}},f_{\setindexdue_{g}} \right]$ is a realization of $\pmgenerder\funcomp g$ for any $g\in \sieve$;\newline
$\genpreder\left[\pmV_{\setindexdue_{\sieve}},f_{\setindexdue_{\sieve}} \right]$ is a realization of $\pmgenerder$ relative to the covering sieve $\sieve$;\newline
$\pmpreder_{\setindexdue_{g}}$ is a decomposition of $\pmgenerder\funcomp g$ for any $g\in \sieve$;\newline
$\pmpreder_{\setindexdue_{\sieve}}$ is a decomposition of $\pmgenerder$ relative to the covering sieve $\sieve$.\newline
If $\sieve$ is the maximal sieve then we only consider skeleton, realization and decomposition relative to arrow $\idobj+\setsymuno+\in \sieve$ and we drop any reference to both $\idobj+\setsymuno+$ and $\sieve$ in the notation.
\end{enumerate}
\end{remark}

In Proposition \ref{carV} below we give a necessary condition for an element belonging to $\genprederspace[\modsym](l)+m+$ to belong to the $\mathbb{R}$-vector subspace $\nullpreder[\modsym](l)+m+$.\newline
We refer to Notation \ref{difrealfunc}-[3], Definitions \ref{domcodfuncont}, \ref{genpointdef}, Remarks \ref{pathinlpreder}, \ref{pathinlgenerder}, Example \ref{ExCzero}.

\begin{proposition}\label{carV}
Fix $l,m \in \mathbb{N}_0$, a $\mathbb{R}$-vector subspace $\modsym$ of $\mathbb{R}^l$, a pre-derivation $\genpreder\in\nullpreder[\modsym](l)+m+$.\newline
Then for any $\gggfw[\gfw] \in\genfquotcccz(m)+1+$, path $\pmtwo$ in $\magtwo$ detecting $\genpreder\left(\gggfw[\gfw]\right)$ with
\begin{equation}
 \exists\delta>0\;:\; \left(\evalcomptwo\left( \pmtwo\left(\unkdue\right)\right)\right)\left(\unkuno\right)\neq\udenunk\quad
\forall \unkdue\in \left(-\delta,\delta\right)\setminus\left\{0\right\}\quad \forall\unkuno\in \dommagtwo\left(\pmtwo\right) \label{esistsemp}
\end{equation}
the continuous function $\pmtwo_{\smint}:\left(-1,1\right)\setminus\left\{0\right\}\rightarrow\Cksp{0}$ given by setting 
\begin{equation*}
\pmtwo_{\smint}\left(\unkdue\right)=\smint_{1}\left(... \smint_{l}\pmtwo\left(\unkdue\right)\right) \hspace{15pt}\forall\unkdue\in\left(-1,1\right)\setminus\left\{0\right\}
\end{equation*}
converges in $\Cksp{0}$ to the constant $0$ function for $\unkdue\rightarrow 0$.
\end{proposition}
\begin{proof}
Relations \eqref{relpreder}, definition of $\nullpreder[\modsym](l)+m+$ in Remark \ref{pathinlgenerder}-[2]entail that there is no loss of generality by assuming that there are a finite set $\setindexuno$, $k_{\elindexuno}, l,m \in \mathbb{N}_0$, $ \elindexotto_{\elindexuno}\in \left\{1,...,l\right\}$, $\vecuno_{\elindexuno}\in \mathbb{R}$, 
$\gggfz[\gfz_{\elindexuno}]\in \genfquotcccz(l)+k_{\elindexuno}+$, $\gggfx[\gfx_{\elindexuno,1}],\gggfx[\gfx_{\elindexuno,2}]\in \genfquotcccz(k_{\elindexuno})+m+$ with $\evalcompquotcontcontzero\left(\gggfx[ \gfx_{1,\elindexuno}]\right)=\evalcompquotcontcontzero\left(\gggfx[\gfx_{2,\elindexuno}]\right)$ for any $\elindexuno \in \setindexuno$ such that 
\begin{equation*}
\genpreder= \underset{\elindexuno\in \setindexuno}{\bigsumgenpreder}\, \vecuno_{\elindexuno}\scalpgenpreder\left(
\genprediff(l)+\gggfx[\gfx_{\elindexuno,1}]+\left( \genpreder\left[\gggfz[\gfz_{\elindexuno}],  \ef_{\elindexotto_{\elindexuno}} \right] \right)\,\sumgenpreder\, -1\scalpgenpreder
\genprediff(l)+\gggfx[\gfx_{\elindexuno,2}]+\left(  \genpreder\left[\gggfz[\gfz_{\elindexuno}],\ef_{\elindexotto _{\elindexuno}}\right]\right)
   \right)\text{.}
\end{equation*}
Set $\hspace{15pt}\gfy_{\elindexuno,0}=\left(\gfw\genfuncomp  \left(\gfx_{\elindexuno,1}\genfuncomp\gfz_{\elindexuno}\right)\right)-\left( \gfw\genfuncomp \left(\gfx_{\elindexuno,2}\genfuncomp\gfz_{\elindexuno}\right)\right) \hspace{15pt}\forall \elindexuno\in \setindexuno$.\newline
Then
\begin{equation*}
\genpreder\left(\gggfw[\gfw]\right)= \underset{\elindexuno\in \setindexuno}{\bigsumquotpoint}\, \vecuno_{\elindexuno}\scalpquotpoint \left( \Fpart_{\elindexotto _{\elindexuno}} \bsfmuquotpoint \gggfy[\gfy_{\elindexuno,0}]\right) \text{.} 
\end{equation*} 
Propositions \ref{pathprop3}, \ref{pathprop4bis}, \ref{pathprop5bis} entail that there is a path $\pmtwo_{\elindexuno}$ in $\genf$ detecting $\gfy_{\elindexuno,0}$ for any $\elindexuno\in \setindexuno$ such that 
\begin{equation*}
\pmtwo= \underset{\elindexuno\in \setindexuno}{\biggenfunsum}\, \vecuno_{\elindexuno}\genfunscalp \left( \Fpart_{\elindexotto _{\elindexuno}} \genfbsfmu \pmtwo_{\elindexuno}\right) \text{.}
\end{equation*} 
Proposition \ref{pathprop7}, Remark \ref{existoperquot}-[1] together entail that there is no loss of generality by assuming that all conditions below hold true:
\begin{flalign*}
&\begin{array}{l}
\gfy_{\elindexuno,0}\notin \magtwoempty \hspace{15pt}\forall \elindexuno\in \setindexuno\text{;}
\end{array}\\[8pt]
&\begin{array}{l}
\evalcomptwo\left(\gfx_{\elindexuno,1}\right)=\evalcomptwo\left(\gfx_{\elindexuno,2}\right) \hspace{15pt}\forall \elindexuno\in \setindexuno\text{;}
\end{array}\\[8pt]
&\begin{array}{l}
\evalcomptwo\left(\gfy_{\elindexuno,0}\right)=\cost<\dommagtwo\left(\gfy_{\elindexuno,0}\right)<>\mathbb{R}>+0+ \hspace{15pt}\forall \elindexuno\in \setindexuno\text{;}
\end{array}\\[8pt]
&\begin{array}{l}
\text{there}\hspace{4pt}\text{is}\hspace{4pt} \delta>0\hspace{4pt}\text{such}\hspace{4pt}\text{that}\hspace{4pt} \pmtwo_{\elindexuno}\left(\unkdue\right)\notin \magtwoempty \hspace{4pt}\text{for}\hspace{4pt}\text{any}\hspace{4pt}\elindexuno\in \setindexuno\text{,}\hspace{4pt}\unkdue\in \left(-\delta,\delta\right)\text{.}
\end{array}
\end{flalign*}
Eventually statement follows by \eqref{derint}-[$(ii)$, $(iii)$], \eqref{intdiffop}-[$(v)$], Proposition \ref{topinDprop}-[5], relations \eqref{(R 10 bis)}, \eqref{(R 16)}.
\end{proof}

In Proposition \ref{algstrderpoint} below we prove that for any not trivial $\mathbb{R}$-vector subspace $\modsym$ of $\mathbb{R}^l$ the $\mathbb{R}$-vector space $\generderspace[\modsym](l)+m+$ is not trivial. We refer to Notations \ref{int}, \ref{realvec}-[2].

\begin{proposition}\label{algstrderpoint}
Fix $l\in \mathbb{N}_0$, $m \in \mathbb{N}_0$, $\mathsf{l}\in \left\{1,...,l\right\}$, a $\mathbb{R}$-vector subspace $\modsym$ of $\mathbb{R}^l$, $\vecuno_{1}, \vecuno_{2}\in\modsym$.
Define the smooth function $h:\mathbb{R}^l\rightarrow \mathbb{R}^m$ by setting
\begin{equation*}
h\left(\unkuno_1,...,\unkuno_l\right)=\left(\unkuno_{\mathsf{l}},0,...,0\right)\hspace{15pt}\forall\left(\unkuno_1,...,\unkuno_l\right)\in\mathbb{R}^l\text{.} 
\end{equation*}
Then\hspace{10pt}$\generder[\argcompl{\genpreder\left[\lininclquot\left(\gczsfx[h]\right),\vecuno_{1}\right]}]=\generder[\argcompl{\genpreder\left[\lininclquot\left(\gczsfx[h]\right),\vecuno_{2}\right]}]\quad \Leftrightarrow\quad\vecuno_{1, \mathsf{l}}=\vecuno_{2,\mathsf{l}}$.
\end{proposition}
\begin{proof}Implication $\Leftarrow$ is straightforward, then we just prove $\Rightarrow$. Relations \eqref{relpreder} entail that statement follows by proving the claim 
\begin{equation}\label{notinvecsp}
\genpreder\left[\lininclquot\left(\gczsfx[h]\right),\vecuno\right]\in\nullpreder[\modsym](l)+m+\quad \Rightarrow\quad \vecuno_{\mathsf{l}}=0\text{.}
\end{equation}
Referring to Notations \ref{realfunc}-[5(c)], \ref{difrealfunc}, we have that Proposition \ref{carV} applied to the evaluation of $\genpreder\left[\lininclquot\left(\gczsfx[h]\right),\vecuno\right]$ at $\lininclquot\left(\gczsfx[\coord{m}{1}]\right)$, assumption $\genpreder\left[\lininclquot\left(\gczsfx[h]\right),\vecuno\right]\in\nullpreder[\modsym](l)+m+$ entail that there is an open neighborhood $\intuno\subseteq\mathbb{R}^l$ of $0\in \mathbb{R}^l$ such that 
\begin{multline}
\left(\smint_{1}\left(... \smint_{l} \left(\smder_{\vecuno_{\mathsf{l}}}\left(h\right)\right) \right)\right)\left(\unkuno_1,...,\unkuno_{2l}\right)=0 \\[4pt]
\forall \left(\unkuno_1,...,\unkuno_{2l}\right)\in \domint\left[...\domint\left[\intuno,l\right],...1\right]\text{.}\label{zerointeg}
\end{multline}
Eventually statement follows since elementary calculus entails 
\begin{multline*}
\left(\smint_{1}\left(... \smint_{l} \left(\smder_{\vecuno_{\mathsf{l}}}\left(h\right)\right) \right)\right)\left(\unkuno_1,...,\unkuno_{2l}\right)=
\vecuno_{\mathsf{l}}\vecprod\underset{\elindexsei=1}{\overset{l}{\prod}} \left(\unkuno_{2\elindexsei}-\unkuno_{2\elindexsei-1}\right) 
\\[4pt] \forall\left(\unkuno_1,...,\unkuno_{2l}\right)\in \domint\left[...\domint\left[\intuno,l\right],...1\right]
\end{multline*}
then \eqref{zerointeg} entail that $\vecuno_{\mathsf{l}}=0$ .
\end{proof}

In Proposition \ref{lindiff} below we give a sufficient condition for a pair of germs of local continuous pointed functions to ensure that the differential of their sum coincides with the sum of their differentials. We refer to Remark \ref{pathinlgenerder}-[6].

\begin{proposition}\label{lindiff}
Fix $m,n \in \mathbb{N}_0$, a $\mathbb{R}$-vector subspace $\modsym$ of $\generderspace(1)+m+$, $\gczsfx, \gczsfy \in \Ckspquotccz{0}(m)+n+$. Assume that both conditions below hold true:
\begin{flalign}
&\begin{array}{l}
\generderspace(1)+m+=\modsym \dirsum 	\ker\left(\gendiff(1)+\gczsfx+\right)\text{;}\label{dirsumDm}
\end{array}\\[8pt]
&\begin{array}{l}
\modsym \subseteq	\ker\left(\gendiff(1)+\gczsfy+\right)\text{.}\label{dentrodirsum}
\end{array}
\end{flalign}
Then $\gendiff(1)+\argcompl{\gczsfx \sumquotpoint\gczsfy}+ = \gendiff(1)+\gczsfx+  \sumquotpoint  \gendiff(1)+\gczsfy+$.
\end{proposition}
\begin{proof} Since
\begin{equation}
\left\{
\begin{array}{l}
\left\{\generder[\argcompl{\genpreder\left[ \gggfz,\ef_1\right]}]\; :\; \gggfz\in \genfquotcccz(1)+m+ \right\}
\hspace{4pt}\text{generates}\hspace{4pt}\generderspace(1)+m+\\[4pt]
\text{as}\hspace{4pt}\text{a}\hspace{4pt}\mathbb{R}\text{-vector}\hspace{4pt}\text{space}
\end{array}\label{genDmRvs}
\right.
\end{equation}
by construction, then statement will follows by proving claim
\begin{multline}
\gendiff(1)+\argcompl{\gczsfx \sumquotpoint \gczsfy}+\left(\generder[\argcompl{\genpreder\left[ \gggfz,\ef_1\right]}]\right) = \gendiff(1)+\gczsfx+\left(\generder[\argcompl{\genpreder\left[ \gggfz,\ef_1\right]}]\right)  \sumquotpoint  \gendiff(1)+\gczsfy+\left(\generder[\argcompl{\genpreder\left[ \gggfz,\ef_1\right]}]\right)\\[4pt]
\forall \gggfz\in \genfquotcccz(1)+m+ \text{.}\label{valutconfr}
\end{multline}
Fix $\gggfz\in \genfquotcccz(1)+m+$. Assumptions \eqref{dirsumDm}, \eqref{dentrodirsum} together with \eqref{genDmRvs} entails that there is no loss of generality by assuming that at least one among the two conditions below holds true:
\begin{flalign}
&\generder[\argcompl{\genpreder\left[ \gggfz,\ef_1\right]}]\in \ker\left(\gendiff(1)+\gczsfx+\right)\text{;}\label{inkerx}\\[6pt]
&\generder[\argcompl{\genpreder\left[ \gggfz,\ef_1\right]}]\in \ker\left(\gendiff(1)+\gczsfy+\right)\text{.}\label{inkery}
\end{flalign} 
We prove claim \eqref{valutconfr} assuming \eqref{inkerx}, the proof of claim \eqref{valutconfr} in case \eqref{inkery} holds true is the same.\newline 
Fix $\gggfx,\gggfy\in \genfquotcccz(1)+m+$ with $\evalcompquotcontcontzero\left(\gggfx\right)=\gczsfx$ and $\evalcompquotcontcontzero\left(\gggfy\right)=\gczsfy$. Chain of equalities
\begin{multline*}
\gendiff(1)+\argcompl{\gczsfx \sumquotpoint \gczsfy}+\left(\generder[\argcompl{\genpreder\left[ \gggfz,\ef_1\right]}]\right)=\\[4pt]
\generder[\argcompl{\genpreder\left[ \left(\gggfx \sumquotpoint \gggfy\right)\compquotpoint\gggfz,\ef_1\right]}]=\generder[\argcompl{\genpreder\left[ \left(\gggfx \compquotpoint\gggfz\right)\sumquotpoint \left(\gggfy\compquotpoint\gggfz\right),\ef_1\right]}]
\end{multline*}
entails that claim \eqref{valutconfr} follows by proving that
\begin{equation}
\generder[\argcompl{\genpreder\left[ \left(\gggfx \compquotpoint\gggfz\right)\sumquotpoint \left(\gggfy\compquotpoint\gggfz\right),\ef_1\right]}]=\generder[\argcompl{\genpreder\left[ \gggfy\compquotpoint\gggfz,\ef_1\right]}]\text{.}\label{piuugualnon}
\end{equation} 
By \eqref{propevaltwo}, Proposition \ref{genpointprop}-[11] and Remark \ref{pathinlgenerder}-[2], equality \eqref{piuugualnon} follows by proving that 
\begin{equation}
\evalcompquotcontcontzero\left(\gggfx \compquotpoint\gggfz\right)=\gczsfx[\cost<\mathbb{R}<>\mathbb{R}^n>+0+]\text{.}\label{zerovalut}
\end{equation}
By assumption \eqref{inkerx}, Proposition \ref{carV} entails that there is $\unkuno\in \mathbb{R}^n$ such that $\evalcompquotcontcontzero\left(\gggfx \compquotpoint\gggfz\right)= \gczsfx[\cost<\mathbb{R}<>\mathbb{R}>+\unkuno+]$. Eventually \eqref{zerovalut} follows since $\gggfx,\gggfz \in \genfquotcccz$ entails that $\unkuno=0$. 
\end{proof}

In Definition \ref{smsubsdef} below we focus attention on generalized $\left(\modsym,l\right)$-derivations of degree $1$ represented by pre-derivations with at least one smooth core.  We refer to Definition \ref{gpdssm}, Proposition \ref{genpointprop}-[10], Remarks \ref{genercontfun}, \ref{existoperquot}-[2], \ref{extremAnCC0}.

\begin{definition}\label{smsubsdef}
Fix $l,m \in \mathbb{N}_0$, a $\mathbb{R}$-vector subspace $\modsym$ of $\mathbb{R}^l$. We denote by $\generderspacesm[\modsym](l)+m+$ the $\mathbb{R}$-vector subspace of $\generderspace[\modsym](l)+m+$ defined by setting
\begin{equation}
\generderspacesm[\modsym](l)+m+=\left\{\generder[\genpreder]\hspace{4pt}:\hspace{4pt}\genpreder\in \genprederspacesm[\modsym](l)+m+\right\}\text{.}\label{gensmgdss}
\end{equation}
Elements belonging to the $\mathbb{R}$-vector space $\generderspacesm[\modsym](l)+m+$ are called $m$-dimensional generalized smooth $\left(\modsym,l\right)$-derivations of degree $1$, or $m$-dimensional generalized smooth $\left(\modsym,l\right)$-tangent vectors.\newline
We set $\inclgdsmgd[\modsym](l)+m+=\incl{\generderspacesm[\modsym](l)+m+}{\generderspace[\modsym](l)+m+}$.
\end{definition}

In Proposition \ref{nattosm} below we study the stability of the $\mathbb{R}$-vector subspace $\generderspacesm[\modsym](l)+m+$ of $\generderspace[\modsym](l)+m+$ with respect to arrows of $\GTopcat$ which are relevant for our aims. We refer to Notations \ref{realvec}-[2], \ref{difrealfunc}-[3], Proposition \ref{genpointprop}-[10], \eqref{pcM}, Remarks \ref{genercontfun}, \ref{existoperquot}-[2], \ref{extremAnCC0}, \ref{pathinlpreder}-[6].

\begin{proposition}\label{nattosm}\mbox{}
\begin{enumerate}
\item Fix $l,m \in \mathbb{N}_0$, a $\mathbb{R}$-vector subspace $\modsym$ of $\mathbb{R}^l$. Then there is one and only one $\mathbb{R}$-linear arrow $\quotsm[\modsym](l)+m+:\generderspacesm[\modsym](l)+m+\rightarrow \mathbb{R}^m$ of $\GTopcat$ defined by setting 
\begin{equation}
\quotsm[\modsym](l)+m+\left(\generder[\genpreder]\right) =
\evspd[\modsym](l)+m+\left(\genpreder\right)\hspace{20pt}
\forall\genpreder\in \genprederspacesm[\modsym](l)+m+\text{.}
\label{eqtrcol}
\end{equation}
\item Fix $m \in \mathbb{N}_0$. Then there is one and only one $\mathbb{R}$-linear arrow
\begin{equation*}
\invquotsm+m+: \mathbb{R}^m\rightarrow\generderspacesm[\mathbb{R}^m](m)+m+
\end{equation*}
of $\GTopcat$ defined by setting 
\begin{equation*}
\invquotsm+m+\left(\vecuno\right)=\generder[\argcompl{\genpreder\left[\gggfx[\argcompl{\idobj+\mathbb{R}^m+}], \vecuno \right]}]\hspace{20pt}\forall \vecuno \in \mathbb{R}^m\text{.}
\end{equation*}
Arrow $\invquotsm+m+$ fulfills
\begin{equation}
\quotsm[\mathbb{R}^m](m)+m+\funcomp\invquotsm+m+=\idobj+\mathbb{R}^m+\text{.}\label{smpr1}
\end{equation}
\item Fix $l,m,n \in \mathbb{N}_0$, a $\mathbb{R}$-vector subspace $\modsym$ of $\mathbb{R}^l$, $\gczsfx  \in \Ckspquot{\infty}(m)+n+$.\newline
Then there is one and only one $\mathbb{R}$-linear arrow  
\begin{equation*}
\gendiffsm[\modsym](l)+\lininclquot\left(\gczsfx\right)+:\generderspacesm[\modsym](l)+m+\rightarrow\generderspacesm[\modsym](l)+n+
\end{equation*}
of $\GTopcat$ defined by setting 
\begin{equation*}
\gendiffsm[\modsym](l)+\lininclquot\left(\gczsfx \right)+\left(\generder[\genpreder]\right)=\generder[ \genprediff[\modsym](l)+\lininclquot\left(\gczsfx \right)+ \left(\genpreder\right)]\hspace{20pt}
\forall\genpreder\in \genprederspacesm[\modsym](l)+m+\text{.}
\end{equation*}
Arrow $\gendiffsm[\modsym](l)+\lininclquot\left(\gczsfx\right)+ $ fulfills both conditions below:
\begin{flalign}
& \gendiff[\modsym](l)+\lininclquot\left(\gczsfx\right)+\funcomp \inclgdsmgd[\modsym](l)+m+=\inclgdsmgd[\modsym](l)+n+\funcomp \gendiffsm[\modsym](l)+ \lininclquot\left(\gczsfx \right)+\text{;}\label{smpr2}\\[6pt]
&\symjacsm<1<>\argcompl{\gczsfx}>\funcomp\quotsm[\modsym](l)+m+ =\quotsm[\modsym](l)+n+ \funcomp \gendiffsm[\modsym](l)+\lininclquot\left(\gczsfx \right)+\text{.}\label{smpr3}
\end{flalign}
\item Fix $l,m,n \in \mathbb{N}_0$, two $\mathbb{R}$-vector subspaces $\modsym_1$, $\modsym_2$ of $\mathbb{R}^l$.\newline
Then there is one and only one $\mathbb{R}$-linear arrow  
\begin{equation*}
\generinclsm<\modsym_1<>\modsym_2>(l)+m+:\generderspacesm[\modsym_1](l)+m+\rightarrow\generderspacesm[\modsym_2](l)+m+
\end{equation*}
of $\GTopcat$ defined by setting 
\begin{equation*}
\generinclsm<\modsym_1<>\modsym_2>(l)+m+\left(\generder[\genpreder]\right)=\generder[ \genpreincl<\modsym_1<>\modsym_2>(l)+m+\left(\genpreder\right)]\hspace{8pt}
\forall\genpreder\in \genprederspacesm[\modsym](l)+m+\text{.}
\end{equation*}
Arrow $\generinclsm<\modsym_1<>\modsym_2>(l)+m+$ fulfills both conditions below:
\begin{flalign}
& \inclgdsmgd[\modsym_2](l)+m+\funcomp \generinclsm<\modsym_1<>\modsym_2>(l)+m+=\generincl<\modsym_1<>\modsym_2>(l)+m+\funcomp \inclgdsmgd[\modsym_1](l)+m+\text{;}\label{smpr7}
\\[6pt]
&\quotsm[\modsym_2](l)+m+\funcomp\generinclsm<\modsym_1<>\modsym_2>(l)+m+=\quotsm[\modsym_1](l)+m+\text{;}\label{smpr8}
\\[6pt]
&\begin{array}{l}
\gendiffsm[\modsym_2](l)+\lininclquot\left(\gczsfx \right)+ \funcomp \generinclsm<\modsym_1<>\modsym_2>(l)+m+=\\
\hspace{26pt}\generinclsm<\modsym_1<>\modsym_2>(l)+n+\funcomp \gendiffsm[\modsym_1](l)+\lininclquot\left(\gczsfx \right)+\hspace{10pt} \forall \gczsfx  \in \Ckspquot{\infty}(m)+n+\text{.}
\label{smpr9}
\end{array} 
\end{flalign}
\end{enumerate}
\end{proposition}
\begin{proof}\mbox{}\newline
\textnormal{\textbf{Proof of statement 1.}}\newline
Statement is achieved by proving that $\quotsm[\modsym](l)+m+$ is well defined. This is done by evaluating right side of \eqref{eqtrcol} at different representatives of $\generder[\genpreder]$. We refer to Proposition \ref{gpdssmprop}.\newline  
Fix $l\in \mathbb{N}_0$, $\vecuno_{1}, \vecuno_{2}\in\mathbb{R}^{l}$, $\gczsfx[f_{1}], \gczsfx[f_{2}] \in \Ckspquot{\infty}(l)+m+$.
Referring to Notations \ref{realfunc}-[1], \ref{difrealfunc}, Proposition \ref{carV} and assumption $\generder[\argcompl{\genpreder\left[ \lininclquot\left(\gczsfx[f_{1}] \right) ,\vecuno_{1}\right]}]= \generder[\argcompl{\genpreder\left[ \lininclquot\left(\gczsfx[f_{2}] \right) ,\vecuno_{2}\right]}]$, finiteness of set $\left\{1,...,l\right\}$ entail that there is an open neighborhood $\intuno\subseteq\mathbb{R}^l$ of $0\in \mathbb{R}^l$ such that
\begin{flalign}
&\begin{array}{l}
f_{1}, f_{2}\in \Cksp{\infty}(\intuno)+\mathbb{R}^{m}+\text{;}\nonumber
\end{array}\\[8pt]
&\left\{
\begin{array}{l}
\left(\smint_{1}\left(... \smint_{l} \left(\smder_{\vecuno_1}\left( f_{1,\mathsf{m}}    \right)\right) \right)\right)\left(\unkuno\right)=\left(\smint_{1}\left(... \smint_{l} \left(\smder_{\vecuno_2}\left( f_{2,\mathsf{m}}    \right)\right) \right)\right)\left(\unkuno\right)\\[4pt]
\text{for}\hspace{4pt}\text{any}\hspace{4pt} \mathsf{m}\in\left\{1,...,m\right\} \text{,}\hspace{4pt} \unkuno\in \domint\left[...\domint\left[\intuno,l\right],...1\right]\text{.}\label{zerointegdue}
\end{array}
\right.
\end{flalign}
Define smooth functions
\begin{flalign*}
&\left\{
\begin{array}{l}
h\in \Cksp{\infty}(\mathbb{R}^{l})+\mathbb{R}^{2l}+	\hspace{4pt}\text{by}\hspace{4pt}\text{setting}\\[4pt]
h_{\elindexsei}\left(\unkuno_1,...,\unkuno_l\right)=\left\{
\begin{array}{ll}
\unkuno_{\elindexsei/2}&\text{if}\hspace{4pt}\elindexsei\in \left\{2\elindexcinque\;:\;\elindexcinque\in \left\{1,...,l\right\}\right\}\text{,}	\\[4pt]
0&\text{if}\hspace{4pt}\elindexsei\in \left\{2\elindexcinque-1\;:\;\elindexcinque\in \left\{1,...,l\right\}\right\}\text{;}
\end{array}
\right.
\end{array}
\right.\\[8pt]
&\left\{
\begin{array}{l}
g_{i,\mathsf{m}}\in \Cksp{\infty}(\intuno)+1+\hspace{4pt}\text{by}\hspace{4pt}\text{setting}\\[4pt]
\hspace{50pt}g_{i,\mathsf{m}}\left(\unkuno\right)= \left(\smint_{1}\left(... \smint_{l} \left(\smder_{\vecuno_i}\left( f_{i,\mathsf{m}}\right)\right) \right)\right)  \left(h\left(\unkuno\right)\right) \hspace{15pt}\forall \unkuno \in\intuno\\[4pt]
\text{for}\hspace{4pt}\text{any}\hspace{4pt} i\in\left\{1,2\right\} \text{,}\hspace{4pt}\mathsf{m}\in\left\{1,...,m\right\} \text{.}
\end{array}
\right.
\end{flalign*}
Eventually statement follows since elementary calculus entails 
\begin{equation*}
\begin{array}{l}
\smder_{\vecuno_i}\left( f_{i,\mathsf{m}}\right)=\smder_l\left(...\smder_1\left(g_{i,\mathsf{m}}\right)...\right)\hspace{4pt} \text{for}\hspace{4pt}\text{any}\hspace{4pt} i\in\left\{1,2\right\} \text{,}\hspace{4pt}\mathsf{m}\in\left\{1,...,m\right\}
\end{array}
\end{equation*}
then \eqref{zerointegdue} entails that  $\smder_{\vecuno_1}\left( f_{1,\mathsf{m}}\right)= \smder_{\vecuno_1}\left( f_{1,\mathsf{m}}\right)$ for any $\mathsf{m}\in\left\{1,...,m\right\}$.\newline
\textnormal{\textbf{Proof of statement 2.}}\ \ Arrow $\invquotsm+m+$ is well defined by \eqref{pcM} and definition of Grothendieck topology on $\GTopcat$. Equality \eqref{smpr1} follow straightforwardly by performing computations.\newline
\textnormal{\textbf{Proof of statement 3.}}\ \ Arrow $\gendiffsm[\modsym](l)+\lininclquot\left(\gczsfx\right)+$ is well defined since arrow $\gendiff[\modsym](l)+\lininclquot\left(\gczsfx \right)+$ is well defined and $\genprediff[\modsym](l)+\lininclquot\left(\gczsfx\right)+$ carries pre-derivations with a smooth core to pre-derivations with a smooth core. \newline
Equalities \eqref{smpr2}, \eqref{smpr3} follow straightforwardly by performing computations.\newline
\textnormal{\textbf{Proof of statement 4.}}\ \ Arrow $\generinclsm<\modsym_1<>\modsym_2>(l)+m+$ is straightforwardly well defined. Equalities \eqref{smpr7}-\eqref{smpr9} follow straightforwardly by performing computations.
\end{proof}

In Proposition \ref{diagRvs} below we prove that the material introduced above can be organized in an exact diagram of $\mathbb{R}$-vector spaces in $\GTopcat$ which is natural with respect to differential of smooth functions.  We refer to Notation \ref{difrealfunc}-[3], Proposition \ref{genpointprop}-[10], Remarks \ref{genercontfun}, \ref{existoperquot}-[2], \ref{extremAnCC0}.

\begin{proposition}\label{diagRvs}\mbox{}
\begin{enumerate}
\item Fix $l,m \in \mathbb{N}_0$, a $\mathbb{R}$-vector subspace $\modsym$ of $\mathbb{R}^l$. Then there is an exact diagram of $\mathbb{R}$-vector spaces in $\GTopcat$
\begin{equation}
\xymatrix{
&&0 \ar[d]\\
0 \ar[r]& \kergenerderspacesm[\modsym](l)+m+\ar[r]^{\kgdssmar[\modsym](l)+m+}&\generderspacesm[\modsym](l)+m+\ar[d]^{\inclgdsmgd[\modsym](l)+m+}\ar[r]^{\quotsm[\modsym](l)+m+}& \mathbb{R}^{m}\ar[r]&0\\
&&\generderspace[\modsym](l)+m+\ar[d]^{\ckgdssmar[\modsym](l)+m+}\\
&&\cokergenerderspace[\modsym](l)+m+\ar[d]\\
&&0}
\label{fig1}
\end{equation}
where: $\kgdssmar[\modsym](l)+m+$ is the kernel of $\quotsm[\modsym](l)+m+$; $\ckgdssmar[\modsym](l)+m+$ is the co-kernel of $\inclgdsmgd[\modsym](l)+m+$.\newline
If $\modsym=\mathbb{R}^m$ then the exact row in diagram \eqref{fig1} splits and splitting is given by $\invquotsm+m+$. 
\item Fix $l, m,n \in \mathbb{N}_0$, a $\mathbb{R}$-vector subspace $\modsym$ of $\mathbb{R}^l$, $ \gczsfx  \in \Ckspquot{\infty}(m)+n+$. Then the diagram below commutes
\begin{equation}
\xy
% prima riga dall'alto
(-36,15)*+{0}="a1",
(-22,15)*+{\scalebox{.7}{$\kergenerderspacesm[\modsym](l)+n+$}}="b1",
(10,15)*+{\scalebox{.7}{$\generderspacesm[\modsym](l)+n+$}}="c1",
(42,15)*+{\scalebox{.7}{$\mathbb{R}^{n}$}}="d1",
(58,15)*+{0}="e1",
% seconda riga dall'alto
(-52,0)*+{0}="a2",
(-38,0)*+{\scalebox{.7}{$\kergenerderspacesm[\modsym](l)+m+$}}="b2",
(-6,0)*+{\scalebox{.7}{$\generderspacesm[\modsym](l)+m+$}}="c2",
(28,0)*+{\scalebox{.7}{$\mathbb{R}^{m}$}}="d2",
(44,0)*+{0}="e2",
% prima colonna da sinistra
(-6,10)*+{0}="cs1",
(-6,-30)*+{\scalebox{.7}{$\generderspace[\modsym](l)+m+$}}="cs3",
(-6,-60)*+{\scalebox{.7}{$\cokergenerderspace[\modsym](l)+m+$}}="cs4",
(-6,-70)*+{0}="cs5",
% seconda colonna da sinistra
(10,25)*+{0}="cd1",
(10,-15)*+{\scalebox{.7}{$\generderspace[\modsym](l)+n+$}}="cd3",
(10,-45)*+{\scalebox{.7}{$\cokergenerderspace[\modsym](l)+n+$}}="cd4",
(10,-55)*+{0}="cd5",
% primo quadrato commutativo
(17,7)*{\scalebox{.7}{A}},
% secondo quadrato commutativo
(2,-7)*{\scalebox{.7}{B}},
% frecce prima riga
\ar @{->} "a1";"b1"
\ar @{->}^<>(.7){\kgdssmar[\modsym](l)+n+} "b1";"c1"
\ar @{->}^<>(.7){\quotsm[\modsym](l)+n+} "c1";"d1"
\ar @{->} "d1";"e1"
% frecce seconda riga
\ar @{->} "a2";"b2"
\ar @{->}^<>(.7){\kgdssmar[\modsym](l)+m+} "b2";"c2"
\ar @{->}^<>(.7){\quotsm[\modsym](l)+m+} "c2";"d2"
\ar @{->} "d2";"e2"
% frecce prima colonna da sinistra
\ar @{->} "cs1";"c2"
\ar @{->}_<>(.6){\inclgdsmgd[\modsym](l)+m+} "c2";"cs3"
\ar @{->}_<>(.6){\ckgdssmar[\modsym](l)+m+} "cs3";"cs4"
\ar @{->} "cs4";"cs5"
% frecce seconda colonna da sinistra
\ar @{->} "cd1";"c1"
\ar @{->}|<>(.5){\hole}^<>(.6){\inclgdsmgd[\modsym](l)+n+} "c1";"cd3"
\ar @{->}^<>(.6){\ckgdssmar[\modsym](l)+n+} "cd3";"cd4"
\ar @{->} "cd4";"cd5"
% frecce diagonali
\ar @{->}|<>(.5){\scalebox{.5}{$\kersqA[\modsym](l)+\lininclquot\left(\gczsfx \right)+$}} "b2";"b1"
\ar @{->}|<>(.5){\scalebox{.5}{$\gendiffsm[\modsym](l)+\lininclquot\left(\gczsfx \right)+$}} "c2";"c1"
\ar @{->}|<>(.5){\scalebox{.5}{$\symjacsm<1<>\argcompl{\gczsfx}>$}} "d2";"d1"
\ar @{->}|<>(.5){\scalebox{.5}{$\gendiff[\modsym](l)+\lininclquot\left(\gczsfx \right)+$}} "cs3";"cd3"
\ar @{->}|<>(.5){\scalebox{.5}{$\cksqB[\modsym](l)+\lininclquot\left(\gczsfx\right)+$}} "cs4";"cd4"
\endxy
\label{fig2}
\end{equation}
where: $\kersqA[\modsym](l)+\lininclquot\left(\gczsfx \right)+$ is the natural arrow induced on kernels by commutative square A; $\cksqB[\modsym](l)+\lininclquot\left(\gczsfx \right)+$ is the natural arrow induced on co-kernels by commutative square B. 
\item Fix $l, m \in \mathbb{N}_0$, two $\mathbb{R}$-vector subspaces $\modsym_1$, $\modsym_2$ of $\mathbb{R}^l$. Then the diagram below commutes 
\begin{equation}
\xy
% prima riga dall'alto
(-36,15)*+{0}="a1",
(-22,15)*+{\scalebox{.7}{$\kergenerderspacesm[\modsym_2](l)+m+$}}="b1",
(10,15)*+{\scalebox{.7}{$\generderspacesm[\modsym_2](l)+m+$}}="c1",
(42,15)*+{\scalebox{.7}{$\mathbb{R}^ {m}$}}="d1",
(58,15)*+{0}="e1",
% seconda riga dall'alto
(-52,0)*+{0}="a2",
(-38,0)*+{\scalebox{.7}{$\kergenerderspacesm[\modsym_1](l)+m+$}}="b2",
(-6,0)*+{\scalebox{.7}{$\generderspacesm[\modsym_1](l)+m+$}}="c2",
(28,0)*+{\scalebox{.7}{$\mathbb{R}^{m}$}}="d2",
(44,0)*+{0}="e2",
% prima colonna da sinistra
(-6,10)*+{0}="cs1",
(-6,-30)*+{\scalebox{.7}{$\generderspace[\modsym_1](l)+m+$}}="cs3",
(-6,-60)*+{\scalebox{.7}{$\cokergenerderspace[\modsym_1](l)+m+$}}="cs4",
(-6,-70)*+{0}="cs5",
% seconda colonna da sinistra
(10,25)*+{0}="cd1",
(10,-15)*+{\scalebox{.7}{$\generderspace[\modsym_2](l)+m+$}}="cd3",
(10,-45)*+{\scalebox{.7}{$\cokergenerderspace[\modsym_2](l)+m+$}}="cd4",
(10,-55)*+{0}="cd5",
% primo quadrato commutativo
(17,7)*{\scalebox{.7}{A}},
% secondo quadrato commutativo
(2,-7)*{\scalebox{.7}{B}},
% frecce prima riga
\ar @{->} "a1";"b1"
\ar @{->}^<>(.7){\kgdssmar[\modsym_2](l)+m+} "b1";"c1"
\ar @{->}^<>(.7){\quotsm[\modsym_2](l)+m+} "c1";"d1"
\ar @{->} "d1";"e1"
% frecce seconda riga
\ar @{->} "a2";"b2"
\ar @{->}^<>(.7){\kgdssmar[\modsym_1](l)+m+} "b2";"c2"
\ar @{->}^<>(.7){\quotsm[\modsym_1](l)+m+} "c2";"d2"
\ar @{->} "d2";"e2"
% frecce prima colonna da sinistra
\ar @{->} "cs1";"c2"
\ar @{->}_<>(.6){\inclgdsmgd[\modsym_1](l)+m+} "c2";"cs3"
\ar @{->}_<>(.6){\ckgdssmar[\modsym_1](l)+m+} "cs3";"cs4"
\ar @{->} "cs4";"cs5"
% frecce seconda colonna da sinistra
\ar @{->} "cd1";"c1"
\ar @{->}|<>(.5){\hole}^<>(.6){\inclgdsmgd[\modsym_2](l)+m+} "c1";"cd3"
\ar @{->}^<>(.6){\ckgdssmar[\modsym_2](l)+m+} "cd3";"cd4"
\ar @{->} "cd4";"cd5"
% frecce diagonali
\ar @{->}|<>(.5){\scalebox{.5}{$\kergenerinclsm<\modsym_1<>\modsym_2>(l)+m+$}} "b2";"b1"
\ar @{->}|<>(.5){\scalebox{.5}{$\generinclsm<\modsym_1<>\modsym_2>(l)+m+$}} "c2";"c1"
\ar @{->}|<>(.5){\scalebox{.5}{$1_{\mathbb{R}^m}$}} "d2";"d1"
\ar @{->}|<>(.5){\scalebox{.5}{$\generincl<\modsym_1<>\modsym_2>(l)+m+$}} "cs3";"cd3"
\ar @{->}|<>(.5){\scalebox{.5}{$\ckgenerincl<\modsym_1<>\modsym_2>(l)+m+$}} "cs4";"cd4"
\endxy
\label{fig4}
\end{equation}
where: $\kergenerinclsm<\modsym_1<>\modsym_2>(l)+m+$ is the natural arrow induced on kernels by commutative square A; $\ckgenerincl<\modsym_1<>\modsym_2>(l)+m+$ is the natural arrow induced on co-kernels by commutative square B. 
\end{enumerate}
\end{proposition}
\begin{proof}
Statement follows by Proposition \ref{nattosm}, Remark \ref{pathinlgenerder}-[6, 7].
\end{proof}

In Remark \ref{diagRvsrem} below we explicitly describe some straightforward consequences of Propositions \ref{nattosm}, \ref{diagRvs}. 

\begin{remark}\label{diagRvsrem}
 Splitting of the exact row of diagram \eqref{fig1} is not natural since in general we have 
$\invquotsm+n+\funcomp\symjacsm<1<>\argcompl{\gczsfx}> \neq \gendiffsm[\mathbb{R}^l](l)+\lininclquot\left(\gczsfx \right)+ \funcomp\invquotsm+m+$.\newline
This happens because pre-differentials change cores of pre-derivations then generalized differentials of degree $1$ change cores of pre-derivations representing generalized derivations on which they acts.  
\end{remark}

In Definition \ref{produnodef}, Propositions \ref{produno}, \ref{produnosm} below we study the behavior of generalized derivation operators of degree $1$ (i.e. generalized tangent vectors) with respect to the sum of dimensions. We refer to Proposition \ref{genpointprop}-[3].

\begin{definition}\label{produnodef}
Fix $l,m_1,m_2\in\mathbb{N}_0$, $i \in \left\{1,2\right\}$, a $\mathbb{R}$-vector subspace $\modsym$ of $\mathbb{R}^l$. We define $\mathbb{R}$-linear arrows of 
$\GTopcat$
\begin{flalign*}
&\left\{
\begin{array}{l}
\natisprodfib[\modsym](l)<m_1+m_2<>m_i>:\generderspace[\modsym](l)+\argcompl{m_1+m_2}+\rightarrow\generderspace[\modsym](l)+m_i+
\hspace{4pt}\text{by}\hspace{4pt}\text{setting:} \\
\hspace{100pt}\natisprodfib[\modsym](l)<m_1+m_2<>m_i>=
\gendiff[\modsym](l)+\lininclquot\left(\gczsfx[\proj<\mathbb{R}^{m_1},\mathbb{R}^{m_2}<>i>] \right)+\text{;}
\end{array}
\right.\\[8pt]
&\left\{
\begin{array}{l}
\natisprodfibinv[\modsym](l)<m_i<>\argcompl{m_1+m_2}>:\generderspace[\modsym](l)+m_i+\rightarrow\generderspace[\modsym](l)+\argcompl{m_1+m_2}+
\hspace{4pt}\text{by}\hspace{4pt}\text{setting:}\\[4pt]
\hspace{5pt}\natisprodfibinv[\modsym](l)<m_i<>m_1+m_2>=\\[4pt]
\hspace{52pt}
\left\{
\begin{array}{ll}
\gendiff[\modsym](l)+\argcompl{\lininclquot\left(\gczsfx[\argcompl{\left(\left(\idobj+\mathbb{R}^{m_1}+,\cost<\mathbb{R}^{m_1}<>\mathbb{R}^{m_2}>+0+\right)\funcomp\Diag{\mathbb{R}^{m_1}}{2}\right)}] \right)}+&\text{if}\;i=1\text{,}\\[4pt]
\gendiff[\modsym](l)+\argcompl{\lininclquot\left(\gczsfx[\argcompl{\left(\left(\cost<\mathbb{R}^{m_2}<>\mathbb{R}^{m_1}>+0+,\idobj+\mathbb{R}^{m_2}+\right)\funcomp\Diag{\mathbb{R}^{m_2}}{2}\right)}] \right)}+&\text{if}\;i=2\text{.}
\end{array}
\right.
\end{array}
\right.
\end{flalign*}
\end{definition}

\begin{proposition}\label{produno}\mbox{}
\begin{enumerate}
\item Fix $l,m_1,m_2\in\mathbb{N}_0$, $i \in \left\{1,2\right\}$, a $\mathbb{R}$-vector subspace $\modsym$ of $\mathbb{R}^l$. Then:
\begin{flalign*}
&\begin{array}{l}
\natisprodfib[\modsym](l)<m_1+m_2<>m_i>\funcomp\natisprodfibinv[\modsym](l)<m_i<>m_1+m_2>=\idobj+\argcompl{\generderspace[\modsym](l)+m_i+}+\text{;}
\end{array}\\[8pt]
&\begin{array}{l}
\natisprodfib[\modsym](l)<m_1+m_2<>m_i>\left(\generderspacesm[\modsym](l)+\argcompl{m_1+m_2}+\right)\subseteq\generderspacesm[\modsym](l)+m_i+\text{;}
\end{array}\\[8pt]
&\begin{array}{l}
\natisprodfibinv[\modsym](l)<m_i<>m_1+m_2>\left(\generderspacesm[\modsym](l)+m_i+\right)\subseteq\generderspacesm[\modsym](l)+\argcompl{m_1+m_2}+
\text{.}
\end{array}
\end{flalign*}
\item Fix $l, m_1,m_2, n_1, n_2\in\mathbb{N}_0$, $i \in \left\{1,2\right\}$, a $\mathbb{R}$-vector subspace $\modsym$ of $\mathbb{R}^l$,  $\gggfx_1 \in \genfquotcccz(m_1)+n_1+$, $\gggfx_2 \in \genfquotcccz(m_2)+n_2+$. Then:
\begin{flalign*}
&\begin{array}{l}
\gendiff[\modsym](l)+\gggfx_1+
\funcomp \natisprodfib[\modsym](l)<m_1+m_2<>m_i>=\\[4pt]
\hspace{136pt}
\natisprodfib[\modsym](l)<n_1+n_2<>n_i>\funcomp
\gendiff[\modsym](l)+\lboundquotpoint \gggfx_1, \gggfx_2\rboundquotpoint+\text{;}
\end{array}\\[8pt]
&\begin{array}{l}
\gendiff[\modsym](l)+\lboundquotpoint \gggfx_1, \gggfx_2\rboundquotpoint+
\funcomp\natisprodfibinv[\modsym](l)<m_i<>m_1+m_2>=\\[4pt]
\hspace{166pt}
\natisprodfibinv[\modsym](l)<n_i<>n_1+n_2>\funcomp 
\gendiff[\modsym](l)+\gggfx_i+
\text{.}
\end{array}
\end{flalign*}
\item Fix $l,m_1,m_2\in\mathbb{N}_0$, $i \in \left\{1,2\right\}$, two $\mathbb{R}$-vector subspaces $\modsym_1$, $\modsym_2$ of $\mathbb{R}^l$. Then:
\begin{flalign*}
&\begin{array}{l}
\natisprodfib[\modsym_2](l)<m_1+m_2<>m_i>\funcomp \generincl<\modsym_1<>\modsym_2>(l)+\argcompl{m_1+m_2}+ =\\[4pt]
 \hspace{137pt}\generincl<\modsym_1<>\modsym_2>(l)+m_i+ \funcomp\natisprodfib[\modsym_1](l)<m_1+m_2<>m_i>\text{;}
\end{array}\\[8pt]
&\begin{array}{l}
\natisprodfibinv[\modsym_2](l)<m_i<>m_1+m_2>\funcomp \generincl<\modsym_1<>\modsym_2>(l)+\argcompl{m_i}+   =\\[4pt]
 \hspace{110pt}   \generincl<\modsym_1<>\modsym_2>(l)+\argcompl{m_1+m_2}+ \funcomp\natisprodfibinv[\modsym_1](l)<m_i<>m_1+m_2>\text{.}
\end{array}
\end{flalign*}
\end{enumerate}
\end{proposition}
\begin{proof} Statement follows by direct computation.
\end{proof}

\begin{remark}\label{smisomdim}\mbox{}
\begin{enumerate}
\item $\natisprodfib[\modsym](l)<m_1+m_2<>m_i>$ extends to generalized tangent vectors the corresponding natural arrow $\natisprodfibsm<m_1+m_2<>m_i>$ which is well known in the smooth setting.\newline
$\natisprodfibinv[\modsym](l)<m_i<>m_1+m_2>$ extends to generalized tangent vectors the corresponding natural arrow $\natisprodfibinvsm<m_i<>m_1+m_2>$ which is well known in the smooth setting.
\item $\natisprodfibsm<m_1+m_2<>m_i>$ and $\natisprodfibinvsm<m_i<>m_1+m_2>$ together define a bi-product diagram where we refer to \cite{SML}\;\;Chap\,8\;\;\textsection\,2. 
\end{enumerate}
\end{remark}

\begin{proposition}\label{produnosm} Fix $l,m_1,m_2\in\mathbb{N}_0$, $i \in \left\{1,2\right\}$, a $\mathbb{R}$-vector subspace $\modsym$ of $\mathbb{R}^l$. Then:
\begin{flalign*}
&\begin{array}{l}
\quotsm[\modsym](l)+m_i+\funcomp\natisprodfib[\modsym](l)<m_1+m_2<>m_i>=\natisprodfibsm<m_1+m_2<>m_i>\funcomp \quotsm[\modsym](l)+\argcompl{m_1+m_2}+\text{;}
\end{array}\\[8pt]
&\begin{array}{l}
\quotsm[\modsym](l)+\argcompl{m_1+m_2}+\funcomp\natisprodfibinv[\modsym](l)<m_i<>m_1+m_2>=
\natisprodfibinvsm<m_1+m_2<>m_i>\funcomp\quotsm[\modsym](l)+m_i+\text{.}
\end{array}
\end{flalign*}
\end{proposition}
\begin{proof} Statement follows by direct computation.
\end{proof}

In Remark \ref{siterem2} below we show the necessity of the Grothendieck topology for spaces defined above.

\begin{remark}\label{siterem2}
Here we show by an example that spaces of derivations of degree $1$ cannot be endowed with a classical topology if we want that Theorem \ref{mainth} makes sense.  In fact we need that for any $m,n \in \mathbb{N}_0$, $\setsymcinque \subseteqdentro\mathbb{R}^m$, $\setsymsei \subseteqdentro\mathbb{R}^n$, continuous function $f:\setsymcinque \rightarrow \setsymsei$ the set function associating to any $\unkuno\in\setsymcinque$ the generalized differential of $f$ at $\unkuno$ must be a continuous function. This is equivalent to ask that functions described in Remark \ref{siterem}-[5] are continuous functions. Any classical topology on spaces of derivations of degree $1$ and compatible with continuity condition of functions described in Remark \ref{siterem}-[5] has the drawback described below.\newline
Fix $l,m \in \mathbb{N}$, $\setsymcinque \subseteqdentro\mathbb{R}^l$, $\elsymcinque_1, \elsymcinque_2 \in \setsymcinque$, a smooth function $f:\setsymcinque \rightarrow \mathbb{R}^m$. Assume that there is an open neighborhood $\setsymcinque_i\subseteq \setsymcinque$ of $\elsymcinque_i$ for any $i \in \left\{1,2\right\}$ such that $\closure\left[\setsymcinque_1,\mathbb{R}^l\right]\cap\closure\left[\setsymcinque_2,\mathbb{R}^l\right]=\udenset$ and $f\funcomp\incl{\setsymcinque_1}{\setsymcinque}=f\funcomp\left(\incl{\setsymcinque_2}{\setsymcinque}\funcomp\trasl[\elsymcinque_2-\elsymcinque_1]\right)$. Fix a path $\pmRn:\left(-1,1\right)\rightarrow \setsymcinque_1$ with $\pmRn\left(0\right)=\elsymcinque_1$ and define:
\begin{flalign*}
&\begin{array}{l}
\text{the set function}\;\pmRn_1:\left(-1,1\right)\rightarrow \setsymcinque\;\text{by setting}\\[4pt]
\pmRn_1\left(\unkdue\right)=
\left\{
\begin{array}{ll}
\pmRn\left(\unkdue\right)&\text{if}\;\unkdue\leq 0\text{,}\\
\left(\trasl[\elsymcinque_2-\elsymcinque_1]\funcomp\pmRn\right)\left(\unkdue\right)&\text{if}\;\unkdue> 0\text{;}
\end{array}
\right.
\end{array}\\[8pt]
&\begin{array}{l} 
\text{the set function}\;\pmdspreder:\left(-1,1\right)\rightarrow \dirsumgenpderspace[\mathbb{R}^l](l)+m+\;\text{by setting}\\[4pt]
\pmdspreder\left(\unkdue\right)=\genpreder\left[
\lininclquot\left(\gczsfx[\argcompl{
\trasl[-f\left(\pmRn\left(\unkdue\right)\right)] \genfuncomp\left(f \genfuncomp\trasl[\pmRn\left(\unkdue\right)] \right)
}] \right)  , \ef_1 \right]\hspace{5pt}\forall \unkdue\in \left(-1,1\right)\text{;}
\end{array}\\[8pt]
&\begin{array}{l}
\text{the set function}\;\pmdspreder_1:\left(-1,1\right)\rightarrow \dirsumgenpderspace[\mathbb{R}^l](l)+m+\;\text{by setting}\\[4pt]
\pmdspreder_1\left(\unkdue\right)=
\genpreder\left[
\lininclquot\left(\gczsfx[\argcompl{
\trasl[-f\left(\pmRn_1\left(\unkdue\right)\right)] \genfuncomp\left(f \genfuncomp\trasl[\pmRn_1\left(\unkdue\right)] \right)
}] \right) , \ef_1 \right]\hspace{5pt}\forall \unkdue\in \left(-1,1\right)\text{.}
\end{array}
\end{flalign*}
Since functions $f \funcomp\trasl[\pmRn\left(\unkdue\right)]$, $f \funcomp\trasl[\pmRn_1\left(\unkdue\right)]$ coincides in a suitable neighborhood of $0$ for any $\unkdue\in \left(-1,1\right)$ we have that set functions $\pmdspreder$ and $\pmdspreder_1$ coincide, then any classical topology on $\dirsumgenpderspace[\mathbb{R}^l](l)+m+$ making $\pmdspreder$ into a continuous function will make $\pmdspreder_1$ a continuous function too.  
Hence, by referring to Definition \ref{defcatD}-[5], any classical topology on $\generderspace[\mathbb{R}^l](l)+m+$ which makes $\quotdirsumdue \funcomp\left(\quotdirsumuno\funcomp\pmdspreder\right)$ into a continuous function also makes $\quotdirsumdue\funcomp\left(\quotdirsumuno\funcomp\pmdspreder_1\right)$ into a continuous function.\newline           
Such a drawback is intrinsic to classical topologies since they are based on subsets and do not allow any control on functions which generate such subsets.
\end{remark}

In Notation \ref{varispder} below we introduce functors, natural arrows, objects and arrows which will be used in forthcoming sections.  We refer to \eqref{bGt2ter}-[$(iv)$], \cite{MM} Chap. III. \textsection\,5.

\begin{notation} \mbox{}\label{varispder}
\begin{enumerate}
\item Fix $l,m \in \mathbb{N}_0$, one among symbols $\genfquot$, $\genfquotcccz$, $\genfquotccczz$ and denote it by $\OAOS$.\newline
We denote by $\OAOS-functor-(l)+m+:\GTopcat\rightarrow \GVecR$ the presheaf defined by setting
\begin{flalign*}
&\begin{array}{l}
\OAOS-functor-(l)+m+\left(\setsymuno\right)=\homF+\GTopcat+\left(\setsymuno,\OAOS(l)+m+\right)\hspace{4pt}\text{for}\hspace{4pt}\text{any}\hspace{4pt}\text{object}\hspace{4pt}\setsymuno\hspace{4pt}\text{of}\hspace{4pt}\GTopcat\text{;}
\end{array}\\[8pt] 
&\begin{array}{l}
\OAOS-funtore-(l)+m+\left(f\right)=\homF+\GTopcat+\left(f,\idobj+\argcompl{\OAOS(l)+m+}+\right)\hspace{4pt}\text{for}\hspace{4pt}\text{any}\hspace{4pt}\text{arrow}\hspace{4pt}f\hspace{4pt}\text{of}\hspace{4pt}\GTopcat\text{.}
\end{array}
\end{flalign*} 
We denote by $\OAOS/+/-funtore-(l)+m+$ the presheaf obtained by applying the plus-construction to $\OAOS-functor-(l)+m+$.\newline
We denote by $\OAOS/++/-funtore-(l)+m+$ the sheaf obtained by applying the plus-construction twice to $\OAOS-functor-(l)+m+$.\newline
For any fixed $n\in \mathbb{N}_0$ we denote by $\OAOS-funtore-(l)+m+\left(n\right)$ the stalk at $0\in \mathbb{R}^n$ of the sheaf $\OAOS/++/-funtore-(l)+m+$.
We emphasize that $\OAOS-funtore-(l)+m+\left(n\right)$ is naturally isomorphic to the stalk of $\OAOS/++/-funtore-(l)+m+$ at any other point $\unkuno\in\mathbb{R}^n$, due to the structure of Euclidean topology of $\mathbb{R}^n$. Any external or internal algebraic operation induced on $\OAOS-funtore-(l)+m+\left(n\right)$ by $\OAOS-functor-(l)+m+$ will be denoted by the same symbol used to denote the corresponding operation on $\OAOS-functor-(l)+m+$, with an abuse of language.\newline
For any fixed $n\in \mathbb{N}_0$, $\intuno\subseteqdentro \mathbb{R}^n$, $\elsymuno \in \OAOS-funtore-(l)+m+\left(\intuno\right)$ we denote by $\grm[\elsymuno]$ the arrow belonging to $\homF+\GTopcat+\left(\intuno,\OAOS-funtore-(l)+m+\left(n\right)\right)$ defined by setting
$\grm[\elsymuno]\left(\unkuno\right)$ to be the germ of $\elsymuno$  at $\unkuno$  for any $\unkuno \in\intuno$.
\item Fix $l,m \in \mathbb{N}_0$, a $\mathbb{R}$-vector subspace $\modsym$ of $\mathbb{R}^l$, one among symbols $\dirsumgenpderspace$, $\nulldirsumgenpderspace$, $\genprederspace$, $\genprederspacesm$, $\nullpreder$, $\nullprederunodue$, $\generderspace$, $\generderspacesm$, $\kergenerderspacesm$, $\cokergenerderspace$ and denote it by $\OAOS$.\newline
We denote by $\OAOS-functor-[\modsym](l)+m+:\GTopcat\rightarrow \GVecR$ the presheaf defined by setting
\begin{flalign*}
&\begin{array}{l}
\OAOS-functor-[\modsym](l)+m+\left(\setsymuno\right)=\homF+\GTopcat+\left(\setsymuno,\OAOS[\modsym](l)+m+\right)\hspace{4pt}\text{for}\hspace{4pt}\text{any}\hspace{4pt}\text{object}\hspace{4pt}\setsymuno\hspace{4pt}\text{of}\hspace{4pt}\GTopcat\text{;}
\end{array}\\[8pt] 
&\begin{array}{l}
\OAOS-funtore-[\modsym](l)+m+\left(f\right)=\homF+\GTopcat+\left(f,\idobj+\argcompl{\OAOS[\modsym](l)+m+}+\right)\hspace{4pt}\text{for}\hspace{4pt}\text{any}\hspace{4pt}\text{arrow}\hspace{4pt}f\hspace{4pt}\text{of}\hspace{4pt}\GTopcat\text{.}
\end{array}
\end{flalign*} 
We denote by $\OAOS/+/-funtore-[\modsym](l)+m+$ the presheaf obtained by applying the plus-construction to $\OAOS-functor-[\modsym](l)+m+$.\newline
We denote by $\OAOS/++/-funtore-[\modsym](l)+m+$ the sheaf obtained by applying the plus-construction twice to $\OAOS-functor-[\modsym](l)+m+$.\newline
For any fixed $n\in \mathbb{N}_0$ we denote by $\OAOS-funtore-[\modsym](l)+m+\left(n\right)$ the stalk at $0\in \mathbb{R}^n$ of the sheaf $\OAOS/++/-funtore-[\modsym](l)+m+$. We emphasize that $\OAOS-funtore-[\modsym](l)+m+\left(n\right)$ is naturally isomorphic to the stalk of $\OAOS/++/-funtore-[\modsym](l)+m+$ at any other point $\unkuno\in\mathbb{R}^n$, due to the structure of Euclidean topology of $\mathbb{R}^n$. Any external or internal algebraic operation induced on $\OAOS-funtore-[\modsym](l)+m+\left(n\right)$ by $\OAOS-functor-[\modsym](l)+m+$ will be denoted by the same symbol used to denote the corresponding operation on $\OAOS-functor-[\modsym](l)+m+$, with an abuse of language.\newline
For any fixed $n\in \mathbb{N}_0$, $\intuno\subseteqdentro \mathbb{R}^n$, $\elsymuno \in \OAOS-funtore-[\modsym](l)+m+\left(\intuno\right)$ we denote by $\grm[\elsymuno]$ the arrow belonging to $\homF+\GTopcat+\left(\intuno,\OAOS-funtore-[\modsym](l)+m+\left(n\right)\right)$ defined by setting
$\grm[\elsymuno]\left(\unkuno\right)$ to be the germ of $\elsymuno$  at $\unkuno$  for any $\unkuno \in\intuno$.\newline
If $\modsym=\mathbb{R}^l$ then symbol $\modsym$ will be dropped in the notation.\newline
We set:
\begin{flalign*} 
&\OAOS+m+=\underset{l \in \mathbb{N}_0}{\bigdirsum}\OAOS(l)+m+\text{;}\\[4pt]
&\OAOS-funtore-+m+=\underset{l \in \mathbb{N}_0}{\bigdirsum}\OAOS-funtore-(l)+m+\text{;}\\[4pt]
&\OAOS/+/-funtore-+m+=\underset{l \in \mathbb{N}_0}{\bigdirsum}\OAOS/+/-funtore-(l)+m+\text{;}\\[4pt]
&\OAOS/++/-funtore-+m+=\underset{l \in \mathbb{N}_0}{\bigdirsum}\OAOS/++/-funtore-(l)+m+\text{;}\\[4pt]
&\OAOS-funtore-+m+\left(n\right)=\underset{l \in \mathbb{N}_0}{\bigdirsum}\OAOS-funtore-(l)+m+\left(n\right)\text{.}
\end{flalign*}
\item Fix $l, m_1,m_2 \in \mathbb{N}_0$, a $\mathbb{R}$-vector subspace $\modsym$ of $\mathbb{R}^l$, $\gggfx\in\genfquotcccz(m_1)+m_2+$, one among symbols $\dirsumgenprediff$, $\genprediff$, $\gendiff$ and denote it by $\OAAS$.\newline
We denote by $\OAAS-trnat-[\modsym](l)+\gggfx+$ the natural arrow defined by setting
\begin{equation*}
\OAAS-trnat-[\modsym](l)+\gggfx+\left(\setsymuno\right)=\homF+\GTopcat+\left(\idobj+\setsymuno+,\OAAS[\modsym](l)+\gggfx+\right)\hspace{4pt}\text{for}\hspace{4pt}\text{any}\hspace{4pt}\text{object}\hspace{4pt}\setsymuno\hspace{4pt}\text{of}\hspace{4pt}\GTopcat\text{.}
\end{equation*} 
We denote by $\OAAS/+/-trnat-[\modsym](l)+\gggfx+$ the natural arrow obtained by applying the plus-construction to $\OAAS-trnat-[\modsym](l)+\gggfx+$.\newline
We denote by $\OAAS/++/-trnat-[\modsym](l)+\gggfx+$ the natural arrow obtained by applying the plus-construction twice to $\OAAS-trnat-[\modsym](l)+\gggfx+$.\newline
For any fixed $n\in \mathbb{N}_0$ we denote by $\OAAS-trnat-[\modsym](l)+\gggfx+\left(n\right)$ the arrow induced on stalks at $0\in \mathbb{R}^n$ by the natural arrow $\OAAS/++/-trnat-[\modsym](l)+\gggfx+$.\newline
Definition \ref{specfunVCC}, Proposition \ref{genpointprop}-[1] entail that $\OAAS-trnat-[\modsym](l)+\gggfx+\left(n\right)$ takes elements represented by at least one special set function to elements represented by at least one special set function.\newline
If $\modsym=\mathbb{R}^l$ then symbol $\modsym$ will be dropped in the notation.\newline
We set:
\begin{flalign*} 
&\OAAS+\gggfx+=\underset{l \in \mathbb{N}_0}{\bigdirsum}\OAAS(l)+\gggfx+\text{;}\\[4pt]
&\OAAS-funtore-+\gggfx+=\underset{l \in \mathbb{N}_0}{\bigdirsum}\OAAS-funtore-(l)+\gggfx+\text{;}\\[4pt]
&\OAAS/+/-funtore-+\gggfx+=\underset{l \in \mathbb{N}_0}{\bigdirsum}\OAAS/+/-funtore-(l)+\gggfx+\text{;}\\[4pt]
&\OAAS/++/-funtore-+\gggfx+=\underset{l \in \mathbb{N}_0}{\bigdirsum}\OAAS/++/-funtore-(l)+\gggfx+\text{;}\\[4pt]
&\OAAS-funtore-+\gggfx+\left(n\right)=\underset{l \in \mathbb{N}_0}{\bigdirsum}\OAAS-funtore-(l)+\gggfx+\left(n\right)\text{.}
\end{flalign*} 
\item Fix $l,m \in \mathbb{N}_0$, a $\mathbb{R}$-vector subspace $\modsym$ of $\mathbb{R}^l$, one among symbols $\quotdirsumuno$, $\quotdirsumdue$, $\inclgdsmgd$, $\kgdssmar$, $\ckgdssmar$ and denote it by $\OAAS$.\newline
We denote by $\OAAS-trnat-[\modsym](l)+m+$ the natural arrow defined by setting
\begin{equation*}
\OAAS-trnat-[\modsym](l)+m+\left(\setsymuno\right)=\homF+\GTopcat+\left(\idobj+\setsymuno+,\OAAS[\modsym](l)+m+\right)\hspace{4pt}\text{for}\hspace{4pt}\text{any}\hspace{4pt}\text{object}\hspace{4pt}\setsymuno\hspace{4pt}\text{of}\hspace{4pt}\GTopcat\text{.}
\end{equation*} 
We denote by $\OAAS/+/-trnat-[\modsym](l)+m+$ the natural arrow obtained by applying the plus-construction to $\OAAS-trnat-[\modsym](l)+m+$.\newline
We denote by $\OAAS/++/-trnat-[\modsym](l)+m+$ the natural arrow obtained by applying the plus-construction twice to $\OAAS-trnat-[\modsym](l)+m+$.\newline
For any fixed $n\in \mathbb{N}_0$ we denote by $\OAAS-trnat-[\modsym](l)+m+\left(n\right)$ the arrow induced on stalks at $0\in \mathbb{R}^n$ by the natural arrow $\OAAS/++/-trnat-[\modsym](l)+m+$.\newline
Arrow $\OAAS-trnat-[\modsym](l)+m+\left(n\right)$ takes elements represented by at least one special set function to elements represented by at least one special set function.\newline
If $\modsym=\mathbb{R}^l$ then symbol $\modsym$ will be dropped in the notation.\newline
We set:
\begin{flalign*} 
&\OAAS+m+=\underset{l \in \mathbb{N}_0}{\bigdirsum}\OAAS(l)+m+\text{;}\\[4pt]
&\OAAS-funtore-+m+=\underset{l \in \mathbb{N}_0}{\bigdirsum}\OAAS-funtore-(l)+m+\text{;}\\[4pt]
&\OAAS/+/-funtore-+m+=\underset{l \in \mathbb{N}_0}{\bigdirsum}\OAAS/+/-funtore-(l)+m+\text{;}\\[4pt]
&\OAAS/++/-funtore-+m+=\underset{l \in \mathbb{N}_0}{\bigdirsum}\OAAS/++/-funtore-(l)+m+\text{;}\\[4pt]
&\OAAS-funtore-+m+\left(n\right)=\underset{l \in \mathbb{N}_0}{\bigdirsum}\OAAS-funtore-(l)+m+\left(n\right)\text{.}
\end{flalign*} 
\item Fix $l, m,n \in \mathbb{N}_0$, a $\mathbb{R}$-vector subspace $\modsym$ of $\mathbb{R}^l$, $\gczsfx \in 
\Ckspquot{\infty}(m)+n+$, one among symbols $\gendiffsm$, $\kersqA$, $\cksqB$ and denote it by $\OAAS$.\newline
We denote by $\OAAS-trnat-[\modsym](l)+\lininclquot\left(\gczsfx\right)+$ the natural arrow defined by setting
\begin{equation*}
\OAAS-trnat-[\modsym](l)+\lininclquot\left(\gczsfx\right)+\left(\setsymuno\right)=\homF+\GTopcat+\left(\idobj+\setsymuno+,\OAAS[\modsym](l)+\lininclquot\left(\gczsfx\right)+\right)\hspace{4pt}\text{for}\hspace{4pt}\text{any}\hspace{4pt}\text{object}\hspace{4pt}\setsymuno\hspace{4pt}\text{of}\hspace{4pt}\GTopcat\text{.}
\end{equation*} 
We denote by $\OAAS/+/-trnat-[\modsym](l)+\lininclquot\left(\gczsfx\right)+$ the natural arrow obtained by applying the plus-construction to $\OAAS-trnat-[\modsym](l)+\lininclquot\left(\gczsfx\right)+$.\newline
We denote by $\OAAS/++/-trnat-[\modsym](l)+\lininclquot\left(\gczsfx\right)+$ the natural arrow obtained by applying the plus-construction twice to $\OAAS-trnat-[\modsym](l)+\lininclquot\left(\gczsfx\right)+$.\newline 
For any fixed $n\in \mathbb{N}_0$ we denote by $\OAAS-trnat-[\modsym](l)+\lininclquot\left(\gczsfx\right)+\left(n\right)$ the arrow induced on stalks at $0\in \mathbb{R}^n$ by the natural arrow $\OAAS/++/-trnat-[\modsym](l)+\lininclquot\left(\gczsfx\right)+$.\newline
Definition \ref{specfunVCC}, Proposition \ref{genpointprop}-[1] entail that $\OAAS-trnat-[\modsym](l)+\lininclquot\left(\gczsfx\right)+\left(n\right)$ takes elements represented by at least one special set function to elements represented by at least one special set function.\newline
If $\modsym=\mathbb{R}^l$ then symbol $\modsym$ will be dropped in the notation.\newline
We set:
\begin{flalign*} 
&\OAAS+\lininclquot\left(\gczsfx\right)+=\underset{l \in \mathbb{N}_0}{\bigdirsum}\OAAS(l)+\lininclquot\left(\gczsfx\right)+\text{;}\\[4pt]
&\OAAS-funtore-+\lininclquot\left(\gczsfx\right)+=\underset{l \in \mathbb{N}_0}{\bigdirsum}\OAAS-funtore-(l)+\lininclquot\left(\gczsfx\right)+\text{;}\\[4pt]
&\OAAS/+/-funtore-+\lininclquot\left(\gczsfx\right)+=\underset{l \in \mathbb{N}_0}{\bigdirsum}\OAAS/+/-funtore-(l)+\lininclquot\left(\gczsfx\right)+\text{;}\\[4pt]
&\OAAS/++/-funtore-+\lininclquot\left(\gczsfx\right)+=\underset{l \in \mathbb{N}_0}{\bigdirsum}\OAAS/++/-funtore-(l)+\lininclquot\left(\gczsfx\right)+\text{;}\\[4pt]
&\OAAS-funtore-+\lininclquot\left(\gczsfx\right)+\left(n\right)=\underset{l \in \mathbb{N}_0}{\bigdirsum}\OAAS-funtore-(l)+\lininclquot\left(\gczsfx\right)+\left(n\right)\text{.}
\end{flalign*}
\item Fix $l,m \in \mathbb{N}_0$, $\mathbb{R}$-vector subspaces $\modsym_1$, $\modsym_2$ of $\mathbb{R}^l$, one among symbols $\generincl$,  $\generinclsm$, $\kergenerinclsm$, $\ckgenerincl$ and denote it by $\OAASdue$.\newline
We denote by $\OAASdue-trnat-<\modsym_1<>\modsym_2>(l)+m+$ the natural arrow defined by setting
\begin{equation*}
\OAASdue-trnat-<\modsym_1<>\modsym_2>(l)+m+\left(\setsymuno\right)=\homF+\GTopcat+\left(\idobj+\setsymuno+,\OAASdue<\modsym_1<>\modsym_2>(l)+m+\right)\hspace{4pt}\text{for}\hspace{4pt}\text{any}\hspace{4pt}\text{object}\hspace{4pt}\setsymuno\hspace{4pt}\text{of}\hspace{4pt}\GTopcat\text{.}
\end{equation*} 
We denote by $\OAASdue/+/-trnat-<\modsym_1<>\modsym_2>(l)+m+$ the natural arrow obtained by applying the plus-construction to $\OAASdue-trnat-<\modsym_1<>\modsym_2>(l)+m+$.\newline
We denote by $\OAASdue/++/-trnat-<\modsym_1<>\modsym_2>(l)+m+$ the natural arrow obtained by applying the plus-construction twice to $\OAASdue-trnat-<\modsym_1<>\modsym_2>(l)+m+$.\newline
For any fixed $n\in \mathbb{N}_0$ we denote by $\OAASdue-trnat-<\modsym_1<>\modsym_2>(l)+m+\left(n\right)$ the arrow induced on stalks at $0\in \mathbb{R}^n$ by the natural arrow $\OAASdue/++/-trnat-<\modsym_1<>\modsym_2>(l)+m+$.\newline
Arrow $\OAASdue-trnat-<\modsym_1<>\modsym_2>(l)+m+\left(n\right)$ takes elements represented by at least one special set function to elements represented by at least one special set function.\newline
\item Fix $l,m, n\in \mathbb{N}_0$, an $\mathbb{R}$-vector subspace $\modsym$ of $\mathbb{R}^l$.\newline 
We denote by $\quotsm-trnat-[\modsym](l)+m+$ the natural arrow defined by setting
\begin{equation*}
\quotsm-trnat-[\modsym](l)+m+\left(\setsymuno\right) =\homF+\GTopcat+\left(\idobj+\setsymuno+,\quotsm[\modsym](l)+m+\right)\hspace{4pt}\text{for}\hspace{4pt}\text{any}\hspace{4pt}\text{object}\hspace{4pt}\setsymuno\hspace{4pt}\text{of}\hspace{4pt}\GTopcat\text{.}
\end{equation*}
We denote by $\quotsm/+/-trnat-[\modsym](l)+m+$ the natural arrow obtained by applying the plus-construction to $\quotsm-trnat-[\modsym](l)+m+$.\newline
We denote by $\quotsm/++/-trnat-[\modsym](l)+m+$ the natural arrow obtained by applying the plus-construction twice to $\quotsm-trnat-[\modsym](l)+m+$.\newline
For any fixed $n\in \mathbb{N}_0$ we denote by $\quotsm-trnat-[\modsym](l)+m+\left(n\right)$ the arrow induced on stalks at $0\in \mathbb{R}^n$ by the natural arrow $\quotsm/++/-trnat-[\modsym](l)+m+$.\newline
If $\modsym=\mathbb{R}^l$ then symbol $\modsym$ will be dropped in the notation.\newline
We set:
\begin{flalign*} 
&\quotsm+m+=\underset{l \in \mathbb{N}_0}{\bigdirsum}\quotsm(l)+m+\text{;}\\[4pt]
&\quotsm-trnat-+m+=\underset{l \in \mathbb{N}_0}{\bigdirsum}\quotsm-trnat-(l)+m+\text{;}\\[4pt]
&\quotsm/+/-trnat-+m+=\underset{l \in \mathbb{N}_0}{\bigdirsum}\quotsm/+/-trnat-(l)+m+\text{;}\\[4pt]
&\quotsm/++/-trnat-+m+=\underset{l \in \mathbb{N}_0}{\bigdirsum}\quotsm/++/-trnat-(l)+m+\text{;}\\[4pt]
&\quotsm-trnat-+m+\left(n\right)=\underset{l \in \mathbb{N}_0}{\bigdirsum}\quotsm-trnat-(l)+m+\left(n\right)\text{.}
\end{flalign*}
\end{enumerate}
\end{notation}

\begin{remark}\label{hogdorem}\mbox{}
\begin{enumerate}
\item We emphasize that if Grothendieck topology $\csGT$ on $\GTopcat$ is subcanonical then functors $\OAOS-functor-[\modsym](l)+m+$, $\OAOS/+/-functor-[\modsym](l)+m+$, $\OAOS/++/-functor-[\modsym](l)+m+$ introduced in Notation \ref{varispder} all define the same sheaf.
\item Definition \ref{alssitesdef}, Remark \ref{alssitesrem} entail that the well known construction of tensor, exterior and symmetric algebras starting from a given vector space, or a given sheaf of vector spaces, specializes to the case of generalized derivations.
More precisely fix one of the functors $\tenalgfun$, $\extalgfun$, $\symalgfun$ and denote it by $\genfunc$. By applying $\genfunc$ to $\mathbb{R}$-vector spaces, arrows, presheaves and natural arrows introduced up to now in this section we get the corresponding $\mathbb{R}$-algebras, arrows and related results in $\GTopcat$. If we fix $k \in \mathbb{N}_0$ and we apply functor $\genfunc_k$ to vector spaces, arrows, presheaves and natural arrows introduced up to now in this section then we get the degree $k$ homogeneous component of algebras, arrows and related results provided by $\genfunc$.\newline
If $\genfunc=\symalgfun$ then we get the algebra of higher degree generalized derivations. If we fix $k \in \mathbb{N}_0$ and we apply functor $\symalgfun_k$ then we get degree $k$ generalized derivations. It is well known that $\symalgfun\left(\mathbb{R}^m\right)$ is the algebra $\symalgtvsm<\segnvar<>m>$ of derivations introduced in Notation \ref{difrealfunc}, then results of Propositions \ref{nattosm}, \ref{diagRvs}, \ref{produnosm}, Remarks \ref{diagRvsrem}, \ref{smisomdim}, provide relations among derivations and generalized derivations.
\item Notation \ref{varispder}-[1] can be repeated word by word both in case of germs of local smooth functions by replacing $\genfquot$, $\genfquotcccz$ respectively by $\Ckspquot{\infty}$, $\Ckspquotccz{\infty}$ and in case of germs of local continuous functions by replacing $\genfquot$, $\genfquotcccz$ respectively by $\Ckspquot{0}$, $\Ckspquotccz{0}$.
\end{enumerate}
\end{remark}

\chapter{Foundations of generalized calculus \label{bascalc}}

In this chapter we introduce the notion of calculus in category $\GTopcat$. Calculus is defined locally, as in category $\Topcat$, where locality is to be understood with respect to Grothendieck topology and derivatives are computed with respect to generalized derivations, then any arrow of $\GTopcat$ turns out to be locally smooth.\\[12pt]

\section{Corresponding arrows}
 
In Definition \ref{isoinddef} below we introduce the notion of corresponding arrows taking values in spaces of generalized germs.\newline
We refer to Definitions \ref{genpointdef}, \ref{genpointdefmnloc}, Proposition \ref{genpointprop}.

\begin{definition}\label{isoinddef}\mbox{} 
\begin{enumerate}
\item Fix $l,m,n \in \mathbb{N}_0$, interval $\intuno\sqsubseteqdentro \mathbb{R}^{m}$, $\unkuno\in\intuno$, $\pmV_1\in\homF+\GTopcat+\left(\intuno,\genfquot(l)+n+\right)$, $\pmV_2\in\homF+\GTopcat+\left(\intuno,\genfquot(\argcompl{l+m})+n+\right)$.\newline
Assume that there are an open neighborhood $\intuno_{\unkuno} \sqsubseteqdentro\intuno$ of $\unkuno\in \intuno$, an open neighborhood $\setsymdue_{\unkuno}$ of $0\in \mathbb{R}^l$, $\gfy_{\unkuno}\in \genf(\argcompl{\setsymdue \times \left(\trasl[-\unkuno]\left(\intuno_{\unkuno}\right)\right)})+\mathbb{R}^{n}+$ fulfilling both conditions below:
\begin{flalign}
&\pmV_1\left(\unktre\right)=\preffuncaap{\gfy}_{\unkuno}\left(0,\unktre-\unkuno\right)\compquotpoint \lininclquot\left(\gczsfx[\argcompl{\pointedincl<\argcompl{\left(\mathbb{R}^{l},0\right),\left(\mathbb{R}^{m},0\right)}<>1>}]\right)\hspace{15pt}\forall \unktre \in \intuno_{\unkuno}\text{;}\label{arrdef}\\[6pt]
&\pmV_2\left(\unktre\right)=\preffuncaap{\gfy}_{\unkuno}\left(0,\unktre-\unkuno\right)\hspace{140pt}\forall \unktre \in \intuno_{\unkuno}\text{.}\label{corrarrdef}
\end{flalign}
We say that: $\pmV_2$ corresponds to $\pmV_1$ at $\unkuno$ or that $\pmV_1$ and $\pmV_2$ are corresponding arrows at $\unkuno$ or that there is a correspondence between $\pmV_1$ and $\pmV_2$ at $\unkuno$; the triplet $\left(\setsymdue_{\unkuno},\intuno_{\unkuno}, \gfy_{\unkuno}\right)$ is a representation of the correspondence between $\pmV_1$ and $\pmV_2$ at $\unkuno$.
\item Fix $l,m,n \in \mathbb{N}_0$, interval $\intuno\sqsubseteqdentro \mathbb{R}^{m}$, arrows $\pmV_1\in\homF+\GTopcat+\left(\intuno,\genfquot(l)+n+\right)$, $\pmV_2\in\homF+\GTopcat+\left(\intuno,\genfquot(\argcompl{l+m})+n+\right)$.\newline
We say that $\pmV_2$ corresponds to $\pmV_1$ or that $\pmV_1$ and $\pmV_2$ are corresponding arrows or that there is a correspondence between $\pmV_1$ and $\pmV_2$ if and only if  $\pmV_2$ corresponds to $\pmV_1$ at $\unkuno$ for any $\unkuno \in \intuno$.\newline
\item Fix $l,m,n \in \mathbb{N}_0$, $\intuno\sqsubseteqdentro \mathbb{R}^{m}$, $\pmV\in\homF+\GTopcat+\left(\intuno,\genfquot(l)+n+\right)$. We define
\begin{equation*}
\ihtlre{\intuno}{l}{n}\left(\pmV\right)=\left\{ \corrarrre[\pmV]\in \homF+\GTopcat+\left(\intuno,\genfquot(l+m)+n+\right)\;:\; \corrarrre[\pmV] \hspace{4pt} \text{corresponds}\hspace{4pt}\text{to}\hspace{4pt}\pmV\right\}\text{.}
\end{equation*}
\end{enumerate}
\end{definition}

In Proposition \ref{dualnot01} below we study the structure of corresponding arrows.\newline  
We refer to Notations  \ref{magtwopartic}-[2], \ref{ins}-[1, 11], Definitions \ref{intdiffmon}, \ref{dualnot0}, Propositions \ref{genpointprop}, \ref{pointspecfun}, Remarks \ref{Cinfdentro}, \ref{siterem}, \ref{siteremdue}. 
 
\begin{proposition}\label{dualnot01}\mbox{}
\begin{enumerate}
\item Fix $l,m,n \in \mathbb{N}_0$, $\intuno\sqsubseteqdentro \mathbb{R}^{m}$, $\pmV\in\homF+\GTopcat+\left(\intuno,\genfquot(l)+n+\right)$, $\corrarrre[\pmV]\in\homF+\GTopcat+\left(\intuno,\genfquot(\argcompl{l+m})+n+\right)$.\newline
Then both statement below hold true:
\begin{flalign}
&\left\{
\begin{array}{l}
\corrarrre[\pmV]\left(\unkuno\right)  \compquotpoint\lininclquot\left(\gczsfx[\argcompl{\pointedincl<\argcompl{\left(\mathbb{R}^{l},0\right),\left(\mathbb{R}^{m},\unktre\right)}<>1>}] \right)=\pmV\left(\unkuno+\unktre\right)\\[4pt]
\text{for}\hspace{4pt}\text{any}\hspace{4pt} \unkuno\in \intuno\text{,}\hspace{4pt}\unktre\in \mathbb{R}^{m}\text{;}
\end{array}\label{starv}
\right.\\[8pt]
&\left\{
\begin{array}{l}
\corrarrre[\pmV]\left(\unkuno\right)  \compquotpoint\lininclquot\left(\gczsfx[\argcompl{\pointedincl<\argcompl{\left(\mathbb{R}^{l},0\right),\left(\mathbb{R}^{m},\unktre_1\vecsum\unktre_2\right)}<>1>}] \right)=\\[4pt]
\hspace{70pt}\corrarrre[\pmV]\left(\unkuno\vecsum\unktre_1\right)\compquotpoint\lininclquot\left(\gczsfx[\argcompl{\pointedincl<\argcompl{\left(\mathbb{R}^{l},0\right),\left(\mathbb{R}^{m},\unktre_2\right)}<>1>}] \right)\\[4pt]
\text{for}\hspace{4pt}\text{any}\hspace{4pt} \unkuno\in \intuno\text{,}\hspace{4pt}\unktre_1,\unktre_2\in \mathbb{R}^{m}\text{.}
\end{array}\label{starstarv}
\right.
\end{flalign} 
\item Fix $l,m,n \in \mathbb{N}_0$, $\intuno\sqsubseteqdentro \mathbb{R}^{m}$, $\gfx \in\genf$, $\pmV\in\homF+\GTopcat+\left(\intuno,\genfquot(l)+n+\right)$, $\corrarrre[\pmV]$ corresponding to $\pmV$. Then $\gfx \compgenfquot\corrarrre[\pmV]$ corresponds to $\gfx \compgenfquot\pmV$.
\item Fix $l_1,l_2,m_1,m_2,n_1,n_2 \in \mathbb{N}_0$, intervals $\intuno_1\sqsubseteqdentro \mathbb{R}^{m_1}$, $\intuno_2\sqsubseteqdentro \mathbb{R}^{m_2}$, $\pmV_1\in\homF+\GTopcat+\left(\intuno_1,\genfquot(l_1)+n_1+\right)$, $\corrarrre[\pmV]_{1}$ corresponding to $\pmV_1$, $\pmV_2\in\homF+\GTopcat+\left(\intuno_2,\genfquot(l_2)+n_2+\right)$, $\corrarrre[\pmV]_{2}$ corresponding to $\pmV_2$.\newline
Define:
\begin{flalign*}
&\left\{
\begin{array}{l}
h:\mathbb{R}^{l_1}\times\mathbb{R}^{l_2}\times\mathbb{R}^{m_1}\times\mathbb{R}^{m_2}\rightarrow\mathbb{R}^{l_1}\times\mathbb{R}^{m_1}\times\mathbb{R}^{l_2}\times\mathbb{R}^{m_2}\hspace{4pt}\text{by}\hspace{4pt} \text{setting}\\[4pt]
h\left(\unkuno_1,\unkuno_2,\unkuno_3,\unkuno_4\right)=\left(\unkuno_1,\unkuno_3,\unkuno_2,\unkuno_3\right)\\[4pt] 
\hspace{95pt}\forall \left(\unkuno_1,\unkuno_2,\unkuno_3,\unkuno_3\right)\in \mathbb{R}^{l_1}\times\mathbb{R}^{l_2}\times\mathbb{R}^{m_1}\times\mathbb{R}^{m_2}\text{;}
\end{array}
\right.\\[8pt]
&\left\{
\begin{array}{l}
\pmV_3\in\homF+\GTopcat+\left(\intuno_1\times\intuno_2,\genfquot(l_1+l_2)+\argcompl{n_1+n_2}+\right) \hspace{4pt}\text{by setting}\\[4pt]
\pmV_3\left(\unkuno_1,\unkuno_2\right)=\lboundquotpoint\pmV_1\left(\unkuno_1\right),\pmV_2\left(\unkuno_2\right)\rboundquotpoint
\hspace{15pt} \forall \left(\unkuno_1,\unkuno_2\right)\in \intuno_1\times\intuno_2\text{;}
\end{array}
\right.\\[8pt]
&\left\{
\begin{array}{l}\corrarrre[\pmV]_3\in\homF+\GTopcat+\left(\intuno_1\times\intuno_2,\genfquot(l_1+l_2+m_1+m_2)+\argcompl{n_1+n_2}+\right) \\[4pt]
\text{by setting}\\[4pt]
\corrarrre[\pmV]_3\left(\unkuno_1,\unkuno_2\right)=\lboundquotpoint\corrarrre[\pmV]_{1}\left(\unkuno_1\right),\corrarrre[\pmV]_{2}\left(\unkuno_2\right)\rboundquotpoint\compquotpoint\lininclquot\left(\gczsfx[h]\right)
\hspace{15pt}\forall \left(\unkuno_1,\unkuno_2\right)\in \intuno_1\times\intuno_2\text{.}
\end{array}
\right.
\end{flalign*} 
Then $\corrarrre[\pmV]_3$ corresponds to $\pmV_3$.
\item Fix $l,m_1,m_2,n \in \mathbb{N}_0$, $\intuno_1\sqsubseteqdentro \mathbb{R}^{m_1}$, $\intuno_2\sqsubseteqdentro \mathbb{R}^{m_2}$, $\unkuno_{1,0}\in \intuno_1$, $\unkuno_{2,0}\in \intuno_2$,  $f \in \homF+\GTopcat+\left(\intuno_1,\intuno_2\right)$ with $f\left(\unkuno_{1,0}\right)=\unkuno_{2,0}$, $\pmV_2\in\homF+\GTopcat+\left(\intuno_2,\genfquotcccz(l)+n+\right)$, an arrow $\corrarr[\pmV]_2$ corresponding to $\pmV_2$, a representation $\left(\setsymdue_{\unkuno_{2,0}},\intuno_{\unkuno_{2,0}}, \gfy_{\unkuno_{2,0}}\right)$ of the correspondence between $\pmV_2$ and $\corrarr[\pmV]_2$ at $\unkuno_{2,0}$, an open neighborhood $\intuno_{\unkuno_{1,0}}$ of $\unkuno_{1,0} \in \intuno_1$ with $f\left(\intuno_{\unkuno_{1,0}}\right)\compacont \intuno_{\unkuno_{2,0}}$, $\gfx \in \genfcontcont(\intuno_{\unkuno_{1,0}})+\mathbb{R}^{m_2}+$ with $\evalcomptwo\left(\gfx\right)=f\funcomp \incl{\intuno_{\unkuno_{1,0}}}{\intuno_{1}}$.\newline
Define:
\begin{flalign*}
&\left\{
\begin{array}{l}
\gfy_{\unkuno_{1,0}} \in \genfcontcont(\setsymdue_{\unkuno_{2,0}}\times\trasl[-\unkuno_{1,0}]\left(\intuno_{\unkuno_{1,0}}\right))+\mathbb{R}^{n}+\hspace{4pt}\text{by}\hspace{4pt}\text{setting}\\[4pt]
\gfy_{\unkuno_{1,0}}=\gfy_{\unkuno_{2,0}} \genfuncomp \lboundgenf \incl{\setsymdue_{\unkuno_{2,0}}}{\mathbb{R}^l}, \trasl[-\unkuno_{2,0}]  \genfuncomp  \gfx  \genfuncomp  \trasl[\unkuno_{1,0}]   \rboundgenf\text{;}
\end{array}
\right.\\[8pt]
&\left\{
\begin{array}{l}
\text{the}\hspace{4pt}\text{arrow}\hspace{4pt}\pmV_{1,0}\in \homF+\GTopcat+\left(\intuno_{1,0},\genfquot(l)+n+\right)\hspace{4pt}\text{by}\hspace{4pt}\text{setting}\\[4pt]
\pmV_{1,0}\left(\unkuno\right)=\preffuncaap{\gfy}_{\unkuno_{1,0}}\left(0,\unkuno-\unkuno_{1,0}\right)\compquotpoint \lininclquot\left(\gczsfx[\argcompl{\pointedincl<\argcompl{\left(\mathbb{R}^{l},0\right),\left(\mathbb{R}^{m_1},0\right)}<>1>}]\right)\hspace{10pt}\forall \unkuno \in \intuno_{\unkuno_{1,0}}\text{;}
\end{array}
\right.\\[8pt]
&\left\{
\begin{array}{l}
\text{the}\hspace{4pt}\text{arrow}\hspace{4pt}\corrarrre[\pmV]_{1,0}\in \homF+\GTopcat+\left(\intuno_{1,0},\genfquot(l+m_1)+n+\right)\hspace{4pt}\text{by}\hspace{4pt}\text{setting}\\[4pt]
\corrarrre[\pmV]_{1,0}\left(\unkuno\right)=\preffuncaap{\gfy}_{\unkuno_{1,0}}\left(0,\unkuno-\unkuno_{1,0}\right)\hspace{15pt}\forall \unkuno \in \intuno_{\unkuno_{1,0}}\text{.}
\end{array}
\right.
\end{flalign*}
Then:
\begin{flalign}
&\begin{array}{l}
\corrarrre[\pmV]_{1,0}\hspace{4pt} \text{corresponds}\hspace{4pt} \text{to}\hspace{4pt} \pmV_{1,0}\text{;}
\end{array}\nonumber\\[8pt]
&\left\{
\begin{array}{l}
\left(\setsymdue_{\unkuno_{2,0}},\intuno_{\unkuno_{1,0}}, \gfy_{\unkuno_{1,0}}\right)\hspace{4pt} \text{is}\hspace{4pt} \text{a}\hspace{4pt} \text{representation}\hspace{4pt} \text{of}\hspace{4pt} \text{the}\hspace{4pt} \text{correspondence}\\[4pt]
\text{between}\hspace{4pt} \pmV_{1,0}\hspace{4pt} \text{and}\hspace{4pt} \corrarrre[\pmV]_{1,0}\hspace{4pt} \text{at}\hspace{4pt} \unkuno_{1,0}\text{;}
\end{array}
\right.\nonumber\\[8pt]
&\begin{array}{l}
\evalcompquotcontcontzero\funcomp\pmV_{1,0}=\evalcompquotcontcontzero\funcomp\left(\pmV_2\funcomp\left( f\funcomp\incl{\intuno_{\unkuno_{1,0}}}{\mathbb{R}^{m_1}} \right)\right)\text{.}\label{ugualeval}
\end{array}
\end{flalign}
\item Fix $l,m,n \in \mathbb{N}_0$, $\bsfm\in\left\{\Fint_{\mathsf{l}}, \Fpart_{\mathsf{l}},\Fp_{\mathsf{l}}, \Fq_{\mathsf{l}},\Fqq_{\mathsf{l}}\;:\; \mathsf{l}\in \left\{1,...,l\right\}\right\}$, $\intuno\sqsubseteqdentro \mathbb{R}^{m}$, $\pmV\in\homF+\GTopcat+\left(\intuno,\genfquot(l)+n+\right)$,
$\corrarrre[\pmV]$ corresponding to $\pmV$.\newline
Then $\bsfm\bsfmuquotpoint\corrarrre[\pmV]$ corresponds to $\bsfm\bsfmuquotpoint\pmV$.
\end{enumerate}
\end{proposition}
\begin{proof}\mbox{}\newline
\textnormal{\textbf{Proof of statement 1.}}\ \ Fix points $\unkuno_0\in \intuno$, $\unktre, \unktre_1, \unktre_2\in \mathbb{R}^m$, 
a representation $\left(\setsymdue_{\unkuno_0},\intuno_{\unkuno_0}, \gfy_{\unkuno_0}\right)$ of the correspondence between $\pmV$ and $\corrarrre[\pmV]$ at $\unkuno_0$.\newline
Statement \eqref{starv} follows by the chain of equalities below
\begin{equation*}
\begin{array}{l}
\corrarrre[\pmV]\left(\unkuno\right)  \compquotpoint\lininclquot\left(\gczsfx[\argcompl{\pointedincl<\argcompl{\left(\mathbb{R}^{l},0\right),\left(\mathbb{R}^{m},\unktre\right)}<>1>}] \right)=\preffuncaap{\gfy}_{\unkuno_0}\left(0,\unkuno-\unkuno_0\vecsum\unktre\right)= \\[4pt]
\preffuncaap{\gfy}_{\unkuno_0}\left(0,\unkuno\vecsum\unktre-\unkuno_0\right)\compquotpoint\lininclquot\left(\gczsfx[\argcompl{\pointedincl<\argcompl{\left(\mathbb{R}^{l},0\right),\left(\mathbb{R}^{m},0\right)}<>1>}] \right)=\pmV\left(\unkuno+\unktre\right)\\[4pt]
\text{for}\hspace{4pt}\text{any}\hspace{4pt} \unkuno\in \intuno\text{,}\hspace{4pt}\unktre\in \mathbb{R}^{m}\text{.}
\end{array}
\end{equation*}
Statement \eqref{starstarv} follows by the chain of equalities below
\begin{equation*}
\begin{array}{l}
\corrarrre[\pmV]\left(\unkuno\right)  \compquotpoint\lininclquot\left(\gczsfx[\argcompl{\pointedincl<\argcompl{\left(\mathbb{R}^{l},0\right),\left(\mathbb{R}^{m},\unktre_1\vecsum\unktre_2\right)}<>1>}] \right)=\preffuncaap{\gfy}_{\unkuno_0}\left(0,\unkuno-\unkuno_0\vecsum\unktre_1\vecsum\unktre_2\right)= \\[4pt]
\preffuncaap{\gfy}_{\unkuno_0}\left(0,\unkuno\vecsum\unktre_1-\unkuno_0\right)\compquotpoint\lininclquot\left(\gczsfx[\argcompl{\pointedincl<\argcompl{\left(\mathbb{R}^{l},0\right),\left(\mathbb{R}^{m},\unktre_2\right)}<>1>}] \right)=\\[4pt]
\corrarrre[\pmV]\left(\unkuno\vecsum\unktre_1\right)\compquotpoint\lininclquot\left(\gczsfx[\argcompl{\pointedincl<\argcompl{\left(\mathbb{R}^{l},0\right),\left(\mathbb{R}^{m},\unktre_2\right)}<>1>}] \right)\\[4pt]
\text{for}\hspace{4pt}\text{any}\hspace{4pt} \unkuno\in \intuno\text{,}\hspace{4pt}\unktre_1,\unktre_2\in \mathbb{R}^{m}\text{.}
\end{array}
\end{equation*} 
\textnormal{\textbf{Proof of statement 2.}}\ \ Fix $\unkuno_0\in \intuno$, a representation $\left(\setsymdue_{\unkuno_0},\intuno_{\unkuno_0}, \gfy_{\unkuno_0}\right)$ of the correspondence between $\pmV$ and $\corrarrre[\pmV]$ at $\unkuno_0$.\newline
Proposition \ref{genpointprop}-[1] entails that
\begin{equation*}
\left(\gfx \compgenfquot\preffuncaap{\gfy}_{\unkuno_0}\right)\left(\vecuno,\unkuno\right)=\preffuncaap{\left(\gfx \genfuncomp \gfy_{\unkuno_0}\right)}\left(\vecuno,\unkuno\right)\hspace{15pt}\forall \left(\vecuno,\unkuno\right)\in\setsymdue \times \left(\trasl[-\unkuno_0]\left(\intuno_0\right)\right)\text{.}
\end{equation*}
Eventually statement follows since:
\begin{equation*}
\begin{array}{ll}
\gfx \compgenfquot\pmV\left(\unkuno\right)=\preffuncaap{\left(\gfx \genfuncomp \gfy_{\unkuno_0}\right)}\left(0,\unkuno-\unkuno_0\right)\compquotpoint \lininclquot\left(\gczsfx[\argcompl{\pointedincl<\argcompl{\left(\mathbb{R}^{l},0\right),\left(\mathbb{R}^{m},0\right)}<>1>}]\right)&\forall \unkuno \in \intuno_0\text{;}\\[6pt]
\gfx \compgenfquot\corrarrre[\pmV]\left(\unkuno\right)=\preffuncaap{\left(\gfx \genfuncomp \gfy_{\unkuno_0}\right)}\left(0,\unkuno-\unkuno_0\right)&\forall \unkuno \in \intuno_0\text{.}
\end{array}
\end{equation*}
\textnormal{\textbf{Proof of statement 3.}}\ \ Fix points $\unkuno_{1,0}\in \intuno_1$, $\unkuno_{2,0}\in \intuno_2$, a representation $\left(\setsymdue_{\unkuno_{1,0}},\intuno_{\unkuno_{1,0}}, \gfy_{\unkuno_{1,0}}\right)$ of the correspondence between $\pmV_1$ and $\corrarrre[\pmV]_1$ at $\unkuno_{1,0}$, a representation $\left(\setsymdue_{\unkuno_{2,0}},\intuno_{\unkuno_{2,0}}, \gfy_{\unkuno_{2,0}}\right)$ of the correspondence between $\pmV_2$ and $\corrarrre[\pmV]_2$ at $\unkuno_{2,0}$.\newline
Define:
\begin{flalign*}
&\begin{array}{l}
\setsymdue_{3}=\setsymdue_{\unkuno_{1,0}}\times\setsymdue_{\unkuno_{2,0}}\text{;}
\end{array}\\[8pt]
&\begin{array}{l}
\intuno_{3}=\intuno_{\unkuno_{1,0}}\times\intuno_{\unkuno_{2,0}}\text{;}
\end{array}\\[8pt]
&\left\{
\begin{array}{l}
\gfy_3\in \genf(\argcompl{\setsymdue_{3}\times \left(\trasl[-\left(\unkuno_{1,0},\unkuno_{2,0}\right)]\left(\intuno_{3}\right)\right)})+\argcompl{\mathbb{R}^{n_1+n_2}}+\hspace{4pt}\text{by}\hspace{4pt}\text{setting}\\[4pt]
\gfy_3=\lboundgenf\gfy_{\unkuno_{1,0}},\gfy_{\unkuno_{2,0}}\rboundgenf\genfuncomp\left(h\genfuncomp \incl{\setsymdue_3\times \left(\trasl[-\left(\unkuno_{1,0},\unkuno_{2,0}\right)]\left(\intuno_{3}\right)\right)}{\mathbb{R}^{l_1+l_2+m_1+m_2}}\right)\text{.}
\end{array}
\right.
\end{flalign*}
Eventually statement follows by proving that $\left(\setsymdue_{3}, \intuno_{3}, \gfy_3\right)$ is a representation of the correspondence between $\pmV_3$ and $\corrarrre[\pmV]_3$ at $\left(\unkuno_{1,0},\unkuno_{2,0}\right)$.\newline
This follows by the two chain of equalities below:
\begin{equation*}
\begin{array}{l}
\pmV_3\left(\unkuno_1,\unkuno_2\right)=\lboundquotpoint\pmV_1\left(\unkuno_1\right),\pmV_2\left(\unkuno_2\right)\rboundquotpoint=\\[4pt]
\lboundquotpoint\preffuncaap{\left(\gfy_{\unkuno_{1,0}}\right)}\left(0,\unkuno_1-\unkuno_{1,0}\right)\compquotpoint \lininclquot\left(\gczsfx[\argcompl{\pointedincl<\argcompl{\left(\mathbb{R}^{l_1},0\right),\left(\mathbb{R}^{m_1},0\right)}<>1>}]\right), \\[4pt]
\hspace{75pt}\preffuncaap{\left(\gfy_{\unkuno_{2,0}}\right)}\left(0,\unkuno_2-\unkuno_{2,0}\right)\compquotpoint \lininclquot\left(\gczsfx[\argcompl{\pointedincl<\argcompl{\left(\mathbb{R}^{l_2},0\right),\left(\mathbb{R}^{m_2},0\right)}<>1>}]\right)\rboundquotpoint=\\[4pt]
\preffuncaap{\lboundgenf\gfy_{\unkuno_{1,0}},\gfy_{\unkuno_{2,0}}\rboundgenf}\left(0,\unkuno_1-\unkuno_{1,0},0,\unkuno_2-\unkuno_{2,0}\right)\compquotpoint\\[4pt]
\hspace{50pt}
\lboundquotpoint \lininclquot\left(\gczsfx[\argcompl{\pointedincl<\argcompl{\left(\mathbb{R}^{l_1},0\right),\left(\mathbb{R}^{m_1},0\right)}<>1>}]\right), \lininclquot\left(\gczsfx[\argcompl{\pointedincl<\argcompl{\left(\mathbb{R}^{l_2},0\right),\left(\mathbb{R}^{m_2},0\right)}<>1>}]\right)\rboundquotpoint=\\[4pt]
\preffuncaap{\gfy_3}\left(0,0,\unkuno_1-\unkuno_{1,0},\unkuno_2-\unkuno_{2,0}\right)\compquotpoint\\[4pt]
\hspace{50pt}
\lininclquot\left(\gczsfx[\argcompl{\pointedincl<\argcompl{\left(\mathbb{R}^{l_1+l_2},0\right),\left(\mathbb{R}^{m_1+m_2},0\right)}<>1>}]\right)
\hspace{15pt}\forall \left(\unkuno_1,\unkuno_2\right)\in \intuno_1\times\intuno_2\text{;}\\[10pt]
\corrarrre[\pmV]_3\left(\unkuno_1,\unkuno_2\right)=\lboundquotpoint\corrarrre[\pmV]_{1}\left(\unkuno_1\right),\corrarrre[\pmV]_{2}\left(\unkuno_2\right)\rboundquotpoint\compquotpoint\lininclquot\left(\gczsfx[h]\right)=\\[4pt]
\lboundquotpoint\preffuncaap{\gfy_{\unkuno_{1,0}}}\left(0,\unkuno_1-\unkuno_{1,0}\right),\preffuncaap{\gfy_{\unkuno_{2,0}}}\left(0,\unkuno_2-\unkuno_{2,0}\right)\rboundquotpoint\compquotpoint\lininclquot\left(\gczsfx[h]\right)=\\[4pt]
\preffuncaap{\lboundgenf\gfy_{\unkuno_{1,0}},\gfy_{\unkuno_{2,0}}\rboundgenf}\left(0,\unkuno_1-\unkuno_{1,0},0,\unkuno_2-\unkuno_{2,0}\right)\compquotpoint\lininclquot\left(\gczsfx[h]\right)=
\\[4pt]
\preffuncaap{\gfy_3}\left(0,0,\unkuno_1-\unkuno_{1,0},\unkuno_2-\unkuno_{2,0}\right)\hspace{15pt}\forall \left(\unkuno_1,\unkuno_2\right)\in \intuno_1\times\intuno_2\text{.}
\end{array}
\end{equation*}
\textnormal{\textbf{Proof of statement 4.}}\ \ Statement follows straightforwardly by direct checking.\newline
\textnormal{\textbf{Proof of statement 5.}}\ \ Fix $\unkuno_0\in \intuno$, a representation $\left(\setsymdue_{\unkuno_0},\intuno_{\unkuno_0}, \gfy_{\unkuno_0}\right)$ of the correspondence between $\pmV$ and $\corrarrre[\pmV]$ at $\unkuno_0$.\newline
Referrign to Notation \ref{realvec}-[6] define:
\begin{flalign*}
&\begin{array}{l}
k=\left\{
\begin{array}{ll}
l&\text{if}\hspace{4pt}\bsfm\in\left\{\Fpart_{\mathsf{l}},\Fp_{\mathsf{l}}\right\}\hspace{4pt}\text{for}\hspace{4pt} \mathsf{l}\in \left\{1,...,l\right\}\text{,}\\[4pt]
l+1&\text{if}\hspace{4pt}\bsfm\in\left\{\Fint_{\mathsf{l}}, \Fq_{\mathsf{l}},\Fqq_{\mathsf{l}}\right\}\hspace{4pt}\text{for}\hspace{4pt} \mathsf{l}\in \left\{1,...,l\right\}\text{;}
\end{array}
\right.
\end{array}\\[8pt]
&\begin{array}{l}
\setsymtre_{\unkuno_0}=\left\{
\begin{array}{ll}
\setsymdue_{\unkuno_0}&\text{if}\hspace{4pt}\bsfm\in\left\{\Fpart_{\mathsf{l}},\Fp_{\mathsf{l}}\right\}\hspace{4pt}\text{for}\hspace{4pt} \mathsf{l}\in \left\{1,...,l\right\}\text{,}\\[4pt]
\domint\left[\setsymdue_{\unkuno_0},\mathsf{l}\right]&\text{if}\hspace{4pt}\bsfm\in\left\{\Fint_{\mathsf{l}}, \Fq_{\mathsf{l}},\Fqq_{\mathsf{l}}\right\}\hspace{4pt}\text{for}\hspace{4pt} \mathsf{l}\in \left\{1,...,l\right\}\text{;}
\end{array}
\right.
\end{array}\\[8pt]
&\begin{array}{l}
\gfz_{\unkuno_0}\in \genf(\argcompl{\setsymtre_{\unkuno_{0}}\times \left(\trasl[-\unkuno_{0}]\left(\intuno_{\unkuno_{0}}\right)\right)})+\mathbb{R}^{n}+\hspace{4pt}\text{by}\hspace{4pt}\text{setting}\hspace{4pt}\gfz_{\unkuno_0}=\bsfm\bsfmuquotpoint\gfy_{\unkuno_0}\text{.}
\end{array}
\end{flalign*}
Then statement follows by proving that  $\left(\setsymtre_{\unkuno_0},\intuno_{\unkuno_0}, \gfz_{\unkuno_0}\right)$ a representation of the correspondence between $\bsfm\bsfmuquotpoint\pmV$ and $\bsfm\bsfmuquotpoint\corrarrre[\pmV]$ at $\unkuno_0$.\newline 
This follows since $\bsfm$ acts on the first $l$ unknowns and by the two chain of equalities below:
\begin{flalign*}
&\begin{array}{l}
\bsfm\bsfmuquotpoint\left(\pmV\left(\unkuno\right)\right)=\bsfm\bsfmuquotpoint\left(\preffuncaap{\gfy}_{\unkuno_0}\left(0,\unkuno-\unkuno_0\right)\compquotpoint \lininclquot\left(\gczsfx[\argcompl{\pointedincl<\argcompl{\left(\mathbb{R}^{l},0\right),\left(\mathbb{R}^{m},0\right)}<>1>}]\right)\right)=\\[4pt]
\left(\bsfm\bsfmuquotpoint\preffuncaap{\gfy}_{\unkuno_0}\right)\left(0,\unkuno-\unkuno_0\right)\compquotpoint \lininclquot\left(\gczsfx[\argcompl{\pointedincl<\argcompl{\left(\mathbb{R}^{k},0\right),\left(\mathbb{R}^{m},0\right)}<>1>}]\right)\hspace{15pt}\forall \unkuno \in \intuno_{\unkuno_0}\text{;}
\end{array}\\[10pt]
&\begin{array}{l}
\bsfm\bsfmuquotpoint\left(\corrarrre[\pmV]\left(\unkuno\right)\right)=\bsfm\bsfmuquotpoint\left(\preffuncaap{\gfy}_{\unkuno_0}\left(0,\unkuno-\unkuno_0\right)\right)=
\left(\bsfm\bsfmuquotpoint\preffuncaap{\gfy}_{\unkuno_0}\right)\left(0,\unkuno-\unkuno_0\right)\hspace{15pt}\forall \unkuno \in \intuno_{\unkuno_0}\text{.}
\end{array}
\end{flalign*}
\end{proof}

In Definition \ref{dualnot1} below we introduce the notion of p-corresponding arrows taking values in spaces of continuous pointed generalized germs.\newline
We refer to Notation \ref{gentopnot}-[7], Definitions \ref{genpointdef}, \ref{genpointdefmnloc}, Proposition \ref{genpointprop}.

\begin{definition}\label{dualnot1} \mbox{} 
\begin{enumerate}
\item Fix $l,m,n \in \mathbb{N}_0$, $\intuno\sqsubseteqdentro \mathbb{R}^{m}$, $\unkuno\in\intuno$, $\pmV_1\in\homF+\GTopcat+\left(\intuno,\genfquotcccz(l)+n+\right)$, $\pmV_2\in\homF+\GTopcat+\left(\intuno,\genfquotcccz(\argcompl{l+m})+n+\right)$.\newline
Assume that there are an open neighborhood $\intuno_{\unkuno} \sqsubseteqdentro\intuno$ of $\unkuno\in \intuno$, an open neighborhood $\setsymdue_{\unkuno}$ of $0\in \mathbb{R}^l$, $\gfy_{\unkuno}\in \genf(\argcompl{\setsymdue \times \left(\trasl[-\unkuno]\left(\intuno_{\unkuno}\right)\right)})+\mathbb{R}^{n}+$ fulfilling both conditions below:
\begin{equation*}
\begin{array}{ll}
\pmV_1\left(\unktre\right)=\ffuncaap{\gfy}_{\unkuno}\left(0,\unktre-\unkuno\right)\compquotpoint \lininclquot\left(\gczsfx[\argcompl{\pointedincl<\argcompl{\left(\mathbb{R}^{l},0\right),\left(\mathbb{R}^{m},0\right)}<>1>}]\right)&\forall \unktre \in \intuno_{\unkuno}\text{;}\\[6pt]
\pmV_2\left(\unktre\right)=\ffuncaap{\gfy}_{\unkuno}\left(0,\unktre-\unkuno\right)&\forall \unktre \in \intuno_{\unkuno}\text{.}
\end{array}
\end{equation*}
We say: $\pmV_2$ p-corresponds to $\pmV_1$ at $\unkuno$ or that $\pmV_1$ and $\pmV_2$ are p-corresponding arrows at $\unkuno$ or that there is a p-correspondence between $\pmV_1$ and $\pmV_2$ at $\unkuno$; the triplet $\left(\setsymdue_{\unkuno},\intuno_{\unkuno}, \gfy_{\unkuno}\right)$ is a representation of the p-correspondence between $\pmV_1$ and $\pmV_2$ at $\unkuno$.
\item Fix $l,m,n \in \mathbb{N}_0$, interval $\intuno\sqsubseteqdentro \mathbb{R}^{m}$, arrows $\pmV_1\in\homF+\GTopcat+\left(\intuno,\genfquotcccz(l)+n+\right)$, $\pmV_2\in\homF+\GTopcat+\left(\intuno,\genfquotcccz(\argcompl{l+m})+n+\right)$.\newline
We say that $\pmV_2$ p-corresponds to $\pmV_1$ or that $\pmV_1$ and $\pmV_2$ are p-corresponding arrows or that there is a p-correspondence between $\pmV_1$ and $\pmV_2$ if and only if  $\pmV_2$ p-corresponds to $\pmV_1$ at $\unkuno$ for any $\unkuno \in \intuno$.
\item Fix $l,m,n \in \mathbb{N}_0$, $\intuno\sqsubseteqdentro \mathbb{R}^{m}$, $\pmV\in\homF+\GTopcat+\left(\intuno,\genfquotcccz(l)+n+\right)$. We define
\begin{equation*}
\ihtl{\intuno}{l}{n}\left(\pmV\right)=\left\{ \corrarr[\pmV]\in \homF+\GTopcat+\left(\intuno,\genfquotcccz(l+m)+n+\right)\;:\; \corrarr[\pmV] \hspace{4pt} \text{p-corresponds}\hspace{4pt}\text{to}\hspace{4pt}\pmV\right\}\text{.}
\end{equation*}
\end{enumerate}
\end{definition}

In Proposition \ref{preimmerfun} below we study the structure of p-corresponding arrows.\newline  
We refer to Notations  \ref{magtwopartic}-[2], \ref{ins}-[1, 11], Definitions \ref{intdiffmon}, \ref{dualnot0}, Propositions \ref{genpointprop}, \ref{pointspecfun}, Remarks \ref{Cinfdentro}, \ref{siterem}, \ref{siteremdue}.

\begin{proposition} \label{preimmerfun}\mbox{} 
\begin{enumerate}
\item Fix $l,m,n \in \mathbb{N}_0$, interval $\intuno\sqsubseteqdentro \mathbb{R}^{m}$, arrows $\pmV\in\homF+\GTopcat+\left(\intuno,\genfquotcccz(l)+n+\right)$, $\corrarr[\pmV]\in\homF+\GTopcat+\left(\intuno,\genfquotcccz(\argcompl{l+m})+n+\right)$.\newline
Then both statement below hold true: 
\begin{flalign}
&\left\{
\begin{array}{l}
\corrarr[\pmV]\left(\unkuno\right)  \compquotpoint\lininclquot\left(\gczsfx[\argcompl{\pointedincl<\argcompl{\left(\mathbb{R}^{l},0\right),\left(\mathbb{R}^{m},\unktre\right)}<>1>}] \right)  \sumquotpoint\\[4pt] 
\hspace{40pt}-1 \scalpquotpoint \left(\corrarr[\pmV]\left(\unkuno\right)  \compquotpoint\lininclquot\left(\gczsfx[\argcompl{\cost<\mathbb{R}^{l}<>\mathbb{R}^{l+m}>+\left(0, \unktre\right)+}]\right)\right)=\pmV\left(\unkuno+\unktre\right)\\[4pt]
\text{for}\hspace{4pt}\text{any}\hspace{4pt} \unkuno\in \intuno\text{,}\hspace{4pt}\unktre\in \mathbb{R}^{m}\text{;}
\end{array}\label{tildv}
\right.\\[8pt]
&\left\{
\begin{array}{l}
\corrarr[\pmV]\left(\unkuno\right)  \compquotpoint\lininclquot\left(\gczsfx[\argcompl{\pointedincl<\argcompl{\left(\mathbb{R}^{l},0\right),\left(\mathbb{R}^{m},\unktre_1\vecsum\unktre_2\right)}<>1>}] \right) \sumquotpoint\\[4pt] 
\hspace{80pt}-1 \scalpquotpoint
\corrarr[\pmV]\left(\unkuno\right)  \compquotpoint\lininclquot\left(\gczsfx[\argcompl{\cost<\mathbb{R}^{l}<>\mathbb{R}^{l+m}>+\left(0, \unktre_1\right)+}]\right)
=\\[4pt]
\corrarr[\pmV]\left(\unkuno\vecsum\unktre_1\right)\compquotpoint\lininclquot\left(\gczsfx[\argcompl{\pointedincl<\argcompl{\left(\mathbb{R}^{l},0\right),\left(\mathbb{R}^{m},\unktre_2\right)}<>1>}] \right)\\[4pt]
\text{for}\hspace{4pt}\text{any}\hspace{4pt} \unkuno\in \intuno\text{,}\hspace{4pt}\unktre_1,\unktre_2\in \mathbb{R}^{m}\text{.}
\end{array}\label{tildtildv}
\right.
\end{flalign} 
\item Fix $l,m,n \in \mathbb{N}_0$, $\intuno\sqsubseteqdentro \mathbb{R}^{m}$, $\gfx \in\genf$, $\pmV\in\homF+\GTopcat+\left(\intuno,\genfquotcccz(l)+n+\right)$, $\corrarr[\pmV]$ p-correspon\-ding to $\pmV$.\newline
Then $\left(\gfx \compgenfquot\corrarr[\pmV]\right) \sumquotpoint
\left(-1\scalpquotpoint 
 \gfx \compgenfquot \left(\corrarr[\pmV]\left(\unkuno\right)\compquotpoint 
\lininclquot\left(\gczsfx[\argcompl{\lboundquotpoint   \cost<\mathbb{R}^l<>\mathbb{R}^l>+0+, \idobj+\mathbb{R}^m+ \rboundquotpoint}]\right)  \right)\right)$ p-correspo\-nds to $\left(\gfx \compgenfquot\pmV\right) \sumquotpoint\left(-1\scalpquotpoint   \gczsfx[\argcompl{\gfx\genfuncomp\cost<\mathbb{R}^{l+m}<>\mathbb{R}^{n}>+0+}] \right)$.
\item Fix $l_1,l_2,m_1,m_2,n \in \mathbb{N}_0$, $\intuno_1\sqsubseteqdentro \mathbb{R}^{m_1}$, $\intuno_2\sqsubseteqdentro \mathbb{R}^{m_2}$,  $\unkuno_{1,0}\in \intuno_1$, $\unkuno_{2,0}\in \intuno_2$, $\pmV_1\in\homF+\GTopcat+\left(\intuno_1,\genfquotcccz(l_1)+n+\right)$, $\corrarrre[\pmV]_{1}\in \homF+\GTopcat+\left(\intuno_1,\genfquotcccz(l_1+m_1)+n+\right)$ corresponding to $\pmV_1$, $\pmV_2\in\homF+\GTopcat+\left(\intuno_2,\genfquotcccz(l_2)+l_1+\right)$, $\corrarr[\pmV]_{2}\in \homF+\GTopcat+\left(\intuno_2,\genfquotcccz(l_2+m_2)+l_1+\right)$ p-corresponding to $\pmV_2$, a representation $\left(\setsymdue_{\unkuno_{1,0}},\intuno_{\unkuno_{1,0}}, \gfy_{\unkuno_{1,0}}\right)$ of the correspondence between $\pmV_1$ and $\corrarr[\pmV]_1$ at $\unkuno_{1,0}$, a representation $\left(\setsymdue_{\unkuno_{2,0}},\intuno_{\unkuno_{2,0}}, \gfy_{\unkuno_{2,0}}\right)$ of the p-correspondence between $\pmV_2$ and $\corrarr[\pmV]_2$ at $\unkuno_{2,0}$.\newline
Define:
\begin{flalign*}
&\left\{
\begin{array}{l}
h:\setsymdue_{\unkuno_{2,0}}\times     \trasl[-\unkuno_{2,0}]\left(\intuno_{\unkuno_{2,0}}\right)\rightarrow\mathbb{R}^{l_2}\times\mathbb{R}^{m_2}\hspace{4pt}\text{by}\hspace{4pt}\text{setting}\\[4pt]
h\left(\unkuno_1,\unkuno_2\right)=\left(0,\unkuno_2\right)\hspace{15pt}\forall\left(\unkuno_1,\unkuno_2\right)\in\setsymdue_{\unkuno_{2,0}}\times     \trasl[-\unkuno_{2,0}]\left(\intuno_{\unkuno_{2,0}}\right)\text{;}
\end{array}
\right.\\[8pt]
&\begin{array}{l}
 \gfy_{\left(\unkuno_{1,0},\unkuno_{2,0}\right)}= \gfy_{\unkuno_{1,0}}\genfuncomp \lboundgenf \gfy_{\unkuno_{2,0}}\genfunsum 
 -1 \genfunscalp\left( \gfy_{\unkuno_{2,0}}\genfuncomp  h    \right)
, \incl{\trasl[-\unkuno_{1,0}]\left(\intuno_{\unkuno_{1,0}}\right)}{\mathbb{R}^{m_1}}\rboundgenf  \text{;}
\end{array}\\[8pt]
&\left\{
\begin{array}{l}
\pmV_3\in\homF+\GTopcat+\left(\intuno_1\times\intuno_2,\genfquot(l_2)+n+\right) \hspace{4pt}\text{by setting}\\[4pt]
\pmV_3\left(\unkuno_1,\unkuno_2\right)=\pmV_1\left(\unkuno_1\right)\compquotpoint\pmV_2\left(\unkuno_2\right) \hspace{15pt} \forall \left(\unkuno_1,\unkuno_2\right)\in \intuno_1\times\intuno_2\text{;}
\end{array}
\right.\\[8pt]
&\left\{
\begin{array}{l}
\text{the}\hspace{4pt}\text{arrow}\hspace{4pt}\corrarrre[\pmV]_{3}\in \homF+\GTopcat+\left(\intuno_{\unkuno_{1,0}}\times\intuno_{\unkuno_{2,0}},\genfquot(l+m_1+m_2)+n+\right)\hspace{4pt}\text{by}\hspace{4pt}\text{setting}\\[4pt]
\hspace{5pt}\corrarrre[\pmV]_{3}\left(\unkuno_1,\unkuno_2\right)=\preffuncaap{\gfy}_{\left(\unkuno_{1,0},\unkuno_{2,0}\right)}\left(0,\left(\unkuno_1,\unkuno_2\right)-\left(\unkuno_{1,0}, \unkuno_{2,0}\right)\right)\\[4pt]
\hspace{175pt}\forall \left(\unkuno_1,\unkuno_2\right) \in \intuno_{\unkuno_{1,0}}\times\intuno_{\unkuno_{2,0}}\text{.}
\end{array}
\right.
\end{flalign*}
Then:
\begin{flalign*}
&\begin{array}{l}
\corrarrre[\pmV]_{3}\hspace{4pt} \text{corresponds}\hspace{4pt} \text{to}\hspace{4pt} \pmV_{3}\funcomp \incl{\intuno_{\unkuno_{1,0}}\times\intuno_{\unkuno_{2,0}}}{\intuno_{\unkuno_{1}}\times\intuno_{\unkuno_{2}}}\text{;}
\end{array}\\[8pt]
&\left\{
\begin{array}{l}
\left(\setsymdue_{\unkuno_{1,0}}\times\setsymdue_{\unkuno_{2,0}},\intuno_{\unkuno_{2,0}}, \gfy_{3}\right)\hspace{4pt} \text{is}\hspace{4pt} \text{a}\hspace{4pt} \text{representation}\hspace{4pt} \text{of}\hspace{4pt} \text{the}\hspace{4pt} \text{correspondence}\\[4pt]
\text{between}\hspace{4pt} \pmV_{3}\hspace{4pt} \text{and}\hspace{4pt} \corrarrre[\pmV]_{3}\funcomp \incl{\intuno_{\unkuno_{1,0}}\times\intuno_{\unkuno_{2,0}}}{\intuno_{\unkuno_{1}}\times\intuno_{\unkuno_{2}}}\hspace{4pt} \text{at}\hspace{4pt} \left(\unkuno_{1,0},\unkuno_{2,0}\right)\text{.}
\end{array}
\right.
\end{flalign*}
\item Fix $l_1,l_2,m_1,m_2,n \in \mathbb{N}_0$, $\intuno_1\sqsubseteqdentro \mathbb{R}^{m_1}$, $\intuno_2\sqsubseteqdentro \mathbb{R}^{m_2}$,  $\unkuno_{1,0}\in \intuno_1$, $\unkuno_{2,0}\in \intuno_2$, $\pmV_1\in\homF+\GTopcat+\left(\intuno_1,\genfquotcccz(l_1)+n+\right)$, $\corrarr[\pmV]_{1}\in \homF+\GTopcat+\left(\intuno_1,\genfquotcccz(l_1+m_1)+n+\right)$ 
p-corresponding to $\pmV_1$, $\pmV_2\in\homF+\GTopcat+\left(\intuno_2,\genfquotcccz(l_2)+l_1+\right)$, $\corrarr[\pmV]_{2}\in \homF+\GTopcat+\left(\intuno_2,\genfquotcccz(l_2+m_2)+l_1+\right)$ p-corresponding to $\pmV_2$, a representation $\left(\setsymdue_{\unkuno_{1,0}},\intuno_{\unkuno_{1,0}}, \gfy_{\unkuno_{1,0}}\right)$ of the p-correspon\-dence between $\pmV_1$ and $\corrarr[\pmV]_1$ at $\unkuno_{1,0}$, a representation $\left(\setsymdue_{\unkuno_{2,0}},\intuno_{\unkuno_{2,0}}, \gfy_{\unkuno_{2,0}}\right)$ of the p-correspondence between $\pmV_2$ and $\corrarr[\pmV]_2$ at $\unkuno_{2,0}$.\newline
Define:
\begin{flalign*}
&\left\{
\begin{array}{l}
h:\setsymdue_{\unkuno_{2,0}}\times     \trasl[-\unkuno_{2,0}]\left(\intuno_{\unkuno_{2,0}}\right)\rightarrow\mathbb{R}^{l_2}\times\mathbb{R}^{m_2}\hspace{4pt}\text{by}\hspace{4pt}\text{setting}\\[4pt]
h\left(\unkuno_1,\unkuno_2\right)=\left(0,\unkuno_2\right)\hspace{15pt}\forall\left(\unkuno_1,\unkuno_2\right)\in\setsymdue_{\unkuno_{2,0}}\times     \trasl[-\unkuno_{2,0}]\left(\intuno_{\unkuno_{2,0}}\right)\text{;}
\end{array}
\right.\\[8pt]
&\begin{array}{l}
 \gfy_{\left(\unkuno_{1,0},\unkuno_{2,0}\right)}= \gfy_{\unkuno_{1,0}}\genfuncomp \lboundgenf \gfy_{\unkuno_{2,0}}\genfunsum 
 -1 \genfunscalp\left( \gfy_{\unkuno_{2,0}}\genfuncomp h \right)
, \incl{\trasl[-\unkuno_{1,0}]\left(\intuno_{\unkuno_{1,0}}\right)}{\mathbb{R}^{m_1}}\rboundgenf  \text{;}
\end{array}\\[8pt]
&\left\{
\begin{array}{l}
\pmV_3\in\homF+\GTopcat+\left(\intuno_1\times\intuno_2,\genfquot(l_2)+n+\right) \hspace{4pt}\text{by setting}\\[4pt]
\pmV_3\left(\unkuno_1,\unkuno_2\right)=\pmV_1\left(\unkuno_1\right)\compquotpoint\pmV_2\left(\unkuno_2\right) \hspace{15pt} \forall \left(\unkuno_1,\unkuno_2\right)\in \intuno_1\times\intuno_2\text{;}
\end{array}
\right.\\[8pt]
&\left\{
\begin{array}{l}
\text{the}\hspace{4pt}\text{arrow}\hspace{4pt}\corrarr[\pmV]_{3}\in \homF+\GTopcat+\left(\intuno_{\unkuno_{1,0}}\times\intuno_{\unkuno_{2,0}},\genfquot(l+m_1+m_2)+n+\right)\hspace{4pt}\text{by}\hspace{4pt}\text{setting}\\[4pt]
\hspace{5pt}\corrarr[\pmV]_{3}\left(\unkuno_1,\unkuno_2\right)=\ffuncaap{\gfy}_{\left(\unkuno_{1,0},\unkuno_{2,0}\right)}\left(0,\left(\unkuno_1,\unkuno_2\right)-\left(\unkuno_{1,0}, \unkuno_{2,0}\right)\right)\\[4pt]
\hspace{175pt}\forall \left(\unkuno_1,\unkuno_2\right) \in \intuno_{\unkuno_{1,0}}\times\intuno_{\unkuno_{2,0}}\text{.}
\end{array}
\right.
\end{flalign*}
Then:
\begin{flalign*}
&\begin{array}{l}
\corrarr[\pmV]_{3}\hspace{4pt} \text{p-corresponds}\hspace{4pt} \text{to}\hspace{4pt} \pmV_{3}\funcomp \incl{\intuno_{\unkuno_{1,0}}\times\intuno_{\unkuno_{2,0}}}{\intuno_{\unkuno_{1}}\times\intuno_{\unkuno_{2}}}\text{;}
\end{array}\\[8pt]
&\left\{
\begin{array}{l}
\left(\setsymdue_{\unkuno_{1,0}}\times\setsymdue_{\unkuno_{2,0}},\intuno_{\unkuno_{2,0}}, \gfy_{3}\right)\hspace{4pt} \text{is}\hspace{4pt} \text{a}\hspace{4pt} \text{representation}\hspace{4pt} \text{of}\hspace{4pt} \text{the}\hspace{4pt} \text{p-correspondence}\\[4pt]
\text{between}\hspace{4pt} \pmV_{3}\hspace{4pt} \text{and}\hspace{4pt} \corrarr[\pmV]_{3}\funcomp \incl{\intuno_{\unkuno_{1,0}}\times\intuno_{\unkuno_{2,0}}}{\intuno_{\unkuno_{1}}\times\intuno_{\unkuno_{2}}}\hspace{4pt} \text{at}\hspace{4pt} \left(\unkuno_{1,0},\unkuno_{2,0}\right)\text{.}
\end{array}
\right.
\end{flalign*}
\item Fix $l_1,l_2,m_1,m_2,n_1,n_2 \in \mathbb{N}_0$,  $\intuno_1\sqsubseteqdentro \mathbb{R}^{m_1}$, $\intuno_2\sqsubseteqdentro \mathbb{R}^{m_2}$, arrows\newline
$\pmV_1\in\homF+\GTopcat+\left(\intuno_1,\genfquotcccz(l_1)+n_1+\right)$, $\corrarr[\pmV]_{1}$ p-corresponding to $\pmV_1$,\newline
$\pmV_2\in\homF+\GTopcat+\left(\intuno_2,\genfquotcccz(l_2)+n_2+\right)$, $\corrarr[\pmV]_{2}$ p-corresponding to $\pmV_2$.\newline
Define:
\begin{flalign*}
&\left\{
\begin{array}{l}
h:\mathbb{R}^{l_1}\times\mathbb{R}^{l_2}\times\mathbb{R}^{m_1}\times\mathbb{R}^{m_2}\rightarrow\mathbb{R}^{l_1}\times\mathbb{R}^{m_1}\times\mathbb{R}^{l_2}\times\mathbb{R}^{m_2}\hspace{4pt}\text{by}\hspace{4pt} \text{setting}\\[4pt]
h\left(\unkuno_1,\unkuno_2,\unkuno_3,\unkuno_4\right)=\left(\unkuno_1,\unkuno_3,\unkuno_2,\unkuno_3\right)\\[4pt] 
\hspace{95pt}\forall \left(\unkuno_1,\unkuno_2,\unkuno_3,\unkuno_3\right)\in \mathbb{R}^{l_1}\times\mathbb{R}^{l_2}\times\mathbb{R}^{m_1}\times\mathbb{R}^{m_2}\text{;}
\end{array}
\right.\\[8pt]
&\left\{
\begin{array}{l}
\pmV_3\in\homF+\GTopcat+\left(\intuno_1\times\intuno_2,\genfquotcccz(l_1+l_2)+\argcompl{n_1+n_2}+\right) \hspace{4pt}\text{by setting}\\[4pt]
\pmV_3\left(\unkuno_1,\unkuno_2\right)=\lboundquotpoint\pmV_1\left(\unkuno_1\right),\pmV_2\left(\unkuno_2\right)\rboundquotpoint
\hspace{15pt} \forall \left(\unkuno_1,\unkuno_2\right)\in \intuno_1\times\intuno_2\text{;}
\end{array}
\right.\\[8pt]
&\left\{
\begin{array}{l}\corrarr[\pmV]_3\in\homF+\GTopcat+\left(\intuno_1\times\intuno_2,\genfquotcccz(l_1+l_2+m_1+m_2)+\argcompl{n_1+n_2}+\right) \\[4pt]
\text{by setting}\\[4pt]
\corrarr[\pmV]_3\left(\unkuno_1,\unkuno_2\right)=\lboundquotpoint\corrarr[\pmV]_{1}\left(\unkuno_1\right),\corrarr[\pmV]_{2}\left(\unkuno_2\right)\rboundquotpoint\compquotpoint\lininclquot\left(\gczsfx[h]\right)
\hspace{15pt}\forall \left(\unkuno_1,\unkuno_2\right)\in \intuno_1\times\intuno_2\text{.}
\end{array}
\right.
\end{flalign*} 
Then $\corrarr[\pmV]_3$ p-corresponds to $\pmV_3$.
\item Fix $l,m_1,m_2,n \in \mathbb{N}_0$, $\intuno_1\sqsubseteqdentro \mathbb{R}^{m_1}$, $\intuno_2\sqsubseteqdentro \mathbb{R}^{m_2}$, $\unkuno_{1,0}\in \intuno_1$, $\unkuno_{2,0}\in \intuno_2$,  $f \in \homF+\GTopcat+\left(\intuno_1,\intuno_2\right)$ with $f\left(\unkuno_{1,0}\right)=\unkuno_{2,0}$, $\pmV_2\in\homF+\GTopcat+\left(\intuno_2,\genfquotcccz(l)+n+\right)$, an arrow $\corrarr[\pmV]_2$ p-corresponding to $\pmV_2$, a representation $\left(\setsymdue_{\unkuno_{2,0}},\intuno_{\unkuno_{2,0}}, \gfy_{\unkuno_{2,0}}\right)$ of the p-correspondence between $\pmV_2$ and $\corrarr[\pmV]_2$ at $\unkuno_{2,0}$, an open neighborhood $\intuno_{\unkuno_{1,0}}$ of $\unkuno_{1,0} \in \intuno_1$ with $f\left(\intuno_{\unkuno_{1,0}}\right)\compacont \intuno_{\unkuno_{2,0}}$, $\gfx \in \genfcontcont(\intuno_{\unkuno_{1,0}})+\mathbb{R}^{m_2}+$ with $\evalcomptwo\left(\gfx\right)=f\funcomp \incl{\intuno_{\unkuno_{1,0}}}{\intuno_{1}}$.\newline
Define:
\begin{flalign*}
&\left\{
\begin{array}{l}
\gfy_{\unkuno_{1,0}} \in \genfcontcont(\setsymdue_{\unkuno_{2,0}}\times\trasl[-\unkuno_{1,0}]\left(\intuno_{\unkuno_{1,0}}\right))+\mathbb{R}^{n}+\hspace{4pt}\text{by}\hspace{4pt}\text{setting}\\[4pt]
\gfy_{\unkuno_{1,0}}=\gfy_{\unkuno_{2,0}} \genfuncomp \lboundgenf \incl{\setsymdue_{\unkuno_{2,0}}}{\mathbb{R}^l},\trasl[-\unkuno_{2,0}]  \genfuncomp  \gfx  \genfuncomp  \trasl[\unkuno_{1,0}]   \rboundgenf\text{;}
\end{array}
\right.\\[8pt]
&\left\{
\begin{array}{l}
\text{the}\hspace{4pt}\text{arrow}\hspace{4pt}\pmV_{1,0}\in \homF+\GTopcat+\left(\intuno_{1,0},\genfquot(l)+n+\right)\hspace{4pt}\text{by}\hspace{4pt}\text{setting}\\[4pt]
\pmV_{1,0}\left(\unkuno\right)=\ffuncaap{\gfy}_{\unkuno_{1,0}}\left(0,\unkuno-\unkuno_{1,0}\right)\compquotpoint \lininclquot\left(\gczsfx[\argcompl{\pointedincl<\argcompl{\left(\mathbb{R}^{l},0\right),\left(\mathbb{R}^{m_1},0\right)}<>1>}]\right)\hspace{10pt}\forall \unkuno \in \intuno_{\unkuno_{1,0}}\text{;}
\end{array}
\right.\\[8pt]
&\left\{
\begin{array}{l}
\text{the}\hspace{4pt}\text{arrow}\hspace{4pt}\corrarr[\pmV]_{1,0}\in \homF+\GTopcat+\left(\intuno_{1,0},\genfquot(l+m_1)+n+\right)\hspace{4pt}\text{by}\hspace{4pt}\text{setting}\\[4pt]
\corrarr[\pmV]_{1,0}\left(\unkuno\right)=\ffuncaap{\gfy}_{\unkuno_{1,0}}\left(0,\unkuno-\unkuno_{1,0}\right)\hspace{15pt}\forall \unkuno \in \intuno_{\unkuno_{1,0}}\text{.}
\end{array}
\right.
\end{flalign*}
Then:
\begin{flalign}
&\begin{array}{l}
\corrarr[\pmV]_{1,0}\hspace{4pt} \text{p-corresponds}\hspace{4pt} \text{to}\hspace{4pt} \pmV_{1,0}\text{;}
\end{array}\nonumber\\[8pt]
&\left\{
\begin{array}{l}
\left(\setsymdue_{\unkuno_{2,0}},\intuno_{\unkuno_{1,0}}, \gfy_{\unkuno_{1,0}}\right)\hspace{4pt} \text{is}\hspace{4pt} \text{a}\hspace{4pt} \text{representation}\hspace{4pt} \text{of}\hspace{4pt} \text{the}\hspace{4pt} \text{p-correspondence}\\[4pt]
\text{between}\hspace{4pt} \pmV_{1,0}\hspace{4pt} \text{and}\hspace{4pt} \corrarr[\pmV]_{1,0}\hspace{4pt} \text{at}\hspace{4pt} \unkuno_{1,0}\text{;}
\end{array}
\right.\nonumber\\[8pt]
&\begin{array}{l}
\evalcompquotcontcontzero\funcomp\pmV_{1,0}=\evalcompquotcontcontzero\funcomp\left(\pmV_2\funcomp\left( f\funcomp\incl{\intuno_{\unkuno_{1,0}}}{\mathbb{R}^{m_1}} \right)\right)\text{.}
\end{array}\label{pugualeval}
\end{flalign}
\item Fix $l,m,n \in \mathbb{N}_0$, $\bsfm\in\left\{ \Fp_{\mathsf{l}}, \Fq_{\mathsf{l}},\Fqq_{\mathsf{l}}\;:\; \mathsf{l}\in \left\{1,...,l\right\}\right\}$, interval $\intuno\sqsubseteqdentro \mathbb{R}^{m}$, arrows $\pmV\in\homF+\GTopcat+\left(\intuno,\genfquotcccz(l)+n+\right)$,
$\corrarr[\pmV]$ p-corresponding to $\pmV$.\newline
Then $\bsfm\bsfmuquotpoint\corrarr[\pmV]$ p-corresponds to $\bsfm\bsfmuquotpoint\pmV$.
\item Fix $l,m,n \in \mathbb{N}_0$, $\mathsf{l}\in \left\{1,...,l\right\}$, $\intuno\sqsubseteqdentro \mathbb{R}^{m}$, $\pmV\in\homF+\GTopcat+\left(\intuno,\genfquotcccz(l)+n+\right)$, $\corrarr[\pmV]$ p-corresponding to $\pmV$.\newline
Then arrow $\Fint_{\mathsf{l}}\bsfmuquotpoint\left(\corrarr[\pmV] \sumquotpoint
-1 \genfunscalp \left( \corrarr[\pmV]\compquotpoint 
\lininclquot\left(\gczsfx[\argcompl{\lboundquotpoint   \cost<\mathbb{R}^l<>\mathbb{R}^l>+0+, \idobj+\mathbb{R}^m+ \rboundquotpoint}]\right)\right)\right)$ p-corresponds to $\Fint_{\mathsf{l}}\bsfmuquotpoint\pmV$.
\end{enumerate}
\end{proposition}
\begin{proof}\mbox{}\newline
\textnormal{\textbf{Proof of statement 1.}}\ \ Fix $\unkuno_0\in \intuno$, a representation $\left(\setsymdue_{\unkuno_0},\intuno_{\unkuno_0}, \gfy_{\unkuno_0}\right)$ of the p-correspondence between $\pmV$ and $\corrarr[\pmV]$ at $\unkuno_0$.\newline
Statement \eqref{tildv} follows by the chain of equalities below
\begin{flalign*}
&\begin{array}{l}
\corrarr[\pmV]\left(\unkuno\right)  \compquotpoint\lininclquot\left(\gczsfx[\argcompl{\pointedincl<\argcompl{\left(\mathbb{R}^{l},0\right),\left(\mathbb{R}^{m},\unktre\right)}<>1>}] \right)= \\[4pt]
\hspace{140pt} 
\ffuncaap{\gfy}_{\unkuno_0}\left(0,\unkuno-\unkuno_0\right)\compquotpoint\lininclquot\left(\gczsfx[\argcompl{\pointedincl<\argcompl{\left(\mathbb{R}^{l},0\right),\left(\mathbb{R}^{m},\unktre\right)}<>1>}] \right)=
\end{array} \\[4pt]
&\begin{array}{l}
\genfquotfun\Big(\prefuncaap{\gfy}_{\unkuno_0}\left(0,\unkuno-\unkuno_0\right) \genfunsum
 -1 \genfunscalp\left(\prefuncaap{\gfy}_{\unkuno_0}\left(0,\unkuno-\unkuno_0\right)       \genfuncomp \cost<\mathbb{R}^{l+m}<>\mathbb{R}^{l+m}>+0+ \right)\Big)\compquotpoint\\[4pt]
\hspace{210pt} \lininclquot\left(\gczsfx[\argcompl{\pointedincl<\argcompl{\left(\mathbb{R}^{l},0\right),\left(\mathbb{R}^{m},\unktre\right)}<>1>}] \right)=
\end{array} \\[4pt]
&\begin{array}{l}
\genfquotfun\Big(\big(\prefuncaap{\gfy}_{\unkuno_0}\left(0,\unkuno-\unkuno_0\right) \genfunsum
 -1 \genfunscalp\left(\prefuncaap{\gfy}_{\unkuno_0}\left(0,\unkuno-\unkuno_0\right)       \genfuncomp \cost<\argcompl{\mathbb{R}^{l+m}}<>\argcompl{\mathbb{R}^{l+m}}>+0+ \right)\big)\genfuncomp\\[4pt]
\hspace{210pt} \pointedincl<\argcompl{\left(\mathbb{R}^{l},0\right),\left(\mathbb{R}^{m},\unktre\right)}<>1>   \Big)=
\end{array}\\[4pt]
&\begin{array}{l}
\genfquotfun\Big(\prefuncaap{\gfy}_{\unkuno_0}\left(0,\unkuno\vecsum \unktre-\unkuno_0\right) \genfuncomp \pointedincl<\argcompl{\left(\mathbb{R}^{l},0\right),\left(\mathbb{R}^{m},0\right)}<>1> \genfunsum\\[4pt]
 -1 \genfunscalp\left(\prefuncaap{\gfy}_{\unkuno_0}\left(0,\unkuno-\unkuno_0\right)       \genfuncomp \cost<\argcompl{\mathbb{R}^{l+m}}<>\argcompl{\mathbb{R}^{l+m}}>+0+ \right)\genfuncomp \pointedincl<\argcompl{\left(\mathbb{R}^{l},0\right),\left(\mathbb{R}^{m},0\right)}<>1>   \Big)=
\end{array} \\[4pt]
&\begin{array}{l}
\genfquotfun\Big(\big(\prefuncaap{\gfy}_{\unkuno_0}\left(0,\unkuno\vecsum\unktre-\unkuno_0\right) \genfunsum\\[4pt]
 -1 \genfunscalp\left(\prefuncaap{\gfy}_{\unkuno_0}\left(0,\unkuno\vecsum\unktre-\unkuno_0\right)       \genfuncomp \cost<\argcompl{\mathbb{R}^{l+m}}<>\argcompl{\mathbb{R}^{l+m}}>+0+ \right)\big)\genfuncomp \pointedincl<\argcompl{\left(\mathbb{R}^{l},0\right),\left(\mathbb{R}^{m},0\right)}<>1> \genfunsum\\[4pt]
\big( 
\prefuncaap{\gfy}_{\unkuno_0}\left(0,\unkuno\vecsum\unktre-\unkuno_0\right)       \genfuncomp \cost<\argcompl{\mathbb{R}^{l+m}}<>\argcompl{\mathbb{R}^{l+m}}>+0+  \genfunsum\\[4pt]
-1 \genfunscalp\left(\prefuncaap{\gfy}_{\unkuno_0}\left(0,\unkuno-\unkuno_0\right)       \genfuncomp \cost<\argcompl{\mathbb{R}^{l+m}}<>\argcompl{\mathbb{R}^{l+m}}>+0+ \right)\big)\genfuncomp \pointedincl<\argcompl{\left(\mathbb{R}^{l},0\right),\left(\mathbb{R}^{m},0\right)}<>1>   \Big)=
\end{array} \\[4pt]
&\begin{array}{l}
\ffuncaap{\gfy}_{\unkuno_0}\left(0,\unkuno\vecsum\unktre-\unkuno_0\right)\compquotpoint\lininclquot\left(\gczsfx[\argcompl{\pointedincl<\argcompl{\left(\mathbb{R}^{l},0\right),\left(\mathbb{R}^{m},0\right)}<>1>}] \right)\sumquotpoint\\[4pt]
\hspace{10pt}\genfquotfun\Big(\big( 
\prefuncaap{\gfy}_{\unkuno_0}\left(0,\unkuno-\unkuno_0\right)  \genfunsum
-1 \genfunscalp\left(\prefuncaap{\gfy}_{\unkuno_0}\left(0,\unkuno-\unkuno_0\right)       \genfuncomp \cost<\argcompl{\mathbb{R}^{l+m}}<>\argcompl{\mathbb{R}^{l+m}}>+0+ \right)\big)\\[4pt]
\hspace{200pt}\genfuncomp  \cost<\mathbb{R}^{l}<>\argcompl{\mathbb{R}^{l+m}}>+\left(0, \unktre\right)+ \Big)=
\end{array} \\[4pt]
&\begin{array}{l}
\pmV\left(\unkuno+\unktre\right) \sumquotpoint  \corrarr[\pmV]\left(\unkuno\right)  \compquotpoint\lininclquot\left(\gczsfx[\argcompl{\cost<\mathbb{R}^{l}<>\argcompl{\mathbb{R}^{l+m}}>+\left(0, \unktre\right)+}]\right)\hspace{10pt}
\text{for}\hspace{4pt}\text{any}\hspace{4pt} \unkuno\in \intuno\text{,}\hspace{4pt}\unktre\in \mathbb{R}^{m}\text{.}
\end{array}
\end{flalign*}
Statement \eqref{tildtildv} follows by the chain of equalities below 
\begin{flalign*}
&\begin{array}{l}
\corrarr[\pmV]\left(\unkuno\right)  \compquotpoint\lininclquot\left(\gczsfx[\argcompl{\pointedincl<\argcompl{\left(\mathbb{R}^{l},0\right),\left(\mathbb{R}^{m},\unktre_1\vecsum \unktre_2\right)}<>1>}] \right)= \\[4pt]
\hspace{100pt} 
\ffuncaap{\gfy}_{\unkuno_0}\left(0,\unkuno-\unkuno_0\right)\compquotpoint\lininclquot\left(\gczsfx[\argcompl{\pointedincl<\argcompl{\left(\mathbb{R}^{l},0\right),\left(\mathbb{R}^{m},\unktre_1\vecsum \unktre_2\right)}<>1>}] \right)=
\end{array} \\[4pt]
&\begin{array}{l}
\genfquotfun\Big(\prefuncaap{\gfy}_{\unkuno_0}\left(0,\unkuno-\unkuno_0\right) \genfunsum
 -1 \genfunscalp\left(\prefuncaap{\gfy}_{\unkuno_0}\left(0,\unkuno-\unkuno_0\right)       \genfuncomp \cost<\argcompl{\mathbb{R}^{l+m}}<>\argcompl{\mathbb{R}^{l+m}}>+ 0+ \right)\Big)\compquotpoint\\[4pt]
\hspace{170pt} \lininclquot\left(\gczsfx[\argcompl{\pointedincl<\argcompl{\left(\mathbb{R}^{l},0\right),\left(\mathbb{R}^{m},\unktre_1\vecsum \unktre_2\right)}<>1>}] \right)=
\end{array} \\[4pt]
&\begin{array}{l}
\genfquotfun\Big(\big(\prefuncaap{\gfy}_{\unkuno_0}\left(0,\unkuno-\unkuno_0\right) \genfunsum
 -1 \genfunscalp\left(\prefuncaap{\gfy}_{\unkuno_0}\left(0,\unkuno-\unkuno_0\right)       \genfuncomp \cost<\argcompl{\mathbb{R}^{l+m}}<>\argcompl{\mathbb{R}^{l+m}}>+ 0+ \right)\big)\genfuncomp\\[4pt]
\hspace{170pt} \pointedincl<\argcompl{\left(\mathbb{R}^{l},0\right),\left(\mathbb{R}^{m},\unktre_1\vecsum \unktre_2\right)}<>1>   \Big)=
\end{array}\\[4pt]
&\begin{array}{l} \\[4pt]
\genfquotfun\Big(\prefuncaap{\gfy}_{\unkuno_0}\left(0,\unkuno\vecsum \unktre_1\vecsum \unktre_2-\unkuno_0\right) \genfuncomp \pointedincl<\argcompl{\left(\mathbb{R}^{l},0\right),\left(\mathbb{R}^{m},0\right)}<>1> \genfunsum\\[4pt]
 -1 \genfunscalp\left(\prefuncaap{\gfy}_{\unkuno_0}\left(0,\unkuno-\unkuno_0\right)       \genfuncomp \cost<\argcompl{\mathbb{R}^{l+m}}<>\argcompl{\mathbb{R}^{l+m}}>+ 0+ \right)\genfuncomp \pointedincl<\argcompl{\left(\mathbb{R}^{l},0\right),\left(\mathbb{R}^{m},0\right)}<>1>   \Big)=
\end{array} \\[4pt]
&\begin{array}{l}
\genfquotfun\Big(\big(\prefuncaap{\gfy}_{\unkuno_0}\left(0,\unkuno\vecsum\unktre_1-\unkuno_0\right) \genfunsum\\[4pt]
 -1 \genfunscalp\left(\prefuncaap{\gfy}_{\unkuno_0}\left(0,\unkuno\vecsum\unktre_1-\unkuno_0\right)       \genfuncomp \cost<\argcompl{\mathbb{R}^{l+m}}<>\argcompl{\mathbb{R}^{l+m}}>+0+\right)\big)\genfuncomp \pointedincl<\argcompl{\left(\mathbb{R}^{l},0\right),\left(\mathbb{R}^{m},\unktre_2}<>1> \right)\genfunsum\\[4pt]
\big( 
\prefuncaap{\gfy}_{\unkuno_0}\left(0,\unkuno\vecsum\unktre_1-\unkuno_0\right)       \genfuncomp \cost<\argcompl{\mathbb{R}^{l+m}}<>\argcompl{\mathbb{R}^{l+m}}>+0+  \genfunsum\\[4pt]
-1 \genfunscalp\left(\prefuncaap{\gfy}_{\unkuno_0}\left(0,\unkuno-\unkuno_0\right)       \genfuncomp \cost<\argcompl{\mathbb{R}^{l+m}}<>\argcompl{\mathbb{R}^{l+m}}>+0+ \right)\big)\genfuncomp \pointedincl<\argcompl{\left(\mathbb{R}^{l},0\right),\left(\mathbb{R}^{m},\unktre_2\right)}<>1>  \Big)=
\end{array} \\[4pt]
&\begin{array}{l}
\ffuncaap{\gfy}_{\unkuno_0}\left(0,\unkuno\vecsum\unktre_1-\unkuno_0\right)\compquotpoint\lininclquot\left(\gczsfx[\argcompl{\pointedincl<\argcompl{\left(\mathbb{R}^{l},0\right),\left(\mathbb{R}^{m},\unktre_2\right)}<>1>}] \right)\sumquotpoint\\[4pt]
\hspace{10pt}\genfquotfun\Big(\big( 
\prefuncaap{\gfy}_{\unkuno_0}\left(0,\unkuno-\unkuno_0\right)  \genfunsum
-1 \genfunscalp\left(\prefuncaap{\gfy}_{\unkuno_0}\left(0,\unkuno-\unkuno_0\right)       \genfuncomp \cost<\argcompl{\mathbb{R}^{l+m}}<>\argcompl{\mathbb{R}^{l+m}}>+0+ \right)\big)\\[4pt]
\hspace{200pt}\genfuncomp  \cost<\mathbb{R}^{l}<>\argcompl{\mathbb{R}^{l+m}}>+\left(0, \unktre_1\right)+  \Big)=
\end{array} \\[4pt]
&\begin{array}{l}
\ffuncaap{\gfy}_{\unkuno_0}\left(0,\unkuno\vecsum\unktre_1-\unkuno_0\right)\compquotpoint\lininclquot\left(\gczsfx[\argcompl{\pointedincl<\argcompl{\left(\mathbb{R}^{l},0\right),\left(\mathbb{R}^{m},\unktre_2\right)}<>1>}] \right)\sumquotpoint\\[4pt]
\hspace{60pt} \corrarr[\pmV]\left(\unkuno\right)  \compquotpoint\lininclquot\left(\gczsfx[\argcompl{\cost<\mathbb{R}^{l}<>\argcompl{\mathbb{R}^{l+m}}>+\left(0, \unktre_1\right)+}]\right)\hspace{10pt}
\text{for}\hspace{4pt}\text{any}\hspace{4pt} \unkuno\in \intuno\text{,}\hspace{4pt}\unktre_1,\unktre_2\in \mathbb{R}^{m}\text{.}
\end{array}
\end{flalign*}
\textnormal{\textbf{Proof of statement 2.}}\ \ Fix $\unkuno_0\in \intuno$, a representation $\left(\setsymdue_{\unkuno_0},\intuno_{\unkuno_0}, \gfy_{\unkuno_0}\right)$ of the p-correspondence between $\pmV$ and $\corrarr[\pmV]$ at $\unkuno_0$.\newline
Define:
\begin{flalign*}
&\left\{
\begin{array}{l}
h:\setsymdue_{\unkuno_0}\times     \trasl[-\unkuno_{0}]\left(\intuno_{\unkuno_0}\right)\rightarrow\mathbb{R}^l\times\mathbb{R}^m\hspace{4pt}\text{by}\hspace{4pt}\text{setting}\\[4pt]
h\left(\unkuno_1,\unkuno_2\right)=\left(0,\unkuno_2\right)\hspace{15pt}\forall\left(\unkuno_1,\unkuno_2\right)\in\setsymdue_{\unkuno_0}\times     \trasl[-\unkuno_{0}]\left(\intuno_{\unkuno_0}\right)\text{;}
\end{array}
\right.\\[8pt]
&\left\{
\begin{array}{l}
\gfz_{\unkuno_0}\in \genf(\argcompl{\setsymdue_{\unkuno_{0}}\times \left(\trasl[-\unkuno_{0}]\left(\intuno_{\unkuno_{0}}\right)\right)})+\mathbb{R}^{n}+\hspace{4pt}\text{by}\hspace{4pt}\text{setting}\\[4pt]
\gfz_{\unkuno_0}=\gfx \compgenfquot\left(\gfy_{\unkuno_0}\genfunsum 
 -1 \genfunscalp\left(\gfy_{\unkuno_0}\genfuncomp h \right)\right)\text{.}
\end{array}
\right.
\end{flalign*}
Statement follows by proving that  $\left(\setsymdue_{\unkuno_0},\intuno_{\unkuno_0}, \gfz_{\unkuno_0}\right)$ is a representation of the correspondence between arrows $\left(\gfx \compgenfquot\pmV\right) \sumquotpoint\left(-1\scalpquotpoint   \gczsfx[\argcompl{\gfx\genfuncomp\cost<\argcompl{\mathbb{R}^{l+m}}<>\mathbb{R}^{n}>+0+}] \right)$ and\newline
$\left(\gfx \compgenfquot\corrarr[\pmV]\right) \sumquotpoint
\left(-1\scalpquotpoint 
 \gfx \compgenfquot \left(\corrarr[\pmV]\left(\unkuno\right)\compquotpoint 
\lininclquot\left(\gczsfx[\argcompl{\lboundquotpoint   \cost<\mathbb{R}^l<>\mathbb{R}^l>+0+, \idobj+\mathbb{R}^m+ \rboundquotpoint}]\right)  \right)\right) $ at $\unkuno_0$.\newline
This is proven by the two chain of equalities below, where equalities $(1)$ hold true by relation \eqref{(S 7 bis)}:
\begin{flalign*}
&\begin{array}{l}
\gfx \compgenfquot\left(\pmV\left(\unkuno\right)\right)=\gfx \compgenfquot\left(\ffuncaap{\gfy}_{\unkuno_0}\left(0,\unkuno-\unkuno_0\right)\compquotpoint \lininclquot\left(\gczsfx[\argcompl{\pointedincl<\argcompl{\left(\mathbb{R}^{l},0\right),\left(\mathbb{R}^{m},0\right)}<>1>}]\right)\right)=
\end{array}\\[4pt]
&\begin{array}{l}
\left(\gfx \compgenfquot\left(\ffuncaap{\gfy}_{\unkuno_0}\left(0,\unkuno-\unkuno_0\right)\right)\right)\compquotpoint \lininclquot\left(\gczsfx[\argcompl{\pointedincl<\argcompl{\left(\mathbb{R}^{l},0\right),\left(\mathbb{R}^{m},0\right)}<>1>}]\right)=
\end{array}\\[4pt]
&\begin{array}{l} 
\bigg(\gfx \compgenfquot\genfquotfun\Big(\gfy_{\unkuno_0} \genfuncomp \trasl[\left(0, \unkuno-\unkuno_0\right)] \genfunsum\\[4pt]
\hspace{50pt} 
 -1 \genfunscalp\left( \gfy_{\unkuno_0} \genfuncomp \trasl[\left(0, \unkuno-\unkuno_0\right)]  \genfuncomp \cost<\argcompl{\mathbb{R}^{l+m}}<>\argcompl{\mathbb{R}^{l+m}}>+0+
\right)\Big)\bigg)\compquotpoint \\[4pt]
\hspace{200pt}\lininclquot\left(\gczsfx[\argcompl{\pointedincl<\argcompl{\left(\mathbb{R}^{l},0\right),\left(\mathbb{R}^{m},0\right)}<>1>}]\right)=
\end{array}\\[4pt]
&\begin{array}{l}
\genfquotfun\Big(\gfx \genfuncomp\big(\gfy_{\unkuno_0} \genfuncomp \trasl[\left(0, \unkuno-\unkuno_0\right)] \genfunsum\\[4pt]
\hspace{50pt} 
 -1 \genfunscalp\left(\gfy_{\unkuno_0} \genfuncomp \trasl[\left(0, \unkuno-\unkuno_0\right)]  \genfuncomp \cost<\argcompl{\mathbb{R}^{l+m}}<>\argcompl{\mathbb{R}^{l+m}}>+0+
\right)\big)\Big)\compquotpoint \\[4pt]
\hspace{200pt}\lininclquot\left(\gczsfx[\argcompl{\pointedincl<\argcompl{\left(\mathbb{R}^{l},0\right),\left(\mathbb{R}^{m},0\right)}<>1>}]\right)\overset{(1)}{=}
\end{array}\\[4pt]
&\begin{array}{l}
\genfquotfun\bigg(\gfx \genfuncomp\big(\gfy_{\unkuno_0} \genfuncomp \trasl[\left(0, \unkuno-\unkuno_0\right)] \genfunsum\\[4pt]
\hspace{75pt} 
 -1 \genfunscalp\left(\gfy_{\unkuno_0} \genfuncomp \trasl[\left(0, \unkuno-\unkuno_0\right)]  \genfuncomp \cost<\argcompl{\mathbb{R}^{l+m}}<>\argcompl{\mathbb{R}^{l+m}}>+0+
\right)\big)\genfunsum\\[4pt]
-1 \genfunscalp\Big(\gfz_{\unkuno_0}\genfuncomp \trasl[\left(0, \unkuno-\unkuno_0\right)]
\genfuncomp\cost<\argcompl{\mathbb{R}^{l+m}}<>\argcompl{\mathbb{R}^{l+m}}>+0+ \Big)\genfunsum\\[4pt]
\Big(\gfz_{\unkuno_0}\genfuncomp \trasl[\left(0, \unkuno-\unkuno_0\right)]
\genfuncomp\cost<\argcompl{\mathbb{R}^{l+m}}<>\argcompl{\mathbb{R}^{l+m}}>+0+ \Big)   \bigg)\compquotpoint \lininclquot\left(\gczsfx[\argcompl{\pointedincl<\argcompl{\left(\mathbb{R}^{l},0\right),\left(\mathbb{R}^{m},0\right)}<>1>}]\right)=
\end{array}\\[4pt]
&\begin{array}{l}
\genfquotfun\bigg(\gfx \genfuncomp\big(\gfy_{\unkuno_0} \genfuncomp \trasl[\left(0, \unkuno-\unkuno_0\right)] \genfunsum\\[4pt]
 -1 \genfunscalp\left(\gfy_{\unkuno_0} \genfuncomp \trasl[\left(0, \unkuno-\unkuno_0\right)]  \genfuncomp \cost<\argcompl{\mathbb{R}^{l+m}}<>\argcompl{\mathbb{R}^{l+m}}>+0+\right)\big)\genfuncomp\pointedincl<\argcompl{\left(\mathbb{R}^{l},0\right),\left(\mathbb{R}^{m},0\right)}<>1>\genfunsum\\[4pt]
-1 \genfunscalp\Big(\gfz_{\unkuno_0}\genfuncomp \trasl[\left(0, \unkuno-\unkuno_0\right)]
\genfuncomp\cost<\argcompl{\mathbb{R}^{l+m}}<>\argcompl{\mathbb{R}^{l+m}}>+0+ \Big)\genfuncomp \pointedincl<\argcompl{\left(\mathbb{R}^{l},0\right),\left(\mathbb{R}^{m},0\right)}<>1>\genfunsum\\[4pt]
\Big(\gfz_{\unkuno_0}\genfuncomp \trasl[\left(0, \unkuno-\unkuno_0\right)]
\genfuncomp\cost<\argcompl{\mathbb{R}^{l+m}}<>\argcompl{\mathbb{R}^{l+m}}>+0+ \Big)\genfuncomp \pointedincl<\argcompl{\left(\mathbb{R}^{l},0\right),\left(\mathbb{R}^{m},0\right)}<>1>   \bigg)\overset{(1)}{=}
\end{array}\\[4pt]
&\begin{array}{l}
\genfquotfun\bigg(
\gfz_{\unkuno_0}\genfuncomp \trasl[\left(0, \unkuno-\unkuno_0\right)]\genfuncomp\pointedincl<\argcompl{\left(\mathbb{R}^{l},0\right),\left(\mathbb{R}^{m},0\right)}<>1>\genfunsum\\[4pt]
-1 \genfunscalp\Big(\gfz_{\unkuno_0}\genfuncomp \trasl[\left(0, \unkuno-\unkuno_0\right)]
\genfuncomp\cost<\argcompl{\mathbb{R}^{l+m}}<>\argcompl{\mathbb{R}^{l+m}}>+0+ \Big)\genfuncomp\pointedincl<\argcompl{\left(\mathbb{R}^{l},0\right),\left(\mathbb{R}^{m},0\right)}<>1>\genfunsum\\[4pt]
\Big(\gfz_{\unkuno_0}\genfuncomp \trasl[\left(0, \unkuno-\unkuno_0\right)]
\genfuncomp\cost<\argcompl{\mathbb{R}^{l+m}}<>\argcompl{\mathbb{R}^{l+m}}>+0+ \Big)\genfuncomp\pointedincl<\argcompl{\left(\mathbb{R}^{l},0\right),\left(\mathbb{R}^{m},0\right)}<>1>  \bigg)=
\end{array}\\[4pt]
&\begin{array}{l}
\genfquotfun\bigg(
\left(\funcaap{\gfz}_{\unkuno_0}\left(0,\unkuno-\unkuno_0\right)\right)
\genfuncomp\pointedincl<\argcompl{\left(\mathbb{R}^{l},0\right),\left(\mathbb{R}^{m},0\right)}<>1>\genfunsum
\gfx\genfuncomp\cost<\argcompl{\mathbb{R}^{l+m}}<>\mathbb{R}^{n}>+0+    \bigg) =
\end{array}\\[4pt]
&\begin{array}{l}
\left(\ffuncaap{\gfz}_{\unkuno_0}\left(0,\unkuno-\unkuno_0\right)\right)\compquotpoint \lininclquot\left(\gczsfx[\argcompl{\pointedincl<\argcompl{\left(\mathbb{R}^{l},0\right),\left(\mathbb{R}^{m},0\right)}<>1>}]\right)
\sumquotpoint\\[4pt]
\hspace{160pt}
\genfquotfun\left(\gfx\genfuncomp\cost<\argcompl{\mathbb{R}^{l+m}}<>\mathbb{R}^{n}>+0+    \right) 
\hspace{15pt}\forall \unkuno \in \intuno_{\unkuno_0}\text{;}
\end{array}\\[20pt]
&\begin{array}{l}
\gfx \compgenfquot\left(\corrarr[\pmV]\left(\unkuno\right)\right)=\gfx \compgenfquot \ffuncaap{\gfy}_{\unkuno_0}\left(0,\unkuno-\unkuno_0\right)\overset{(1)}{=}
\end{array}\\[4pt]
&\begin{array}{l} 
\gfx \compgenfquot\genfquotfun\Big(\gfy_{\unkuno_0} \genfuncomp \trasl[\left(0, \unkuno-\unkuno_0\right)] \genfunsum\\[4pt]
 \hspace{20pt}-1 \genfunscalp\left(\gfy_{\unkuno_0}\genfuncomp  \left(h \funcomp \trasl[\left(0, \unkuno-\unkuno_0\right)]\right)\right)\genfunsum
 \left(\gfy_{\unkuno_0}\genfuncomp \left(h  \funcomp \trasl[\left(0, \unkuno-\unkuno_0\right)]\right)\right)\genfunsum \\[4pt]
\hspace{100pt}  -1 \genfunscalp\left( \gfy_{\unkuno_0} \genfuncomp \trasl[\left(0, \unkuno-\unkuno_0\right)]  \genfuncomp \cost<\argcompl{\mathbb{R}^{l+m}}<>\argcompl{\mathbb{R}^{l+m}}>+0+
\right)\Big)=
\end{array}\\[4pt]
&\begin{array}{l} 
\gfx \compgenfquot\bigg(\genfquotfun\Big(\gfy_{\unkuno_0} \genfuncomp \trasl[\left(0, \unkuno-\unkuno_0\right)] \genfunsum
 -1 \genfunscalp\left(\gfy_{\unkuno_0}\genfuncomp \left(h\funcomp \trasl[\left(0, \unkuno-\unkuno_0\right)]\right)\right)\Big) \sumquotpoint\\[4pt]
 \genfquotfun\Big(\left(\gfy_{\unkuno_0}\genfuncomp \left(h  \funcomp \trasl[\left(0, \unkuno-\unkuno_0\right)]\right)\right)\genfunsum \\[4pt]
\hspace{100pt}  -1 \genfunscalp\left( \gfy_{\unkuno_0} \genfuncomp \trasl[\left(0, \unkuno-\unkuno_0\right)]  \genfuncomp \cost<\argcompl{\mathbb{R}^{l+m}}<>\argcompl{\mathbb{R}^{l+m}}>+0+
\right)\Big)\bigg)=
\end{array}\\[4pt]
&\begin{array}{l} 
\gfx \compgenfquot\bigg(
\genfquotfun
\Big(
\big(
\gfy_{\unkuno_0} \genfuncomp \trasl[\left(0, \unkuno-\unkuno_0\right)] \genfunsum
 -1 \genfunscalp\left(\gfy_{\unkuno_0}\genfuncomp \left(h \funcomp \trasl[\left(0, \unkuno-\unkuno_0\right)]\right)\right)
\big) \genfunsum\\[4pt]
 -1 \genfunscalp
\big(
   \big(
		\gfy_{\unkuno_0}\genfuncomp \trasl[\left(0, \unkuno-\unkuno_0\right)] \genfunsum\\[4pt]
\hspace{50pt}-1 \genfunscalp\left(\gfy_{\unkuno_0}\genfuncomp \left(h \funcomp \trasl[\left(0, \unkuno-\unkuno_0\right)]\right)\right) 
\big)\genfuncomp \left(\cost<\argcompl{\mathbb{R}^{l+m}}<>\argcompl{\mathbb{R}^{l+m}}>+0+ \right)   
\big)
\Big)
\bigg)\sumquotpoint\\[4pt]
 \gfx \compgenfquot \left(\corrarr[\pmV]\left(\unkuno\right)\compquotpoint 
\lininclquot\left(\gczsfx[\argcompl{\contsplbound  \cost<\mathbb{R}^l<>\mathbb{R}^l>+0+, \idobj+\mathbb{R}^m+ \contsprbound}]\right)  \right) =
\end{array}\\[4pt]
&\begin{array}{l} 
 \gfx \compgenfquot\genfquotfun\left(\ffuncaap{\gfz}_{\unkuno_0}\left(0,\unkuno-\unkuno_0\right)\right)\sumquotpoint
 \gfx \compgenfquot \left(\corrarr[\pmV]\left(\unkuno\right)\compquotpoint 
\lininclquot\left(\gczsfx[\argcompl{\contsplbound   \cost<\mathbb{R}^l<>\mathbb{R}^l>+0+, \idobj+\mathbb{R}^m+ \contsprbound}]\right)  \right)
\hspace{15pt}\forall \unkuno \in \intuno_{\unkuno_0}\text{.}
\end{array}
\end{flalign*}
\textnormal{\textbf{Proof of statement 3.}}\ \ Statement follows straightforwardly by direct checking.\newline
\textnormal{\textbf{Proof of statement 4.}}\ \ Statement follows straightforwardly by direct checking.\newline
\textnormal{\textbf{Proof of statement 5.}}\ \ Fix points $\unkuno_{1,0}\in \intuno_1$, $\unkuno_{2,0}\in \intuno_1$, a representation $\left(\setsymdue_{\unkuno_{1,0}},\intuno_{\unkuno_{1,0}}, \gfy_{\unkuno_{1,0}}\right)$ of the correspondence between $\pmV_1$ and $\corrarr[\pmV]_1$ at $\unkuno_{1,0}$, a representation $\left(\setsymdue_{\unkuno_{2,0}},\intuno_{\unkuno_{2,0}}, \gfy_{\unkuno_{2,0}}\right)$ of the correspondence between $\pmV_2$ and $\corrarr[\pmV]_2$ at $\unkuno_{2,0}$.\newline
Define:
\begin{flalign*}
&\begin{array}{l}
\setsymdue_{\left(\unkuno_{1,0},\unkuno_{2,0}\right)}=\setsymdue_{\unkuno_{1,0}}\times\setsymdue_{\unkuno_{2,0}}\text{;}
\end{array}\\[6pt]
&\begin{array}{l}
\intuno_{\left(\unkuno_{1,0},\unkuno_{2,0}\right)}=\intuno_{\unkuno_{1,0}}\times\intuno_{\unkuno_{2,0}}\text{;}
\end{array}\\[6pt]
&\begin{array}{l}
 \gfy_{\left(\unkuno_{1,0},\unkuno_{2,0}\right)}=\lboundgenf \gfy_{\unkuno_{1,0}}, \gfy_{\unkuno_{2,0}} \rboundgenf  \genfuncomp\\[4pt]
\hspace{40pt} \big(h\funcomp \incl{\setsymdue_{\left(\unkuno_{1,0},\unkuno_{2,0}\right)}\times \trasl[-\left(\unkuno_{1,0},\unkuno_{2,0}\right)]\left(\intuno_{\left(\unkuno_{1,0},\unkuno_{2,0}\right)}\right)}{\mathbb{R}^{l_1+l_2+m_1+m_2}}   \big)\text{.}
\end{array}
\end{flalign*}
Eventually statement follows since Proposition \ref{genpointprop}-[3] entails that the triplet $\left(\setsymdue_{\left(\unkuno_{1,0},\unkuno_{2,0}\right)}, \intuno_{\left(\unkuno_{1,0},\unkuno_{2,0}\right)} , \gfy_{\left(\unkuno_{1,0},\unkuno_{2,0}\right)}\right)$ is a representation of the correspondence between $\pmV_3$ and $\corrarr[\pmV]_3$ at $\left(\unkuno_{1,0},\unkuno_{2,0}\right)$.\newline
\textnormal{\textbf{Proof of statement 6.}}\ \  Statement follows straightforwardly by direct checking.\newline
\textnormal{\textbf{Proof of statement 7.}}\ \ Fix $\unkuno_0\in \intuno$, a representation $\left(\setsymdue_{\unkuno_0},\intuno_{\unkuno_0}, \gfy_{\unkuno_0}\right)$ of the correspondence between $\pmV$ and $\corrarr[\pmV]$ at $\unkuno_0$.\newline
Referrign to Notation \ref{realvec}-[6] define:
\begin{flalign*}
&\begin{array}{l}
k=\left\{
\begin{array}{ll}
l&\text{if}\hspace{4pt}\bsfm=\Fp_{\mathsf{l}}\hspace{4pt}\text{for}\hspace{4pt} \mathsf{l}\in \left\{1,...,l\right\}\text{,}\\[4pt]
l+1&\text{if}\hspace{4pt}\bsfm\in\left\{\Fq_{\mathsf{l}},\Fqq_{\mathsf{l}}\right\}\hspace{4pt}\text{for}\hspace{4pt} \mathsf{l}\in \left\{1,...,l\right\}\text{;}
\end{array}
\right.
\end{array}\\[8pt]
&\begin{array}{l}
\setsymtre_{\unkuno_0}=\left\{
\begin{array}{ll}
\setsymdue_{\unkuno_0}&\text{if}\hspace{4pt}\bsfm=\Fp_{\mathsf{l}}\hspace{4pt}\text{for}\hspace{4pt} \mathsf{l}\in \left\{1,...,l\right\}\text{,}\\[4pt]
\domint\left[\setsymdue_{\unkuno_0},\mathsf{l}\right]&\text{if}\hspace{4pt}\bsfm\in\left\{ \Fq_{\mathsf{l}},\Fqq_{\mathsf{l}}\right\}\hspace{4pt}\text{for}\hspace{4pt} \mathsf{l}\in \left\{1,...,l\right\}\text{;}
\end{array}
\right.
\end{array}\\[8pt]
&\begin{array}{l}
\gfz_{\unkuno_0}\in \genf(\argcompl{\setsymtre_{\unkuno_{0}}\times \left(\trasl[-\unkuno_{0}]\left(\intuno_{\unkuno_{0}}\right)\right)})+\mathbb{R}^{n}+\hspace{4pt}\text{by}\hspace{4pt}\text{setting}\hspace{4pt}
\gfz_{\unkuno_0}=\bsfm\genfbsfmu\gfy_{\unkuno_0}\text{.}
\end{array}
\end{flalign*}
Referring to Definition \ref{dualnot0} and performing computations we have 
\begin{equation}
\bsfm\genfbsfmu\left(\funcaap{\gfy}_{\unkuno_0}\left(0,\unkuno-\unkuno_0\right)\right)=\funcaap{\left(\bsfm\genfbsfmu\gfz_{\unkuno_0}\right)}\left(0,\unkuno-\unkuno_0\right)\hspace{15pt}\forall \unkuno \in \intuno_{\unkuno_0} \text{.}\label{pfcdf1}
\end{equation}
Eventually statement follows since \eqref{pfcdf1} entails that $\left(\setsymtre_{\unkuno_0},\intuno_{\unkuno_0}, \gfz_{\unkuno_0}\right)$ is a representation of the correspondence between $\bsfm\bsfmuquotpoint\pmV$ and $\bsfm\bsfmuquotpoint\corrarr[\pmV]$ at $\unkuno_0$.\newline
\textnormal{\textbf{Proof of statement 8.}}\ \  Fix $\unkuno_0\in \intuno$, a representation $\left(\setsymdue_{\unkuno_0},\intuno_{\unkuno_0}, \gfy_{\unkuno_0}\right)$ of the correspondence between $\pmV$ and $\corrarr[\pmV]$ at $\unkuno_0$.\newline
Referrign to Notation \ref{realvec}-[6] define:
\begin{flalign*}
&\left\{
\begin{array}{l}
h:\setsymdue_{\unkuno_0}\times     \trasl[-\unkuno_{0}]\left(\intuno_{\unkuno_0}\right)\rightarrow\mathbb{R}^l\times\mathbb{R}^m\hspace{4pt}\text{by}\hspace{4pt}\text{setting}\\[4pt]
h\left(\unkuno_1,\unkuno_2\right)=\left(0,\unkuno_2\right)\hspace{15pt}\forall\left(\unkuno_1,\unkuno_2\right)\in\setsymdue_{\unkuno_0}\times     \trasl[-\unkuno_{0}]\left(\intuno_{\unkuno_0}\right)\text{;}
\end{array}
\right.\\[8pt]
&\begin{array}{l}
\setsymtre_{\unkuno_0}=
\domint\left[\setsymdue_{\unkuno_0},\mathsf{l}\right]\text{;}
\end{array}\\[8pt]
&\left\{
\begin{array}{l}
\gfz_{\unkuno_0}\in \genf(\argcompl{\setsymtre_{\unkuno_{0}}\times \left(\trasl[-\unkuno_{0}]\left(\intuno_{\unkuno_{0}}\right)\right)})+\mathbb{R}^{n}+\hspace{4pt}\text{by}\hspace{4pt}\text{setting}\\[4pt]
\gfz_{\unkuno_0}=\Fint_{\mathsf{l}}\genfbsfmu\left(\gfy_{\unkuno_0}\genfunsum 
 -1 \genfunscalp\left(\gfy_{\unkuno_0}\genfuncomp h   \right)\right)\text{.}
\end{array}
\right.
\end{flalign*}
Statement follows by proving that  $\left(\setsymtre_{\unkuno_0},\intuno_{\unkuno_0}, \gfz_{\unkuno_0}\right)$ is a representation of the correspondence between $\Fint_{\mathsf{l}}\bsfmuquotpoint\pmV$ and $\Fint_{\mathsf{l}}\bsfmuquotpoint\left(\corrarr[\pmV] \sumquotpoint
-1 \genfunscalp \left( \corrarr[\pmV]\compquotpoint 
\lininclquot\left(\gczsfx[\argcompl{\lboundquotpoint   \cost<\mathbb{R}^l<>\mathbb{R}^l>+0+, \idobj+\mathbb{R}^m+ \rboundquotpoint}]\right)\right)\right)$ at $\unkuno_0$.\newline 
This is proven by the two chain of equalities below, where equalities (1) hold true by relation \eqref{(R 9 ter)}, equalities $(2)$ hold true by relation \eqref{(S 7 bis)}, equalities (3) hold true by relation \eqref{(R 9.3)}:
\begin{flalign*}
&\begin{array}{l}
\Fint_{\mathsf{l}}\bsfmuquotpoint\left(\pmV\left(\unkuno\right)\right)=\Fint_{\mathsf{l}}\bsfmuquotpoint\left(\ffuncaap{\gfy}_{\unkuno_0}\left(0,\unkuno-\unkuno_0\right)\compquotpoint \lininclquot\left(\gczsfx[\argcompl{\pointedincl<\argcompl{\left(\mathbb{R}^{l},0\right),\left(\mathbb{R}^{m},0\right)}<>1>}]\right)\right)\overset{1}{=}
\end{array}\\[4pt]
&\begin{array}{l}
\left(\Fint_{\mathsf{l}}\bsfmuquotpoint\left(\ffuncaap{\gfy}_{\unkuno_0}\left(0,\unkuno-\unkuno_0\right)\right)\right)\compquotpoint \lininclquot\left(\gczsfx[\argcompl{\pointedincl<\argcompl{\left(\mathbb{R}^{l+1},0\right),\left(\mathbb{R}^{m},0\right)}<>1>}]\right)=
\end{array}\\[4pt]
&\begin{array}{l} 
\bigg(\Fint_{\mathsf{l}}\bsfmuquotpoint\genfquotfun\Big(\gfy_{\unkuno_0} \genfuncomp \trasl[\left(0, \unkuno-\unkuno_0\right)] \genfunsum\\[4pt]
\hspace{50pt} 
 -1 \genfunscalp\left( \gfy_{\unkuno_0} \genfuncomp \trasl[\left(0, \unkuno-\unkuno_0\right)]  \genfuncomp \cost<\argcompl{\mathbb{R}^{l+m}}<>\argcompl{\mathbb{R}^{l+m}}>+0+
\right)\Big)\bigg)\compquotpoint \\[4pt]
\hspace{200pt}\lininclquot\left(\gczsfx[\argcompl{\pointedincl<\argcompl{\left(\mathbb{R}^{l+1},0\right),\left(\mathbb{R}^{m},0\right)}<>1>}]\right)=
\end{array}\\[4pt]
&\begin{array}{l}
\genfquotfun\Big(\Fint_{\mathsf{l}}\genfbsfmu\left(\gfy_{\unkuno_0} \genfuncomp \trasl[\left(0, \unkuno-\unkuno_0\right)]\right) \genfunsum\\[4pt]
\hspace{50pt} 
 -1 \genfunscalp\left( \Fint_{\mathsf{l}}\genfbsfmu\left(\gfy_{\unkuno_0} \genfuncomp \trasl[\left(0, \unkuno-\unkuno_0\right)]  \genfuncomp \cost<\argcompl{\mathbb{R}^{l+m}}<>\argcompl{\mathbb{R}^{l+m}}>+0+
\right)\right)\Big)\compquotpoint \\[4pt]
\hspace{200pt}\lininclquot\left(\gczsfx[\argcompl{\pointedincl<\argcompl{\left(\mathbb{R}^{l+1},0\right),\left(\mathbb{R}^{m},0\right)}<>1>}]\right)\overset{2}{=}
\end{array}\\[4pt]
&\begin{array}{l}
\genfquotfun\Big(\Fint_{\mathsf{l}}\genfbsfmu\left(\gfy_{\unkuno_0} \genfuncomp \trasl[\left(0, \unkuno-\unkuno_0\right)]\right) \genfunsum\\[4pt]
-1 \genfunscalp\left( \Fint_{\mathsf{l}}\genfbsfmu\left(\gfy_{\unkuno_0} \genfuncomp\left( h\genfuncomp \trasl[\left(0, \unkuno-\unkuno_0\right)]\right)
\right)\right)\genfunsum \Fint_{\mathsf{l}}\genfbsfmu\left(\gfy_{\unkuno_0} \genfuncomp \left( h\genfuncomp \trasl[\left(0, \unkuno-\unkuno_0\right)]\right)
\right)\genfunsum\\[4pt]
-1 \genfunscalp\left( \Fint_{\mathsf{l}}\genfbsfmu\left(\gfy_{\unkuno_0} \genfuncomp \trasl[\left(0, \unkuno-\unkuno_0\right)]  \genfuncomp \cost<\argcompl{\mathbb{R}^{l+m}}<>\argcompl{\mathbb{R}^{l+m}}>+0+
\right)\right)\Big)\compquotpoint \\[4pt]
\hspace{200pt}
\lininclquot\left(\gczsfx[\argcompl{\pointedincl<\argcompl{\left(\mathbb{R}^{l+1},0\right),\left(\mathbb{R}^{m},0\right)}<>1>}]\right)=
\end{array}\\[4pt]
&\begin{array}{l}
\bigg(\genfquotfun\Big(\Fint_{\mathsf{l}}\genfbsfmu\left(\gfy_{\unkuno_0} \genfuncomp \trasl[\left(0, \unkuno-\unkuno_0\right)] \genfunsum
-1 \genfunscalp\left( \Fint_{\mathsf{l}}\genfbsfmu\left(\left(\gfy_{\unkuno_0} \genfuncomp h\right)\genfuncomp \trasl[\left(0, \unkuno-\unkuno_0\right)]\right)
\right)\right)\Big)
\sumquotpoint\\[4pt]
 \genfquotfun\Big(\Fint_{\mathsf{l}}\genfbsfmu\left(\left(\gfy_{\unkuno_0} \genfuncomp   \trasl[\left(0, \unkuno-\unkuno_0\right)]\right)\genfuncomp h
\right)\genfunsum\\[4pt]
\hspace{25pt}
-1 \genfunscalp\left( \Fint_{\mathsf{l}}\genfbsfmu\left(\left(\gfy_{\unkuno_0} \genfuncomp \trasl[\left(0, \unkuno-\unkuno_0\right)]  \genfuncomp \cost<\argcompl{\mathbb{R}^{l+m}}<>\argcompl{\mathbb{R}^{l+m}}>+0+\right)\genfuncomp h
\right)\right)\Big)\bigg)\compquotpoint \\[4pt]
\hspace{200pt}
\lininclquot\left(\gczsfx[\argcompl{\pointedincl<\argcompl{\left(\mathbb{R}^{l+1},0\right),\left(\mathbb{R}^{m},0\right)}<>1>}]\right)=
\end{array}\\[4pt]
&\begin{array}{l}
\bigg(\genfquotfun\Big(\left(\Fint_{\mathsf{l}}\genfbsfmu\left(\gfy_{\unkuno_0}  \genfunsum
-1 \genfunscalp\left(\gfy_{\unkuno_0} \genfuncomp h\right)\right)\right)
\genfuncomp \trasl[\left(0, \unkuno-\unkuno_0\right)]\Big)
\sumquotpoint\\[4pt]
\genfquotfun\Big(
\Fint_{\mathsf{l}}\genfbsfmu\big(\left(\gfy_{\unkuno_0} \genfuncomp   \trasl[\left(0, \unkuno-\unkuno_0\right)]\right)\genfuncomp h
\genfunsum\\[4pt]
\hspace{25pt}
-1 \genfunscalp \left(\left(\gfy_{\unkuno_0} \genfuncomp \trasl[\left(0, \unkuno-\unkuno_0\right)]  \genfuncomp \cost<\argcompl{\mathbb{R}^{l+m}}<>\argcompl{\mathbb{R}^{l+m}}>+0+\right)\genfuncomp h
\right)\big)
\Big)
\bigg)\compquotpoint\\[4pt]
\hspace{200pt} 
\lininclquot\left(\gczsfx[\argcompl{\pointedincl<\argcompl{\left(\mathbb{R}^{l+1},0\right),\left(\mathbb{R}^{m},0\right)}<>1>}]\right)=
\end{array}\\[4pt]
&\begin{array}{l}
\left(\genfquotfun\left(\prefuncaap{\gfz}_{\unkuno_0}\left(\unkuno-\unkuno_0\right)\right)\right)\compquotpoint 
\lininclquot\left(\gczsfx[\argcompl{\pointedincl<\argcompl{\left(\mathbb{R}^{l+1},0\right),\left(\mathbb{R}^{m},0\right)}<>1>}]\right)
\sumquotpoint\\[4pt]
\bigg(\Fint_{\mathsf{l}}\bsfmuquotpoint\genfquotfun\Big(\big(\left(\gfy_{\unkuno_0} \genfuncomp   \trasl[\left(0, \unkuno-\unkuno_0\right)]\right)\genfunsum\\[4pt]
\hspace{30pt}
-1 \genfunscalp\left(\gfy_{\unkuno_0} \genfuncomp \trasl[\left(0, \unkuno-\unkuno_0\right)]  \genfuncomp \cost<\argcompl{\mathbb{R}^{l+m}}<>\argcompl{\mathbb{R}^{l+m}}>+0+\right)
\big)\genfuncomp h \Big)\bigg)\compquotpoint\\[4pt]
\hspace{190pt}  
\lininclquot\left(\gczsfx[\argcompl{\pointedincl<\argcompl{\left(\mathbb{R}^{l+1},0\right),\left(\mathbb{R}^{m},0\right)}<>1>}]\right)
\overset{3}{=}
\end{array}\\[4pt]
&\begin{array}{l}
\left(\genfquotfun\left(\prefuncaap{\gfz}_{\unkuno_0}\left(\unkuno-\unkuno_0\right) \genfunsum
-1 \genfunscalp\left(\prefuncaap{\gfz}_{\unkuno_0}\left(\unkuno-\unkuno_0\right) \genfuncomp \cost<\argcompl{\mathbb{R}^{l+1+m}}<>\argcompl{\mathbb{R}^{l+1+m}}>+0+\right)    \right)\right)\compquotpoint \\[4pt]
\lininclquot\left(\gczsfx[\argcompl{\pointedincl<\argcompl{\left(\mathbb{R}^{l+1},0\right),\left(\mathbb{R}^{m},0\right)}<>1>}]\right)\sumquotpoint
\Fint_{\mathsf{l}}\bsfmuquotpoint\genfquotfun\Big(\big(
\funcaap{\gfy}_{\unkuno_0}\left(\unkuno-\unkuno_0\right)\genfuncomp h
\big)\genfuncomp \pointedincl<\argcompl{\left(\mathbb{R}^{l},0\right),\left(\mathbb{R}^{m},0\right)}<>1> \Big) =
\end{array}\\[4pt]
&\begin{array}{l}
\left(\ffuncaap{\gfz}_{\unkuno_0}\left(\unkuno-\unkuno_0\right)\right)\compquotpoint 
\lininclquot\left(\gczsfx[\argcompl{\pointedincl<\argcompl{\left(\mathbb{R}^{l+1},0\right),\left(\mathbb{R}^{m},0\right)}<>1>}]\right)\sumquotpoint
\Fint_{\mathsf{l}}\bsfmuquotpoint\genfquotfun\Big(0 \genfunscalp \left(\gfy_{\unkuno_0} \genfuncomp \trasl[\left(0, \unkuno-\unkuno_0\right)] \right) \Big)=
\end{array}\\[4pt]
&\begin{array}{l}
\left(\ffuncaap{\gfz}_{\unkuno_0}\left(\unkuno-\unkuno_0\right)\right)\compquotpoint 
\lininclquot\left(\gczsfx[\argcompl{\pointedincl<\argcompl{\left(\mathbb{R}^{l+1},0\right),\left(\mathbb{R}^{m},0\right)}<>1>}]\right)\hspace{15pt}\forall \unkuno \in \intuno_{\unkuno_0}\text{;}
\end{array}\\[20pt]
&\begin{array}{l}
\Fint_{\mathsf{l}}\bsfmuquotpoint\left(\corrarrre[\pmV]\left(\unkuno\right)\right)=\Fint_{\mathsf{l}}\bsfmuquotpoint\left(\ffuncaap{\gfy}_{\unkuno_0}\left(0,\unkuno-\unkuno_0\right)\right)=
\end{array}\\[4pt]
&\begin{array}{l} 
\Fint_{\mathsf{l}}\bsfmuquotpoint\genfquotfun\Big(\gfy_{\unkuno_0} \genfuncomp \trasl[\left(0, \unkuno-\unkuno_0\right)] \genfunsum
 -1 \genfunscalp\left( \gfy_{\unkuno_0} \genfuncomp \trasl[\left(0, \unkuno-\unkuno_0\right)]  \genfuncomp \cost<\argcompl{\mathbb{R}^{l+m}}<>\argcompl{\mathbb{R}^{l+m}}>+0+
\right)\Big)= 
\end{array}\\[4pt]
&\begin{array}{l}
\genfquotfun\Big(\Fint_{\mathsf{l}}\genfbsfmu\left(\gfy_{\unkuno_0} \genfuncomp \trasl[\left(0, \unkuno-\unkuno_0\right)]\right) \genfunsum 
 -1 \genfunscalp\left( \Fint_{\mathsf{l}}\genfbsfmu\left(\gfy_{\unkuno_0} \genfuncomp \trasl[\left(0, \unkuno-\unkuno_0\right)]  \genfuncomp \cost<\argcompl{\mathbb{R}^{l+m}}<>\argcompl{\mathbb{R}^{l+m}}>+0+
\right)\right)\Big)\overset{2}{=} 
\end{array}\\[4pt]
&\begin{array}{l}
\genfquotfun\Big(\Fint_{\mathsf{l}}\genfbsfmu\left(\gfy_{\unkuno_0} \genfuncomp \trasl[\left(0, \unkuno-\unkuno_0\right)]\right) \genfunsum\\[4pt]
-1 \genfunscalp\left( \Fint_{\mathsf{l}}\genfbsfmu\left(\gfy_{\unkuno_0} \genfuncomp\left( h\funcomp \trasl[\left(0, \unkuno-\unkuno_0\right)]\right)
\right)\right)\genfunsum \Fint_{\mathsf{l}}\genfbsfmu\left(\gfy_{\unkuno_0} \genfuncomp \left( h\funcomp \trasl[\left(0, \unkuno-\unkuno_0\right)]\right)
\right)\genfunsum\\[4pt]
\hspace{115pt}
-1 \genfunscalp\left( \Fint_{\mathsf{l}}\genfbsfmu\left(\gfy_{\unkuno_0} \genfuncomp \trasl[\left(0, \unkuno-\unkuno_0\right)]  \genfuncomp \cost<\argcompl{\mathbb{R}^{l+m}}<>\argcompl{\mathbb{R}^{l+m}}>+0+
\right)\right)\Big)=
\end{array}\\[4pt]
&\begin{array}{l}
\genfquotfun\Big(\Fint_{\mathsf{l}}\genfbsfmu\left(\gfy_{\unkuno_0} \genfuncomp \trasl[\left(0, \unkuno-\unkuno_0\right)] \genfunsum
-1 \genfunscalp\left( \Fint_{\mathsf{l}}\genfbsfmu\left(\left(\gfy_{\unkuno_0} \genfuncomp h\right)\genfuncomp \trasl[\left(0, \unkuno-\unkuno_0\right)]\right)
\right)\right)\Big)
\sumquotpoint\\[4pt]
 \genfquotfun\Big(\Fint_{\mathsf{l}}\genfbsfmu\left(\left(\gfy_{\unkuno_0} \genfuncomp   \trasl[\left(0, \unkuno-\unkuno_0\right)]\right)\genfuncomp h
\right)\genfunsum\\[4pt]
\hspace{85pt}
-1 \genfunscalp\left( \Fint_{\mathsf{l}}\genfbsfmu\left(\left(\gfy_{\unkuno_0} \genfuncomp \trasl[\left(0, \unkuno-\unkuno_0\right)]  \genfuncomp \cost<\argcompl{\mathbb{R}^{l+m}}<>\argcompl{\mathbb{R}^{l+m}}>+0+\right)\genfuncomp h
\right)\right)\Big)=
\end{array}\\[4pt]
&\begin{array}{l}
\genfquotfun\Big(\left(\Fint_{\mathsf{l}}\genfbsfmu\left(\gfy_{\unkuno_0}  \genfunsum
-1 \genfunscalp\left(\gfy_{\unkuno_0} \genfuncomp h\right)\right)\right)
\genfuncomp \trasl[\left(0, \unkuno-\unkuno_0\right)]\Big)
\sumquotpoint\\[4pt]
\genfquotfun\Big(
\Fint_{\mathsf{l}}\genfbsfmu\big(\left(\gfy_{\unkuno_0} \genfuncomp   \trasl[\left(0, \unkuno-\unkuno_0\right)]\right)\genfuncomp h
\genfunsum\\[4pt]
\hspace{110pt}
-1 \genfunscalp \left(\left(\gfy_{\unkuno_0} \genfuncomp \trasl[\left(0, \unkuno-\unkuno_0\right)]  \genfuncomp \cost<\argcompl{\mathbb{R}^{l+m}}<>\argcompl{\mathbb{R}^{l+m}}>+0+\right)\genfuncomp h
\right)\big)
\Big)=
\end{array}\\[4pt]
&\begin{array}{l}
\genfquotfun\left(\prefuncaap{\gfz}_{\unkuno_0}\left(\unkuno-\unkuno_0\right)\right)
\sumquotpoint
\Fint_{\mathsf{l}}\bsfmuquotpoint\genfquotfun\Big(\big(\left(\gfy_{\unkuno_0} \genfuncomp   \trasl[\left(0, \unkuno-\unkuno_0\right)]\right)\genfunsum\\[4pt]
\hspace{120pt}
-1 \genfunscalp\left(\gfy_{\unkuno_0} \genfuncomp \trasl[\left(0, \unkuno-\unkuno_0\right)]  \genfuncomp \cost<\argcompl{\mathbb{R}^{l+m}}<>\argcompl{\mathbb{R}^{l+m}}>+0+\right)
\big)\genfuncomp h \Big)
\overset{3}{=}
\end{array}\\[4pt]
&\begin{array}{l}
\genfquotfun\left(\prefuncaap{\gfz}_{\unkuno_0}\left(\unkuno\right) \genfunsum
-1 \genfunscalp\left(\prefuncaap{\gfz}_{\unkuno_0}\left(\unkuno\right) \genfuncomp \cost<\argcompl{\mathbb{R}^{l+1+m}}<>\argcompl{\mathbb{R}^{l+1+m}}>+0+\right)    \right)\sumquotpoint\\[4pt]
 \hspace{230pt}\Fint_{\mathsf{l}}\bsfmuquotpoint\genfquotfun\left(\funcaap{\gfy}_{\unkuno_0}\left(\unkuno-\unkuno_0\right)\genfuncomp h \right) =
\end{array}\\[4pt]
&\begin{array}{l}
\ffuncaap{\gfz}_{\unkuno_0}\left(\unkuno-\unkuno_0\right)\sumquotpoint
\Fint_{\mathsf{l}}\bsfmuquotpoint\left( \corrarr[\pmV]\compquotpoint 
\lininclquot\left(\gczsfx[\argcompl{\lboundquotpoint   \cost<\mathbb{R}^l<>\mathbb{R}^l>+0+, \idobj+\mathbb{R}^m+ \rboundquotpoint}]\right)\right)
\hspace{15pt}\forall \unkuno \in \intuno_{\unkuno_0}\text{.}
\end{array}
\end{flalign*}
\end{proof}

\section{Higher order generalized differential operators on germs of generalized functions}

In this section we study how germs of generalized derivations compose. We refer to Notation \ref{varispder},  Definition \ref{defcatD}-[5],  Remark \ref{noM}.\\[12pt]

In Definition \ref{dirsumpreclosdef} below we introduce sets which will be used to index local fields of generalized derivations to be iterated in order to get higher order generalized differential operators. We refer to Notation \ref{int}.

\begin{definition}\label{dirsumpreclosdef}\mbox{}
\begin{enumerate} 
\item Fix $\ls\in \mathbb{N}$, $\ddgtv, \dgtv, \uedgtv\in \mathbb{N}_0$, $\ddgtv<seq<=\left(\ddgtv>\lss>\right)_{\lss=1}^{\ls}\subseteq{\mathbb{N}_0}^{\ls}$, $\ouol<seq<=\left(\ouol>\lss>\right)_{\lss=1}^{\ls}\subseteq\left\{0,1\right\}^{\ls}$. Assume that $\ddgtv=\underset{\lss=1}{\overset{\ls}{\sum}}\,\ddgtv>\lss>$.\newline
We say that $\left(\ls, \ddgtv, \dgtv, \uedgtv,\ddgtv<seq<, \ouol<seq<\right)$ is an admissible six-tuple.
\item Fix an admissible six-tuple $\left(\ls, \ddgtv, \dgtv, \uedgtv,\ddgtv<seq<, \ouol<seq<\right)$, $\intunosk\subseteq\mathbb{R}^{\ddgtv+\dgtv}$, $\unkunosk\in\mathbb{R}^{\ddgtv+\dgtv}$.\newline
We set:
\begin{flalign*}
&\begin{array}{l}
\intunosk>\lss, \ddgtv<seq<>=\proj<\left\{\mathbb{R}^{\ddgtv>1>},...,\mathbb{R}^{\ddgtv>\ls>},\mathbb{R}^{\dgtv}\right\}<>\lss>\left(\intuno\right)\hspace{15pt}\forall\lss\in \left\{1,...,\ls\right\}\text{;}
\end{array}\\[8pt]
&\begin{array}{l}
\intunosk>\ddgtv>=\proj<\left\{\mathbb{R}^{\ddgtv},\mathbb{R}^{\dgtv}\right\}<>1>\left(\intuno\right)\text{;}
\end{array}\\[8pt]
&\begin{array}{l}
\intunosk>\dgtv>=\proj<\left\{\mathbb{R}^{\ddgtv},\mathbb{R}^{\dgtv}\right\}<>2>\left(\intuno\right)\text{;}
\end{array}\\[8pt]
&\begin{array}{l}
\unkunosk>\lss, \ddgtv<seq<>=\proj<\left\{\mathbb{R}^{\ddgtv>1>},...,\mathbb{R}^{\ddgtv>\ls>},\mathbb{R}^{\dgtv}\right\}<>\lss>\left(\unkuno\right)\hspace{15pt}\forall\lss\in \left\{1,...,\ls\right\}\text{;}
\end{array}\\[8pt]
&\begin{array}{l}
\unkunosk>\ddgtv>=\proj<\left\{\mathbb{R}^{\ddgtv},\mathbb{R}^{\dgtv}\right\}<>1>\left(\unkuno\right)\text{;}
\end{array}\\[8pt]
&\begin{array}{l}
\unkunosk>\dgtv>=\proj<\left\{\mathbb{R}^{\ddgtv},\mathbb{R}^{\dgtv}\right\}<>2>\left(\unkuno\right)\text{.}
\end{array}
\end{flalign*}
\item We define  
\begin{multline*}
\AQ=
\Big\{\left\{\left(\ls>\antelss>, \ddgtv>\antelss>, \dgtv>\antelss>, \uedgtv>\antelss>,\ddgtv<seq<>\antelss>, \ouol<seq<>\antelss>\right)\right\}_{\antelss=1}^{\antels}\;:\;\antels\in \mathbb{N}\text{,}\hspace{4pt}
 \left(\ls>\antelss>, \ddgtv>\antelss>, \dgtv>\antelss>, \uedgtv>\antelss>,\ddgtv<seq<>\antelss>, \ouol<seq<>\antelss>\right)\hspace{4pt}\text{is}\\[4pt]
\text{an}\hspace{4pt}\text{admissible}\hspace{4pt}
\text{six-tuple}\hspace{4pt}\text{for}\hspace{4pt}\text{any}\hspace{4pt}\antelss\in \left\{1,...,\antels\right\}\Big\}\text{.}
\end{multline*}
Fix $\ls, \antels\in \mathbb{N}$, $\ddgtv, \dgtv,\uedgtv\in \mathbb{N}_0$. We define  
\begin{multline*}
\AQ/\ls/-\antels-+\argcompl{\ddgtv,\dgtv}+<\uedgtv<=
\Big\{\left\{\left(\ls, \ddgtv, \dgtv, \uedgtv,\ddgtv<seq<>\antelss>, \ouol<seq<>\antelss>\right)\right\}_{\antelss=1}^{\antels}\;:\;
 \left(\ls, \ddgtv, \dgtv, \uedgtv,\ddgtv<seq<>\antelss>, \ouol<seq<>\antelss>\right)\hspace{4pt}\text{is}\hspace{4pt}\text{an}\\[4pt]
\text{admissible}\hspace{4pt}
\text{six-tuple}\hspace{4pt}\text{for}\hspace{4pt}\text{any}\hspace{4pt}\antelss\in \left\{1,...,\antels\right\}\Big\}\text{.}
\end{multline*}
Fix $ \antels\in \mathbb{N}$, $\ddgtv, \dgtv,\uedgtv\in \mathbb{N}_0$. We define  
$\AQ-\antels-+\argcompl{\ddgtv,\dgtv}+<\uedgtv<=\underset{\ls\in \mathbb{N}}{\bigcup}\AQ/\ls/-\antels-+\argcompl{\ddgtv,\dgtv}+<\uedgtv<$.\newline
Fix $ \antels\in \mathbb{N}$, $ \dgtv,\uedgtv\in \mathbb{N}_0$. We define  
$\AQ-\antels-+\argcompl{\segnvar,\dgtv}+<\uedgtv<=\underset{\ddgtv\in \mathbb{N}_0}{\bigcup}\AQ-\antels-+\argcompl{\ddgtv,\dgtv}+<\uedgtv<$.\newline
Fix $ \antels\in \mathbb{N}$, $ \uedgtv\in \mathbb{N}_0$. We define  
$\AQ-\antels-<\uedgtv<= \underset{\dgtv\in \mathbb{N}_0}{\bigcup}\AQ-\antels-+\argcompl{\segnvar,\dgtv}+<\uedgtv<$.\newline
Fix $ \antels\in \mathbb{N}_0$. We define  
$\AQ-\antels-= \underset{\uedgtv\in \mathbb{N}_0}{\bigcup} \AQ-\antels-<\uedgtv<$.\newline
Fix $ \uedgtv\in \mathbb{N}_0$. We define  
$\AQ<\uedgtv<= \underset{\antels\in \mathbb{N}}{\bigcup}\AQ-\antels-<\uedgtv<$.\newline
Fix $\ls, \antels\in \mathbb{N}$, $\tddmsk,\uedgtv\in \mathbb{N}_0$. We define 
\begin{multline*}
\AQ/\ls/-\antels-+\tddmsk+<\uedgtv<=
\Big\{\left\{\left(\ls, \ddgtv>\antelss>, \dgtv>\antelss>, \uedgtv,\ddgtv<seq<>\antelss>, \ouol<seq<>\antelss>\right)\right\}_{\antelss=1}^{\antels}\;:\;
 \left(\ls, \ddgtv, \dgtv, \uedgtv,\ddgtv<seq<>\antelss>, \ouol<seq<>\antelss>\right)\hspace{4pt}\text{is}\hspace{4pt}\text{an}\\[4pt]
\text{admissible}\hspace{4pt}\text{six-tuple}\hspace{4pt}\text{with}\hspace{4pt}\ddgtv>\antelss>+ \dgtv>\antelss>=\tddmsk\hspace{4pt}\text{for}\hspace{4pt}\text{any}\hspace{4pt}\antelss\in \left\{1,...,\antels\right\}\Big\}\text{.}
\end{multline*}
Fix $\antels\in \mathbb{N}$, $\tddmsk,\uedgtv\in \mathbb{N}_0$. We define  
$\AQ-\antels-+\tddmsk+<\uedgtv<= \underset{\ls\in \mathbb{N}}{\bigcup} \AQ/\ls/-\antels-+\tddmsk+<\uedgtv<$.\newline
Fix $\antels\in \mathbb{N}$, $\tddmsk\in \mathbb{N}_0$. We define  
$\AQ-\antels-+\tddmsk+= \underset{\uedgtv\in \mathbb{N}_0}{\bigcup} \AQ-\antels-+\tddmsk+<\uedgtv<$.\newline
Fix $\tddmsk,\uedgtv\in \mathbb{N}_0$. We define  $\AQ+\tddmsk+<\uedgtv<=  \underset{\antels\in \mathbb{N}}{\bigcup} \AQ-\antels-+\tddmsk+<\uedgtv< $.\newline
Fix $\tddmsk\in \mathbb{N}_0$. We define  
$\AQ+\tddmsk+= \underset{\uedgtv\in \mathbb{N}_0}{\bigcup} \AQ+\tddmsk+<\uedgtv<$.
\item Fix $\dgtv>1>,\dgtv>2>,\uedgtv\in \mathbb{N}_0$. We define the set function  
\begin{equation*}
\consec-\dgtv>1>-(\dgtv>2>)+\uedgtv+\,:\AQ-1-+\argcompl{\segnvar,\dgtv>1>}+<\argcompl{\dgtv>2>+\uedgtv}<\times \AQ-1-+\argcompl{\segnvar,\dgtv>2>}+<\uedgtv<\rightarrow\AQ-1-+\argcompl{\segnvar,\dgtv>1>+\dgtv>2>}+<\uedgtv<
\end{equation*}
by setting:
\begin{equation*}
\begin{array}{l}
\left(\ls>1>, \ddgtv>1>, \dgtv>1>, \uedgtv>1>, \ddgtv<seq<>1>,\ouol<seq<>1>\right)\consec-\dgtv>1>-(\dgtv_1)+\uedgtv+
\left(\ls>2>, \ddgtv>2>, \dgtv>2>, \uedgtv>2>, \ddgtv<seq<>2>,\ouol<seq<>2>\right)
=\left(\ls>3>, \ddgtv>3>, \dgtv>3>, \uedgtv>3>, \ddgtv<seq<>3>,\ouol<seq<>3>\right)
\\[6pt] 
\hspace{15pt}\forall\left(\left(\ls>1>, \ddgtv>1>, \dgtv>1>, \uedgtv>1>, \ddgtv<seq<>1>,\ouol<seq<>1>\right),
\left(\ls>2>, \ddgtv>2>, \dgtv>2>, \uedgtv>2>, \ddgtv<seq<>2>,\ouol<seq<>2>\right)\right)\in \AQ-1-+\argcompl{\segnvar,\dgtv>1>}+<\argcompl{\uedgtv>1>}<\times \AQ-1-+\argcompl{\segnvar,\dgtv>2>}+<\argcompl{\uedgtv>2>}<
\end{array}
\end{equation*}
where: 
\begin{flalign*}
&\begin{array}{l}
\ls>3>=\ls>1>+\ls>2>\text{;}
\end{array}\\[10pt]
&\begin{array}{l}
\ddgtv>3>=\ddgtv>1>+\ddgtv>2>\text{;}
\end{array}\\[10pt]
&\begin{array}{l}
\dgtv>3>=\dgtv>1>+\dgtv>2>\text{;}
\end{array}\\[10pt]
&\begin{array}{l}
\uedgtv>3>=\uedgtv>2>\text{;}
\end{array}\\[10pt]
& \begin{array}{l}
\ddgtv<seq<>3>=\left( \ddgtv>3,\lss> \right)_{\lss=1}^{\ls>3>}\hspace{4pt}\text{is}\hspace{4pt}\text{given}\hspace{4pt}\text{by}\\[4pt]
\hspace{90pt}\ddgtv>3,\lss>=\left\{
\begin{array}{ll}
\ddgtv>1,\lss>&\text{if}\hspace{4pt}\lss\in\left\{1,..., \ls>1>\right\}\text{,}\\[4pt]
\ddgtv>\argcompl{2,\lss-\ls>1>}>&\text{if}\hspace{4pt}\lss\in\left\{\ls>1>+1,...,\ls>1>+\ls>2>\right\}\text{;}
\end{array}
\right.
\end{array}\\[10pt]
& \begin{array}{l}
\ouol<seq<>3>=\left(  \ouol>3,\lss>    \right)_{\lss=1}^{\ls>3>}\hspace{4pt}\text{is}\hspace{4pt}\text{given}\hspace{4pt}\text{by}\\[4pt]
\hspace{90pt}\ouol>3,\lss>=\left\{
\begin{array}{ll}
\ouol>1,\lss>&\text{if}\hspace{4pt}\lss\in\left\{1,..., \ls>1>\right\}\text{,}\\[4pt]
\ouol>\argcompl{2,\lss-\ls>1>}>&\text{if}\hspace{4pt}\lss\in\left\{\ls>1>+1,...,\ls>1>+\ls>2>\right\}\text{.}
\end{array}
\right.
\end{array}
\end{flalign*}
The three-graded family of set functions $\left\{\consec-\dgtv>1>-(\dgtv>2>)+\uedgtv+\right\}_{\left(\dgtv>1>,\dgtv>2>, \uedgtv\right)\in \mathbb{N}_0\times\mathbb{N}_0\times\mathbb{N}_0}$ defines an internal operation $\consec$ on the bi-graded set $\left\{\AQ-1-+\argcompl{\segnvar,\dgtv}+<\uedgtv<\right\}_{\left(\dgtv,\uedgtv\right)\in \mathbb{N}_0\times\mathbb{N}_0}$.
\item We define the set function  
\begin{equation*}
\anteconsec\,:\AQ\times \AQ\rightarrow\AQ
\end{equation*}
by setting: 
\begin{equation*}
\begin{array}{l}
\left\{\left(\ls>1,\antelss>, \ddgtv>1,\antelss>, \dgtv>1,\antelss>, \uedgtv>1,\antelss>,\ddgtv<seq<>1,\antelss>, \ouol<seq<>1\antelss>\right)\right\}_{\antelss=1}^{\antels>1>}\anteconsec\\[4pt]
\hspace{75pt}\left\{\left(\ls>2,\antelss>, \ddgtv>2,\antelss>, \dgtv>2,\antelss>, \uedgtv>2,\antelss>,\ddgtv<seq<>2,\antelss>, \ouol<seq<>2\antelss>\right)\right\}_{\antelss=1}^{\antels>2>}=\\[4pt]
\hspace{150pt}\left\{\left(\ls>3,\antelss>, \ddgtv>3,\antelss>, \dgtv>3,\antelss>, \uedgtv>3,\antelss>,\ddgtv<seq<>3,\antelss>, \ouol<seq<>3\antelss>\right)\right\}_{\antelss=1}^{\antels>1>+\antels>2>}\\[4pt]
\text{for}\hspace{4pt} \text{any}\\[4pt]
\left\{\left(\ls>1,\antelss>, \ddgtv>1,\antelss>, \dgtv>1,\antelss>, \uedgtv>1,\antelss>,\ddgtv<seq<>1,\antelss>, \ouol<seq<>1\antelss>\right)\right\}_{\antelss=1}^{\antels>1>},
\left\{\left(\ls>2,\antelss>, \ddgtv>2,\antelss>, \dgtv>2,\antelss>, \uedgtv>2,\antelss>,\ddgtv<seq<>2,\antelss>, \ouol<seq<>2\antelss>\right)\right\}_{\antelss=1}^{\antels>2>}\in  \AQ\text{,}
\end{array}
\end{equation*}
where:
\begin{flalign*}
&\begin{array}{l}
\ls>3, \antelss>=\left\{
\begin{array}{ll}
\ls>1, \antelss>&\text{if}\hspace{4pt}\antelss\in \left\{1,...,\antels>1>\right\}\text{,}\\[4pt]
\ls>\argcompl{2, \antelss-\antels>1>}>&\text{if}\hspace{4pt}\antelss\in \left\{\antels>1>+1,...,\antels>1>+\antels>2>\right\}\text{;}
\end{array}
\right.
\end{array}\\[10pt]
&\begin{array}{l}
\ddgtv>3, \antelss>=\left\{
\begin{array}{ll}
\ddgtv>1, \antelss>&\text{if}\hspace{4pt}\antelss\in \left\{1,...,\antels>1>\right\}\text{,}\\[4pt]
\ddgtv>\argcompl{2, \antelss-\antels>1>}>&\text{if}\hspace{4pt}\antelss\in \left\{\antels>1>+1,...,\antels>1>+\antels>2>\right\}\text{;}
\end{array}
\right.
\end{array}\\[10pt]
&\begin{array}{l}
\dgtv>3, \antelss>=\left\{
\begin{array}{ll}
\dgtv>1, \antelss>&\text{if}\hspace{4pt}\antelss\in \left\{1,...,\antels>1>\right\}\text{,}\\[4pt]
\dgtv>\argcompl{2, \antelss-\antels>1>}>&\text{if}\hspace{4pt}\antelss\in \left\{\antels>1>+1,...,\antels>1>+\antels>2>\right\}\text{;}
\end{array}
\right.
\end{array}\\[10pt]
&\begin{array}{l}
\uedgtv>3, \antelss>=\left\{
\begin{array}{ll}
\uedgtv>1, \antelss>&\text{if}\hspace{4pt}\antelss\in \left\{1,...,\antels>1>\right\}\text{,}\\[4pt]
\uedgtv>\argcompl{2, \antelss-\antels>1>}>&\text{if}\hspace{4pt}\antelss\in \left\{\antels>1>+1,...,\antels>1>+\antels>2>\right\}\text{;}
\end{array}
\right.
\end{array}\\[10pt]
&\begin{array}{l}
\ddgtv<seq<>3, \antelss>=\left\{
\begin{array}{ll}
\ddgtv<seq<>1, \antelss>&\text{if}\hspace{4pt}\antelss\in \left\{1,...,\antels>1>\right\}\text{,}\\[4pt]
\ddgtv<seq<>\argcompl{2, \antelss-\antels>1>}>&\text{if}\hspace{4pt}\antelss\in \left\{\antels>1>+1,...,\antels>1>+\antels>2>\right\}\text{;}
\end{array}
\right.
\end{array}\\[10pt]
&\begin{array}{l}
\ouol<seq<>3, \antelss>=\left\{
\begin{array}{ll}
\ouol<seq<>1, \antelss>&\text{if}\hspace{4pt}\antelss\in \left\{1,...,\antels>1>\right\}\text{,}\\[4pt]
\ouol<seq<>\argcompl{2, \antelss-\antels>1>}>&\text{if}\hspace{4pt}\antelss\in \left\{\antels>1>+1,...,\antels>1>+\antels>2>\right\}\text{;}
\end{array}
\right.
\end{array}
\end{flalign*}
\end{enumerate}
\end{definition}

In Definition \ref{ChDiR1} below we algebraically formalize the composition of germs of generalized derivations.
We refer to Notations \ref{int}, \ref{ins}-[10, 11, 12], \ref{alg}-[12, 15], Definition \ref{dirsumpreclosdef}, Remarks \ref{siterem}-[3], \ref{alssitesrem}, \ref{pathinlgenerder}-[9].

\begin{definition}\label{ChDiR1}\mbox{}
\begin{enumerate}
\item We define: 
\begin{flalign*}
&\left\{\begin{array}{l}
\predcR[\left\{\left(\ls>\antelss>, \ddgtv>\antelss>, \dgtv>\antelss>, \uedgtv>\antelss>,\ddgtv<seq<>\antelss>, \ouol<seq<>\antelss>\right)\right\}_{\antelss=1}^{\antels}]=\overset{\antels}{\underset{\antelss=1}{\bigtensprodR}}\left(\overset{\ls>\antelss>}{\underset{\lss=1}{\bigtensprodR}} \,\generderspace-funtore-(\argcompl{\ddgtv>\antelss,\lss>})+\argcompl{\dgtv>\antelss>+\uedgtv>\antelss>}+\left(\dgtv>\antelss>\right)\right)\\[12pt]
\text{for}\hspace{4pt}\text{any}\hspace{4pt}\antelss\in \mathbb{N}\text{,}\hspace{4pt}\left\{\left(\ls>\antelss>, \ddgtv>\antelss>, \dgtv>\antelss>, \uedgtv>\antelss>,\ddgtv<seq<>\antelss>, \ouol<seq<>\antelss>\right)\right\}_{\antelss=1}^{\antels}\in\AQ-\antels-\text{;}
\end{array}
\right.\\[8pt]
&\begin{array}{l}
\predcR-\antels-+\argcompl{\ddgtv,\dgtv}+<\uedgtv<=
\underset{\elAQ\in\AQ-\antels-+\argcompl{\ddgtv,\dgtv}+<\uedgtv<
}{\bigdirsum}\predcR[\elAQ]\hspace{4pt}\text{for}\hspace{4pt}\text{any}\hspace{4pt}\antels\in\mathbb{N}  \text{,}\hspace{4pt}\ddgtv,\dgtv,\uedgtv\in \mathbb{N}_0\text{;}
\end{array}\\[8pt]
&\begin{array}{l}
\predcR-\antels-+\argcompl{\segnvar,\dgtv}+<\uedgtv<=
\left\{\predcR-\antels-+\argcompl{\ddgtv,\dgtv}+<\uedgtv<\right\}_{\ddgtv\in \mathbb{N}_0}\hspace{4pt}\text{for}\hspace{4pt}\text{any}\hspace{4pt}\antels\in\mathbb{N}  \text{,}\hspace{4pt}\dgtv,\uedgtv\in \mathbb{N}_0\text{;}
\end{array}\\[8pt]
&\begin{array}{l}
\predcR-\antels-+\argcompl{\segnvar,\segnvar}+<\uedgtv<=
\left\{\predcR-\antels-+\argcompl{\segnvar,\dgtv}+<\uedgtv<\right\}_{\left(\ddgtv,\dgtv\right)\in \mathbb{N}_0\times\mathbb{N}_0}\hspace{4pt}\text{for}\hspace{4pt}\text{any}\hspace{4pt}\antels\in\mathbb{N}  \text{,}\hspace{4pt}\uedgtv\in \mathbb{N}_0\text{;}
\end{array}\\[8pt]
&\begin{array}{l}
\predcR-\antels-+\argcompl{\segnvar,\segnvar}+<\segnvar<=
\left\{\predcR-\antels-+\argcompl{\ddgtv,\dgtv}+<\uedgtv<\right\}_{\left(\ddgtv,\dgtv,\uedgtv\right)\in \mathbb{N}_0\times\mathbb{N}_0\times\mathbb{N}_0}\hspace{4pt}\text{for}\hspace{4pt}\text{any}\hspace{4pt}\antels\in\mathbb{N}\text{;}
\end{array}\\[8pt]
&\begin{array}{l}
\predcR+\argcompl{\ddgtv,\dgtv}+<\uedgtv<=
\underset{\elAQ\in\AQ+\argcompl{\ddgtv,\dgtv}+<\uedgtv<
}{\bigdirsum}\predcR[\elAQ]\hspace{4pt}\text{for}\hspace{4pt}\text{any}\hspace{4pt}\ddgtv,\dgtv,\uedgtv\in \mathbb{N}_0\text{;}
\end{array}\\[8pt]
&\begin{array}{l}
\predcR+\argcompl{\segnvar,\dgtv}+<\uedgtv<=
\left\{\predcR+\argcompl{\ddgtv,\dgtv}+<\uedgtv<\right\}_{\ddgtv\in \mathbb{N}_0}\hspace{4pt}\text{for}\hspace{4pt}\text{any}\hspace{4pt}\antels\in\mathbb{N}  \text{,}\hspace{4pt}\dgtv,\uedgtv\in \mathbb{N}_0\text{;}
\end{array}\\[8pt]
&\begin{array}{l}
\predcR+\argcompl{\segnvar,\segnvar}+<\uedgtv<=
\left\{\predcR+\argcompl{\ddgtv,\dgtv}+<\uedgtv<\right\}_{\left(\ddgtv,\dgtv\right)\in \mathbb{N}_0\times\mathbb{N}_0}\hspace{4pt}\text{for}\hspace{4pt}\text{any}\hspace{4pt}\antels\in\mathbb{N}  \text{,}\hspace{4pt}\uedgtv\in \mathbb{N}_0\text{;}
\end{array}\\[8pt]
&\begin{array}{l}
\predcR+\argcompl{\segnvar,\segnvar}+<\segnvar<=
\left\{\predcR+\argcompl{\segnvar,\dgtv}+<\uedgtv<\right\}_{\left(\ddgtv,\dgtv,\uedgtv\right)\in \mathbb{N}_0\times\mathbb{N}_0\times\mathbb{N}_0}\text{;}
\end{array}\\[8pt]
&\begin{array}{l}
\predcR-\antels-+\tddmsk+<\uedgtv<=
\underset{\elAQ\in\AQ-\antels-+\tddmsk+<\uedgtv<
}{\bigdirsum}\predcR[\elAQ]\hspace{4pt}\text{for}\hspace{4pt}\text{any}\hspace{4pt}\antels\in\mathbb{N}  \text{,}\hspace{4pt}\tddmsk,\uedgtv\in \mathbb{N}_0\text{;}
\end{array}\\[8pt]
&\begin{array}{l}
\predcR+\tddmsk+<\uedgtv<=
\underset{\elAQ\in\AQ+\tddmsk+<\uedgtv<
}{\bigdirsum}\predcR[\elAQ]\hspace{4pt}\text{for}\hspace{4pt}\text{any}\hspace{4pt}\tddmsk,\uedgtv\in \mathbb{N}_0\text{;}
\end{array}\\[8pt]
&\begin{array}{l}
\predcR+\segnvar+<\uedgtv<=\left\{\predcR+\tddmsk+<\uedgtv<\right\}_{\tddmsk\in \mathbb{N}_0}
\hspace{4pt}\text{for}\hspace{4pt}\text{any}\hspace{4pt}\uedgtv\in \mathbb{N}_0\text{;}
\end{array}
\end{flalign*}
With an abuse of language, for any $\left\{\left(\ls, \ddgtv, \dgtv, \uedgtv,\ddgtv<seq<, \ouol<seq<\right)\right\}\in\AQ-1-$ we will write $\predcR[\ls, \ddgtv, \dgtv, \uedgtv,\ddgtv<seq<, \ouol<seq<]$ in place of $\predcR[\left\{\left(\ls, \ddgtv, \dgtv, \uedgtv,\ddgtv<seq<, \ouol<seq<\right)\right\}]$.\newline
Sum and scalar product of all $\mathbb{R}$-modules above are denoted by $\predcRsum$ and $\predcRscalp$ respectively.\newline
Elements belonging to $\mathbb{R}$-modules above are denoted by $\predcRel$.\newline
Fix $\left(\ls, \ddgtv, \dgtv, \uedgtv,\ddgtv<seq<, \ouol<seq<\right)\in \AQ-1-$, $\predcRel_{\lss}\in \generderspace-funtore-(\ddgtv>\lss>)+\uedgtv+\left(\dgtv\right)$ for any $\lss\in \left\{1,...,\ls\right\}$. To emphasize that we consider the element $\overset{\ls>\antelss>}{\underset{\lss=1}{\bigtensprodR}}\predcRel_{\lss}$ belonging to $\predcR[\ls, \ddgtv, \dgtv, \uedgtv,\ddgtv<seq<, \ouol<seq<]$ we write $\overset{\ls>\antelss>}{\underset{\lss=1}{\bigtensprodR}}\left(\predcRel_{\lss}\right)_{\ouol>\lss>}$.  
\item We define $\mathbb{R}$-bilinear arrows of $\GTopcat$
\begin{multline*}
\predcRprod<\elAQ_1<>\elAQ_2>:\predcR[\elAQ_1]\times \predcR[\elAQ_2]\rightarrow \predcR[\elAQ_1\consec \elAQ_2]\\[4pt]
\dgtv>1>,\dgtv>2>, \uedgtv\in \mathbb{N}_0\text{,}\hspace{4pt}\elAQ_1\in\AQ-1-+\argcompl{\segnvar,\dgtv>1>}+<\argcompl{\dgtv>2>+\uedgtv}<\text{,}\hspace{4pt}\elAQ_2\in\AQ-1-+\argcompl{\segnvar,\dgtv>2>}+<\uedgtv<
\end{multline*}
on generators as follows.\newline
Fix $\dgtv_1,\dgtv_2, \uedgtv\in \mathbb{N}_0$, 
\begin{flalign*}
&\begin{array}{l}
\left(\ls>1>, \ddgtv>1>, \dgtv>1>, \dgtv>2>+\uedgtv,\ddgtv<seq<>1>,\ouol<seq<>1>\right)\in\AQ-1-+\argcompl{\segnvar,\dgtv>1>}+<\argcompl{\dgtv>2>+\uedgtv}<\text{,}
\end{array}\\[6pt]
&\begin{array}{l}
\left(\ls>2>, \ddgtv>2>, \dgtv>2>, \uedgtv,\ddgtv<seq<>2>,\ouol<seq<>2>\right)\in\AQ-1-+\argcompl{\segnvar,\dgtv>2>}+<\uedgtv<\text{,}
\end{array}
\end{flalign*} 
an open neighborhood $\intunosk>1>$ of $0\in \mathbb{R}^{\ddgtv>1>+\dgtv>1>}$, an open neighborhood $\intunosk>2>$ of $0\in \mathbb{R}^{\ddgtv>2>+\dgtv>2>}$,
\begin{flalign*}
&\begin{array}{l}
\left\{\left(\pmV_{1, \lss},f_{1,\lss}\right)\right\}_{\lss\in \left\{1,...,\ls>1>\right\}}\in\\[4pt]
\hspace{20pt}
\left(\homF+\GTopcat+\left(\intunosk>\argcompl{1,\dgtv>1>}>,\genfquotcccz(\ddgtv>1,\lss>)+\argcompl{\dgtv>1>+\dgtv>2>+\uedgtv}+\right)\cap\gensieve{\ssf{\genfquot}}\right)\times\homF+\GTopcat+\left(\intunosk>\argcompl{1,\dgtv>1>}>\mathbb{R}^{\ddgtv>1,\lss>}\right)\text{,}
\end{array}\\[6pt]
&\begin{array}{l}
\left\{\left(\pmV_{2, \lss},f_{2,\lss}\right)\right\}_{\lss\in \left\{1,...,\ls>2>\right\}}\in\\[4pt]
\hspace{40pt}\left(\homF+\GTopcat+\left(\intunosk>\argcompl{2,\dgtv>2>}>,\genfquotcccz(\ddgtv>2,\lss>)+\argcompl{\dgtv>2>+\uedgtv}+\right)\cap\gensieve{\ssf{\genfquot}}\right)\times\homF+\GTopcat+\left(\intunosk>\argcompl{2,\dgtv>2>}>,\mathbb{R}^{\ddgtv>2,\lss>}\right)\text{.}
\end{array}
\end{flalign*}
Define:
\begin{flalign*}
&\begin{array}{l}
\elAQ_{1}=\left(\ls>1>, \ddgtv>1>, \dgtv>1>, \dgtv>2>+\uedgtv,\ddgtv<seq<>1>,\ouol<seq<>1>\right)\text{;}
\end{array}\\[6pt]
&\begin{array}{l}
\elAQ_{2}=\left(\ls>2>, \ddgtv>2>, \dgtv>2>, \uedgtv,\ddgtv<seq<>2>,\ouol<seq<>2>\right)\text{;}
\end{array}\\[6pt]
&\begin{array}{l}
\pmVdue_{1, \lss}=\pmV_{1, \lss}\funcomp\proj<\argcompl{\intunosk>\argcompl{1,\dgtv>1>}>,\intunosk>\argcompl{2,\dgtv>2>}>}<>1>\hspace{15pt}\forall\lss\in\left\{1,...,\ls>1>\right\}\text{;}
\end{array}\\[6pt]
&\begin{array}{l}
g_{1, \lss}=f_{1, \lss}\funcomp\proj<\argcompl{\intunosk>\argcompl{1,\dgtv>1>}>,\intunosk>\argcompl{2,\dgtv>2>}>}<>1>\hspace{15pt}
\forall\lss\in\left\{1,...,\ls>1>\right\}\text{;}
\end{array}\\[6pt]
&\begin{array}{l}
\pmVdue_{2, \lss}=\lininclquot\left(\pointedincl<\left\{\left(\mathbb{R}^{\dgtv>1>},0\right),\left(\mathbb{R}^{\dgtv>2>+\uedgtv},0\right)\right\}<>2>\right)
\compquotpoint \\[4pt]
\hspace{100pt}
\left(\pmV_{2, \lss}\funcomp \proj<\argcompl{\intunosk>\argcompl{1,\dgtv>1>}>,\intunosk>\argcompl{2,\dgtv>2>}>}<>2>\right) \hspace{15pt}\forall\lss\in\left\{1,...,\ls>2>\right\}\text{;}
\end{array}\\[6pt]
&\begin{array}{l}
g_{2, \lss}=f_{2, \lss}\funcomp \proj<\argcompl{\intunosk>\argcompl{1,\dgtv>1>}>,\intunosk>\argcompl{2,\dgtv>2>}>}<>2>\hspace{15pt}
\forall\lss\in\left\{1,...,\ls>2>\right\}\text{;}
\end{array}\\[6pt]
&\begin{array}{l} 
\predcRel_{1}=\overset{\ls>1>}{\underset{\lss=1}{\bigtensprodR}}\gdgerm[\quotdirsumdue\funcomp\left(\quotdirsumuno\funcomp\left( \genpreder\left[\pmV_{1, \lss},g_{1,\lss} \right]\right)\right)]\left(0\right)\text{;}
\end{array}\\[6pt]
&\begin{array}{l} 
\predcRel_{2}=\overset{\ls>2>}{\underset{\lss=1}{\bigtensprodR}}\gdgerm[\quotdirsumdue\funcomp\left(\quotdirsumuno\funcomp\left( \genpreder\left[\pmV_{2, \lss},g_{2,\lss} \right]\right)\right)]\left(0\right)\text{;}
\end{array}\\[6pt]
&\begin{array}{l}
\hspace{5pt}\predcRel_{3} =
\left(\overset{\ls>1>}{\underset{\lss=1}{\bigtensprodR}}\gdgerm[\quotdirsumdue\funcomp\left(\quotdirsumuno\funcomp\left( \genpreder\left[\pmVdue_{1, \lss},g_{1,\lss} \right]\right)\right)]\left(0\right)\right)
\bigtensprodR\\[4pt]
\hspace{130pt}
\left(\overset{\ls>2>}{\underset{\lss=1}{\bigtensprodR}}\gdgerm[\quotdirsumdue\funcomp\left(\quotdirsumuno\funcomp\left( \genpreder\left[\pmVdue_{2, \lss},g_{2,\lss} \right]\right)\right)]\left(0\right)\right)\text{;}
\end{array}\\[6pt]
&\begin{array}{l}
\predcRel_{1}\predcRprod<\elAQ_1<>\elAQ_2> \predcRel_{2}=\predcRel_{3}\text{.}
\end{array}
\end{flalign*} 
For any $\dgtv>1>,\dgtv>2>,\uedgtv\in \mathbb{N}_0$ the bi-graded family of arrows 
\begin{equation*}
\left\{\predcRprod<\elAQ_1<>\elAQ_2>\right\}_{\left(\elAQ_1,\elAQ_2\right)\in
\AQ-1-+\argcompl{\segnvar,\dgtv>1>}+<\argcompl{\dgtv>2>+\uedgtv}<\times\AQ-1-+\argcompl{\segnvar,\dgtv>2>}+<\uedgtv<}
\end{equation*} 
defines a $\mathbb{R}$-bilinear arrow of $\GTopcat$
\begin{equation*}
\predcRprod<\argcompl{\dgtv>1>}<(\dgtv>2>)>\uedgtv>:\predcR-1-+\argcompl{\segnvar,\dgtv>1>}+<\argcompl{\dgtv>2>+\uedgtv}<\times\predcR-1-+\argcompl{\segnvar,\dgtv>2>}+<\uedgtv<\rightarrow \predcR-1-+\argcompl{\segnvar,\dgtv>1>+\dgtv>2>}+<\uedgtv<\text{.}
\end{equation*}
The three-graded family of arrows $\left\{\predcRprod<\argcompl{\dgtv>1>}<(\dgtv>2>)>\uedgtv>\right\}_{\left(\dgtv>1>,\dgtv>2>,\uedgtv\right)\in\mathbb{N}_0\times\mathbb{N}_0\times\mathbb{N}_0}$ defines a $\mathbb{R}$-bilinear associative graded product $\predcRprod$ on the three-graded $\mathbb{R}$-module $\predcR-1-+\argcompl{\segnvar,\segnvar}+<\segnvar<$.\newline
The four-tuple $\left(\predcR-1-+\argcompl{\segnvar,\segnvar}+<\segnvar<, \predcRsum, \predcRprod, \predcRscalp\right)$ is a three-graded $\mathbb{R}$-algebra.
\item
We define $\mathbb{R}$-bilinear arrows of $\GTopcat$
\begin{equation*}
\antedcRprod<\elAQ_1<>\elAQ_2>:\predcR[\elAQ_1]\times \predcR[\elAQ_2]\rightarrow \predcR[\elAQ_1\anteconsec \elAQ_2]\hspace{15pt}
\elAQ_1,\elAQ_2\in\AQ
\end{equation*}
by setting 
\begin{equation*}
\predcRel_1 \antedcRprod<\elAQ_1<>\elAQ_2> \predcRel_2= \predcRel_1 \tensprodR \predcRel_2 \hspace{15pt}\forall \left(\predcRel_1, \predcRel_2\right)\in\predcR[\elAQ_1]\times \predcR[\elAQ_2]\text{.} 
\end{equation*}
For any $\tddmsk, \uedgtv\in \mathbb{N}_0$ the bi-graded family of arrows
\begin{equation*}
\left\{\antedcRprod<\elAQ_1<>\elAQ_2>\right\}_{\left(\elAQ_1,\elAQ_2\right)\in
\AQ+\tddmsk+<\uedgtv<\times\AQ+\tddmsk+<\uedgtv<}
\end{equation*}
defines a $\mathbb{R}$-bilinear arrow of $\GTopcat$
\begin{equation*}
\antedcRprod|\tddmsk,\uedgtv|:\predcR+\tddmsk+<\uedgtv<\times\predcR+\tddmsk+<\uedgtv<\rightarrow \predcR+\tddmsk+<\uedgtv<\text{.}
\end{equation*}
The graded family of arrows $\left\{\antedcRprod|\tddmsk,\uedgtv|\right\}_{\tddmsk\in\mathbb{N}_0}$ defines a $\mathbb{R}$-bilinear associative graded product $\antedcRprod$ on the graded $\mathbb{R}$-module $\predcR<\uedgtv<$.\newline
The four-tuple $\left(\predcR+\segnvar+<\uedgtv<, \predcRsum, \antedcRprod, \predcRscalp\right)$ is a graded $\mathbb{R}$-algebra.
\end{enumerate}
\end{definition}

\begin{remark}\label{idealemodulo}
Fix $\uedgtv\in\mathbb{N}_0$. Then $\predcR-1-+\argcompl{\segnvar, \segnvar}+<\uedgtv<$ is a proper left ideal of the $\mathbb{R}$-algebra $\left(\predcR-1-+\argcompl{\segnvar, \segnvar}+<\segnvar<, \predcRsum, \predcRprod, \predcRscalp\right)$. Whenever needed we consider $\predcR-1-+\argcompl{\segnvar, \segnvar}+<\uedgtv<$ as a left $\left(\predcR-1-+\argcompl{\segnvar, \segnvar}+<\segnvar<, \predcRsum, \predcRprod, \predcRscalp\right)$-module.
\end{remark}

In Definition \ref{raevpunfun} below we introduce a right action of germs of generalized derivations on compositions of germs of generalized derivations.

\begin{definition}\label{raevpunfun}
We define bilinear arrows of $\GTopcat$
\begin{equation*}
\begin{array}{l} 
\PRELOCAT:\predcR-\antels-+\argcompl{\ddgtv,\dgtv}+<\uedgtv< \times \generderspace-funtore-(\ddgtv>0>)+\dgtv+\left(\dgtv>0>\right)
\rightarrow \predcR-\antels-+\argcompl{\ddgtv>0>+\ddgtv,\dgtv>0>+\dgtv}+<\uedgtv<\hspace{15pt} \antels\in\mathbb{N}\text{,}\hspace{4pt} \ddgtv>0>, \ddgtv,\dgtv>0>,\dgtv,\uedgtv \in \mathbb{N}_0\text{.}
\end{array}
\end{equation*}
recursively as follows.
\begin{enumerate}
\item Fix $\left(\ls, \ddgtv, \dgtv, \uedgtv,\ddgtv<seq<,\ouol<seq<\right)\in \AQ-1-+\argcompl{\ddgtv,\dgtv}+<\uedgtv<$, $\predcRel\in\predcR[\argcompl{\ls,\ddgtv, \dgtv, \uedgtv,\ddgtv<seq<,\ouol<seq<}] $, $\predcRel_0\in \generderspace-funtore-(\ddgtv>0>)+\dgtv+\left(\dgtv>0>\right)$.\newline 
Set
\begin{flalign*}
&\begin{array}{l}
\ddgtv<seq<>0>=\left\{  \ddgtv>0,\lss>   \right\}_{\lss=1}^{\ls+1}\hspace{4pt}\text{by}\hspace{4pt}
\ddgtv>0,\lss>=\left\{
\begin{array}{ll}
 \ddgtv>0>&\text{if}\hspace{4pt}\lss=1\text{,}\\[4pt]
 \ddgtv>\lss-1>&\text{if}\hspace{4pt}\lss\in\left\{2,...,\ls+1\right\}\text{;}
\end{array}
\right.
\end{array}\\[8pt]
&\begin{array}{l}
\ouol<seq<>0>=\left\{  \ouol>0,\lss> \right\}_{\lss=1}^{\ls+1}\hspace{4pt}\text{by}\hspace{4pt}
\hspace{5pt}\ouol>0,\lss> =\left\{
\begin{array}{ll}
0&\text{if}\hspace{4pt}\lss=1\text{,}\\[4pt]
\ouol>\lss-1>&\text{if}\hspace{4pt}\lss\in\left\{2,...,\ls+1\right\}\text{.}
\end{array}
\right.
\end{array}
\end{flalign*}
Define $\PRELOCAT-\predcRel-+\argcompl{\predcRel_0}+ \in \predcR[\argcompl{\ls+1, \ddgtv>0>+\ddgtv, \dgtv>0>+\dgtv, \uedgtv,\ddgtv<seq<>0>,\ouol<seq<>0>}]$ by setting
\begin{equation*}
\PRELOCAT-\predcRel-+\predcRel_0+=\left(\left(\gendiff(\ddgtv>0>)+\lininclquot\left(\gczsfx[
\pointedincl<\left\{\left(\mathbb{R}^{\dgtv},0\right),\left(\mathbb{R}^{\uedgtv},0\right)\right\}<>1>
]\right)+\left(\dgtv>0>\right)\right)  \left(\predcRel_{0}\right)\right) \predcRprod\left(\predcRel\right)\text{.}
\end{equation*}
\item Fix $\antels>1>,\antels>2> \in \mathbb{N}$ with $\antels= \antels>1> +\antels>2>$, $\predcRel_1\in\predcR-\argcompl{\antels>1>}-+\argcompl{\ddgtv,\dgtv}+<\uedgtv<$, $\predcRel_2\in\predcR-\argcompl{\antels>2>}-+\argcompl{\ddgtv,\dgtv}+<\uedgtv<$, $\predcRel_0\in \generderspace-funtore-(\ddgtv>0>)+\dgtv+\left(\dgtv>0>\right)$. Define
\begin{equation*}
\PRELOCAT-\argcompl{\left(\predcRel_1\antedcRprod\predcRel_2\right)}-+\predcRel_0+\;=\;\left(\PRELOCAT-\predcRel_1-+\predcRel_0+\right) \; \antedcRprod \;\left(\PRELOCAT-\predcRel_2-+\predcRel_0+ \right)\text{.}
\end{equation*}
\end{enumerate}
\end{definition}

In Definition \ref{laevpunfun} below we introduce a left action of germs of generalized derivations on compositions of germs of generalized derivations.

\begin{definition}\label{laevpunfun}
We define bilinear arrows of $\GTopcat$
\begin{equation*}
\begin{array}{l} 
\PREDER:\generderspace-funtore-(\ddgtv>0>)+\dgtv+\left(\dgtv>0>\right)   \times  \predcR-\antels-+\argcompl{\ddgtv,\dgtv}+<\uedgtv<
\rightarrow \predcR-\antels-+\argcompl{\ddgtv>0>+\ddgtv,\dgtv>0>+\dgtv}+<\uedgtv<\hspace{15pt} \antels\in\mathbb{N}\text{,}\hspace{4pt} \ddgtv>0>, \ddgtv,\dgtv>0>,\dgtv,\uedgtv \in \mathbb{N}_0\text{.}
\end{array}
\end{equation*}
recursively as follows.
\begin{enumerate}
\item Fix $\left(\ls, \ddgtv, \dgtv, \uedgtv,\ddgtv<seq<,\ouol<seq<\right)\in \AQ-1-+\argcompl{\ddgtv,\dgtv}+<\uedgtv<$, $\predcRel\in\predcR[\argcompl{\ls,\ddgtv, \dgtv, \uedgtv,\ddgtv<seq<,\ouol<seq<}] $, $\predcRel_0\in \generderspace-funtore-(\ddgtv>0>)+\dgtv+\left(\dgtv>0>\right)$.\newline 
Set
\begin{flalign*}
&\begin{array}{l}
\ddgtv<seq<>0>=\left\{  \ddgtv>0,\lss>   \right\}_{\lss=1}^{\ls+1}\hspace{4pt}\text{by}\hspace{4pt}
\ddgtv>0,\lss>=\left\{
\begin{array}{ll}
 \ddgtv>0>&\text{if}\hspace{4pt}\lss=1\text{,}\\[4pt]
 \ddgtv>\lss-1>&\text{if}\hspace{4pt}\lss\in\left\{2,...,\ls+1\right\}\text{;}
\end{array}
\right.
\end{array}\\[8pt]
&\begin{array}{l}
\ouol<seq<>0>=\left\{  \ouol>0,\lss> \right\}_{\lss=1}^{\ls+1}\hspace{4pt}\text{by}\hspace{4pt}
\hspace{5pt}\ouol>0,\lss> =\left\{
\begin{array}{ll}
1&\text{if}\hspace{4pt}\lss=1\text{,}\\[4pt]
\ouol>\lss-1>&\text{if}\hspace{4pt}\lss\in\left\{2,...,\ls+1\right\}\text{.}
\end{array}
\right.
\end{array}
\end{flalign*}
Define $\PREDER-\predcRel-+\argcompl{\predcRel_0}+ \in \predcR[\argcompl{\ls+1, \ddgtv>0>+\ddgtv, \dgtv>0>+\dgtv, \uedgtv,\ddgtv<seq<>0>,\ouol<seq<>0>}]$ by setting
\begin{equation*}
\PREDER-\predcRel-+\predcRel_0+=\left(\left(\gendiff(\ddgtv>0>)+\lininclquot\left(\gczsfx[
\pointedincl<\left\{\left(\mathbb{R}^{\dgtv},0\right),\left(\mathbb{R}^{\uedgtv},0\right)\right\}<>1>
]\right)+\left(\dgtv>0>\right)\right)  \left(\predcRel_{0}\right)\right) \predcRprod\left(\predcRel\right)\text{.}
\end{equation*}
\item Fix $\antels>1>,\antels>2> \in \mathbb{N}$ with $\antels= \antels>1> +\antels>2>$, $\predcRel_1\in\predcR-\argcompl{\antels>1>}-+\argcompl{\ddgtv,\dgtv}+<\uedgtv<$, $\predcRel_2\in\predcR-\argcompl{\antels>2>}-+\argcompl{\ddgtv,\dgtv}+<\uedgtv<$, $\predcRel_0\in \generderspace-funtore-(\ddgtv>0>)+\dgtv+\left(\dgtv>0>\right)$. Define
\begin{multline*}
\PREDER-\argcompl{\predcRel_1\antedcRprod\predcRel_2}-+\predcRel_0+\;=\\[4pt]
\left(\PRELOCAT-\predcRel_1-+\predcRel_0+\right) \; \antedcRprod \;\left(\PREDER-\predcRel_2-+\predcRel_0+ \right)\;\predcRsum  \;\left(\PREDER-\predcRel_1-+\predcRel_0+\right) \; \antedcRprod \;\left(\PRELOCAT-\predcRel_2-+\predcRel_0+ \right)\text{.}
\end{multline*}
\end{enumerate}
\end{definition}

In Definitions \ref{corrrapp} below we introduce the notion of skeleton for higher order generalized differential operators, this extends to $\predcR$ the notion of skeleton given in Remark \ref{coeqdef}-[10].\newline
We refer to Notations \ref{varispder}, \ref{ins}-[10], \ref{realvec}-[4], Remark \ref{pathinlgenerder}-[8, 9].

\begin{definition}\label{corrrapp}  We define the notion of skeleton of elements belonging to $\predcR$ recursively as follows.
\begin{enumerate}
\item Fix an element $\left(\ls, \ddgtv, \dgtv, \uedgtv,\ddgtv<seq<,\ouol<seq<\right)\in\AQ-1-$, an open neighborhood $\domsk$ of $0\in\mathbb{R}^{\ddgtv+\dgtv}$, $\left(\pmV_{\lss},f_{\lss}\right)\in\left(\homF+\GTopcat+\left( \domsk>\dgtv>   ,\genfquotcccz(\ddgtv>\lss>)+\argcompl{\dgtv+\uedgtv}+\right)\cap\gensieve{\ssf{\genfquot}}\right)\times\homF+\GTopcat+\left( \domsk>\dgtv>   ,\mathbb{R}^{\ddgtv>\lss>}\right)$, $\pmgenerder_{\lss}\in \generderspace-functor-(\ddgtv>\lss>)+\argcompl{\dgtv+\uedgtv}+\left(\domsk>\dgtv>\right)$, $\predcRel_{\lss}\in \generderspace-functor-(\ddgtv>\lss>)+\argcompl{\dgtv+\uedgtv}+\left(\dgtv\right)$ for any $\lss\in \left\{1,...,\ls\right\}$.\newline
Assume that all conditions below are fulfilled:
\begin{flalign}
&\begin{array}{l}
\ihtl{\domsk>\dgtv>}{\ddgtv>\lss>}{\dgtv+\uedgtv}\left(\pmV_{\lss}\right)  \neq \udenset\hspace{15pt}\forall\lss\in \left\{1,...,\ls\right\}\text{;}
\end{array}\label{nonvuotoihtl}\\[8pt]
&\begin{array}{l} 
\pmgenerder_{\lss}=\quotdirsumdue\funcomp\left(\quotdirsumuno\funcomp\genpreder\left[\pmV_{\lss},f_{\lss} \right]\right)\hspace{15pt}\forall\lss\in \left\{1,...,\ls\right\}\text{;}
\end{array}\\[8pt]
&\begin{array}{l}
\predcRel_{\lss}= \gdgerm[\pmgenerder_{\lss}]\left(0\right)\hspace{15pt}\forall\lss\in \left\{1,...,\ls\right\}\text{.}
\end{array}
\end{flalign}
Set:
\begin{flalign*}
&\begin{array}{l}
\left( \pmV,f,\ouol<seq<\right)=
\left(\left( \pmV_{\lss},f_{\lss}\right),\ouol<seq<>\lss>\right)_{\lss=1}^{\ls}\text{;}
\end{array}\\[6pt]
&\begin{array}{l}
\predcRel=\overset{\ls}{\underset{\lss=1}{\bigtensprodR}} \, \predcRel_{\lss}\text{.}
\end{array}
\end{flalign*}
We say that:
\begin{flalign}
&\begin{array}{l}
\left(\predcRel_{\lss}\right)_{\lss=1}^{\ls} \hspace{4pt}\text{is}\hspace{4pt}\text{a}\hspace{4pt}\text{decomposition}\hspace{4pt}\text{vector}\hspace{4pt}\text{of}\hspace{4pt}\predcRel\text{;}
\end{array}\nonumber\\[6pt]
&\begin{array}{l}
\domsk\hspace{4pt}\text{is}\hspace{4pt}\text{a}\hspace{4pt}\text{domain}\hspace{4pt}\text{of}\hspace{4pt}\predcRel\text{;}\label{skfdomdue}
\end{array}\\[6pt]
&\begin{array}{l}
\left( \pmV,f,\ouol<seq<\right)\hspace{4pt}\text{is}\hspace{4pt}\text{a}\hspace{4pt}\text{skeleton}\hspace{4pt}\text{of}\hspace{4pt}\predcRel\text{;}
\label{skfdomuno}
\end{array}\\[6pt]
&\begin{array}{l}
\left(\pmgenerder_{\lss}\right)_{\lss=1}^{\ls} \hspace{4pt}\text{is}\hspace{4pt}\text{a}\hspace{4pt}\text{representative}\hspace{4pt}\text{vector}\hspace{4pt}\text{of}\hspace{4pt}\predcRel\hspace{4pt}\text{related}\hspace{4pt}\text{to}\hspace{4pt}\left( \pmV,f,\ouol<seq<\right)\text{;}\nonumber
\end{array}\\[6pt]
&\begin{array}{l}
\overset{\ls}{\underset{\lss=1}{\bigtensprodR}} \,    \pmgenerder_{\lss}\hspace{4pt}\text{is}\hspace{4pt}\text{a}\hspace{4pt}\text{representative}\hspace{4pt}\text{of}\hspace{4pt}\predcRel\hspace{4pt}\text{related}\hspace{4pt}\text{to}\hspace{4pt}\left( \pmV,f,\ouol<seq<\right)\text{.}
\end{array}\nonumber
\end{flalign}
With an abuse of language we summarize \eqref{skfdomdue}, \eqref{skfdomuno} by saying that $\left( \pmV,f,\ouol<seq<\right)$ is a skeleton of $\predcRel$ with domain $\domsk$.
\item Fix a finite set $\setindexuno$, $\tddmsk\in \mathbb{N}_0$, an open neighborhood $\domsk$ of $0\in \mathbb{R}^{\tddmsk}$, $\left( \ls>\elindexuno>, \ddgtv>\elindexuno>, \dgtv>\elindexuno>, \uedgtv>\elindexuno>,\ddgtv<seq<>\elindexuno>,\ouol<seq<>\elindexuno>\right)\in\AQ-1-+\tddmsk+$,    
$\predcRel_{\elindexuno}\in\predcR[\ls>\elindexuno>, \ddgtv>\elindexuno>, \dgtv>\elindexuno>, \uedgtv>\elindexuno>,\ddgtv<seq<>\elindexuno>,\ouol<seq<>\elindexuno>]$, 
a skeleton $\left( \pmV_{\elindexuno},f_{\elindexuno},\ouol<seq<>\elindexuno>\right)$ of $\predcRel_{\elindexuno}$ with domain $\domsk$, a representative vector $\left(\pmgenerder_{\elindexuno,\lss}\right)_{\lss=1}^{\ls>\elindexuno>}$ of $\predcRel_{\elindexuno}$ related to $\left( \pmV_{\elindexuno},f_{\elindexuno},\ouol<seq<>\elindexuno>\right)$, a decomposition $\left(\predcRel_{\elindexuno,\lss}\right)_{\lss=1}^{\ls>\elindexuno>}$ of $\predcRel_{\elindexuno}$ related to $\left( \pmV_{\elindexuno},f_{\elindexuno},\ouol<seq<>\elindexuno>\right)$ for any $\elindexuno\in \setindexuno$.\newline
Set $\left( \pmV,f,\ouol<seq<\right)=
\left(\left( \pmV_{\elindexuno},f_{\elindexuno},\ouol<seq<>\elindexuno>\right)\right)_{\elindexuno \in\setindexuno}$.\newline
We say that:
\begin{flalign}
&\begin{array}{l}
\left(\left(\predcRel_{\elindexuno,\lss}\right)_{\lss=1}^{\ls>\elindexuno>}\right)_{\elindexuno \in\setindexuno} \hspace{4pt}\text{is}\hspace{4pt}\text{a}\hspace{4pt}\text{decomposition}\hspace{4pt}\text{vector}\hspace{4pt}\text{of}\hspace{4pt}\left(\predcRel_{\elindexuno}\right)_{\elindexuno\in\setindexuno}\text{;}
\end{array}\nonumber\\[6pt]
&\begin{array}{l}
\domsk\hspace{4pt}\text{is}\hspace{4pt}\text{a}\hspace{4pt}\text{domain}\hspace{4pt}\text{of}\hspace{4pt}\left(\predcRel_{\elindexuno}\right)_{\elindexuno\in\setindexuno}\text{;}\label{Iskfdomdue}
\end{array}\\[6pt]
&\begin{array}{l}
\left( \pmV,f,\ouol<seq<\right)\hspace{4pt}\text{is}\hspace{4pt}\text{a}\hspace{4pt}\text{skeleton}\hspace{4pt}\text{of}\hspace{4pt} \left(\predcRel_{\elindexuno}\right)_{\elindexuno\in\setindexuno}\text{;}
\end{array}\label{Iskfdomuno}\\[6pt]
&\left\{
\begin{array}{l}
\left(\left(\pmgenerder_{\elindexuno,\lss}\right)_{\lss=1}^{\ls>\elindexuno>}\right)_{\elindexuno \in\setindexuno} \hspace{4pt}\text{is}\hspace{4pt}\text{a}\hspace{4pt}\text{representative}\hspace{4pt}\text{vector}\hspace{4pt}\text{of}\hspace{4pt}\left(\predcRel_{\elindexuno}\right)_{\elindexuno\in\setindexuno}\hspace{4pt}\text{related}\\[4pt]\text{to}\hspace{4pt}\left( \pmV,f,\ouol<seq<\right)\text{.}
\end{array}
\right.\nonumber
\end{flalign}
With an abuse of language we summarize \eqref{Iskfdomdue}, \eqref{Iskfdomuno} by saying that $\left( \pmV,f,\ouol<seq<\right)$ is a skeleton of a $\left(\predcRel_{\elindexuno}\right)_{\elindexuno\in\setindexuno}$ with domain $\domsk$.
\item Fix a finite set $\setindexuno$, $\tddmsk\in \mathbb{N}_0$, an open neighborhood $\domsk$ of $0\in \mathbb{R}^{\tddmsk}$, $\vecdue_{\elindexuno}\in \mathbb{R}$,
\begin{flalign*}
&\begin{array}{l}
\left\{\left(\ls>\elindexuno,\antelss>, \ddgtv>\elindexuno,\antelss>, \dgtv>\elindexuno,\antelss>, \uedgtv>\elindexuno,\antelss>,\ddgtv<seq<>\elindexuno,\antelss>, \ouol<seq<>\elindexuno,\antelss>\right)\right\}_{\antelss=1}^{\antels>\elindexuno>}\in\AQ+\tddmsk+\text{,}
\end{array}\\[8pt]
&\begin{array}{l}
\predcRel_{\elindexuno,\antelss}\in\predcR[\ls>\elindexuno,\antelss>, \ddgtv>\elindexuno,\antelss>, \dgtv>\elindexuno,\antelss>, \uedgtv>\elindexuno,\antelss>,\ddgtv<seq<>\elindexuno,\antelss>,\ouol<seq<>\elindexuno,\antelss>]
\end{array}
\end{flalign*}
for any $\elindexuno \in  \setindexuno$, $\antelss\in\left\{1,...,\antels>\elindexuno>\right\}$, a skeleton $\left( \pmV,f,\ouol<seq<\right)$ of $\left(\left(\predcRel_{\elindexuno,\antelss}\right)_{\antelss=1}^{\antels>\elindexuno>}\right)_{\elindexuno \in  \setindexuno}$.\newline
Set $\predcRel= \underset{\elindexuno \in  \setindexuno}{\bigdirsum}\, \vecdue_{\elindexuno} \predcRscalp\left(\overset{\antels>\elindexuno>}{\underset{\antelss=1}{\bigtensprodR}} \predcRel_{\elindexuno,\antelss}\right)$.\newline
We say that:
\begin{flalign*}
&\begin{array}{l}
\left(\left(\predcRel_{\elindexuno,\antelss}\right)_{\antelss=1}^{\antels>\elindexuno>}\right)_{\elindexuno \in  \setindexuno} \hspace{4pt}\text{is}\hspace{4pt}\text{a}\hspace{4pt}\text{decomposition}\hspace{4pt}\text{vector}\hspace{4pt}\text{of}\hspace{4pt}\predcRel\text{;}
\end{array}\nonumber\\[6pt]
&\begin{array}{l}
\left( \pmV,f,\ouol<seq<\right) \hspace{4pt}\text{is}\hspace{4pt}\text{a}\hspace{4pt}\text{skeleton}\hspace{4pt}\text{of}\hspace{4pt}\predcRel\hspace{4pt}\text{with}\hspace{4pt}\text{domain}\hspace{4pt}\domsk\text{.}
\end{array}
\end{flalign*}
\end{enumerate}
\end{definition}

\begin{remark}\label{spiegafinlens}\mbox{}
\begin{enumerate}
\item Condition \eqref{nonvuotoihtl} does not affect generality because of Proposition \ref{preimmerfun}-[2, 3, 4, 5, 7, 8] and Remark \ref{pathinlgenerder}-[2] together with \eqref{pugualeval}.
\item Existence of skeletons of finite sets of germs follows by conditions \eqref{bGt7bis}, \eqref{pcM} in definition of Grothendieck topology of the site $\GTopcat$.
\end{enumerate}
\end{remark}

In Definition \ref{corrrapp2} below we introduce the p-corresponding skeleton of elements belonging to $\predcR$.

\begin{definition}\label{corrrapp2} We define the notion of p-corresponding skeleton of elements belonging to $\predcR$, recursively as follows:
\begin{enumerate}
\item Fix an element $\left(\ls, \ddgtv, \dgtv, \uedgtv,\ddgtv<seq<,\ouol<seq<\right)\in\AQ-1-$, a germ $\predcRel\in\predcR[\argcompl{\ls, \ddgtv, \dgtv, \uedgtv,\ddgtv<seq<,\ouol<seq<}]$,
a skeleton $\left( \pmV,f,\ouol<seq<\right)$ of $\predcRel$ with domain $\domsk$, an arrow $\corrarr[\pmV]_{\lss,\elindexdue}$ p-corresponding to $\pmV_{\lss,\elindexdue}$, $\gfx_{\lss,\elindexdue}\in \genf$ with $\evalcomptwo\left(\gfx_{\lss,\elindexdue}\right)=f_{\lss,\elindexdue}$ for any $\lss\in \left\{1,...,\ls\right\}$, $\elindexdue \in \setindexdue_{\lss}$.\newline
Set $\left( \corrarr[\pmV],\gfx,\ouol<seq<\right)=\left(\left( \corrarr[\pmV]_{\lss},\gfx_{\lss}\right),\ouol<seq<>\lss>\right)_{\lss=1}^{\ls}$.\newline
We say that $\left( \corrarr[\pmV],\gfx,\ouol<seq<\right)$ is a p-corresponding skeleton of $\predcRel$ related to $\left( \pmV,\gfx,\ouol<seq<\right)$.
\item Fix a finite set $\setindexuno$, $\tddmsk\in \mathbb{N}_0$, an open neighborhood $\domsk\sqsubseteq \mathbb{R}^{\tddmsk}$ of $0\in \mathbb{R}^{\tddmsk}$,  $\left(
\ls>\elindexuno>, \ddgtv>\elindexuno>, \dgtv>\elindexuno>, \uedgtv>\elindexuno>,\ddgtv<seq<>\elindexuno>,\ouol<seq<>\elindexuno>\right)\in\AQ-1-+\tddmsk+$, $\predcRel_{\elindexuno}\in\predcR[\ls>\elindexuno>, \ddgtv>\elindexuno>, \dgtv>\elindexuno>, \uedgtv>\elindexuno>,\ddgtv<seq<>\elindexuno>,\ouol<seq<>\elindexuno>]$, a skeleton $\left( \pmV,f, \ouol<seq<\right)$ of $\left(\predcRel_{\elindexuno}\right)_{\elindexuno\in\setindexuno}$ with domain $\domsk$, a p-corresponding skeleton $\left( \corrarr[\pmV]_{\elindexuno},\gfx_{\elindexuno},\ouol<seq<>\elindexuno>\right)$ of $\predcRel_{\elindexuno}$ related to $\left( \pmV_{\elindexuno},f_{\elindexuno},\ouol<seq<>\elindexuno>\right)$ for any $\elindexuno\in \setindexuno$.\newline
Set $\left( \corrarr[\pmV],\gfx,\ouol<seq<\right)=
\left(\left( \corrarr[\pmV]_{\elindexuno},\gfx_{\elindexuno},\ouol<seq<\right)\right)_{\elindexuno\in \setindexuno}$.\newline
We say that $\left( \corrarr[\pmV],\gfx,\ouol<seq<\right)$ is a p-corresponding skeleton of $\left(\predcRel_{\elindexuno}\right)_{\elindexuno\in\setindexuno}$ related to $\left( \pmV,f,\ouol<seq<\right)$.
\item Fix a finite set $\setindexuno$, $\tddmsk\in \mathbb{N}_0$, $\vecdue_{\elindexuno}\in \mathbb{R}$,
$\left\{\left(\ls>\elindexuno,\antelss>, \ddgtv>\elindexuno,\antelss>, \dgtv>\elindexuno,\antelss>, \uedgtv>\elindexuno,\antelss>,\ddgtv<seq<>\elindexuno,\antelss>, \ouol<seq<>\elindexuno,\antelss>\right)\right\}_{\antelss=1}^{\antels>\elindexuno>}\in\AQ+\tddmsk+$,
$\predcRel_{\elindexuno,\antelss}\in\predcR[\ls>\elindexuno,\antelss>, \ddgtv>\elindexuno,\antelss>, \dgtv>\elindexuno,\antelss>, \uedgtv>\elindexuno,\antelss>,\ddgtv<seq<>\elindexuno,\antelss>,\ouol<seq<>\elindexuno,\antelss>]$ for any $\elindexuno \in  \setindexuno$, $\antelss\in\left\{1,...,\antels>\elindexuno>\right\}$, a skeleton $\left( \pmV,f,\ouol<seq<\right)$ of $\left(\left(\predcRel_{\elindexuno,\antelss}\right)_{\antelss=1}^{\antels>\elindexuno>}\right)_{\elindexuno \in  \setindexuno}$, a p-corresponding skeleton $\left( \corrarr[\pmV],\gfx,\ouol<seq<\right)$ of $\left(\left(\predcRel_{\elindexuno,\antelss}\right)_{\antelss=1}^{\antels>\elindexuno>}\right)_{\elindexuno \in  \setindexuno}$ related to $\left( \pmV,f,\ouol<seq<\right)$.\newline
Set $\predcRel= \underset{\elindexuno \in  \setindexuno}{\bigdirsum}\, \vecdue_{\elindexuno} \predcRscalp\left(\overset{\antels>\elindexuno>}{\underset{\antelss=1}{\bigtensprodR}} \predcRel_{\elindexuno,\antelss}\right)$.\newline
We say that $\left( \corrarr[\pmV],\gfx,\ouol<seq<\right)$ is a p-corresponding skeleton of $\predcRel$ related to $\left( \pmV,f,\ouol<seq<\right)$.
\end{enumerate}
\end{definition}

In Definition \ref{geneval} below we introduce evaluation of elements belonging to $\predcR $. We refer to Notation \ref{ins}-[10, 11], Definitions \ref{genpointdef}, \ref{genderdef}, Proposition \ref{genpointprop}, Remark \ref{siterem}-[3].

\begin{definition}\label{geneval}\mbox{}
We define the notion of evaluation of elements belonging to $\predcR$ recursively as follows:
\begin{enumerate}
\item Fix $\left(\ls, \ddgtv, \dgtv, \uedgtv,\ddgtv<seq<,\ouol<seq<\right)\in\AQ-1-$, a germ $\predcRel\in\predcR[\argcompl{\ls, \ddgtv, \dgtv, \uedgtv,\ddgtv<seq<,\ouol<seq<}]$,
a skeleton $\left( \pmV,f,\ouol<seq<\right)$ of $\predcRel$ with domain $\domsk$, a p-corresponding skeleton $\left( \corrarr[\pmV],\gfx,\ouol<seq<\right)$ of $\predcRel$ related to $\left( \pmV,f,\ouol<seq<\right)$.\newline
We define the arrow of $\GTopcat$
\begin{equation*}
\evpdf[\left( \corrarr[\pmV],\gfx, \ouol<seq<\right)] :\domsk\times\genfquot(\argcompl{\uedgtv})+1+ \rightarrow \genfquot(\argcompl{\ddgtv+\dgtv+\uedgtv})+1+
\end{equation*}
 recursively as follows.\newline
We define the arrow 
\begin{equation*}
\evpdf<\ls<[\left( \corrarr[\pmV],\gfx, \ouol<seq<\right)] :
\domsk\times\genfquot(\argcompl{\uedgtv})+1+ \rightarrow \genfquot(\argcompl{\ddgtv+\dgtv+\uedgtv})+1+
\end{equation*}
of $\GTopcat$ by setting
\begin{equation*}
\begin{array}{l}
\evpdf<\ls<[\left( \corrarr[\pmV],\gfx, \ouol<seq<\right)] \left(\unkuno,\gggfw\right)=
\\[10pt] 
\left\{ 
\begin{array}{ll}
\left(\gggfw\compquotpoint \left(\corrarr[\pmV]_{1,\ls}\left(\unkuno\right)\right) \right)\compquotpoint \lininclquot\left(
\gczsfx[ \proj<\left\{\mathbb{R}^{\ddgtv>\ls>+\dgtv+\uedgtv},\mathbb{R}^{\dgtv}\right\}<>1>  \funcomp h_{\ls}]\right)&  \text{if}\hspace{4pt}\ouol>\ls>= 0\text{,}\\[10pt]
\left(\underset{\elindexddgtv=1}{\overset{\ddgtv>\ls>}{\bigsumquotpoint}}\,
\left( \Fpart_{\elindexddgtv} \bsfmuquotpoint \left(\gggfw\compquotpoint \left(\corrarr[\pmV]_{1,\ls}\left(\unkuno\right)\right)\right)\right) \sqmultquotpoint  \gczsfx[\coord{\ddgtv>\ls>}{\elindexddgtv} \genfuncomp \gfx_{\ls}] \right)\compquotpoint
 \lininclquot\left(
\gczsfx[h_{\ls}]\right)& \text{if}\hspace{4pt}\ddgtv>\ls>\neq 0\hspace{4pt}\text{and}\hspace{4pt}\ouol>\ls>= 1\text{,}\\[10pt]
\left(\left(  \gggfw\compquotpoint \left(\corrarr[\pmV]_{1,\ls}\left(\unkuno\right)\right)\right) 
\sqmultquotpoint  \gczsfx[\cost<\argcompl{\mathbb{R}^{\ddgtv>\ls>}}<>\mathbb{R}>+0+ \genfuncomp \gfx_{\ls}] \right)\compquotpoint \lininclquot\left(
\gczsfx[h_{\ls}]\right)&  \text{if}\hspace{4pt}\ddgtv>\ls>= 0\hspace{4pt}\text{and}\hspace{4pt}\ouol>\ls>= 1
\end{array}
\right.\\[35pt]
\forall\left(\unkuno,\gggfw\right)\in  \domsk \times\genfquot(\uedgtv)+1+
\end{array}
\end{equation*}
where 
\begin{flalign*}
&\begin{array}{l}
\corrarr[\pmV]_{1,\ls}= \left(\lininclquot\left(\gczsfx[
\proj<\left\{\mathbb{R}^{\dgtv},\mathbb{R}^{\uedgtv}\right\}<>2>
]\right)\compquotpoint
\corrarr[\pmV]_{\ls}\right)\sqsumquotpoint   \lininclquot\left(
\gczsfx[\idobj+\mathbb{R}^{\uedgtv}+]\right)\text{;}
\end{array}\\[8pt]
&\left\{
\begin{array}{l}
h_{\ls}:\mathbb{R}^{\ddgtv}\times\mathbb{R}^{\dgtv}\times\mathbb{R}^{\uedgtv} \rightarrow  \mathbb{R}^{\ddgtv>\ls>}\times\mathbb{R}^{\dgtv}\times\mathbb{R}^{\uedgtv}\times\mathbb{R}^{\dgtv}\hspace{4pt}
\text{is}\hspace{4pt}\text{given}\hspace{4pt}\text{by}\\[4pt]
h_{\ls}\left(\unkuno_1,\unkuno_2,\unkuno_3\right)=\left(\proj<\left\{\mathbb{R}^{\ddgtv>1>},...,\mathbb{R}^{\ddgtv>\ls>}\right\}<>\ls> \left(\unkuno_1\right)     ,\unkuno_2,\unkuno_3,\unkuno_2\right)\\[4pt]
\hspace{150pt}\forall\left(\unkuno_1,\unkuno_2,\unkuno_3\right)\in\mathbb{R}^{\ddgtv}\times\mathbb{R}^{\dgtv}\times\mathbb{R}^{\uedgtv} \text{.}
\end{array}
\right.
\end{flalign*}
For any $\lss\in \left\{1,...,\ls-1\right\}$ we define the arrow
\begin{equation*}
\evpdf<\lss<[\left( \corrarr[\pmV],\gfx, \ouol<seq<\right)] :
\domsk\times\genfquot(\argcompl{\uedgtv})+1+ \rightarrow \genfquot(\argcompl{\ddgtv+\dgtv+\uedgtv})+1+
\end{equation*}
of $\GTopcat$ by setting
\begin{equation*}
\begin{array}{l}
\evpdf<\lss<[\left( \corrarr[\pmV],\gfx, \ouol<seq<\right)] \left(\unkuno,\gggfw\right)=\\[10pt] 
\left\{ 
\begin{array}{ll}
\left(  \gggfu\left(\unkuno\right)\compquotpoint \left(\corrarr[\pmV]_{1,\lss}\left(\unkuno\right)\right)\right)\compquotpoint \lininclquot\left(\gczsfx[
\proj<\left\{\mathbb{R}^{\ddgtv>\ls>+\dgtv+\ddgtv+\dgtv+\uedgtv},\mathbb{R}^{\dgtv}\right\}<>1>  \funcomp
h_{\lss}]\right)&  \text{if}\hspace{4pt}\ouol>\ls>= 0\\[10pt]
\left(\underset{\elindexddgtv=1}{\overset{\ddgtv>\lss>}{\bigsumquotpoint}}\,
\left( \Fpart_{\elindexddgtv} \bsfmuquotpoint \left(\gggfu\left(\unkuno\right)\compquotpoint \left(\corrarr[\pmV]_{1,\lss}\left(\unkuno\right)\right)\right)\right) \sqmultquotpoint  \gczsfx[\coord{\ddgtv>\lss>}{\elindexddgtv} \genfuncomp \gfx_{\lss}] \right)\compquotpoint \lininclquot\left(\gczsfx[h_{\lss}]\right)&  \text{if}\hspace{4pt}\ddgtv>\lss>\neq 0\hspace{4pt}\text{and}\hspace{4pt}\ouol>\ls>= 1\text{,}\\[10pt]
\left(\left(  \gggfu\left(\unkuno\right)\compquotpoint \left(\corrarr[\pmV]_{1,\lss,\elindexdue}\left(\unkuno\right)\right)\right) \sqmultquotpoint  \gczsfx[\cost<\argcompl{\mathbb{R}^{\ddgtv>\lss>}}<>\mathbb{R}>+0+ \genfuncomp \gfx_{\lss}]
\right)\compquotpoint \lininclquot\left(\gczsfx[h_{\lss}]\right)&  \text{if}\hspace{4pt}\ddgtv>\lss>= 0\hspace{4pt}\text{and}\hspace{4pt}\ouol>\ls>= 1
\end{array}
\right.\\[25pt]
\forall\left(\unkuno,\gggfw\right)\in  \domsk \times\genfquot(\uedgtv)+1+ 
\end{array}
\end{equation*}
where 
\begin{flalign*}
&\begin{array}{l}
\gggfu\left(\unkuno\right)=\evpdf<\lss+1<[\left( \corrarr[\pmV],\gfx, \ouol<seq<\right)] \left(\unkuno,\gggfw\right)\text{;}
\end{array}\nonumber\\[8pt]
&\begin{array}{l}
\corrarr[\pmV]_{1,\lss}= \left(\lininclquot\left(\gczsfx[\argcompl{\pointedincl<\left\{\left(\mathbb{R}^{  \ddgtv},0\right), \left(\mathbb{R}^{\dgtv+\uedgtv},0\right)\right\}<>2>}]\right)\compquotpoint\corrarr[\pmV]_{\lss}\right)\sqsumquotpoint  \lininclquot\left(\gczsfx[\idobj+\mathbb{R}^{\ddgtv+\dgtv+\uedgtv}+]\right)\text{;}
\end{array}\\[8pt]
&\left\{
\begin{array}{l}
h_{\lss}:\mathbb{R}^{\ddgtv}\times\mathbb{R}^{\dgtv}\times\mathbb{R}^{\uedgtv} \rightarrow\mathbb{R}^{\ddgtv>\lss>}\times\mathbb{R}^{\dgtv}\times\mathbb{R}^{\ddgtv}\times\mathbb{R}^{\dgtv}\times\mathbb{R}^{\uedgtv}\times\mathbb{R}^{\dgtv}\hspace{4pt}
\text{is}\hspace{4pt}\text{given}\hspace{4pt}\text{by}\\[4pt]
h_{\ls}\left(\unkuno_1,\unkuno_2,\unkuno_3\right)=\left(\proj<\left\{\mathbb{R}^{\ddgtv>1>},...,\mathbb{R}^{\ddgtv>\ls>}\right\}<>\lss> \left(\unkuno_1\right)     ,\unkuno_2,\unkuno_1,\unkuno_2,\unkuno_3,\unkuno_2\right)\\[4pt]
\hspace{150pt}\forall\left(\unkuno_1,\unkuno_2,\unkuno_3\right)\in\mathbb{R}^{\ddgtv}\times\mathbb{R}^{\dgtv}\times\mathbb{R}^{\uedgtv} \text{.}
\end{array}
\right.
\end{flalign*}
Eventually we set $\evpdf[\left( \corrarr[\pmV],\gfx, \ouol<seq<\right)]=\evpdf<1<[\left( \corrarr[\pmV],\gfx, \ouol<seq<\right)]$.
\item Fix a finite set $\setindexuno$, $\tddmsk\in \mathbb{N}_0$, an open neighborhood $\domsk\sqsubseteq \mathbb{R}^{\tddmsk}$ of $0\in \mathbb{R}^{\ddgtv+\dgtv}$, $\left(\ls>\elindexuno>, \ddgtv>\elindexuno>, \dgtv>\elindexuno>, \uedgtv>\elindexuno>,\ddgtv<seq<>\elindexuno>,\ouol<seq<>\elindexuno>\right)\in\AQ-1-+\tddmsk+$, $\predcRel_{\elindexuno}\in\predcR[\ls>\elindexuno>, \ddgtv>\elindexuno>, \dgtv>\elindexuno>, \uedgtv>\elindexuno>,\ddgtv<seq<>\elindexuno>,\ouol<seq<>\elindexuno>]$ for any $\elindexuno\in \setindexuno$, a skeleton $\left( \pmV,f, \ouol<seq<\right)$ of $\left(\predcRel_{\elindexuno}\right)_{\elindexuno\in\setindexuno}$ with domain $\domsk$, a p-corresponding skeleton $\left( \corrarr[\pmV]_{\elindexuno},\gfx_{\elindexuno},\ouol<seq<\right)$ of $\left(\predcRel_{\elindexuno}\right)_{\elindexuno\in\setindexuno}$ related to $\left( \pmV,f,\ouol<seq<\right)$.\newline
We define the arrow   
\begin{equation*}
\evpdf[\argcompl{\left( \corrarr[\pmV],\gfx,\ouol<seq<\right)}] :\domsk\times\left(\underset{\elindexuno\in\setindexuno}{\prod}\,\genfquot(\uedgtv_{\elindexuno})+1+\right) \rightarrow 
\underset{\elindexuno\in\setindexuno}{\prod}\,\genfquot(\tddmsk +\uedgtv_{\elindexuno})+1+
\end{equation*}
of $\GTopcat$ by setting  
\begin{multline*}
\evpdf[\argcompl{\left( \corrarr[\pmV],\gfx,\ouol<seq<\right)}]\left(\unkuno,\left(\gggfw_{\elindexuno}\right)_{\elindexuno\in\setindexuno}\right)=
\left(
\evpdf[\argcompl{\left( \corrarr[\pmV]_{\elindexuno},\gfx_{\elindexuno},\ouol>\elindexuno>\right)}]
\left(\unkuno,\gggfw_{\elindexuno}\right)\right)_{\elindexuno\in\setindexuno}\\[4pt]
\forall \left(\unkuno,\left(\gggfw_{\elindexuno}\right)_{\elindexuno\in\setindexuno}\right)\in\domsk\times\left(\underset{\elindexuno\in\setindexuno}{\prod}\genfquot(\uedgtv_{\elindexuno})+1+\right)
\text{.}
\end{multline*}
We say that $\evpdf[\argcompl{\left( \corrarr[\pmV],\gfx,\ouol<seq<\right)}]$ is the evaluation of $\left(\predcRel_{\elindexuno}\right)_{\elindexuno\in\setindexuno}$ related to $\left( \corrarr[\pmV],\gfx,\ouol<seq<\right)$, or simply an evaluation of $\left(\predcRel_{\elindexuno}\right)_{\elindexuno\in\setindexuno}$whenever there is no risk of confusion.
\item Fix a finite set $\setindexuno$, $\tddmsk,\uedgtv\in \mathbb{N}_0$, $\vecdue_{\elindexuno}\in \mathbb{R}$,
$\left\{\left(\ls>\elindexuno,\antelss>, \ddgtv>\elindexuno,\antelss>, \dgtv>\elindexuno,\antelss>, \uedgtv,\ddgtv<seq<>\elindexuno,\antelss>, \ouol<seq<>\elindexuno,\antelss>\right)\right\}_{\antelss=1}^{\antels>\elindexuno>}\in\AQ+\tddmsk+$,
$\predcRel_{\elindexuno,\antelss}\in\predcR[\ls>\elindexuno,\antelss>, \ddgtv>\elindexuno,\antelss>, \dgtv>\elindexuno,\antelss>, \uedgtv,\ddgtv<seq<>\elindexuno,\antelss>,\ouol<seq<>\elindexuno,\antelss>]$ for any $\elindexuno \in  \setindexuno$, $\antelss\in\left\{1,...,\antels>\elindexuno>\right\}$, a skeleton $\left( \pmV,f,\ouol<seq<\right)$ of $\left(\left(\predcRel_{\elindexuno,\antelss}\right)_{\antelss=1}^{\antels>\elindexuno>}\right)_{\elindexuno \in  \setindexuno}$ with domain $\domsk$, a p-corresponding skeleton $\left( \corrarr[\pmV],\gfx,\ouol<seq<\right)$ of $\left(\left(\predcRel_{\elindexuno,\antelss}\right)_{\antelss=1}^{\antels>\elindexuno>}\right)_{\elindexuno \in  \setindexuno}$ related to $\left( \pmV,f,\ouol<seq<\right)$.\newline
Set $\predcRel= \underset{\elindexuno \in  \setindexuno}{\bigdirsum}\, \vecdue_{\elindexuno} \predcRscalp\left(\overset{\antels>\elindexuno>}{\underset{\antelss=1}{\bigtensprodR}} \predcRel_{\elindexuno,\antelss}\right)$.\newline
We define the arrow of $\GTopcat$ 
\begin{equation*}
\evpdf[\argcompl{\left( \corrarr[\pmV],\gfx,\ouol<seq<\right)}](\predcRel) :\domsk\times\genfquot(\uedgtv)+1+ \rightarrow \genfquot(\argcompl{\tddmsk+1})+1+
\end{equation*}
by setting
\begin{equation*}
\begin{array}{l}
\evpdf[\argcompl{\left( \corrarr[\pmV],\gfx,\ouol<seq<\right)}](\predcRel)\left(\unkuno,\gggfw\right)=
\left(\underset{\elindexuno\in\setindexuno}{\bigsumquotpoint}\,\vecdue_{\elindexuno} \scalpquotpoint\left(\underset{\antelss=1}{\overset{\antels>\elindexuno>}{\bigmultquotpoint}}
\evpdf[\argcompl{\left( \corrarr[\pmV]_{\elindexuno,\antelss},\gfx_{\elindexuno,\antelss},\ouol<seq<>\elindexuno,\antelss>\right)}]
\left(\unkuno,\gggfw\right)\right)\right)\\[4pt]
\hspace{197pt}\forall \left(\unkuno,\gggfw\right)\in\domsk\times\genfquot(\uedgtv)+1+\text{.}
\end{array}
\end{equation*}
We say that $\evpdf[\argcompl{\left( \corrarr[\pmV],\gfx,\ouol<seq<\right)}](\predcRel)$ is the evaluation of $\predcRel$ related to $\left( \corrarr[\pmV],\gfx,\ouol<seq<\right)$, or simply an evaluation of $\predcRel$ whenever there is no risk of confusion.\newline
Fix $\left(\unkuno,\gggfw\right)\in\domsk\times\genfquot(\uedgtv)+1+$, a path $\pmtwo$ in $\genf$ detecting $\evpdf[\argcompl{\left( \corrarr[\pmV],\gfx,\ouol<seq<\right)}](\predcRel)\left(\unkuno,\gggfw\right)$. We say that $\pmtwo$ is a path in $\genf$ detecting $\predcRel$ at $\left(\unkuno,\gggfw\right)$.
\end{enumerate}
\end{definition}

In Definition \ref{nullpredcRdef} below we introduce the notion of vanishing elements in $\predcR+\segnvar+<\uedgtv<$.  We refer to Remark \ref{siterem}-[3].

\begin{definition}\label{nullpredcRdef}\mbox{}
\begin{enumerate}
\item Fix $\uedgtv\in \mathbb{N}_0$, $\predcRel\in\predcR-1-+\argcompl{\segnvar,\segnvar}+<\uedgtv<$, $\gggfw\in\genfquot(\uedgtv)+1+$, a skeleton $\left( \pmV,f,\ouol<seq<\right)$ of $\predcRel$ with domain $\domsk$, a p-corresponding skeleton $\left( \corrarr[\pmV],\gfx,\ouol<seq<\right)$ of $\predcRel$ related to $\left( \pmV,f,\ouol<seq<\right)$ such that
\begin{equation*}
\evpdf[\left( \corrarr[\pmV],\gfx,\ouol<seq<\right)](\predcRel)\left(\unkuno,\gggfw\right)=0\scalpquotpoint\evpdf[\left( \corrarr[\pmV],\gfx,\ouol<seq<\right)](\predcRel)\left(\unkuno,\gggfw\right)\hspace{15pt}\forall\unkuno\in\domsk\text{.}
\end{equation*}
We say: $\predcRel$ vanishes at $\gggfw$ for $\left( \corrarr[\pmV],\gfx,\ouol<seq<\right)$; $\left( \corrarr[\pmV],\gfx,\ouol<seq<\right)$ realizes the vanishing of $\predcRel$ at $\gggfw$. With an abuse of language we say that $\predcRel$ is a vanishing element at $\gggfw$ whenever there is no need to specify the p-corresponding skeleton realizing the vanishing.\newline
We set:
\begin{flalign*}
&\begin{array}{l}
\nullpredcR-1-+\argcompl{\ddgtv,\dgtv}+<\uedgtv<(\gggfw)=\left\{\predcRel\in\predcR-1-+\argcompl{\ddgtv,\dgtv}+<\uedgtv<\;:\;\predcRel\hspace{4pt}\text{is}\hspace{4pt}\text{a}\hspace{4pt}\text{vanishing}\hspace{4pt}\text{element}\hspace{4pt}\text{at}\hspace{4pt}\gggfw\right\}\\[4pt]
\hspace{156pt}\text{for}\hspace{4pt}\text{any}\hspace{4pt}\ddgtv,\dgtv,\uedgtv \in \mathbb{N}_0 \text{,}\hspace{4pt} \gggfw\in\genfquot(\uedgtv)+1+\text{;}
\end{array}\\[6pt]
&\begin{array}{l}
\nullpredcR-1-+\argcompl{\segnvar,\dgtv}+<\uedgtv<(\gggfw)=\left\{\nullpredcR-1-+\argcompl{\ddgtv,\dgtv}+<\uedgtv<(\gggfw)\right\}_{\ddgtv\in\mathbb{N}_0 }
\hspace{7pt}\text{for}\hspace{4pt}\text{any}\hspace{4pt}\dgtv,\uedgtv \in \mathbb{N}_0 \text{,}\hspace{4pt} \gggfw\in\genfquot(\uedgtv)+1+\text{;}
\end{array}\\[6pt]
&\begin{array}{l}
\nullpredcR-1-+\argcompl{\segnvar,\segnvar}+<\uedgtv<(\gggfw)=\left\{\nullpredcR-1-+\argcompl{\ddgtv,\dgtv}+<\uedgtv<(\gggfw)\right\}_{\left(\ddgtv,\dgtv\right)\in\mathbb{N}_0 \times\mathbb{N}_0 }
\hspace{17pt}\text{for}\hspace{4pt}\text{any}\hspace{4pt}\uedgtv \in \mathbb{N}_0 \text{,}\hspace{4pt} \gggfw\in\genfquot(\uedgtv)+1+\text{;}
\end{array}\\[6pt]
&\begin{array}{l}
\nullpredcR-1-+\argcompl{\segnvar,\segnvar}+<\segnvar<(\gggfw)=\left\{\nullpredcR-1-+\argcompl{\ddgtv,\dgtv}+<\uedgtv<(\gggfw)\right\}_{\left(\ddgtv,\dgtv,\uedgtv\right)\in\mathbb{N}_0 \times\mathbb{N}_0 \times\mathbb{N}_0}
\hspace{17pt}\text{for}\hspace{4pt}\text{any}\hspace{4pt} \gggfw\in\genfquot(\uedgtv)+1+\text{.}
\end{array}
\end{flalign*}
If $\uedgtv=1$ and $\gggfw  =\lininclquot\left(\gczsfx[\idobj+\mathbb{R}+ ]\right)$ then we drop  $\lininclquot\left(\gczsfx[\idobj+\mathbb{R}+ ]\right)$ from the notation.
\item Fix $\uedgtv\in \mathbb{N}_0$, $\predcRel\in\predcR+\segnvar+<\uedgtv<$, $\gggfw\in\genfquot(\uedgtv)+1+$, a skeleton $\left( \pmV,f,\ouol<seq<\right)$ of $\predcRel$ with domain $\domsk$, a p-corresponding skeleton $\left( \corrarr[\pmV],\gfx,\ouol<seq<\right)$ of $\predcRel$ related to $\left( \pmV,\gfx,\ouol<seq<\right)$ such that
\begin{equation*}
\evpdf[\left( \corrarr[\pmV],\gfx,\ouol<seq<\right)](\predcRel)\left(\unkuno,\gggfw\right)=0\scalpquotpoint\evpdf[\left( \corrarr[\pmV],\gfx,\ouol<seq<\right)](\predcRel)\left(\unkuno,\gggfw\right)\hspace{15pt}\forall\unkuno\in\domsk\text{.}
\end{equation*}
We say: $\predcRel$ vanishes at $\gggfw$ for $\left( \corrarr[\pmV],\gfx,\ouol<seq<\right)$; $\left( \corrarr[\pmV],\gfx,\ouol<seq<\right)$ realizes the vanishing of $\predcRel$ at $\gggfw$. With an abuse of language we say that $\predcRel$ is a vanishing element at $\gggfw$ whenever there is no need to specify the p-corresponding skeleton realizing the vanishing.\newline
We set:
\begin{flalign*}
&\begin{array}{l}
\nullpredcR+\tddmsk+<\uedgtv<(\gggfw)=\left\{\predcRel\in\predcR+\tddmsk+<\uedgtv<\;:\;\predcRel\hspace{4pt}\text{is}\hspace{4pt}\text{a}\hspace{4pt}\text{vanishing}\hspace{4pt}\text{element}\hspace{4pt}\text{at}\hspace{4pt}\gggfw\right\}\hspace{15pt}\forall \tddmsk \in \mathbb{N}_0\text{;}
\end{array}\\[6pt]
&\begin{array}{l}
\nullpredcR+\segnvar+<\uedgtv<(\gggfw)=\left\{\nullpredcR+\tddmsk+<\uedgtv<(\gggfw)\right\}_{\tddmsk\in\mathbb{N}_0}\text{.}
\end{array}
\end{flalign*}
If $\uedgtv=1$ and $\gggfw  =\lininclquot\left(\gczsfx[\idobj+\mathbb{R}+ ]\right)$ then we drop  $\lininclquot\left(\gczsfx[\idobj+\mathbb{R}+ ]\right)$ from the notation, whenever there is no risk of confusion.
\end{enumerate}
\end{definition}

In Proposition \ref{preprevalideal} below we study the behavior of vanishing elements belonging to $\left(\predcR-1-+\argcompl{\segnvar, \segnvar}+<\segnvar<, \predcRsum, \predcRprod, \predcRscalp\right)$.

\begin{proposition}\label{preprevalideal}\mbox{}
\begin{enumerate}
\item Fix $\ddgtv, \dgtv,\uedgtv\in \mathbb{N}_0$, $\vecdue\in \mathbb{R}$, $\gggfw\in\genfquot(n)+1+$, $\predcRel \in \nullpredcR-1-+\argcompl{\ddgtv,\dgtv}+<\uedgtv<(\gggfw)$.\newline
Then $\vecdue \predcRscalp\predcRel\in \nullpredcR-1-+\argcompl{\ddgtv,\dgtv}+<\uedgtv<(\gggfw)$.  
\item Fix $\ddgtv, \dgtv,\uedgtv\in \mathbb{N}_0$, $\vecuno\in \mathbb{R}$, $\gggfw\in\genfquot(n)+1+$, $\predcRel_1,\predcRel_2 \in \nullpredcR-1-+\argcompl{\ddgtv,\dgtv}+<\uedgtv<(\gggfw)$.\newline
Then $\predcRel_1 \predcRsum\predcRel_2\in \nullpredcR-1-+\argcompl{\ddgtv,\dgtv}+<\uedgtv<(\gggfw)$. 
\item Fix $\ddgtv>1>,\ddgtv>2>,\dgtv>1>, \dgtv>2>,\uedgtv\in \mathbb{N}_0$, $\gggfw\in\genfquot(n)+1+$, $\predcRel_1 \in \predcR-1-+\argcompl{\ddgtv>1>,\dgtv>1>}+<\argcompl{\dgtv>2>+\uedgtv}<$, 
$\predcRel_2 \in\nullpredcR-1-+\argcompl{\ddgtv>2>,\dgtv>2>}+<\uedgtv<(\gggfw)$.\newline
Then $\predcRel_1 \predcRprod\predcRel_2\in \nullpredcR-1-+\argcompl{\ddgtv>1>+\ddgtv>2>,\dgtv>1>+\dgtv>2>}+<\uedgtv<(\gggfw)$.
\item Fix $\ddgtv>0>, \ddgtv, \dgtv>0>,\dgtv,\uedgtv\in \mathbb{N}_0$, $\gggfw\in\genfquot(n)+1+$, $\predcRel \in \nullpredcR-1-+\argcompl{\ddgtv,\dgtv}+<\uedgtv<(\gggfw)$, $\predcRel_0 \in\generderspace-funtore-(\ddgtv>0>)+\dgtv+\left(\dgtv>0>\right)$.\newline
Then $ \PREDER-\predcRel-+\predcRel_0+, \PRELOCAT-\predcRel-+\predcRel_0+\in \nullpredcR-1-+\argcompl{\ddgtv>0>+\ddgtv,\dgtv>0>+\dgtv}+<\uedgtv<(\gggfw)$.
\item Fix $\uedgtv\in \mathbb{N}_0$, $\gggfw\in\genfquot(n)+1+$.\newline
Then $\nullpredcR-1-+\argcompl{\segnvar,\segnvar}+<\uedgtv<(\gggfw)$ is a left ideal of the algebra
$\left(\predcR-1-+\argcompl{\segnvar, \segnvar}+<\segnvar<, \predcRsum, \predcRprod, \predcRscalp\right)$.
\end{enumerate}
\end{proposition}
\begin{proof}\mbox{}\newline
\textnormal{\textbf{Proof of statement 1.}}\ \ Fix a skeleton $\left( \pmV,f, \ouol<seq<\right)$ of $\predcRel$ with domain $\domsk$, a p-corresponding skeleton $\left( \corrarr[\pmV],\gfx, \ouol<seq<\right)$ of $\predcRel$ related to $\left( \pmV,f, \ouol<seq<\right)$ and realizing the vanishing of $\predcRel$ at $\gggfw$.\newline
Then: $\left( \pmV,f, \ouol<seq<\right)$ is a skeleton of $\vecdue \predcRscalp\predcRel$ with domain $\domsk$, by Definition \ref{corrrapp}-[3]; $\left( \corrarr[\pmV],\gfx, \ouol<seq<\right)$ is a p-corresponding skeleton  of $\vecdue \predcRscalp\predcRel$ related to $\left( \pmV,f, \ouol<seq<\right)$, by Definition \ref{corrrapp2}-[3].\newline
Eventually $\left( \corrarr[\pmV],\gfx, \ouol<seq<\right)$ realizes the vanishing of $\vecdue \predcRscalp\predcRel$ at $\gggfw$, by Definition \ref{geneval}-[3].\newline
\textnormal{\textbf{Proof of statement 2.}}\ \ Fix a decomposition vector $\left(\predcRel_{\elindexsei,\elindexuno}\right)_{\elindexuno\in \setindexuno_{\elindexsei}}$ of $\predcRel_{\elindexsei}$, a skeleton $\left( \pmV_{\elindexsei},f_{\elindexsei}, \ouol<seq<>\elindexsei>\right)$ of $\left(\predcRel_{\elindexsei,\elindexuno}\right)_{\elindexuno\in \setindexuno_{\elindexsei}}$ with domain $\domsk_{\elindexsei}$, a p-corresponding skeleton $\left( \corrarr[\pmV]_{\elindexsei},\gfx_{\elindexsei}, \ouol<seq<>\elindexsei>\right)$ of $\predcRel_{\elindexsei}$ related to $\left( \pmV_{\elindexsei},f_{\elindexsei}, \ouol<seq<>\elindexsei>\right)$ and realizing the vanishing of $\predcRel_{\elindexsei}$ at $\gggfw$ for any $\elindexsei\in\left\{1,2\right\}$.\newline
Locality of the argument entails that there is no loss of generality by assuming that there is an open neighborhood $\domsk$ of $0 \in \mathbb{R}^{\ddgtv+\dgtv}$ with $\domsk=\domsk_{1}=\domsk_{2}$.\newline
Set:
\begin{flalign*}
&\begin{array}{l}
\setindexuno \hspace{4pt}\text{the}\hspace{4pt}\text{disjoint}\hspace{4pt}\text{union}\hspace{4pt}\text{of}\hspace{4pt}\text{sets}\hspace{4pt}\setindexuno_1\text{,}\hspace{4pt}\setindexuno_2\text{;}
\end{array}\\[6pt]
&\begin{array}{l}
\left(\predcRel_{\elindexuno}\right)_{\elindexuno\in \setindexuno}\hspace{4pt}\text{by}\hspace{4pt} {\predcRel_{\elindexuno}}=\left\{
\begin{array}{ll}
\predcRel_{1,\elindexuno}
&\text{if}\hspace{4pt}\elindexuno\in \setindexuno_1\text{,}\\[4pt]
\predcRel_{2,\elindexuno}
&\text{if}\hspace{4pt}\elindexuno\in \setindexuno_2\text{;}
\end{array}
\right.
\end{array}\\[6pt]
&\begin{array}{l}
\left( \pmV,f, \ouol<seq<\right)=\left(\left( \pmV_{\elindexuno},f_{\elindexuno}, \ouol<seq<>\elindexuno>\right)\right)_{\elindexuno\in \setindexuno}\hspace{4pt}\text{where}\hspace{4pt} \left( \pmV_{\elindexuno},f_{\elindexuno}, \ouol<seq<>\elindexuno>\right)=\left\{
\begin{array}{ll}
\left( \pmV_{1,\elindexuno},f_{1,\elindexuno}, \ouol<seq<>1,\elindexuno>\right)&\text{if}\hspace{4pt}\elindexuno\in \setindexuno_1\text{,}\\[4pt]
\left( \pmV_{2,\elindexuno},f_{2,\elindexuno}, \ouol<seq<>2,\elindexuno>\right)&\text{if}\hspace{4pt}\elindexuno\in \setindexuno_2\text{;}
\end{array}
\right.
\end{array}\\[6pt]
&\begin{array}{l}
\left( \corrarr[\pmV],\gfx, \ouol<seq<\right)=\left(\left( \corrarr[\pmV]_{\elindexuno},\gfx_{\elindexuno}, \ouol<seq<>\elindexuno>\right)\right)_{\elindexuno\in \setindexuno}\hspace{4pt}\text{where}\hspace{4pt} \left( \corrarr[\pmV]_{\elindexuno},\gfx_{\elindexuno}, \ouol<seq<>\elindexuno>\right)=\left\{
\begin{array}{ll}
\left( \corrarr[\pmV]_{1,\elindexuno},\gfx_{1,\elindexuno}, \ouol<seq<>1,\elindexuno>\right)&\text{if}\hspace{4pt}\elindexuno\in \setindexuno_1\text{,}\\[4pt]
\left( \corrarr[\pmV]_{2,\elindexuno},\gfx_{2,\elindexuno}, \ouol<seq<>2,\elindexuno>\right)&\text{if}\hspace{4pt}\elindexuno\in \setindexuno_2\text{.}
\end{array}
\right.
\end{array}
\end{flalign*}
Then: $\left(\predcRel_{\elindexuno}\right)_{\elindexuno\in \setindexuno}$ is a decomposition vector of $\predcRel_1 \predcRsum\predcRel_2$; $\left( \pmV,f, \ouol<seq<\right)$ is a skeleton of $\left(\predcRel_{\elindexuno}\right)_{\elindexuno\in \setindexuno}$ with domain $\domsk$; 
$\left( \corrarr[\pmV],\gfx, \ouol<seq<\right)$ is a p-corresponding skeleton of $\left(\predcRel_{\elindexuno}\right)_{\elindexuno\in \setindexuno}$ related to 
$\left( \pmV,f, \ouol<seq<\right)$.\newline
Eventually $\left( \corrarr[\pmV],\gfx, \ouol<seq<\right)$ realizes the vanishing of $\predcRel_1 \predcRsum\predcRel_2$ at $\gggfw$, by construction.\newline
\textnormal{\textbf{Proof of statement 3.}}\ \ Fix a decomposition vector $\left(\predcRel_{\elindexsei,\elindexuno}  \right)_{\elindexuno\in \setindexuno_{\elindexsei}}$ of $\predcRel_{\elindexsei}$, a skeleton $\left( \pmV_{\elindexsei},f_{\elindexsei}, \ouol<seq<>\elindexsei>\right)$ of $\left(\predcRel_{\elindexsei,\elindexuno}\right)_{\elindexuno\in \setindexuno_{\elindexsei}}$ with domain $\domsk_{\elindexsei}$, a p-corresponding skeleton $\left( \corrarr[\pmV]_{\elindexsei},\gfx_{\elindexsei}, \ouol<seq<>\elindexsei>\right)$ of $\predcRel_{\elindexsei}$ related to $\left( \pmV_{\elindexsei},f_{\elindexsei}, \ouol<seq<>\elindexsei>\right)$ for any $\elindexsei\in\left\{1,2\right\}$. Assume that $\left( \corrarr[\pmV]_{2},\gfx_{2}, \ouol<seq<>2>\right)$ realizes the vanishing of $\predcRel_{2}$ at $\gggfw$.\newline
Define:
\begin{flalign*}
&\begin{array}{l}
\setindexuno =\setindexuno_1\times\setindexuno_2\text{;}
\end{array}\\[8pt]
&\begin{array}{l}
\left(\predcRel_{\left(\elindexuno_1,\elindexuno_2\right)}\right)_{\left(\elindexuno_1,\elindexuno_2\right)\in \setindexuno}\hspace{4pt}\text{by}\hspace{15pt} \predcRel_{\left(\elindexuno_1,\elindexuno_2\right)}= \predcRel_{\elindexuno_1}\predcRprod   \predcRel_{\elindexuno_2} \hspace{15pt}\forall \left(\elindexuno_1,\elindexuno_2\right)\in \setindexuno\text{;}
\end{array}\\[8pt] 
&\begin{array}{l}
\pmV_{\left(\elindexuno_1,\elindexuno_2\right),\lss}=\pmV_{1,\elindexuno_1,\lss}\funcomp\proj<\argcompl{\domsk>\argcompl{1,\dgtv>1>}>,\domsk>\argcompl{2,\dgtv>2>}>}<>1>\hspace{4pt}
\text{for}\hspace{4pt}\text{any}\hspace{4pt}\left(\elindexuno_1,\elindexuno_2\right) \in\setindexuno\text{,}\hspace{4pt}\lss\in\left\{1,...,\ls>\argcompl{1,\elindexuno_1}>\right\}\text{;}
\end{array}\\[8pt]
&\begin{array}{l}
f_{\left(\elindexuno_1,\elindexuno_2\right),\lss}=f_{1,\elindexuno_1,\lss}\funcomp\proj<\argcompl{\domsk>\argcompl{1,\dgtv>1>}>,\domsk>\argcompl{2,\dgtv>2>}>}<>1>\hspace{4pt}
\text{for}\hspace{4pt}\text{any}\hspace{4pt}\left(\elindexuno_1,\elindexuno_2\right) \in\setindexuno\text{,}\hspace{4pt}\lss\in\left\{1,...,\ls>\argcompl{1,\elindexuno_1}>\right\}\text{;}
\end{array}\\[8pt]
&\begin{array}{l}
\ouol>\left(\elindexuno_1,\elindexuno_2\right),\lss>=\ouol>1,\elindexuno_1,\lss>\hspace{4pt}
\text{for}\hspace{4pt}\text{any}\hspace{4pt}\left(\elindexuno_1,\elindexuno_2\right) \in\setindexuno\text{,}\hspace{4pt}\lss\in\left\{1,...,\ls>\argcompl{1,\elindexuno_1}>\right\}\text{;}
\end{array}\\[8pt]
&\left\{ 
\begin{array}{l}
\pmV_{\left(\elindexuno_1,\elindexuno_2\right),\lss}=\lininclquot\left(\pointedincl<\left\{\left(\mathbb{R}^{\dgtv>1>},0\right),\left(\mathbb{R}^{\dgtv>2>+\uedgtv},0\right)\right\}<>2>\right)\compquotpoint
\left(\pmV_{2, \elindexuno_2, \lss-\ls>\argcompl{1,\elindexuno_1}>}\funcomp
\proj<\argcompl{\domsk>\argcompl{1,\dgtv>1>}>,\domsk>\argcompl{2,\dgtv>2>}>}<>2>\right)\\[4pt]
\text{for}\hspace{4pt}\text{any}\hspace{4pt}\left(\elindexuno_1,\elindexuno_2\right) \in\setindexuno\text{,}\hspace{4pt}\lss\in\left\{\ls>\argcompl{1,\elindexuno_1}>+1,...,\ls>\argcompl{2,\elindexuno_2}>\right\}\text{;}
\end{array}
\right.\\[8pt]
&\left\{
\begin{array}{l}
f_{\left(\elindexuno_1,\elindexuno_2\right),\lss}=f_{2, \elindexuno_2, \lss-\ls>\argcompl{1,\elindexuno_1}>} \funcomp
\proj<\argcompl{\domsk>\argcompl{1,\dgtv>1>}>,\domsk>\argcompl{2,\dgtv>2>}>}<>2>  \\[4pt]
\text{for}\hspace{4pt}\text{any}\hspace{4pt}\left(\elindexuno_1,\elindexuno_2\right) \in\setindexuno\text{,}\hspace{4pt}\lss\in\left\{\ls>\argcompl{1,\elindexuno_1}>+1,...,\ls>\argcompl{2,\elindexuno_2}>\right\}\text{;}
\end{array}
\right.\\[8pt]
&\begin{array}{l}
\ouol>\left(\elindexuno_1,\elindexuno_2\right),\lss>=\ouol>\argcompl{2, \elindexuno_2, \lss-\ls>\argcompl{1,\elindexuno_1}>}>\hspace{4pt}
\text{for}\hspace{4pt}\text{any}\hspace{4pt}\left(\elindexuno_1,\elindexuno_2\right) \in\setindexuno\text{,}\hspace{4pt}\lss\in\left\{\ls>\argcompl{1,\elindexuno_1}>+1,...,\ls>\argcompl{2,\elindexuno_2}>\right\}\text{;}
\end{array}\\[8pt]
&\begin{array}{l}
\left( \pmV_{\left(\elindexuno_1,\elindexuno_2\right)},f_{\left(\elindexuno_1,\elindexuno_2\right)}, \ouol<seq<>\left(\elindexuno_1,\elindexuno_2\right)>\right)=\left(\left( \pmV_{\left(\elindexuno_1,\elindexuno_2\right),\lss},f_{\left(\elindexuno_1,\elindexuno_2\right),\lss}\right), \ouol<seq<>\left(\elindexuno_1,\elindexuno_2\right),\lss>\right)_{\lss=1}^{\ls>1>+\ls>2>}\hspace{15pt}\forall \left(\elindexuno_1,\elindexuno_2\right)\in \setindexuno\text{;}
\end{array}\\[8pt]
&\begin{array}{l}
\left( \pmV,f, \ouol<seq<\right)=\left(\left( \pmV_{\left(\elindexuno_1,\elindexuno_2\right)},f_{\left(\elindexuno_1,\elindexuno_2\right)}, \ouol<seq<>\left(\elindexuno_1,\elindexuno_2\right)>\right)\right)_{\left(\elindexuno_1,\elindexuno_2\right)\in \setindexuno}\text{;}
\end{array}\\[8pt]
&\left\{ 
\begin{array}{l}
 \corrarr[\pmV]_{\left(\elindexuno_1,\elindexuno_2\right),\lss}\left(\unkuno_{1}, \unkuno_{2}\right)=\corrarr[\pmV]_{1,\elindexuno_1,\lss}\left(\unkuno_{1}\right)  \compquotpoint  \lininclquot\left(
\gczsfx[\proj<\left\{\mathbb{R}^{\ddgtv>1, \elindexuno,\lss>},\mathbb{R}^{\dgtv>1>},\mathbb{R}^{\dgtv>2>}\right\}<>\left\{1,2\right\}>]\right)\\[4pt]
\text{for}\hspace{4pt}\text{any}\hspace{4pt}\left(\elindexuno_1,\elindexuno_2\right) \in\setindexuno\text{,}\hspace{4pt}\lss\in\left\{1,...,\ls>\argcompl{1,\elindexuno_1}>\right\}\text{,}\hspace{4pt}\left(\unkuno_{1}, \unkuno_{2}\right)\in \domsk>\argcompl{1,\dgtv>1>}>\times\domsk>\argcompl{2,\dgtv>2>}>\text{;}
\end{array}
\right.\\[8pt]
&\begin{array}{l}
\gfx_{\left(\elindexuno_1,\elindexuno_2\right),\lss}=\gfx_{1,\elindexuno_1,\lss}\compquotpoint\proj<\argcompl{\domsk>\argcompl{1,\dgtv>1>}>,\domsk>\argcompl{2,\dgtv>2>}>}<>1>\hspace{4pt}
\text{for}\hspace{4pt}\text{any}\hspace{4pt}\left(\elindexuno_1,\elindexuno_2\right) \in\setindexuno\text{,}\hspace{4pt}\lss\in\left\{1,...,\ls>\argcompl{1,\elindexuno_1}>\right\}\text{;}
\end{array}\\[8pt]
&\left\{
\begin{array}{l}
\corrarr[\pmV]_{\left(\elindexuno_1,\elindexuno_2\right),\lss}\left(\unkuno_{1}, \unkuno_{2}\right)=
\lininclquot\left(\gczsfx[\argcompl{\pointedincl<\left\{\left(\mathbb{R}^{\dgtv>1>},0\right),\left(\mathbb{R}^{\dgtv>2>+\uedgtv},0\right)\right\}<>2>}]\right)\compquotpoint\\[4pt]
\hspace{80pt}\Big( \corrarr[\pmV]_{2, \elindexuno_2, \lss-\ls>\argcompl{1,\elindexuno_1}>}\left(\unkuno_{2}\right)  \compquotpoint   
\lininclquot\left(
\gczsfx[\proj<\left\{\mathbb{R}^{\ddgtv>2, \elindexuno,\lss>},\mathbb{R}^{\dgtv>1>},\mathbb{R}^{\dgtv>2>}\right\}<>\left\{1,3\right\}>]\right)\Big)\\[4pt]
\text{for}\hspace{4pt}\text{any}\hspace{4pt}\left(\elindexuno_1,\elindexuno_2\right) \in\setindexuno\text{,}\hspace{4pt}\lss\in\left\{\ls>\argcompl{1,\elindexuno_1}>+1,...,\ls>\argcompl{2,\elindexuno_2}>\right\}\text{,}\hspace{4pt}\forall\left(\unkuno_{1}, \unkuno_{2}\right)\in \domsk>\argcompl{1,\dgtv>1>}>\times\domsk>\argcompl{2,\dgtv>2>}>\text{;}
\end{array}
\right.\\[8pt]
&\left\{
\begin{array}{l}
\gfx_{\left(\elindexuno_1,\elindexuno_2\right),\lss}=\gfx_{2, \elindexuno_2, \lss-\ls>\argcompl{1,\elindexuno_1}>} \compquotpoint
\proj<\argcompl{\domsk>\argcompl{1,\dgtv>1>}>,\domsk>\argcompl{2,\dgtv>2>}>}<>2> \\[4pt]
\text{for}\hspace{4pt}\text{any}\hspace{4pt}\left(\elindexuno_1,\elindexuno_2\right) \in\setindexuno\text{,}\hspace{4pt}\lss\in\left\{\ls>\argcompl{1,\elindexuno_1}>+1,...,\ls>\argcompl{2,\elindexuno_2}>\right\}\text{;}
\end{array}
\right.\\[8pt]
&\begin{array}{l}
\left( \corrarr[\pmV]_{\left(\elindexuno_1,\elindexuno_2\right)},\gfx_{\left(\elindexuno_1,\elindexuno_2\right)}, \ouol<seq<>\left(\elindexuno_1,\elindexuno_2\right)>\right)=\left(\left( \corrarr[\pmV]_{\left(\elindexuno_1,\elindexuno_2\right),\lss},\gfx_{\left(\elindexuno_1,\elindexuno_2\right),\lss}\right), \ouol<seq<>\left(\elindexuno_1,\elindexuno_2\right),\lss>\right)_{\lss=1}^{\ls>1>+\ls>2>}\hspace{15pt}\forall \left(\elindexuno_1,\elindexuno_2\right)\in \setindexuno\text{;}
\end{array}\\[8pt]
&\begin{array}{l}
\left( \corrarr[\pmV],\gfx, \ouol<seq<\right)=\left(\left( \corrarr[\pmV]_{\left(\elindexuno_1,\elindexuno_2\right)},\gfx_{\left(\elindexuno_1,\elindexuno_2\right)}, \ouol<seq<>\left(\elindexuno_1,\elindexuno_2\right)>\right)\right)_{\left(\elindexuno_1,\elindexuno_2\right)\in \setindexuno}\text{;}
\end{array}\\[8pt]
&\left\{
\begin{array}{l}
h: \mathbb{R}^{\ddgtv>1>}\times\mathbb{R}^{\ddgtv>2>}\times\mathbb{R}^{\dgtv>1>}\times\mathbb{R}^{\dgtv>2>}\times\mathbb{R}^{\uedgtv}\rightarrow \mathbb{R}^{\ddgtv>1>}\times\mathbb{R}^{\dgtv>1>}\times\mathbb{R}^{\ddgtv>2>}\times\mathbb{R}^{\dgtv>2>}\times\mathbb{R}^{\uedgtv}\\[4pt]
\text{by}\hspace{4pt}\text{setting}\\[4pt]
h\left(\unkuno_1, \unkuno_2, \unkuno_3, \unkuno_4, \unkuno_5\right)=\left(\unkuno_1, \unkuno_3, \unkuno_2, \unkuno_4, \unkuno_5\right)\\[4pt]
\hspace{60pt}\forall \left(\unkuno_1, \unkuno_2, \unkuno_3, \unkuno_4, \unkuno_5\right)\in\mathbb{R}^{\ddgtv>1>}\times\mathbb{R}^{\ddgtv>2>}\times\mathbb{R}^{\dgtv>1>}\times\mathbb{R}^{\dgtv>2>}\times\mathbb{R}^{\uedgtv}\text{;}
\end{array}
\right.\\[8pt] 
&\left\{
\begin{array}{l}
{_1 \corrarr[\pmV]}_{\elindexuno,\lss}=\lininclquot\left(\gczsfx[
\pointedincl<\left\{\left(\mathbb{R}^{\dgtv>1>},0\right), \left(\mathbb{R}^{\ddgtv>2>},0\right), \left(\mathbb{R}^{\dgtv>2>},0\right), \left(\mathbb{R}^{\uedgtv},0\right)\right\}<>\left\{1,3,4\right\}>
]\right) \compquotpoint 
\corrarr[\pmV]_{1,\elindexuno,\lss}
\\[4pt]
\text{for}\hspace{4pt}\text{any}\hspace{4pt}\elindexuno\in\setindexuno_1\text{,}\hspace{4pt}\lss\in\left\{1,...,\ls>\argcompl{1,\elindexuno}>\right\}\text{;}
\end{array}
\right.\\[8pt]
&\begin{array}{l}
\left(  {_1 \corrarr[\pmV]}_{\elindexuno}, \gfx_{\elindexuno}\right)=
\left(
\left( {_1 \corrarr[\pmV]}_{\elindexuno,\lss}, \gfx_{\elindexuno,\lss}\right)
\right)_{\lss\in\left\{1,...,\ls>\argcompl{1,\elindexuno}>\right\}}\hspace{4pt}
\text{for}\hspace{4pt}\text{any}\hspace{4pt}\elindexuno\in\setindexuno_1\text{;}
\end{array}\\[8pt]
&\begin{array}{l}
\left( {_1 \corrarr[\pmV]}, \gfx\right)=\left\{
\left(  {_1 \corrarr[\pmV]}_{\elindexuno},  \gfx_{\elindexuno}\right)\right\}_{\elindexuno\in\setindexuno_1}\text{.}
\end{array}
\end{flalign*}
Definition \ref{geneval}-[3] entails that
\begin{multline*}
\evpdf[\argcompl{\left( \corrarr[\pmV],\gfx,\ouol<seq<\right)}](\predcRel_1\predcRprod \predcRel_2)\left(\unkuno_1,\unkuno_2,\gggfw\right)=\\[4pt]
\evpdf[\argcompl{\left( {_1 \corrarr[\pmV]},\gfx,\ouol<seq<>1>\right)}](\predcRel_1)\Bigg(\unkuno_1,
\evpdf[\argcompl{\left( \corrarr[\pmV]_{2},\gfx_{2},\ouol<seq<>2>\right)}](\predcRel_2)\left(\unkuno_2,\gggfw\right)
\Bigg)\compquotpoint\lininclquot\left(\gczsfx[h]\right)\\[4pt]
\forall \left(\unkuno_1,\unkuno_2\right)\in \domsk>1>\times\domsk>2>\text{.} 
\end{multline*}
Eventually statement follows since by assumption we have
\begin{equation*}
\evpdf[\argcompl{\left( \corrarr[\pmV]_{2},\gfx_{2},\ouol<seq<>2>\right)}](\predcRel_2)\left(\unkuno_2,\gggfw\right)=
0\scalpquotpoint\evpdf[\argcompl{\left( \corrarr[\pmV]_{2},\gfx_{2},\ouol<seq<>2>\right)}](\predcRel_2)\left(\unkuno_2,\gggfw\right)\hspace{15pt}\forall \unkuno_2\in \domsk>2>\text{.} 
\end{equation*}
\textnormal{\textbf{Proof of statement 4.}}\ \ Statement follows straightforwardly by Definitions \ref{raevpunfun}, \ref{laevpunfun}, statements 2, 3.\newline
\textnormal{\textbf{Proof of statement 5.}}\ \ Statement follows by statements 1, 2, 3.  
\end{proof}

In Proposition \ref{anteprevalideal} below we study the behavior of $\nullpredcR-1-+\argcompl{\segnvar}+<1<(\gggfw)$ with respect to the algebraic structure of $\left(\predcR+\segnvar+<\uedgtv<, \predcRsum, \antedcRprod, \predcRscalp\right)$. We refer to Notations \ref{int}, \ref{realfunc}-[5(c)]

\begin{proposition}\label{anteprevalideal}\mbox{}
\begin{enumerate}
\item Fix $\tddmsk, \uedgtv\in \mathbb{N}_0$, $\vecdue\in \mathbb{R}$, $\gggfw\in\genfquot(\uedgtv)+1+$, $\predcRel \in \nullpredcR+\tddmsk+<\uedgtv<(\gggfw)$. 
Then $\vecdue \predcRscalp\predcRel\in\nullpredcR+\tddmsk+<\uedgtv<(\gggfw)$.  
\item Fix $\tddmsk, \uedgtv\in \mathbb{N}_0$, $\gggfw\in\genfquot(\uedgtv)+1+$, $\predcRel_1, \predcRel_2 \in \nullpredcR+\tddmsk+<\uedgtv<(\gggfw)$.
Then $\predcRel_1 \predcRsum\predcRel_2 \in \nullpredcR+\tddmsk+<\uedgtv<(\gggfw)$. 
\item Fix $\tddmsk, \uedgtv\in \mathbb{N}_0$, $\gggfw\in\genfquot(\uedgtv)+1+$, $\predcRel_1, \predcRel_2 \in \predcR+\tddmsk+<\uedgtv<$. Assume that $\predcRel_1$ or $\predcRel_2$ is a vanishing element at $\gggfw$. Then $\predcRel_1\antedcRprod\predcRel_2 \in \nullpredcR+\tddmsk+<\uedgtv<(\gggfw)$.
\item Fix $ \uedgtv\in \mathbb{N}_0$, $\gggfw\in\genfquot(\uedgtv)+1+$.\newline
Then $\nullpredcR+\segnvar+<\uedgtv<(\gggfw)$ is an ideal of the algebra $\left(\predcR+\segnvar+<\uedgtv<, \predcRsum, \antedcRprod, \predcRscalp\right)$.
\item Fix $ \uedgtv\in \mathbb{N}_0$, $\elindexuedgtv\in \left\{1,..., \overline{\uedgtv}\right\}$. Then $\nullpredcR+\segnvar+<\uedgtv<(\argcompl{\lininclquot\left(\gczsfx[\coord{\uedgtv}{\elindexuedgtv} ]\right)})$ is a proper ideal of the algebra
$\left(\predcR+\segnvar+<1<, \predcRsum, \antedcRprod, \predcRscalp\right)$.
\end{enumerate}
\end{proposition}
\begin{proof}\mbox{}\newline
Statements 1, 2, 3, 4 straightforwardly follow by Definitions \ref{ChDiR1}, \ref{nullpredcRdef}-[2].\newline
Statement 5 follows by proving the claim 
\begin{equation}
\left\{
\begin{array}{l}
\text{Fix}\hspace{4pt}\dgtv\in \mathbb{N}_0\text{.}\hspace{4pt}\text{Then}\hspace{4pt}\text{there}\hspace{4pt}\text{is}\hspace{4pt}\predcRel\in \predcR+\argcompl{\uedgtv+\dgtv}+<\uedgtv<\hspace{4pt}\text{such}\hspace{4pt}\text{that}\hspace{4pt}\text{for}\hspace{4pt}\text{any}\\[4pt]
\text{skeleton}\hspace{4pt}\left( \pmV,f,\ouol<seq<\right)\hspace{4pt}
\text{of}\hspace{4pt}\predcRel\hspace{4pt}
\text{any}\hspace{4pt}\text{p-corresponding}\hspace{4pt}\text{skeleton}\\[4pt]
\left( \corrarr[\pmV],\gfx,\ouol<seq<\right)\hspace{4pt}\text{related}
\hspace{4pt}\text{to}\hspace{4pt}\left( \pmV,f,\ouol<seq<\right) \hspace{4pt}
\text{does}\hspace{4pt}\text{not}\hspace{4pt}\text{realize}\hspace{4pt}\text{the}\hspace{4pt}\text{vanishing}\hspace{4pt}\text{of}\hspace{4pt}\predcRel\\[4pt]
\text{at}\hspace{4pt}\lininclquot\left(\gczsfx[\coord{\uedgtv}{\elindexuedgtv}]\right)\text{.}
\end{array}
\right.\label{maivano}
\end{equation}
Referring to Notations \ref{varispder}-[2], \ref{realvec}-[2], Remarks \ref{pathinlpreder}-[2], \ref{pathinlgenerder}-[2], \eqref{pcM} we define $\predcRel=\gdgerm[\argcompl{\quotdirsumdue\funcomp\quotdirsumuno\funcomp\cost<\mathbb{R}^{\dgtv}<>\argcompl{\dirsumgenpderspace(\uedgtv)+\uedgtv+}>+\argcompl{\genpreder\left[\lininclquot\left(\gczsfx[\idobj+\mathbb{R}^{\uedgtv}+]\right),  \ef_{\elindexuedgtv}\right]}+}]$. \newline
By Definitions \ref{isoinddef}, \ref{dualnot1} of corresponding arrows and of p-corresponding arrows respectively there is one and only one evaluation of $\predcRel$ at $\lininclquot\left(\gczsfx[\coord{\uedgtv}{\elindexuedgtv} ]\right)$, namely the arrow $\cost<\mathbb{R}^{\uedgtv+\dgtv}<>\argcompl{\genfquot(\argcompl{\uedgtv+\dgtv+\uedgtv})+1+}>+\argcompl{ \lininclquot\left(\gczsfx[\argcompl{\cost<\argcompl{\mathbb{R}^{\uedgtv+\dgtv+\uedgtv}}<>\mathbb{R}>+1+}]\right)}+$.\newline
Eventually $\predcRel$ straightforwardly verifies claim \eqref{maivano}.
\end{proof}

In Definition \ref{connid} below we point the attention on a class of elements in $\predcR+\segnvar+<1<$ which will be crucial in development of calculus in $\GTopcat$.

\begin{definition} \label{connid} 
We define sets
\begin{equation*}
\gennullwpredcR+\tddmsk+\subseteq \predcR+\tddmsk+<1< \hspace{15pt} \tddmsk\in \mathbb{N}_0
\end{equation*}
recursively as follows.
\begin{enumerate}
\item Fix  $\ddgtv,\dgtv \in \mathbb{N}_0$, $\vecdue_1,\vecdue_2\in \mathbb{R}$, $\vecuno\in \mathbb{R}^{\ddgtv}$, an open neighborhood $\domsk\subseteq \mathbb{R}^{\dgtv}$ of $0\in\mathbb{R}^{\dgtv}$, $f_1,f_2\in \Cksp{0}(\domsk)+\mathbb{R}+$, $\pmV\in \genfquotcccz(\ddgtv)+\dgtv+$.\newline
Define:
\begin{flalign*}
&\begin{array}{l}
\generder=\quotdirsumdue\left(\quotdirsumuno\left( \genpreder\left[\pmV,\vecuno \right]\right)\right)\text{;}
\end{array}\\[8pt]
&\left\{
\begin{array}{l}
\text{the}\hspace{4pt}\text{arrow}\hspace{4pt}\pmgenerderdue_{\elindexquattro}:   \domsk   \rightarrow\generderspace(\ddgtv)+\argcompl{\dgtv+1}+\hspace{4pt}\text{by}\hspace{4pt}\text{setting}\\[4pt]
\pmgenerderdue_{\elindexquattro}\left(\elsymtre\right)=
\gendiff(\ddgtv)+\argcompl{\gczsfx[\ffuncaap{\left(\pointedincl<\left\{\left(\mathbb{R}^{\dgtv},0\right),\left(\mathbb{R},0\right)\right\}<>2>\funcomp f_{\elindexquattro}\right)}\left(\elsymtre\right)]}+ \left(\generder\right)  \hspace{15pt}\forall   \elsymtre\in\domsk\\[4pt]
\text{for}\hspace{4pt}\text{any}\hspace{4pt}\elindexquattro\in \left\{1,2\right\}\text{;}
\end{array}
\right.\\[8pt]
&\left\{
\begin{array}{l}
\text{the}\hspace{4pt}\text{arrow}\hspace{4pt}\pmgenerderdue_{3}:   \domsk   \rightarrow\generderspace(\ddgtv)+\argcompl{\dgtv+1}+\hspace{4pt}\text{by}\hspace{4pt}\text{setting}\\[4pt]
\pmgenerderdue_{3}\left(\elsymtre\right)=\\[4pt]
\hspace{5pt}\gendiff(\ddgtv)+\argcompl{\gczsfx[ \ffuncaap{\left(\pointedincl<\left\{\left(\mathbb{R}^{\dgtv},0\right),\left(\mathbb{R},0\right)\right\}<>2>\funcomp \left(\vecdue_1\funscalp f_{1} \funsum \vecdue_2\funscalp  f_{2}\right)\right)}\left(\elsymtre\right)]}+ \left(\generder\right)\\[4pt]
\hspace{240pt}\forall\elsymtre\in\domsk\text{.}
\end{array}
\right.
\end{flalign*}
Then:
\begin{flalign*}
&\begin{array}{l}
\left(\gdgerm[\pmgenerderdue_3]\left(0\right)\right)_0 \predcRsum -1 \predcRscalp
\left(\vecdue_1\predcRscalp\left(\gdgerm[\pmgenerderdue_1]\left(0\right)\right)_0 \predcRsum  \vecdue_2\predcRscalp\left(\gdgerm[\pmgenerderdue_2]\left(0\right)\right)_0 \right)\in \gennullwpredcR+\argcompl{\ddgtv+\dgtv}+\text{;}
\end{array}\\[8pt]
&\begin{array}{l}
\left(\gdgerm[\pmgenerderdue_3]\left(0\right)\right)_1 \predcRsum -1 \predcRscalp
\left(\vecdue_1\predcRscalp\left(\gdgerm[\pmgenerderdue_1]\left(0\right)\right)_1 \predcRsum  \vecdue_2\predcRscalp\left(\gdgerm[\pmgenerderdue_2]\left(0\right)\right)_1 \right)\in \gennullwpredcR+\argcompl{\ddgtv+\dgtv}+\text{.}
\end{array}
\end{flalign*}
\item Fix  $\ddgtv,\dgtv \in \mathbb{N}_0$, $\vecuno\in \mathbb{R}^{\ddgtv}$, an open neighborhood $\domsk\subseteq \mathbb{R}^{\dgtv}$ of $0\in\mathbb{R}^{\dgtv}$, $f_1,f_2\in \Cksp{0}(\domsk)+\mathbb{R}+$, $\pmV\in \genfquotcccz(\ddgtv)+\dgtv+$.\newline
Define:
\begin{flalign*}
&\begin{array}{l}
\generder=\quotdirsumdue\left(\quotdirsumuno\left( \genpreder\left[\pmV,\vecuno \right]\right)\right)\text{;}
\end{array}\\[8pt]
&\left\{
\begin{array}{l}
\text{the}\hspace{4pt}\text{arrow}\hspace{4pt}\pmgenerderdue_{\elindexquattro}:   \domsk   \rightarrow\generderspace(\ddgtv)+\argcompl{\dgtv+1}+\hspace{4pt}\text{by}\hspace{4pt}\text{setting}\\[4pt]
\pmgenerderdue_{\elindexquattro}\left(\elsymtre\right)=
\gendiff(\ddgtv)+\argcompl{\gczsfx[\ffuncaap{\left(\pointedincl<\left\{\left(\mathbb{R}^{\dgtv},0\right),\left(\mathbb{R},0\right)\right\}<>2>\funcomp
f_{\elindexquattro}\right)}\left(\elsymtre\right)]}+ \left(\generder\right)  \hspace{15pt}\forall   \elsymtre\in\domsk\\[4pt]
\text{for}\hspace{4pt}\text{any}\hspace{4pt}\elindexquattro\in \left\{1,2\right\}\text{;}
\end{array}
\right.\\[8pt]
&\left\{
\begin{array}{l}
\text{the}\hspace{4pt}\text{arrow}\hspace{4pt}\pmgenerderdue_{3}:   \domsk   \rightarrow\generderspace(\ddgtv)+\argcompl{\dgtv+1}+\hspace{4pt}\text{by}\hspace{4pt}\text{setting}\\[4pt]
\pmgenerderdue_{3}\left(\elsymtre\right)=
\gendiff(\ddgtv)+\argcompl{\gczsfx[ \ffuncaap{\pointedincl<\left\{\left(\mathbb{R}^{\dgtv},0\right),\left(\mathbb{R},0\right)\right\}<>2>]\funcomp\left( f_{1} \funmult   f_{2}\right)}\left(\elsymtre\right)]}+ \left(\generder\right)\hspace{15pt}\forall\elsymtre\in\domsk\text{.}
\end{array}
\right.
\end{flalign*}
Then
\begin{multline*}
\left(\gdgerm[\pmgenerderdue_3]\left(0\right)\right)_1 \predcRsum -1 \predcRscalp
\Big(  \left(\gdgerm[\pmgenerderdue_1]\left(0\right)\right)_0  \antedcRprod\left(\gdgerm[\pmgenerderdue_2]\left(0\right)\right)_1 \predcRsum\\[4pt]  \left(\gdgerm[\pmgenerderdue_1]\left(0\right)\right)_1 \antedcRprod\left(\gdgerm[\pmgenerderdue_2]\left(0\right)\right)_0 \Big)\in \gennullwpredcR+\argcompl{\ddgtv+\dgtv}+\text{.}
\end{multline*}
\end{enumerate}
\end{definition}

In Proposition \ref{connidprop} below we prove that elements introduced in  Definition \ref{connid} belongs to $\nullpredcR+\segnvar+<1<(\argcompl{\lininclquot\left(\gczsfx[\idobj+\mathbb{R}+ ]\right)})$. We refer to Remark \ref{siterem}-[3].

\begin{proposition} \label{connidprop}
Fix $\tddmsk\in \mathbb{N}_0$. Then $\gennullwpredcR+\tddmsk+\subseteq \nullpredcR+\tddmsk+<1<(\argcompl{\lininclquot\left(\gczsfx[\idobj+\mathbb{R}+ ]\right)})$.
\end{proposition}
\begin{proof} We prove the claim 
\begin{equation}
\left\{
\begin{array}{l}
\text{Fix}\hspace{4pt}\text{a}\hspace{4pt}\text{finite}\hspace{4pt}\text{set}\hspace{4pt}\setindexuno\text{,}\hspace{4pt} \left\{\predcRel_{\elindexuno}\right\}_{\elindexuno \in \setindexuno}\subseteq \gennullwpredcR+\tddmsk+\text{.}\hspace{4pt}
\text{Then}\hspace{4pt}\text{there}\hspace{4pt}\text{is}\hspace{4pt}\text{a}\hspace{4pt}\text{skeleton}\\[4pt]\left( \pmV,\vecuno,\ouol<seq<\right)\hspace{4pt}\text{of}\hspace{4pt}\left(\predcRel_{\elindexuno}\right)_{\elindexuno \in \setindexuno}\hspace{4pt} 
\text{realizing}\hspace{4pt}\text{the}\hspace{4pt}\text{vanishing}\hspace{4pt}\text{of}\hspace{4pt}\predcRel_{\elindexuno} \hspace{4pt}\text{at}\hspace{4pt}
\lininclquot\left(\gczsfx[\idobj+\mathbb{R}+]\right)\hspace{4pt}\text{for}\\[4pt]
\text{any}\hspace{4pt}\elindexuno\in \setindexuno  \text{.}
\end{array}
\right.\label{Xèzero}
\end{equation}
Referring to Notation \ref{ins}-[14], Remark \ref{siteremdue} by construction of $\gennullwpredcR+\tddmsk+$ locality of the argument entails that there are a partition $\left\{\setindexuno_0, \setindexuno_1, \setindexuno_2\right\}$ of $\setindexuno$, $\ddgtv>\elindexuno>,\uedgtv>\elindexuno>\in \mathbb{N}_0$,  $\vecuno_{\elindexuno}\!\in\! \mathbb{R}^{\ddgtv>\elindexuno>}$, $f_{1,\elindexuno},f_{2,\elindexuno}\!\in\!  \Cksp{0}(\domsk>\argcompl{\uedgtv>\elindexuno>}>)+\mathbb{R}+$, $\gggfy_{\elindexuno}\!\in\! \genfquotcccz(\ddgtv>\elindexuno>)+\uedgtv>\elindexuno>+$, $\pmgenerderdue_{1,\elindexuno}, \pmgenerderdue_{2,\elindexuno}, \pmgenerderdue_{3,\elindexuno}\!\in\! \homF+\GTopcat+\left( \domsk>\argcompl{\uedgtv>\elindexuno>}>  ,\generderspace(\ddgtv>\elindexuno>)+\argcompl{\uedgtv>\elindexuno>+1}+\right)$, $\generder_{\elindexuno}\in \generderspace(\ddgtv>\elindexuno>)+\argcompl{\uedgtv>\elindexuno>}+$, 
$\predcRel_{1,\elindexuno}, \predcRel_{2,\elindexuno}, \predcRel_{3,\elindexuno}\in \generderspace-funtore-(\argcompl{\ddgtv>\elindexuno>})+1+\left(\uedgtv>\elindexuno>\right)$, for any $\elindexuno\in\setindexuno$, $\vecdue_{1,\elindexuno},\vecdue_{2,\elindexuno}\in \mathbb{R}$ for any $\elindexuno\in\setindexuno_0\cup\setindexuno_1$ fulfilling all conditions below:
\begin{flalign*}
&\begin{array}{l}
\predcRel_{\elindexuno}=\\
\hspace{10pt}\left\{
\begin{array}{ll}
\left(\predcRel_{3,\elindexuno}\right)_0\;\predcRsum\;-1\predcRscalp\left(\vecdue_1\predcRscalp\left(\predcRel_{1,\elindexuno}
\right)_0\predcRsum\vecdue_2\predcRscalp\left(\predcRel_{2,\elindexuno}\right)_0\right)&\text{if}\hspace{4pt}\elindexuno\in\setindexuno_0\text{,}\\[8pt]
\left(\predcRel_{3,\elindexuno}\right)_1\;\predcRsum\;-1\predcRscalp\left(\vecdue_1\predcRscalp\left(\predcRel_{1,\elindexuno}
\right)_1\predcRsum\vecdue_2\predcRscalp\left(\predcRel_{2,\elindexuno}\right)_1\right)&\text{if}\hspace{4pt}\elindexuno\in\setindexuno_1\text{,}\\[8pt]
\left(\predcRel_{3,\elindexuno}\right)_1\;\predcRsum\;-1\predcRscalp\Big(\left(\predcRel_{1,\elindexuno}\right)_0\antedcRprod\left(\predcRel_{2,\elindexuno}\right)_1\predcRsum 
\left(\predcRel_{1,\elindexuno}\right)_1\antedcRprod\left(\predcRel_{2,\elindexuno}\right)_0\Big)&\text{if}\hspace{4pt}\elindexuno\in\setindexuno_2\text{;}
\end{array}
\right.
\end{array}\\[8pt]
&\begin{array}{l}
\generder_{\elindexuno}=\quotdirsumdue\left(\quotdirsumuno\left( \genpreder\left[\gggfy_{\elindexuno},\vecuno_{\elindexuno} \right]\right)\right)
\hspace{15pt} \forall \elindexuno \in\setindexuno\text{;}
\end{array}\\[8pt]
&\begin{array}{l}
\predcRel_{\elindexsei,\elindexuno}=\gdgerm[\pmgenerderdue_{\elindexsei,\elindexuno}]\left(0\right)
\hspace{4pt}\text{for}\hspace{4pt}\text{any}\hspace{4pt}\elindexsei	\in \left\{1,2,3\right\} \text{,}\hspace{4pt}\elindexuno \in\setindexuno\text{;}
\end{array}\\[8pt]
&\begin{array}{l}
\pmgenerderdue_{\elindexsei,\elindexuno}\left(\elsymtre\right)=
\gendiff(\ddgtv>0>)+\argcompl{\gczsfx[\ffuncaap{\left(\pointedincl<\left\{\left(\mathbb{R}^{\uedgtv>\elindexuno>},0\right),\left(\mathbb{R},0\right)\right\}<>2>\funcomp
f_{\elindexsei,\elindexuno}\right)}\left(\elsymtre\right)]}+ \left(\generder_{\elindexuno}\right)  \\[4pt]
\hspace{160pt}\text{for}\hspace{4pt}\text{any}\hspace{4pt}\elindexsei\in \left\{1,2\right\}\text{,}\hspace{4pt}\elindexuno\in\setindexuno\text{,}\hspace{4pt}
\elsymtre\in\domsk>\argcompl{\uedgtv>\elindexuno>}>\text{;}
\end{array}\\[8pt]
&\begin{array}{l}
\pmgenerderdue_{3,\elindexuno}\left(\elsymtre\right)=\\[4pt]
\hspace{10pt}\left\{
\begin{array}{l}
\gendiff(\ddgtv>0>)+\argcompl{\gczsfx[ \ffuncaap{\left(\pointedincl<\left\{\left(\mathbb{R}^{\uedgtv>\elindexuno>},0\right),\left(\mathbb{R},0\right)\right\}<>2>\funcomp \left(\vecdue_{1,\elindexuno}\funscalp f_{1,\elindexuno} \funsum \vecdue_{2,\elindexuno}\funscalp  f_{2,\elindexuno}\right)\right)}\left(\elsymtre\right)]}+ \left(\generder_{\elindexuno}\right)\\[4pt]
\hspace{160pt}\text{for}\hspace{4pt}\text{any}\hspace{4pt}   \elindexuno \in\setindexuno_0\cup\setindexuno_1\text{,}\hspace{4pt}
\elsymtre\in\domsk>\argcompl{\uedgtv>\elindexuno>}>\text{,}\\[8pt]
\gendiff(\ddgtv>0>)+\argcompl{\gczsfx[ \ffuncaap{\pointedincl<\left\{\left(\mathbb{R}^{\uedgtv>\elindexuno>},0\right),\left(\mathbb{R},0\right)\right\}<>2>]\funcomp\left( f_{1,\elindexuno} \funmult   f_{2,\elindexuno}\right)}\left(\elsymtre\right)]}+ \left(\generder_{\elindexuno}\right)\\[4pt]
\hspace{160pt}\text{for}\hspace{4pt}\text{any}\hspace{4pt}   \elindexuno \in\setindexuno_2 \text{,}\hspace{4pt}
\elsymtre\in\domsk>\argcompl{\uedgtv>\elindexuno>}> \text{.} 
\end{array}
\right.
\end{array}
\end{flalign*}
Fix $\gfx_{1,\elindexuno},\gfx_{2,\elindexuno}\in  \genfcontcont(\domsk>\argcompl{\uedgtv>\elindexuno>}>)+\mathbb{R}+$ for any $\elindexuno\in \setindexuno$ fulfilling condition below:
\begin{equation*}
\evalcomptwo\left(\gfx_{\elindexsei,\elindexuno}\right)=f_{\elindexsei,\elindexuno}\hspace{4pt}\text{for}\hspace{4pt}\text{any}\hspace{4pt}\elindexsei\left\{1,2\right\}\text{,}\hspace{4pt}\elindexuno\in \setindexuno\text{.}
\end{equation*}
Define the arrows belonging to $\homF+\GTopcat+\left( \domsk>\argcompl{\uedgtv>\elindextre>}>  ,\genfquotcccz(\ddgtv>\elindextre>)+\argcompl{\uedgtv>\elindextre>+1}+\right)$ by setting: 
\begin{flalign*}
&\begin{array}{l}
\pmV_{\elindexsei,\elindexuno}\left(\unkuno\right)=\ffuncaap{\gfx_{\elindexsei,\elindexuno}}\left(\unkuno\right) \compquotpoint \gggfy_{\elindexuno} \hspace{4pt}\text{for}\hspace{4pt}\text{any}\hspace{4pt} \elindexsei\in\left\{1,2\right\}\text{,}\hspace{4pt}\elindexuno\in \setindexuno\text{,}\hspace{4pt}\unkuno\in \domsk>\argcompl{\uedgtv>\elindexuno>}>\text{;}
\end{array}\\[8pt]
&\begin{array}{l}
\pmV_{3,\elindexuno}\left(\unkuno\right)=\vecdue_{1,\elindexuno}\scalpquotpoint\pmV_{1,\elindexuno}\left(\unkuno\right)\sumquotpoint\vecdue_{2,\elindexuno}\scalpquotpoint\pmV_{2,\elindexuno}\left(\unkuno\right)
 \hspace{4pt}\text{for}\hspace{4pt}\text{any}\hspace{4pt} \elindexuno\in \setindexuno_0\cup\setindexuno_1\text{,}\hspace{4pt}\unkuno\in \domsk>\argcompl{\uedgtv>\elindexuno>}>\text{;}
\end{array}\\[8pt]
&\begin{array}{l}
\pmV_{3,\elindexuno}\left(\unkuno\right)=\pmV_{1,\elindexuno}\left(\unkuno\right)\multquotpoint \pmV_{2,\elindexuno}\left(\unkuno\right)
 \hspace{4pt}\text{for}\hspace{4pt}\text{any}\hspace{4pt} \elindexuno\in \setindexuno_2\text{,}\hspace{4pt}\unkuno\in \domsk>\argcompl{\uedgtv>\elindexuno>}>\text{.}
\end{array}
\end{flalign*}
Then Proposition \ref{preimmerfun}-[4], locality of the argument together entail that 
\begin{equation*}
\begin{array}{l}
\ihtl{\domsk>\argcompl{\uedgtv>\elindexuno>}>}{\ddgtv>\elindexuno>}{\uedgtv>\elindexuno>}\left(\pmV_{\elindexsei,\elindexuno}\right)  \neq \udenset \hspace{4pt}\text{for}\hspace{4pt}\text{any}\hspace{4pt}\elindexsei\in\left\{1,2\right\}\text{,}\hspace{4pt}\elindexuno\in \setindexuno\text{.}
\end{array}
\end{equation*}
Hence:
\begin{flalign*}
&\left\{
\begin{array}{l}
\left(\lininclquot\left(\gczsfx[\pointedincl<\left\{\left(\mathbb{R}^{\uedgtv>\elindexuno>},0\right),\left(\mathbb{R},0\right)\right\}<>2>]\right)\compgenfquot \pmV_{\elindexsei,\elindexuno},\vecuno_{\elindexuno},0\right)\hspace{4pt}\text{is}\hspace{4pt}\text{a}\hspace{4pt}\text{skeleton}\hspace{4pt}\text{of}\hspace{4pt}\left(\predcRel_{\elindexsei,\elindexuno}\right)_0 \\[4pt]\text{with}\hspace{4pt}\text{domain}\hspace{4pt}\domsk\hspace{4pt}
\text{for}\hspace{3pt}\text{any}\hspace{4pt}
\elindexsei\in\left\{1,2\right\}\text{,}\hspace{4pt}\elindexuno\in\setindexuno_0\cup\setindexuno_2\text{;}
\end{array}
\right.\\[8pt]
&\left\{
\begin{array}{l}
\left(\lininclquot\left(\gczsfx[\pointedincl<\left\{\left(\mathbb{R}^{\uedgtv>\elindexuno>},0\right),\left(\mathbb{R},0\right)\right\}<>2>]\right)\compgenfquot \pmV_{\elindexsei,\elindexuno},\vecuno_{\elindexuno},1\right)\hspace{4pt}\text{is}\hspace{4pt}\text{a}\hspace{4pt}\text{skeleton}\hspace{4pt}\text{of}\hspace{4pt}\left(\predcRel_{\elindexsei,\elindexuno}\right)_1 \\[4pt]\text{with}\hspace{4pt}\text{domain}\hspace{4pt}\domsk\hspace{4pt}
\text{for}\hspace{3pt}\text{any}\hspace{4pt}
\elindexsei\in\left\{1,2\right\}\text{,}\hspace{4pt}\elindexuno\in\setindexuno_1\cup\setindexuno_2\text{;}
\end{array}
\right.\\[8pt]
&\left\{
\begin{array}{l}
\left(\lininclquot\left(\gczsfx[\pointedincl<\left\{\left(\mathbb{R}^{\uedgtv>\elindexuno>},0\right),\left(\mathbb{R},0\right)\right\}<>2>]\right)\compgenfquot \pmV_{3,\elindexuno},\vecuno_{\elindexuno},0\right)\hspace{4pt}\text{is}\hspace{4pt}
\text{a}\hspace{4pt}\text{skeleton}\hspace{4pt}\text{of}\hspace{4pt}\left(\predcRel_{3,\elindexuno}\right)_0\\[4pt] 
\text{with}\hspace{4pt}\text{domain}\hspace{4pt}\domsk\hspace{4pt}
\text{for}\hspace{3pt}\text{any}\hspace{4pt}\elindexuno\in\setindexuno_0\text{;}
\end{array}
\right.\\[8pt]
&\left\{
\begin{array}{l}
\left(\lininclquot\left(\gczsfx[\pointedincl<\left\{\left(\mathbb{R}^{\uedgtv>\elindexuno>},0\right),\left(\mathbb{R},0\right)\right\}<>2>]\right)\compgenfquot \pmV_{3,\elindexuno},\vecuno_{\elindexuno},1\right)\hspace{4pt}\text{is}\hspace{4pt}
\text{a}\hspace{4pt}\text{skeleton}\hspace{4pt}\text{of}\hspace{4pt}\left(\predcRel_{3,\elindexuno}\right)_1 \\[4pt]
\text{with}\hspace{4pt}\text{domain}\hspace{4pt}\domsk\hspace{4pt}
\text{for}\hspace{3pt}\text{any}\hspace{4pt}\elindexuno\in\setindexuno_1\text{;}
\end{array}
\right.\\[8pt]
&\left\{
\begin{array}{l}
\left(\lininclquot\left(\gczsfx[\pointedincl<\left\{\left(\mathbb{R}^{\uedgtv>\elindexuno>},0\right),\left(\mathbb{R},0\right)\right\}<>2>]\right)\compgenfquot \pmV_{3,\elindexuno},\vecuno_{\elindexuno},1\right)\hspace{4pt}\text{is}\hspace{4pt}
\text{a}\hspace{4pt}\text{skeleton}\hspace{4pt}\text{of}\hspace{4pt}\left(\predcRel_{3,\elindexuno}\right)_1\\[4pt]
\text{with}\hspace{4pt}\text{domain}\hspace{4pt}\domsk\hspace{4pt}
\text{for}\hspace{3pt}\text{any}\hspace{4pt}\elindexuno\in\setindexuno_2\text{.}
\end{array}
\right.
\end{flalign*}
Fix $\corrarr[\pmV]_{\elindexsei,\elindexuno}\in  \homF+\GTopcat+\left( \domsk>\argcompl{\uedgtv>\elindexuno>}>  ,\genfquotcccz(\tddmsk)+\argcompl{\uedgtv>\elindexuno>+1}+\right)$ for any $\elindexsei\in\left\{1,2\right\}$, $\elindextre\in \setindexuno$.\newline
Define the arrows belonging to $\homF+\GTopcat+\left( \domsk>\argcompl{\uedgtv>\elindexuno>}>  ,\genfquotcccz(\tddmsk)+\argcompl{\uedgtv>\elindexuno>+1}+\right)$ by setting: 
\begin{flalign*}
&\begin{array}{l}
\corrarr[\pmV]_{3,\elindexuno}\left(\unkuno\right)=\vecdue_{1,\elindexuno}\scalpquotpoint\corrarr[\pmV]_{1,\elindexuno}\left(\unkuno\right)\sumquotpoint\vecdue_{2,\elindexuno}\scalpquotpoint\corrarr[\pmV]_{2,\elindexuno}\left(\unkuno\right)
 \hspace{4pt}\text{for}\hspace{4pt}\text{any}\hspace{4pt} \elindexuno\in \setindexuno_0\cup\setindexuno_1\text{,}\hspace{4pt}\unkuno\in \domsk>\argcompl{\uedgtv>\elindexuno>}>\text{;}
\end{array}\\[8pt]
&\begin{array}{l}
\corrarr[\pmV]_{3,\elindexuno}\left(\unkuno\right)=\corrarr[\pmV]_{1,\elindexuno}\left(\unkuno\right)\multquotpoint \corrarr[\pmV]_{2,\elindexuno}\left(\unkuno\right)
 \hspace{4pt}\text{for}\hspace{4pt}\text{any}\hspace{4pt} \elindexuno\in \setindexuno_2\text{,}\hspace{4pt}\unkuno\in \domsk>\argcompl{\uedgtv>\elindexuno>}>\text{.}
\end{array}
\end{flalign*}
Then, referring to Proposition \ref{genpointprop}, we have: 
\begin{flalign*}
&\left\{
\begin{array}{l}
\left(\lininclquot\left(\gczsfx[\pointedincl<\left\{\left(\mathbb{R}^{\uedgtv>\elindexuno>},0\right),\left(\mathbb{R},0\right)\right\}<>2>]\right)\compgenfquot \corrarr[\pmV]_{\elindexsei,\elindexuno},\vecuno_{\elindexuno},0\right)\hspace{4pt}\text{is}\hspace{4pt}\text{a}\hspace{4pt}\text{p-corresponding}\hspace{4pt}\text{skeleton} \\[4pt]
\text{of}\hspace{4pt}\left(\predcRel_{\elindexsei,\elindexuno}\right)_0\hspace{4pt}\text{related}\hspace{4pt}\text{to}\hspace{4pt}
\left(\lininclquot\left(\gczsfx[\pointedincl<\left\{\left(\mathbb{R}^{\uedgtv>\elindexuno>},0\right),\left(\mathbb{R},0\right)\right\}<>2>]\right)\compgenfquot \pmV_{\elindexsei,\elindexuno},\vecuno_{\elindexuno},0\right) \hspace{4pt}
\text{for}\hspace{3pt}\text{any}\\[4pt]
\elindexsei\in\left\{1,2\right\}\text{,}\hspace{4pt}\elindexuno\in\setindexuno_0\cup\setindexuno_2\text{;}
\end{array}
\right.\\[8pt]
&\left\{
\begin{array}{l}
\left(\lininclquot\left(\gczsfx[\pointedincl<\left\{\left(\mathbb{R}^{\uedgtv>\elindexuno>},0\right),\left(\mathbb{R},0\right)\right\}<>2>]\right)\compgenfquot \corrarr[\pmV]_{\elindexsei,\elindexuno},\vecuno_{\elindexuno},1\right)\hspace{4pt}\text{is}\hspace{4pt}\text{a}\hspace{4pt}\text{p-corresponding}\hspace{4pt}\text{skeleton}\\[4pt]
\text{of}\hspace{4pt}\left(\predcRel_{\elindexsei,\elindexuno}\right)_1 \hspace{4pt}\text{related}\hspace{4pt}\text{to}\hspace{4pt}
\left(\lininclquot\left(\gczsfx[\pointedincl<\left\{\left(\mathbb{R}^{\uedgtv>\elindexuno>},0\right),\left(\mathbb{R},0\right)\right\}<>2>]\right)\compgenfquot \pmV_{\elindexsei,\elindexuno},\vecuno_{\elindexuno},1\right)\hspace{4pt}\text{for}\hspace{3pt}\text{any}\\[4pt]
\elindexsei\in\left\{1,2\right\}\text{,}\hspace{4pt}\elindexuno\in\setindexuno_1\cup\setindexuno_2\text{;}
\end{array}
\right.\\[8pt]
&\left\{
\begin{array}{l}
\left(\lininclquot\left(\gczsfx[\pointedincl<\left\{\left(\mathbb{R}^{\uedgtv>\elindexuno>},0\right),\left(\mathbb{R},0\right)\right\}<>2>]\right)\compgenfquot \corrarr[\pmV]_{3,\elindexuno},\vecuno_{\elindexuno},0\right)\hspace{4pt}\text{is}\hspace{4pt}
\text{a}\hspace{4pt}\text{p-corresponding}\hspace{4pt}\text{skeleton}\\[4pt]
\text{of}\hspace{4pt}\left(\predcRel_{3,\elindexuno}\right)_0
\hspace{4pt}\text{related}\hspace{4pt}\text{to}\hspace{4pt}
\left(\lininclquot\left(\gczsfx[\pointedincl<\left\{\left(\mathbb{R}^{\uedgtv>\elindexuno>},0\right),\left(\mathbb{R},0\right)\right\}<>2>]\right)\compgenfquot \pmV_{3,\elindexuno},\vecuno_{\elindexuno},0\right)\hspace{4pt}
\text{for}\hspace{3pt}\text{any}\\[4pt] 
\elindexuno\in\setindexuno_0\text{;}
\end{array}
\right.\\[8pt]
&\left\{
\begin{array}{l}
\left(\lininclquot\left(\gczsfx[\pointedincl<\left\{\left(\mathbb{R}^{\uedgtv>\elindexuno>},0\right),\left(\mathbb{R},0\right)\right\}<>2>]\right)\compgenfquot \corrarr[\pmV]_{3,\elindexuno},\vecuno_{\elindexuno},1\right)\hspace{4pt}\text{is}\hspace{4pt}
\text{a}\hspace{4pt}\text{p-corresponding}\hspace{4pt}\text{skeleton}\\[4pt]
\text{of}\hspace{4pt}\left(\predcRel_{3,\elindexuno}\right)_1 \hspace{4pt}\text{related}\hspace{4pt}\text{to}\hspace{4pt}
\left(\lininclquot\left(\gczsfx[\pointedincl<\left\{\left(\mathbb{R}^{\uedgtv>\elindexuno>},0\right),\left(\mathbb{R},0\right)\right\}<>2>]\right)\compgenfquot \pmV_{3,\elindexuno},\vecuno_{\elindexuno},1\right)\hspace{4pt}\text{for}\hspace{3pt}\text{any}\\[4pt]
\elindexuno\in\setindexuno_1\text{;}
\end{array}
\right.\\[8pt]
&\left\{
\begin{array}{l}
\left(\lininclquot\left(\gczsfx[\pointedincl<\left\{\left(\mathbb{R}^{\uedgtv>\elindexuno>},0\right),\left(\mathbb{R},0\right)\right\}<>2>]\right)\compgenfquot \corrarr[\pmV]_{3,\elindexuno},\vecuno_{\elindexuno},1\right)\hspace{4pt}\text{is}\hspace{4pt}
\text{a}\hspace{4pt}\text{p-corresponding}\hspace{4pt}\text{skeleton}\\[4pt]
\text{of}\hspace{4pt}\left(\predcRel_{3,\elindexuno}\right)_1 \hspace{4pt}\text{related}\hspace{4pt}\text{to}\hspace{4pt}
\left(\lininclquot\left(\gczsfx[\pointedincl<\left\{\left(\mathbb{R}^{\uedgtv>\elindexuno>},0\right),\left(\mathbb{R},0\right)\right\}<>2>]\right)\compgenfquot \pmV_{3,\elindexuno},\vecuno_{\elindexuno},1\right)\hspace{4pt}
\text{for}\hspace{3pt}\text{any}\\[4pt]
\elindexuno\in\setindexuno_2\text{.}
\end{array}
\right.
\end{flalign*}
Define arrows belonging to $\homF+\GTopcat+\left( \domsk,\genfquot(\tddmsk+1)+1+\right)$ by setting:
\begin{flalign*}
&\left\{
\begin{array}{l}
\pmVdue_{\elindexsei,\elindexuno,0}\left(\unkuno\right)=\evpdf[\argcompl{\left(\lininclquot\left(\gczsfx[\pointedincl<\left\{\left(\mathbb{R}^{\uedgtv>\elindexuno>},0\right),\left(\mathbb{R},0\right)\right\}<>2>]\right)\compgenfquot \corrarr[\pmV]_{\elindexsei,\elindexuno} ,\vecuno_{\elindexuno},0\right)}](\predcRel_{\elindexsei,\elindexuno})\left(\unkuno,\lininclquot\left(\gczsfx[\idobj+\mathbb{R}+]\right)\right)\\[4pt]
\text{for}\hspace{4pt}\text{any}\hspace{4pt}\elindexsei\in\left\{1,2\right\}\text{,}\hspace{4pt}\elindexuno\in\setindexuno_0\cup\setindexuno_2\text{,}\hspace{4pt}\unkuno\in \domsk\text{;}
\end{array}
\right.\\[8pt]
&\left\{
\begin{array}{l}
\pmVdue_{\elindexsei,\elindexuno,1}\left(\unkuno\right)=\evpdf[\argcompl{\left(\lininclquot\left(\gczsfx[\pointedincl<\left\{\left(\mathbb{R}^{\uedgtv>\elindexuno>},0\right),\left(\mathbb{R},0\right)\right\}<>2>]\right)\compgenfquot \corrarr[\pmV]_{\elindexsei,\elindexuno} ,\vecuno_{\elindexuno},1\right)}](\predcRel_{\elindexsei,\elindexuno})\left(\unkuno,\lininclquot\left(\gczsfx[\idobj+\mathbb{R}+]\right)\right)\\[4pt]
\text{for}\hspace{4pt}\text{any}\hspace{4pt}\elindexsei\in\left\{1,2\right\}\text{,}\hspace{4pt}\elindexuno\in\setindexuno_1\cup\setindexuno_2\text{,}\hspace{4pt}\unkuno\in \domsk\text{;}
\end{array}
\right.\\[8pt]
&\left\{
\begin{array}{l}
\pmVdue_{3,\elindexuno}\left(\unkuno\right)=\evpdf[\argcompl{\left(\lininclquot\left(\gczsfx[\pointedincl<\left\{\left(\mathbb{R}^{\uedgtv>\elindexuno>},0\right),\left(\mathbb{R},0\right)\right\}<>2>]\right)\compgenfquot \corrarr[\pmV]_{3,\elindexuno} ,\vecuno_{\elindexuno},0\right)}](\predcRel_{3,\elindexuno})\left(\unkuno,\lininclquot\left(\gczsfx[\idobj+\mathbb{R}+]\right)\right)\\[4pt]
\text{for}\hspace{4pt}\text{any}\hspace{4pt}\elindexuno\in\setindexuno_0\text{,}\hspace{4pt}\unkuno\in \domsk\text{;}
\end{array}
\right.\\[8pt]
&\left\{
\begin{array}{l}
\pmVdue_{3,\elindexuno}\left(\unkuno\right)=\evpdf[\argcompl{\left(\lininclquot\left(\gczsfx[\pointedincl<\left\{\left(\mathbb{R}^{\uedgtv>\elindexuno>},0\right),\left(\mathbb{R},0\right)\right\}<>2>]\right)\compgenfquot \corrarr[\pmV]_{3,\elindexuno} ,\vecuno_{\elindexuno},1\right)}](\predcRel_{3,\elindexuno})\left(\unkuno,\lininclquot\left(\gczsfx[\idobj+\mathbb{R}+]\right)\right)\\[4pt]
\text{for}\hspace{4pt}\text{any}\hspace{4pt}\elindexuno\in\setindexuno_1\cup\setindexuno_2\text{,}\hspace{4pt}\unkuno\in \domsk\text{.}
\end{array}
\right.
\end{flalign*}
Eventually claim \eqref{Xèzero} follows since, by referring Definition \ref{genpointdef} and performing computations, we have that: 
\begin{flalign*}
&\begin{array}{l}
\pmVdue_{3,\elindexuno}\sumquotpoint-1\scalpquotpoint\left(\vecdue_{1,\elindexuno}\scalpquotpoint\pmVdue_{1,\elindexuno,0}\sumquotpoint \vecdue_{2,\elindexuno}\scalpquotpoint\pmVdue_{2,\elindexuno,0}\right) =\zeroquot\hspace{15pt}\forall \elindexuno\in \elindexuno_0 \text{;}
\end{array} \\[8pt]
&\begin{array}{l}
\pmVdue_{3,\elindexuno}\sumquotpoint-1\scalpquotpoint\left(\vecdue_{1,\elindexuno}\scalpquotpoint\pmVdue_{1,\elindexuno,1}\sumquotpoint \vecdue_{2,\elindexuno}\scalpquotpoint\pmVdue_{2,\elindexuno,1}\right)=\zeroquot\hspace{15pt}\forall \elindexuno\in \elindexuno_1 \text{;}  
\end{array} \\[8pt]
&\begin{array}{l}
 \pmVdue_{3,\elindexuno}\sumquotpoint-1\scalpquotpoint\left(\pmVdue_{1,\elindexuno,0}\multquotpoint \pmVdue_{2,\elindexuno,1}\sumquotpoint \pmVdue_{1,\elindexuno,1}\multquotpoint \pmVdue_{2,\elindexuno,0}\right)=\zeroquot\hspace{15pt}\forall \elindexuno\in \elindexuno_2 \text{.}
\end{array}
\end{flalign*}
\end{proof}

\section{Differential closure of $\mathbb{R}$}

Motivated by Proposition \ref{connidprop}, in this section we explicitly construct the smallest $\mathbb{R}$-algebra of $\GTopcat$ containing $\mathbb{R}$ as a subfield, such that any arrow $f$ of $\GTopcat$ to such algebra is differentiable infinitely many times in the generalized setting, with classical smooth differential if $f$ itself is classically smooth.\newline
We refer to Notations \ref{magtwopartic}, \ref{varispder}, \ref{ins}-[1], Definition \ref{defcatD}-[5], Remark \ref{noM}, \ref{specfunrem}, \ref{siterem}-[3, 6], \ref{pathinlgenerder}-[1, 6].\\[12pt]

In Definition \ref{dirsumpreclosdefdiuni} below we focus our attention on suitable admissible six-tuples. We refer to Definition \ref{dirsumpreclosdef}.

\begin{definition}\label{dirsumpreclosdefdiuni}\mbox{}
We define  
\begin{multline*}
\AQdiuni=
\Big\{\left\{\left(\ls>\antelss>, \ddgtv>\antelss>, \dgtv>\antelss>, \uedgtv>\antelss>,\ddgtv<seq<>\antelss>, \ouol<seq<>\antelss>\right)\right\}_{\antelss=1}^{\antels}\in \AQ\;:\hspace{15pt} \ddgtv>\antelss>=1\hspace{15pt}\forall \antelss\in \left\{1,...,\antels\right\}\Big\}\text{.}
\end{multline*}
Fix $\ls, \antels\in \mathbb{N}$, $\ddgtv, \dgtv,\uedgtv\in \mathbb{N}_0$. We define  $
\AQdiuni/\ls/-\antels-+\argcompl{\ddgtv,\dgtv}+<\uedgtv<= \AQ/\ls/-\antels-+\argcompl{\ddgtv,\dgtv}+<\uedgtv<$.\newline
Fix $ \antels\in \mathbb{N}$, $\ddgtv, \dgtv,\uedgtv\in \mathbb{N}_0$. We define  
$\AQdiuni-\antels-+\argcompl{\ddgtv,\dgtv}+<\uedgtv<=\AQ-\antels-+\argcompl{\ddgtv,\dgtv}+<\uedgtv<\cap\AQdiuni$.\newline
Fix $ \antels\in \mathbb{N}$, $ \dgtv,\uedgtv\in \mathbb{N}_0$. We define  
$\AQdiuni-\antels-+\argcompl{\segnvar,\dgtv}+<\uedgtv<=\AQ-\antels-+\argcompl{\segnvar,\dgtv}+<\uedgtv<\cap\AQdiuni$.\newline
Fix $ \antels\in \mathbb{N}$, $ \uedgtv\in \mathbb{N}_0$. We define  
$\AQdiuni-\antels-<\uedgtv<=\AQ-\antels-<\uedgtv<\cap\AQdiuni$.\newline
Fix $ \antels\in \mathbb{N}_0$. We define  
$\AQdiuni-\antels-= \AQ-\antels-\cap\AQdiuni$.\newline
Fix $ \uedgtv\in \mathbb{N}_0$. We define  
$\AQdiuni<\uedgtv<= \AQ<\uedgtv<\cap\AQdiuni$.\newline
Fix $\ls, \antels\in \mathbb{N}$, $\tddmsk,\uedgtv\in \mathbb{N}_0$. We define $\AQdiuni/\ls/-\antels-+\tddmsk+<\uedgtv<=\AQ/\ls/-\antels-+\tddmsk+<\uedgtv<\cap\AQdiuni$.\newline
Fix $\antels\in \mathbb{N}$, $\tddmsk,\uedgtv\in \mathbb{N}_0$. We define  
$\AQdiuni-\antels-+\tddmsk+<\uedgtv<= \AQ-\antels-+\tddmsk+<\uedgtv<\cap\AQdiuni$.\newline
Fix $\antels\in \mathbb{N}$, $\tddmsk\in \mathbb{N}_0$. We define  
$\AQdiuni-\antels-+\tddmsk+= \AQ-\antels-+\tddmsk+\cap\AQdiuni$.\newline
Fix $\tddmsk,\uedgtv\in \mathbb{N}_0$. We define  $\AQdiuni+\tddmsk+<\uedgtv<= \AQ+\tddmsk+<\uedgtv<\cap\AQdiuni$.\newline
Fix $\tddmsk\in \mathbb{N}_0$. We define  
$\AQdiuni+\tddmsk+= \AQ+\tddmsk+\cap\AQdiuni$.
\end{definition}

\begin{remark}\label{dirsumpreclosremdiuni}
Operations introduced in Definition \ref{dirsumpreclosdef}-[4, 5] induce operations on sets introduced in Definition \ref{dirsumpreclosdefdiuni}. Induced operations will be denoted  by the same symbols. 
\end{remark}

In Definition \ref{ChDiR1diuni} below we focus our attention on suitable admissible six-tuples. We refer to Definition \ref{ChDiR1}.

\begin{definition}\label{ChDiR1diuni}\mbox{}
We define: 
\begin{flalign*}
&\left\{\begin{array}{l}
\predcRdiuni[\left\{\left(\ls>\antelss>, \ddgtv>\antelss>, \dgtv>\antelss>, \uedgtv>\antelss>,\ddgtv<seq<>\antelss>, \ouol<seq<>\antelss>\right)\right\}_{\antelss=1}^{\antels}]=\overset{\antels}{\underset{\antelss=1}{\bigtensprodR}}\left(\overset{\ls>\antelss>}{\underset{\lss=1}{\bigtensprodR}} \,\generderspace-funtore-(\argcompl{\ddgtv>\antelss,\lss>})+\argcompl{\dgtv>\antelss>+\uedgtv>\antelss>}+\left(\dgtv>\antelss>\right)\right)\\[12pt]
\text{for}\hspace{4pt}\text{any}\hspace{4pt}\antelss\in \mathbb{N}\text{,}\hspace{4pt}\left\{\left(\ls>\antelss>, \ddgtv>\antelss>, \dgtv>\antelss>, \uedgtv>\antelss>,\ddgtv<seq<>\antelss>, \ouol<seq<>\antelss>\right)\right\}_{\antelss=1}^{\antels}\in\AQdiuni-\antels-\text{;}
\end{array}
\right.\\[8pt]
&\begin{array}{l}
\predcRdiuni-\antels-+\argcompl{\ddgtv,\dgtv}+<\uedgtv<=
\underset{\elAQ\in\AQdiuni-\antels-+\argcompl{\ddgtv,\dgtv}+<\uedgtv<
}{\bigdirsum}\predcRdiuni[\elAQ]\hspace{4pt}\text{for}\hspace{4pt}\text{any}\hspace{4pt}\antels\in\mathbb{N}  \text{,}\hspace{4pt}\ddgtv,\dgtv,\uedgtv\in \mathbb{N}_0\text{;}
\end{array}\\[8pt]
&\begin{array}{l}
\predcRdiuni-\antels-+\argcompl{\segnvar,\dgtv}+<\uedgtv<=
\left\{\predcRdiuni-\antels-+\argcompl{\ddgtv,\dgtv}+<\uedgtv<\right\}_{\ddgtv\in \mathbb{N}_0}\hspace{4pt}\text{for}\hspace{4pt}\text{any}\hspace{4pt}\antels\in\mathbb{N}  \text{,}\hspace{4pt}\dgtv,\uedgtv\in \mathbb{N}_0\text{;}
\end{array}\\[8pt]
&\begin{array}{l}
\predcRdiuni-\antels-+\argcompl{\segnvar,\segnvar}+<\uedgtv<=
\left\{\predcRdiuni-\antels-+\argcompl{\segnvar,\dgtv}+<\uedgtv<\right\}_{\left(\ddgtv,\dgtv\right)\in \mathbb{N}_0\times\mathbb{N}_0}\hspace{4pt}\text{for}\hspace{4pt}\text{any}\hspace{4pt}\antels\in\mathbb{N}  \text{,}\hspace{4pt}\uedgtv\in \mathbb{N}_0\text{;}
\end{array}\\[8pt]
&\begin{array}{l}
\predcRdiuni-\antels-+\argcompl{\segnvar,\segnvar}+<\segnvar<=
\left\{\predcRdiuni-\antels-+\argcompl{\ddgtv,\dgtv}+<\uedgtv<\right\}_{\left(\ddgtv,\dgtv,\uedgtv\right)\in \mathbb{N}_0\times\mathbb{N}_0\times\mathbb{N}_0}\hspace{4pt}\text{for}\hspace{4pt}\text{any}\hspace{4pt}\antels\in\mathbb{N}\text{;}
\end{array}\\[8pt]
&\begin{array}{l}
\predcRdiuni+\argcompl{\ddgtv,\dgtv}+<\uedgtv<=
\underset{\elAQ\in\AQdiuni+\argcompl{\ddgtv,\dgtv}+<\uedgtv<
}{\bigdirsum}\predcRdiuni[\elAQ]\hspace{4pt}\text{for}\hspace{4pt}\text{any}\hspace{4pt}\ddgtv,\dgtv,\uedgtv\in \mathbb{N}_0\text{;}
\end{array}\\[8pt]
&\begin{array}{l}
\predcRdiuni+\argcompl{\segnvar,\dgtv}+<\uedgtv<=
\left\{\predcRdiuni+\argcompl{\ddgtv,\dgtv}+<\uedgtv<\right\}_{\ddgtv\in \mathbb{N}_0}\hspace{4pt}\text{for}\hspace{4pt}\text{any}\hspace{4pt}\antels\in\mathbb{N}  \text{,}\hspace{4pt}\dgtv,\uedgtv\in \mathbb{N}_0\text{;}
\end{array}\\[8pt]
&\begin{array}{l}
\predcRdiuni+\argcompl{\segnvar,\segnvar}+<\uedgtv<=
\left\{\predcRdiuni+\argcompl{\ddgtv,\dgtv}+<\uedgtv<\right\}_{\left(\ddgtv,\dgtv\right)\in \mathbb{N}_0\times\mathbb{N}_0}\hspace{4pt}\text{for}\hspace{4pt}\text{any}\hspace{4pt}\antels\in\mathbb{N}  \text{,}\hspace{4pt}\uedgtv\in \mathbb{N}_0\text{;}
\end{array}\\[8pt]
&\begin{array}{l}
\predcRdiuni+\argcompl{\segnvar,\segnvar}+<\segnvar<=
\left\{\predcRdiuni+\argcompl{\segnvar,\dgtv}+<\uedgtv<\right\}_{\left(\ddgtv,\dgtv,\uedgtv\right)\in \mathbb{N}_0\times\mathbb{N}_0\times\mathbb{N}_0}\text{;}
\end{array}\\[8pt]
&\begin{array}{l}
\predcRdiuni-\antels-+\tddmsk+<\uedgtv<=
\underset{\elAQ\in\AQdiuni-\antels-+\tddmsk+<\uedgtv<
}{\bigdirsum}\predcRdiuni[\elAQ]\hspace{4pt}\text{for}\hspace{4pt}\text{any}\hspace{4pt}\antels\in\mathbb{N}  \text{,}\hspace{4pt}\tddmsk,\uedgtv\in \mathbb{N}_0\text{;}
\end{array}\\[8pt]
&\begin{array}{l}
\predcRdiuni+\tddmsk+<\uedgtv<=
\underset{\elAQ\in\AQdiuni+\tddmsk+<\uedgtv<
}{\bigdirsum}\predcRdiuni[\elAQ]\hspace{4pt}\text{for}\hspace{4pt}\text{any}\hspace{4pt}\tddmsk,\uedgtv\in \mathbb{N}_0\text{;}
\end{array}\\[8pt]
&\begin{array}{l}
\predcRdiuni+\segnvar+<\uedgtv<=\left\{\predcRdiuni+\tddmsk+<\uedgtv<\right\}_{\tddmsk\in \mathbb{N}_0}
\hspace{4pt}\text{for}\hspace{4pt}\text{any}\hspace{4pt}\uedgtv\in \mathbb{N}_0\text{;}
\end{array}
\end{flalign*}
With an abuse of language, for any $\left\{\left(\ls, \ddgtv, \dgtv, \uedgtv,\ddgtv<seq<, \ouol<seq<\right)\right\}\in\AQdiuni-1-$ we will write $\predcRdiuni[\ls, \ddgtv, \dgtv, \uedgtv,\ddgtv<seq<, \ouol<seq<]$ in place of $\predcRdiuni[\left\{\left(\ls, \ddgtv, \dgtv, \uedgtv,\ddgtv<seq<, \ouol<seq<\right)\right\}]$.
\end{definition}

\begin{remark}\label{ChDiR1remdiuni}
Operations introduced in Definitions \ref{ChDiR1}-[2, 3], \ref{raevpunfun}, \ref{laevpunfun} where we fix $\ddgtv>0>=1$ induce operations on sets introduced in Definition \ref{dirsumpreclosdefdiuni}. Induced operations will be denoted  by the same symbols. 
\end{remark}

In Definition \ref{IDeR} below we define the differential closure of $\mathbb{R}$.\newline
We refer to Notation \ref{alg}-[16].

\begin{definition} \label{IDeR}
We define
\begin{flalign*}
&\left\{
\begin{array}{l}
\text{the}\hspace{7pt}\mathbb{R}\text{-algebra}\hspace{7pt}\left(\dcR+\tddmsk+, \dcRsum, \dcRprod, \dcRscalp\right)\hspace{7pt}\text{as}\hspace{7pt}\text{the}\hspace{7pt}\text{quotient}\hspace{7pt}\text{of}\hspace{7pt}\text{the}\hspace{7pt}\mathbb{R}\text{-algebra}\\[4pt]
\left(\predcRdiuni+\tddmsk+<1<, \predcRsum, \antedcRprod, \predcRscalp\right)\hspace{4pt}\text{by}\hspace{4pt}\text{its}\hspace{4pt}\text{ideal}\hspace{4pt}\nullpredcR+\tddmsk+<1<(\argcompl{\lininclquot\left(\gczsfx[\idobj+\mathbb{R}+ ]\right)})\cap\predcRdiuni+\tddmsk+<1< \text{,}\hspace{4pt} \text{for}\hspace{4pt}\text{any}\hspace{4pt}   \tddmsk \in \mathbb{N}_0\text{;}
\end{array}
\right.\\[8pt]
&\begin{array}{l}
\predcRquotfun+\tddmsk+ =   \predcRdiuni+\tddmsk+<1< \rightarrow \dcR+\tddmsk+ \hspace{4pt}\text{the}\hspace{4pt}\text{quotient}\hspace{4pt}\text{function,}\hspace{4pt} \text{for}\hspace{4pt}\text{any}\hspace{4pt}   \tddmsk \in \mathbb{N}_0\text{;}
\end{array}\\[8pt]
&\begin{array}{l}
\left(\dcR, \dcRsum, \dcRprod, \dcRscalp\right)=\Aug \left(
\left\{\left(\dcR+\tddmsk+, \dcRsum, \dcRprod, \dcRscalp\right)\right\}_{\tddmsk\in \mathbb{N}_0}\right)\text{.}
\end{array}
\end{flalign*}
\end{definition}

In Definition \ref{evpunfun} below we describe the localization of $\dcR$-valued arrows of $\GTopcat$ at a given generalized tangent vector.\newline
We refer to Notation \ref{catfun}-[8], Definition \ref{ChDiR1}-[1], Remarks \ref{genercontfun}-[2], \ref{existoperquot}-[2], \ref{hogdorem}-[3].

\begin{definition}\label{evpunfun}
We define bilinear arrows of $\GTopcat$
\begin{equation*}
\begin{array}{l} 
\LOCAT:\homF+\GTopcat+\left(\intuno_2, \dcR\right)\times 
\generderspace-funtore-(1)+\uedgtv>2>+\left(\intuno_1\right)\rightarrow 
\homF+\GTopcat+\left(\intuno_1\times\intuno_2,\dcR\right)\\[4pt]
\uedgtv>1>,\uedgtv>2> \in \mathbb{N}_0\text{,}\hspace{4pt} \intuno_1\subseteqdentro \mathbb{R}^{\uedgtv>1>}\text{,}\hspace{4pt} \intuno_2\subseteqdentro \mathbb{R}^{\uedgtv>2>}\text{.}
\end{array}
\end{equation*}
recursively as follows.
\begin{enumerate}
\item Fix $\pmgenerder\in \generderspace-funtore-(1)+\uedgtv>2>+\left(\intuno_1\right)$, $f\in \homF+\GTopcat+\left(\intuno_2,\mathbb{R}\right)$. Define 
\begin{equation*}
\begin{array}{l}
\LOCAT-f-+\pmgenerder+\left(\unkuno_1,\unkuno_2\right)=\left(\gdgerm[
\gendiff(1)+\argcompl{\gczsfx[
\ffuncaap{\left(\pointedincl<\left\{\left(\mathbb{R}^{\uedgtv>2>},0\right),\left(\mathbb{R},0\right)\right\}<>2>\funcomp f\right)}\left(\unkuno_2\right)]}+ \funcomp\pmgenerder
]\left(\unkuno_1\right)\right)_0\\[4pt]
\hspace{230pt}\forall \left(\unkuno_{1},\unkuno_{2}\right)\in\intuno_1\times\intuno_2
\text{.}
\end{array}
\end{equation*}
\item Fix $\ddgtv, \dgtv \in \mathbb{N}_0$, $\left(\ls, \ddgtv, \dgtv, 1,\ddgtv<seq<,\ouol<seq<\right)\in \AQ-1-+\argcompl{\ddgtv,\dgtv}+<1<$, $\pmgenerder\in \generderspace-funtore-(1)+\uedgtv>2>+\left(\intuno_1\right)$, the arrow $\pmgenerderdue\in \homF+\GTopcat+\left(\intuno_2,  \predcRdiuni[\ls, \ddgtv, \dgtv, 1,\ddgtv<seq<,\ouol<seq<]\right)$.\newline 
Define the arrow $\LOCAT-\argcompl{\left(\predcRquotfun+\argcompl{\ddgtv+\dgtv}+\funcomp\pmgenerderdue\right)}-+\pmgenerder+\in \homF+\GTopcat+\left(\intuno_1\times\intuno_2, \dcR\right)$ by pointwise arguing as follows.\newline
Fix $\left(\unkuno_{1},\unkuno_{2}\right)\in\intuno_1\times\intuno_2$. By construction of $\predcRdiuni[\argcompl{\ls, \ddgtv, \dgtv, 1,\ddgtv<seq<,\ouol<seq<}]$ and of Grothen\-dieck topology of $\GTopcat$ in Section \ref{Grtopss} there are open neighborhoods $\intuno_{\unkuno_2}\sqsubseteq\intuno_2$ of $\unkuno_{2}$, $\setsymtre\sqsubseteq\mathbb{R}^{\ddgtv+\dgtv}$ of $0\in\mathbb{R}^{\ddgtv+\dgtv}$, an arrow $\pmgenerderdue_{\lss}\in\homF+\GTopcat+\left(\intuno_{\unkuno_2} ,\generderspace-funtore-(\ddgtv>\lss>)+\argcompl{\dgtv+1}+\left(\setsymtre_{\dgtv}\right)\right)$ for any $\lss\in\left\{1,...,\ls\right\}$ fulfilling condition below 
\begin{equation*}
\pmgenerderdue\left(\unkuno\right)=\overset{\ls}{\underset{\lss=1}{\bigtensprodR}} \,\gdgerm[\pmgenerderdue_{\lss}\left(\unkuno\right)]\left(0\right)\hspace{20pt}\forall\unkuno\in\intuno_{\unkuno_2}\text{.}
\end{equation*}
Define:
\begin{flalign*}
&\begin{array}{l}
\predcRel_{\unkuno_1}=\gdgerm[\gendiff(1)+\gczsfx[
\pointedincl<\left\{\left(\mathbb{R}^{\dgtv},0\right),\left(\mathbb{R}^{\uedgtv>2>},0\right)\right\}<>2>
]+  \funcomp \pmgenerder]\left(\unkuno_1\right)\text{;}
\end{array}\\[8pt]
&\left\{
\begin{array}{l}
\text{the}\hspace{4pt}\text{arrow}\hspace{4pt}\swinvadj{\pmgenerderdue}_{\lss}:\setsymtre_{\dgtv}\times\intuno_{\unkuno_2}\rightarrow\generderspace(\ddgtv>\lss>)+\argcompl{\dgtv+1}+\hspace{4pt}\text{to}\hspace{4pt}\text{be}\hspace{4pt}\text{the}\hspace{4pt}\text{adjoint}\hspace{4pt}\text{to}\hspace{4pt}
\pmgenerderdue_{\lss}\\[4pt]
\text{for}\hspace{4pt}\text{any}\hspace{4pt}\lss\in\left\{1,...,\ls\right\}\text{;}
\end{array}
\right.\\[8pt]
&\begin{array}{l}
\predcRel_{\lss,\unkuno_2}=
\gdgerm[
\gendiff(\ddgtv>\lss>)+\gczsfx[
\pointedincl<\left\{\left(\mathbb{R}^{\dgtv},0\right),\left(\mathbb{R}^{\uedgtv>2>},0\right),\left(\mathbb{R},0\right)\right\}<>\left\{1,3\right\}>
]+ \funcomp \swinvadj{\pmgenerderdue}_{\lss}] \left(0,\unkuno_2 \right)\text{;}
\end{array}\\[8pt]
&\begin{array}{l}
\LOCAT-\argcompl{\left(\predcRquotfun+\argcompl{\ddgtv+\dgtv}+\funcomp\pmgenerderdue\right)}-+\pmgenerder+\left(\unkuno_1,\unkuno_2\right)=
\predcRquotfun+\argcompl{\ddgtv+\dgtv}+\left(
\PRELOCAT-\argcompl{\left(\overset{\ls}{\underset{\lss=1}{\bigtensprodR}} \,\predcRel_{\lss,\unkuno_2}\right)} -+\predcRel_{\unkuno_1}+\right)
\text{.}
\end{array}
\end{flalign*}
\item Fix $\pmgenerder\in \generderspace-funtore-(1)+\uedgtv>2>+\left(\intuno_1\right)$, $f_1,f_2\in \homF+\GTopcat+\left(\intuno_2,\dcR\right)$. Define
\begin{multline*}
\LOCAT-\argcompl{\left(f_1 \dcRprod f_2\right)}-+\pmgenerder+\left(\unkuno_1,\unkuno_2\right)\;=\;\left(\LOCAT-\argcompl{\left(f_1\left(\unkuno_2\right)\right)}-+\gdgerm[\pmgenerder]\left(\unkuno_1\right)+\right) \; \dcRprod \;\left(\LOCAT-\argcompl{\left(f_2\left(\unkuno_2\right)\right)}-+\gdgerm[\pmgenerder]\left(\unkuno_1\right)+ \right)\\[4pt]\forall \left(\unkuno_1,\unkuno_2\right)\in \intuno_1\times\intuno_2\text{.}
\end{multline*}
\end{enumerate}
\end{definition}

In Definition \ref{connsufun} below we define the derivation of $\predcRdiuni<1<$-valued arrows of $\GTopcat$ with respect to generalized tangent vectors.\newline
We refer to Notation \ref{catfun}-[8], Definition \ref{ChDiR1}-[1], Remarks \ref{genercontfun}-[2], \ref{existoperquot}-[2], \ref{hogdorem}-[3].

\begin{definition}\label{connsufun}
We define bilinear arrows of $\GTopcat$
\begin{equation*}
\begin{array}{l} 
\DER:\generderspace-funtore-(1)+\uedgtv>2>+\left(\intuno_1\right)\times\homF+\GTopcat+\left(\intuno_2, \dcR\right) 
\rightarrow \homF+\GTopcat+\left(\intuno_1\times\intuno_2,\dcR\right)\\[4pt]
\uedgtv>1>,\uedgtv>2> \in \mathbb{N}_0\text{,}\hspace{4pt} \intuno_1\subseteqdentro \mathbb{R}^{\uedgtv>1>}\text{,}\hspace{4pt} \intuno_2\subseteqdentro \mathbb{R}^{\uedgtv>2>}\text{.}
\end{array}
\end{equation*}
recursively as follows.
\begin{enumerate}
\item Fix $\pmgenerder\in \generderspace-funtore-(1)+\uedgtv>2>+\left(\intuno_1\right)$, $f\in \homF+\GTopcat+\left(\intuno_2,\mathbb{R}\right)$. Define 
\begin{equation*}
\begin{array}{l}
\DER-f-+\pmgenerder+\left(\unkuno_1,\unkuno_2\right)=\left(\gdgerm[
\gendiff(1)+\argcompl{\gczsfx[
\ffuncaap{\left(\pointedincl<\left\{\left(\mathbb{R}^{\uedgtv>2>},0\right),\left(\mathbb{R},0\right)\right\}<>2>\funcomp f\right)}\left(\unkuno_2\right)]}+ \funcomp\pmgenerder
]\left(\unkuno_1\right)\right)_1\\[4pt]
\hspace{220pt}\forall \left(\unkuno_{1},\unkuno_{2}\right)\in\intuno_1\times\intuno_2
\text{.}
\end{array}
\end{equation*}
\item Fix $\ddgtv, \dgtv \in \mathbb{N}_0$, $\left(\ls, \ddgtv, \dgtv, 1,\ddgtv<seq<,\ouol<seq<\right)\in \AQ-1-+\argcompl{\ddgtv,\dgtv}+<1<$, $\pmgenerder\in \generderspace-funtore-(1)+\uedgtv>2>+\left(\intuno_1\right)$, the arrow $\pmgenerderdue\in \homF+\GTopcat+\left(\intuno_2,  \predcRdiuni[\ls, \ddgtv, \dgtv, 1,\ddgtv<seq<,\ouol<seq<]\right)$.\newline 
Define the arrow $\DER-\argcompl{\left(\predcRquotfun+\argcompl{\ddgtv+\dgtv}+\funcomp\pmgenerderdue\right)}-+\pmgenerder+\in \homF+\GTopcat+\left(\intuno_1\times\intuno_2, \dcR\right)$ by pointwise arguing as follows.\newline
Fix $\left(\unkuno_{1},\unkuno_{2}\right)\in\intuno_1\times\intuno_2$. By construction of $\predcRdiuni[\argcompl{\ls, \ddgtv, \dgtv, 1,\ddgtv<seq<,\ouol<seq<}]$ and of Grothend\-ieck topology of $\GTopcat$ in Section \ref{Grtopss} there are open neighborhoods $\intuno_{\unkuno_2}\sqsubseteq\intuno_2$ of $\unkuno_{2}$, $\setsymtre\sqsubseteq\mathbb{R}^{\ddgtv+\dgtv}$ of $0\in\mathbb{R}^{\ddgtv+\dgtv}$, an arrow $\pmgenerderdue_{\lss}\in\homF+\GTopcat+\left(\intuno_{\unkuno_2} ,\generderspace-funtore-(\ddgtv>\lss>)+\argcompl{\dgtv+1}+\left(\setsymtre_{\dgtv}\right)\right)$ for any $\lss\in\left\{1,...,\ls\right\}$ fulfilling condition below 
\begin{equation*}
\pmgenerderdue\left(\unkuno\right)=\overset{\ls}{\underset{\lss=1}{\bigtensprodR}} \,\gdgerm[\pmgenerderdue_{\lss}\left(\unkuno\right)]\left(0\right)\hspace{20pt}\forall\unkuno\in\intuno_{\unkuno_2}\text{.}
\end{equation*}
Define:
\begin{flalign*}
&\begin{array}{l}
\predcRel_{\unkuno_1}=\gdgerm[\gendiff(1)+\gczsfx[
\pointedincl<\left\{\left(\mathbb{R}^{\dgtv},0\right),\left(\mathbb{R}^{\uedgtv>2>},0\right)\right\}<>2>
]+  \funcomp \pmgenerder]\left(\unkuno_1\right)\text{;}
\end{array}\\[8pt]
&\left\{
\begin{array}{l}
\text{the}\hspace{4pt}\text{arrow}\hspace{4pt}\swinvadj{\pmgenerderdue}_{\lss}:\setsymtre_{\dgtv}\times\intuno_{\unkuno_2}\rightarrow\generderspace(\ddgtv>\lss>)+\argcompl{\dgtv+1}+\hspace{4pt}\text{to}\hspace{4pt}\text{be}\hspace{4pt}\text{the}\hspace{4pt}\text{adjoint}\hspace{4pt}\text{to}\hspace{4pt}
\pmgenerderdue_{\lss}\\[4pt]
\text{for}\hspace{4pt}\text{any}\hspace{4pt}\lss\in\left\{1,...,\ls\right\}\text{;}
\end{array}
\right.\\[8pt]
&\begin{array}{l}
\predcRel_{\lss,\unkuno_2}=
\gdgerm[
\gendiff(\ddgtv>\lss>)+\gczsfx[
\pointedincl<\left\{\left(\mathbb{R}^{\dgtv},0\right),\left(\mathbb{R}^{\uedgtv>2>},0\right),\left(\mathbb{R},0\right)\right\}<>\left\{1,3\right\}>
]+ \funcomp \swinvadj{\pmgenerderdue}_{\lss}] \left(0,\unkuno_2 \right)\text{;}
\end{array}\\[8pt]
&\begin{array}{l}
\DER-\argcompl{\left(\predcRquotfun+\argcompl{\ddgtv+\dgtv}+\funcomp\pmgenerderdue\right)}-+\pmgenerder+\left(\unkuno_1,\unkuno_2\right)=
\predcRquotfun+\argcompl{\ddgtv+\dgtv}+\left(
\DER-\argcompl{\left(\overset{\ls}{\underset{\lss=1}{\bigtensprodR}} \,\predcRel_{\lss,\unkuno_2}\right)} -+\predcRel_{\unkuno_1}+\right)
\text{.}
\end{array}
\end{flalign*}
\item Fix $\pmgenerder\in \generderspace-funtore-(1)+\uedgtv>2>+\left(\intuno_1\right)$, $f_1,f_2\in \homF+\GTopcat+\left(\intuno_2,\dcR\right)$. Define
\begin{equation*}
\begin{array}{l}
\left(\DER-\argcompl{f_1 \dcRprod f_2}-+\pmgenerder+\right)\left(\unkuno_1,\unkuno_2\right)=\\[4pt]
\hspace{20pt}\left(\LOCAT-\argcompl{\left(f_1\left(\unkuno_2\right)\right)}-+\gdgerm[\pmgenerder]\left(\unkuno_1\right)+\right)  \dcRprod \left(\DER-\argcompl{f_2\left(\unkuno_2\right)}-+\gdgerm[\pmgenerder]\left(\unkuno_1\right)+ \right) \dcRsum\\[4pt]
\hspace{120pt}\left(\DER-\argcompl{f_1\left(\unkuno_2\right)}-+\gdgerm[\pmgenerder]\left(\unkuno_1\right)+\right) \dcRprod \left(\LOCAT-\argcompl{\left(f_2\left(\unkuno_2\right)\right)}-+\gdgerm[\pmgenerder]\left(\unkuno_1\right)+ \right)\\[4pt]
\forall \left(\unkuno_1,\unkuno_2\right)\in \intuno_1\times\intuno_2\text{.}
\end{array}
\end{equation*}
\end{enumerate}
\end{definition}

\begin{remark}
By construction we have that the field $\mathbb{R}$ of real numbers is contained in $\dcR$ as a subfield.\newline
Definitions \ref{connid}, \ref{connsufun}, Propositions \ref{preprevalideal}-[4], \ref{connidprop} entail that any arrow $f$ of $\GTopcat$ to $\dcR$ is differentiable infinitely many times in the generalized setting, moreover if $f$ is smooth (i.e. factorizes through the natural inclusion $\mathbb{R}\rightarrow \dcR$ by a smooth function) then any derivative of $f$ with respect a to vector field in spaces introduced in Definition \ref{smsubsdef} is still smooth.
\end{remark}

\part{Applications \label{PII}}

\chapter{Generalized tangent bundle on topological manifolds\label{App}}

Classical differential geometry extends calculus from Euclidean vector spaces to smooth manifolds. This is done by gluing together derivation operators defined locally in calculus on Euclidean vector spaces. More precisely there is a functor $\smoothfib$ defined on the category of local charts of smooth manifolds and smooth functions which associates to any object the corresponding trivial bundle whose fibre is the $\mathbb{R}$-vector space of derivations of degree $1$, and to any smooth function between manifolds the corresponding differential. Algebraic structure of functor $\smoothfib$ translates in algebraic terms all relations which hold true in calculus among derivation operators.\newline We extend $\smoothfib$ to a functor $\genfib$ defined on the category local charts of topological manifolds and continuous functions, by preserving algebraic structure. We refer to Notation \ref{varispder}, \ref{realvec}-[5].\\[12pt]

In Definition \ref{settings} below we introduce categories and functors which are involved in statement and proof of Theorem \ref{mainth}. We refer to Notations \ref{catfun}, \ref{ins}-[10], \ref{gentopnot}-[3], \ref{realvec}-[3, 5, 7], \ref{difrealfunc}-[3], Definitions \ref{defcatD}-[5], \ref{alssitesdef}, Remark \ref{alssitesrem}.   

\begin{definition}\label{settings}\mbox{}
\begin{enumerate}
\item We denote by $\left(\contchart,\Gtcch\right)$ the site of continuous local charts defined as follows.\newline
Objects of $\contchart$ are pairs $\left(\intuno,m\right)$ where $m\in \mathbb{N}_0$, $\intuno \opsspint \mathbb{R}^m$.\newline
Arrows of $\contchart$ are continuous functions $f:\intuno\rightarrow \intdue$, where $\left(\intuno,m\right)$, $\left(\intdue,n\right)$ are objects of $\contchart$.\newline
Grothendieck topology $\Gtcch$ is given by selecting for any object $\left(\intuno,m\right)$ of $\contchart$ all sieves containing at least one covering of $\intuno$ by open sets.
\item We denote by $\Tmcat$ the category of topological manifolds and continuous functions.
\item We denote by $\left(\smoothchart,\Gtsch\right)$ the site of smooth local charts defined as follows.\newline
Objects of $\smoothchart$ are pairs  $\left(\intuno,m\right)$ where $m\in \mathbb{N}_0$, $\intuno \opsspint \mathbb{R}^m$.\newline
We emphasize that objects of $\smoothchart$ and objects of $\contchart$ are the same.\newline
Arrows of $\smoothchart$ are smooth functions $f:\intuno\rightarrow \intdue$ where $\left(\intuno,m\right)$, $\left(\intdue,n\right)$ are objects of $\smoothchart$.\newline
Grothendieck topology $\Gtsch$ is given by selecting for any object $\left(\intuno,m\right)$ of $\smoothchart$ all sieves containing at least one covering of $\intuno$ by open sets.
\item We denote by $\Smcat$ the category of smooth manifolds and smooth functions.
\item We denote by $\ntbf:\GTopcat\rightarrow \GTopcat$ the functor defined by setting 
\begin{flalign*}
&\ntbf\left(\intuno,m\right)=\intuno \times\left\{0\right\}\hspace{20pt}\text{for any object}\;\left(\intuno,m\right)\;\text{of}\; \GTopcat\text{,}\\[6pt]
&\ntbf\left(f\right)=f\times\idobj+\left\{0\right\}+\hspace{33pt}
\text{for any arrow}\; f:\left(\intuno,m\right)\rightarrow \left(\intdue,n\right) \;\text{of}\;  \GTopcat\text{.}
\end{flalign*}
We say that $\ntbf$ is the trivial null bundle functor.\newline
With an abuse of language we denote by the same symbol $\ntbf$ and by the same name both $\ntbf\funcomp\inclfun{\contchart}{\GTopcat}$ and $\ntbf\funcomp\inclfun{\smoothchart}{\GTopcat}$. 
\item We denote by $\smoothfib:\smoothchart\rightarrow \GTopcat$ the functor defined by setting 
\begin{flalign}
&\begin{array}{l}
\smoothfib\left(\intuno,m\right)=\intuno \times\mathbb{R}^m\hspace{20pt}\text{for any object}\;\left(\intuno,m\right)\;\text{of}\; \smoothchart\text{,}
\end{array}\\[16pt]
&\left\{
\begin{array}{l}
\left(\smoothfib\left(f\right)\right)\left(\unkuno,\vecuno\right)=\left(f\left(\unkuno\right),\symjacsm<1<
>f>[\unkuno]\left(\vecuno\right)\right)\quad \forall \unkuno \in\intuno\;\;\forall \vecuno \in\mathbb{R}^m\\[4pt]
\text{for any arrow}\; f:\left(\intuno,m\right)\rightarrow \left(\intdue,n\right) \;\text{of}\;  \smoothchart\text{.}
\end{array}
\right.\label{crossf1}
\end{flalign}
We say that $\smoothfib$ is the trivial $\mathbb{R}^m$-bundle functor.
\item We denote by $\prntsm:\smoothfib\rightarrow \inclfun{\smoothchart}{\GTopcat}$ the natural arrow defined by setting $\prntsm_{\left(\intuno,m\right)}=\proj<\intuno,\mathbb{R}^m<>1>$.
\end{enumerate}
\end{definition}

Theorem \ref{mainth} below is the main result of this work. We refer to Notations \ref{realvec}-[7], \ref{realfunc}-[6], Definitions \ref{defcatD}-[5], \ref{alssitesdef}, Proposition \ref{pointspecfun}, \ref{alssitesrem}, \ref{hogdorem}-[2],  Remarks \ref{siteremdue}, \ref{pathinlgenerder}-[6].

\begin{theorem}\label{mainth}\mbox{}
\begin{enumerate}
\item There is a functor $\genfib:\contchart \rightarrow \GTopcat$ defined as follows.\\[8pt]
Fix an object $\left(\intuno,m\right)$ of $\contchart$. We define the object $\genfib\left(\intuno,m\right)$ by setting 
\begin{equation*}
\genfib\left(\intuno,m\right)=\intuno\times \generderspace(1)+m+\text{.}
\end{equation*}
Fix an arrow $f:\left(\intuno,m\right)\rightarrow \left(\intdue,n\right)$ of $\contchart$.\newline 
We define the arrow $\genfib\left(f\right): \genfib\left(\intuno,m\right)\rightarrow \genfib\left(\intdue,n\right)$ by setting
\begin{equation*}
\left(\genfib\left(f\right)\right)\left(\unkuno,\elsymuno\right)=\left(
f\left(\unkuno\right), 
\gendiff(1)+\gggfx[\argcompl{
\funcaap{\left(f\funcomp\attfun\left[\intuno\right]\right)}\left(\unkuno\right)}]+\left(\elsymuno\right)\right)\hspace{10pt}
\forall\left(\unkuno,\elsymuno\right)\in\intuno\times\generderspace(1)+m+
\text{.}
\end{equation*}
We say that $\genfib$ is the trivial generalized vector bundle functor. 
\item There is a functor $\genfibsm:\smoothchart \rightarrow  \GTopcat$ defined as follows.\\[8pt]
Fix an object $\left(\intuno,m\right)$ of $\smoothchart$. We define the object $\genfibsm\left(\intuno,m\right)$ by setting
\begin{equation*}
\genfibsm\left(\intuno,m\right)=\intuno\times \generderspacesm(1)+m+\text{.}
\end{equation*}
Fix an arrow $f:\left(\intuno,m\right)\rightarrow \left(\intdue,n\right)$ of $\smoothchart$.\newline 
We define the arrow $\genfibsm\left(f\right): \genfibsm\left(\intuno,m\right)\rightarrow \genfibsm\left(\intdue,n\right)$ by setting
\begin{multline*}
\left(\genfibsm\left(f\right)\right)\left(\unkuno,\elsymuno\right)=\\
\left(
f\left(\unkuno\right),
\gendiffsm(1)+ 
\gggfx[\argcompl{
\linincltwo\left(\funcaap{\left(f\funcomp\attfun\left[\intuno\right]\right)}\left(\unkuno\right)\right)}]+\left(\elsymuno\right)\right)\hspace{10pt}
\forall\left(\unkuno,\elsymuno\right)\in\intuno\times\generderspacesm(1)+m+
\text{.}
\end{multline*}
We say that $\genfibsm$ is the trivial smooth generalized vector bundle functor.
\item There is a functor $\kergenfibsm:\smoothchart \rightarrow  \GTopcat$ defined as follows.\\[8pt]
Fix an object $\left(\intuno,m\right)$ of $\smoothchart$. We define the object $\kergenfibsm\left(\intuno,m\right)$ by setting
\begin{equation*}
\kergenfibsm\left(\intuno,m\right)=\intuno\times \kergenerderspacesm(1)+m+\text{.}
\end{equation*}
Fix an arrow $f:\left(\intuno,m\right)\rightarrow \left(\intdue,n\right)$ of $\smoothchart$.\newline 
We define the arrow $\kergenfibsm\left(f\right): \kergenfibsm\left(\intuno,m\right)\rightarrow \kergenfibsm\left(\intdue,n\right)$ by setting
\begin{multline*}
\left(\kergenfibsm\left(f\right)\right)\left(\unkuno,\elsymuno\right)=\\
\left(
f\left(\unkuno\right),
\kersqA(1)+\gggfx[\argcompl{
\linincltwo\left(\funcaap{\left(f\funcomp\attfun\left[\intuno\right]\right)}\left(\unkuno\right)\right)}]+\left(\elsymuno\right)\right)\hspace{10pt}
\forall\left(\unkuno,\elsymuno\right)\in\intuno\times\kergenerderspacesm(1)+m+
\text{.}
\end{multline*}
\item There is a functor $\cokergenfib:\smoothchart \rightarrow  \GTopcat$ defined as follows.\\[8pt]
Fix an object $\left(\intuno,m\right)$ of $\smoothchart$. We define the object $\cokergenbunfunc\left(\intuno,m\right)$ by setting
\begin{equation*}
\cokergenfib\left(\intuno,m\right)=\intuno\times\cokergenerderspace(1)+m+\text{.}
\end{equation*}
Fix an arrow $f:\left(\intuno,m\right)\rightarrow \left(\intdue,n\right)$ of $\smoothchart$.\newline 
We define the arrow $\cokergenfib\left(f\right): \cokergenfib\left(\intuno,m\right)\rightarrow \cokergenfib\left(\intdue,n\right)$ by setting
\begin{multline*}
\left(\cokergenfib\left(f\right)\right)\left(\unkuno,\elsymuno\right)=\\
\left(
f\left(\unkuno\right),
\cksqB(1)+ 
\gggfx[\argcompl{
\linincltwo\left(\funcaap{\left(f\funcomp\attfun\left[\intuno\right]\right)}\left(\unkuno\right)\right)}]+\left(\elsymuno\right)\right)\hspace{10pt}
\forall\left(\unkuno,\elsymuno\right)\in\intuno\times\cokergenerderspace(1)+m+
\text{.}
\end{multline*}
\item There are natural arrows
\begin{flalign}
&\left\{
\begin{array}{l}
\kerptbf : \kergenfibsm\rightarrow\genfibsm \;\text{defined by setting}\\[4pt]
\kerptbf\left(\intuno,m\right)=\idobj+\intuno+\times \kgdssmar+m+ \\[4pt]
\text{for any object}\;\left(\intuno,m\right)\;\text{of}\;\smoothchart\text{;}
\end{array} 
\right.\\[8pt]
&\left\{
\begin{array}{l}
\ptbf : \genfibsm\rightarrow \smoothfib\;\text{defined by setting}\\[4pt]
\ptbf\left(\intuno,m\right)=\idobj+\intuno+\times\quotsm+m+\\[4pt]
\text{for any object}\;\left(\intuno,m\right)\;\text{of}\;\smoothchart\text{;}
\end{array} 
\right.\\[8pt]
&\left\{
\begin{array}{l}
\itbf :\genfibsm \rightarrow\genfib \funcomp\inclfun{\smoothchart}{\contchart}\;\text{defined by setting}\\[4pt]
\itbf\left(\intuno,m\right)= \idobj+\intuno+\times \inclgdsmgd+m+\\[4pt]
\text{for any object}\;\left(\intuno,m\right)\;\text{of}\;\smoothchart\text{;}
\end{array}
\right.\\[8pt]
&\left\{
\begin{array}{l}
\ckitbf :\genfib \funcomp\inclfun{\smoothchart}{\contchart}\rightarrow \cokergenfib  \;\text{defined by setting}\\[4pt]
\ckitbf\left(\intuno,m\right)= \idobj+\intuno+\times \ckgdssmar+m+\\[4pt]
\text{for any object}\;\left(\intuno,m\right)\;\text{of}\;\smoothchart\text{;}
\end{array}
\right.\\[8pt]
&\left\{
\begin{array}{l}
\prntgf:\genfib\rightarrow \inclfun{\contchart}{\GTopcat}\;\text{defined by setting}\\[4pt]%mono
\prntgf\left(\intuno,m\right)=\proj<\intuno,\generderspace+m+<>1>\\[4pt]
\text{for any object}\;\left(\intuno,m\right)\;\text{of}\;\contchart\text{;}
\end{array}\label{prnatcch}
\right.\\[8pt]
&\left\{
\begin{array}{l}
\ckprntgf:\cokergenfib\rightarrow \inclfun{\smoothchart}{\GTopcat}\;\text{defined by setting}\\[4pt]
\prntgf\left(\intuno,m\right)=\proj<\intuno,\cokergenerderspace+m+<>1>\\[4pt]
\text{for any object}\;\left(\intuno,m\right)\;\text{of}\;\smoothchart\text{;}
\end{array}
\right.\\[8pt]
&\left\{
\begin{array}{l}
\prntgfsm:\genfibsm \rightarrow \inclfun{\smoothchart}{\GTopcat}\;\text{defined by setting}\\[4pt]%mono
\prntgfsm\left(\intuno,m\right)=\proj<\intuno,\generderspacesm+m+<>1>\\[4pt]
\text{for any object}\;\left(\intuno,m\right)\;\text{of}\;\smoothchart\text{.}
\end{array}
\right.\\[8pt]
&\left\{
\begin{array}{l}
\prkergenfibsm:\kergenfibsm \rightarrow \inclfun{\smoothchart}{\GTopcat}\;\text{defined by setting}\\[4pt]
\prkergenfibsm\left(\intuno,m\right)=\proj<\intuno,\kergenerderspacesm+m+<>1>\\[4pt]
\text{for any object}\;\left(\intuno,m\right)\;\text{of}\;\smoothchart\text{.}
\end{array}
\right.
\end{flalign} 
\item There is an exact diagram of functors
\begin{equation*}
\xymatrix{
&&0 \ar[d]\\
0 \ar[r]& \kergenfibsm\ar[r]^{\kerptbf}&\genfibsm\ar[d]^{\itbf}\ar[r]^{\ptbf}&\smoothfib\ar[r]&0\\
&&\genfib\funcomp\inclfun{\smoothchart}{\contchart}\ar[d]^{\ckitbf}\\
&&\cokergenfib\ar[d]\\
&&0}
\end{equation*}
here exactness is to be understood on fibres.
\end{enumerate} 
\end{theorem}
\begin{proof}
Statement follow straightforwardly by definition of site $\left(\GTopcat, \csGT\right)$ in Definition \ref{defcatD}-[5], Propositions \ref{nattosm}, \ref{diagRvs}, Remark \ref{pathinlgenerder}-[6].
\end{proof}

In Proposition \ref{sitevart} below we prove that any manifold defines a sub-site of $\left(\contchart,\Gtcch\right)$ or of $\left(\smoothchart,\Gtsch\right)$ depending on its smoothness.

\begin{proposition}\label{sitevart}
Fix $m \in\mathbb{N}_0$, an $m$-dimensional manifold $\left(\setsymuno, \atlas\right)$. Assume that both condition below are fulfilled: $\atlas=\left\{\setsymnove_{\elsymcinque}\overset{\cl_{\elsymcinque}}{\rightarrow}\intuno_{\elsymcinque}\;:\;\elsymcinque\in \setsymcinque\right\}$ is the maximal atlas; $\intuno_{\elsymcinque}\opsspint\mathbb{R}^m$ for any $\elsymcinque\in \setsymcinque$.\newline
Then there is a subsite $\left(\Catmn\left[\setsymuno, \atlas\right], \Gtmn\left[\setsymuno, \atlas\right]\right)$ of $\left(\contchart,\Gtcch\right)$ defined as follows.\newline
For any subset $\setsymsei\subseteq \setsymcinque$, $\elsymsei_1,\elsymsei_2\in \setsymsei$, $\intuno \opsspint \cl_{\elsymsei_1}\left(\underset{\elsymsei\in \setsymsei}{\bigcap}\setsymnove_{\elsymsei}\right)$  we have that pair $\left(\intuno,m\right)$ is an object of $\Catmn\left[\setsymuno, \atlas\right]$ and $\cl_{\elsymsei_2} \funcomp\left(\cl^{-1}_{\elsymsei_1}\funcomp\incl{\intuno}{\intuno_{\elsymsei_1}}\right) $ is an arrow of $\Catmn\left[\setsymuno, \atlas\right]$.\newline
Grothendieck topology $\Gtmn\left[\setsymuno, \atlas\right]$ is given by selecting for any object $\left(\intuno,m\right)$ of $\Catmn\left[\setsymuno, \atlas\right]$ all sieves containing at least one covering of $\intuno$ by open sets.\newline
If $\left(\setsymuno, \atlas\right)$ is a smooth manifold then $\left(\Catmn\left[\setsymuno, \atlas\right], \Gtmn\left[\setsymuno, \atlas\right]\right)$ is a subsite of $\left(\smoothchart,\Gtsch\right)$.\newline
We say that: \newline
$\Catmn\left[\setsymuno, \atlas\right]$ is the underlying category of the manifold $\left(\setsymuno, \atlas\right)$;\newline
$\Gtmn\left[\setsymuno, \atlas\right]$ is the underlying Grothendieck topology of the manifold $\left(\setsymuno, \atlas\right)$;\newline
$\left(\Catmn\left[\setsymuno, \atlas\right], \Gtmn\left[\setsymuno, \atlas\right]\right)$ is the underlying site of the manifold $\left(\setsymuno, \atlas\right)$.\newline
We drop any reference to the manifold $\left(\setsymuno, \atlas\right)$ from the notation whenever no confusion is possible.
\end{proposition}
\begin{proof}
Since the manifold $\left(\setsymuno, \atlas\right)$ is fixed then we drop any reference to it from the notation.\newline 
We prove that $\Catmn$ is a category, other parts of the statement follows straightforwardly.\newline
Fix an object $\left(\intuno,m\right)$ of $\Catmn$. We prove that the arrow $\idobj+\intuno+$ is an arrow of $\Catmn$ by arguing as follows. Choose $\elsymcinque \in \setsymcinque$ with $\intuno\opsspint\intuno_{\elsymcinque}$. Set $\cl=\cl_{\elsymcinque}\funcomp\incl{\cl_{\elsymcinque}^{-1}\left(\intuno\right)}{\setsymnove_{\elsymcinque}}$. Since $\atlas$ is a maximal atlas there is $\elsymcinque_0\in \setsymcinque$ with $\cl_{\elsymcinque_0}=\cl$. Statement follows by taking $\setsymsei=\left\{\elsymcinque_0\right\}$.\newline
Fix two arrows $\left(\intuno_1, m\right)\overset{\arrmn_1}{\rightarrow}\left(\intuno_2, m\right)$, $\left(\intuno_2, m\right)\overset{\arrmn_2}{\rightarrow}\left(\intuno_3, m\right)$ of $\Catmn$. Assume that for any $i \in \left\{1,2\right\}$ there are $\setsymsei_i\subseteq \setsymcinque$, $\elsymsei_{i,1},\elsymsei_{i,2}\in \setsymsei_i$ with:\newline
\centerline{$\intuno_i\opsspint\cl_{\elsymsei_{i,1}}\left(\underset{\elsymsei\in \setsymsei_i}{\bigcap}\setsymnove_{\elsymsei}\right)$;\hspace{20pt}$\cl_{\elsymsei_{i,2}} \funcomp\left(\cl^{-1}_{\elsymsei_{i,1}}\funcomp\incl{\intuno_i}{\intuno_{\elsymsei_{i,1}}}\right)=\arrmn_i$.}\newline
We prove that there are $\setsymsei\subseteq \setsymcinque$, $\elsymsei_{1},\elsymsei_{2}\in \setsymsei$ with:\newline
\centerline{$\intuno_1\opsspint\cl_{\elsymsei_{1}}\left(\underset{\elsymsei\in \setsymsei}{\bigcap}\setsymnove_{\elsymsei}\right)$;\hspace{20pt}$\cl_{\elsymsei_{2}} \funcomp\left(\cl^{-1}_{\elsymsei_{1}}\funcomp\incl{\intuno_1}{\intuno_{\elsymsei_{1}}}\right)=\arrmn_2\funcomp\arrmn_1$.}\newline
Set: $\cl_1=\cl_{\elsymsei_{1,1}}\funcomp\incl{\cl_{\elsymsei_{1,1}}^{-1}\left(\intuno_1\right)}{\setsymnove_{\elsymsei_{1,1}}}$; $\cl_2=\arrmn_2\funcomp\cl_{\elsymsei_{1,2}}$. Since $\atlas$ is a maximal atlas there are $\elsymcinque_1,\elsymcinque_2 \in \setsymcinque$ with $\cl_{\elsymcinque_i}=\cl_i$ for any $i \in\ \left\{1,2\right\}$. Statement follows by taking $\setsymsei=\left\{\elsymcinque_1,\elsymcinque_2\right\}$.
\end{proof}

\begin{remark}\label{costrbun}\mbox{}
\begin{enumerate}
\item We denote by $\smoothbunfunc$ the functor associating to any smooth manifold the bundle whose fibre is the topological $\mathbb{R}$-vector space of all derivation of degree $1$ (i.e. tangent vectors), and to any smooth function between smooth manifolds the corresponding differential. Functor $\smoothbunfunc$ is constructed from $\smoothfib$ as co-limit of the natural arrow $\prntsm$ introduced in Definition \ref{settings}-[7]. For any fixed smooth manifold the co-limit is computed by taking the restriction of $\prntsm$ to the underlying category of the chosen manifold (see Proposition \ref{sitevart}).\newline
We denote by $\genbunfunc$ the functor associating to any topological manifold the bundle whose fibre is the topological $\mathbb{R}$-vector space of generalized derivations of degree $1$, and to any continuous function between topological manifolds the corresponding generalized differential. Functor $\genbunfunc$ is constructed from $\genfib$ as co-limit of the natural arrow $\prntgf$ introduced in \eqref{prnatcch}. For any fixed manifold the co-limit is computed by taking the restriction of $\prntgf$ to the underlying category of the chosen manifold (see Proposition \ref{sitevart}).\newline
We denote by $\smgenbunfunc$ the functor associating to any smooth manifold the bundle whose fibre is the topological $\mathbb{R}$-vector space of all generalized smooth derivations of degree $1$, and to any smooth function between smooth manifolds the corresponding generalized differential. Functor $\smgenbunfunc$ is constructed from $\genfibsm$ as co-limit of natural arrow $\prntgfsm$.For any fixed smooth manifold the co-limit is computed by taking the restriction of $\prntgfsm$ to the underlying category of the chosen manifold (see Proposition \ref{sitevart}).\newline
We denote by $\kergenbunfuncsm$ the functor constructed from $\kergenfibsm$ as co-limit of natural arrow $\prkergenfibsm$. For any fixed smooth manifold the co-limit is computed by taking the restriction of $\prkergenfibsm$ to the underlying category of the chosen manifold (see Proposition \ref{sitevart}).\newline
We denote by $\cokergenbunfunc$ the functor constructed from $\cokergenfib$ as co-limit of natural arrow $\ckprntgf$. For any fixed smooth manifold the co-limit is computed by taking the restriction of $\ckprntgf$ to the underlying category of the chosen manifold (see Proposition \ref{sitevart}).\newline
There is an exact diagram of functors
\begin{equation*}
\xymatrix{
&&0 \ar[d]\\
0 \ar[r]& \kergenbunfuncsm\ar[r]^{\kerptbb}&\smgenbunfunc\ar[d]^{\itbb}\ar[r]^{\ptbb}&\smoothbunfunc\ar[r]&0\\
&&\genbunfunc\funcomp\inclfun{\Smcat}{\Tmcat}\ar[d]^{\ckitbb}\\
&&\cokergenbunfunc\ar[d]\\
&&0}
\end{equation*}
where here exactness is to be understood on fibres.
\item Behavior of functors $\smoothbunfunc$, $\genbunfunc$, $\smgenbunfunc$ $\kergenbunfuncsm$, $\cokergenbunfunc$ with respect to direct sum, tensor product, hom-functor of bundles can be deduced by general vector bundle theory and Proposition \ref{produno}.\newline
Arrows $\natisprodfibsm<m_1<>m_2>$, $\natisprodfibinvsm<m_1<>m_2>$ induce natural arrows of bundles which we denote by $\natisprodsm$, $\natisprodinvsm$ respectively. Arrows $\natisprodfib<m_1<>m_2>$, $\natisprodfibinv<m_1<>m_2>$ induce natural arrows of bundles which we denote by $\natisprod$, $\natisprodinv$ respectively.\newline
Commutativity of diagrams involving arrows $\ptbb$, $\itbb$, $\kerptbb$, $\ckitbb$, $\natisprod$, $\natisprodinv$, $\natisprodsm$, $\natisprodinvsm$ can be deduced by Propositions \ref{produno}, \ref{produnosm}.\newline
We emphasize that: $\natisprod$ is the left inverse of $\natisprodinv$; $\natisprodsm$ and $\natisprodinvsm$ are mutually inverse isomorphisms.\end{enumerate}
\end{remark}

\chapter{Combing hairy spheres\label{chssec}}
It is well known that nowhere vanishing tangent vector fields exist on odd dimensional spheres and do not exist on even dimensional spheres (\cite{JWM2}). This is true if we consider classical tangent vectors fields, that is sections of the classical tangent bundle. In this chapter we exhibit a nowhere vanishing generalized tangent vector field on spheres by considering a section of the generalized tangent bundle. The image of such generalized tangent vector field through $\ptbb_1$ vanishes somewhere whenever $n$ is an even integer by \cite{JWM2}.\newline
We refer to Notation \ref{realvec}-[2], Chapter \ref{App}.

\begin{definition}\label{fundefsp}
Fix $n \in \mathbb{N}$. We set:
\begin{flalign*}
& \begin{array}{l}
\opdisc\left[n\right]=\left\{\unkuno\in\mathbb{R}^n\;:\; \normeul\unkuno\normeur<1\right\}\text{;}
\end{array} \\[6pt]
&\begin{array}{l}
\circfiss\left[n\right]=\left\{\unkuno\in\mathbb{R}^n\;:\; \normeul\unkuno\normeur=\frac{2}{3}\right\}\text{;}
\end{array}\\[6pt]
&\begin{array}{l}
 \sph\left[n\right]=\left\{\left(\unkuno, \unkuno_{n+1}\right)\in\mathbb{R}^{n+1}\;:\; \normeul\left(\unkuno, \unkuno_{n+1}\right)\normeur=1 \right\} \text{;}
\end{array}\\[6pt]
& \begin{array}{l}
\upcal\left[n\right]=\left\{\left(\unkuno, \unkuno_{n+1}\right)\in\mathbb{R}^{n+1}\;:\; \normeul\left(\unkuno, \unkuno_{n+1}\right)\normeur=1 \;\text{and}\;\unkuno_{n+1}> -\frac{\sqrt{2}}{2}\right\}\text{;}
\end{array}\\[6pt]
& \begin{array}{l}
\uncal\left[n\right]=\left\{\left(\unkuno, \unkuno_{n+1}\right)\in\mathbb{R}^{n+1}\;:\; \normeul\left(\unkuno, \unkuno_{n+1}\right)\normeur=1 \;\text{and}\;\unkuno_{n+1}<\frac{\sqrt{2}}{2}\right\}\text{;}
\end{array}\\[6pt]
&\begin{array}{l}
\mappn\left[n\right]:\opdisc\left[n\right]\rightarrow \mathbb{R}^{n+1}\;\text{by setting}\\[4pt]
 \hspace{30pt}\mappn\left[n\right]\left(\unkuno\right)=\left( \frac{\sin\left(\frac{3}{4}\pi\normeul\unkuno\normeur\right)}{\normeul\unkuno\normeur}\unkuno ,\cos\left(\frac{3}{4}\pi\normeul\unkuno\normeur\right)
\right) \hspace{20pt}\forall \unkuno \in \opdisc\left[n\right]\text{;}
\end{array}\\[8pt]
&\begin{array}{l}
\mapps\left[n\right]:\opdisc\left[n\right]\rightarrow \mathbb{R}^{n+1}\;\text{by setting}\\[4pt]
 \hspace{30pt}\mapps\left[n\right]\left(\unkuno\right)=\left( \frac{\sin\left(\frac{3}{4}\pi\normeul\unkuno\normeur\right)}{\normeul\unkuno\normeur}\unkuno ,-\cos\left(\frac{3}{4}\pi\normeul\unkuno\normeur\right)
\right) \hspace{29pt}\forall \unkuno \in \opdisc\left[n\right]\text{;}
\end{array}\\[6pt]
& \begin{array}{l}
h:\mathbb{R}\rightarrow \mathbb{R}^n\hspace{4pt}\text{by}\hspace{4pt}\text{setting}\hspace{10pt}
h\left(\unkuno\right)=\left(\unkuno,0\right)\hspace{10pt}\forall \unkuno \in \mathbb{R}\text{;}
\end{array}\\[6pt]
& \begin{array}{l}
\cvup\left[n\right]: \opdisc\left[n\right]\rightarrow  \generderspacesm+n+\;\text{defined by setting}\\[4pt]
\hspace{30pt}\cvup\left[n\right]\left(\unkuno\right)= \generder[\argcompl{\genpreder\left[\lininclquot\left(\gczsfx[h]\right),\ef_1\right]}]\qquad \forall \unkuno \in \opdisc\left[n\right]\text{.}
\end{array}
\end{flalign*}
\end{definition}
     
In Remark \ref{chiarfundefsp} below we list some elementary facts about object introduced in Definition \ref{fundefsp}. We refer to Notation \ref{difrealfunc}, Definition \ref{smsubsdef}, Proposition \ref{nattosm}. 

\begin{remark}\label{chiarfundefsp}\mbox{}
\begin{enumerate}
\item Function $\mappn$ fulfills all conditions below:
\begin{flalign}
&\begin{array}{l}
\mappn\;\text{is a diffeomorphism between}\;\opdisc\left[n\right]\text{and}\;\upcal\left[n\right]\text{;}\label{omeopn}
\end{array}\\[6pt]
&\begin{array}{l}\mappn\;\text{maps}\; 0 \in \opdisc\left[n\right]\;\text{to}\;\left(0,1\right)\text{;}
\end{array}\nonumber\\[6pt]
&\begin{array}{l}
\mappn\;\text{maps}\; \circfiss\left[n\right]\;\text{to the equator of}\;\sph\left[n\right]\text{;}\nonumber
\end{array}
\end{flalign}
\item Function $\mapps$ fulfills all conditions below:
\begin{flalign}
&\begin{array}{l}
\mapps\;\text{is a diffeomorphism between}\;\opdisc\left[n\right]\text{and}\;\uncal\left[n\right]\text{;}\label{omeops}
\end{array}\\[6pt]
&\begin{array}{l}
\mapps \;\text{maps}\; 0 \in \opdisc\left[n\right]\;\text{to}\;\left(0,-1\right)\text{;}\nonumber
\end{array}\\[6pt]
&\begin{array}{l}
\mapps\;\text{maps}\; \circfiss\left[n\right]\;\text{to the equator of}\;\sph\left[n\right]\text{;}\nonumber
\end{array}
\end{flalign}
\item $\cvup\left[n\right]\left(\unkuno\right)\neq 0 \qquad \forall \unkuno \in \opdisc\left[n\right]$.
\end{enumerate}
\end{remark}

In Proposition \ref{prolcv} below we prove that the nowhere vanishing generalized tangent vector field defined by $\cvup$ through $\mappn$ to the upper emisphere can be completed to a nowhere vanishing generalized tangent vector field to the lower emisphere.   

\begin{proposition}\label{prolcv}
Fix $n \in \mathbb{N}$. Then there are $\varepsilon>0$, a continuous function 
\begin{equation*}
\cvun\left[n\right]: \opdisc\left[n\right]\rightarrow \generderspacesm+n+
\end{equation*}
fulfilling all conditions below:
\begin{flalign}
&\left\{
\begin{array}{l}
\gendiff+ \gczsfx[\argcompl{
\trasl[-\unkdue]\funcomp\left(
\mapps^{-1}\funcomp\left(\mappn\funcomp
\trasl[\unkuno]\right)\right)
 }]+\left(\cvup\left[n\right]\left(\unkuno\right)\right)=
\cvun\left[n\right]\left(\unkdue\right)\\
\text{for any}\quad \unkuno,\unkdue \in \opdisc\left[n\right]\quad\text{with}\quad \normeul \unkdue\normeur>\frac{2}{3}-\varepsilon\text{,}\quad\mappn\left(\unkuno\right)=\mapps\left(\unkdue\right)\text{;}
\end{array}
\right.\label{idsusfer}\\[8pt]
&\begin{array}{l}
\cvun\left[n\right]\left(\unkuno\right)\neq 0\hspace{15pt}\forall \unkuno \in  \opdisc\left[n\right]\text{;}\label{nonzerocond}
\end{array}\\[8pt]
&\begin{array}{l}
\quotsm+m+\left(\cvun\left(\unkuno\right)\right)=0\;\text{entails}\; \normeul\unkuno\normeur=\frac{1}{2}\text{.}\label{chsc2}
\end{array}
\end{flalign}
\end{proposition}
\begin{proof} For any fixed $a,b\in \left(0,+\infty\right)$ with $a<b$ we set:
\begin{flalign}
&\opdisc\left[n,a\right]=\left\{\vecuno\in\mathbb{R}^n\;:\; \normeul\vecuno\normeur<a\right\}\text{;}\nonumber\\[6pt]
&\sph\left[n,a\right]=\left\{\vecuno\in\mathbb{R}^{n+1}\;:\; \normeul\vecuno\normeur=a\right\}\text{;}\label{sphnot}\\[6pt]
&\cor\left[n,a,b\right]=\left\{\vecuno\in \mathbb{R}^n\;:\;a<\normeul\vecuno\normeur<b\right\}\text{.}\nonumber
\end{flalign}
Choose $\varepsilon\in\left(0,\frac{1}{6}\right)$, $\preisocordue\in \Cksp{\infty}(\mathbb{R})+\mathbb{R}+$ fulfilling all conditions below
\begin{flalign}
&\begin{array}{l}
\preisocordue\left(\unkdue\right)=\unkdue-\frac{4}{3}\hspace{30pt}\forall \unkdue \in \left(-\infty, \frac{1}{3}+\varepsilon\right)\text{;}\label{traslsucor}
\end{array}\\[6pt]
&\begin{array}{l}
\preisocordue\left(\frac{1}{2}\right)=0\text{;}\label{zerocor}
\end{array}\\[6pt]
&\begin{array}{l}
\preisocordue\left(\unkdue\right)=\unkdue\hspace{49pt}\forall \unkdue \in \left(\frac{2}{3}-\varepsilon,+\infty\right)\text{;}\label{idsucor}
\end{array}\\[6pt]
&\begin{array}{l}
\preisocordue'\left(\unkdue\right)> 0\hspace{45pt}\forall \unkdue \in\mathbb{R}\text{.}\label{derpiu}
\end{array}
\end{flalign}
Define the smooth function $\isocordue:\cor\left[n,\frac{1}{3},\frac{2}{3}\right]\rightarrow\mathbb{R}^n$ by setting 
\begin{equation*}
\isocordue\left(\vecuno\right)=\preisocordue\left(\normeul\vecuno\normeur\right)\frac{\vecuno}{\normeul\vecuno\normeur}\qquad \forall \vecuno \in \cor\left[n,\frac{1}{3},\frac{2}{3}\right]\text{.}
\end{equation*}
We have that:
\begin{flalign}
&\begin{array}{l}
\isocordue\funcomp \incl{\cor\left[n,\frac{1}{3},\frac{1}{2}\right]}{\cor\left[n,\frac{1}{3},\frac{2}{3}\right]}\;\text{is a diffeomorphism between} \\[4pt]
\cor\left[n,\frac{1}{3},\frac{1}{2}\right]\;\text{and}\; \cor\left[n,0,1\right]\text{;}
\end{array}\label{diffuno}\\[8pt]
&\begin{array}{l}
\isocordue\funcomp \incl{\cor\left[n,\frac{1}{2},\frac{2}{3}\right]}{\cor\left[n,\frac{1}{3},\frac{2}{3}\right]}\;\text{is a diffeomorphism between}\\[4pt]
\cor\left[n,\frac{1}{2},\frac{2}{3}\right]\; \text{and}\; \cor\left[n,0,\frac{2}{3}\right]\text{.}\label{diffdue}
\end{array}
\end{flalign}
Fix $\unkdue \in \cor\left[n,\frac{1}{3},1\right]$, $\delta\left[\unkdue\right] \in \left(0, \min\left\{\normeul\unkdue\normeur-\frac{1}{3}, 1-\normeul\unkdue\normeur\right\}\right)$. Define the set function $\corcv\left[\unkdue\right]:\opdisc\left[n,\delta\left[\unkdue\right]\right]\rightarrow \mathbb{R}^n$ by setting
\begin{equation*}
\corcv\left[\unkdue\right]\left(\vecuno\right) =\isocordue\left(\left(\mapps^{-1}\funcomp\mappn\right)\left(\vecuno+\left(\mappn^{-1}\funcomp\mapps\right)\left(\unkdue\right)\right)\right)-\isocordue\left(\unkdue\right)\hspace{15pt} \forall \vecuno \in \opdisc\left[n,\delta\left[\unkdue\right]\right]\text{.}
\end{equation*}
Define the set function $\cvun\left[n\right]$ by setting 
\begin{equation} 
\cvun\left(\unkdue\right)=\left\{
\begin{array}{lll}
(i)&\generder[\argcompl{\genpreder\left[\lininclquot\left(\gczsfx[\argcompl{\corcv\left[\unkdue\right]\funcomp h}]\right),\ef_1\right] }]& \forall \unkdue \in \cor\left[n,\frac{1}{3},1\right]\text{.}\\[6pt]
(ii)&\generder[\argcompl{\genpreder\left[\lininclquot\left(\gczsfx[\argcompl{-h}]\right),\ef_1\right]}]& \forall \unkdue \in \opdisc\left[n,\frac{1}{3}+\varepsilon\right]\text{.}
\end{array}
\right.\label{defcvun}
\end{equation}
By direct computation we have 
\begin{flalign*}
&\left(\mapps^{-1}\funcomp\mappn\right)\left(\unkuno\right)=\left(\frac{4}{3}-\normeul\unkuno\normeur\right)\frac{\unkuno}{\normeul\unkuno\normeur}\hspace{20pt}\forall \unkuno \in \cor\left[n,\frac{1}{3},1\right]\text{,}\\[6pt]
&\left(\mappn^{-1}\funcomp\mapps\right)\left(\unkdue\right)=\left(\frac{4}{3}-\normeul\unkdue\normeur\right)\frac{\unkdue}{\normeul\unkdue\normeur}\hspace{20pt}\forall \unkdue \in \cor\left[n,\frac{1}{3},1\right]\text{,}
\end{flalign*}
then \eqref{traslsucor} entails that $\cvun$ is well defined since \eqref{defcvun}-[(i)] and \eqref{defcvun}-[(ii)] coincide on $\cor\left[n,\frac{1}{3},\frac{1}{3}+\varepsilon\right]$.\newline
By Definitions \ref{specfunVCC}, \ref{defcatD}-[5] we have that $\cvun$ is a continuous function.\newline
By \eqref{omeopn}, \eqref{omeops}, \eqref{idsucor} we have that \eqref{idsusfer} holds true.\newline
By \eqref{omeopn}, \eqref{omeops} and since $\isocordue$ is a smooth function by construction we have that  
\begin{equation*}
\cvun\left(\unkdue\right)\in  \generderspacesm+n+ \quad \forall \unkdue \in \opdisc\left[n\right]\text{.}
\end{equation*}
By \eqref{omeopn}, \eqref{omeops}, \eqref{diffuno}, \eqref{diffdue} and since $\sph\left[n-1,\frac{1}{2}\right]$ has empty interior in $\opdisc\left[n,1\right]$ we have that \eqref{nonzerocond} holds true.\newline 
By \eqref{omeopn}, \eqref{omeops}, \eqref{diffuno}, \eqref{diffdue} we have that \eqref{chsc2} holds true.
\end{proof}

We refer to Remark \ref{costrbun}-[1, 3].

\begin{theorem}\label{chsth}
Fix $n \in \mathbb{N}$. There is at least one nowhere vanishing generalized tangent vector field $ \cv : \sph\left[n\right] \rightarrow \smgenbunfunc_1\left(\sph\left[n\right]\right)$. 
\end{theorem}
\begin{proof}
We define set function $\cv$ by setting
\begin{equation*}
\cv\left(\unkuno,\unkuno_{n+1}\right)=\left\{
\begin{array}{ll}
\left(\left(\unkuno, \unkuno_{n+1}\right),\cvup\left[n\right]\left(\mappn^{-1}\left[n\right]\left(\unkuno, \unkuno_{n+1}\right)\right)\right)&\text{if}\;\unkuno_{n+1}>\cos\left(\frac{\pi}{2}+\frac{3}{4}\pi\varepsilon\right)\text{,}\\[4pt]
\left(\left(\unkuno, \unkuno_{n+1}\right),\cvun\left[n\right]\left(\mapps^{-1}\left[n\right]\left(\unkuno, \unkuno_{n+1}\right)\right)\right)&\text{if}\;\unkuno_{n+1}<\frac{\sqrt{2} 
}{2}\text{.}\\[4pt]
\end{array}
\right.
\end{equation*} 
Proposition \ref{prolcv} entails that $\cv$ is well defined and continuous. By construction $\cv$ is a section of the bundle $\smgenbunfunc_1\left(\sph\left[n\right]\right)\rightarrow\sph\left[n\right]$. Eventually Remark \ref{chiarfundefsp}-[3], \eqref{nonzerocond} entail that $\cv$ is a nowhere vanishing vector field.
\end{proof}

In Remark \ref{nullproj} below we show how existence of nowhere vanishing generalized tangent vector fields on even dimensional spheres proved in Theorem \ref{chsth} is compatible with \cite{JWM2}. We refer to Notation \ref{difrealfunc}-[3],\eqref{sphnot}.

\begin{remark}\label{nullproj}
If $n$ is an even integer then Theorem 1 of \cite{JWM2} entails that there is at least one $\unkuno \in \sph\left[n\right]$ fulfilling 
$\left(\ptbb_1\left(\sph\left[n\right]\right)\right)\left( \cv\left(\unkuno\right)\right)=0$.\newline
If we consider the special case of the generalized tangent vector field defined in Theorem \ref{chsth} then $\left(\ptbb_1\left(\sph\left[n\right]\right)\right)\left( \cv\left(\mapps\left(\unkdue\right)\right)\right)=0$ if and only if $\unkdue\in \sph\left[n-1,\frac{1}{2}\right]$ and vector $\symjacsm<1<>[\argcompl{\mapps^{-1}\funcomp\mappn}>[\argcompl{\mappn^{-1}\left(\mapps\left(\unkdue\right)\right)}]\left(\quotsm+m+\left(\cvup\left[n\right]\left(\mappn^{-1}\left(\mapps\left(\unkdue\right)\right)\right)\right)\right)$ is tangent to $\sph\left[n-1,\frac{1}{2}\right]$ at $\unkdue$.
\end{remark}

\chapter{Questions\label{Qstns}}

In this chapter we list some questions that the theory developed in this work raises in various areas. These questions are hints for future researches.

 \section{Mathematics}

\begin{question}\label{mq1}
Do results achieved in this work hold true in a topos with a natural numbers object?\newline
We think this is true since we have used only categorical arguments in our constructions and proofs, however a careful exam is needed in order to confirm our hypothesis.   
\end{question}

\section{Physics}
\begin{question}\label{fq1}
Which is the physical meaning of the dimensions of the tangent space to the universe exceeding the classical ones?\newline
From now on we refer to these dimensions as to exotic dimensions, or exotic directions.   
\end{question}

\begin{question}\label{fq2}
Which is the real meaning of speed?\newline
Speed $V(t)$ is the differential of the function $P(t)$ representing the trajectory of the point $P$ moving in a space $S$ in an interval $T$ of time. Since $T$ is a one dimensional smooth manifold then classically we have that $V(t)$ is representable by a vector belonging to the tangent space to $S$ at $P(t)$. In the generalized setting such representation of $V(t)$ fails since the tangent space to $T$ is not a one dimensional $\mathbb{R}$-vector space. Then in the generalized setting $V(t)$ must be considered in its entirety as a continuous $\mathbb{R}$-linear function.\newline
Clearly this question is strictly related to Question \ref{fq1}.
\end{question}

\begin{question}\label{fq3}
Are there vector fields which propagates along exotic directions?\newline 
The same question also arises for tensor fields in general.   
\end{question}

\begin{question}\label{fq4}
Is it possible to build instruments which are able to measure phenomena propagating along exotic directions?    
\end{question}

\begin{question}\label{fq5}
Could some still unexplained phenomena be the manifestations of kinds of interactions occurring along exotic directions?\newline
An example could be the dark energy. We present our hypothesis about how dark energy could be hidden in the formula of kinetic energy. We refer to Question \ref{fq2} for both notation and physics considerations.\newline
Kinetic energy $E(t)=\frac{1}{2}m V(t)^2$ of a material point $P$ depends on the mass $m$ of $P$ and on speed $V(t)$ of $P$.\newline
Classically, for any instant $t$, we have that $V(t)$ is representable by a vector belonging to the tangent space to $S$ at $P(t)$ which is a finite dimensional $\mathbb{R}$-vector space.\newline
In the new generalized setting, for any instant $t$, speed $V(t)$ is a continuous $\mathbb{R}$-linear function. Then classical speed is only a single component of speed when considered in its entirety. This entails that kinetic energy computed taking in account only the classical speed turns out to be only a fraction of the total energy. We conjecture that the remaining part is the dark energy.
 
\end{question}

\begin{question}\label{fq6}
Could some measures that underlie physical theories be affected by errors due to the fact that instruments available up to now presumably do not detect phenomena occurring along exotic directions?\newline
An example could be the black body radiation which is the base of quantum mechanics.   
\end{question}

\begin{question}\label{fq7}
Could the discretization of the universe suggested by quantum mechanics be a distortion of reality due to the fact that we are observing and measuring only that part of reality which comes to us through that section of the tangent space to the universe corresponding to the classical tangent bundle?
\end{question}

\part{Appendix \label{PIII}}

\chapter{Relations in free magma $\left(\fremag, \compone\right)$\label{relations}}

In this chapter we list those relations in $\fremag$ which are used to prove results in this work.          
We refer to Notations \ref{ins}, \ref{realvec}, \ref{realfunc}, \ref{difrealfunc}, Proposition \ref{linincl},  Remark \ref{rprop}, Examples \ref{ExCinf}, \ref{ExCzero}.
\begin{flalign}
&\left\{\!
\begin{array}{l}
\lboundone \fmx_{0,1},...,\fmx_{0,n}\rboundone \!\!\relsymb\!\! \lboundone \fmx_{1,1},...,\fmx_{1,l_1},...,\fmx_{i,1},...,\fmx_{i,l_i},...,\fmx_{n,1},...,\fmx_{n,l_n}\rboundone\vspace{4pt}\\
\forall \left\{\left\{\fmx_{i,j}\right\}_{j=1}^{l_i}\right\}_{i=1}^n  \subseteq \fremag \quad\text{with}\vspace{4pt}\\
\hspace{114pt} \fmx_{0,i}=\lboundone \fmx_{i,1},...,\fmx_{i,l_i}\rboundone \quad \forall i \in \{1,...,n\}\text{;} 
\end{array}\label{(R 5)}
\right.\\[8pt]
&\left\{
\begin{array}{l}
\lboundone \fmx_1,...,\fmx_n\rboundone \relsymb \lboundone \fmy_1,...,\fmy_n\rboundone \vspace{4pt}\\
\forall \fmx_1,...,\fmx_n, \fmy_1,...,\fmy_n \in \fremag\quad\text{with}\vspace{4pt}\\
\hspace{118pt}\fmx_i\relsymb \fmy_i\quad\text{or}\quad \fmx_i=\fmy_i\quad \forall i\in\{1,...,n\}\text{;}
\end{array}\label{(R 4 bis)}
\right.\\[8pt]
&\begin{array}{l}
\left\{
\begin{array}{l}
\lboundone \fmx_1,...,\fmx_n\rboundone \compone \lboundone \fmy_1,...,\fmy_n\rboundone  \relsymb  \lboundone \fmx_1\compone  \fmy_1,..., \fmx_n\compone  \fmy_n\rboundone \vspace{4pt}\\
\forall \fmx_1,...,\fmx_n, \fmy_1,...,\fmy_n \in \fremag\quad\text{with}\vspace{4pt}\\
\hspace{118pt}\dommagone\left(\fmx_i\right)\subseteqdentro \codmagone\left(\fmy_i\right) \quad \forall i\in\left\{1,...,n\right\}\text{;}
\end{array}
\right.\label{(S 0)}
\end{array}\\[8pt]
&\left\{
\begin{array}{l}
\fmx\relsymb \lboundone\proj<\left\{\mathbb{R}^{l_i}\right\}_{i=1}^n<>1> \compone  \fmx,...,\proj<\left\{\mathbb{R}^{l_i}\right\}_{i=1}^n<>n> \compone \fmx\rboundone\compone  \Diag{\dommagone\left(\fmx\right)}{n}\\
\forall l_1,...,l_n \in \mathbb{N}_0\text{,}\quad\forall \fmx \in \fremag\quad\text{with}\quad
\codmagone\left(\fmx\right)=\underset{i=1}{\overset{n}{\prod}}\mathbb{R}^{l_i}\text{;}
\end{array}\label{(R 1)}
\right.\\[8pt]
&\left\{
\begin{array}{l}
\proj<\codmagone\left(\fmx_1\right),\codmagone\left(\fmx_2\right)<>1> \compone \lboundone \fmx_1, \fmx_2\rboundone \relsymb
\fmx_1 \compone  \proj<\dommagone\left(\fmx_1\right),\dommagone\left(\fmx_2\right)<>1>\\[4pt]
\forall \fmx_1 \in \fremag\text{,}\quad \forall \fmx_2 \in \fremagsmoothsmooth \text{;}
\end{array}\label{(R 1 bis)}
\right.\\[8pt]
&\left\{
\begin{array}{l}
\lboundone \fmx_1,...,\fmx_n\rboundone \compone \Diag{\dommagone\left(\fmx_1\right)}{n} \relsymb \Diag{\codmagone\left(\fmx_1\right)}{n} \compone  \fmx_1 \\[4pt]
\forall \left(\fmx_1,...,\fmx_n\right) \in \fremag^n \quad \text{with}\quad\fmx_1=...=\fmx_n\text{;}
\end{array}\label{(R 2)}
\right.\\[8pt]
&\left\{
\begin{array}{l}
 \lboundone \fmx_1,\fmx_2\rboundone\relsymb\\[4pt]
\hspace{48pt} \left(\switch<\codmagone\left(\fmx_1\right),\codmagone\left(\fmx_2\right)<\compone\lboundone \fmx_2,\fmx_1\rboundone\right)\compone \switch<\dommagone\left(\fmx_1\right),\dommagone\left(\fmx_2\right)<\vspace{4pt}\\
\forall \fmx_1,\fmx_2 \in\fremag\text{;}
\end{array}\label{(R 3)}
\right.\\[8pt]
&\begin{array}{l}
\left\{
\begin{array}{l}
\left( \left(\vecprod\left[1,n\right] \compone\lboundone\cost<\dommagone\left(\fmx\right)<>\mathbb{R}>+a+,\fmx\rboundone\right)\compone \Diag{\dommagone\left(\fmx\right)}{2}\right)\compone \fmy \relsymb\\[4pt]
\hspace{68pt}\left(\vecprod\left[1,n\right]\compone\lboundone\cost<\dommagone\left(\fmy\right)<>\mathbb{R}>+a+,\fmx\compone \fmy \rboundone\right)\compone \Diag{\dommagone\left(\fmy\right)}{2}\vspace{4pt}\\
\forall \fmx,\fmy \in \fremag\text{;}
\end{array}
\right.\label{(S 3)}
\end{array}\\[8pt]
& \fmx\relsymb \idobj+\codmagone\left(\fmx\right)+  \compone   \fmx \qquad \forall \fmx \in \fremag\text{;}\label{(R 7)}\\[8pt]
& \left\{
\begin{array}{l}
\fmx\relsymb \vecsum<n<>2> \compone  \lboundone \fmy, \vecprod\left[1,n\right]\compone \lboundone\cost<\dommagone\left(\fmx\right)<>\mathbb{R}>+0+, \fmz  \rboundone \rboundone\compone \Diag{\dommagone\left(\fmx\right)}{3}\\[4pt]
\text{if}\hspace{4pt}\text{and}\hspace{4pt}\text{only}\hspace{4pt}\text{if}\hspace{4pt}\fmx\relsymb \fmy\\[4pt]
\forall n \in \mathbb{N}_0\text{,}\hspace{4pt}\forall \fmx, \fmy,\fmz \in \fremag \hspace{4pt}\text{with}\hspace{4pt} \occset_{\fmx}=\occset_{\fmy}\text{,}\hspace{4pt}
 \dommagone\left(\fmx\right)=\dommagone\left(\fmy\right)=
\dommagone\left(\fmz\right)\text{,}\\[4pt]
 \codmagone\left(\fmx\right)=\codmagone\left(\fmy\right)=
\codmagone\left(\fmz\right)=\mathbb{R}^n\text{;}  
\end{array}\label{(S 7 bis)}
\right.\\[8pt]
& \fmx\relsymb \fmx  \compone \incl{\dommagone\left(\fmx\right)}{\mathbb{R}^m}  \qquad
\forall m \in \mathbb{N}_0\text{,}\quad\forall \fmx \in \fremag\quad\text{with}\quad \dommagone\left(\fmx\right)\subseteqdentro \mathbb{R}^m\text{;}\label{(R 7 ter)}\\[8pt]
& \left\{
\begin{array}{l}
\left(\fmx\compone \fmy \right) \compone   \fmz \relsymb \fmx\compone\left(\fmy\compone \fmz\right) \vspace{4pt}\\
\forall \fmx,\fmz \in \fremag\quad \forall \fmy \in \fremagsmoothsmooth \quad
\text{with}\quad
\Ima\left(\evalcompone\left(\fmy\right)\right)\subseteq \dommagone\left(\fmx\right)\text{;}
\end{array}\label{(R 7 quater)}
\right.\\[8pt]
& \left\{
\begin{array}{l}
\left(\fmx\compone \fmy \right) \compone   \fmz\relsymb  \fmx\compone\left(\fmy\compone \fmz\right)\vspace{4pt}\\
\forall \fmx,\fmy,\fmz \in \fremag \quad\text{with}\quad
\codmagone\left(\fmy\right)= \dommagone\left(\fmx\right)\text{;}
\end{array}\label{(R 7 penta)}
\right.\\[8pt]
& \fmx\relsymb\bsfuno \acmone  \fmx \qquad \forall \fmx \in \fremag \text{;}\label{(R 9)}\\[8pt]
& \bsfm \acmone (\bsfn \acmone \fmx)\relsymb(\bsfm\,\bsfn)\acmone \fmx  \qquad \forall \bsfm,\bsfn \in \bsfM\text{,}\quad \forall \fmx \in \fremag\text{;}
\label{(R 9.1)}\\[8pt]
& \bsfm \acmone \fmx \relsymb \bsfm\acmone \fmy  \qquad \forall \bsfm \in \bsfM\quad \forall \fmx,\fmy \in \fremag\quad \text{with}\quad \fmx \relsymb  \fmy\text{;}
\label{(R 9.2)}\\[8pt]
& \left\{
\begin{array}{l}
\left(\Fint_i \acmone \fmx\right) \compone  h \relsymb \vecprod\left[1,n\right]\compone \lboundone\cost<\dommagone\left(\fmx\right)<>\mathbb{R}>+0+, \fmx  \rboundone\compone \Diag{\dommagone\left(\fmx\right)}{2} \\[4pt]
\text{where}\hspace{4pt} h:\dommagone\left(\fmx\right)\rightarrow \domint\left[\dommagone\left(\fmx\right),i\right]\hspace{4pt}\text{is}\hspace{4pt}\text{defined}\hspace{4pt}\text{by}\hspace{4pt}\text{setting}\\[4pt]
h\left(\unkuno_1,...,\unkuno_m\right)=\\[4pt]
\hspace{25pt}\left(\unkuno_1,...,\unkuno_{i-1}, \unkuno_i, \unkuno_i, \unkuno_{i+1},...,\unkuno_m\right)\quad\forall\left(\unkuno_1,...,\unkuno_m\right)\in \dommagone\left(\fmx\right) \\[4pt]
\forall m\in \mathbb{N}_0\text{,}\quad \forall \fmx \in \fremag\quad\text{with}\quad i\leq m\text{,}\quad\dommagone\left(\fmx\right)\subseteqdentro\mathbb{R}^m\text{;}
\end{array}\label{(R 9.3)}
\right.\\[8pt]
& \left\{
\begin{array}{l}
\Fint_i \acmone \left(\fmx \compone \incl{\intuno}{\mathbb{R}^m}\right)\relsymb
\left(\Fint_i \acmone \fmx\right) \compone \incl{\domint\left[\intuno,i\right]}{\mathbb{R}^{m+1}}\\[4pt]
\forall m\in \mathbb{N}_0\text{,}\quad \forall \intuno\subseteqdentro \mathbb{R}^m\text{,}\quad  \forall \fmx \in \fremag\quad \text{with}\quad \intuno\subseteqdentro\dommagone\left(\fmx\right)\subseteqdentro \mathbb{R}^m\text{;}
\end{array}\label{(R 9 bis)}
\right.\\[8pt]
& \left\{
\begin{array}{l}
\Fint_i \acmone \left(\fmx \compone \pointedincl<\argcompl{\left(\intuno_1,0\right),\left(\intuno_2,0\right)}<>1> \right)\relsymb\\[4pt]
\hspace{105pt}
\left(\Fint_i \acmone \fmx\right) \compone \pointedincl<\argcompl{\left(\domint\left[\intuno_1,i\right],0\right),\left(\intuno_2,0\right)}<>1> \\[4pt]
\forall m_1,m_2\in \mathbb{N}_0\text{,}\quad \forall \intuno_1\subseteqdentro \mathbb{R}^{m_1}\text{,}\quad \forall \intuno_2\subseteqdentro \mathbb{R}^{m_2}\text{,}\quad \forall \fmx \in \fremag\quad \text{with}\\[4pt]
i\leq m_1\text{,}\quad \left(0,0\right)\in \intuno_1\times\intuno_2\subseteqdentro\dommagone\left(\fmx\right)\text{;}
\end{array}\label{(R 9 ter)}
\right.\\[8pt] 
& \left\{\! 
\begin{array}{l}
\Fint_i \acmone \fmx\relsymb\lboundone \fmy_1,...,\fmy_n\rboundone\compone  \Diag{\domint\left[\dommagone\left(\fmx\right),i\right]}{n}\quad \text{where}\\[6pt]
\fmy_{\mathsf{n}}=\left\{
\begin{array}{ll}
\Big(\funmult\left[m_{\mathsf{n}},1\right]\compone \big(\fmx_{\mathsf{n}}\compone\proj<\left\{\intuno_{i,j}\right\}_{j=1}^n<>\mathsf{n}>,\\
\funsum \left[1\right]\compone\left(\coord{L_n}{i+1},\funminus\left[1\right]\compone\coord{L_n}{i}\right) \big)\Big)\compone \\
\hspace{160pt}\Diag{\domint\left[\dommagone\left(\fmx\right),i\right]}{2}\\
\hspace{144pt}\text{if}\;i<L_{\mathsf{n}-1}\;\text{or}\;L_{\mathsf{n}}<i\leq L_n\text{,}\\[6pt]
\left(\Fint_{i-\left(L_{\mathsf{n}-1}\right)}\acmone x_{\mathsf{n}}\right)\compone 
\proj<\left\{\intuno_{i,j}\right\}_{j=1}^n<>\mathsf{n}>\\
\hspace{144pt}\text{if}\;L_{\mathsf{n}-1}< i\leq L_{\mathsf{n}}\text{,}\\[6pt]
\Big(\funmult\left[m_j,1\right]\compone \big(\fmx_{\mathsf{n}}\compone\proj<\left\{\intuno_{i,j}\right\}_{j=1}^{n+1}<>\mathsf{n}>,\\
\funsum \left[1\right]\compone\left(\coord{L_n}{i+1},\funminus\left[1\right]\compone\coord{L_n}{i}\right) \big)\Big)\compone \\
\hspace{160pt}\Diag{\domint\left[\dommagone\left(\fmx\right),i\right]}{2}\\
\hspace{144pt}\text{if}\;L_{n}<i\text{,}
\end{array}
\right.\\[6pt]
\intuno_{i,j}=\left\{
\begin{array}{ll}
\dommagone\left(\fmx_j\right)&\text{if}\;1\leq j \leq n\;\text{and}\;i<L_{j-1}\;\text{or}\;L_{j}<i\leq L_{n}\text{,}\\[4pt]
\domint\left[\dommagone\left(\fmx_j\right),i-\left(L_{j-1}\right)\right]&\text{if}\;1\leq j \leq n\;\text{and}\;L_{j-1}< i\leq L_{j}\text{,}\\[4pt]
\mathbb{R}^{i-L_n+1}&\text{if}\;j=n+1\text{,}\;L_{n}<i\text{,}
\end{array}\right.\\[6pt]
L_0=0\text{,}\hspace{4pt} L_j=\underset{k=1}{\overset{j}{\sum}}\,l_k \hspace{4pt}\forall j \in \{1,...,n\}\text{,}\\[4pt]
\forall i \in \mathbb{N}\text{,}\quad \forall l_1,...,l_n,m_1,...,m_n \in \mathbb{N}_0\text{,} \quad\forall \fmx,\fmx_1,...,\fmx_n \in \fremag \quad  \text{with}\vspace{4pt}\\
\fmx=\lboundone \fmx_1,...,\fmx_n \rboundone\text{,}\hspace{7pt} \dommagone\left(\fmx_j\right)\subseteqdentro \mathbb{R}^{l_j}\text{,}\hspace{7pt} \codmagone\left(\fmx_j\right)=\mathbb{R}^{m_j} \hspace{7pt} \forall j \in \{1,...,n\}   \text{;}
\end{array}\label{(R 10)}
\right.\\[8pt]
&\left\{  
\begin{array}{l}
\left(\Fint_i \acmone \fmx\right)\compone\left(\switch<\mathbb{R}^{m+1}<>\sigma_i^{i+1}>\compone\incl{\domint\left[\dommagone\left(\fmx\right),i\right]}{\mathbb{R}^{m+1}} \right)\relsymb\\[4pt]
\hspace{209pt} \vecminus[n]\compone \left(\Fint_i \acmone \fmx\right)\\[4pt]
\text{where}\;\sigma_i^{i+1}\;\text{is the permutation of}\;\left\{1,...,m+1\right\}\;\text{such that}\\[4pt]
\sigma_i^{i+1}\left(j\right)=\left\{
\begin{array}{ll}
j&\text{if}\;j\notin\left\{i,i+1\right\}\\
i& \text{if}\;j=i+1\\
i+1& \text{if}\;j=i\text{,}
\end{array}
\right.\\[4pt]
\forall m,i \in \mathbb{N}\text{,}\hspace{7pt} \forall n \in \mathbb{N}_0 \hspace{7pt}\text{with}\hspace{7pt} i\leq m\text{,}\hspace{7pt}
\dommagone\left(\fmx\right)\subseteqdentro \mathbb{R}^m\text{,}\hspace{7pt} \codmagone\left(\fmx\right)=\mathbb{R}^n\text{;}
\end{array}\label{(R 10.1)}
\right.\\[8pt]
& \left\{
\begin{array}{ll}
\Fint_i \acmone \left( \left(\vecsum<n<>2>\compone \lboundone \fmx_1,\fmx_2\rboundone\right)\compone\Diag{\dommagone\left(\fmx_1\right)}{2}\right)\relsymb\vspace{4pt}\\
\hspace{67pt}\left( \vecsum<n<>2>\compone \lboundone \Fint_i \acmone \fmx_1,\Fint_i \acmone \fmx_2\rboundone \right)\compone\Diag{\domint\left[\dommagone\left(\fmx_1\right),i\right]}{2}\vspace{4pt}\\
\forall n \in \mathbb{N}_0\text{,}\quad \forall i \in \mathbb{N} \text{,}\quad\forall \fmx_1,\fmx_2 \in \fremag \quad  \text{with}\vspace{4pt}\\
\dommagone\left(\fmx_1\right)=\dommagone\left(\fmx_2\right) \text{,}\;\;\codmagone\left(\fmx_1\right)=\codmagone\left(\fmx_2\right)=\mathbb{R}^n\text{;}
\end{array}\label{(R 10 bis)}
\right.\\[8pt]
& \left\{
\begin{array}{l}
\Fpart_i \acmone \fmx\relsymb\lboundone \fmy_1,...,\fmy_n\rboundone\qquad \text{where}\vspace{4pt}\\
	\fmy_j=\left\{
\begin{array}{ll}
\cost<\codmagone\left(\fmx_j\right)<>\codmagone\left(\fmx_j\right)>+0+\compone \fmx_j& \text{if}\; M_{j}<i\;\text{or}\;i\leq M_{j-1}\vspace{4pt}\\
\Fpart_{\left(i-M_{j-1}\right)}\acmone \fmx_j& \text{if}\; M_{j-1}<i\leq M_j\text{,}
\end{array}\right.\vspace{4pt}\\
M_0=0\text{,}\quad M_j=\underset{k=1}{\overset{j}{\sum}}\,m_k \quad\forall j \in \{1,...,n\}\text{,}\vspace{4pt}\\
\forall i,n \in \mathbb{N}\text{,}\quad \forall m_1,...,m_n \in \mathbb{N}_0 \text{,}\quad\forall \fmx, \fmx_1,...,\fmx_n \in \fremag \quad \text{with}\vspace{4pt}\\
\fmx=\lboundone \fmx_1,...,\fmx_n\rboundone\text{,}\quad \dommagone\left(\fmx_j\right)\subseteqdentro \mathbb{R}^{m_j}\quad \forall j \in \{1,...,n\}   \text{;}
\end{array}\label{(R 11)}
\right.\\[8pt]
& \left\{
\begin{array}{l} 
\Fpart_i \acmone \left(\fmx_1 \compone \fmx_2\right) \relsymb \left(\vecsum<n<>m>\compone\lboundone \fmy_1,...,\fmy_m\rboundone \right)\compone \Diag{\dommagone\left(\fmx_2\right)}{m}\\[4pt]
\text{where}\\[4pt]
\fmy_k = \left(\vecprod\left[n,1\right]\compone\lboundone \left(\Fpart_k \acmone \fmx_1\right)\compone \fmx_2,\left( \Fp_k\Fpart_i \right)\acmone \fmx_2 \rboundone\right)\compone \Diag{\dommagone\left(\fmx_2\right)}{2}\\
\forall l,m,n \in \mathbb{N}_0\text{,}\quad \forall i\in \mathbb{N}\text{,}\quad\forall 	\fmx_1, \fmx_2 \in \fremag \quad \text{with}\vspace{4pt}\\
 \dommagone\left(\fmx_1\right) \subseteqdentro \mathbb{R}^m\text{,}\quad\codmagone\left(\fmx_1\right) = \mathbb{R}^n\text{,}\quad \dommagone\left(\fmx_2\right)\subseteqdentro \mathbb{R}^l\text{,}\quad \codmagone\left(\fmx_2\right)= \mathbb{R}^m \text{;}
\end{array}\label{(R 12)}
\right.\\[8pt]
& \left\{
\begin{array}{l}
\Fpart_i\acmone \fmx \relsymb\left(\vecsum<n<>2>\compone \lboundone \Fpart_i\acmone \fmx, \cost<\mathbb{R}^n<>\mathbb{R}^n>+0+\compone \fmx  \rboundone
\right)\compone \Diag{\dommagone\left(\fmx\right)}{2}\\[4pt]
\forall n\in \mathbb{N}_0\text{,}\quad \forall i\in \mathbb{N}\text{,}\quad \forall \fmx \in \fremag\quad \text{with} \quad\codmagone\left(\fmx\right)=\mathbb{R}^n \text{;}  
\end{array}\label{(R 12.1)}
\right.\\[8pt]
& \left\{
\begin{array}{l}
\Fp_i \acmone \fmx\relsymb 
\left\{
\begin{array}{ll}
\fmx& \text{if}\;n=0\vspace{4pt}\\
\coord{n}{i}\compone x& \text{if}\;i\leq n\vspace{4pt}\\ 
\cost<\codmagone\left(\fmx\right)<>\mathbb{R}>+0+ \compone \fmx& \text{if}\;i>n>0
\end{array}
\right.\\[4pt]
\forall n \in \mathbb{N}_0\text{,}\quad\forall i \in \mathbb{N}
\text{,}\quad
\forall \fmx \in \fremag\quad\text{with}\quad
\codmagone\left(\fmx\right)=\mathbb{R}^n\text{;}
\end{array}\label{(R 13)}
\right.\\[8pt]
& \left\{
\begin{array}{l}
\Fq_i \acmone \fmx\relsymb 
\left\{
\begin{array}{ll}
\fmx \compone \left(\proje{m}{i+1}\compone \incl{\domint\left[\dommagone\left(\fmx\right),i\right]}{\mathbb{R}^{m+1}}\right)& \text{if}\; i \leq m\\[4pt]
\fmx \compone \proj<\dommagone\left(x\right), \mathbb{R}^{i-m+1}<>1> & \text{if}\; i > m
\end{array}
\right.\vspace{4pt}\\
\forall m \in \mathbb{N}_0\text{,}\quad\forall i \in \mathbb{N}\text{,}\quad
\forall \fmx \in \fremag\quad\text{with}\quad\dommagone\left(\fmx\right)\subseteqdentro \mathbb{R}^m\text{;}
\end{array}\label{(R 14)}
\right.\\[8pt]
& \left\{
\begin{array}{l}
\Fqq_i \acmone \fmx\relsymb\vspace{4pt}\\ 
\left\{
\begin{array}{ll}
\left(\vecminus[n]\compone \fmx\right) \compone\left( \proje{m}{i}\compone \incl{\domint\left[\dommagone\left(\fmx\right),i\right]}{\mathbb{R}^{m+1}}\right)& \text{if}\; i \leq m\vspace{4pt}\\
\left(\vecminus[n]\compone \fmx\right) \compone  \proj<\dommagone\left(\fmx\right),\mathbb{R}^{i-m+1}<>1> & \text{if}\; i > m
\end{array}
\right.\vspace{4pt}\\
\forall m \in \mathbb{N}_0\text{,}\quad\forall i \in \mathbb{N}
\text{,}\quad
\forall \fmx \in \fremag\quad\text{with}\quad
\quad\dommagone\left(\fmx\right)\subseteqdentro \mathbb{R}^m\text{;}
\end{array}\label{(R 15)}
\right.\\[8pt]
&\left\{
\begin{array}{ll}
\Fq_i \acmone \fmx\relsymb
\left(\vecsum<n<>2>\compone \lboundone\left(\Fint_i\Fpart_i\right) \acmone \fmx ,\vecminus[n]\compone\left(\Fqq_i \acmone \fmx\right)\rboundone\right)\compone\\[4pt]
\hspace{200pt} \Diag{\domint\left[\dommagone\left(\fmx\right),i\right]}{2}\\[4pt] 
\forall n \in \mathbb{N}_0\text{,}\quad\forall i \in \mathbb{N}
\text{,} \quad
 \forall \fmx\in \fremag \quad \text{with}\quad
 \codmagone\left(\fmx\right)=\mathbb{R}^n\text{;}
\end{array} \label{(R 16)}
\right.\\[8pt]
& \left\{
\begin{array}{ll} 
\Fint_i\acmone\left( \left(\vecprod\left[m,n\right]\compone\lboundone \fmx_1, \Fq_i \acmone 	\fmx_2\rboundone\right)\compone\Diag{\dommagone\left(\fmx_1\right)}{2}    \right) \relsymb \vspace{4pt}\\
\hspace{12pt}\left(\vecprod\left[m,n\right]\compone\lboundone\Fint_i \acmone \fmx_1,\Fq_i \acmone \left(\Fq_i \acmone \fmx_2\right)\rboundone\right)\compone\Diag{\domint\left[\domint\left[\dommagone\left(\fmx_2\right),i\right],i\right]}{2}\vspace{4pt}\\
\forall m,n \in \mathbb{N}_0\text{,}\quad\forall i \in \mathbb{N}
\text{,} \quad 
 \forall \fmx_1,\fmx_2 \in \fremag\quad \text{with}\vspace{4pt}\\
\dommagone\left(\fmx_1\right)=\domint\left[\dommagone\left(\fmx_2\right),i\right]\text{,}\quad \codmagone\left(\fmx_1\right)=\mathbb{R}^m\text{,}\quad
\codmagone\left(\fmx_2\right)=\mathbb{R}^n\text{;}
\end{array}\label{(R 16.1)} 
\right.\\[8pt]
& \left\{
\begin{array}{ll}
\Fint_{i+1}\acmone\left(\left( \vecprod\left[m,n\right]\compone\lboundone \fmx_1, \Fqq_i \acmone \fmx_2\rboundone\right)\compone\Diag{\dommagone\left(\fmx_1\right)}{2}    \right) \relsymb \vspace{4pt}\\
\hspace{20pt}\left(\vecprod\left[m,n\right]\compone\lboundone\Fint_{i+}\acmone \fmx_1,\Fqq_i \acmone \left( \vecminus[n]\compone \left(\Fqq_i\acmone \fmx_2\right)\right)\rboundone\right)\compone\vspace{4pt}\\
\hspace{170pt}\Diag{\domint\left[\domint\left[\dommagone\left(\fmx_2\right),i\right],i\right]}{2}\vspace{4pt}\\
\forall m,n \in \mathbb{N}_0\text{,}\quad\forall i \in \mathbb{N}
\text{,} \quad
 \forall \fmx_1,\fmx_2 \in \fremag \quad \text{with}\vspace{4pt}\\
\dommagone\left(\fmx_1\right)=\domint\left[\dommagone\left(\fmx_2\right),i\right]\text{,}\quad \codmagone\left(\fmx_1\right)=\mathbb{R}^m\text{,}\quad \codmagone\left(\fmx_2\right)=\mathbb{R}^n\text{;}
\end{array}\label{(R 16.2)}
\right.\\[8pt]
& \left\{
\begin{array}{ll}
\left(\Fq_i\Fint_i\right)\acmone \fmx \relsymb\vspace{4pt}\\
\hspace{7pt}\left(  \vecsum<n<>2>\compone\lboundone\left(\Fq_{i+1}\Fint_i\right)\acmone \fmx ,\left(\Fqq_{i+1}\Fint_i\right)\acmone \fmx\rboundone\right)\compone\Diag{\domint\left[\domint\left[\dommagone\left(\fmx\right),i\right],i\right]}{2}\vspace{4pt}\\
\forall n \in \mathbb{N}_0\text{,}\quad\forall i \in \mathbb{N}
\text{,} \quad
 \forall \fmx \in \fremag \quad \text{with}\quad \codmagone\left(\fmx\right)=\mathbb{R}^n\text{;}
\end{array}\label{(R 16.3)} 
\right.\\[8pt]
& \left\{
\begin{array}{ll}
\left(\Fqq_i\Fint_i\right)\acmone \fmx \relsymb \vecminus[n]\compone \left(\left(\Fq_{i+1} \Fint_i\right)\acmone \fmx\right)\vspace{4pt}\\
\forall n \in \mathbb{N}_0\text{,}\quad\forall i \in \mathbb{N}
\text{,} \quad
 \forall \fmx \in \fremag \quad \text{with}\quad \codmagone\left(\fmx\right)=\mathbb{R}^n\text{;}
\end{array}\label{(R 16.4)} 
\right.\\[8pt]
& \left\{
\begin{array}{ll}
\left(\Fqq_i\Fq_i\right)\acmone \fmx \relsymb \vecminus[n]\compone \left(\left(\Fq_{i+1} \Fq_i\right)\acmone \fmx\right)\vspace{4pt}\\
\forall n \in \mathbb{N}_0\text{,}\quad\forall i \in \mathbb{N}
\text{,} \quad
 \forall \fmx \in \fremag \quad \text{with}\quad \codmagone(x)=\mathbb{R}^n\text{;}
\end{array}\label{(R 16.5)}
\right.\\[8pt]
& \fmx \relsymb \lininclone\left(\evalcompone\left( \fmx\right)\right)\qquad \forall \fmx \in \fremagsmoothsmooth \text{;}\label{(R 17)}\\[8pt]
& \bsfm \acmone \fmx \relsymb \lininclone(\bsfm \smospbsfmu \fmx)\qquad \forall \bsfm \in \bsfM\quad\forall \fmx \in \Cksp{\infty}\setminus \left\{\smospempty\right\} \text{.}\label{(R 17 bis)}
\end{flalign}

\chapter{Notation \label{notat}}
We set up the notation which is used throughout the work.

\begin{notation}\label{catfun}
We establish notation relative to categories and functors, we refer to \cite{SML} for details.
\begin{enumerate}
\item We denote by $\idobj+\objuno+$ the identity arrow of and object $\objuno$ of a category. We denote by $\idobj+\Catuno+$ the identity functor of a category. We denote by $\comparrfun$ the composition of arrows in a category and of functors between categories. There is not abuse of language since categories are object and functors are arrows of the category of categories.
\item Fix a category $\Catuno$. We denote by $\Catuno^{op}$ the opposite category of $\Catuno$.
\item Fix a category $\Catuno$. We denote by $\homF+\Catuno+:\Catuno^{op}\times\Catuno\rightarrow\Setcat$ the hom-functor.
For any object $\objuno_0$ of $\Catuno$ we denote by
\begin{flalign*}
&\left\{
\begin{array}{l}
\homF+\Catuno+>\objuno_0>:\Catuno^{op}\rightarrow\Setcat\;\text{the functor defined by setting:}\\[4pt]
\hspace{5pt}\homF+\Catuno+>\objuno_0>\left(\objuno\right)=\homF+\Catuno+\left(\objuno,\objuno_0\right)\;\text{for any object}\;\objuno\;\text{of}\;\Catuno^{op}\text{,}\\[4pt]
\hspace{5pt}\homF+\Catuno+>\objuno_0>\left(\arruno\right)=\homF+\Catuno+\left(\arruno,\idobj+\objuno_0+\right)\;\text{for any arrow}\;\arruno\;\text{of}\;\Catuno^{op}\text{;}
\end{array}
\right.\\[8pt]
&\left\{
\begin{array}{l}
\homF+\Catuno+<\objuno_0<:\Catuno\rightarrow\Setcat\;\text{the functor defined by setting:}\\[4pt]
\hspace{5pt}\homF+\Catuno+<\objuno_0<\left(\objuno\right)=\homF+\Catuno+\left(\objuno_0,\objuno\right)\;\text{for any object}\;\objuno\;\text{of}\;\Catuno\text{,}\\[4pt]
\hspace{5pt}\homF+\Catuno+<\objuno_0<\left(\arruno\right)=\homF+\Catuno+\left(\idobj+\objuno_0+,\arruno\right)\;\text{for any arrow}\;\arruno\;\text{of}\;\Catuno\text{.}
\end{array}
\right.
\end{flalign*}
\item Fix two categories $\Catuno$, $\Catdue$, an object $\objdue$ of 
$\Catdue$. We denote by $\costfun<\Catuno<>\Catdue>+\objdue+$ the constant functor defined by setting $\costfun<\Catuno<>\Catdue>+\objdue+\left(\objuno\right)=\objdue$ for any object $\objuno$ of $\Catuno$ and $\costfun<\Catuno<>\Catdue>+\objdue+\left(\arruno\right)=\idobj+\objdue+$ for any arrow $\arruno$ of $\Catuno$.
\item Fix a category $\Catuno$.\newline
We denote by $\Diagcat:\Catuno\rightarrow\Catuno\times\Catuno$ the diagonal functor.\newline
We denote by $\switchcat:\Catuno\times\Catuno\rightarrow\Catuno\times\Catuno$ the switch functor. 
\item Fix two categories $\Catuno$, $\Catdue$. If $\Catuno$ is a subcategory of $\Catdue$ then we denote by $\inclfun{\Catuno}{\Catdue}$ the inclusion functor.
\item We denote an adjunction by the $4$-tuple $\left(\Catsei, \Catsette, \genfuncuno, \genfuncdue, \adjnattrasf \right)$ where: $\Catsei$, $\Catsette$ are categories; $\genfuncuno:\Catsei\rightarrow\Catsette$, $\genfuncdue:\Catsette\rightarrow\Catsei$ are functors; $\adjnattrasf$ is the natural isomorphism, explicitly $\adjnattrasf_{\objsei,\objsette}:\homF+\Catsette+\left(\genfuncuno\left(\objsei\right),\objsette\right)\rightarrow\homF+\Catsei+\left(\objsei,\genfuncdue\left(\objsette\right)\right)$ for any object $\objsei$ of $\Catsei$, $\objsette$ of $\Catsette$. 
\item We denote by $\adjphins$ the natural isomorphism realizing the adjunction between functors $\times$ and $\homF+\Setcat+$.
Explicitly we have
\begin{equation*}
\adjphins_{\setsymuno_1,\setsymuno_2,\setsymuno_3}:\homF+\Setcat+\left(\setsymuno_1\times\setsymuno_2,\setsymuno_3\right)\rightarrow\homF+\Setcat+\left(\setsymuno_1,\homF+\Setcat+\left(\setsymuno_2,\setsymuno_3\right)\right)
\end{equation*}
for any triplet of objects $\setsymuno_1$, $\setsymuno_2$, $\setsymuno_3$ of $\Setcat$.\newline
Fix objects $\setsymuno_1$ ,$\setsymuno_2$, $\setsymuno_3$ of $\Setcat$, a set function $f\in \homF+\Setcat+\left(\setsymuno_1\times\setsymuno_2,\setsymuno_3\right)$. We set:\hspace{20pt}$\diradj{f}=\adjphins_{\setsymuno_1,\setsymuno_2,\setsymuno_3}\left(f\right)$;\hspace{20pt}$\swdiradj{f}=\adjphins_{\setsymuno_2,\setsymuno_1,\setsymuno_3}\left(f\funcomp \switch<\setsymuno_2,\setsymuno_1<\right)$.\newline
Fix objects $\setsymuno_1$ ,$\setsymuno_2$, $\setsymuno_3$ of $\Setcat$, a set function $f\in \homF+\Setcat+\left(\setsymuno_1,\homF+\Setcat+\left(\setsymuno_2,\setsymuno_3\right)\right)$. We set:\hspace{20pt}$\invadj{f}=\adjphins^{-1}_{\setsymuno_1,\setsymuno_2,\setsymuno_3}\left(f\right)$;\hspace{20pt}$\swinvadj{f}=\adjphins^{-1}_{\setsymuno_1,\setsymuno_2,\setsymuno_3}\left(f\right)\funcomp \switch<\setsymuno_2,\setsymuno_1<$.
\item Fix a site $\left(\Catuno, \Grtop\right)$, an object $\objuno$ of $\Catuno$, a set $\setsymnove$ of arrows of $\Catuno$ with co-domain $\objuno$.\newline
We say that $\setsymnove$ is a base component set for $\Grtop\left(\objuno\right)$ if and only if any sieve $\sieve$ on $\objuno$ with $\setsymnove\subseteq\sieve$ belongs to $\Grtop\left(\objuno\right)$.\newline
Set 
 $\gensieve{\setsymnove}=\left\{\arruno_1\comparrfun\arruno_2\;:\; \objuno_2\overset{\arruno_1}{\rightarrow}\objuno \hspace{4pt}\text{belonging}\hspace{4pt}\text{to}\hspace{4pt}\setsymnove\text{,}\;\objuno_3\overset{\arruno_2}{\rightarrow}\objuno_2\hspace{4pt}\text{arrow}\hspace{4pt}\text{of}\hspace{4pt}\Catuno\right\}$.\newline
If $\setsymnove$ is a base component set for $\Grtop\left(\objuno\right)$ then we say that $\gensieve{\setsymnove}$ is the covering sieve of $\objuno$ generated by $\setsymnove$.
\item Fix a site $\left(\Catuno, \Grtop\right)$, an arrow  $\objuno_1\overset{\arruno}{\rightarrow}\objuno_2$ of $\Catuno$, a covering sieve $\sieve\in \Grtop\left(\objuno_2\right)$. We set denote by $\pbsieve[\arruno](\sieve)$ the covering sieve of $\objuno_1$ obtained by pulling-back $\sieve$ along $\arruno$.  
\item Fix a site $\left(\Catuno, \Grtop\right)$, two objects $\objuno_1$, $\objuno_2$ of $\Catuno$, $\arruno_1, \arruno_2 \in \homF+\Catuno+\left(\objuno_1,\objuno_2\right)$.\newline
Assume that there is covering sieve $\sieve\in \Grtop\left(\objuno_1\right)$ such that 
$\arruno_1\comparrfun \arrdue =  \arruno_2\comparrfun \arrdue$ for any $\arrdue$ belonging to $\sieve$.\newline
Then we say that $\arruno_1$ is locally equal to $\arruno_2$ and we write $\arruno_1 \eqsite[\Grtop] \arruno_2$; we drop symbol $\Grtop$ whenever no confusion is possible.
\end{enumerate}
\end{notation}

\begin{notation}\label{int}
We establish the notation relative to integer numbers.\newline 
Let $\mathbb{N}_{0}=\{0,1,2,...\}$, $\mathbb{N}=\{1,2,...\}$. For any  $n\in \mathbb{N}_0$ we set $\overline{n}=\mbox{\textnormal{max}}\{1,n\}$.
\end{notation}

\begin{notation}\label{ins}We establish the notation relative to sets and set functions.
\begin{enumerate}
\item We denote the empty set by $\udenset$ and the empty function by $\emptysetfun$.\newline
Fix two sets $\setsymuno$, $\setsymdue$. The empty function holds between them $\setsymuno$ and $\setsymdue$. If $\setsymuno=\udenset$ or $\setsymdue=\udenset$ then $\emptysetfun$ is the only one function between them, in particular $\idobj+\udenset+=\emptysetfun$.\newline
We denote the single point set by $\sielset$ and its unique element by $\sielel$ as well, with a abuse of language. 
\item Fix a finite set $\setsymuno$. We denote by $\cardset\left(\setsymuno\right)$ the cardinality of $\setsymuno$.
\item Fix a not empty set $\setindexuno$, a family of sets $\left\{\setsymuno_{\elindexuno}\right\}_{\elindexuno\in\setindexuno}$.\newline
If $\setsymuno_{\elindexuno}=\sielset$ for any $\elindexuno$ belonging to a subset $\setindexuno_0$ of $\setindexuno$ then $\underset{\elindexuno\in\setindexuno}{\prod}\setsymuno_{\elindexuno}=\underset{\elindexuno\in\setindexuno\setminus\setindexuno_0}{\prod}\setsymuno_{\elindexuno}$.\newline
If $\setsymuno_{\elindexuno}=\udenset$ for any $\elindexuno$ belonging to a not empty subset $\setindexuno_0$ of $\setindexuno$ then
 $\underset{\elindexuno\in\setindexuno}{\prod}\setsymuno_{\elindexuno}=\udenset$.
\item Fix two sets $\setsymcinque$, $\setsymuno$. We set 
\begin{equation*}
\setsymuno^{\setsymcinque} = \left\{
\begin{array}{ll}
\text{the single point set}& \text{if}\; \setsymcinque=\udenset\text{,}\\[4pt]
\underset{\elsymcinque\in\setsymcinque}{\prod}\setsymuno&\text{if}\;  \setsymcinque\neq\udenset\text{.} 
\end{array}
\right. 
\end{equation*} 
If $\setsymcinque$ is a finite set then we write $\setsymuno^{\cardset\left(\setsymcinque\right)}$ in place of $\setsymuno^{\setsymcinque}$.
\item Fix a set function $f:\setsymuno \rightarrow \setsymdue$ not necessarily everywhere defined, a set $\setsymtre\subseteq \setsymuno$, a set $\setsymquattro \subseteq \setsymdue$, $\elsymuno_0\in\setsymuno$, $\elsymdue_0\in\setsymdue$.\newline
We write $f\left(\elsymuno_0\right)\neq\udenunk$ to mean that $f$ is well defined for $\elsymuno$, we write $f\left(\elsymuno_0\right)=\udenunk$ to mean that $f$ is not defined for $\elsymuno_0$,  we write $f\left(\elsymuno_0\right)\neq\elsymdue_0$ to mean that $f\left(\elsymuno_0\right)\neq\udenunk$ and $f\left(\elsymuno_0\right)\in \setsymdue\setminus\left\{\elsymdue_0\right\}$.\newline
We adopt the following notation:
\begin{flalign*} 
&\Dom\left(f\right)=\setsymuno\text{,}\hspace{113pt}\Cod\left(f\right)=\setsymdue\text{,}\\[6pt]
&\Cae\left(f\right)=\left\{\elsymuno\in \setsymuno\;:\;f\left(\elsymuno\right)\neq\udenunk  \right\}\text{,}\hspace{24pt}\Ima\left(f\right)=\left\{f\left(\elsymuno\right),\; \elsymuno\in \setsymuno \right\}\text{,}\\[4pt]
&f^{-1}\left(\setsymquattro\right)=\left\{\elsymuno\in \setsymuno:\; f\left(\elsymuno\right)\in \setsymquattro \right\}\text{,}\hspace{20pt} f\left(\setsymtre\right)=\left\{f\left(\elsymtre\right),\; \elsymtre\in \setsymtre\cap \Cae\left(f\right) \right\}\text{.}
\end{flalign*}
If $f=\emptysetfun$ then $\Dom\left(f\right)$ and $\Cod\left(f\right)$ are not well defined since the empty function holds between any pair of sets.    
\item Fix a set $\setsymcinque$, a family $\left\{f_{\elsymcinque}\right\}_{\elsymcinque\in\setsymcinque}$ of set functions with $f_{\elsymcinque}:\setsymuno_{\elsymcinque} \rightarrow \setsymdue_{\elsymcinque}$.
We denote by $\underset{\elsymcinque\in \setsymcinque}{\prod}f_{\elsymcinque} : \underset{\elsymcinque\in \setsymcinque}{\prod}\setsymuno_{\elsymcinque} \rightarrow \underset{\elsymcinque\in \setsymcinque}{\prod}\setsymdue_{\elsymcinque}$ the set function defined by setting 
\begin{equation*}
\underset{\elsymcinque\in \setsymcinque}{\prod}f_{\elsymcinque}\left( \underset{\elsymcinque\in \setsymcinque}{\prod}\elsymuno_{\elsymcinque}\right)=\underset{\elsymcinque\in \setsymcinque}{\prod}f_{\elsymcinque}\left(\elsymuno_{\elsymcinque}\right)
 \quad \forall\underset{\elsymcinque\in \setsymcinque}{\prod}\elsymuno_{\elsymcinque} \in\underset{\elsymcinque\in \setsymcinque}{\prod}\setsymuno_{\elsymcinque}\text{.}
\end{equation*} 
If there is a set function $f:\setsymuno\rightarrow \setsymdue$ with $\setsymuno=\setsymuno_{\elsymcinque}$, $\setsymdue=\setsymdue_{\elsymcinque}$, $f_{\elsymcinque}=f$ for any $\elsymcinque\in \setsymcinque$ then we denote $\underset{\elsymcinque\in \setsymcinque}{\prod}f_{\elsymcinque}$ by $f^{\setsymcinque}$ if $\setsymcinque$ is an infinite set, or by $f^{\cardset\left(\setsymcinque\right)}$ if $\setsymcinque$ is a finite set.
\item Fix two sets $\setsymuno$, $\setsymdue$ with $\setsymdue\neq\udenset$, an element $\elsymdue \in \setsymuno$. We denote by $\cost<\setsymuno<>\setsymdue>+\elsymdue+:\setsymuno \rightarrow \setsymdue$ the constant function $\cost<\setsymuno<>\setsymdue>+\elsymdue+\left(\elsymuno\right)=\elsymdue$ for any $\elsymuno\in \setsymuno$. If $\setsymuno=\udenset$ then $\cost<\setsymuno<>\setsymdue>+\elsymdue+=\emptysetfun$.
\item Fix two sets $\setsymuno$, $\setsymdue$. We denote by $\Diag{\setsymuno}{\setsymdue}:\setsymuno\rightarrow \setsymuno^{\setsymdue}$ the diagonal function. If $\setsymdue$ is a finite set then we replace $\setsymdue$ by the integer $\cardset\left(\setsymdue\right)$ in the notation. If $\cardset\left(\setsymdue\right)=1$ then $\Diag{\setsymuno}{1}=\idobj+\setsymuno+$. If $\setsymuno=\udenset$ then $\Diag{\setsymuno}{\setsymdue}=\emptysetfun$. 
\item Fix two sets $\setsymuno$, $\setsymdue$ with $\setsymuno \subseteq \setsymdue$. We denote by $\incl{\setsymuno}{\setsymdue}:\setsymuno \rightarrow \setsymdue$ the inclusion function. If $\setsymuno=\setsymdue$ then $\incl{\setsymuno}{\setsymdue}= \idobj+\setsymuno+$. If $\setsymuno=\udenset$ or $\setsymdue=\udenset$ then $\incl{\setsymuno}{\setsymdue}=\emptysetfun$.
\item Fix a not empty set $\setindexuno$, $\setindexuno_0\subseteq\setindexuno$, a family of sets $\left\{\setsymuno_{\elindexuno}\right\}_{\elindexuno\in\setindexuno}$.\newline
If $\setindexuno_0=\udenset$ then we set $\proj<\left\{\setsymuno_{\elindexuno}\right\}_{\elindexuno\in\setindexuno}<>\setindexuno_0>=\cost<\underset{\elindexuno\in\setindexuno}{\prod}\setsymuno_{\elindexuno}<>\sielset>+\sielel+$.\newline
If $\setindexuno_0\neq\udenset$ then we denote by $\proj<\left\{\setsymuno_{\elindexuno}\right\}_{\elindexuno\in\setindexuno}<>\setindexuno_0>:\underset{\elindexuno\in\setindexuno}{\prod}\setsymuno_{\elindexuno}\rightarrow \underset{\elindexuno\in\setindexuno_0}{\prod}\setsymuno_{\elindexuno}$ the function defined by setting
\begin{equation*}
\proj<\left\{\setsymuno_{\elindexuno}\right\}_{\elindexuno\in\setindexuno}<>\setindexuno_0>\left(\left(\elsymuno_{\elindexuno}\right)_{\elindexuno\in\setindexuno}\right)=\left(\elsymuno_{\elindexuno}\right)_{\elindexuno\in\setindexuno_0}\hspace{15pt}\forall \left(\elsymuno_{\elindexuno}\right)_{\elindexuno\in\setindexuno}\in \underset{\elindexuno\in\setindexuno}{\prod}\setsymuno_{\elindexuno}\text{.}
\end{equation*} 
If $\setindexuno=\setindexuno_0$ then $\proj<\left\{\setsymuno_{\elindexuno}\right\}_{\elindexuno\in\setindexuno}<>\setindexuno_0>= \idobj+\argcompl{\underset{\elindexuno\in\setindexuno}{\prod}\setsymuno_{\elindexuno}}+$.\newline
If there is $\elindexuno_0\in\setindexuno$ with $\setindexuno_0=\left\{\elindexuno_0\right\}$ then we write $\proj<\left\{\setsymuno_{\elindexuno}\right\}_{\elindexuno\in\setindexuno}<>\elindexuno_0>$ in place of $\proj<\left\{\setsymuno_{\elindexuno}\right\}_{\elindexuno\in\setindexuno}<>\argcompl{\left\{\elindexuno_0\right\}}>$, with an abuse of language.\newline
If there is $\elindexuno_1\in\setindexuno$ with $\setsymuno_{\elindexuno_1}=\udenset$ then by Notation \ref{ins}-[3] we have $\proj<\left\{\setsymuno_{\elindexuno}\right\}_{\elindexuno\in\setindexuno}<>\setindexuno_0>=\emptysetfun$ for any $\setindexuno_0\subseteq\setindexuno$.\newline
With an abuse of language we write $\underset{\elindexuno\in\setindexuno}{\prod}\setsymuno_{\elindexuno}$ in place of $\left\{\setsymuno_{\elindexuno}\right\}_{\elindexuno\in\setindexuno}$ or we drop any reference to $\left\{\setsymuno_{\elindexuno}\right\}_{\elindexuno\in\setindexuno}$ or to $\setindexuno_0$ whenever no confusion is possible.
\item Fix a not empty set $\setindexuno$, $\setindexuno_0\subseteq\setindexuno$, a family of pairs $\left\{\left(\setsymuno_{\elindexuno}, \elsymuno_{0,\elindexuno}\right)\right\}_{\elindexuno\in\setindexuno}$ where $\setsymuno_{\elindexuno}$ is a set and $\elsymuno_{\elindexuno}\in \setindexuno_{\elindexuno}$ for any $\elindexuno \in \setindexuno$.\newline
If $\setindexuno_0=\udenset$ then we set $\pointedincl<\left\{\left(\setsymuno_{\elindexuno}, \elsymuno_{0,\elindexuno}\right)\right\}_{\elindexuno\in\setindexuno}<>\setindexuno_0>=\cost<\sielset<>\underset{\elindexuno\in\setindexuno}{\prod}\setsymuno_{\elindexuno}>+\left(\elsymuno_{0,\elindexuno}\right)_{\elindexuno\in\setindexuno}+$.\newline
If $\setindexuno_0\neq\udenset$ then we denote by $\pointedincl<\left\{\left(\setsymuno_{\elindexuno}, \elsymuno_{\elindexuno}\right)\right\}_{\elindexuno\in\setindexuno}<>\setindexuno_0>:\underset{\elindexuno\in \setindexuno_0}{\prod}  \setsymuno_{\elindexuno}\rightarrow \underset{\elindexuno\in \setindexuno}{\prod}  \setsymuno_{\elindexuno}$ the set function defined by setting
\begin{multline*}
\pointedincl<\left\{\left(\setsymuno_{\elindexuno}, \elsymuno_{\elindexuno}\right)\right\}_{\elindexuno\in\setindexuno}<>\setindexuno_0>\left(\left(\elsymuno_{\elindexuno}\right)_{\elindexuno\in\setindexuno}\right)=\left(\elsymdue_{\elindexuno}\right)_{\elindexuno\in\setindexuno}\\
 \text{where}\hspace{4pt}
\elsymdue_{\elindexuno}=\left\{
\begin{array}{ll}
\elsymuno_{\elindexuno}&\text{if}\hspace{4pt}\elindexuno\in\setindexuno_0  \text{,}\\[4pt]
\elsymuno_{0,\elindexuno}&\text{if}\hspace{4pt}\elindexuno\notin\setindexuno_0  \\[4pt]
\end{array}
\right.
\hspace{20pt}\forall \elsymuno\in \setsymuno_{\elindexuno_0}\text{.}
\end{multline*}
If $\setindexuno=\setindexuno_0$ then $\pointedincl<\left\{\left(\setsymuno_{\elindexuno}, \elsymuno_{\elindexuno}\right)\right\}_{\elindexuno\in\setindexuno}<>\setindexuno_0>= \idobj+\argcompl{\underset{\elindexuno\in\setindexuno}{\prod}\setsymuno_{\elindexuno}}+$.\newline
If there is $\elindexuno_0\in\setindexuno$ with $\setindexuno_0=\left\{\elindexuno_0\right\}$ then we write $\pointedincl<\left\{\left(\setsymuno_{\elindexuno}, \elsymuno_{\elindexuno}\right)\right\}_{\elindexuno\in\setindexuno}<>\elindexuno_0>$ in place of $\pointedincl<\left\{\left(\setsymuno_{\elindexuno}, \elsymuno_{\elindexuno}\right)\right\}_{\elindexuno\in\setindexuno}<>\left\{\elindexuno_0\right\}>$, with an abuse of language.\newline
If there is $\elindexuno_1\in\setindexuno$ with $\setsymuno_{\elindexuno_1}=\udenset$ then by Notation \ref{ins}-[3] we have $\pointedincl<\left\{\left(\setsymuno_{\elindexuno}, \elsymuno_{\elindexuno}\right)\right\}_{\elindexuno\in\setindexuno}<>\setindexuno_0>=\emptysetfun$ for any $\setindexuno_0\subseteq\setindexuno$.\newline
With an abuse of language we replace pairs $\left(\setsymuno_{\elindexuno}, \elsymuno_{\elindexuno}\right)$ by sets $\setsymuno_{\elindexuno}$, or we drop any reference to $\left\{\left(\setsymuno_{\elindexuno}, \elsymuno_{\elindexuno}\right)\right\}_{\elindexuno\in\setindexuno}$ or to $\setindexuno_0$ whenever no confusion is possible.
\item  Fix $m \in \mathbb{N}$, an ordered family of sets $\left\{\setsymuno_1,...,\setsymuno_m\right\}$ and a permutation $\perm$ of the set $\left\{1,...,m\right\}$.\newline
If $\setsymuno_{\mathsf{m}}=\udenset$ for some ${\mathsf{m}} \in \left\{1,...,m\right\}$ then we set $\switch<\setsymuno_1,...,\setsymuno_m<>\perm>=\emptysetfun$.\newline
If $\setsymuno_{\mathsf{m}}\neq \udenset$ for any ${\mathsf{m}} \in \left\{1,...,m\right\}$ then 
\begin{equation}
\switch<\setsymuno_1,...,\setsymuno_m<>\perm>:\overset{m}{\underset{{\mathsf{m}}=1}{\prod}}\setsymuno_{\mathsf{m}} \rightarrow \overset{m}{\underset{{\mathsf{m}}=1}{\prod}}\setsymuno_{\perm\left({\mathsf{m}}\right)}
\end{equation}
is the set function defined by setting 
\begin{equation*}
\switch<\setsymuno_1,...,\setsymuno_m<>\perm>\left(\elsymuno_1,...,\elsymuno_m\right)=\left(\elsymuno_{\perm\left(1\right)},...,\elsymuno_{\perm\left(m\right)}\right) \hspace{10pt} \forall \left(\elsymuno_1,...,\elsymuno_m\right)\in \overset{m}{\underset{{\mathsf{m}}=1}{\prod}}\setsymuno_{\mathsf{m}}\text{.} 
\end{equation*}
If $\perm\left({\mathsf{m}}\right)={\mathsf{m}}$ for any ${\mathsf{m}} \in \left\{1,...,m\right\}$ then $\switch<\setsymuno_1,...,\setsymuno_m<>\perm>$ is the identity on $\overset{m}{\underset{{\mathsf{m}}=1}{\prod}}\setsymuno_{\mathsf{m}}$.\newline
If $m=2$ and $\perm\left(1\right)=2$ and $\perm\left(2\right)=1$ then we drop $\perm$ in writing $\switch<\setsymuno_1,\setsymuno_2<>\perm>$.\newline
If $\setsymuno_1=...=\setsymuno_{m}$ then we the list $\setsymuno_1,...,\setsymuno_m$ will be replaced by $\setsymuno_1$.
\item Fix a set $\setsymuno$, an equivalence relation $\approx $ on $\setsymuno$. We denote by $\frac{\setsymuno}{\approx}$ the quotient set.
\item Fix a set $\setsymuno$. We denote by $\partset{\setsymuno}$ the set of all subsets of $\setsymuno$. We set $\finpartset{\setsymuno}=\left\{\setsymdue\in \partset{\setsymuno}\;:\;\setsymdue\;\text{is a finite set}\right\}$.\newline
Fix $\partiz\subseteq \partset{\setsymuno}$.\newline
We say that $\partiz$ is a partition of $\setsymuno$ if and only if all conditions below are fulfilled: $\udenset\notin \partiz$; $\setsymdue\cap \setsymtre=\udenset$ for any $\setsymdue, \setsymtre\in \partiz$ with $\setsymdue\neq \setsymtre$; $\setsymuno=\underset{\setsymdue\in\partiz}{\bigcup}\setsymdue$.\newline
We say that $\partiz$ is a finite partition of $\setsymuno$ if $\partiz$ is a finite set.\newline
Fix $ \setsymcinque,\setsymsei\subseteq \partset{\setsymuno}$, $\setsymquattro\subseteq \setsymuno$.\newline
We say that $\setsymsei$ is a family of $\setsymcinque$-sets if and only if $\setsymsei\subseteq\setsymcinque$.\newline
We say that $\setsymsei$ is a $\setsymcinque$-cover of $\setsymquattro$ if and only if both conditions below are fulfilled: 
$\setsymsei$ is a family of $\setsymcinque$-sets; $\setsymquattro=\underset{\setsymdue\in\setsymsei}{\bigcup}\,\setsymdue$.\newline
We say that $\setsymsei$ is a sieve-type family of $\setsymcinque$-sets if and only if both conditions below are fulfilled: 
$\setsymsei$ is a family of $\setsymcinque$-sets; for any $\setsymdue\in\setsymsei$, $\setsymtre \in \setsymcinque$ with $\setsymtre \in \setsymdue$ we have that $\setsymtre \in \setsymsei$.\newline
\item We denote by $\Setcat$ the category of all sets and set functions.
\end{enumerate}
\end{notation}

\begin{notation}\label{alg}
We establish notation relative to general algebra, we refer to books of basic algebra for details.
\begin{enumerate}
\item We call magma a pair $\left(\setsymuno,\bot\right)$ where $\setsymuno$ is a set and $\bot:\setsymuno\times \setsymuno \rightarrow \setsymuno$ is a set function.\newline
Consider two magmas $\left(\setsymuno,\bot\right)$, $\left(\setsymdue,\top\right)$ and a set function $f:\setsymuno \rightarrow \setsymdue$. We say that $f$ is a function of magmas if and only if
\begin{equation*}
f\left(\bot\left(\elsymuno_1,\elsymuno_2\right)\right)=\top\left(f\left(\elsymuno_1\right), f\left(\elsymuno_2\right)\right)\qquad \forall \elsymuno_1,\elsymuno_2 \in \setsymuno\text{,}
\end{equation*}
we use the compact notation $f:\left(\setsymuno,\bot\right)\rightarrow \left(\setsymdue,\top\right)$ to mean that $f$ is a function of magmas.
\item Fix two sets $\setsymuno$, $\setsymdue$, two set functions $\bot:\setsymdue \times\setsymuno \rightarrow \setsymuno$, $\top:\setsymuno \times \setsymuno \rightarrow \setsymuno$, a relation (of any kind) $\prec$ on $\setsymuno$.\newline
We say that $\prec$ is compatible with $\bot$ if and only if
\begin{equation*}
\elsymdue\bot \elsymuno_1 \prec \elsymdue \bot \elsymuno_2 \qquad \forall \elsymdue \in \setsymdue\quad\forall \elsymuno_1,\elsymuno_2 \in \setsymuno \qquad \text{with}\quad \elsymuno_1\prec \elsymuno_2\text{.}
\end{equation*}
We say that $\prec$ is compatible with $\top$ if and only if 
\begin{equation*}
\elsymuno_1\top \elsymuno_2 \prec \elsymuno_3 \top \elsymuno_4 \qquad \forall \elsymuno_1,\elsymuno_2,\elsymuno_3,\elsymuno_4 \in X \quad \text{with}\quad \elsymuno_1\prec \elsymuno_3\text{,}\quad \elsymuno_2\prec \elsymuno_4 \text{.}
\end{equation*}
\item  Fix a commutative ring $\ringsym$ with unit. A weak $\ringsym$-algebra is a four-tuple $\left(\setsymuno,\algplus, \algper, \algsclp\right)$ where
\begin{equation*}
\algplus : \setsymuno\times \setsymuno \rightarrow \setsymuno\text{,}\qquad \algper : \setsymuno\times \setsymuno \rightarrow \setsymuno\text{,}\qquad \algsclp : \ringsym\times \setsymuno \rightarrow \setsymuno
\end{equation*}
are set functions fulfilling all conditions below: 
\begin{flalign}
&\left\{
\begin{array}{ll}
\left(\elsymuno\algplus \elsymdue\right)\algplus \elsymtre=\elsymuno\algplus \left(\elsymdue\algplus \elsymtre\right)& \forall \elsymuno,\elsymdue,\elsymtre\in \setsymuno\text{,}\\
\elsymuno\algplus \elsymdue=\elsymdue\algplus \elsymuno& \forall \elsymuno,\elsymdue\in \setsymuno\text{,}\\
\exists \neutsumsym\in \setsymuno:\;\; \elsymuno\algplus \neutsumsym=\elsymuno & \forall \elsymuno\in \setsymuno\text{,} 
\end{array}
\right.\label{plussa}\\[8pt]
&\left\{
\begin{array}{ll}
\left(\elsymuno\algper \elsymdue\right)\algper \elsymtre=\elsymuno\algper \left(\elsymdue\algper \elsymtre\right)&\forall \elsymuno,\elsymdue,\elsymtre\in \setsymuno\text{,}\\
\exists \neutmultsym\in \setsymuno:\;\; \elsymuno\algper \neutmultsym=\neutmultsym\algper \elsymuno=\elsymuno & \forall \elsymuno\in \setsymuno\text{,}\\
\end{array}
\right.\nonumber\\[8pt]
&\left\{
\begin{array}{ll}
\left(\elringsymuno+\elringsymdue\right)\algsclp \elsymuno=\left(\elringsymuno\algsclp \elsymuno\right) \algplus \left(\elringsymdue\algsclp \elsymuno\right)&\forall \elringsymuno, \elringsymdue\in \ringsym\quad  \forall \elsymuno\in \setsymuno \text{,}\\
\elringsymuno\algsclp\left(\elringsymdue\algsclp \elsymuno\right)=\left(\elringsymuno\elringsymdue\right)\algsclp \elsymuno& \forall \elringsymuno, \elringsymdue\in \ringsym\quad \forall \elsymuno\in \setsymuno \text{,}\\
\unoring\algsclp \elsymuno=\elsymuno& \forall \elsymuno\in \setsymuno\text{,}\\
\elringsymuno\algsclp\left(\elsymuno\algplus  \elsymdue\right)=\left(\elringsymuno\algsclp \elsymuno\right)\algplus\left( \elringsymuno\algsclp  \elsymdue\right) & \forall \elringsymuno\in \ringsym\quad\forall \elsymuno, \elsymdue\in \setsymuno \text{,}\\
\elringsymuno\algsclp\left(\elsymuno\algper \elsymdue\right)=\left(\elringsymuno\algsclp \elsymuno\right)\algper  \elsymdue&\forall \elringsymuno\in \ringsym\quad \forall \elsymuno, \elsymdue\in\setsymuno \text{,}\\
\elringsymuno\algsclp\left(\elsymuno\algper  \elsymdue\right)= \elsymuno\algper \left(\elringsymuno\algsclp  \elsymdue\right) &  \forall \elringsymuno\in \ringsym\quad \forall \elsymuno, \elsymdue\in \setsymuno\text{.}
\end{array}
\right.\label{scalpsa}
\end{flalign}
We call: $\algplus$ sum; $\algper$ product; $\algsclp$ scalar product.\newline 
The weak $\ringsym$-algebra $\left(\setsymuno,\algplus, \algper, \algsclp\right)$ is commutative if and only if
\begin{equation}
\elsymuno\algper \elsymdue=\elsymdue\algper \elsymuno\qquad \forall \elsymuno,\elsymdue\in \setsymuno\text{.}\label{comalg}
\end{equation}
\item Fix a commutative ring $\ringsym$ with unit, a weak $\ringsym$-algebra $\left(\setsymuno,\algplus, \algper,\algsclp\right)$ and a subset $\setsymdue\subseteq \setsymuno$.\newline
We say that the triplet $\left(\setsymdue,\algplus, \algper, \algsclp\right)$ is a sub weak $\ringsym$-algebra of the weak $\ringsym$-algebra $\left(\setsymuno,\algplus,\algper, \algsclp\right)$ if and only if
\begin{equation}
\left\{
\begin{array}{ll}
\neutsumsym\in \setsymdue\text{,}\\
\elsymuno\algplus \elsymdue,  \elsymuno\algper\elsymdue\in \setsymdue & \forall \elsymuno,\elsymdue \in \setsymdue\text{,}\\
\elringsymuno\algsclp \elsymuno\in \setsymdue & \forall \elringsymuno\in \ringsym\quad \forall \elsymuno \in \setsymdue\text{.}
\end{array}
\right.\label{semalg}
\end{equation}
We say that $\setsymdue$ is an ideal of $\left(\setsymuno,\algplus,\algper, \algsclp\right)$ if and only if
\begin{equation}
\left\{
\begin{array}{ll}
\elsymuno\algplus \elsymdue\in \setsymdue & \forall \elsymuno,\elsymdue \in \setsymdue\text{,}\\
\elsymuno\algper \elsymdue, \elsymdue\algper \elsymuno\in \setsymdue & \forall \elsymuno\in \setsymuno\quad \forall \elsymdue \in \setsymdue\text{,}\\
\elringsymuno\algsclp \elsymdue\in \setsymdue & \forall \elringsymuno\in \ringsym\quad \forall \elsymdue \in \setsymdue\text{.}
\end{array}
\right.
\end{equation}
\item Fix a commutative ring $\ringsym$ with unit, a weak $\ringsym$-algebra $\left(\setsymuno,\algplus, \algper, \algsclp\right)$, a relation $\prec$ on $\setsymuno$.\newline
We say that $\prec$ is compatible with the weak $\ringsym$-algebra structure of the weak $\ringsym$-algebra $\left(\setsymuno,\algplus, \algper, \algsclp\right)$ if and only if $\prec$ is compatible with both $\algplus$ and $\algper$ and $\algsclp$.
\item  Fix a commutative ring $\ringsym$ with unit. A weak $\ringsym$-module is a triplet $\left(\setsymuno,\algplus, \algsclp\right)$ where \hspace{20pt}$\algplus : \setsymuno\times \setsymuno \rightarrow \setsymuno\text{,}\qquad \algsclp : \ringsym\times \setsymuno \rightarrow \setsymuno$\newline
are set functions fulfilling \eqref{plussa} and all conditions among \eqref{scalpsa} which don't refer to $\algper$.\newline
The weak $\ringsym$-module $\left(\setsymuno,\algplus, \algsclp\right)$ is commutative if and only if \eqref{comalg} holds true.
\item Fix a commutative ring $\ringsym$ with unit, a weak $\ringsym$-module $\left(\setsymuno,\algplus, \algsclp\right)$ and a subset $\setsymdue\subseteq \setsymuno$.\newline
We say that the triplet $\left(\setsymdue,\algplus,  \algsclp\right)$ is a sub weak $\ringsym$-module of $\left(\setsymuno,\algplus, \algsclp\right)$ if and only if $\left(\setsymdue,\algplus,  \algsclp\right)$ fulfills all conditions among \eqref{semalg} which don't refer to $\algper$.
\item Fix a commutative ring $\ringsym$ with unit, a weak $\ringsym$-module $\left(\setsymuno,\algplus, \algsclp\right)$, a relation $\prec$ on $\setsymuno$.\newline
We say that $\prec$ is compatible with the weak $\ringsym$-module structure of the weak $\ringsym$-module $\left(\setsymuno,\algplus, \algsclp\right)$ if and only if $\prec$ is compatible with both $\algplus$ and $\algsclp$.
\item Fix a set $\setindexuno$, a commutative ring $\ringsym$, a weak $\ringsym$-module $\setsymuno_{\elindexuno}$ for any $\elindexuno\in \setindexuno$.\newline 
Referring to Notation \ref{ins}-[11] we set $\pointedincl<\right\{\setsymuno_{\elindexuno}\left\}_{\elindexuno\in \setindexuno}<>\elindexuno_0>=\pointedincl<\right\{\left(\setsymuno_{1},0\right)\left\}_{\elindexuno\in \setindexuno}<>\elindexuno_0>$ for any $\elindexuno_0\in \setindexuno$.
\item We denote by $\VecR$ the category of all $\mathbb{R}$-vector spaces and functions between $\mathbb{R}$-vector spaces.
\item We denote by $\AlgR$ the category of all $\mathbb{R}$-algebras and functions between $\mathbb{R}$-algebras.
\item We denote by $\tenalgfun:\VecR\rightarrow\AlgR$ the functor associating to any $\mathbb{R}$-vector space the corresponding  
tensor $\mathbb{R}$-algebra and to any function between $\mathbb{R}$-vector spaces the corresponding function of tensor $\mathbb{R}$-algebras. It is well known that $\tenalgfun$ is an exact functor, where exactness follows since $\mathbb{R}$ is a field.\newline
For any object $\setsymuno$ of $\VecR$ we denote by $\unoT[\setsymuno]$ the unity of $\tenalgfun\left(\setsymuno\right)$.
\item We denote by $\extalgfun:\VecR\rightarrow\AlgR$ the functor associating to any $\mathbb{R}$-vector space the corresponding  
exterior $\mathbb{R}$-algebra and to any function between $\mathbb{R}$-vector spaces the corresponding function of exterior $\mathbb{R}$-algebras. It is well known that $\extalgfun$ is an exact functor, where exactness follows since $\mathbb{R}$ is a field.
\item We denote by $\symalgfun:\VecR\rightarrow\AlgR$ the functor associating to any $\mathbb{R}$-vector space the corresponding  
symmetric $\mathbb{R}$-algebra and to any function between $\mathbb{R}$-vector spaces the corresponding function of symmetric $\mathbb{R}$-algebras. It is well known that $\symalgfun$ is an exact functor, where exactness follows since $\mathbb{R}$ is a field.
\item Fix $k\in \mathbb{N}_0$, one of the functors $\tenalgfun$, $\extalgfun$, $\symalgfun$ and denote it by $\genfunc$.\newline
We denote by $\genfunc[k]:\VecR\rightarrow \VecR$ the functor given by the degree $k$ homogeneous component of $\genfunc$. It is well known that $\genfunc[k]$ is an exact functor, where exactness follows since $\mathbb{R}$ is a field.\newline
With an abuse of language, whenever no confusion is possible, we denote $\genfunc[k]\left(\setsymuno\right)$ by $\genfunc\left(\setsymuno\right)_k$ for any $\mathbb{R}$-vector space $\setsymuno$.
\item We denote by $\Aug:\AlgR\rightarrow\AlgR$ the functor associating to any $\mathbb{R}$-algebra the corresponding  
unitary $\mathbb{R}$-algebra and to any function between $\mathbb{R}$-algebras the corresponding function of unitary $\mathbb{R}$-algebras.\newline
For any object $\setsymuno$ of $\AlgR$ we denote by $\unoA[\setsymuno]$ the unity of $\Aug\left(\setsymuno\right)$, by $\inclunoA[\setsymuno]:\setsymuno \rightarrow \Aug\left(\setsymuno\right)$ the canonical inclusion.\newline
If $\setsymdue$ is an ideal or a proper ideal of $\setsymuno$ then $\inclunoA[\setsymuno]\left(\setsymdue\right)$ is respectively an ideal or a proper ideal of $\Aug\left(\setsymuno\right)$.\newline
For any $\setsymdue \subseteq \setsymuno$ we denote $\inclunoA[\setsymuno]\left(\setsymdue\right)$ again by $\setsymdue$, with an abuse of language.
\end{enumerate}
\end{notation}

\begin{notation}\label{gentopnot}We establish notation relative to general topology, we refer to any book of basic topology for details.
\begin{enumerate}
\item A topological space is usually denoted by a pair $\left(\setsymuno, \topo\right)$, where $\setsymuno$ is the underlying set and $\topo$ is the chosen topology. In this work, if not otherwise specified, we consider only one topology on any fixed set $\setsymuno$, for this reason we drop the symbol $\topo$ from the notation and we simply refer to the topological space $\setsymuno$.
\item Fix a set $\setsymuno$, a not empty family of topological spaces $\left\{\left(\setsymdue_i, \topo_i\right)\right\}_{i \in I}$ (we are not assuming $\left(\setsymdue_i, \topo_i\right)\neq \left(\setsymdue_j, \topo_j\right)$ for $i\neq j$), a family $\left\{f_i\right\}_{i \in I}$ of set functions with $f_i:\setsymuno \rightarrow \setsymdue_i$ for any $i \in I$.\newline
We define topology $\topo$ on $\setsymuno$ by setting
\begin{flalign*}
&\openset\in \topo\qquad \Leftrightarrow \qquad \exists \; \text{a set}\;J\qquad \exists \left\{\openset_{i,j}\right\}_{i\in I, j\in J} \qquad\text{such that}\\[4pt]
&\openset_{i,j}\in \topo_i \quad \forall i \in I\quad \forall j\in J\qquad\text{and}\qquad \openset=\underset{j\in J}{\bigcup}\left(\underset{i\in I}{\bigcap} \,f_i^{-1}\left(\openset_{i,j}\right)\right)\text{.}
\end{flalign*}
We say that $\topo$ is the initial topology on $\setsymuno$ with respect to $\left\{f_i\right\}_{i \in I}$.\newline
We say that $\setsymdue$ is endowed with the initial topology with respect to $\left\{f_i\right\}_{i \in I}$.
\item Fix a topological space $\left(\setsymuno, \topo\right)$, a set $\setsymdue \subseteq \setsymuno$, a set $\setsymtre \subseteq \setsymdue$. We say that the initial topology on $\setsymdue$ with respect to the single function family $\left\{\incl{\setsymdue}{\setsymuno}\right\}$ is the subspace topology of $\setsymdue$ with respect to $\left(\setsymuno, \topo\right)$. We say that $\setsymdue$ is endowed with the subspace topology with respect to $\left(\setsymuno, \topo\right)$, or simply that $\setsymdue$ is endowed with the subspace topology when no confusion is possible. We say that a set $\setsymtre$ is open in $\setsymdue$ if and only if $\setsymtre$ is open with respect to the subspace topology of $\setsymdue$.
\item Fix a not empty family of topological spaces $\left\{\left(\setsymuno_i, \topo_i\right)\right\}_{i\in I}$. We denote by $\underset{i\in I}{\prod}\topo_i$ the product topology on $\underset{i\in I}{\prod}\setsymuno_i$. If not otherwise stated we always assume that the product of a family of topological spaces is endowed with the product topology.  
\item Fix a set $\setsymdue$, a not empty family of topological spaces $\left\{\left(\setsymuno_i, \topo_i\right)\right\}_{i \in I}$ (we are not assuming $\left(\setsymuno_i, \topo_i\right)\neq \left(\setsymuno_j, \topo_j\right)$ for $i\neq j$), a family $\left\{f_i\right\}_{i \in I}$ of set functions with $f_i:\setsymuno_i \rightarrow \setsymdue$ for any $i \in I$.\newline
We define topology $\topo$ on $\setsymdue$ by setting
\begin{equation*}
\openset\in \topo \Leftrightarrow f_i^{-1}\left(\openset\right)\in \topo_i \quad \forall i \in I\text{.}
\end{equation*}
We say that $\topo$ is the final topology on $\setsymdue$ with respect to $\left\{f_i\right\}_{i \in I}$ with $f_i:\setsymuno_i \rightarrow \setsymdue$.\newline
We say that $\setsymdue$ is endowed with the final topology with respect to $\left\{f_i\right\}_{i \in I}$ with $f_i:\setsymuno_i \rightarrow \setsymdue$.
\item Fix a topological space $\left(\setsymuno, \topo\right)$, an equivalence relation $\approx$ in $\setsymuno$. Set $\setsymdue= \setsymuno\setminus \approx $ the quotient space, $\setquotfun:\setsymuno\rightarrow \setsymdue$ the quotient function. We say that the final topology on $\setsymdue$ with respect to $\left\{\setquotfun\right\}$ is the quotient topology of $\setsymdue$ with respect to $\setquotfun$. We say that $\setsymdue$ is endowed with the quotient topology with respect to $\setquotfun$.
\item Fix a topological space $\left(\setsymuno, \topo\right)$, a set $\setsymdue \subseteq \setsymuno$, $\elsymuno \in \setsymdue$.\newline
We say that $\setsymdue$ is compactly contained in $\setsymuno$, and we write $\setsymdue \compacont \setsymuno$, to mean that exists a compact subset $\setsymtre$ of $\setsymuno$ such that $\setsymdue\subseteq \setsymtre$.\newline
We denote by $\interior\left[\setsymdue,\left(\setsymuno, \topo\right) \right]$ the interior of $\setsymdue$ in $\left(\setsymuno, \topo\right)$. We adopt the simpler notation $\interior\left[\setsymdue,\setsymuno \right]$, whenever no confusion is possible on the topology of $\setsymuno$.\newline
We denote by $\closure\left[\setsymdue,\left(\setsymuno, \topo\right)\right]$ the closure of $\setsymdue$ in $\left(\setsymuno, \topo\right)$. We adopt the simpler notation $\closure\left[\setsymdue,\setsymuno \right]$, whenever no confusion is possible on the topology of $\setsymuno$.\newline
A set $\setsymdue_1$ is a neighborhood of $\elsymuno$ in $\setsymdue$ if and only if there is $\setsymdue_2\in \topo$ with $\elsymuno \in \setsymdue \cap \setsymdue_2 \subseteq \setsymdue_1$. 
\item Fix a topological space $\left(\setsymuno, \topo\right)$, a set $\setsymdue \subseteq \setsymuno$, $\elsymuno \in \setsymuno$. We write: $\setsymdue \subseteqdentro \setsymuno$ if and only if $\setsymdue \in \topo$; $\elsymuno\dentro \setsymdue$ if and only if $\elsymuno \in \interior\left[\setsymdue,\left(\setsymuno, \topo\right) \right]$. We say that $\elsymuno$ is adherent to $\setsymdue$, or that $\elsymuno$ is in the adherence of $\setsymdue$, if and only if for any neighborhood $\setsymtre$ of $\elsymuno$ we have $\setsymdue \cap\left(\setsymtre\setminus\left\{\elsymuno\right\}\right)\neq\udenset$.
\item Fix a topological space $\left(\setsymuno, \topo\right)$, a set function $f:\left(-1,1\right)\rightarrow \setsymuno$. We say that $f$ is a $\topo$-path in $\setsymuno$ if and only if $f$ is continuous with respect to the Euclidean topology on $\left(-1,1\right)$ and to $\topo$ on $\setsymuno$. We say simply that $f$ is a path in $\setsymuno$, by dropping the reference to $\topo$, whenever no confusion is possible.
\item Fix a topological space $\left(\setsymuno, \topo\right)$, a set $\setsymquattro\subseteq \setsymuno$. Referring to Notation \ref{ins}-[14] and choosing $\setsymcinque=\topo$ we define the notion of family of open sets, open cover of $\setsymquattro$, sieve-type family of open sets.
\item We denote by $\Topcat$ the category of all topological spaces and continuous functions. Category $\Topcat$ is naturally endowed with the Grothendieck topology $\GrTc$ induced by the topological structure of its objects. 
\item Fix a topological space $\left(\setsymuno, \topo\right)$. We denote by $\SiteSpTop(\setsymuno, \topo)$ the site whose objects are all elements of $\topo$, arrows are inclusions of open sets and Grothendieck topology is given by selecting for any $U\in \topo$ all sieves containing at least one covering of $U$ by open sets. With an abuse of language we drop the symbol $\topo$ whenever no confusion is possible.
\item  Fix a topological space $\left(\setsymuno, \topo\right)$. We say that $\left(\setsymuno, \topo\right)$ is locally compact if for any $\elsymuno \in\setsymuno$ and any $\setsymsei\in\topo$ with $\elsymuno \in \setsymsei$ there are a compact set $K$, $\setsymsette\in\topo$ with $\elsymuno \in \setsymsette \subseteq K \subseteq \setsymsei$.
\item  Fix a topological space $\left(\setsymuno, \topo\right)$. We say that $\left(\setsymuno, \topo\right)$ fulfills the first numerability axiom if and only if for any point $\elsymuno\in \setsymuno$ there is a countable set $\setsymtre_{\elsymuno}\subseteq\topo$ such that for any $\setsymcinque\in \topo$ with $\elsymuno\in \setsymcinque$ there are $\setsymsei_1,\setsymsei_2\in \setsymtre$ with $\setsymsei_1\subseteq\setsymcinque\subseteq\setsymsei_2$.
\end{enumerate}
\end{notation}

\begin{notation}\label{realvec} 
Let $\mathbb{R}$ denotes the field of real numbers, we establish notation relative to vector spaces over $\mathbb{R}$.
\begin{enumerate}
\item We define $\realpartset=\left\{\intuno\;:\;\exists m \in \mathbb{N}_0\;\text{with}\; \intuno\subseteq \mathbb{R}^m\right\}$.
\item Fix $m\in \mathbb{N}_0$.\newline
If $m=0$ then we set $\mathbb{R}^{0}=\{0\}$. \newline
If $m\geq 1$ then we denote by $\mathbb{R}^m$ the real vector space of dimension $m$ and by $\left\{\ef_1,...,\ef_m \right\}$ the standard orthonormal Euclidean frame.\newline
The usual notation is adopted to denote sum and scalar product; inner product is denoted by $\innpreul\;\innpreur_{m}$, the associated norm is denoted by $\normeul\;\normeur_{m} $, in both cases the index $m$ is omitted whenever there is no risk of confusion.\newline
Elements belonging to $\mathbb{R}^m$ are denoted by $\unkuno$ or $\unkdue$ when we consider them just points belonging to the set $\mathbb{R}^m$, or by $\vecuno$ when we consider them as vectors belonging to $\mathbb{R}^m$.\newline
We also use the extended notation $\left(\unkuno_1,...,\unkuno_{\overline{m}}\right)$,  $\left(\unkdue_1,...,\unkdue_{\overline{m}}\right)$, $\left(\vecuno_1,...,\vecuno_{\overline{m}}\right)$ when we want to emphasize components with respect to the Euclidean frame, in particular:\newline
if $m=0$ then the extended notation becomes $\unkuno_1$, $\unkdue_1$ or $\vecuno_1$ respectively; \newline
if $m>0$ then the extended notation becomes $\left(\unkuno_1,...,\unkuno_m\right)$, $\left(\unkdue_1,...,\unkdue_m\right)$ or $\left(\vecuno_1,...,\vecuno_m\right)$ respectively.
\item Fix $m \in \mathbb{N}_0$, $\intuno\subseteq\mathbb{R}^m$. Refer to Notation \ref{gentopnot}-[1, 3].\newline
The real vector space $\mathbb{R}^m$ is endowed with the Euclidean topology.\newline
The subset $\intuno$ of $\mathbb{R}^m$ is endowed with the subspace topology with respect to $\mathbb{R}^m$.
\item Fix $m\in \mathbb{N}_0$, $\intuno, \intdue \subseteq \mathbb{R}^m$.\newline
If $m=0$ then we write $\intdue \sqsubseteq \intuno$ if and only if $\intdue \subseteq \intuno$ and $\intdue = \mathbb{R}^0$.\newline 
If $m = 1$ then we write $\intdue \sqsubseteq \intuno$ if and only if $\intdue \subseteq \intuno$ and one among the following mutually exclusive cases occurs
\begin{equation*}
\begin{array}{l}
\intdue=\mathbb{R} \text{,} \vspace{4pt}\\
\exists \unkuno \in \mathbb{R}\,:\; \intdue=\left(-\infty,\unkuno\right] \text{,} \vspace{4pt}\\
\exists \unkuno \in \mathbb{R}\,:\; \intdue=\left(-\infty,\unkuno\right) \text{,} \vspace{4pt}\\
\exists \unkuno \in \mathbb{R}\,:\;\intdue=\left[\unkuno, +\infty\right) \text{,} \vspace{4pt}\\
\exists \unkuno \in \mathbb{R}\,:\;\intdue=\left(\unkuno, +\infty\right) \text{,} \vspace{4pt}\\
\exists \unkuno_1,\unkuno_2 \in \mathbb{R}\,:\;\;\unkuno_1<\unkuno_2\; \;\text{and}\;\;\intdue=\left(\unkuno_1, \unkuno_2\right]\text{,} \vspace{4pt}\\
\exists \unkuno_1,\unkuno_2 \in \mathbb{R}\,:\;\;\unkuno_1<\unkuno_2\; \;\text{and}\;\;\intdue=\left[\unkuno_1, \unkuno_2\right)\text{,} \vspace{4pt}\\
\exists \unkuno_1,\unkuno_2 \in \mathbb{R}\,:\;\;\unkuno_1<\unkuno_2\; \;\text{and}\;\;\intdue=\left(\unkuno_1, \unkuno_2\right)\text{.}
\end{array}
\end{equation*}
If $m \geq 1$ then we write $\intdue \sqsubseteq \intuno$ if and only if
\begin{equation*}
\begin{array}{l}
\intdue \subseteq \intuno\qquad \text{and}\qquad  \exists\, \intdue_1,...,\intdue_m \sqsubseteq \mathbb{R}\;:\quad  \intdue=\overset{m}{\underset{{\mathsf{m}}=1}{\prod}} \intdue_{\mathsf{m}}\text{.}
\end{array}
\end{equation*}
We write $\intdue \sqsubseteqdentro \intuno$ if and only if $\intdue \sqsubseteq \intuno$ and there is $\setsymcinque\subseteqdentro\mathbb{R}^m$ with $\intdue= \intuno \cap \setsymcinque$.\newline
We say that: $\intdue$ is an interval if and only if $\intdue \sqsubseteq\mathbb{R}^m$; $\intdue$ is an open interval if and only if $\intdue \sqsubseteqdentro\mathbb{R}^m$.
\item Fix $m\in \mathbb{N}_0$, $\intuno, \intdue \subseteq \mathbb{R}^m$.\newline
If $m=0$ then we write $\intdue\opsspint\intuno$ if and only if $\intdue\subseteq\intuno$ and $\intdue=\mathbb{R}^0$.\newline
If $m\geq 1$ then we write $\intdue\opsspint\intuno$ if and only if  $\intdue\subseteq\intuno$ and $\intdue$ is an open set of $ \left[0,+\infty\right)\times\mathbb{R}^{m-1}$.   
\item Fix $m \in \mathbb{N}_0$, $i \in \mathbb{N}$, $\intuno\subseteq \mathbb{R}^m$.\newline
If $i \leq m$ then we define 
\begin{multline*}
\domint\left[\intuno,i\right]=
\big\{\unkuno\in\mathbb{R}^{m+1}\;:\; \left(\unkuno_1,...,\unkuno_{i-1},\unkdue,\unkuno_{i+2},...,\unkuno_{m+1}\right)\in \intuno\\ \forall \unkdue \in \left[\min\left\{\unkuno_{i}, \unkuno_{i+1}\right\}, \max\left\{\unkuno_{i}, \unkuno_{i+1}\right\}\right]\big\}\text{.}
\end{multline*}
If $i > m$ then we define $\domint\left[\intuno,i\right]= \intuno \times \mathbb{R}^{i-m+1}$.
\item Fix $m \in \mathbb{N}_0$, $\intuno\opsspint \mathbb{R}^m$. We associate to $\intuno$ a suitable open set $\att\left[\intuno\right] \subseteqdentro \mathbb{R}^m$ and a suitable set function $\attfun:\att\left[\intuno\right] \rightarrow \intuno$.\newline
If $\intuno$ is an open set of $\mathbb{R}^m$ then we define: $\att\left[\intuno\right]=\intuno$; $\attfun\left[\intuno\right]=\idobj+\intuno+$.\newline
If $m\geq 1$ and $\intuno$ is an open set of $\left[0,+\infty\right)\times\mathbb{R}^{m-1}$ and $\intuno$ is not an open set of $\mathbb{R}^m$ then we define: 
\begin{flalign*}
&\att\left[\intuno\right]=\left\{\left(\unkuno,\unkuno_1,...,\unkuno_{\overline{m}}\right)\in \mathbb{R}^m\;:\;\left(|\unkuno|,\unkuno_1,...,\unkuno_{\overline{m}}\right)\in \intuno\right\}\text{;}\\[6pt]
&\attfun\left[\intuno\right]\left(\unkuno,\unkuno_1,...,\unkuno_{\overline{m}}\right)=\left(|\unkuno|,\unkuno_1,...,\unkuno_{\overline{m}}\right) \hspace{10pt} \forall\left(\unkuno,\unkuno_1,...,\unkuno_{\overline{m}}\right)\in  \att\left[\intuno\right]\text{.}
\end{flalign*}
\item Fix $m \in \mathbb{N}_0$. We denote by $\eumea\left[m\right]$ the Lebesgue measure on the Borel $\sigma$-algebra of $\mathbb{R}^m$ endowed with euclidean topology. We drop $m$ from the notation whenever no confusion is possible.
\item Fix $m\in\mathbb{N}_0$, a set $\setsymquattro\subseteq \mathbb{R}^m$. Referring to Notation \ref{ins}-[14] and setting $\setsymcinque=\left\{\intuno\;:\;\intuno\sqsubseteqdentro\mathbb{R}^m\right\}$ we define the notion of family of open intervals, open intervals cover of $\setsymquattro$, sieve-type family of open intervals.
\end{enumerate}
\end{notation}

\begin{remark}\label{rprop}We state some basic properties of sets introduced in  Notation \ref{realvec}-[4, 6].
\begin{enumerate}
\item Fix $n \in \mathbb{N}$, a finite not empty set $\{m_1,...,m_n\}\subseteq \mathbb{N}$, subsets $\intuno_{\mathsf{n}}, \intdue_{\mathsf{n}} \subseteq \mathbb{R}^{m_{\mathsf{n}}}$ for ${\mathsf{n}}\in\{1,...,n\}$. Assume $ \intuno_{\mathsf{n}} \sqsubseteq \intdue_{\mathsf{n}}\quad \forall {\mathsf{n}}\in\{1,...,n\}$.\newline
Then: 
$\overset{n}{\underset{{\mathsf{n}}=1}{\prod}}\, \intuno_{\mathsf{n}} \sqsubseteq \overset{n}{\underset{{\mathsf{n}}=1}{\prod}}\, \intdue_{\mathsf{n}} \subseteq \mathbb{R}^{m_1 + ... + m_n}$.
\item Fix  $m\subseteq \mathbb{N}_0$, $n \in \mathbb{N}$, subsets $\intuno_1,...,\intuno_n , \intdue_1,...,\intdue_n, \subseteq \mathbb{R}^{m}$. Assume that all conditions below are fulfilled:
\begin{equation*}
	\intuno_{\mathsf{n}} \sqsubseteq \intdue_{\mathsf{n}}\quad \forall {\mathsf{n}}\in \{1,...,n\}\text{,}\hspace{30pt} \overset{n}{\underset{{\mathsf{n}}=1}{\bigcap}}\, \interior\left[\intuno_{\mathsf{n}},\mathbb{R}^m\right]\sqsubseteq\mathbb{R}^m\text{.}
\end{equation*}
Then 
$ \overset{n}{\underset{{\mathsf{n}}=1}{\bigcap}}\, \intuno_{\mathsf{n}} \sqsubseteq \overset{n}{\underset{{\mathsf{n}}=1}{\bigcap}}\, \intdue_{\mathsf{n}} \subseteq \mathbb{R}^m$.
\item Fix $i,n \in \mathbb{N}$, a finite not empty set $\{m_1,...,m_n\}\subseteq \mathbb{N}$, subsets $\intuno_{\mathsf{n}} \subseteq \mathbb{R}^{m_{\mathsf{n}}}$ for ${\mathsf{n}}\in\{1,...,n\}$. 
Set $\deridxuno_0=0$, $\deridxuno_{\mathsf{n}}=\overset{{\mathsf{n}}}{\underset{k=1}{\sum}}\, m_k$ for any ${\mathsf{n}}\in\{1,...,n\}$.\newline 
If $i\leq \deridxuno_n$ then $\domint\left[\overset{n}{\underset{{\mathsf{n}}=1}{\prod}}\, \intuno_{\mathsf{n}},i\right] = \left(\overset{n}{\underset{{\mathsf{n}}=1}{\prod}}\, \intdue_{\mathsf{n}}\right)$ where
\begin{equation*}
\intdue_{\mathsf{n}}=\left\{
\begin{array}{ll}
\intuno_{\mathsf{n}}&\text{if}\; \deridxuno_{\mathsf{n}} < i \;\text{or}\;i\leq \deridxuno_{{\mathsf{n}}-1}\text{,}\\[4pt]\\
\domint\left[\intuno_{\mathsf{n}},i-\deridxuno_{{\mathsf{n}}-1}\right]&\text{if}\; \deridxuno_{{\mathsf{n}}-1} < i\leq \deridxuno_{\mathsf{n}}\text{.}
\end{array}
\right.
\end{equation*}
If $i> \deridxuno_n$ then $\domint\left[\overset{n}{\underset{{\mathsf{n}}=1}{\prod}}\, \intuno_{\mathsf{n}},i\right] = \left(\overset{n}{\underset{{\mathsf{n}}=1}{\prod}}\, \intuno_{\mathsf{n}}\right)\times \mathbb{R}^{i-\deridxuno_n+1}$.\newline
\item Fix $m \in \mathbb{N}_0$, $i,j \in \mathbb{N}$, $\intuno, \intdue\subseteq \mathbb{R}^{m}$. Then:
\begin{flalign*}
&\text{if}\;\intuno\subseteq\intdue\;\text{then}\;\domint\left[\intuno,i\right]\subseteq \domint\left[\intdue,i\right]\text{;}\\[6pt]
&\text{if}\;\intuno\subseteqdentro\intdue\;\text{then}\;\domint\left[\intuno,i\right]\subseteqdentro \domint\left[\intdue,i\right]\text{;}\\[6pt]
&\domint\left[\domint\left[\intuno,i\right],j\right]=\domint\left[\domint\left[\intuno,j\!-\!1\right],i\right]\hspace{20pt} \text{if}\;\;i < j\text{;}\\[6pt]
&\domint\left[\domint\left[\intuno,i\right],j\right]=\domint\left[\domint\left[\intuno,j\right],i\!+\!1\right]\hspace{20pt} \text{if}\;\;i \geq j\text{.}
\end{flalign*}
\end{enumerate}
\end{remark}

\begin{notation}\label{realfunc}
We fix the notation relative to real functions.
\begin{enumerate}
\item Fix $m,n\in \mathbb{N}_{0}$, $\intuno\subseteq \mathbb{R}^{m}$, $f:\intuno\rightarrow \mathbb{R}^{n}$.\newline
If $n=0$ or $n=1$ then $f_{1}=f$ is always understood.\newline
If $n\geq 2$ then we denote by $f_{\mathsf{n}}$ the ${\mathsf{n}}$-th component of $f$ for ${\mathsf{n}}\in\{1,...,n\}$.\newline
We set $\dimCod\left(f\right)=m$ and $\dimDom\left(f\right)=n$.
\item Fix $m,n \in \mathbb{N}_0$. Refer to Notation \ref{ins}-[7]. Assume: $\intuno \subseteq \mathbb{R}^m$; $\intdue \subseteq \mathbb{R}^n$. Then we set: $\Dom\left(\cost<\intuno<>\intdue>+\unkdue+\right)=\intuno$, $\mathsf{Cod}\left(\cost<\intuno<>\intdue>+\unkdue+\right)=\mathbb{R}^{n}$.\newline
Fix $m \in \mathbb{N}_0$, $n \in \mathbb{N}$. Refer to Notation \ref{ins}-[8]. Assume $\intuno \subseteq \mathbb{R}^m$. Then we set: $\Dom\left(\Diag{\intuno}{n}\right)=\intuno$, $\mathsf{Cod}\left(\Diag{\intuno}{n}\right)=\mathbb{R}^{mn}$.\newline
Fix $m \in \mathbb{N}_0$. Refer to Notation \ref{ins}-[9]. Assume $\intuno \subseteq \intdue \subseteq \mathbb{R}^m$. Then we set: $\Dom\left(\incl{\intuno}{\intdue}\right)=\intuno$, $\mathsf{Cod}\left(\incl{\intuno}{\intdue}\right)=\mathbb{R}^{m}$.\newline
Fix $m_1,...,m_n \in \mathbb{N}_0$. Refer to Notation \ref{ins}-[10]. Assume $\intuno_{\mathsf{n}} \subseteq \mathbb{R}^{m_{\mathsf{n}}}$ for any ${\mathsf{n}} \in \left\{1,...,n\right\}$. Then we set:
\begin{flalign*}
&\Dom\left(\proj<\left\{\intuno_{\mathsf{n}}\right\}_{\mathsf{n}=1}^n<>\mathsf{n}>\right)=\underset{{\mathsf{n}}=1}{\overset{n}{\prod}}\, \intuno_{\mathsf{n}}\text{;}\\[4pt]
&\mathsf{Cod}\left(\proj<\left\{\intuno_{\mathsf{n}}\right\}_{\mathsf{n}=1}^n<>\mathsf{n}>\right)=\mathbb{R}^{m_{\mathsf{n}}} \hspace{15pt} \forall {\mathsf{n}} \in \left\{1,...,n\right\}\text{.}
\end{flalign*}
Fix $m_1,...,m_n \in \mathbb{N}_0$. Refer to Notation \ref{ins}-[12]. Assume $\intuno_{\mathsf{n}} \subseteq \mathbb{R}^{m_{\mathsf{n}}}$ for any ${\mathsf{n}} \in \left\{1,...,n\right\}$. Then we set:
\begin{flalign*}
&\Dom\left(\switch<\intuno_1,...,\intuno_n<>\sigma>\right)=\underset{{\mathsf{n}}=1}{\overset{n}{\prod}}\, \intuno_{\mathsf{n}}\text{;}\\[4pt]
&\mathsf{Cod}\left(\switch<\intuno_1,...,\intuno_n<>\sigma>\right)=\mathbb{R}^{m_1+...+m_n}\text{.}
\end{flalign*}
\item Fix $m\in \mathbb{N}_{0}$, $n \in \mathbb{N}$. Then:
\begin{enumerate}
\item $\vecsum<m<>n>:\left(\mathbb{R}^m\right)^n\rightarrow \mathbb{R}^m$ denotes the vector sum of $n$ vectors of dimension $m$;\newline
\item $\vecminus\left[m\right]:\mathbb{R}^m\rightarrow \mathbb{R}^m$ denotes the vector sign change.
\end{enumerate}
If there is no risk of confusion then we drop both $m$ and $n$ from the notation, with an abuse of language.
\item Fix $m,n\in \mathbb{N}_{0}$. If $m=n=0$ then set $l=0$, else set $l=\overline{m}\,\overline{n}$.\newline
We denote by $\vecprod\left[m,n\right]:\mathbb{R}^{m}\times \mathbb{R}^{n}\rightarrow \mathbb{R}^{l}$ the function defined by setting
\begin{equation*}
\vecprod\left[m,n\right]\left(\unkuno,\unkdue\right)=\left(\unkuno_1\unkdue_1,...,\unkuno_1\unkdue_{\overline{n}},...,\unkuno_{\overline{m}}\unkdue_1,...,\unkuno_{\overline{m}}\unkdue_{\overline{n}}\right)\qquad   \forall \left(\unkuno,\unkdue\right)\in \mathbb{R}^{m}\times \mathbb{R}^{n}\text{.}
\end{equation*}
If $m=1$ or $n=1$ and there is no risk of confusion then we drop both $m$ and $n$ from the notation, with an abuse of language.
\item Fix $i,m\in \mathbb{N}$ with $i\leq m$. Then:
\begin{enumerate}
\item $\sectn[m,i]:	\mathbb{R}^{m+1}\rightarrow \mathbb{R}^m$ is the function defined by setting
\begin{equation*}
\left(\sectn[m,i]\left(\unkuno\right)\right)_{\mathsf{m}} =
\left\{
\begin{array}{ll}
\unkuno_{\mathsf{m}} &\text{if}\; {\mathsf{m}} <i\vspace{4pt}\\
0&  \text{if}\; {\mathsf{m}} = i\vspace{4pt}\\
\unkuno_{{\mathsf{m}}+1}& \text{if}\; {\mathsf{m}}>i
\end{array}  \right.\qquad \forall \unkuno \in \mathbb{R}^{m+1}\text{,}\;\;{\mathsf{m}}\in\{1,...,m\}\text{;}
\end{equation*}
\item $\proje{m}{i}:	\mathbb{R}^{m+1}\rightarrow \mathbb{R}^m$ is the function defined by setting
\begin{equation*}
\left(\proje{m}{i}\left(\unkuno\right)\right)_{\mathsf{m}}=
\left\{
\begin{array}{ll}
\unkuno_{\mathsf{m}} &\text{if}\; {\mathsf{m}} <i \vspace{4pt}\\
\unkuno_{{\mathsf{m}}+1}& \text{if}\; {\mathsf{m}}\geq i
\end{array}  \right.\qquad \forall \unkuno \in \mathbb{R}^{m+1}\text{,}\;\;{\mathsf{m}}\in\{1,...,m\}\text{;}
\end{equation*}
\item $\coord{m}{i}:	\mathbb{R}^m\rightarrow \mathbb{R}$ is the function defined by setting
\begin{equation*}
\coord{m}{i}\left(\unkuno\right)=\unkuno_i \qquad \forall \unkuno \in \mathbb{R}^m\text{,}
\end{equation*}
in case 
$\setsymuno_{\mathsf{m}}=\mathbb{R}$ for any ${\mathsf{m}}\in\{1,...,m\}$ function $\coord{m}{i}$ coincides with function $\proj<\left\{\setsymuno_{\mathsf{m}}\right\}_{\mathsf{m}=1}^m<>i>$
 introduced in Notation \ref{ins}-[10].   
\end{enumerate}
\item Fix $m\in \mathbb{N}_0$, $t \in\mathbb{R}^m$.\newline
Then $\trasl[t]:\mathbb{R}^m\rightarrow \mathbb{R}^m$ is the function defined by setting 
\begin{equation*}
\trasl[t]\left(s\right)=s+t\qquad\forall s \in \mathbb{R}^m\text{.}
\end{equation*}
\end{enumerate}
\end{notation}

\begin{notation}\label{difrealfunc}We establish the notation relative to differentiable real functions. We refer to Notation \ref{gentopnot}-[8].
\begin{enumerate}
\item Fix $m,n\in \mathbb{N}_{0}$, $l\in \mathbb{%
N}_{0}\cup \left\{+\infty \right\}$, $\intuno\subseteqdentro\mathbb{R}^m$.\newline
If $\intuno = \udenset$ then we set $\Cksp{l}(\intuno)+\mathbb{R}^{n}+=\{\udenset\}$.\newline
If $m=0$ and $\intuno \neq \udenset$ then we set $\Cksp{l}(\intuno)+\mathbb{R}^{n}+=\mathbb{R}^{n}$.\newline
If $m\geq 1$ and $\intuno \neq \udenset$ then we denote by $\Cksp{l}(\intuno)+\mathbb{R}^{n}+$ the set of all everywhere defined functions $\intuno\rightarrow \mathbb{R}^{n}$ with continuous differential up to degree $l$. We emphasize that $\Cksp{l}(\intuno)+\mathbb{R}^0+=\{\cost<\intuno<>\mathbb{R}^0>+0+\}$.
\item Fix $m,n\in \mathbb{N}_{0}$, $l\in \mathbb{%
N}_{0}\cup \left\{+\infty \right\}$, $\intuno \subseteqdentro \mathbb{R}^m$.\newline
If $\intuno = \udenset$ then $\Cksp{l}(\intuno)+\mathbb{R}^{n}+=\{\emptysetfun\}$ is endowed with the trivial topology.\newline
If $\intuno \neq \udenset$ then $\Cksp{l}(\intuno)+\mathbb{R}^{n}+$ is endowed with topology defined in \cite{RUD1} 1.44, 1.46.\newline 
We emphasize that $\Cksp{l}(\intuno)+\mathbb{R}^{n}+$ is path connected and complete. 
\item Fix $i,j,k, m\in \mathbb{N}$, $n \in \mathbb{N}_{0}$, $l \in \mathbb{N}_0\cup \left\{+\infty \right\}$, $\vecuno\in \mathbb{R}^m$, $\intuno \subseteqdentro \mathbb{R}^m$, $f\in \Cksp{l}(\intuno)+\mathbb{R}^n+$, $g\in \Cksp{l}(\intuno)+\mathbb{R}^n+$. Refer to Notations \ref{int}, \ref{ins}-[5], \ref{realvec}-[2; 6], \ref{realfunc}-[5(a)].\newline
Assume that both the following conditions are fulfilled: $i\leq \overline{n}$; $j\leq m$.\newline
We define set function $\smint_j f:\domint\left[\intuno,j\right]\rightarrow\mathbb{R}^n$ by setting: 
\begin{multline*}
\left(\smint_j f\right)_i\left(\unkuno\right) =
 \int_{\unkuno_j}^{\unkuno_{j+1}} \!f_i \left(\sectn\left[m,j\right]\left(\unkuno\right)+\intvar \ef_j\right) d\intvar \\
\hspace{29pt} \forall \unkuno \in \domint\left[\intuno,j\right]\quad \forall i \in \left\{1,...,\overline{n}\right\}\text{.}
\end{multline*}
It is well known that $\smint_j f\in \Cksp{l+1}(\domint\left[\intuno,j\right])+\mathbb{R}^n+$.\newline 
If $l\geq 1$ then we define set function $\smder_{\vecuno}f:\intuno\rightarrow\mathbb{R}^n$ by setting
\begin{equation*}
\left(\smder_{\vecuno}f\right)_{i}\left(\unkuno\right)=
\frac{\partial f_{i}}{\partial \vecuno}\left(\unkuno\right)\quad \forall\unkuno\in\intuno\quad \forall i \in \left\{1,...,\overline{n}\right\}\text{.}
\end{equation*}
It is well known that $\smder_{\vecuno}f\in \Cksp{l-1}(\intuno)+\mathbb{R}^n+$.\newline
If $\vecuno=\ef_j$ then we denote $\smder_{\vecuno}$ by $\smder_{j}$.\newline
For any fixed $\unkuno\in \intuno$ we denote by $\symjacsm<1<>f>[\unkuno]$ the jacobian of $f$ at $\unkuno$. It is well known that $\symjacsm<1<>f>[\unkuno]$ depends on $f$ locally at $\unkuno$, then we write $\gczsfx[f]$ in place of $f$ where $\gczsfx[f]$ is the germ represented by $f$ at $\unkuno$. When $\unkuno$ is fixed we drop it from the notation.\newline
The $\mathbb{R}$-vector space $\symalgtvsm<k<>m>$ of all derivations of degree $k$ of real valued smooth functions with $m$ unknowns, $m$-dimensional derivations of degree $k$ for short, is identified with the $\mathbb{R}$-vector space $\mathbb{R}^{\frac{(k+m-1)!}{k!(m-1)!}}$ by associating monomials $\underset{\mathsf{m}=1}{\overset{m}{\prod}}\left(\smder_{\mathsf{m}}\right)^{i_{\mathsf{m}}}$ for any $\left(i_1,...,i_m\right)\in \left(\mathbb{N}_0\right)^m$ with $\underset{\mathsf{m}=1}{\overset{m}{\sum}}{i_{\mathsf{m}}}=k$, ordered by lexicographical criterion, to vectors belonging to the standard orthonormal Euclidean frame of $\mathbb{R}^{\frac{(k+m-1)!}{k!(m-1)!}}$.\newline
We set $\symalgtvsm<0<>m>=\mathbb{R}$ where the identity on smooth functions is identified with $1 \in \mathbb{R}$.\newline 
The graded topological $\mathbb{R}$-vector space $\symalgtvsm<\segnvar<>m>=\left\{\symalgtvsm<k<>m>\right\}_{k \in \mathbb{N}_0}$ endowed with composition of derivatives becomes a topological commutative $\mathbb{R}$-algebra with unit. It is the symmetric algebra generated by the topological $\mathbb{R}$-vector space $\mathbb{R}^m$.\newline
For any fixed $\unkuno\in \intuno$, we denote by $\symjacsm<k<>f>[\unkuno]$ the continuous $\mathbb{R}$-linear map induced by $f$ at $\unkuno$ on the space of all $m$-dimensional derivations of degree $k$.\newline
We set $\symjacsm<0<>f>[\unkuno]=\idobj+\mathbb{R}+$.\newline
The graded continuous $\mathbb{R}$-linear function $\symjacsm<\segnvar<>f>[\unkuno]=\left\{\symjacsm<k<>f>[\unkuno]\right\}_{k \in \mathbb{N}_0}$ is a continuous function between topological commutative $\mathbb{R}$-algebras. It is the arrow between symmetric algebras generated by the continuous $\mathbb{R}$-linear function $\symjacsm<1<>f>[\unkuno]$.\newline
Functions $\symjacsm<k<>f>[\unkuno]$ depends on $f$ locally at $\unkuno$, then we write $\gczsfx[f]$ in place of $f$ where $\gczsfx[f]$ is the germ represented by $f$ at $\unkuno$. When $\unkuno$ is fixed we drop it from the notation.  
\end{enumerate}
\end{notation}

\end{document}